\documentclass[12pt]{article} 

\usepackage{amsmath,amsfonts,amssymb,mathtools,amsthm,epsfig}
\usepackage{mathrsfs}
\usepackage{comment}

\usepackage[all]{xy}
\usepackage{psfrag}

\usepackage{arydshln}
\usepackage{indentfirst}
\usepackage{hyperref}
\hypersetup{colorlinks=true,urlcolor=blue,linkcolor=blue,citecolor=blue}

\usepackage{enumitem}

\DeclarePairedDelimiter\floor{\lfloor}{\rfloor}

\newcommand{\op}[1]{\operatorname{#1}}

\theoremstyle{plain}

\newtheorem{global-theorem}{Theorem}
\newtheorem{theorem}{Theorem}[section]
\newtheorem{lemma}[theorem]{Lemma}

\newtheorem{corollary}[theorem]{Corollary}
\newtheorem{conjecture}[theorem]{Conjecture}
\newtheorem{definition}[theorem]{Definition}
\newtheorem{proposition}[theorem]{Proposition}

\newtheorem{claim}[theorem]{Claim}

\theoremstyle{definition}

\newtheorem{hypothesis}[theorem]{Hypothesis}
\newtheorem{remark}[theorem]{Remark}

\newtheorem{problem}[theorem]{Problem}

\newtheorem{subclaim}[subsubsection]{Claim}
\newtheorem{subcorollary}[subsubsection]{Corollary}

\DeclareFontFamily{U}{rsf}{}
\DeclareFontShape{U}{rsf}{m}{n}{
  <5> <6> rsfs5 <7> <8> <9> rsfs7 <10->  rsfs10}{}
\DeclareMathAlphabet{\mathscr}{U}{rsf}{m}{n}
\newcommand{\mycal}[1]{\mathscr{#1}}

\newcommand{\lesone}[6]{
\xymatrix{     
 0 \ar[r] & {#1} \ar[r]  &  {#2} \ar[r]  &  {#3} 
\ar@{->}`r/10pt[d] `[l] `^dl[dlll]  `^dr/10pt[dll]    [dll] \\
 &  {#4} \ar[r] & {#5} \ar[r] & {#6} \ar[r] & 0 }
}

\renewenvironment{description}[1][0pt]
  {\list{}{\labelwidth=0pt \leftmargin=#1
   }}
  {\endlist}

\DeclareMathAlphabet{\mathsfit}{T1}{\sfdefault}{\mddefault}{\sldefault}
\SetMathAlphabet{\mathsfit}{bold}{T1}{\sfdefault}{\bfdefault}{\sldefault}

\newcommand{\cc}{{\mathbb C}}
\newcommand{\pp}{{\mathbb P}}
\newcommand{\rr}{{\mathbb R}}
\newcommand{\qq}{{\mathbb Q}}
\newcommand{\zz}{{\mathbb Z}}

\newcommand{\Gm}{{\mathbb G}_{\rm m}}
\newcommand{\aaaa}{{\mathbb A}}
\newcommand{\ttt}{{\mathbb T}}

\newcommand{\Ff}{{\mathcal F}}
\newcommand{\Oo}{{\mathcal O}}
\newcommand{\Uu}{{\mathcal U}}

\newcommand{\Vv}{{\mathcal V}}

\newcommand{\Ww}{{\mathcal W}}

\newcommand{\Bb}{{\mathcal B}}

\newcommand{\Gg}{{\mathcal G}}
\newcommand{\Ll}{{\mathcal L}}
\newcommand{\Ee}{{\mathcal E}}

\newcommand{\srF}{{\mathscr F}}
\newcommand{\srD}{{\mathscr D}}

\newcommand{\srX}{{\mathscr X}}

\newcommand{\srL}{{\mathscr L}}

\newcommand{\Dol}{{\mathsf{Dol}^{\op{par}}_{L^{2}}}}

\newcommand{\sq}{{\mathsf{sq}}}

\newcommand{\Chat}{{\widehat{C} }}
\newcommand{\pihat}{\widehat{\pi}}
\newcommand{\sqhat}{\widehat{\sq}}

\newcommand{\Cbar}{{\overline{C} }}
\newcommand{\Wbar}{{\overline{W}}}

\newcommand{\Ctilde}{{\widetilde{C} }}

\newcommand{\Wob}{{\rm Wob}}

\newcommand{\Sym}{{\rm Sym}}
\newcommand{\he}{\text{\bfseries\sf h}}

\newcommand{\ba}{\mathsf{a}}
\newcommand{\fD}{\boldsymbol{\mathsf{f}}}

\newcommand{\Residual}{{\operatorname{{\sf Residual}}}}
\newcommand{\Residualhat}{\widehat{\Residual}}

\newcommand{\bigHecke}{\mathcal{H}}
\newcommand{\bigHeckebar}{{\overline{\mathcal{H}}}}
\newcommand{\bigHeckehat}{{\widehat{\mathcal{H}}}}
\newcommand{\Prym}{{\mathcal{P}}}

\newcommand{\blo}{\boldsymbol{\varepsilon}}
\newcommand{\hyp}{\mathsf{h}}
\newcommand{\pr}{\op{pr}}
\newcommand{\pw}{\boldsymbol{\mathsf{p}}}
\newcommand{\Kum}{\mathsf{Kum}}
\newcommand{\hi}{\boldsymbol{\iota}}
\newcommand{\brb}{\boldsymbol{\mathsf{b}}}
\newcommand{\brx}{\boldsymbol{\mathsf{x}}}
\newcommand{\bry}{\boldsymbol{\mathsf{y}}}
\newcommand{\brz}{\boldsymbol{\mathsf{z}}}
\newcommand{\Higgs}{\mathsf{Higgs}}
\newcommand{\hit}{\boldsymbol{h}}
\newcommand{\trans}{\boldsymbol{t}}
\newcommand{\emb}{\boldsymbol{\imath}}
\newcommand{\embj}{\boldsymbol{\jmath}}
\newcommand{\pwtilde}{\widetilde{\pw}}

\newcommand{\Lprym}{\mathfrak{L}}
\newcommand{\LY}{\mycal{L}}
\newcommand{\thetaprym}{\boldsymbol{\xi}}
\newcommand{\modalpha}{\boldsymbol{\alpha}}
\newcommand{\ExY}{\boldsymbol{\mathit{E}}}
\newcommand{\FxY}{\boldsymbol{\mathit{F}}}
\newcommand{\RxY}{\boldsymbol{\mathit{R}}}
\newcommand{\bigExc}{\boldsymbol{\mathit{Exc}}}
\newcommand{\btheta}{\boldsymbol{\theta}}
\newcommand{\bTheta}{\boldsymbol{\Theta}}
\newcommand{\thetatilde}{\widetilde{\btheta}}
\newcommand{\Thetatilde}{\widetilde{\bTheta}}
\newcommand{\qfrak}{\boldsymbol{\mathfrak{q}}}
\newcommand{\specN}{\boldsymbol{\mathit{N}}}
\newcommand{\pzo}{d}
\newcommand{\qzo}{b}
\newcommand{\add}{\mathbf{add}}
\newcommand{\summ}{\mathsf{sum}}
\newcommand{\diff}{\mathsf{diff}}
\newcommand{\trope}{\mathsf{Trope}}
\newcommand{\conic}{\mathfrak{C}}

\newcommand{\Fix}{\mathfrak{Q}}
\newcommand{\cstar}{\cc^{\times}}
\newcommand{\fullf}{\boldsymbol{\mathfrak{f}}}
\newcommand{\unif}{\boldsymbol{\mathsf{x}}}
\newcommand{\flow}{\mathbf{fl}}

\newcommand{\ql}{\mathsf{ql}}
\newcommand{\bT}{\boldsymbol{\mathcal{T}}}
\newcommand{\bN}{\boldsymbol{\mathcal{N}}}
\newcommand{\bPP}{\boldsymbol{\mathfrak{p}}}
\newcommand{\bQQ}{\boldsymbol{\mathfrak{q}}}
\newcommand{\lphi}{\boldsymbol{\mathfrak{i}}}
\newcommand{\lalpha}{\boldsymbol{\mathfrak{j}}}
\newcommand{\lA}{\mathbf{A}}
\newcommand{\lO}{\boldsymbol{\mathsf{O}}}
\newcommand{\ldelta}{\boldsymbol{\mathfrak{U}}}
\newcommand{\emm}{\mathsfit{m}}
\newcommand{\ellfrak}{\boldsymbol{\mathfrak{l}}}
\newcommand{\lquad}{\mathsfit{quad}}
\newcommand{\lC}{\mathsfit{Conic}}
\newcommand{\lv}{\mathsfit{v}}
\newcommand{\rul}{\boldsymbol{\mathsf{rul}}}
\newcommand{\bnu}{\boldsymbol{\nu}}
\newcommand{\cJ}{\boldsymbol{\mathcal{J}}}
\newcommand{\cKum}{\boldsymbol{\mathcal{K}um}}
\newcommand{\eKK}{\mathsfit{K}}
\newcommand{\ekk}{\mathsfit{k}}
\newcommand{\blowb}{\mathsfit{b}}
\newcommand{\bii}{\mathsfit{i}}
\newcommand{\bLL}{\mathsfit{L}}

\newcommand{\KumKthree}{\boldsymbol{\mathit{K}}}
\newcommand{\Heisen}{\mathsf{Heisen}}

\def\punkt{\refstepcounter{subsubsection}
           \noindent{\bfseries \thesubsubsection.\ }}

\title{Twistor Hecke eigensheaves in 
genus $2$}
\author{R. Donagi, T. Pantev and C. Simpson}
\date{}

\begin{document}
\maketitle

\begin{abstract}
  Following the strategy outlined in \cite{DP1,DonagiPantev} for
  bundles of rank $2$ on a smooth projective curve of genus $2$, we
  construct flat connections over the moduli of stable bundles, with
  singularities along the wobbly locus. We verify that the
  associated $\srD$-modules are Hecke eigensheaves. The local systems
  are constructed by the nonabelian Hodge correspondence from Higgs
  bundles. The spectral varieties of the Higgs bundles are the Hitchin
  fibers corresponding to the Hecke eigenvalues.
\end{abstract}

\setcounter{tocdepth}{2}
\tableofcontents


\section{Introduction}
\label{chapter-intro}

\subsection{Geometric Langlands}
\label{sec-gl}

Suppose $C$ is a smooth projective curve over the complex numbers. Let
$\srX$ be a moduli stack of vector bundles or principal $G$-bundles
for a reductive group $G$.  The open substack of semistable bundles
has a good coarse moduli space $X$.  The {\em geometric Langlands
  program} predicts the existence of certain perverse sheaves on
$\srX$. These are characterized by a property called the {\em Hecke
  eigensheaf property} with respect to a local system $\Lambda$ on
$C$, and it is predicted that to a given $\Lambda$ there is a unique
Hecke eigensheaf. In turn, a perverse sheaf on $\srX$ leads to a local
system on an open subset of $X$. So, we have a construction going from
a local system $\Lambda$ on the curve $C$ to a local system on an open
subset of $X$.

Our purpose in this paper, pursuing the program of the first two
authors, is to look at this construction from the viewpoint of the
Hitchin equations.  These equations lead to what is variously known as
the {\em nonabelian Hodge correspondence} or {\em Kobayashi-Hitchin
  correspondence} relating local systems to Hitchin pairs or {\em
  Higgs bundles}. These geometric objects showing up on the other side
of the correspondence are in many respects more tangible than local
systems, in particular they are associated to spectral data including
a spectral cover which is typically a ramified covering of the
base. We would like to understand the relationship between the
spectral data of the Higgs bundle corresponding to the input local
system $\Lambda$, and the spectral data corresponding to the resulting
output local system over an open subset of $X$.

Before getting to a more detailed look at what we do, let's recall that
the geometric Langlands program originated as a geometrization of the
{\em Langlands program}, where automorphic functions are categorified to
data of perverse or constructible sheaves. In positive characteristic
for $\ell$-adic sheaves, this categorification is based on the
function-sheaf correspondence, where a constructible sheaf defines a
function on the set of finite field-valued points by associating to
each point the trace of its Frobenius action on the stalk of the
sheaf. Due to several authors
\cite{Drinfeld-ICM,Drinfeld,Drinfeld-LOMI,LaumonDuke,
  Laumon95,BD-Hitchin,Gaitsgory,GaitsgoryBourbaki}, the geometric
Langlands program was then carried from positive characteristic to
characteristic zero, where we no longer have a function-sheaf
correspondence, but the resulting statement at the sheaf level still
makes sense. Somewhat surprisingly, over $\mathbb{C}$, or more
generally over local fields of characteristic zero, there is an
analytic function-theoretic version of the Langlands correspondence
originally envisioned in \cite{langlands,frenkelC} and further
developed in \cite{efk1}. The recent works
\cite{efk2,efk3,efk4,bk-survey} have made exciting advances in
understanding and proving instances of this analytic statement. The
same works also make connections with the categorified geometric
version of the Langlands conjecture and we expect them to ultimately
have a direct relation with our Hodge theoretic approach. Exploring
this analytic picture and the expected relation is a very interesting
question which, unfortunately,  is beyond the scope of the present paper.

A main player in both the function-theoretic and the categorified
Langlands conjecture over $\mathbb{C}$  is the algebra of {\em Hecke
  correspondences}. These act on the moduli stack of principal
$G$-bundles over a Riemann surface $C$. This is most easily understood
for $G=GL_2$. Viewed as a multivalued function, the Hecke
correspondence at a point $t\in C$ takes a rank $2$ bundle $E$ to the
sum of its {\em Hecke transforms} at $t$: these are the elementary
transforms $E'$ fitting into exact sequences
$$
0 \rightarrow E' \rightarrow E \rightarrow \cc_t \rightarrow 0
$$
where $\cc_t$ is the skyscraper sheaf of length $1$ at $t$. The set
of Hecke transforms of $E$ is parametrized by the set of rank $1$
quotients $E_t \rightarrow \cc_t$, which is a $\pp^1$.  Thus, a
bundle $E$ is sent to a formal sum of a $\pp^1$'s-worth of new bundles
$E'$. In the function world, one just takes the sum over the discrete
set of points of $\pp^1$.  In the sheaf-theoretical viewpoint, the
formal sum is replaced by the cohomology of a sheaf over $\pp^1$.

The correspondence depends on the choice of point\footnote{This needs
to be modified in case we look at bundles of fixed determinant as we
do in the present paper---the parametrizing data is then a point in a
covering $\Cbar$ of $C$.}  $t\in C$. Thus, letting $X$ denote the moduli
space of bundles, we obtain the Hecke correspondence
 \[
\xymatrix@M+0.25pc{
  & \bigHecke \ar[dl]_-{p}  \ar[dr]^-{q}  & \\
X & & X\times C .
} 
\]
If $\mycal{F}_B$ is a perverse sheaf on $X$, its {\em Hecke transform}
is the perverse sheaf $Rq_{*}(p^{*} \mycal{F}_B)$ on $X\times
C$.  The value of the Hecke transform $\bigHecke (t)$ at a point $t$
is the restriction of this on $X\times \{ t\}$.

In a formal viewpoint, the Hecke operations at different points
$t,t'\in C$ commute. The classical theory therefore views the whole
algebra of Hecke operations as an algebra of commuting operators, and
it becomes natural to look for a common diagonalization of these
operators.

In the sheaf-theoretical viewpoint, it means that we are looking for
{\em Hecke eigensheaves} $\mycal{F}_B$ on $X$, corresponding to {\em Hecke
  eigenvalues} $\Lambda_{B}$ that are perverse sheaves on $C$. The
eigenvalue equation, saying in naive terms that the Hecke operation
$\bigHecke (t)$ multiplies $\mycal{F}_B$ by the eigenvalue that
depends on $t$, is written as
$$
Rq_{*}(p^{*} \mycal{F}_B ) \cong \mycal{F}_B \boxtimes \Lambda_{B} .
$$ One of the main tasks of the geometric Langlands program is to
construct, for a given eigenvalue $\Lambda_{B}$ (that's really the
sheaf-theoretical version of the notion of ``collection of eigenvalues
one for each Hecke operation $\bigHecke (t)$''), a Hecke eigensheaf
$\mycal{F}_B$ corresponding to this eigenvalue. The ``de Rham''
version of the geometric Langlands conjecture a Hecke eigenvalue is a
flat bundle $\Lambda_{dR}$ on $C$ while a corresponding Hecke
eigensheaf is a $\srD$-module $\mycal{F}_{dR}$ on $X$ and these fit better
with parametrizations for the next paragraph.

The translation of Langlands' conjecture from the automorphic world to
the geometric context predicts that there will be a unique such
eigensheaf for each eigenvalue $\Lambda_{dR}$. Furthermore, it states
that as a function of $\Lambda_{dR}$, the eigensheaf varies in a
coherent way. Namely, rather than starting with an individual flat
bundle $\Lambda_{dR}$ on $C$, we could start with a combination of
these in the form of a coherent sheaf on the moduli space $Flat$ of
$\Lambda_{dR}$'s. The (naive) geometric Langlands conjecture predicts
that to such a coherent sheaf should be associated a unique
$\srD$-module $\mycal{F}_{dR}$ on the moduli $Bun$ bundles on $C$, and
that this correspondence should set up a duality between complexes of
coherent sheaves on $Flat$ and complexes of $\srD$-modules on $Bun$.

The notion of {\em Langlands dual group} enters here: if we use a
reductive group $G$ to speak of the moduli stack $Bun(G)$ of principal
$G$-bundles on $C$, then we need to take the Langlands dual group
${}^LG$ and look at coherent sheaves on the moduli $Flat({}^LG)$ of
flat ${}^LG$ bundles on $C$.

Much important progress has been made on establishing the geometric
Langlands correspondence.  Drinfeld was the first to make a
construction of Hecke eigensheaves for the group $GL_2$.  His article
\cite{Drinfeld} changed over from the sheaf viewpoint to the function
viewpoint somewhere in the middle, so it really constructs Hecke
eigenfunctions. Laumon   \cite{Laumon95,LaumonDuke} formalized and
generalized this construction, yielding a solution of the geometric
Langlands problem for $GL_2$, and a conjectural framework for
$GL_n$. Gaitsgory gave an alternative proof in his thesis 
\cite{GaitsgoryThesis}. Then Lafforgue   \cite{LafforgueChtoucas}
proved it for $GL_n$ in the number-theoretical context, and
Frenkel-Gaitsgory-Vilonen proved the geometric version for $GL_n$
\cite{FrenkelGaitsgoryVilonen}. In a first noncompact case, Arinkin
\cite{Arinkin}   treated the case of parabolic bundles on $\pp^1$
with $4$ singular points.

We will be looking in detail at the non-abelian Hodge theory approach
to Drinfeld's original construction in Section \ref{chapter-drinfeld},
and will show how its Dolbeault version can be understandood via
abelianization and the spectral cover construction. The constructions
for higher rank can probably also be recast in terms that would be
more familiar to geometers, although this is bound to contain a
certain level of complication and we don't attempt it here.

It has remained, at least until fairly recently, elusive how one would
attack the problem for general groups $G$. One may take note of recent
progress such as
\cite{Beraldo-Chen,Beraldo,Beraldo-spectral,Beraldo-Whittaker,
  RozenblyumHecke,FaegermanRaskin,Faegerman-qt}
in the de Rham setting and
\cite{BenZviNadlerBetti,agkrrv1,agkrrv2,agkrrv3} in the Betti setting.

It has also been understood, in the preceding period, that the simple
description we have tried to approximate above is not adequate to
describe a duality between the two sides of the geometric Langlands
correspondence, and indeed that the two categories have different
properties making it so that they couldn't be equivalent. Therefore,
the categories need to be modified. The work of Gaitsgory and
Rozenblyum on Ind-coherent sheaves
\cite{Gaitsgory-indcoh,GaitsgoryRozenblyum} and of Arinkin, Gaitsgory,
and collaborators on nilpotent singular supports
\cite{ArinkinGaitsgory-nilp,agkrrv2,agkrrv3} aims to solve these
problems. They mostly have to do with parts of the categories that are
supported on or close to locations in the moduli spaces where various
kinds of singularities occur. Therefore, in looking for a global
understanding of at least some part of the correspondence, we will
ignore these subtleties.

Recently, Gaitsgory and a group of co-authors have announced a full proof
of the geometric Langlands conjecture for all groups, 
with first drafts available \cite{GLC}. 

Something that has not, in our view, been sufficiently emphasized in
previous works in this area is the fact that the geometric Langlands
program predicts something very specific about the topology and
geometry of moduli spaces of vector bundles on curves.  This question
was however raised in Sawin's MathOverFlow post \cite{Sawin}.  For one
thing, the moduli stack is a pretty wild beast, being only locally of
finite type. However, it contains an open substack that is the moduli
of semistable bundles. This substack is in turn close to being a
projective variety, in that its coarse moduli space is the projective
moduli space $X$ of $S$-equivalence classes of semistable bundles,
whose points are in $1:1$ correspondence with the polystable bundles.

Therefore, for a perverse sheaf on $Bun(G)$ we will have an open
subset $X^{\rm vs}$ on which the perverse sheaf is a locally constant
sheaf, i.e. a representation of $\pi_1(X^{\rm vs})$. It was known
early on by Laumon that if the perverse sheaf in question is a Hecke
eigensheaf, the corresponding open subset $X^{\rm vs}$ has an explicit
description as the moduli space of {\em very stable} bundles: a very
stable bundle is one that does not admit a nonzero nilpotent Higgs
field.

So, and in spite of the oversimplifications in the above presentation,
the abstract geometric Langlands correspondence predicts something
concrete and easily understandable: that to a perverse sheaf
$\Lambda$, let us say itself a local system on $C$, there should be
naturally associated a local system over $X^{\rm vs}$. It is this
construction that we would like to study in the present paper.

There are several motivations for the viewpoint we adapt here. The
first was a small detail immediately noticed by the first author, in
the Manin volume of the {\em Duke Mathematical Journal}.  Hitchin's
article on the moduli space of Higgs bundles \cite{Hitchin-spectral}
and Laumon's article on the geometric Langlands correspondence
\cite{LaumonDuke} were both in this volume. In Remarque 5.5.2, Laumon
states that Deligne, in unpublished communication, calculated the
multiplicity of the zero-section in the characteristic cycle of the
$\srD$-module: this is the same as the rank of the local system over
$X^{\rm vs}$. The answer (for $SL_2$-bundles on a curve of genus $g$)
was $2^{3g-3}$.  This number was the same as the degree of the map
from a general fiber of the Hitchin system to the moduli space of
bundles.

The general fiber of the Hitchin system is a subvariety in what is,
basically, the cotangent bundle of $X$ (the necessary modifications to
that statement will be the subject of discussion later).

In light of the nonabelian Hodge correspondence, started in Hitchin's
paper \cite{Hitchin-selfd} and developed in the highest level of
generality by the third author \cite{SimpsonHiggs} and Mochizuki
\cite{Mochizuki-kh1,Mochizuki-kh2} some time later, it would look
natural to think of the general fiber of the Hitchin system, a
subvariety of $T^{\vee}X$ of degree $2^{3g-3}$ over the base $X$, as a
good candidate for being the spectral variety of the geometric
Langlands local system of rank $2^{3g-3}$.  This idea became the
conjecture of the first two authors \cite{DP1}, and is the essence of
what we will be trying to do here.

Some early ideas in this direction were contained in a letter from
Hausel to Hitchin (unpublished), and in Faltings' talk at Deligne's 61th
birthday conference.

A next element of motivation explains more precisely how this should
be organized. This is known as {\em electric-magnetic duality} in the
work of Kapustin and Witten \cite{KapustinWitten} or the {\em
  classical limit of geometric Langlands} \cite{DP2012,DP1}.  Hausel
and Thaddeus \cite{HauselThaddeus,Hausel} view it as a form of mirror
symmetry.  

The Hitchin moduli spaces of principal Higgs bundles for the groups
$G$ and ${}^LG$ fit into a diagram
$$
\xymatrix@M+0.3pc@C-1pc@R-1pc{
  \Higgs(G) \ar[dr]_-{\hit} & & \Higgs({}^{L}G) \ar[dl]^-{{}^{L}\hit} \\
 & \Bb & 
}
$$
where:

\begin{theorem} 
\label{electromagnetic}
The base affine spaces of the  Hitchin fibration ${\bf B}$ for 
$G$ and ${}^LG$ are naturally isomorphic. 
Furthermore, the two Hitchin maps 
$\hit$ and ${}^L\hit$ are generically dual SYZ-type torus fibrations.
\end{theorem}

We refer to \cite{HauselThaddeus,KapustinWitten,DP2012} for the proof.

The prediction of  the de Rham version of the geometric Langlands
correspondence may be stated as in the following several paragraphs.
Let $X_G$ be the coarse
moduli space of semistable $G$-bundles on $C$.  Then given an
${}^LG$-flat bundle $\Lambda_{dR}$ on $C$, there should be a
$\srD_X$-module $\mycal{F}_{dR}$ on $X_G$ that is (reasonably
approximately\footnote{ The statement for the coarse moduli space
$X_G$ is not exactly the same as the stack-theoretical statement,
because the extension from a local system on $X^{\rm vs}_G$ to a
perverse sheaf on $X_G$ is a different object than the extension to a
perverse sheaf on $Bun_G$; we can expect some modifications to the
Hecke eigensheaf property. Such a distinction does not seem to
intervene at the level discussed in the present paper.  })  a Hecke
eigensheaf for $\Lambda_{dR}$.

Through the non-abelian Hodge theorem the eigenvalue ${}^{L}G$-flat
bundle $\Lambda_{dR}$ corresponds to a point $\Lambda_{Dol}$ of
$\Higgs ({}^LG)$, in the fiber ${}^L\hit^{-1}(\mathsf{b})$ over a
point $\mathsf{b}\in \Bb$ in the Hitchin base that we will suppose to
be general.

The fiber $\Prym := \hit^{-1}(\mathsf{b})\subset \Higgs $ is the dual
torus of ${}^L\hit^{-1}(\mathsf{b})$, so the Fourier-Mukai transform
for the dual torus fibrations says that the point $\Lambda_{Dol}$ (and
hence the point $\Lambda_{dR}$) corresponds to a line bundle $\Lprym$
over $\Prym$.

The moduli stack of $G$-Higgs bundles is isomorphic to the cotangent
stack of $Bun(G)$. For moduli spaces, this is no longer true globally,
but it remains true birationally.  The set of very stable points
$X^{\rm vs}$ of the coarse moduli space of $G$-bundles is smooth, and
its cotangent space $T^{\vee}(X^{\rm vs})$ is the subspace $\Higgs
^{\rm vs}$ of $G$-Higgs bundles whose underlying  bundle is very
stable.

The first two authors propose to use $\Prym$, viewed in a birational
sense, after blowing up, as a subvariety $Y$ of the (logarithmic)
cotangent bundle of $X$, to be a spectral cover of $X$.  And to use
$\Lprym$ as input spectral datum to construct a logarithmic parabolic
Higgs sheaf $\mycal{F}_{Dol}$ on $X$ that, under the non-abelian Hodge
correspondenc \cite{Mochizuki-kh1,Mochizuki-kh2}, should correspond to
the Hecke eigensheaf $\srD_X$-module $\mycal{F}_{dR}$.

The divisor of singularities, complement of Laumon's open set $X^{\rm
  vs}$ of very stable points, was termed
the {\em wobbly divisor}
$\Wob \subset X$ in \cite{DP1}. Points of $\Wob$ are semistable
bundles such that $H^0(\op{ad}(E)\otimes \omega_C)$ has a nonzero
nilpotent element.\footnote{It turns out, looking into the 
study of the classical {\em quadric line complex} which is our moduli
space $X_1$, that the wobbly locus in that case was known 
as the union of what are called {\em special lines} \cite[page 792]{GH}.
}

Under this birational transformation thhe degree $0$ line bundle
$\Lprym$ is modified to a line bundle $\LY$ on $Y$, and parabolic
structure is supposed to be added, along divisors lying over the
wobbly divisor.

Here is a first, albeit incomplete, formulation of the conjecture of
the first two authors.

\begin{conjecture}[\cite{DP1,DonagiPantev}]
\label{conj:DP}
Suppose $\Lambda_{dR}$ is an ${}^LG$ local system on $C$ corresponding
to a point $\Lambda_{Dol}$ in the Hitchin fiber of $\Higgs ({}^LG)$
over a general point $b\in \Bb$ in the Hitchin base, and let
$\Lprym$ be the dual line bundle on the Hitchin fiber $\Prym$ of
$\Higgs (G)$ over the same point $b$.  Then the Hecke eigensheaf
$\mycal{F}_{dR}$ for eigenvalue $\Lambda_{dR}$ corresponds, via Mochizuki's
Kobayashi-Hitchin correspondence, to a parabolic logarithmic Higgs
bundle $\mycal{F}_{Dol , \bullet}$ on $(X_G,\Wob )$. The spectral data
over the open subset $X^{\rm vs}_G=X_{G}-\Wob$ consists of the pullback
open subset $\Prym ^{\rm vs}$ considered as a subvariety of $T^{\vee}
(X^{\rm vs}_G)$, with the restriction of $\Lprym$ as spectral line
bundle.
\end{conjecture}

The nonabelian Hodge correspondence is fundamentally global, involving
the solution of Hermite-Yang-Mills-Higgs equations minimizing a
Yang-Mills functional over the manifold. In particular, the
specification of the Higgs bundle by its spectral data over an open
subset $X_{G}^{\rm vs}$ does not uniquely specify how it might be extended
to $(X_{G},\Wob )$. Such an extension involves choosing an extension of
the bundle, and possibly a parabolic structure.

On the smooth parts of the  wobbly divisor $\Wob = X_{G} -
X_{G}^{\rm vs}$, one might be able to formulate a more precise conjecture
spelling out how the parabolic structure is supposed to look.  This
was done for the case of the root stack of $\pp^1$ with five singular
points in \cite{DonagiPantev}.

One of the
objectives of the study we do in the present paper is to continue 
the investigation of what type of behavior to expect for this structure. 

The wobbly divisor will, in general, have singularities that are more
complicated than normal crossings. Therefore, in order to apply the
general theory \cite{Mochizuki-kh1}, we need to blow up to resolve
those singularities. Because of the Bogomolov-Gieseker inequality, the
way to extend the parabolic structure---that will have been defined in
codimension $1$---is unique. It is determined by the condition of
minimizing $c_2^{\rm par}$ (or equivalently maximizing ${\rm
  ch}_2^{\rm par}$) subject to the constraint $c_1^{\rm par} =
0$. However, no recipe is currently known for doing this. In the
present paper, we are faced with two different situations of this
type, and we adopt two different strategies for attaining the minimum,
and showing that it is $c_2^{\rm par} = 0$ that leads to a harmonic
bundle and hence a local system.

The long-term hope is that we could understand these processes in a
uniform way for all groups $G$, obtaining a uniform construction of
Hecke eigensheaves at least for generic initial eigenvalue data
$\Lambda$.

\subsection{The case of genus \texorpdfstring{$2$}{2} and rank
  \texorpdfstring{$2$}{2}}
\label{sec-g2r2}

In this paper, we take a much less lofty goal: to understand how this
works in the case where $G$ is $PSL_2$ and the curve $C$ has genus
$2$. At the end of the paper, we will prove a comparison with Drinfeld's
construction. This proof in Chapter \ref{chapter-drinfeld} will
actually prove Conjecture \ref{conj:DP} for $SL_2$ on curves of
arbitrary genus, but it does not immediately give the further
information about parabolic structures. That seems attainable and
could be the topic of future work.

Let us now look more carefully at what is to be done.  The reader may
refer to Section \ref{sec-mt} below for concise statements of the main
theorems.

Here are a few introductory notations.  The curve $C$ is a smooth
projective curve of genus $g=2$. It is hyperelliptic, with degree $2$
map $\hyp : C\rightarrow \pp^1$ and hyperelliptic involution
$\hi_C$. Out of the $6$ ramification or \emph{\bfseries Weierstrass}
points, fix one of them denoted $\pw$ to be used as a basepoint
(e.g.\ for Abel-Jacobi maps) throughout the paper.

One of the first
main inputs is the beautiful classical theory of moduli spaces 
of semistable bundles of rank $2$ and fixed determinant.

These were investigated extensively by Tyurin \cite{Tyurin}, 
Narasimhan and Ramanan
\cite{NR}, Desale and Ramanan \cite{DesaleRamanan}, 
Newstead \cite{Newstead}, and then further by Previato-van
Geemen \cite{PreviatoVanGeemen}, Beauville
\cite{Beauville-genus2and3}, Heu-Loray
\cite{Heu,HeuLoray-flat,HeuLoray-Higgs}, Pal-Pauly
\cite{PalPauly}, etc. 

Let $X$ be a moduli space of polystable vector bundles of rank $2$ on
$C$ with a fixed determinant line bundle. Up to isomorphism given by
tensoring with a line bundle, there are two possibilities depending on
the parity of the degree. The following explicit descriptions are
provided by Narasimhan and Ramanan \cite{NR}:
\begin{itemize}

\item 
$X_0$ is the moduli space of bundles of degree $0$ with determinant
  $\Oo _C$, and $X_0 \cong \pp^3$;

\item
$X_1$ is the moduli space of bundles of degree $1$ with determinant
  $\Oo _C(\pw )$, and $X_1 =Q_1\cap Q_2 \subset \pp^5$ is the
  intersection of two quadrics.
\end{itemize}

The explicit descriptions were generalized to bundles of rank $2$ over
higher genus hyperelliptic curves by Desale-Ramanan
\cite{DesaleRamanan}.

In even degree $X_0$ is only a coarse moduli space: it is known that a
universal family does not exist even over any Zariski open set.  For
bundles of odd degree, a universal family exists, and heuristically it
seems to be significantly easier explicitly to describe points of
$X_1$.

Let $\Higgs_0$ resp.\ $\Higgs _1$ denote the moduli space of Higgs
bundles of the same rank, degree and determinant.  These were studied
early on by Previato and van Geemen \cite{PreviatoVanGeemen}.

For comments that apply to both cases, and to some extent for other
curves and other groups, we will just write $X$ and $\Higgs$ for either
of the spaces.  For example we can say in that $\dim (X) = 3$ and
$\dim (\Higgs )=6$.

The \emph{\bfseries Hitchin fibration} is the map $\Higgs
\stackrel{\hit}{\longrightarrow} \Bb = \aaaa^3$ sending a Higgs bundle
to the moduli point of its associated spectral curve $\Ctilde \subset
T^{\vee}C$.  The \emph{\bfseries nilpotent cone} has $X$ as its principal
component:
$$
\hit^{-1}(0) = X \cup \; (\mbox{other components}).
$$
The other components usually correspond to the components of the
wobbly locus $\Wob$  \cite{PalPauly,PeonNieto,HauselHitchin}, in
the sense that the other components of $\hit^{-1}(0)$ intersect $X$
along components of $\Wob$. Our case of $\Higgs _0$ is an exception:
there, $\Wob_0$ has an extra component that does not correpond to an
additional component of the nilpotent cone.

Indeed, the wobbly locus of $X_0$ decomposes as 
$$
\Wob _0 = \Kum \cup  \left[\cup_{\kappa \in \mathsf{Spin}(C)}
\trope_{\kappa}\right],
$$
where $\Kum\subset \pp^3$ is the Kummer surface associated to $C$,
and $\trope_{\kappa}$ are the sixteen \emph{\bfseries trope planes},
which are naturally labeled by the set $\mathsf{Spin}(C)$ 
of theta characteristics $\kappa$ of
$C$. The trope planes do correspond to the extra components of the
nilpotent cone. But the Kummer surface does not, rather it consists
entirely of singular points in $\Higgs_0$ that could be viewed as some
kind of infinitesimally nearby components of $\hit^{-1}(0) $.

The Kummer surface itself has $16$ nodes and these form, with the $16$
trope planes, the famous \emph{\bfseries Kummer $16_6$ configuration} 
\cite{GH,BeauvilleSurfaces,Keum,Dolgachev200}.  Each plane passes
through $6$ points and each point is contained in $6$ planes.

In $X_1$, the wobbly locus is a singular surface whose normalization is
$$
\pp^1 \times \overline{C} \rightarrow \Wob _1.
$$
This map comes about in the following way: the curve $\overline{C}$, a
$16$-sheeted etale cover of $C$, sits inside $\Higgs _1$ as the second
fixed point locus of the $\cstar$ action; the downward direction
of the $\cstar$ action is forms a projectively trivial bundle over
$\overline{C}$ and the limits of downward orbits are the points of
$\Wob _1$. In particular, $\Wob _1$ does equal the locus where the
single extra component of the nilpotent cone meets $X_1$.

Several authors have recently investigated more closely the structure
of the structure of the nilpotent cone and the wobbly locus in general
situations: Bozec \cite{Bozec-nilp}, Gothen-Z\'{u}\~{n}iga-Rojas
\cite{Gothen-nilp}, Pal-Pauly \cite{PalPauly}, Pe\'on-Nieto
\cite{PeonNieto,FGOP}, Zelaci \cite{Zelaci}, Hellmann
\cite{Helmann-nilp}, Hausel-Hitchin \cite{HauselHitchin}.

The fiber of the Hitchin fibration over a general point $\mathsf{b}\in
\aaaa^3$ of the Hitchin base is a Prym variety $\Prym :=
\hit^{-1}(\mathsf{b})$.  The point $\mathsf{b}$ corresponds to a
spectral curve
$$
\pi : \widetilde{C} \rightarrow C, 
$$
which is a $2:1$ covering ramified over $4$ points.  The abelian
variety $\Prym $ is the Prym variety of line bundles on
$\widetilde{C}$ whose norm down to $C$ is the line bundle $\Oo _C(2\pw
)$ for $X_0$ or $\Oo _C(3\pw )$ for $X_1$. The shift by $\Oo _C(2\pw )
\cong \omega _C$ is due to the fact that $\pi _{*}(\Oo _{\Ctilde})
\cong \Oo _C \oplus \omega _C^{-1}$ so the determinant line bundle of
the rank $2$ vector bundle is $\omega _C^{-1}$ times the norm line
bundle.

Over the open subset of very stable points, 
$$
\xymatrix@M+0.3pc@R-0.2pc@C-0.2pc{
\Higgs ^{\rm vs} \cong  T^{*} X^{\rm vs} \ar[d]^-{\fullf} \\
 X^{\rm vs}} 
$$
The identification between $\Higgs ^{\rm vs}$ and $T^{\vee} (X^{\rm vs})$ 
comes from Serre duality wherein the cotangent space
$$
T^{\vee}(X^{\rm vs})_E = H^1(\text{End}^0(E))^{*} \cong
H^0(\text{End}^0(E)\otimes \omega_C)
$$
identifies with the space of Higgs fields on a given bundle $E\in
X^{\rm vs}$.

On $\Prym ^{\rm vs}:= \Prym \cap \Higgs ^{\rm vs}$ we obtain a map $f^{\rm vs} 
:\Prym ^{\rm vs} \rightarrow X^{\rm vs}$ which, in our case, is a
finite $8:1$ ramified covering.

In order to go from the given spectral data over the open subset
$X^{\rm vs}$ to completed objects over $X$, the first task is to blow
up $\Prym$ in order to resolve the rational map $\Prym \dashrightarrow
X$, given by $f^{\rm vs}$ on $\Prym ^{\rm vs}$, into a morphism.

This might be complicated in the general setting.  In rough terms, the
projection from our hoped-for spectral variety $\Prym ^{\rm vs}
\rightarrow X^{\rm vs}$ should extend to a finite map $Y
\stackrel{f}{\longrightarrow} X$ where $Y\rightarrow \Prym $ is some
kind of a blowing-up of $\Prym $ along the locus $\Fix\cap \Prym$ for
the complement $\Fix := \Higgs -\Higgs ^{\rm s}$ of the set of Higgs
bundles whose underlying bundle is stable (that's a little bigger than
$\Higgs^{\rm vs}$). The subvariety $\Fix$ is the union of the
collection of incoming directions to the higher level fixed point sets
of the $\cstar$ action.

In general, we don't know how to describe the appropriate blow-up. It
will probably contain a sequence of blow-ups along various
subvarieties within $\Fix$. Luckily, in the case of rank $2$ bundles
on a genus $2$ curve, it suffices to blow up once and there are
explicit descriptions:
\begin{itemize}
\item 
In $\Prym_{3} = \{ L \in \op{Jac}^{3}(\Ctilde) \, | \,
\mathsf{Nm}_{\pi}(L) = \mathcal{O}_{C}(3\pw)\}$ we have a smooth curve
$\widehat{C} := \widetilde{C}\times _C \overline{C}$ and $Y_1$ is the
blow-up of $\Prym_{3}$ along this curve. Let $\ExY _1\subset Y_1$
denote the exceptional divisor, mapping to $\Wob _1$. We will denote
the blow-up maps by $\blo$ with indices if necessary. 
\item
In $\Prym _{2} = \{ L \in \op{Jac}^{3}(\Ctilde) \, | \,
\mathsf{Nm}_{\pi}(L) = \mathcal{O}_{C}(2\pw)\}$ we should blow up $16$
points; the exceptional locus $\ExY_{0} = \sqcup_{\kappa \in
  \mathsf{Spin}(C)} \ExY_{0,\kappa}$ is a union of $16$ planes in $Y_0$
labeled by theta characteristics of $C$ and mapping to the
corresponding trope planes in $X_0$.
\end{itemize}

\

\noindent
The inclusion 
$\Prym ^{\rm vs} \hookrightarrow T^{\vee} X^{\rm vs}$ 
extends to an inclusion
$$
Y\hookrightarrow T^{\vee}X(\log \Wob )
$$
well-defined at least up to and including codimension $1$. The
essential reason for this is that the birational map from $Y$ to
$T^{\vee}X$ is given by a global $1$-form $\modalpha$ on $Y$, termed the
\emph{\bfseries spectral $1$-form} that is in fact pulled back from a
$1$-form on $\Prym$ associated to the tautological $1$-form of the
original spectral covering $\Ctilde/C$.  The section over $Y$ of the
pullback of $T^{\vee}X$ can have poles at points where $Y/X$ is
ramified and $\modalpha$ does not vanish on the vertical directions of
the ramification, but in this case the poles are logarithmic.

The \emph{\bfseries spectral correspondence}
\cite{BNR,ron-spectral,ron.eyal-spectral} says that from a spectral
covering contained in the cotangent bundle, together with a line
bundle, we obtain a Higgs sheaf on the base. For parabolic Higgs
sheaves, the spectral covering should be inside the logarithmic
cotangent bundle with logarithmic poles along $\Wob$. In our
situation, the inclusion of $Y$ into the logarithmic cotangent bundle
insures that the logarithmic property holds along the smooth points of
$\Wob$. Understanding what happens at singular points of $\Wob$ is one
of the main technical difficulties.

The image of $Y\hookrightarrow T^{\vee}X(\log \Wob )$ is the spectral
covering predicted by the program of the first two authors
\cite{DP1,DP2012,DonagiPantev}, as stated in Conjecture \ref{conj:DP}
above.  The proposal for realizing this subvariety as a spectral
covering associated to a local system, is to apply the nonabelian
Hodge correspondence initiated by Hitchin \cite{Hitchin-selfd}, and
taken to the case of open varieties of higher dimension by Mochizuki
\cite{Mochizuki-D1, Mochizuki-D1,Mochizuki-kh1, Mochizuki-kh2}, in
order to get a harmonic bundle and therefore a corresponding local
system on $X^{\rm vs}$.

This requires resolving the singularities of $\Wob$ to get a normal
crossings divisor, choosing a line bundle over $Y$ and then over the
resolution of singularities, and putting an appropriate parabolic
structure on the resolution in order to get a logarithmic Higgs field
with respect to a normal crossings divisor.

Once these are done, one needs vanishing of the parabolic Chern
classes $c_1^{\rm par}$ and $c_2^{\rm par}$. The first two authors
conjecture that placing a well-adjusted parabolic structure over the
wobbly locus should insure this vanishing.  This was done for the case
of $\pp^1$ with $5$ marked points in \cite{DonagiPantev}.  That paper
involves a nontrivial interaction between parabolic structures on the
curve $C$ and these parabolic structures over $(X,\Wob )$.

In the present paper devoted to the same the program for bundles of
rank $2$ on smooth projective curves of genus $2$, the parabolic
weight structure is much simpler.  Indeed, we find that a parabolic
structure, with weight $\alpha = 1/2$, is needed for $X_1$; and that
for $X_0$ no parabolic structure is needed over the trope planes or
the Kummer surface.

The Chern class calculations will prescribe the numerical equivalence
class of the line bundle $\LY$ necessary to obtain a parabolic Higgs
sheaf $\mycal{F} := f _{*}\LY$ on $X$ with logarithmic structure
along the smooth locus of $\Wob $ such that the appropriately
calculated $c_1^{\rm par}$ and $c_2^{\rm par}$ vanish.

Over $X_1$, the parabolic structure is localized upstairs near the
exceptional divisor $\ExY_1$.  The map $f$ has ramification index of
$2$ at a general point of $\ExY_1 $, and that gives us the sub-bundle
of the direct image that we will use to put the parabolic structure.

The standard restriction theorems for the Kobayashi-Hitchin
correspondence \cite{Mochizuki-kh1} say that it suffices to work down
to codimension $2$ in $X$.  This convention will be adopted throughout
the paper.  An understanding of the codimension $2$ pieces, namely the
singular locus in codimension $1$ of $\Wob$, is needed in order to
verify the vanishing of the second Chern class.  In our case, up to
codimension two in $X$, $\Wob_1$ has normal crossings and curves along
which it is cuspidal, whereas $\Wob_0$ has normal crossings and
curves along which there is a tacnode (the trope planes meet the
Kummer surface in double conics).

Some special strategies are employed in Chapters \ref{chapter-d1} and
\ref{chapter-d0} to obtain the calculation of parabolic Chern classes
at these singularities, helped by the fact that the parabolic weights
were limited to $0, 1/2$. This promises to be a stumbling block in the
general situation.

The conclusion of the Chern class calculations is that the line bundle
over $Y_0$ should be of the form
$$
\LY_0 =
\blo _0^{*} (\Lprym _0)\otimes f_0^{*}\Oo _{X_0}(2)\otimes
\Oo _{Y_0}(\ExY _0)
$$
where $\Lprym_0$ is a degree $0$ line bundle on the Prym variety $\Prym_{2}$, and
$\blo _0 : Y_0 \rightarrow \Prym _{2}$ is the blow-up map. 

Similarly,  
over $Y_1$ it should have the form 
$$
\LY_1 = \blo _1^{*} (\Lprym _1)\otimes f_1^{*}\Oo _{X_1}(1)
$$
where $\Lprym_1$ is a degree $0$ line bundle on the Prym variety
$\Prym_{3}$ and again $\blo _1$ is the blow-up map. 

These degree $0$ line bundles correspond to the line
bundle $\specN$ over $\Ctilde$ that is input in the geometric
Langlands picture, i.e. the spectral line bundle defining the
eigenvalue Higgs bundle $\Lambda_{Dol}$ on $C$. Showing that these are
the appropriate adjustments in view of the Hecke eigensheaf condition
is a significant part of our work.

The modular spectral coverings $Y_{0}$ and $Y_{1}$ are irreducible so
the constructed Higgs bundles $\mycal{F}_{Dol,0}$ and $\mycal{F}_{Dol,1}$ are
automatically stable.  Applying Mochizuki's Kobayashi-Hitchin
correspondence \cite{Mochizuki-kh1} yields a corresponding flat bundles
$\mycal{F}_{dR,0}$ and $\mycal{F}_{dR,1}$ on $X_{0}^{\rm vs}$ and $X_{1}^{\rm vs}$.

From this construction we derive immediately some basic properties of
these flat bundles:
\begin{itemize}
\item
For $(X_1,\Wob _1)$ the monodromy transformation in $\mycal{F}_{dR,1}$
around a general point of $\Wob _1$ is of order $2$ with two
eigenvalues of $-1$ and six eigenvalues of $1$.
\item
For $(X_0,\Wob _0)$ the monodromy transformations in
$\mycal{F}_{dR,0}$ around general points of the trope planes are
transvections, whereas the monodromy around a general point of $\Kum$
is a direct sum of four transvections.
\end{itemize}

These flat bundles are supposed to be the ones given by the geometric
Langlands correspondence.  In particular, we would like to show that
they satisfy the Hecke eigensheaf property. We do not address the
general problem of uniqueness of Hecke eigensheaves directly, however
it will be shown that the ones produced by our construction agree
with the ones constructed by Drinfeld-Laumon.

The Hecke property requires computing the higher direct image of a
harmonic bundle in the Dolbeault framework. For this we use the
formalism and results from \cite{dirim}. Some further work is needed
here to generalize and adapt the computational machinery of
\cite{dirim} to the specific setup given by the Hecke
correspondences. This is done in Chapter \ref{chapter-pushforward},
and for the reader's convenience we provide a summary of the end
results that will be applied to the Hecke computations, in Section
\ref{pushforward-statements}.

The general setting for the pushforward calculation is a map
$H\rightarrow S$ from a surface to a curve, but let us look at how it
comes about in our application.  As we are dealing with bundles having
a fixed determinant, the Hecke correspondence is parametrized by the
$16:1$ etale covering $\Cbar \rightarrow C$ parametrizing points $t\in
C$ plus a square-root of $\Oo _C(t-\pw)$. Going from $X_i$ ($i=0,1$)
to the other space $X_j$ ($j=1,0$) this is
\[
\xymatrix@M+0.25pc{
  & \bigHeckebar \ar[dl]_-{p}  \ar[dr]^-{q}  & \\
X_i & & X_j \times \Cbar .
} 
\]
Choosing a point $a\in \Cbar$ and restricting to a line $S:=\ell
\subset X_j \times \{ a\}$, let $H=H_{\ell}$ be the pullback of $\ell$
in the Hecke correspondence $\bigHeckebar$. There is a parabolic Higgs
bundle $(\mycal{E},\varphi )$ on $H$ coming by pullback from our
constructed Higgs bundle $\mycal{F}_{Dol,i}$ on $X_i$.  One defines
\cite{dirim} the $H/S$ relative \emph{\bfseries Dolbeault complex} of
$(\mycal{E},\varphi )$ by setting
$$
\Dol (H/S , \mycal{E},\varphi ):=
\left[
W_0 \mycal{E}\stackrel{\varphi _{H/S}}{\longrightarrow} 
W_{-2}\mycal{E}\otimes \Omega ^1_{H/S}
\right]
$$
over $H$.  Here, $W_{k}\mycal{E}$ is the subsheaf of sections of
$\mycal{E}$ whose restriction to the horizontal part of the
parabolic divisor in $H$ lies in the $k$-th piece of the
monodromy weight filtration of ${\rm res}(\varphi )$ on the parabolic
weight $0$ part.

To avoid cluttering the notation we will continue to write $q : H \to S$
for the restriction of the map $q :\bigHeckebar \to X_{j}\times \Cbar$. 
The Dolbeault higher direct image vector bundle on $S$ is
$$
(\mycal{E},\varphi )\mapsto {\mathcal F}:= {\bf R}^1q_{*}
\Dol (H/S , \mycal{E},\varphi ),
$$
and the global Higgs field upstairs leads to a Higgs field $\phi$
on this $\mathcal{F}$

For sufficiently general Higgs bundles, as will be the case in this
instance, the Dolbeault higher direct image takes a particularly nice
form, with the cohomology along each fiber being localized at a finite
set of points corresponding to the zeros of the relative Higgs
field. This family of finite sets gives the spectral variety for the
higher direct image.

There is a natural subscheme ${\rm Crit}$ of the projectivization $\pp
(\mycal{E} / H)$, the \emph{\bfseries relative critical locus}
consisting of the zeros of the relative Higgs field.

A point of ${\rm Crit}$ corresponds to a point $z\in H$ and a vector
$e\in \mycal{E} $ such that $\varphi (e)$ projects to zero in the
$\mycal{E}$-valued relative differentials $\mycal{E} \otimes \Omega
^1_{H/S}$.

In the case where ${\rm Crit} / S$ is finite (i.e. the zeros of the
relative Higgs field are isolated in each fiber), then ${\rm Crit}
\rightarrow S$ will be the spectral cover of the Dolbeault higher
direct image $R^1q_{*}(\mycal{E} , \varphi )$.  Over each point in the
base, this statement says that the cohomology of the fiber localizes
at the zeros of the Higgs field and decomposes naturally as a direct
sum indexed by these zeros.  It is a form of Witten's Morse theory
\cite{Witten-Morse}. The proof is that the relative Dolbeault complex
becomes quasiisomorphic to the cokernel sheaf that is supported on
${\rm Crit}$.

The critical locus description matches up with the \emph{\bfseries
abelianized Hecke variety} defined in general by the first two
authors:
$$
\bigHeckehat ^{\rm ab} := \{ (L, L',\alpha ,A), \;\;
L\stackrel{\alpha}{\rightarrow} L' \}
$$
where $L$ and $L'$ are line bundles over the spectral curve
$\widetilde{C}$ such that $U=\pi _{*} (L)$ has determinant $\Oo _C$
and $U'=\pi _{*}(L')$ has determinant $\Oo _C(t)$, and $A$ is a
square-root of $\Oo_C(t-\pw)$. Thus $U'\otimes A^{-1}$ has determinant
$\Oo _C(\pw )$.  An alternate formulation is
$$
\bigHeckehat ^{\rm ab} := \{ (L, \tilde{t} ,A), \;\; \tilde{t}\in
\widetilde{C} \}
$$
where we put $L':= L(\tilde{t})$ and require $A^{\otimes 2} = \Oo
_C(t-\pw)$ for $t:= \pi (\tilde{t})$.  The abelianized Hecke variety
played a main role in the paper \cite{DP2012} on the classical limit.

The abelianized Hecke correspondence maps to the usual Hecke
correspondence: \linebreak
$\bigHeckehat ^{\rm ab} \rightarrow \bigHeckebar$.
Indeed, an abelianized Hecke transform between line bundles induces a
usual Hecke transform between the rank $2$ bundles they induce by
pushforward from $\Ctilde \rightarrow C$.

In the situation of our application,  $H\hookrightarrow \bigHeckebar$. 
The Higgs bundle $\mycal{E}$ on $H$ is the pullback of the constructed
Higgs bundle $\mycal{F}_i$ that has spectral variety $Y_i/X_i$.

The relative critical locus then identifies with the abelianized
Hecke: ${\rm Crit} = \bigHeckehat^{\rm ab} |_H$.  The basic idea is to
see a point of the full Hecke correspondence is a pair consisting of a
rank two bundle and a rank $1$ subspace over a point of $C$.  When the
bundle comes from a point of $\Prym ^{\rm vs}$, that point is a line
bundle $L$ on $\widetilde{C}$ and the bundle over $C$ is the
pushforward. At a general point this pushforward has two distinguished
directions coming from the two sheets of $\widetilde{C}$.  One shows
that our point in $\bigHeckebar$ is in the zero set of the Higgs field
exactly when the rank $1$ quotient of the bundle goes in one of the
distinguished directions. Hence, the Hecke point comes from
$(L,\tilde{t},A )$, that is to say a point of $\bigHeckehat ^{\rm
  ab}$.

The pushforward computations applied in our case show that the higher
direct image from $H$ to $S=\ell$ is a Higgs bundle having
$\bigHeckehat^{\rm ab} |_H$ as its spectral variety.

If we move back from the restriction to a line in $X$ to looking at
the direct image from $\bigHeckebar$ to $X_j\times \Cbar$, we see that
the higher direct image Higgs bundle has spectral variety
$\bigHeckehat ^{\rm ab}$.  Birationally,
$$
\bigHeckehat ^{\rm ab} \cong Y_j\times \Ctilde \longrightarrow X_j
\times \Cbar .
$$
This is the statement needed on the level of spectral varieties to get
the Hecke eigensheaf property saying that the higher direct image has
the form $\mycal{F} \boxtimes \Lambda$.  To prove the full eigensheaf
property, one needs to identify the spectral line bundle, and to deal
with parabolic structures and various singularities.

This gives the basic idea of the proof of the Hecke property, although
more technical discussion is needed in order to get the statement
precisely.

\subsection{Main theorems}
\label{sec-mt}

This section contains consolidated statements of the main theorems.
Recall that $C$ is a smooth projective curve of genus $2$.  Suppose
given a rank $2$ flat bundle $\Lambda_{dR}$ over $C$, corresponding to
a Higgs bundle $\Lambda_{Dol} = (E,\theta)$ with spectral curve
$\Ctilde \hookrightarrow T^{\vee}C$. Suppose that the spectral curve
is corresponds to a general point in the Hitchin base.More precisely
we will assume that $\Ctilde$ is smooth and unramified over any of the
Weierstrass points of $C$

Use indices $i=0,1$ to indicate the moduli spaces $X_0,X_1$. Let
$X_i^{\circ}:= X_i - \Wob _i^{\rm sing}$ and $\Wob _i^{\circ}:= \Wob_i
- \Wob _i^{\rm sing}$. Let $Y_i^{\circ}$ be the preimage of
$X_i^{\circ}$ in $Y_i$.

\

\begin{theorem}
\label{mt:construction}
There is a tame purely imaginary harmonic bundle over $X_i - \Wob_i$,
corresponding to a pure twistor $\srD$-module whose Dolbeault fiber is
a parabolic logarithmic Higgs bundle $\mycal{F}_{i,Dol,\bullet}$ on
$(X_i^{\circ}, \Wob _i^{\circ})$, such that the parabolic weights
along $\Wob _i^{\circ}$ are trivial for $i=0$ and $0,1/2$ for
$i=1$. The spectral data for $\mycal{F}_{i, Dol,\bullet}$ consist of the
spectral covering $Y_i^{\circ}$ together with a spectral line bundle
$\LY _i$ defined on $Y_i$ as follows:
\begin{itemize}

\item 
  $\LY_0 = \blo _0^{*} (\Lprym _0)\otimes
  f_0^{*}\Oo _{X_0}(2)\otimes \Oo _{Y_0}(\ExY _0)$

\item 
$\LY_1 = \blo _1^{*} (\Lprym _1)\otimes f_1^{*}\Oo _{X_1}(1)$,
\end{itemize}
where $\Lprym _0$ and $\Lprym _1$ are line bundles on the Prym variety
associated to the spectral line bundle of $\Lambda_{Dol} = (E,\theta )$, and the
spectral $1$-form is given by the tautological form from the inclusion
$$
Y_i^{\rm vs} \cong \Prym _i^{\rm vs} 
\hookrightarrow \Higgs _i^{\rm vs} \cong T^{*}(X^{\rm vs}).
$$ The parabolic structure is given from the spectral data by setting
$\mycal{F}_{i, Dol, 0} =f_{i*} \LY_i$ for $i =0,1$ and, in addition,
in case $i=1$, setting $\mycal{F}_{1,Dol,1/2} = f_{1*} \LY_i(\ExY _1)$.
\end{theorem}

\

\noindent
In the twistor $\pp^1$ we can consider the de Rham point $\lambda =
1$; let $\mycal{F}_{i,dR}$ be the $\srD$-module associated to the
fiber over $\lambda = 1$. Let $\mycal{F}_{i,B}$ be the perverse sheaf
corresponding to this $\srD$-module. Similarly, we will write
$\Lambda_{dR}$ and $\Lambda_{B}$ for the flat bundle and local system on $C$
corresponding to $\Lambda_{Dol}$.

We recall that the Hecke operations are defined for points in the
curve $\Cbar$ that maps to $C$ by a $16:1$ etale covering $\sq : \Cbar
\to C$, so the big Hecke operation goes from sheaves on $X_0$ to
sheaves on $X_1 \times \Cbar$ and vice-versa.

\

\begin{theorem}
\label{mt:Hecke}
The pair of perverse sheaves over $X_0 \sqcup X_1$ is a Hecke
eigensheaf with Hecke eigenvalue $\Lambda$ in the sense that the big
Hecke operation applied to $\mycal{F}_{0,B}$ is $\mycal{F}_{1,B}
\boxtimes \sq ^{*}\Lambda_{B}$ and the big Hecke operation applied to
$\mycal{F}_{1,B}$ is $\mycal{F}_{0,B} \boxtimes \sq ^{*} \Lambda_{B}$.
\end{theorem}

\

Drinfeld \cite{Drinfeld} used a Radon transform to construct Hecke
eigensheaves for rank $2$ local systems $\Lambda_{B}$ on smooth compact
curves of any genus. His construction was put into a more geometric
form by Laumon \cite{Laumon95}.

\

\begin{theorem}
\label{mt:conj-Drinfeld}
The purely imaginary tame harmonic bundles associated to the Hecke
eigensheaves constructed by Drinfeld in rank $2$ have Dolbeault fiber,
i.e. parabolic logarithmic Higgs bundles, that satisfies Conjecture
\ref{conj:DP} for compact curves of any genus.
\end{theorem}

\

\begin{theorem}
\label{mt:comparison-Drinfeld}
For a curve $C$ of genus $2$, Drinfeld's Hecke eigensheaves on $X_0$
and $X_1$ coincide with the $\mycal{F}_{0,B}$ and $\mycal{F}_{1,B}$
that come from the Higgs bundles constructed in Theorem
\ref{mt:construction}.
\end{theorem}

\

\subsection{Structure of the paper}
\label{sec-str}

The case we consider in this paper, moduli of rank $2$ bundles on a curve of genus $2$,
has very classical roots. The basic geometry involved may be viewed as coming from
the expression of $X_1$ as the intersection of a pencil of quadrics in $\pp^5$. This
viewpoint is recalled and developed in Chapter \ref{synthetic}, and we prove some
of the properties needed, notably concerning the lines in $X_1$. Hecke correspondences
make their appearance here from a synthetic point of view. 

In Chapter \ref{chapter-general} we introduce the basic notations of the modular
approach in a
more complete way than was done in the introduction, and discuss
several different types of general considerations that will be used
later in the discussion. These range from Chern class calculations
(\ref{sec-Chern}) to the geometry of the $\cstar$ flow on the Hitchin
moduli space (\ref{orbits}, \ref{logprop}), and include a discussion
of the nonabelian Hodge correspondence (\ref{subsec-nah}) and
parabolic structures (\ref{subsec-parabolic}) to be used in the basic
construction.  We state in (\ref{pushforward-statements}) the results
on Dolbeault higher direct images in a form that will be most useful
for calculating the Hecke correspondences; their proofs are deferred
to Chapter \ref{chapter-pushforward}.

In Chapter \ref{chapter-d1} we construct the parabolic logarithmic
Higgs bundle for our candidate Hecke eigensheaf over the moduli space
$X_1$ of bundles of degree $1$.  This involves a precise description
of the wobbly locus $\Wob_1$ paired with a description of the blown-up
Prym that forms the spectral covering of $X_{1}$, with exceptional
divisor $\ExY_1$ above $\Wob _1$.  The main technical work is to
arrive at a calculation of the parabolic Chern class, given that $\Wob
_1$ has cuspidal (as well as normal-crossings) singularities in
codimension $2$ in $X_{1}$.  The technique used here is to pass to a
finite covering of Kawamata type \cite{Kawamata-cover}, with smooth
total space and having ramification of order two along $\Wob _1$. This
works because the parabolic structure to be used here has weights
$0,1/2$ so it goes away upon pullback to the finite cover.  It
therefore does not matter that the inverse image of $\Wob _1$ has a
triple point where the cusp used to be: the Higgs bundle extends
smoothly across the divisor and we can just compute its Chern class.

In Chapter \ref{chapter-d0} we construct the parabolic logarithmic
Higgs bundle for our candidate Hecke eigensheaf over the moduli space
$X_0$ of bundles of degree $0$.  This again involves a precise
description of the wobbly locus $\Wob _0$, which turns out to be
Kummer's $16_6$ configuration combining the Kummer surface with $16$
nodes in $X_0\cong \pp^3$, with the $16$ trope planes meeting the
surface along trope conics that transversally make tacnodes.  Once
again, the main problem is how to compute the parabolic Chern class
contributions from the singularities. In this case, at smooth points
of $\Wob _0$ there is no parabolic structure, but rather the Higgs
field has nonzero nilpotent residues. The corresponding monodromies of
the local system are unipotent. However, a naive extension of the
bundle across the tacnodes has Higgs field that isn't logarithmic on a
resolution. Our technique in this case, different from the case of
degree $1$, is to resolve the singularities of the tacnodes by two
blow-ups, and then put an appropriate parabolic structure over the
exceptional divisors. In other words, even if the Higgs bundle does not
have parabolic structure along the smooth points of $\Wob _0$, it does
have a `hidden' parabolic structure inside the tacnodes, and indeed
the resulting local system will have nontrivial monodromy eigenvalues
around components of the exceptional divisors. The good parabolic
structure was found by trial and error using some computer
calculations. Those were really bad so we don't reproduce them here,
rather we just state what is the good parabolic structure, verify that
it makes the Higgs field logarithmic, and verify that it yields
vanishing of the parabolic Chern classes.

From Chapters \ref{chapter-d1} and \ref{chapter-d0} we thus obtain the
constructions of local systems, and their associated purely imaginary
tame harmonic bundles, on $X_1^{\rm vs}$ and $X_0^{\rm vs}$,
completing the proof of Theorem \ref{mt:construction}.  The remainder
of the paper is devoted to verifying the Hecke eigensheaf property
plus a few other things.

In Chapter \ref{chapter-hecke} we introduce the general setup of the
Hecke correspondence in the modular viewpoint.

In Chapter \ref{chapter-abelianized} we introduce the \emph{\bfseries
abelianized Hecke} variety, which is the main player in the proof of
the Hecke property. In this chapter, we give a first approach by
showing that the Hecke property holds at the level of spectral data
using the abelianized Hecke.  This part involves consideration of the
``big Hecke correspondence'' which is the total space of the family
parametrized by points of the curve $\Cbar$.  Consideration of the
Hecke correspondences at a single point is done in the next two
chapters.

In Chapter \ref{chapter-heckex0x1} we fix a point $a\in \Cbar$ and
show that the Hecke correspondence at the point $a$ takes the
constructed Higgs bundle on $X_0$ to the constructed Higgs bundle on
$X_1$. This is done by restricting to a general line, and applying the
pushforward statements given in Subsection
\ref{pushforward-statements} and to be proven in Chapter
\ref{chapter-pushforward}. This section contains a subtle point about
apparent singularities: the Hecke pushforward morphism seems to be
singular along an additional subvariety of $X_1$, namely a Kummer K3
surface, known classically, and that depends on the point $a$.  We
need to show that the higher direct image harmonic bundle does not
really have singularities there. This requires the discussion of
Subsection \ref{t5points}.

In Chapter \ref{chapter-heckex1x0} for the fixed point $a\in \Cbar$ we
show that the Hecke correspondence at $a$ takes the constructed Higgs
bundle on $X_1$ to the constructed Higgs bundle on $X_0$. Again the
pushforward statements are applied after restricting to a general
line. In this direction, the difficulty with apparent singularities
does not occur.

In Chapter \ref{chapter-bighecke} we go back to the big Hecke
correspondence, and use the results of the previous chapters to prove
the Hecke eigensheaf property.  This completes the proof of Theorem
\ref{mt:Hecke}.

In Chapter \ref{chapter-thirdconstruction} we propose a third
construction of a parabolic logarithmic Higgs bundle on $X_0$, having
trivial parabolic structure and nilpotent residues over the trope
planes, but with parabolic structure of weights $0,1/2$ over the
Kummer subvariety. We posit that this should be associated to a
$PGL_2$ local system of odd degree, but that is not treated here.

In Chapter \ref{chapter-pushforward} we proceed in several stages to
apply the $L^2$ Dolbeault direct image formulas of \cite{dirim} to
obtain the pushforward statements required for computation of the
Hecke transforms in our cases. The first step is to extend the general
theory to the case when the parabolic divisors can have multiplicity.
Then we consider the case where the horizontal divisor has a simple
ramification point, by blowing up. Then we treat the case of ``points
of type \ref{type5}'' that is needed to show that the higher direct
image does not have singularities at the apparent singularities that
show up in the $(X_0\rightarrow X_1)$ direction of the Hecke
operation. This chapter completes the proof of the pushforward
statements made in Subsection \ref{pushforward-statements}.

In Chapter \ref{chapter-drinfeld} we change gears and consider
Drinfeld's Radon transform construction.  We show how the Dolbeault
direct image technology can apply to gain information about the Higgs
bundles associated to the perverse sheaves constructed by Drinfeld. In
particular, we obtain the spectral coverings of these Higgs bundles,
and prove Conjecture \ref{conj:DP} for them, i.e.  Theorem
\ref{mt:conj-Drinfeld}. We also show that the Hecke eigensheaves that
we construct coincide with those constructed by Drinfeld, which is
Theorem \ref{mt:comparison-Drinfeld}.  At the end of this chapter, we
show that one can get an explicit description of the Hecke eigensheaf
over $X_1$.

\subsection{Acknowledgements}
\label{sec-ack}

We would like to thank the many colleagues who have responded with
comments and suggestions on the occasion of talks we have given on
this work. In particular we would like to thank Dima Arinkin, Tamas
Hausel, Nigel Hitchin, Jochen Heinloth, Ludmil Katzarkov, Frank
Loray, Ana Pe\'{o}n-Nieto, Vivek Shende, and Edward Witten for
illuminating conversations and insightful comments.

During the preparation of this work Ron Donagi was supported in part
by NSF grants DMS 1603526, DMS 2001673, by NSF FRG grant DMS 2244978,
and by Simons HMS Collaboration grant \# 390287. Tony Pantev was
supported in part by NSF grant DMS 1901876, by NSF/BSF grant
DMS-2200914, NSF FRG grant DMS 2244978, and by Simons HMS
Collaboration grant \# 347070. Carlos Simpson was supported by a grant
from the Institute for Advanced Study, by the Agence Nationale de la
Recherche program 3ia Côte d'Azur ANR-19-P3IA-0002, and by the
European Research Council Horizons 2020 grant 670624 (Mai Gehrke's
DuaLL project).  This work was completed during a visit to the
University of Pennsylvania supported by the Simons Collaboration on
Homological Mirror Symmetry through grants \# 390287 and \# 347070.

\section{Synthetic approach}  \label{synthetic}

A special advantage in studying the Geometric Langlands Conjecture for
curves of genus $2$ is that we can utilize two independent geometric
approaches to the problem.  We can study our various objects modularly,
in term of their interpretation as moduli spaces of bundles with
various decorations, or synthetically, in terms of the geometry of the
intersection of two quadrics. In this section we follow the synthetic
approach. In the rest of this work we follow the modular approach. The
modular techniques we develop are considerably more difficult, but
they are less dependent on specifics of our particular situation, so
it is more likely that they can be extended to more general cases of
GLC.

\

In this section we will describe synthetically the  geometric
correspondences between most fundamental objects of interest in this
paper - the base locus $X_{1}$ of a general pencil of quadrics in
$\pp^{5}$ and a ruling $X_{0}$ of a general quadric in the pencil. We
will also see synthetically how this geometry ties up with the genus
two curve $C$ parametrizing the rulings of the quadrics in the pencil
and with the geometry  of the Jacobian of $C$.

The significance of the synthetic considerations stems from the fact
that, through the works of Narasimhan-Ramanan \cite{NR} and Newstead
\cite{Newstead}, the objects $X_{1}$ and $X_{0}$ are identified with
the moduli of semistable rank $2$ bundles of fixed determinant of degree $1$
and $0$ on the curve $C$. Later on we will refine these
identifications and will in particular show that synthetic
constructions of this section reproduce the Hecke correspondences and
their discriminant (wobbly) loci in moduli. This will give us easy
geometric proofs of most of the important geometric properties needed
for the analysis of the Hecke eigensheaf condition.

\

\subsection{Pencils of quadrics in \texorpdfstring{$\pp^{5}$}{P5}}

We will study the geometry of a general pencil of quadrics in
$\pp^{5}$ by analyzing the families of linear subspaces contained in
these quadrics.

\

\punkt {\bfseries Notation} \ Our notation in this section is: 

\smallskip

\begin{description}[2pc]
\item[$X=X_1 = \cap _{x \in \pp^1} Q_x \subset \pp^5$] is the smooth
  intersection of a generic pencil of quadrics in $\pp^5$.
\item[$\ell \subset X$] is a line in $X$.
\item[$v \in X$] is a point in $X$ . 
\item[$\mathsf{Grass}(3,6)$] is the Grassmannian of projective planes $\Pi
  \subset \pp^5$.
\end{description}

\

\noindent
The universal ruling is $\mycal{R} := \left\{(x,\Pi) \in \pp^1\times
\mathsf{Grass}(3,6) \, \left| \, \Pi \subset Q_x\right.\right\}$.  The
fiber of $\op{pr}_1: \mycal{R} \to \pp^1$ over a general $x \in \pp^1$
has 2 components, the rulings of the quadric $Q_x$. We can thus
consider the Stein factorization of $\op{pr}_1: \mycal{R} \to \pp^1$:
\[
\xymatrix@1@C+0.1pc@M+0.3pc{
  \mycal{R}  \ar[r]^-{\rul} \ar@/_1pc/[rr]_-{\op{pr}_{1}}
  & C \ar[r]^-{\hyp_{C}} & \pp^1. 
}
\]
Then $C$ is a hyperelliptic curve of genus $2$, with $6$ Weierstrass
points $\pw_i$, which are the ramification points of the hyperelliptic
map $\hyp_C$. They correspond to the rulings on the $6$ singular quadrics
in the pencil.  We fix one of them, labelled $\pw \in C$.

The fibers of $\rul$, parametrizing planes in a given ruling $t \in
C$, are spinor varieties, isomorphic to $\pp^3$. The subvariety
parametrizing planes in ruling $t$ that pass through a specified point
$v \in X$ is a line in this $\pp^3$.  Through a specified line $\ell
\subset X$ there is a unique plane in each ruling.

We choose coordinates on $\pp^5$ that are adapted to our pencil,
meaning that the $6$ coordinate points are the vertices of the $6$ quadric
cones in our pencil. All the quadrics $Q_x$ become simultaneously
diagonalized. If we take the matrix of one of the quadrics to be the
identity and that of another quadric in the pencil to be
$diag(\brx_1,\dots,\brx_6)$, the equation of $C$ becomes $y^2 = \Pi _{i=1}^6
(x-\brx_i)$, and $\hyp_C(\pw_i) = \brx_i$.

The modular approach studies bundles on this
genus $2$ hyperelliptic curve $C$.  When we wish to compare the
synthetic and modular approaches, the line $\ell$ will correspond to a
line bundle $L \in \op{Pic}^0(C)$, and the point $v$ will correspond to a
rank $2$ vector bundle $V$ on $C$ with determinant $\det(V) \cong
\Oo_C(\pw)$.

\subsection{Basic geometry of the lines}

Let $\lA$ denote the variety of lines $\ell \subset X$ in $X$. 

\

\begin{lemma} \label{4lines}
Through any point $v \in X$ there are four lines $\ell \in \lA$ (counting
multiplicities).
\end{lemma}
\begin{proof} For a point $v \in X$ let
  $\ttt_{v}X \cong \pp^{3} \subset \pp^{5}$ denote the projective
  tangent space to $X$ at $v$.  A line $\ell \in \pp^5$ passing
  through $v$ is in $X$ iff it is in $X \cap \ttt_{v}X$. But $X \cap
  \ttt_{v} X$ is a curve in $\ttt_{v}X \cong \pp^{3}$ that has degree
  $4$ and is a cone with vertex $v$, so it consists of $4$ lines.
\end{proof}

\

\noindent In any ruling of a quadric $Q_x$ there is a unique plane
$\Pi$ containing any $\ell \subset Q_x$, and in particular any $\ell
\subset X$. In the latter case, the intersection $\Pi \cap X$ consists
of two lines, the given $\ell$ and another line $\ell' \subset X$.  We
get a natural map
\[
\lphi: C \times \lA \to \lA
\]
sending $(t,\ell) \mapsto \ell'$, where $t \in C$ labels a ruling
$\mycal{R}_{t}$ of the quadric $Q_{\hyp_{C}(t)}$, $\ell \subset X$ is a
line, $\Pi$ is the plane in $\mycal{R}_{t}$ containing $\ell$, and
$\ell'$ is the other line in $\Pi \cap X$.  Fixing $\ell \in \lA$ we get
a map $\lphi_{\ell} : C \to \lA$.  Fixing $t \in C$ we get an involution
$\lphi_{t}: \lA \to \lA$.  We will be particularly interested in
$\lphi_{\pw}$, where $\pw \in C$ is our chosen Weierstrass point.  The
action of $C \times C$ on $\lA$ sending $\ell \to \lphi_{t_{1}} \circ
\lphi_{t_{2}}(\ell)$ is easily seen \cite{donagi-quadrics} to descend to
$\op{Sym}^2(C)$ and further to the Jacobian $\mathsf{J} =
\op{Jac}(C)=\op{Pic}^0(C) \cong \op{Pic}^2(C)$, where the last
identification sends $0$ to the canonical bundle $\omega_C =
\Oo_C(2\pw)$. The resulting map
\[
\mathsf{J} \times \lA \to \lA
\]
turns $\lA$ into a $\mathsf{J}$-torsor. The choice of $\ell \in \lA$
thus gives an isomorphism
\[
\lalpha_{\ell} : \mathsf{J}  \stackrel{\cong}{\to} \lA
\]
sending the origin to $\ell$ and restricting to $\lphi_{\ell}$ on $C$,
embedded in $\mathsf{J}$ via Abel-Jacobi with $\pw$ mapped to the
origin.  Note that if $\ell, \emm \in \lA$ are two lines in $X$, then
$\emm = \lalpha_{\ell}(M)$ for a unique $M \in \mathsf{J}$, and  
\[
\lalpha_{\emm} = \lalpha_{\ell}  \circ \trans_M,
\] 
where $\trans_M : \mathsf{J} \to \mathsf{J}$ is translation by $M$: 
$\trans_M(N) := M \otimes N$.

\

\noindent
As above, given a smooth point $z$ of a subvariety $Z \subset
\pp^{N}$ in a projective space, we let $\ttt_{z}Z \cong \pp^{\dim Z}
\subset \pp^{N}$ denote the projective tangent space at $z$, i.e. the linear
subspace of $\pp^{N}$ containing $z$ and pointing in the direction of the
usual tangent space $T_{z}Z$.

\

\begin{lemma} \label{16lines}
There are 16 lines $\lO \in \lA$ that are fixed by the involution
$\lphi_{\pw}$:
\[
\lO = \lphi_{\pw} \lO.
\]
\end{lemma}
\begin{proof}
We start with an initial $\ell \in \lA$ and will see how to modify it to be
a fixed point.  Conjugating $\lphi_{\pw}$ by $\lalpha_{\ell}$ we get an
involution
\[
\lalpha_{\ell}^{-1} \circ \lphi_{\pw} \circ \lalpha_{\ell}: \mathsf{J}
\to \mathsf{J}.
\]
This sends $L \in \mathsf{J}$ to $\ldelta_{\ell} \otimes L^{-1}$,
where $\ldelta_{\ell} := \lalpha_{\ell}^{-1}(\lphi_{\pw} \ell) \in
\mathsf{J}$.  This involution of $\mathsf{J}$ has $16$ fixed points,
namely the square roots $M$ of $\ldelta_{\ell}$.  Replacing our initial $\ell$
by $\lO := \lalpha_{\ell}(M)$ gives a new line $O$ such that $\lalpha_{\lO} =
\lalpha_{\ell} \circ \trans_M$, so
\[
\lalpha_{\lO}^{-1} \circ \lphi_{\pw} \circ \lalpha_{\lO} = \trans_M^{-1} \circ
\lalpha_{\ell}^{-1} \circ \lphi_{\pw} \circ \lalpha_{\ell} \circ \trans_M
\]
is just inversion, so $\lO = \lphi_{\pw} \lO$ as desired.  (The
existence of $16$ $\lO$'s with this property was shown in
\cite{donagi-quadrics} via a Schubert cycle calculation.)
\end{proof}

\

\noindent
If we choose one of these fixed points $\lO\in \lA$ as origin, then
the composition
\[
(\lalpha_{\lO})^{-1} \circ \lphi \circ (\op{id} \times \lalpha_{\lO}) :
C \times \mathsf{J} \to \mathsf{J}
\]
becomes the Abel-Jacobi map
\[
(t,L) \mapsto L(t-\pw),
\]
so its restriction:
\[
\lalpha_{\lO}^{-1} \circ \lphi_{\lO}: C \to \lA \to \mathsf{J}
\]
becomes the Abel-Jacobi map
\[
t \mapsto \Oo_{C}(t-\pw),
\]
and the involutions 
\[
\lalpha_{\lO}^{-1} \circ \lphi_{r} \circ \lalpha_{\lO}: \mathsf{J} \to \mathsf{J}
\]
\[
L \mapsto L^{-1}(t-{\pw}).
\]
so as we have seen, $\lalpha_{\lO}^{-1} \circ \lphi_{\pw} \circ
become:
\lalpha_{\lO}: \mathsf{J} \to \mathsf{J}$ is just inversion.

\

\noindent
For the modular/synthetic dictionary, we set
\[
\ell = \lalpha_{\lO}(L), \qquad L= \lalpha_{\lO}^{-1}(\ell).
\]
Note that while $\lA$ depends only on $X$, the isomorphism
$\lalpha_{\lO}$ depends on the auxiliary choices of the Weierstrass
point $\pw$ and the origin $\lO$.

\

\

\punkt {\bfseries Line incidence and special lines} \
Given $\ell\in \lA$, consider the family 
\[
I_{\ell} \subset \lA
\] 
of lines $\emm \in \lA$ that intersect $\ell$.
By definition, this is the closure (in $\lA$) of:
\[
I^{\circ}_{\ell}  := \left\{\emm \in \lA \,  \left| \,
\emm \neq \ell, \emm \cap \ell \neq \varnothing  \right.\right\}.
\]
For general $\ell$, $I^{\circ}_{\ell} = I_{\ell}$ is closed, so we do not
need the closure. We say that the line $\ell$ is \emph{\bfseries
special} if it ``intersects itself'', in the sense that $\ell \in
I_{\ell}$, so $I^{\circ}_{\ell} \neq I_{\ell}$ (Compare with \cite[page 792]{GH}).

Our isomorphism
$\lalpha_{\lO} : J \stackrel{\cong}{\to} \lA$
identifies $I_{\ell}$ with the theta divisor
\[
\Theta_{\ell} = \Theta_{L} = \left\{ M \in \mathsf{J} \, \left| \,
h^0(L \otimes M(\pw) >0 \right.\right\} = \left\{\Oo_{C}(t-\pw)
\otimes L^{-1}) \, \left| t \in C\right.\right\} \subset \mathsf{J},
\]
where $\ell = \lalpha_{\lO}(L)$.

$I_{\ell}$ is the image of our map  $\lphi_{\ell}: C \to
\lA$, and in fact $\lphi_{\ell}$ induces an isomorphism \linebreak
$\lphi_{\ell}
: C \stackrel{\cong}{\to} I_{\ell}$.  The composed map $
(\lalpha_{\lO})^{-1}\circ \lphi_{\ell} : C \to \Theta_{L}$ is just the
Abel-Jacobi map:
\[
t \mapsto (L(\pw))^{-1}(t).
\]
Consider the curve
\[
\Cbar := \left\{ \, (L,t) \ \left| \  L^2 \cong \Oo(t-\pw) \,
\right.\right\} \subset
\mathsf{J} \times C.
\]
 The second projection realizes $\Cbar$ as the $16$-sheeted cover $\sq:
 \Cbar \to C$ induced from the doubling map $\mathsf{J} \to
 \mathsf{J}, \ \ L \mapsto L^2$ and the Abel-Jacobi map $C \to
 \mathsf{J}, \ \ t \mapsto \Oo_{C}(t-\pw)$.  The first projection
 identifies $\Cbar$ with its image in $\mathsf{J}$. Composing with the
 isomorphism $\lalpha_{\lO}$, we identify $\Cbar$ as a subvariety of
 $\lA$. In summary, this gives the following

 \

\begin{corollary}
  The line $\ell$ is special iff $L = \lalpha_{\lO}^{-1}(\ell) \in \Cbar
  \subset \mathsf{J}$, iff $L^2(\pw)$ is effective.
\end{corollary}

\

\noindent
From now on we will think of $\Cbar$ as the subvariety of $\lA$
parametrizing special lines. It is intrinsic to $X$, independent of
choices of $\pw$ and $\lO$.

A central object in our study is the Wobbly locus in $X$. We define it
to be:
\[
\Wob := \bigcup_{\ell \in \Cbar}\  \ell \ \subset \ X,
\]
i.e. the union of all the special lines.

\

\punkt {\bfseries Trigonal bundles} \ Via the map $\lphi_{\ell} : C
\to I_{\ell}$, we have identified $C$ with the family $I_{\ell}$ of
lines $\emm$ meeting a given line $\ell$, while $\lalpha_{\lO}$
identifies $I_{\ell}$ with $\Theta_{\ell}$. The composition
$\lalpha_{\lO}\circ \lphi_{\ell} : C \to \Theta_{\ell} $ sends $t
\mapsto \Oo_{C}(t-\pw) \otimes L^{-1}$, where again $L =
\lalpha_{\lO}^{-1}(\ell)$.

\

\begin{lemma} Let $L = \lalpha_{\lO}^{-1}(\ell)$. 
The map $\mathsf{a}^{\circ}=\mathsf{a}_{\ell}^{\circ}:
I_{\ell}^{\circ} \to \ell$ sending $\emm \to \emm \cap \ell$ extends
to $\mathsf{a} = \mathsf{a}_{\ell}: C \cong I_{\ell} \to \ell$, where
it is given by the linear system of sections of the degree $3$ line
bundle $L^{2}(3\pw)$.
\end{lemma}
\begin{proof}
We have a morphism $X \to \op{Pic}^{0}(C)$ sending $v$ to
$\otimes_{i=1}^4 L_i$, where the $L_i$ are the line bundles
corresponding to the $4$ lines $\ell_i$ through $v$, i.e. $L_i =
\lalpha_{\lO}^{-1} (\ell_i)$. Since $X$ is unirational (in fact,
rational), this map must be constant. Our choice of origin assures us
that this constant value of this map is the origin in
$\op{Pic}^{0}(C)$. So for each point $v$ of $\ell$, the sum of the
three lines other than $\ell$ through $v$ (i.e. the product of the
corresponding line bundles) must be $L^{-1} \in \mathsf{J}$. The
identification of $\Theta_{\ell}$ with $C$ involves a translation by
$L(\pw)$. Suppose that the four lines $\{ \ell_{i}\}_{i=1}^{4}$
through $v \in \ell$ are labeled so that $\ell_{4} = \ell$.  So when we
convert the three $L_i$ to points $t_i \in C$, we see that their sum,
in $\op{Pic}^3(C)$, is:
\[
\Oo_C (\Sigma_{i=1}^3 t_i) = L^{-1} \otimes (L(\pw)^{\otimes 3}) = L^2(3\pw)
\]
as claimed.  
\end{proof}

\

\noindent
There are two possibilities: $L(3\pw)$ could be base-point free and
therefore give a genuine trigonal map, or it could have a base
point. This happens iff $L(3\pw) = \omega_C(t)$ for some $t \in C$, or
equivalently iff $L^2 \cong \Oo_{C}(t-\pw)$. We immediately get

\

\begin{corollary}
The trigonal bundle $L^{2}(3\pw)$ has a base point iff $\ell$ is special,
or equivalently $(L,t)$ belongs to $\Cbar$.
 \end{corollary}

\

\noindent
In the case of a special $\ell$, the map $\mathsf{a}_{\ell}: C
\dashrightarrow \ell$ sending $\emm \to \emm \cap \ell$ is still given
by the linear system of sections of $L^{2}(3\pw)$, but now this linear
system has the base point $t$, so it is only a rational map. That
means that the trigonal curve is now reducible, $C \cup \pp^{1}$ with
$t \in C$ glued to $\hyp_C(t) \in \pp^1$, where $\hyp_C: C \to \pp^1$
is the hyperelliptic double cover. The rational map
$\mathsf{a}_{\ell}$ lifts to a morphism
$\widetilde{\mathsf{a}}_{\ell}: C \cup \pp^1 \to \ell$. When
restricted to $\pp^1$ this gives a natural isomorphism
$\mathsf{a}^{\circ}_{\ell}: \pp^1 \to \ell$, and its restriction to
$C$ agrees with $\hyp_C : C \to \pp^1$ followed by this isomorphism
$\mathsf{a}^{\circ}_{\ell} : \pp^1 \to \ell$.  In this case the $3$
lines meeting $\ell$ at any $v \in \ell$ consist of the fixed line
$\ell$, plus the moving hyperelliptic pair $\hyp_{C}^{-1}(v)$, or more
precisely, their images under $\lphi_{\ell}$. So the $4$ lines through
any such $v$ consist of twice $\ell$, plus the pair
$\hyp_{C}^{-1}(v)$.

Going back to a general $\ell$, the trigonal map $\mathsf{a} : C \to
\ell$ has $8$ branch points $\brb_i \in \ell$, with
$\mathsf{a}^{-1}(\brb_i) =2r_i + s_i$ consisting of the ramification
point $r_i$ and one other point $s_i$. The corresponding $4$ lines
through such a $\brb_i$ are then $\ell, \emm_{i}$, and $\ell_i$ occurring
with multiplicity $2$. Here $\ell_i$ corresponds to the degree $0$ line
bundle $L^{-1}(r_i -\pw)$, and similarly $\emm_{i}$ corresponds to the
degree $0$ line bundle $L^{-1}(s_i -\pw)$. Here $\ell$ is arbitrary, but the
$\ell_i$ are in $\Cbar$. In fact, the $8$ ramification points are the
intersection of $\Cbar$ with an appropriate translate of the $\Theta$
divisor.

\

\

\punkt {\bfseries Patterns of lines} \ \label{sssec:lpatterns}
We have now proved
everything we need in order to describe all possible patterns of the $4$
lines through some $v \in X$. These are summarized in the follwoing:

\

\begin{theorem}{\label {patterns}}
For a point $v$ in a line $\ell \in \lA$, the possible patterns of lines
through $v$ are:
\begin{itemize}
\item $2\ell + \ell_{1} + \ell_{2}$: \ $\ell$ counts twice among the
  $4$ lines through $v$ iff $\ell \in I_{\ell}$, iff $\ell$ is
  special, iff the trigonal line bundle $L^{2}(3\pw)$ has a base
  point.
\item $\ell + 2\ell_{1} + \ell_{2}$: \ One of the lines in
  $(\mathsf{a}_{\ell})^{-1}(v)$ counts twice iff $v$ is a branch point
  of $\mathsf{a}_{\ell} :C \to \ell$.
\item $2\ell+2\ell_{1}$: \ The fiber $(\mathsf{a}_{\ell})^{-1}(v)$
  consists of two lines $\ell,\ell_{1}$, each counted twice, iff the
  lines $\ell,\ell_{1} \subset X$ intersect at the point $v$, which is
  then a singular point of $\Wob$. This occurs iff
  $\ell_{1}=\lalpha_{\ell}(\pw_i) = \lphi(\pw_i, \ell)$ where $\pw_i
  := (\lalpha_{\ell})^{-1}(\ell_{1})$ is one of the 6 Weierstrass
  points of the hyperelliptic $C$, and $v=\hyp_C(\pw_i)$ is its image.
\item $3\ell+\ell_{1}$: \ The line $\ell$ counts $3$ times among the
  $4$ lines through $v$ iff $\ell$ is special, corresponding to $(L,t)
  \in \Cbar$, and $v=\hyp_C(t)$.
\item $\ell+3\ell_{1}$: \ $(\mathsf{a}_{\ell})^{-1}(v)$ consists of a
  single line $\ell_{1}$ counted $3$ times iff
  $(\lalpha_{\ell})^{-1}(\ell_{1})$ is a total ramification point of
  $\mathsf{a}_{\ell} : C \to \ell$, and $v = \hyp_C
  ((\lalpha_{\ell})^{-1}(\ell_{1}))$ is the corresponding total branch
  point.
\item $4\ell$: \ The lines that count $4$ times in a fiber are the
  $16$ origin-candidates $\lO$ from Lemma~\ref{16lines}.
\end{itemize}
\end{theorem}

\

\

\noindent
Note in particular that pattern $2\ell+ \ell_{1}+ \ell_{2}$ is
independent of $v \in \ell$. We can prove this directly:

\

\begin{lemma}
If $\ell$ counts twice through some $v \in \ell$, it does so for every
$v \in \ell$.
\end{lemma}
\begin{proof}
Given $\ell \subset X = \cap_{x \in \pp^1} Q_x$, let $\pp = \pp^3_{\ell}$ be
the quotient projective space $\pp^5 /\ell$, parametrizing planes in
$\pp^5$ through $\ell$, and let $\pp^{\vee}$ denote its dual,
parametrizing hyperplanes in $\pp^5$ through $\ell$. Let the line
$\ellfrak_v$ denote the projection of $\ttt_v X$ to $\pp$, and
$\ellfrak^{\vee}_v$ the dual line in $\pp^{\vee}$.

Consider the morphism $\ell \times \pp^1 \to \pp^{\vee}$ sending
$(v,x) \mapsto \ttt_v Q_x$. It is well-defined (each $v$ is a
non-singular point of each $Q_x$) and linear in each of its arguments
$v,x$, i.e. it is given by a base-point free sublinear system of
$\Oo_{\pp^{1}\times \pp^{\vee}}(1,1)$. There are only two
possibilities: the map can be an embedding, identifying $\ell \times
\pp^1$ with a smooth quadric $\lquad^{\vee} \subset \pp^\vee$, or it
can be $2$-to-$1$ onto a plane $\Pi \subset \pp^\vee$, with branch
locus a conic $\lC \subset \pp^\vee$. In the first case, the lines
$\ellfrak^{\vee}_v$ form a ruling of $\lquad^{\vee}$, and dually the
lines $\ellfrak_v$ form a ruling of the dual quadric $\lquad \subset
\pp$. In the second case, the lines $\ellfrak^{\vee}_v$ are the
tangents in $\Pi$ to $\lC$, and dually the lines $\ellfrak_v$ form the
ruling of a quadratic cone $\lquad \subset \pp$.

Now the line $\ell$ counts more than once among the $4$ lines in $X$
through $v \in \ell$ if and only if the intersection $\ttt_ X \cap X$ is
singular at every $v' \in \ell$. This happens if and only if $\ttt_v X$
is not transversal to $\ttt_{v'} X$, if and only if the lines
$\ellfrak_v, \ellfrak_{v'}$ intersect in $\pp$. This happens if and only
if our quadric $\lquad \subset \pp$ is singular, i.e. the second case
above. But that means that all the lines $\ellfrak_{v'}$ intersect
each other (at the vertex of $\lquad$).
\end{proof}

\

\subsection{The wobbly divisor} 

We defined $\Wob$ as the union in
$X$ of the special lines. For a general line $\ell \subset X$, the
family $I_{\ell}$ of lines intersecting $\ell$ is given by $I_{\ell}
\cong C$, which comes with a trigonal map $\mathsf{a}_{\ell} :
I_{\ell} \to \ell$. According to the first pattern in
Theorem~\ref{patterns}, the $8$ branch points of $\mathsf{a}_{\ell}$
are the intersection $\ell \cap \Wob$. In particular, the class of
$\Wob$ is $8H$, where $H$ is the positive generator of $\op{Pic}(X)
\cong \zz$.

Consider the surface
\[
D := \Cbar \times \pp^1 = \left\{ \, (\ell,x) \ \left| \ {\ell
  \ \text{special}},\, x \in \pp^1 \, \right. \right\} \cong \left\{
\, (\ell,v) \ \left| \ {\ell \ \text{special}}, \, v \in \ell
\right.\right\},
\]
where $v = \mathsf{a}^{\circ}_{\ell}(x)$, using the canonical
identification $\mathsf{a}^{\circ}_{\ell} : \pp^1 \to \ell$ of each
special $\ell$ with the hyperelliptic $\pp^1$.

Let $\Gamma \subset D$ be the graph of the map $\Cbar \to \pp^1$
sending $(A,t) \to \hyp_C (t)$. For each of the $6$ Weierstrass points
$\pw_i$, let $\Gamma_i := \Cbar \times \{\brx_i\} \subset D =\Cbar
\times \pp^1$, where $\brx_i := \hyp_C(p_i)$. The involution
$\lphi_{\pw_i}: \lA \to \lA$ acts on $\Gamma_i$. In terms of the
$\pw$-dependent identification of $\Cbar$ with its image in
$\mathsf{J}$, this involution sends $A \to A^{-1}(\pw_i-\pw)$. So it
has $16$ fixed points, and the quotient $\Gamma_i / \lphi_{\pw_i}$ has
genus $5$.

\

\begin{theorem}{\label{wobbly}}
The map $\bnu: D \to X$ sending $(\ell,v) \to v$ is a finite morphism and
maps $D$ birationally onto $\Wob$, so it gives its normalization. It
is an isomorphism away from $\Gamma$ and the $\Gamma_i$. The map $\bnu$
is $2$-to-$1$ on $\Gamma_i$, with $16$ fixed points. $D$ has normal
crossings along the genus $5$ curve $\Gamma_i / \lphi_{\pw_i}$, except
at the 16 fixed points. The fixed points correspond to pattern
$4\ell$, and the other points of $\Gamma_i / \lphi_{\pw_i}$ correspond
to pattern $2\ell+2\ell_{1}$. The restriction of $\bnu$ to $\Gamma$ is
an embedding, and the image parametrizes pattern $3\ell + \ell_{1}$.
But $\bnu$ is not immersive along $\Gamma$ - the image $\Wob = \bnu(D)$
has a curve of cusps along it.
\end{theorem}
\begin{proof}
The two distinct points $(\ell,v), (\emm,v) \in D$ map to the same $v \in
X$ iff the special lines $\ell,\emm$ intersect at $v \in X$. This means
that $\emm = \lphi_{r}(\ell)$ for some $r \in C$, According to pattern
$2\ell+2\emm$ in Theorem \ref{patterns}, this occurs iff $\emm =
\lphi_{\pw_i}(\ell)$ for one of the Weierstrass points $\pw_i \in C$. This
means that we are on one of the $\Gamma_i$, and $\ell,\emm$ are related by
the involution $\lphi_{\pw_i}$.

\

Pattern $3\ell + \emm$ in Theorem \ref{patterns} shows that $\Gamma$
is the inverse image under $\bnu$ of the locus of $(\ell,v)$ such that
$(\ell,v)$ counts $3$ times among the lines through $v$. Finally we
claim that special line $\ell \subset X$ is tangent to $\bnu(\Gamma)
\subset X$ at the point $v \in \ell$ at which the patern of lines in
$X$ is $3\ell + \emm$. We will leave it to the reader to verify this
tangency in the synthetic language but we will give a modular proof of
this statement in section~\ref{app-proof-cusp}.
\end{proof}

\subsection{Synthetic correspondences and rigidification}
\label{ssec:synth.corr}

\

\punkt {\bfseries The  connection} \ 
The universal ruling $\mycal{R}$ is a $\pp^3$ bundle $\rul: \mycal{R}
\to C$ over our hyperelliptic curve $C$.

\begin{lemma}
\label{Rbar}
The pullback 
$\overline{\mycal{R}} := \mycal{R} \times _C \Cbar$
is a product,
$\overline{\mycal{R}} \cong \pp^3 \times \Cbar$.
\end{lemma}
\begin{proof}
We work with a coordinate system adapted to our pencil, so the $Q_x$
become diagonal matrices.  For $x, x_0 \in \pp^1 \setminus \{\brx_1,
\dots, \brx_6\}$, let $S:=S_{x_0,x}$ be a diagonal matrix such that
$S^2 = Q_{x_0}^{-1} Q_{x}$. Then $S$ acts on $\pp^5$, taking $Q_{x_0}$
to $Q_{x}$.  The set of these square roots $S_{x_0,x}$ has cardinality
$2^6$, and forms a torsor under the action of the group $G:=(\zz
/2)^6$.  If $t,t_0 \in C$ lie above $x, x_0$, i.e. $\hyp_C(t)=x,
\ \hyp_C(t_0)=x_0$, then half of these $S$'s take the ruling
$\mycal{R} _{t_0}$ of $Q_{x_0}$ to the ruling $\mycal{R} _{t}$ of
$Q_{x}$, while the other half take $\mycal{R} _{t_0}$ to the other
ruling of $Q_{x}$.  The set of square roots $S_{x_0,x}$ taking $t_0$
to $t$ is a torsor under the subgroup $G_0 \cong (\zz /2)^5 \subset G$
which is the kernel of the sum map $(\zz /2)^6 \to \zz /2$.  Further,
we have the diagonal embedding $\zz/2 \to G_0$. Its non-zero element
exchanges $S$ with $-S$, so it does not affect the transformation of
$\pp^5$.  The quotient $G_0/(\zz/2)$ is canonically identified with
$\mathsf{J}[2]$, the $2$-torsion subgroup of the Jacobian of $C$.  We
see that the family of quadrics $\{Q_{h_C(t)} \ | \ t \in C \setminus
\{\pw_1, \dots, \pw_6\} \ \}$ has a flat connection with monodromy
$J[2]$, and it becomes a product when pulled back to $\Cbar$ minus the
inverse image of the Weierstrass points.  The same holds therefore for
$\mycal{R}$. But since $\rho: \mycal{R} \to C$ is a $\pp^3$ bundle
(with no degenerations) over $C$, including the Weierstrass points,
the flat connection extends to all of $\mycal{R}$ ,with the same
monodromy $\mathsf{J}[2]$, and therefore the pullback to $\Cbar$ is a
product as claimed.
\end{proof}

\

A fancier way to understand this is by identifying the $\pp^3$ fiber
of $\rul : \mycal{R} \to C$ over $t \in C$ with a translated
``$2\Theta$'' linear system on $\mathsf{J}$. The right translation is
seen to be by a square root of $\Oo(t-\pw)$, which gives rise to the
$\mathsf{J}[2]$-monodromy and to the cover $\Cbar \to C$. The
projective connection we described here explicitly becomes a special
case of the general theory of theta groups and their actions, due to
Kummer, Heisenberg, Mumford etc. This approach is discussed in section
\ref{sec-heisenberg}.

\

\punkt {\bfseries The Hecke correspondence} \ We will encounter
several versions of the Hecke correspondence. The basic one is the
incidence correspondence:
\[
\bigHecke := \{ (v,\Pi) \ | \ v \in \Pi \} \subset X \times \mycal{R} 
\]
between points of $X$ and planes contained in some quadric through
$X$. Fixing $t \in C$, i.e. fixing the quadric and one of its rulings,
we have
\[
\bigHecke(t) := \bigHecke \cap (X \times \rul^{-1}(t)) = \{ (v,\Pi)
\ | \ v \in \Pi, \rul(\Pi)=t \} \ \subset \ X \times \rul^{-1}(t)
\cong X \times \pp^3.
\]
The fiber of $\bigHecke$ over $\Pi \in \mycal{R}$ (or the fiber of
$\bigHecke(t)$ over $\Pi \in \rul^{-1}(t)$) is the conic $X \cap \Pi$,
which could be smooth, a pair of intersecting lines, or a double
line. The fiber of $\bigHecke(t)$ over $v \in X$ is a line in $
\rul^{-1}(t) \cong \pp^3$. The fiber of $\bigHecke$ over $v \in X$ is
therefore a $\pp^1$-bundle $\pp E_v$ over $C$. (In the modular
approach, points $v \in X$ will correspond to rank-2 vector bundles
$E_v$ on $C$ whose determinant is $\Oo_C(\pw)$, and our $\pp E_v$ will
be the projectivization of $E_v$.)

Recall from Lemma \ref{Rbar} that the cover $ \overline{\mycal{R}} =
\mycal{R} \times _C \Cbar$ is a product, $\overline{\mycal{R}} \cong
\pp^3 \times \Cbar$.  A rigidified  Hecke correspondence is obtained by
pulling back to $\overline{\mycal{R}}$ via $\overline{\mycal{R}} 
\to \mycal{R}$:
\[
\bigHeckebar := (X\times \overline{\mycal{R}})\times_{X\times
  \mycal{R}} \bigHecke \ \subset \ X \times \overline{\mycal{R}} \cong
X \times \pp^3 \times \Cbar.
\]
From the modular point of view, we will be interested in two moduli
spaces of bundles: \ $X_1 =X$ parametrizes rank $2$ bundles on $C$ with
determinant $\Oo_C(\pw)$, while $X_0 \cong \pp^3$ parametrizes rank $2$
bundles on $C$ with determinant $\Oo_C$. The big Hecke correspondence
can then be understood as a subvariety:
\[
\bigHeckebar   \subset X_1 \times X_0 \times \Cbar.
\]
The group $\mathsf{J}[2]$ acts on all three factors. The action of
$\mathsf{J}[2] = G_0/ (\zz/2)$ on $X_1=X$ was described in the proof
of lemma \ref{Rbar}. Let $\mathsf{M}_1 := [X_1/\mathsf{J}[2]]$ denote
the quotient.  Since $X_0 \cong \pp^3$ was defined as the family of
flat sections of $ \overline{\mycal{R}} \to \Cbar$, we get an action
of $\mathsf{J}[2]$ on $X_0$ induced from its action on
$\overline{\mycal{R}}$.  Let $\mathsf{M}_0$ denote the quotient $[X_0/
  \mathsf{J}[2]]$. The action of $\mathsf{J}[2]$ on $\Cbar$ is the
familiar one, with quotient $C$. The Hecke correspondence is
compatible with these actions, so we get another version of the  Hecke
correspondence
\[
\mathbf{Hecke} \subset \mathsf{M}_1 \times \mathsf{M}_0 \times C.
\]
This is the Hecke correspondence for the group $\pp GL(2)$.

\

\

\subsection{The quadric line complex}
\label{sec-linecomplex}

Starting with the intersection
of quadrics $X$ we retrieve the curve $C$, its cover $\Cbar$, the
universal ruling $\mycal{R}$, and its cover $\overline{\mycal{R}} \to
\Cbar$.  We can therefore recover `the' $\pp^3$ as the family of flat
sections of $\overline{\mycal{R}} \to \Cbar$. We can also go in the
opposite direction. 

\

The Grassmannian $\mathsf{Grass}(2,4)$ of lines in $\pp^3 = \pp(V)$
can be identified, via the Pl\"{u}cker embedding, as a smooth quadric
in $\pp^5 = \pp (\wedge^2 V)$.  All smooth quadrics in $\pp^5$ are
isomorphic by a projective transformation, and the automorphisms of
$\pp^3$ correspond to automorphisms of $\pp^5$ preserving the quadric
and each of its rulings, via the group isomorphism $\pp GL(4)
\cong \pp SO(6)$. Thus, any smooth quadric with a specified ruling may be
viewed as being the Grassmanian $\mathsf{Grass}(2,4)$ up to adjustment
of the identification $\pp^5 = \pp (\wedge^2 V)$, in a way that is
unique up to projective transformations of $\pp^3$.

The Grassmanian has two rulings: a point $p \in \pp^3$ determines the
plane $\Pi = \Pi_p \subset \mathsf{Grass}(2,4)$ of all lines in
$\pp^3$ through $p$, while a plane $P \subset \pp^3$ determines the
plane $\Pi = \Pi_P \subset \mathsf{Grass}(2,4)$ of all lines in $P$.
Each of these rulings may be viewed as being specified, once chosen
the isomorphism of a certain quadric with the Grassmanian.

If we choose another quadric meeting $\mathsf{Grass}(2,4)$
transversally, then these two quadrics span a pencil giving the
situation at the start of this chapter. The intersection of
$\mathsf{Grass}(2,4)$ with the other quadric is known as the
\emph{\bfseries quadric line complex}. In the classical terminology cf
\cite{GH}, the ``line complex'' means the family of lines in $\pp^3$
parametrized by this quadric section of $\mathsf{Grass}(2,4)$, the
\emph{\bfseries ``lines of the quadric line complex''}.

Thus the quadric line complex is the intersection $X$ of a general pencil of
quadrics in $\pp^5$.  It carries a bit more information though:
specifying our $X$ is equivalent to specifying the hyperelliptic curve
$C$, while specifying a quadric line complex is equivalent to
specifying the hyperelliptic curve $C$ together with one
non-Weierstrass point $t \in C$: the point $\hyp_C(t) \in \pp^1$
corresponds to the quadric we identify as the Grassmannian, the point
$t \in C$ corresponds to its ruling by planes $\Pi_p$ for points $p
\in \pp^3$.

In the other direction, saying that $X$ is the quadric line complex
depends on the choice of $t\in C$ or indeed on its lifting to
$(A,t)\in \Cbar$, because the correspondence between points of $X$ and
lines in $\pp^3$ is exactly the Hecke correspondence depending on $t$,
with furthermore the identification between the ruling and a fixed
$\pp^3$ being dependent upon the lifting to $\Cbar$.

\subsection{Hecke curves and Kummer surfaces} 
\label{Hecke-curves-Kummer}

Consider the incidence:
\[
\cJ = \left\{ \, (\ell, \Pi)
\ | \ \ell \subset \Pi \, \right\}
\ \subset \ \lA \times \mycal{R}.
\]
We have an isomorphism
\[
C \times \lA \stackrel{\cong}{\to} \cJ
\]
sending $(t,\ell) \mapsto (\ell, \Pi)$, where $\Pi \subset
Q_{\hyp_{C}(t)}$ is the plane spanned by $\ell$ and $ \lphi(t,\ell)$.
The involution $\ell \mapsto \lphi(t,\ell)$ acts on
$\cJ$, with quotient:
\[
\cKum := \left\{ \, \Pi \in \mycal{R} \ | \ \Pi \cap X
\ \text{contains a line} \, \right\}.
\]
We have maps
\[
\xymatrix@1@M+0.5pc@C+0.7pc{
  C \times \lA \ \cong \
  \cJ \ar[r]^-{2:1} & \cKum \ar@{^{(}->}[r] & \mycal{R},
  }
\]
or, if we fix $t \in C$:
\[
\lA \ \stackrel{2:1}{\to} \ \Kum \ \hookrightarrow \ \pp^3.
\]
Here $\Kum$ is the Kummer surface, image of $\lA \cong \mathsf{J}$ by
the $2\Theta$ linear system, which embeds $\Kum = \lA/ \lphi_t \cong
\mathsf{J}/ {\pm 1} $ into $\pp^3$ as a singular hypersurface of degree $4$.  
Note that the Kummer surface $\Kum
\subset \pp^3$ has $16$ nodes, images of the $16$ fixed points of the
involution $\lphi_t : \lA \to \lA$. These are the points where the fiber
$\ell \cup \lphi_t (\ell)$ becomes a double line.

A basic
fact is that $\Kum$ is isomorphic to its dual \cite{GH,Keum}.

\

Given a point $p\in \pp^3=\rul^{-1}(t)$, the fiber of the Hecke correspondence over
$p$ is a subset of $X$ that was denoted by
$X_p$ in \cite{GH}. In the terminology of the quadric line complex  
it is the set of lines in the line complex that pass through $p$.

\

\begin{corollary} 
\label{Hecke-curve-trichotomy}
Fix $t \in C$ and denote $\pp^3:= \rul^{-1}(t)$. 
\begin{itemize}
\item If $p$ is a
  point of $\pp^3 \setminus \Kum$, then $X_p$ is a smooth conic $\Pi \cap X$.
\item If $p$ is a smooth point of
  $\Kum$, then $X_p$ is a union of two lines $ \ell \cup \lphi(t,\ell)$ that touch.
\item If $p$ is one of the $16$ nodes of
  $\Kum$, then $X_p$ is a double line $2\ell$.
\end{itemize}
\end{corollary}

\

This was stated in \cite[pp 762-763]{GH}.

\

\noindent
The locus of points in $X$ that are intersection points of the two
lines, for the second case of $X_p$, has closure that is a surface
denoted $\Sigma \subset X$ in \cite{GH}.  This is seen to be the
\emph{\bfseries Kummer K3 surface} obtained by resolving the $16$ nodes of
$\Kum \subset \pp^3$ (indeed, over the nodes we get not a single point
but a line of points in $\Sigma$ corresponding to the full double
line).

The K3 surface will show up in our situation as a locus of apparent
singularities for the Hecke transform from $X_0 = \pp^3$ to $X_1=X$. For the
main local systems constructed here, it turns out that the
singularities along $\Sigma$ are removable. On the other hand for the
third construction in Chapter \ref{chapter-thirdconstruction}, the
Hecke transform will have singularities along $\Sigma$. We refer the
reader to that chapter for more discussion.

In the synthetic picture, the trope planes may be characterized as
planes in $\pp^3$ that contain $6$ of the $16$ nodes. There are $16$ of
these. They are also the planes that correspond to the $16$ nodes of
the dual Kummer surface.  The trope planes meet the Kummer surface in
plane conics counted with multiplicity two; these conics are
characterized also by passing through the $6$ nodes that define the
plane.

The $16_6$ property, saying that each plane contains $6$ nodes (by
definition) and each node is contained in $6$ planes, has been well-known since the 19th century.
Modern discussions can be found in
\cite{GH,BeauvilleSurfaces,Keum,Dolgachev200}. 
\

\

\section{General considerations}
 \label{chapter-general}

Throughout this paper, we consider a smooth projective curve $C$ of
genus $2$. It is therefore hyperelliptic, with hyperelliptic
involution denoted $\hi_{C} : C\rightarrow C$, and the quotient
projection denoted by $\hyp_{C} : C\rightarrow \pp^1$. Let
$\pw_1,\ldots , \pw_6\in C$ denote the Weierstrass points,
i.e. the ramification points for the map $\hyp_{C}$, and let
$\brx_1,\ldots , \brx_6\in \pp ^1$ denote the corresponding branch
points. For brevity we will choose one of the Weierstrass points and
drop the index, writing $\pw :=\pw_1$.  Thus the canonical line
bundle of $C$ is given by $\omega_C=\Oo _C(2\pw)$.

If $x\in C$ we will usually denote by $x':= \hi_{C}(x)$ its image by
the hyperelliptic involution; then $\Oo _C(x+x')\cong \Oo _C(2\pw)=
\omega_C$.

Fix a line bundle $\mathbf{d}$ to be either $\Oo _C$ or $\Oo _C(\pw)$.
We will look at the moduli space $X$ of polystable rank $2$ vector
bundles $E$ provided with an isomorphism $\det(E)\cong \mathbf{d}$. The
choice of $\mathbf{d}$, that is to say degree $0$ or $1$, will be made
according to section in the paper---we prefer not to overload the
notations by indexing on this choice. When it becomes necessary to
distinguish them, we will denote the two moduli spaces by $X_0$ and
$X_1$ and similarly for their associated constructions (e.g. $Y_0$,
$Y_1$).

By Narasimhan-Ramanan \cite{NR}, for degree $0$ we have
$X=\pp^3$ and for degree $1$ we have $X\subset \pp^5$ is a smooth
complete intersection of two quadrics. These descriptions will be
developed in more detail below.

Let $\Higgs$ denote the coarse moduli space of semistable Higgs
bundles $(E,\theta)$ of rank $2$ with isomorphism $\det(E)\cong
\mathbf{d}$ and satisfying $\op{tr} \theta =0$.  In particular,
$X\subset \Higgs$ is the subset of Higgs bundles with Higgs field
equal to zero.  The {\em\bfseries Hitchin map} on rank $2$ Higgs
bundles was defined by Hitchin \cite{Hitchin-selfd,Hitchin-spectral}
and is hiven by
\[
\hit : \Higgs \rightarrow \Bb := H^0(\omega_C^{\otimes 2}), \quad
(E,\theta) \mapsto \det(\theta)
\]
The space of quadratic differentials $\Bb := H^0(\omega_C^{\otimes
  2})$ is called the {\em\bfseries Hitchin base} and in our case is a
complex vector space of dimension $3$.

We will write $T^{\vee}C$ for the total space of
$\omega_{C}$. $\boldsymbol{\pi} : T^{\vee}C \to C$ for the natural
projection, and $\lambda \in
H^{0}(T^{\vee}C,\boldsymbol{\pi}^{*}\omega_{C})$ for the tautological
section.  For any quadratic differential $\mathsf{b} \in \Bb$, the
associated {\em\bfseries spectral curve} is defined as the curve
$\Ctilde_{\mathsf{b}} \subset T^{\vee}C$ given by the equation
\[
\Ctilde_{\mathsf{b}} \, : \, \lambda^{2}  - \boldsymbol{\pi}^{*}\mathsf{b} = 0.
\]
In other words, the spectal curve correcsponding to $\mathsf{b}$ is
the zero divisor of the holomorphic section $\lambda^{2} -
\boldsymbol{\pi}^{*}\mathsf{b} \in
H^{0}(T^{\vee}C,\boldsymbol{\pi}^{*}\omega_{C})$.  By definition, the
spectral curve associated to a Higgs bundle $(E,\theta)$ is the
spectral curve corresponding to $\hit(E,\theta)$. It is given explicitly
the equation
\[
\det(\lambda\cdot \op{id} - \boldsymbol{\pi}^{*}\theta) =
\lambda^{2} - \det(\theta)  = 0,
\]
where $\lambda\cdot \op{id} - \boldsymbol{\pi}^{*}\theta$ is viewed as a map
\[
\lambda\cdot \op{id} - \boldsymbol{\pi}^{*}\theta : \, 
\boldsymbol{\pi}^{*}E \longrightarrow \boldsymbol{\pi}^{*}E\otimes
\boldsymbol{\pi}^{*}\omega_{C}.
\]

For future reference we record the following

\begin{proposition} \label{prop:hitchin.map} Let $\Higgs$ be the moduli
  space of rank $2$ Higgs bundles with fixed determinant $\mathbf{d}$.  Then
\begin{itemize}
\item[(a)] The Hitchin map $\hit : \Higgs \to \Bb$ is proper and  surjective.
\item[(b)] For any $(E,\theta)\in \Higgs$ the spectral curve
  $\Ctilde \subset T^{\vee}C$ of $\theta$ is a  curve of arithmetic genus $5$
  and the projection 
$$
\pi = \boldsymbol{\pi}_{|\Ctilde}  : \Ctilde \rightarrow C
$$ is the degree $2$ covering branched over the bicanonical divisor
$\det(\theta) = 0$. For a general $(E,\theta)$ the spectral curve
$\Ctilde$ is smooth and connected.
\item[(c)] For any $\mathsf{b} \in \Bb$, the branch divisor of the
  associated spectral curve $\pi : \Ctilde_{\mathsf{b}} \to C$ is the
  $\hyp_{C}$- pullback of a degree $2$ effective divisor $\bry + \brz$ on
  $\mathbb{P}^{1}$. For a general $\mathsf{b}$ the branch divisor
  $\mathsf{zero}(\mathsf{b}) = \hyp_{C}^{-1}(\bry + \brz) = \tilde{y} +
  \tilde{y}' + \tilde{z} + \tilde{z}'$ consists of two
  $\hi_{C}$-conjugate pairs of points in $C$.
\item[(d)] For general $\mathsf{b} \in \Bb$ the Hitchin fiber
  $\hit^{-1}(\mathsf{b})$ is identified with the $3$-dimensional abelian
  variety $\Prym$ of line bundles $L$ on $\Ctilde $ such that
  $\det(\pi_{*}(L))\cong \mathbf{d}$, or equivalently
  $\mathsf{Nm}_{\pi}(L) = \omega_{C}\otimes \mathbf{d}$. The Higgs bundle
  corresponding to $L$ is $(E,\theta) = (\pi _{*}L,
  \pi_{*}(\lambda\otimes (-)))$.
\end{itemize}
\end{proposition}
\begin{proof}
Parts (a), (b), and (d) are standard and are proven in many classical
sources, e.g.  \cite{Hitchin-selfd,Hitchin-spectral,BNR}. For part (c)
we only need to note that $\omega_{C} = \hyp_{C}^{*}\mathcal{O}(1)$,
and that $\hyp_{C}^{*} : H^{0}(\mathbb{P}^{1},\mathcal{O}(2)) \to
H^{0}(C,\omega_{C}^{\otimes 2})$ is injective, Since both spaces are
$3$-dimensional this shows that the pullback map is an isomorphism,
and so any quadratic differential on $C$ is a pullback from an
effective degree $2$ divisor on $\mathbb{P}^{1}$. This proves (c) and
completes the proof of the proposition.
\end{proof}

\

Throughout the paper, we will fix a point $\mathsf{b}$ and the
corresponding spectral curve $\Ctilde_{\mathsf{b}}$. Since
$\mathsf{b}$ will be fixed we will drop the subscript from the
notation and will simply write $\Ctilde$ finstead of
$\Ctilde_{\mathsf{b}}$. We will also fix $\Lambda$ the rank $2$ local
system corresponding to a rank $2$ Higgs bundle with trivial
determinant $(E,\theta )$.

\subsection{Curves}
\label{subsec-curves}

Let us make the following notations and definitions. 

\begin{description}
\item[$C$] is a curve of genus $g(C) = 2$ with hyperelliptic 
  map $\hyp_{C} : C \rightarrow \pp^1$. 

\item[$\pw\in C$]  is a fixed Weierstrass point.

\item[$\Cbar$] $:= \left\{ (A,t) \ \left| A\in \op{Jac}^0(C),\;\;
  t\in C, \;\; A^{\otimes 2} = \Oo _C (t-\pw) \right.\right\}$.

  Note that the map
  \[
  \sq : \Cbar \to C, \quad \sq(A,t) = t
  \]
  is a $16$-sheeted
  \'{e}tale cover of $C$, while the map
  \[
  \emb_{\Cbar} : \Cbar
  \hookrightarrow \op{Jac}^{0}(C), \quad \emb_{\Cbar}(A,t) = A
  \]
  is a closed embedding. In particular
  $\Cbar$ is a smooth connected curve of genus $g(\Cbar) = 17$.

\item[$\pi : \Ctilde \rightarrow C$] is a fixed spectral curve,
  assumed to be smooth. We also assume 
  that the spectral curve has two branches over $\pw$. 
  
\item[$\Chat := \Ctilde \times _C \Cbar$] is a smooth connected curve
  of genus $g(\Chat) = 65$ which fits in the fiber square
  \[
  \xymatrix{
    \Chat \ar[r]^-{\pihat} \ar[d]_-{\sqhat} & \Cbar \ar[d]^-{\sq} \\
    \Ctilde \ar[r]_-{\pi} & C
   }
  \]
\end{description}

Suppose $(E,\theta )$ is a Higgs bundle with trivial
determinant on $C$, i.e. $\det (E) \cong \Oo _C$ and $tr(\theta )= 0$.
Our main condition is that $\Ctilde$ is the spectral cover of
$(E,\theta)$, this means that there is a line bundle $L$ on $\Ctilde$
such that $E = \pi_{*}L$ and $\theta = \pi_{*}(\lambda\otimes(-))$.

The choice of $(E,\theta)$ will determine
(via the construction we are going to do,
using its spectral data
$(\Ctilde,L)$) two tame parabolic Higgs bundles of rank $8$,
\begin{equation} \label{eq:eigensheaves}
( \srF _{0,\bullet}, \Phi _0) / X_0, \;\;\;\; ( \srF _{1,\bullet},
\Phi _1) / X_1.
\end{equation}
Our goal is to give a detailed construction of these Higgs bundles and
to show that they satisfy the Hecke eigensheaf property on
$X_{0}\sqcup X_{1}$.

\subsection{Moduli spaces} \label{ssec:moduli}

The coarse moduli spaces of semistable bundles will be denoted by $X$
with, if necessary, a subscript depending on the degree. These have
coverings destined to become the modular spectral covers of the parabolic
Higgs bundles we are going to construct. The notations are as follows.

\begin{description}
\item[$X_0$] is the moduli space of rank $2$ bundles $F$ with $\det
  (F) = \Oo _C$.  We have $X_0 \cong \pp^3$ \cite{NR}.
\item[$X_1$] is the moduli space of rank $2$ bundles $E$ with $\det
  (E) = \Oo _C(p)$.  We have \cite{NR,Newstead}
$$
X_1 \subset \pp ^5, \;\;\;\; X_1 = \bigcap _{x\in \pp^1} Q_x,
$$
where $\{Q_{x}\}_{x \in \mathbb{P}^{!}}$ is a pencil of quadrics in
$\mathbb{P}^{5}$ parametrized by the hyperelliptic $\mathbb{P}^{1}$ of
$C$, with a discriminant divisor being exactly the ramification divisor
of the hyperelliptic map $\hyp_{C} : C \to \mathbb{P}^{1}$. 
\item[$\Prym_2$] is the degree two Prym for $\Ctilde$, that is to say
$$
  \Prym_2 := \left\{\left. L \in \op{Jac}^2(\Ctilde )\;\; \right| \;\;
  \mathsf{Nm}_{\pi}(L)
  \cong \omega_C \right\} .
  $$
This is identified with the Hitchin fiber over $[\lambda : \Ctilde
    \hookrightarrow T^{\vee}C] \in \Bb$ (corresponding point in the
  Hitchin base) for the $SL(2)$ Hitchin fibration $\hit : \Higgs_{0}
  \to \Bb$.
\item[$\Prym_3$] is the degree three Prym:
$$ 
  \Prym_3 := \left\{\left. L \in \op{Jac}^3(\Ctilde)\;\; \right|
  \;\; \mathsf{Nm}_{\pi}(L)
  \cong \omega _C (\pw) \right\} .
  $$
It is the Hitchin fiber over $[\lambda : \Ctilde \hookrightarrow
    T^{\vee}C] \in \Bb$ for the odd degree Hitchin fibration $\hit :
  \Higgs_{1} \to \Bb$ (Higgs bundles whose determinant Higgs bundle is
  $(\Oo _C(\pw), 0)$).
\item[$Y_0$] is the blow-up of $\Prym_2$ at the pullbacks of the theta
  characteristics on  $C$. More precisely, consider the un(semi)stable locus
  in $\Prym_{2}$,
  $$
\begin{aligned}
  \Prym_{2}^{\rm unss} 
  & = \left\{ L \in \Prym_{2} \;\; \left| \;\; \pi_{*}L
    \text{ is un(semi)stable } \right.\right\} \\ 
& = \left\{  \pi^{*}\kappa \in \Prym_{2}
    \;\; \left| \;\; \kappa \in \op{Jac}^{1}(C), \text{ s.t. }
    \kappa^{\otimes 2} = \omega_{C} \right.\right\}
    & = \pi^{*}\mathsf{Spin}(C)
\end{aligned}
$$
which has $16$ points, and $\blo_{0}: Y_{0} \rightarrow \Prym_{2}$ is the
blow-up of $\Prym_{2}^{\rm unss} \cong \mathsf{Spin}(C)$.  
\item[$\ExY_0\subset Y_0$] denotes
the exceptional divisor, it has $16$ connected components
$$
\ExY_0 = \sqcup _{\kappa} \ExY_{0,\kappa}.
$$
\item[$f_0$] is a morphism $f_0 : Y_0 \rightarrow X_0$ such that if
  $L\in \Prym_2$ is a point not on $\Prym_2^{\rm unss} $ we have
  $f_0(L)= \pi _{*}(L)$. In Theoroem~\ref{qresolves} the morphism
  $f_{0}$ is constructed as the minimal resolution of the rational map
  $\pi_{*}(-) : \Prym_{2} \dashrightarrow X_{0}$. From the
  construction it follows that for a general spectral curve $f_{0} :
  Y_{0} \to X_{0}$ is finite and by Lemma~\ref{deg8c1} and
  Remark~\ref{rem:deg8X0} has degree $8$.
  \item[$Y_1$] is the blow-up of $\Prym_3$ in the un(semi)stable
    ($=$unstable) locus 
$$
\Prym_3^{\rm unss} = \left\{ L \in \Prym_{3} \;\; \left| \;\; \pi_{*}L
\text{ is un(semi)stable } (=\text{is unstable})
\right.\right\} \cong \Chat.
$$
To see that the unstable locus is related to $\Chat$, suppose we have
a line bundle $L \in \Prym_{3}$ such that $\pi_{*}L$ is unstable. Then
$\pi_{*}L$ will have a destabilizing line subbundle $M \subset
\pi_{*}L$ of degree one. By adjunction we get an injective map of
locally free rank one sheaves $\pi^{*} M \hookrightarrow L$ and for
degree reasons we must have a short exact sequence
\[
0 \to \pi^{*}M \to
L \to \mathcal{O}_{\tilde{t}} \to 0
\]
for some point $\tilde{t} \in
\Ctilde$. Writing $M = A^{-1}(\pw)$ for some line bundle of degree
zero, we see that $\mathsf{Nm}_{\pi}(L)\otimes \mathcal{O}_{C}(-
\pi(\tilde{t})) = \mathsf{Nm}_{\pi}(\pi^{*}A^{-1}(\pw)) =
A^{-2}(2\pw)$. Since $\mathsf{Nm}_{\pi}(L) = \mathcal{O}_{C}(3\pw)$ we
have that $A^{\otimes 2}(\pw) = \mathcal{O}_{C}(\pi(\tilde{t}))$. Thus
$(A,\tilde{t}) \in \Chat$ and   $L = \pi^{*}(A^{-1}(\pw))(\tilde{t})$. 
In Lemma~\ref{lemma:embedChat} below  we show that the map $\Chat \to
\Prym_{3}$, $(A,\tilde{t}) \mapsto \pi^{*}(A^{-1}(\pw))(\tilde{t})$ is
a closed embedding which gives the identification $\Prym_{3}^{\rm
  unss} = \Chat$.

Again we will write $\blo_1:Y_1\rightarrow \Prym_3$ for the blow-up
map.

\item[$f_1$] is a morphism $f_1 : Y_1 \rightarrow X_1$ such that if
  $L\in \Prym_3$ is a point not on $\Prym_3^{\rm unss}$ we have \linebreak 
  $f_1(L)= \pi _{*}L$. In Theoroem~\ref{qresolves} the morphism
  $f_{1}$ is constructed as the minimal resolution of the rational map
  $\pi_{*}(-) : \Prym_{3} \dashrightarrow X_{1}$ and by construction is finite for a general spectral cover. By Lemma~\ref{deg8c1}  $f_{1}$ has  degree $8$.
\item[${\ExY}_1\subset Y_1$] denotes the exceptional divisor, it is a
  $\mathbb{P}^{1}$-bundle
$$
  \ExY_1 = \mathbb{P}(N_{\Chat/\Prym_{3}}) \rightarrow \Chat. 
$$
In Theoroem~\ref{qresolves} we will see that in fact
$$
\ExY_1\cong  \Chat \times \pp^1,
$$
where the second factor  is naturally idenitified with the hyperelliptic
$\pp^1$ of $C$, that is it is identified with the projective line
$\mathbb{P}(H^{0}(C,\omega_{C})^{\vee})$.  We have a diagram
$$
\xymatrix{
\ExY_1 \ar[d] \ar@{^{(}->}[r] & Y_1 \ar[d] \\
\Cbar \times \pp^1 \ar[r] & X_1 
}
$$
and in section~\ref{chapter-d1} we will see that the map on the
bottom factors as
$$
\Cbar \times \pp^1 \rightarrow \Wob _1  \hookrightarrow  X_1 
$$
where the surface $\Wob_{1}$ is the wobbly divisor discussed in the
next section, and the map $\Cbar \times \pp^1 \rightarrow \Wob _1$ is
the normalization.

To understand the map 
\begin{equation} \label{eq:wobnorm}
\Cbar\times \pp^{1} = \Cbar\times \mathbb{P}(H^{0}(C,\omega_{C})^{\vee})
\to X_{1}.
\end{equation}
in more concrete terms fix points $(A,t) \in \Cbar$ and $x \in
\mathbb{P}(H^{0}(C,\omega_{C})^{\vee})$. Any non-split extension of
$A^{-1}(\pw)$ by $A$ will be a stable rank two bundle with determinant
$\mathcal{O}_{C}(\pw)$, i.e. will give us a point in $X_{1}$. However
the space of extensions of $A^{-1}(\pw)$ by $A$ is canonically
identified with the space $H^{0}(C,\omega_{C})^{\vee}$. Indeed we have
\begin{equation} \label{eq:exttriv}
\begin{aligned}
\op{Ext}^{1}(A^{-1}(\pw),A) & = H^{1}(C,A^{\otimes 2}(-\pw)) =
H^{1}(C,\mathcal{O}_{C}(t-2\pw)) \\ & =
H^{1}(C,\mathcal{O}_{C}(-t')) = H^{1}(C,\mathcal{O}_{C}) \\
& =
H^{0}(C,\omega_{C})^{\vee},
\end{aligned}
\end{equation}
where $t'$ is the hyperelliptic conjugate of $t$, and the
identification $H^{1}(C,\mathcal{O}_{C}(-t')) =
H^{1}(C,\mathcal{O}_{C})$ is induced by the natural inclusion of
sheaves $\mathcal{O}_{C}(-t') \subset \mathcal{O}_{C}$.

With this picture in mind we can now describe the map
\eqref{eq:wobnorm}. It sends a point $((A,t),x) \in \Cbar\times
\mathbb{P}(H^{0}(C,\omega_{C})^{\vee})$ to the vector bundle $E \in
X_{1}$, where $E$ is defined as the unique up to isomorphism extension
\[
0 \to A \to E \to A^{-1}(\pw) \to 0
\]
which corresponds to the extension class $x \in
\mathbb{P}(H^{0}(C,\omega_{C})^{\vee}) \cong
\mathbb{P}(\op{Ext}^{1}(A^{-1}(\pw),A))$ under the identification
\eqref{eq:exttriv}.

\

Finally, note that the curve $\Cbar$ embeds in $\Cbar\times \pp^1$ as
the graph of the  composition
\[
\xymatrix@1@M+0.5pc{
\Cbar \ar[r]^-{\sq} & C \ar[r]^-{\hyp_{C}} & \pp^1,
  }
\]
with $\hyp_{C}$ denoting the hyperelliptic map. It is easy to check
that restricting the map \eqref{eq:wobnorm} to this embedded copy of
$\Cbar$ yields a closed embedding
\[
\xymatrix@M+0.5pc{
  \ql : \hspace{-3.5pc}  & \Cbar \ar@{^{(}->}[r]  & X_{1} \ \subset
  \mathbb{P}^{5}}
\]
of $\Cbar$ in the intersection of two quadrics $X_{1}$.
\end{description}

\

\noindent
We conclude this subsection with the promised check that the curve
$\Chat = \Ctilde\times_{C} \Cbar$ embeds in $\Prym_{3}$.

\

\begin{lemma} \label{lemma:embedChat}
  The map
\begin{equation} \label{eq:ChatinY1}
\emb_{\Chat} : \Chat \to \Prym_{3}, \quad (A,\tilde{t}) \mapsto
\pi^{*}(A^{-1}(\pw))\otimes \mathcal{O}_{\Ctilde}(\tilde{t}).
\end{equation}
  embeds the curve $\Chat$ inside
$\Prym_{3}$. In particular $\ExY_{1} = \mathbb{P}(N_{\Chat/\Prym_{3}})$.
\end{lemma}
\begin{proof}
  Suppose $(A_{1},\tilde{t}_{1}), (A_{2},\tilde{t}_{2}) \in \Chat$. If
  these two points map to the same point in $\Prym_{3}$, then we will
  have
  \begin{equation} \label{eq:equalpts}
    (\pi^{*}A_{1}^{-1})(\tilde{t}_{1}) =  (\pi^{*}A_{2}^{-1})(\tilde{t}_{2}).
  \end{equation}
Squaring this identity and using the fact that
  $A_{i}^{-2}(\pi(\tilde{t}_{i})) = \mathcal{O}_{C}(\pw)$ for $i = 1,2$  we get 
\[
\begin{aligned}
  \mathcal{O}_{\Ctilde} & = (\pi^{*}A_{2}^{2})( - 2\tilde{t}_{2})
  \otimes (\pi^{*}A_{1}^{-2})(2\tilde{t}_{1}) \\
& = 
  \left(\pi^{*}(A_{2}^{2}(-\pi(\tilde{t}_{2}))\right)
  \left(\pi^{*}(\pi(\tilde{t}_{2}) - 2\tilde{t}_{2}\right)
  \otimes \left(\pi^{*}A_{1}^{-2}(\pi(\tilde{t}_{1})\right)
  \left(-\pi^{*}(\pi(\tilde{t}_{1}) + 2\tilde{t}_{1}\right) \\
  & =
  \mathcal{O}_{\Ctilde}(\tilde{t}_{1} + \tau(\tilde{t}_{2}) -
  \tau(\tilde{t}_{1}) - \tilde{t}_{2}),
\end{aligned}
\]
where, as usual, $\tau : \Ctilde \to \Ctilde$ denotes the covering
involution for the map $\pi : \Ctilde \to C$.

Thus either $\tilde{t}_{1} + \tau(\tilde{t}_{2})$ and
$\tau(\tilde{t}_{1}) + \tilde{t}_{2}$ are equal as divisors or
$\tilde{t}_{1} + \tau(\tilde{t}_{2})$ and $\tau(\tilde{t}_{1}) +
\tilde{t}_{2}$ span a $g^{1}_{2}$ linear system on $\Ctilde$, that is
$\tilde{t}_{1} + \tau(\tilde{t}_{2})$ and $\tau(\tilde{t}_{1}) +
\tilde{t}_{2}$ are disjoint divisors in the hyperelliptic linear
system on $\Ctilde$.  This gives the following possibilities

\

\noindent
{\bfseries Case 1.} \ We have $\tilde{t}_{1} = \tilde{t}_{2}$. In this
case the equality \eqref{eq:equalpts} implies $\pi^{*}A_{1} =
\pi^{*}A_{2}$. Since $\pi$ is ramified, the pullback $\pi^{*} :
\op{Jac}^{0}(C) \to \op{Jac}^{0}(\Ctilde)$ is injective, and so $A_{1}
= A_{2}$.

\

\noindent
{\bfseries Case 2.} \ We have $\tilde{t}_{1} = \tau(\tilde{t}_{1})$
and $\tilde{t}_{2} = \tau(\tilde{t}_{2})$, and $\tilde{t}_{1} \neq
\tilde{t}_{2}$, i.e.  $\tilde{t}_{1}$ and $\tilde{t}_{2}$ are two
distinct ramification points of $\pi : \Ctilde \to C$. This 
violates the identity \eqref{eq:equalpts} which 
implies that
\[
\mathcal{O}_{\Ctilde}(\tilde{t}_{1} - \tilde{t}_{2}) =
\pi^{*}\left(A_{1}\otimes A_{2}^{-1}\right).
\]
However the line bundle $\mathcal{O}_{\Ctilde}(\tilde{t}_{1} -
\tilde{t}_{2})$ can not be a pullback of a line bundle from $C$.
Indeed, every line bundle $\mycal{L}$ on $\Ctilde$ which is a pullback
from a line bundle on $C$ admits a $\tau$-equivariant structure which
acts trivially on the fibers of $\mycal{L}$ at all ramification points
of $\pi : \Ctilde \to C$. On the other hand, if $\tilde{r}$ is a fixed
point of $\tau$, then $\tau$ preserves the ideal subsheaf
$\mathcal{O}_{\Ctilde}(-\tilde{r})$ and so the locally free sheaf
$\mathcal{O}_{\Ctilde}(-\tilde{r})$ is equipped with a canonical
$\tau$-equivariant structure in which $\tau$ acts as multiplication by
$(-1)$ on the fiber of $\mathcal{O}_{\Ctilde}(-\tilde{r})$ at
$\tilde{r}$ and acts as multiplication by $1$ on the of
$\mathcal{O}_{\Ctilde}(-\tilde{r})$ at any fixed point different from
$\tilde{r}$. By duality $\mathcal{O}_{\Ctilde}(\tilde{r})$ has an
equivariant structure with the same exact property and so by tensoring
we see that for two distinct ramification points $\tilde{t}_{1}$ and
$\tilde{t}_{2}$ we see that $\mathcal{O}_{\Ctilde}(\tilde{t}_{1} -
\tilde{t}_{2})$ has a $\tau$-equivariant structure in which $\tau$
acts by multiplication by $(-1)$ on the fibers at $\tilde{t}_{1}$ and
$\tilde{t}_{2}$ and by multiplication by $1$ on the fibers at the
other two ramification points. But the only other $\tau$-equivariant
structure on this line bundle will be obtained by multiplying the
given equivariant structure by the sign charater $\langle \tau \rangle
\to \mathbb{C}^{\times}$, $\tau \mapsto -1$. In either of these
structures $\tau$ acts non-trivially on the fibers of
$\mathcal{O}_{\Ctilde}(\tilde{t}_{1} - \tilde{t}_{2})$ at two of the
ramification points, and so this line bundle is not a pullback of a
line bundle on $C$. This shows that in this case equation
\eqref{eq:equalpts} does not have a solution.

\

\noindent
{\bfseries Case 3.} \ We have that $\tilde{t}_{1} +
\tau(\tilde{t}_{2})$ and $\tilde{t}_{2} + \tau(\tilde{t}_{1})$ are
disjoint divisors in the hyperelliptic linear system on $\Ctilde$. To
analyze this case better we first recall the basic diagram governing
the geometry of the curve $\Ctilde$. To construct the spectral curve
$\Ctilde$ we start with a quadratic differential $\beta \in
H^{0}(C,\omega_{C}^{\otimes 2})$ having simple zeroes and we take
$\Ctilde \subset \op{tot}(\omega_{C})$ to be the unique double cover
of $C$ branched at the zeroes of $\beta$. Let $\hyp_{C} : C \to
\mathbb{P}^{1}$ be the hyperelliptic map. Since $\omega_{C} =
\hyp_{C}^{*} \mathcal{O}_{\mathbb{P}^{1}}(1)$ and it is easy to see
that every section of $\omega_{C}^{\otimes 2}$ is a pullback of a
unique section in $\mathcal{O}_{\mathbb{P}^{1}}(2)$.  Thus the divisor
of $\beta$ is the pullback of a degree two divisor $\bry + \brz$ in
$\mathbb{P}^{1}$, where $\bry, \brz \in \mathbb{P}^{1}$ are two distinct
points, neither of which is a branch point of $\hyp_{C}$. In
particular, the degree $4$ cover $\hyp_{C}\circ\pi : \Ctilde \to
\mathbb{P}^{1}$ factors as $\Ctilde \to \mathbb{P}^{1} \to
\mathbb{P}^{1}$, where $\mathbb{P}^{1} \to \mathbb{P}^{1}$ is the
double cover branched at $\bry + \brz$, and $\Ctilde \to \mathbb{P}^{1}$ is
the hyperelliptic map on $\Ctilde$ which we will denote by
$\hyp_{\Ctilde}$. This also implies that the covering involution $\tau
: \Ctilde \to \Ctilde$ for $\pi : \Ctilde \to C$ commutes with the
hyperelliptic involution $\sigma : \Ctilde \to \Ctilde$ and that the
composition $\rho = \sigma\circ\tau$ is a fixed point free involution
with quotient $D = \Ctilde/\langle \rho \rangle$ which is a smooth
hyperelliptic curve of genus $3$. All these data can be organized in
the commutative diagram
\begin{equation} \label{eq:hyp.diagram}
  \begin{minipage}[c]{4in}
    \[
\xymatrix@M+1pc{
  & \Ctilde \ar[ld]_-{\hyp_{\Ctilde}} \ar[rd]^-{\pi} \ar[d]_-{\mathsf{a}} & \\
  \mathbb{P}^{1} \ar[rd]_-{\hyp_{\mathbb{P}^{1}}}  & D
  \ar[d]_-{\hyp_{D}} &
  C \ar[ld]^-{\hyp_{C}} \\
  & \mathbb{P}^{1} &
 }
\]
\end{minipage}
\end{equation}
where the maps $\hyp_{\Ctilde}$, $\mathsf{a}$, and $\pi$ have covering
involutions $\sigma$, $\rho$, and $\tau$ respectively.

Therefore the divisor $\tilde{t}_{1} + \tau(\tilde{t}_{2})$ is a
fiber of the map $\hyp_{\Ctilde}$ if and only if $\tilde{t}_{1} =
\sigma\circ \tau (\tilde{t}_{2}) = \rho(\tilde{t}_{2})$. So in this case the
identity \eqref{eq:equalpts}  becomes
\begin{equation} \label{eq:case3identity} 
\mathcal{O}_{\Ctilde}(\tilde{t}_{1} - \rho(\tilde{t}_{1})) =
\pi^{*}(A_{1}\otimes A_{2}^{-1}).
\end{equation}
Note that the line bundle $\mathcal{O}_{\Ctilde}(\tilde{t}_{1} -
\rho(\tilde{t}_{1}))$ is a $\rho$-antiinvariant line bundle on
$\Ctilde$, i.e. as a point in $\op{Jac}^{0}(\Ctilde)$ it belongs to the
subgroup
\[
\ker\left[\xymatrix@1{\op{Jac}^{0}(\Ctilde)
    \ar[r]^-{1+\sigma} & \op{Jac}^{0}(\Ctilde)}\right] =
\ker\left[\xymatrix{\op{Jac}^{0}(\Ctilde)
    \ar[r]^-{\mathsf{Nm}_{\mathsf{a}}} &  \op{Jac}^{0}(C)}\right] \ \subset
  \op{Jac}^{0}(\Ctilde).
  \]
 From the classical theory of Prym varieties
\cite{Mumford-Prym} of \'{e}tale double covers, it is known that $\ker
(1 + \rho) = \ker \mathsf{Nm}_{\mathsf{a}}$ is a disconnected abelian
subgroup with two connected components, and a connected component of
the identity equal to the (degree zero) Prym variety for the pair
$(\Ctilde,D)$:
\[
\mathsf{Prym}(\Ctilde,D) =
\op{im}\left[\xymatrix@1{\op{Jac}^{0}(\Ctilde) \ar[r]^-{1-\sigma} &
    \op{Jac}^{0}(\Ctilde)}\right].
\]
Next consider the Abel-Prym map
\[
\mathfrak{ap}_{(\Ctilde,D)}^{k} : \op{Sym}^{k}\Ctilde \longrightarrow
\ker \mathsf{Nm}_{\mathsf{a}} \subset \op{Jac}^{0}(\Ctilde),
\ \mathfrak{d} \mapsto \mathcal{O}_{\Ctilde}(\mathfrak{d} -
\rho(\mathfrak{d})).
\]
As explained in \cite{Mumford-Prym} the image of this map lands in the
identity component
\[
\mathsf{Prym}(\Ctilde,D) = \left(\ker
\mathsf{Nm}_{\mathsf{a}}\right)_{o} \subset \ker
\mathsf{Nm}_{\mathsf{a}}
\]
if and only if $k$ is even. Thus
the image of Abel-Prym map
\[
\mathfrak{ap}_{(\Ctilde,D)}^{1} : \Ctilde \longrightarrow
\ker \mathsf{Nm}_{\mathsf{a}} \subset \op{Jac}^{0}(\Ctilde),
\ \tilde{t} \mapsto \mathcal{O}_{\Ctilde}(\tilde{t} -
\rho(\tilde{t})),
\]
is contained in the non-neutral component of $\ker
\mathsf{Nm}_{\mathsf{a}}$ and so the curve
$\mathfrak{ap}_{(\Ctilde,D)}^{1}(\Ctilde) \subset
\op{Jac}^{0}(\Ctilde)$ is disjoint\footnote{In our setting we can see
this fact directly without appealing to Mumford's parity analysis
\cite{Mumford-Prym}. Indeed, by definition the curve $D$ is the double
cover of $\mathbb{P}^{1}$ branched at the $8$ points which are the
union of the $6$ branch points of $\hyp_{C} : C \to \mathbb{P}^{1}$
and the two points $\bry, \brz \in \mathbb{P}^{1}$. Let $\mathsf{Y}$ denote the
ramification point of $\hyp_{D} : D \to \mathbb{P}^{1}$, sitting over
$\bry \in \mathbb{P}^{1}$. Then $\mathsf{a}^{-1}(\mathsf{Y})$ consists of two
distinct points $\tilde{y}, \rho(\tilde{y})$, which are both ramification points
of $\pi : \Ctilde \to C$. But in {\bfseries Case 2.} above we argued
that the line bundle
\[
\mathcal{O}_{\Ctilde}(\tilde{y} - \rho(\tilde{y})) 
\]
cannot be a pullback from $\op{Jac}^{0}(C)$. Thus the point
$\mathfrak{ap}_{(\Ctilde,D)}^{1}(\tilde{t})$ must belong to the
non-neutral component of $\ker \mathsf{Nm}_{\mathsf{a}}$. Since
$\mathfrak{ap}_{(\Ctilde,D)}^{1}(\Ctilde)$ is connected, this implies
that $\mathfrak{ap}_{(\Ctilde,D)}^{1}(\Ctilde)$ is entirely contained
in the 
non-neutral component.
} from the two dimensional abelian subvariety
$\mathsf{Prym}(\Ctilde,D) \subset \op{Jac}^{0}(\Ctilde)$ inside
$\op{Jac}^{0}(\Ctilde)$. But the configuration \eqref{eq:hyp.diagram}
of double covers implies that
\[
\mathsf{Prym}(\Ctilde,D) = \pi^{*} \op{Jac}^{0}(C)  \
\subset \ \op{Jac}^{0}(\Ctilde),
\]
and hence the equation \eqref{eq:case3identity} has no solution.

This completes the analysis of the third case and shows that the map
\eqref{eq:ChatinY1} is injective. To check that the map
\eqref{eq:ChatinY1} also separates tangent directions, consider the
degree six version of the Prym for $(\Ctilde,C)$:
\[
\Prym_{6} = \left\{\left. M \in \op{Jac}^{6}(\Ctilde) \, \right| \,
\mathsf{Nm}_{\pi}(M) = \omega_{C}^{\otimes 3} \right\}.
\]
We have a natural multiplication-by-2 map
\[
\mathsf{mult}_{2} : \Prym_{3} \to \Prym_{6}, \quad M \mapsto M^{\otimes 2}, 
\]
which fits in a commutative diagram
\begin{equation} \label{eq:immerse}
\xymatrix{
\Chat \ar[r] \ar[d]_-{\sqhat} & \Prym_{3}
  \ar[d]^-{\mathsf{mult}_{2}} \\
 \Ctilde \ar[r] & \Prym_{6} 
}
\end{equation}
where the top horizontal map is the map \eqref{eq:ChatinY1}, while the
bottim horizontal map $\Ctilde \to \Prym_{6}$ is given by
$\tilde{t} \mapsto  \pi^{*}(\omega_{C}(\pw))\otimes
\mathcal{O}_{\Ctilde}(\tilde{t} -
\tau(\tilde{t}))$.

Now observe that for any $\tilde{t} \in \Ctilde$ the divisor
$\tilde{t} + \sigma(\tilde{t})$ is in the hyperelliptic linear system
on $\Ctilde$. By the same token $\rho(\tilde{t}) + \sigma\circ\rho(\tilde{t}) =
\rho(\tilde{t}) + \tau(\tilde{t})$ is in the hyperelliptic linear system on
$\Ctilde$. Thus we have a linear equivalence 
\[
\tilde{t} + \sigma(\tilde{t})  \sim \rho(\tilde{t}) + \tau(\tilde{t})
\]
and so
\[
\mathcal{O}_{\Ctilde}(\tilde{t} - \tau(\tilde{t})) \cong
\mathcal{O}_{\Ctilde}(\rho(\tilde{t}) - \sigma(\tilde{t})) =
\mathcal{O}_{\Ctilde}(\rho(\tilde{t}) - \tau(\rho(\tilde{t}))).
\]
Hence the map $\Ctilde \to \Prym_{6}$ factors through $D$. The induced
map $D \to \op{Jac}^{0}(\Ctilde)$ is given by $s \mapsto
a^{*}\mathcal{O}_{D}(s)\otimes
\hyp_{\Ctilde}\mathcal{O}_{\mathbb{P}^{1}}(-1)$ for all $s \in D$. This shows that
the map $\Ctilde \to \Prym_{6}$ factors as
\[
\xymatrix@1@M+1pc{ \Ctilde \ar[r]^-{\mathsf{a}} & D
  \ar[r]^-{\mathfrak{aj}} & \op{Jac}^{1}(D) \ar[rrr]^-{a^{*}(-)\otimes
    \hyp_{\Ctilde}\mathcal{O}_{\mathbb{P}^{1}}(-1)\otimes
    \pi^{*}(\omega_{C}(\pw))} & & & \Prym_{6}.}
\]
The map $\op{Jac}^{1}(D) \to \Prym_{6}$ is an \'{e}tale double cover
which is isomorphic to the quotient of $\op{Jac}^{1}(D)$ by the the translation
action of the $2$-torsion line bundle defining the cover $\Ctilde \to
D$ \cite{Mumford-Prym}. From the diagram \eqref{eq:hyp.diagram} we see
that this $2$-torsion line bundle is given explicitly as
$\mathcal{O}_{D}(\mathsf{Y} - \mathsf{Z})$, where $\mathsf{Y}$ and
$\mathsf{Z}$ are the two Weierstrass points in $D$ that sit over $\bry, \brz
\in \mathbb{P}^{1}$. Since the Abel-Jacobi map $\mathfrak{aj} : D \to
\op{Jac}^{1}(D)$ is a closed embedding, this shows that the map $D \to
\Prym_{6}$ embeds $D - \{\mathsf{Y},\mathsf{Z}\}$ and glues
$\mathsf{Y}$ and $\mathsf{Z}$ into a node. If we write $\underline{D}$
for the iimage of this nodal curve in $\Prym_{6}$ we get that the map
$\Ctilde \to \underline{D} \subset \Prym_{6}$ is injective on tangent
spaces. But the map $\sqhat : \Chat \to \Ctilde$ is \'{e}tale, and so
the composition $\Chat \to \Ctilde \to \underline{D} \subset
\Prym_{6}$ is injective on tangent spaces.  From the diagram
\eqref{eq:immerse} we see that this map also factors as $\Chat \to
\Prym_{3} \to \Prym_{6}$ and so the $\Chat \to \Prym_{3}$ is injective
on tangent spaces. This completes the proof of the lemma.
\end{proof}

\

\subsection{The wobbly locus}
\label{subsec-wobbly}

We recall that Laumon defines the notion of {\em\bfseries very stable}
vector bundle as one that does not admit a non-zero nilpotent Higgs
field. Such bundles are automatically stable
\cite{LaumonVeryStable}. The first two authors therefore introduced
the complementary notion of {\em\bfseries wobbly} as a semistable vector bundle
that isn't very stable. The {\em\bfseries wobbly locus} $\Wob \subset X$ thus
consists of those polystable vector bundles that admit a nonzero
nilpotent Higgs field.

With indexation on the degree our notations for the wobbly loci are:
$$
\Wob_0 \subset X_0, \;\;\;\; \Wob_1 \subset X_1.
$$
These will be the supports for the parabolic structures and
logarithmic poles of the Higgs fields for the Higgs bundles
\eqref{eq:eigensheaves}.

In section~\ref{chapter-d0} we show that the wobbly locus $\Wob_{0}$
is a divisor in the $2$-theta space of $C$ comprising several familiar
players in the classical geometry of the quadric line complex.
Specifically the divisor $\Wob_{0}$ in $X_{0} \cong \mathbb{P}^{3}$
has $17$ irreducible components: the quartic Kummer surface $\Kum
\subset \mathbb{P}^{3}$ and its  $16$ trope planes $\trope_{\kappa}$,
labeled by the theta characteristics of $C$. Thus
\[
\Wob_{0} = \Kum \cup  \left[\cup_{\kappa \in \mathsf{Spin}(C)}
\trope_{\kappa}\right].
\]
From the classical $16_{6}$ configuration in the quadratic line
complex \cite{GH} we know that each trope plane
$\trope_{\kappa}$ is tangent to the Kummer surface $\Kum$ along
a trope conic $\conic_{\kappa}$.  Hence $\Wob_{0}$ fails to be
normal crossing along the trope conics
$\left\{\conic_{\kappa}\right\}_{\kappa\in
  \mathsf{Spin}(C)}$.

It is also easy to characterize the components of the wobbly divisor
from the moduli point of view. The Kummer surface $\Kum$ parametrizes
$S$-equivalence classes of semistable bundles with polystable
representatives of the form $\mathfrak{a}\oplus \mathfrak{a}^{-1}$ for
$\mathfrak{a} \in \op{Jac}^{0}(C)$. The trope plane $\trope_{\kappa}$
can be canonically identified with
$\mathbb{P}(H^{1}(C,\kappa^{-\otimes 2})) =
\mathbb{P}(H^{1}(C,\omega_{C}^{-1}))$ and parametrizes all bundles
that can be realized as non-split extensions $0 \to \kappa^{-1} \to F
\to \kappa \to 0$.

Similarly, in section~\ref{chapter-d1} we show that $\Wob_{1} \subset
X_{1}$ is an irreducible divisor.  Specifically, as mentioned in the
previous section, $\Wob_{1}$ parametrizes all vector bundles $E \in
X_{1}$ which arise as non-split extensions $0 \to A \to E \to
A^{-1}(\pw) \to 0$ for some point $(A,t) \in \Cbar$. In fact we will
see that $\Wob_{1}$ is the tangent developable of the map $\ql :
\Cbar \to X_{1} \subset \mathbb{P}^{5}$, i.e. the union of all
projective tangent lines to points in the curve $\ql(\Cbar)$. This
shows that $\Wob_{1}$ is also a non normal crossings divisor in
$X_{1}$. It has the curve $\ql(\Cbar)$ as its curve of cusps.

\subsection{Hecke correspondences}
\label{subsec-Hecke}

The {\em\bfseries big Hecke correspondence} $\bigHeckebar \rightarrow
X_1\times X_0 \times \Cbar$ is the moduli of quadruples
$(E,F,(A,t),\beta)$ such that $(E,F,(A,t))\in X_1\times X_0 \times
\Cbar$ and $\beta : F\otimes A^{-1} \to E$ is a map that fits in a
short exact sequence
$$
\xymatrix@1@M+0.5pc{
  0\ar[r] & F\otimes A^{-1} \ar[r]^-{\beta}
  &   E \ar[r] &  \cc_{t} \ar[r] & 0.
  }
$$
We will view $\bigHeckebar$ either as a correspondence from $X_{1}$
to $X_{0}\times \Cbar$ or as a correspondence from $X_{0}$ to
$X_{1}\times \Cbar$. We will label  the respective projections as 
\begin{equation} \label{eq:pqdb} 
\quad   \xymatrix@M+0.5pc@-0.5pc{
    & \bigHeckebar \ar[dl]_-{p} \ar[dr]^-{q} & \\
    X_{1} & & X_{0}\times \Cbar
} \qquad
\xymatrix@M+0.5pc@-0.5pc{
    & \bigHeckebar \ar[dl]_-{\pzo} \ar[dr]^-{\qzo} & \\
    X_{0} & & X_{1}\times \Cbar
    }
\end{equation}
where $p = \op{pr}_{X_{1}}$, $\pzo = \op{pr}_{X_{0}}$, $q =
\left(\op{pr}_{X_{0}},\op{pr}_{\Cbar}\right)$, and $\qzo =
\left(\op{pr}_{X_{0}},\op{pr}_{\Cbar}\right)$. In fact the moduli
$\bigHeckebar$ is a subvariety in the triple product $X_{0}\times
X_{1}\times \Cbar$.  To see this, note that when $\beta$ exists for a
given point $(E,F,(A,t)) \in X_{0}\times X_{1}\times \Cbar$, then it
is unique up to scale: $\Oo_{C}(t)$, and hence $t$, is determined as
the ratio of determinants, and then for this given $t$, the Hecke
fiber is a line in one direction, a conic in the other which excludes
any self-intersections.

\

\noindent
In some parts of the discussion we will look at Hecke transformations
supported at a single point.  In those situations we will write
$\bigHeckebar (a)=\op{pr}_{\Cbar}^{-1}(\{ a \} )$ for the preimage of
$a= (A,t) \in \Cbar$ in $\bigHeckebar$ and will view $\bigHeckebar(a)$
as a correspondence between $X_{0}$ and $X_{1}$.

\

The (big) Hecke correspondences appeared in Chapter \ref{synthetic}, Sections  \ref{ssec:synth.corr} 
and \ref{Hecke-curves-Kummer}, from the synthetic viewpoint.
We will review the comparison between the synthetic and modular viewpoints in Section \ref{synth-comparison}. 

\

The {\em\bfseries big abelianized Hecke correspondence} $\bigHeckehat^{\rm
  ab}$ is the blow-up of
$Y_0 \times \Chat$ along a copy of $\Chat \times \Chat$.
As we will see in section~\ref{ssec:abelinized.in.context} the 
map
$$
\xymatrix@M+0.5pc@R-2.5pc{
\emb_{\Chat\times\Chat} : \hspace{-3.5pc} & \Chat\times \Chat \ar@{^{(}->}[r] &
  \Prym_{2}\times \Chat
  \\
& ((A_{1},\tilde{t}_{1}),(A_{2},\tilde{t}_{2})) \ar[r] &
  \left(\pi^{*}\left(A_{2}\otimes A_{1}^{\vee}(\pw)\right)(\tilde{t}_{1}
  - \tilde{t}_{2}), (A_{2},\tilde{t}_{2})\right), }
$$
is a closed embedding. 
This embedding clearly preserves the second projections to $\Chat$ and
the strict transform of $\Chat \times \Chat$ in $Y_{0}\times \Chat$ is
isomorphic copy of $\Chat \times \Chat$ embedded as a subvariety of
$Y_0 \times \Chat$. The blow-up of $Y_{0}\times \Chat$ along this copy
of $\Chat \times \Chat$ is the big abelianized Hecke correspondence
$\bigHeckehat^{\rm ab}$.

The variety $\bigHeckehat^{\rm ab}$ maps to $Y_{0}\times Y_{1} \times
\Chat$ and thus gives correspondences
\[
  \begin{minipage}[c]{2.5in}
 \[
\xymatrix@M+0.25pc{
  & \bigHeckehat^{\rm ab} \ar[dl]_-{p^{\rm ab}}  \ar[dr]^-{q^{\rm ab}}  & \\
 Y_{1} & & Y_{0}\times \Chat
} 
\]
  \end{minipage} \qquad
\begin{minipage}[c]{2.5in}
 \[
\xymatrix@M+0.25pc{
  & \bigHeckehat^{\rm ab} \ar[dl]_-{\pzo^{\rm ab}}  \ar[dr]^-{\qzo^{\rm ab}}  & \\
 Y_{0} & & Y_{1}\times \Chat
} 
\]
\end{minipage} 
\]
There is also  map 
$$
g: \bigHeckehat^{\rm ab} \rightarrow \bigHeckebar, 
$$
such that altogether we get a map of abelianized and usual Hecke
correspondences.
The formulas for the maps $p^{\rm ab}$, $q^{\rm ab}$, $\pzo^{\rm ab}$, $\qzo^{\rm
  ab}$, and $g$  are given in section~\ref{ssec:abelinized.in.context}.

Again, most of the time, we will work with a single fiber of the
abelianized Hecke correspondence. If $\tilde{a} = (A,\tilde{t}) \in \Chat$ is a
fixed point we will focus on the fiber $\bigHeckehat^{\rm
  ab}(\tilde{a}) = \op{pr}_{\Chat}^{-1}(\tilde{a})$. 
  

\

\subsection{Spectral line bundles}
\label{subsec-spectral}

In the abelianization strategy to the Hecke eigensheaf problem we use
spectral data on $C$ to describe the eigenvalue Higgs bundle
$(E,\theta)$ and spectral data on $X_{0}$ and $X_{1}$ to describe the
eigensheaf parabolic Higgs bundles $(\mycal{F}_{0,\bullet},\Phi_{0})$
and $(\mycal{F}_{1,\bullet},\Phi_{1})$. The spectral covers for these
spectral data are $\pi : \Ctilde \to C$ and $f_{0} : Y_{0} \to X_{0}$
and $f_{1} : Y_{1} \to X_{1}$ respectively.

Suppose $(\Ctilde \subset T^{\vee}C,\specN)$ is the spectral data for
$(E,\theta)$, i.e. $(E,\theta) = \left(\pi_{*}\specN, \pi_{*}\left(
\lambda\otimes (-)\right)\right)$, where $\lambda : \Ctilde \to
T^{\vee}C$ is the embedding, and $\specN \in \Prym_{2}$ is the
spectral line bundle.  We will use $\specN$ to construct the
{\em\bfseries spectral line bundles} $\srL _0$ on $Y_0$, and $\srL _1$
on $Y_1$ which will define our parabolic eigensheaf Higgs bundles.

An appropriately normalized Fourier-Mukai transform on the Jacobian of
$\Ctilde$ (see section~\ref{ssec:abelinized.in.context} for the
precise details of the normalization) sends the skyscraper sheaf
$\mathcal{O}_{\specN} \in
D^{b}_{\op{coh}}(\op{Jac}^{2}(\Ctilde))$ to a line bundle
on $\op{Jac}^{2}(\Ctilde)$ with a vanishing first Chern class. We
denote the restriction of this line bundle to $\Prym_{2} \subset
\op{Jac}^{2}(\Ctilde)$ by $\Lprym_{0}$ and define a spectral line bundle
$\LY_{0}$ on the modular spectral cover $f_{0} : Y_{0} \to X_{0}$ by setting
$$
\LY_0 = \left(\blo_{0}^{*}\Lprym_{0}\right)(\ExY_{0})\otimes f_{0}^{*}
\mathcal{O}_{X_{0}}(2) \ \text{in} \ \op{Pic}(Y_{0}).
$$
Using the modular spectral data $(Y_{0},\LY_{0})$ we can now define
a meromorphic Higgs bundle on $X_{0}$ by 
$$
\left(\srF_{0,0},\Phi_{0}\right) :=
\left(f_{0*}\LY_{0},f_{0*}(\modalpha_{0}\otimes (-))\right),
$$
where $\modalpha_{0} : Y_{0} \to T^{\vee}_{X_{0}}(\log \Wob_{0})$ is
the tautological map, defined away from the tacnodes of $\Wob_{0}$.
We will later provide this with a parabolic structure to build a
parabolic rank $8$ bundle $\srF _{0,\bullet}$. Away from the non
normal crossing codimension two strata of $\Wob_0$, the parabolic
structures are trivial; the bundle needs to have a parabolic structure
on a blown-up version of $X_0$ at the tacnodes. The geometry of the
blow-up and the construction of the parabolic structure will be explained
in section~\ref{chapter-d0} and section~\ref{chapter-heckex1x0}.

In this analysis we will also see that the Higgs field of $\Phi_{0} :
\srF_{0,\bullet} \to \srF_{0,\bullet}\otimes \Omega^{1}_{X_{0}}(\log
\Wob_{0})$ has logarithmic poles on both the trope planes and the
Kummer surface.  The residues are nilpotent and have Jordan blocks as
follows: two Jordan blocks of size $2$ and four of size $1$ over the
tropes, and four Jordan blocks of size $2$ over the Kummer.

To describe the spectral line bundle $\LY_{1}$ on $Y_{1}$ it is
convenient to fix a point $\pwtilde\in \Ctilde$ lying over
$\pw$. Recall the assumption that the spectral curve has two branches
over $\pw$. 

The choice of $\pwtilde$ gives an isomorphism
$\trans_{-\pwtilde} : \Prym_{3}
\stackrel{\cong}{\rightarrow} \Prym_{2}$ sending $L$ to $L(-\pwtilde)$
and we can pullback $\Lprym_{0}$ by this isomorphism to define a line
bundle $\Lprym_{1} = \trans_{-\pwtilde}^{*} \Lprym_{0}$ on $\Prym_{3}$
with vanishing first Chern class. With this notation we now set 
$$
\LY_{1}  := \blo_{1}^{*}\Lprym_{1}\otimes \mathcal{O}_{X_{1}}(1),
$$
and 
$$ \left(\srF_{1,0},\Phi_{1}\right) := \left(
f_{1*}\LY_{1},f_{1*}\left(\modalpha_{1}\otimes (-)\right)\right).
$$
In section~\ref{chapter-d1} this is going to be given a parabolic
structure over $\Wob_1$, with parabolic weights $0$ and $1/2$, such
that the associated graded bundle of grading $1/2$ has rank $2$.
This means that
$$
\srF_{1,s} = \srF_1\;\;\; \mbox{for} \;\;\; s \in [0,1/2) .
$$

\subsection{Orbits of the \texorpdfstring{$\cstar$}{C*}-action}
\label{orbits}

An important aspect of the geometry of the Hitchin moduli space is the
$\cstar$-action.  In this section, $\Higgs$ will denote a moduli space
of Higgs bundles of some degree and fixed determinant. In our
situation, it means $\Higgs = \Higgs _0$ or $\Higgs _1$. Some of this
discussion takes place in more general settings so the notation is
left non-specific when possible.

An element $z\in \cstar$ sends $(E,\theta )\in \Higgs$ to
$(E,z\theta)$.  The fixed-point locus $\Higgs ^{\cstar}$
decomposes
$$
\Higgs ^{\cstar} = \Higgs ^{\cstar,\rm u} \sqcup \Higgs ^{\cstar,\rm nu}
$$ into the ``unitary piece'' that is just the moduli of bundles (=
moduli of bundles equipped with zero Higgs fields), i.e.  $\Higgs
^{\cstar,\rm u} \cong X$, and the disjoint union of remaining pieces
that are moduli spaces of polystable Hodge bundles with nonzero Higgs
field.

We recall (see \cite{SimpsonHiggs}) that a Higgs bundle is called a
{\em\bfseries Hodge bundle} if $E=\bigoplus E^p$ and $\theta : E^p
\rightarrow E^{p-1}\otimes \omega _C$. These correspond under the
nonabelian Hodge correspondence to complex variations of Hodge
structure. In our case, we consider bundles of rank $2$, so a
non-unitary Hodge bundle is a direct sum of two line bundles
$$
E\cong L^1 \oplus L^0, \;\;\;\; \theta : L^1 \rightarrow L^0\otimes \omega _C.
$$
We remark that a Hodge bundle with a nonzero Higgs field that is
semistable but not stable could be $S$-equivalent to a polystable
unitary Hodge bundle, that is a to a Hodge bundle with a zero Higgs
field. This will happen in our case if $L^p$ are both line bundles of
degree $0$. Those don't count as points in $\Higgs ^{\cstar,\rm nu}$
since the moduli space parametrizes $S$-equivalence classes and the
polystable representative would be the same bundle with trivial Higgs
field.

\

\begin{proposition}
\label{descrip-fixed}
For the case of rank $2$ Higgs bundles on a curve of genus $2$
considered in this paper, the non-unitary fixed point loci in
$\Higgs_0$ and $\Higgs_1$ are described as follows:
\begin{itemize}

\item[{\rm\bfseries (a)}]
$\Higgs _0^{\cstar,\rm nu}$ is a disjoint union of $16$ points
  parametrizing the uniformizing Higgs bundles
  $\{(E_{\kappa},\theta_{\kappa})\}_{\kappa\in \mathsf{Spin}(C)}$,
  where $E_{\kappa} = \kappa\oplus \kappa^{-1}$, and $\theta_{\kappa}
  : \kappa \to \kappa^{-1}\otimes \omega_{C} = \kappa$ is an
  isomorphism;

\item[{\rm\bfseries (b)}]
$\Higgs _1^{\cstar,\rm nu}$ is connected, isomorphic to the curve $\Cbar$ that
is a $16$-sheeted \'{e}tale covering of $C$.

\end{itemize}
\end{proposition}
\begin{proof}
Suppose $E\cong L^1 \oplus L^0$ is a Hodge bundle. Polystability with
nontrivial Higgs field implies (since we are in the rank $2$ case)
that the Higgs bundle is in fact stable, and $(L^0,0)$ is a sub-Higgs
bundle, so $\deg(L^0) < (\deg E)/2$.  If $E$ has degree $0$ with
determinant $\Oo _C$ then $L^0 = (L^1)^{\vee}$, $\deg L^{0} < 0$, and
the only possibility is that $\deg(L^{0}) = -1$, the Higgs field
$\theta : L^{1} \to L^{0}\otimes \omega_{C}$ is an isomorphism, and
hence $L^1=\kappa$ is one of the $16$ square-roots of the canonical
bundle. If $E$ has degree $1$ with determinant $\Oo _C(\pw )$ then $E
\cong L \oplus L^{-1}(\pw )$, $\deg L \leq 0$, and $\theta :
L^{-1}(\pw)\to L\otimes \omega_{C}$ is a non-zero map.
Therefore we must have $\deg L = 0$ and so the Higgs
field $\theta \in H^0(C,L ^{\otimes 2}(\pw ))$ has a single zero at some 
point $t\in C$, and we get $L^{\otimes -2} \cong \Oo _C(\pw -t)$. This
corresponds to a point of the connected curve $\Cbar$.
\end{proof}

\

The Hitchin fibration $\hit : \Higgs \rightarrow \Bb = \aaaa^N$ is
equivariant for the $\cstar$-action, where the weights of the action
on the base are strictly positive, determined by the degrees of the
invariant polynomials on the Lie algebra or, equivalently, by the
expression in terms of sections of powers of $\omega_C$. In the
present case, the Hitchin base is simply $H^0(C,\omega_C^{\otimes 2})$
so there is a single weight $2$.

Properness of the Hitchin fibration implies that if $y\in \Higgs$ is
any point, then the limit $\lim _{z\rightarrow 0} (zy)$ exists and is
a $\cstar$-fixed point. This fixed point is a unitary Higgs bundle,
i.e. a semistable vector bundle with zero Higgs field, if and only if
$y=(E,\theta )$ with semistable underlying vector bundle $E$, and in
this case $\lim _{z\rightarrow 0} (zy) = (E,0)$.  Let $\Higgs^{\rm su}
\subset \Higgs^{\rm seu}\subset \Higgs$ denote the open subsets of
Higgs bundles whose underlying vector bundle is stable or semistable
respectively. The limiting construction provides a regular map
$\Higgs^{\rm seu} \rightarrow X$ to the moduli space of semistable
bundles.

Over the open subset of stable bundles, the map that we will write as
$\Higgs^{\rm su}\rightarrow X^{\rm s}$ may be identified with the
cotangent bundle
$$
\Higgs^{\rm su} \cong T^{\vee} (X^{\rm s}),
$$
compatibly with $\cstar$ actions, the one on the right being the
scaling action on the total space of the cotangent bundle. This
identification preserves the symplectic structure
\cite{Hitchin-spectral}, in particular the fiber over any point of
$X^{\rm s}$ is a Lagrangian subspace.
 
Let $X^{\rm vs}\subset X$ be the open subset of very stable bundles,
and let $\Higgs^{\rm vs}$ be the open subset of Higgs bundles whose
underlying bundle is very stable.  We have $\Higgs^{\rm vs} \subset
\Higgs^{\rm su}$ and also $\Higgs ^{\rm vs} \cong T^{\vee} (X^{\rm
  vs})$.

The limiting map over the open subset $\Higgs^{\rm seu}$ provides a
rational map
$$
\Higgs \dashrightarrow X.
$$
One may think of this in terms of {\em\bfseries broken
  orbits}. Those are defined as maximal $\cstar$-invariant subsets of
the form
$$
Z_0\cup Z_1 \cup \cdots \cup Z_k \subset \Higgs 
$$
such that $Z_0 \cong \aaaa^1$ and $Z_i\cong \pp^1$ for $i=1,\ldots
,k$, provided with points $I_i,O_i\in Z_i$ (except $I_0$) such that
$O_i = I_{i+1}$, where $O_i$ corresponds to the origin and $I_i$ to
$\infty$ in $Z_i\cong \pp^1$ (or just the origin for $Z_0\cong
\aaaa^1$), such that everything is compatible with the $\cstar$
action. The input point $y$ corresponding to $1\in \aaaa^1=Z_0$ is
mapped by the correspondence of the rational map, to the output point
$O_k\in \Higgs ^{\cstar,\rm u}$.  The intermediate points
$O_i=I_{i+1}$ are non-unitary fixed points. The data needed to
determine such a broken orbit consists of fixing the downward or
outgoing direction (see below) of the next orbit $Z_{i+1}$ whenever
the limiting point $O_i$ is a non-unitary fixed point. The process
stops when we get to a unitary fixed point.

If we start at a general point $y_{\varepsilon}$ nearby to $y$, the
orbit of $y_{\varepsilon}$ limiting to a unitary fixed point, will be
near to a broken orbit starting at $y$.

Let $\Fix\subset \Higgs$ be the complement of $\Higgs ^{\rm seu}$. It is
the set of points $y$ such that $\lim _{z\rightarrow 0} (zy) \in
\Higgs ^{\cstar,\rm nu}$. The limiting map provides a constructible
map from $\Fix$ to the non-unitary fixed point set. In general, this will
not be a regular map. However, in our case there are not very many
non-unitary fixed points so broken orbits have at most one
break. Points of $\Fix$ are those having one break which is a single
well-defined point of $\Higgs ^{\cstar,\rm nu}$.

\begin{lemma}
In our case, the map $y\mapsto \lim _{z\rightarrow 0} (zy)$ is a
regular map from $\Fix$ to the fixed-point set $\Higgs ^{\cstar,\rm
  nu}$.
\end{lemma}
\begin{proof}
In the general situation the map can be non-regular if there is a
broken orbit with two non-unitary fixed points joined by an
orbit. However, in view of the description of
Proposition~\ref{descrip-fixed}, this does not happen in our case.
\end{proof}

In order to understand how broken orbits can work, or equivalently
what their nearby orbits could look like, we need to consider the
local picture of $\Higgs$ at a fixed point. Assume we are at a fixed
point $y$ that is a stable Higgs bundle, so it is a smooth point of
$\Higgs$. The general theory of Bialynicki-Birula \cite{BB} tells us
that the fixed point set is smooth at $y$. The $\cstar$-action on
$\Higgs$ determines an action of $\cstar$ on its tangent space at the
fixed point $y$.

\begin{lemma}
\label{vhs}
The tangent space decomposes into pieces according to weights of the
$\cstar$-action:
$$
T_y(\Higgs ) = \bigoplus T_y(\Higgs )^p
$$
such that $T_y(\Higgs )^p$ is Serre-dual to $T_y(\Higgs )^{1-p}$
and in particular they have the same dimension. The tangent space of
the fixed point set is $T_y(\Higgs )^0$. In our case, at a
non-unitary fixed point, the weights that occur are $-1$ and $2$ in
the case of $\Higgs_0$, and $-1,0,1,2$ in the case of $\Higgs_1$.
\end{lemma}
\begin{proof}
The weights can be understood by thinking of the Hodge bundle
associated to $y$ as being associated to a complex variation of Hodge
structure $V$. The hermitian structure of $V$ (given by a flat but
indefinite hermitian form) leads to a real structure on the VHS
$\text{End}^0(V)$ of trace-free endomorphisms.  Thus \cite{Zucker} the
first cohomology group has a real Hodge structure of weight $1$:
$$
T _y(\Higgs ) \cong H^1(C,\text{End}^0(V)) = \bigoplus _{p+q=1}
T^{p,q}
$$
with $T^{p,q} = \overline{T^{q,p}})$.  This complex conjugation may
be identified, using the real structure, with Serre duality at the
level of Dolbeault cohomology.  The $\cstar$-action acts by weight $p$
on the piece $T^{p,1-p}$.

More concretely if we write $E=\bigoplus E^p$ for the general case, then
$$
T_y(\Higgs ) ^p = T^{p,1-p} = 
{\mathbb H}^{1}\left[ \bigoplus _q
\text{Hom}(E^q, E^{p+q})\stackrel{[-,\theta ]}{\longrightarrow} 
\bigoplus _q
Hom(E^q, E^{p+q-1})
\otimes \omega _C 
\right] 
$$
(where one should also include the trace-free condition but we
didn't put that in so as not to complicate the notation). The complex
is Serre-dual to itself with a shift by $1$.  The tangent space of the
fixed point set is the set of fixed points in the tangent space, as
may be seen by the linearization result to be mentioned shortly
below. In our particular case, the VHS at a fixed point has only two
adjacent Hodge weights, so the VHS of trace-free endomorphisms has
Hodge weights $(-1,1)$, $(0,0)$ and $(1,-1)$.  Thus, the Hodge
structure on $\mathbb{H}^1$ has weights obtained by adding $(1,0)$ and $(0,1)$,
namely $(-1,2)$, $(0,1)$, $(1,0)$ and $(2,-1)$. The $\cstar$-weights
are therefore $-1,0,1,2$. The piece of weight $1$ is complex conjugate
or Serre-dual to the piece of weight $0$, and this is the tangent
space of the fixed point set. Thus, in our case of $\Higgs_0$ this
vanishes whereas for $\Higgs_1$ it has dimension $1$.  We see that
the piece of weight $-1$ has dimension $3$ for $\Higgs _0$ and
dimension $2$ for $\Higgs_1$.
\end{proof}

\

Define the {\em\bfseries incoming} or {\em\bfseries upward} directions
to be the directions on which the weight is $>0$, and the
{\em\bfseries outgoing} or {\em\bfseries downward} directions to be
those on which the weight is $<0$.  The general theory \cite{BB} of
$\cstar$ actions on complex manifolds tells us that there is an
analytic open neighborhood around a fixed point $z$ that is
isomorphic, as a complex manifold with $\cstar$-action, to a
neighborhood of the origin in the tangent space $T_y(\Higgs )$.  In
particular, the fixed points, and incoming and outgoing manifolds are
smooth and identified with the corresponding subspaces of the tangent
space. The incoming manifold is the local subset of points whose limit
as $z\rightarrow 0$ is equal to the fixed point, whereas the outgoing
manifold is the local subset of points whose $z\rightarrow \infty$ is
equal to the fixed point.

The terminologies ``upward'' and ``downward'' come from Hitchin's
picture \cite{Hitchin-selfd} of the energy, or $L^2$-norm of the Higgs
field, as a Morse function on $\Higgs$. We note that the orbit
structure of the $\cstar$-action has been identified, by work of
Collier, Wentworth, Wilkin and others, as being the same as the structure of
flow lines for the gradient of the Morse function
\cite{Wentworth,Wilkin}.  It is because of this picture that we call
the incoming directions ``upward'' (they are the directions where one
goes upwards by following the gradient of the energy function) and the
outgoing directions ``downward'' (they go downwards following the
gradient of the energy function).  The flow stops at the ``bottom''
which is the unitary fixed point set, the moduli space of bundles with
zero Higgs field which is obviously the minimum for the $L^2$ norm.

From Lemma \ref{vhs} we see that the space of incoming directions
(weights $p>0$) is Serre dual to the tangent of the fixed points plus
the space of outgoing directions (weights $p\leq 0$). Thus, the space
of incoming directions has $1/2$ the dimension. In fact, each of the
fibers of the projection over a fixed point, is a Lagrangian subspace
of the symplectic moduli space. The space of outgoing directions, plus
the fixed point directions, together form the tangent space of the
manifold of points that flow out of somewhere in that fixed point
set. This is locally (at a general point of the fixed point set) one
of the components of the nilpotent cone in $\Higgs$, also a Lagrangian
subspace.

In general, the space of directions that are outgoing from a given
fixed point will have dimension smaller than half, although it is also
half the dimension if the fixed point set is $0$-dimensional as is our
case in $\Higgs_0$. Similarly, the total space $\Fix$ of directions
incoming to the local piece of the fixed point set, will generally
have dimension greater than half, but again it is half if the fixed
point set is $0$-dimensional.

In general one will be considering a Prym variety $\Prym \subset
\Higgs$ that is a general fiber of the Hitchin map. The rational map
$\Higgs \dashrightarrow X$ therefore provides a rational map $\Prym
\dashrightarrow X$, and to resolve it one should blow up $\Fix\cap \Prym$
to get $Y\rightarrow \Prym$. The dimension of $\Fix\cap \Prym$ is equal
to the dimension of the fixed point locus.

In our cases, for $\Higgs_0$ the dimension of $\Fix$ is $3$ and it is
Lagrangian, while for $\Higgs _1$ the dimension of $\Fix$ is $4$ and it
fibers over $\Cbar$ with fibers that are Lagrangian.  The locus $\Fix\cap
\Prym$ to be blown up has dimension $0$ for $\Prym \subset \Higgs _0$
and dimension $1$ for $\Prym \subset \Higgs_1$.

One can understand from the local description at the fixed point that
a single blowing-up will be sufficient in our case, yielding a
morphism $Y\rightarrow X$.  We will denote the exceptional divisor by
$\ExY\subset Y$.  In principle, one can use the description of the
$\cstar$ action to understand the resulting ramification of the map
$Y\rightarrow X$ along $\ExY$. We have chosen instead to give more
direct proofs (Lemma \ref{tropes}, Proposition \ref{descrip-ram})
that the map is simply ramified along the exceptional locus.

In the general case when the $\cstar$-action can have a wider range of
weights, and when there can be multiply broken flow lines,
understanding the birational transformation needed to resolve the
rational map, and understanding the resulting ramification, seem to be
difficult questions.

\subsection{The logarithmic property}
\label{logprop}

As before we will write $X$ for the coarse moduli space of semistable
rank two bundles on $C$ with fixed determinant $\mathbf{d}$ and
$\Higgs$ for the moduli space of semistable Higgs bundles with fixed
determinant $\mathbf{d}$. Recall that $\Higgs$ is equipped with an
algebraic $\mathbb{C}^{\times}$-action, scaling the Higgs fields,
i.e. $z \in \mathbb{C}^{\times}$ acts by $(E,\theta) \mapsto
(E,z\theta)$.  Note that with this definition the Hitchin map $\hit :
\Higgs \to \Bb$ becomes $\mathbb{C}^{\times}$-equivariant once we
equip $\Bb = H^{0}(\omega_{C}^{\otimes 2})$ with a scaling action of
degree $2$.

We will also keep using the notation $\Higgs^{\cstar,\text{nu}}
\subset \Higgs$ for the union of components of the fixed point locus
that are disjoint from $X$, and let $\Fix\subset \Higgs$ denote the
incoming variety to $\Higgs^{\cstar,\text{nu}}$, that is
$$
\Fix = \left\{ (E,\theta ) \in \Higgs \ \left| \  \lim _{z\rightarrow 0}
(E,z\theta) \in \Higgs^{\cstar,\text{nu}}\right.\right\}.
$$
In the degree $1$ case, $\Higgs^{\cstar,\text{nu}}=\Cbar $ is an
irreducible curve that is a $16$-sheeted etale covering of $C$,
whereas in the degree $0$ case $\Higgs^{\cstar,\text{nu}}$ consists of
$16$ distinct points.

In our situation, the fixed point locus has only one higher level,
that is to say that the limit for $z\rightarrow 0$ along any outgoing
direction at $\Higgs^{\cstar,\text{nu}}$ lies in $X$. This will no
longer be the case for higher genus or higher rank. It allows for some
simplification.

\

\begin{theorem}
\label{qresolves}
The subvariety $\Fix$ is smooth and contained in the smooth locus of
$\Higgs$. It coincides with the locus of polystable Higgs bundles
whose underlying vector bundle is un(semi)stable.
\begin{itemize}
\item[(a)] Let
$\widetilde{\Higgs}$ denote the blow-up of $\Higgs$ along $\Fix$. Then
this resolves the projection map to $X$, in other words the rationally
defined map $\Higgs \dashrightarrow X$ extends to a morphism $\fullf :
\widetilde{\Higgs} \rightarrow X$.
\item[(b)] $\widetilde{\hit} : \widetilde{\Higgs} \to \Bb$ denote the
  composition of the blow-up morphism $\widetilde{\Higgs} \to \Higgs$
  with the Hitchin map. Let $\mathsf{b} \in \Bb$ be a point
  corresponding to a smooth spectral curve $\Ctilde$ and let $Y =
  \widetilde{\hit}^{-1}(\mathsf{b}) \subset \widetilde{\Higgs}$ be
  the fiber of $\widetilde{\hit}$. Then $Y$ is a smooth compact
  threefold which is a blow-up of of the Hitchin fiber
  $\Prym = \hit^{-1}(\mathsf{b})$ in a smooth center. The induced map $f :=
  \fullf_{|Y} : Y \to X$ is finite and maps each connected component
  of the exceptional divisor of the blow-up $Y \to \Prym$ onto some
  irreducible component of the wobly divisor $\Wob \subset X$.
\end{itemize}
\end{theorem}
\begin{proof}
Note that $\Higgs$ is smooth in the degree $1$ case. In the degree $0$
case, the singular locus consists of the reducible Higgs bundles.  A
point of $\Higgs^{\cstar,\text{nu}}$ is a Hodge bundle, for which the
underlying bundle is of the form $\kappa \oplus \kappa^{-1}$ with
$\kappa^{\otimes 2}= \omega_C$ and for which the Higgs field sends
$\kappa^{-1}$ to zero, and maps $\kappa$ isomorphically to
$\kappa^{-1}\otimes \omega_C$.  Such Higgs bundles are stable. The
stable locus in $\Higgs$ is $\cstar$-invariant and contains a
neighborhood of $\Higgs^{\cstar,\text{nu}}$ so it contains $\Fix$.  In
the previous section we checked that $\Fix$ itself is smooth as it is
a rank three affine bundle over $\Cbar$ in the degree one case and is
the disjoint union of $16$ copies of a three dimensional affine space
in the case of degree zero.  This proves the first statement.

For the second statement we will first deal with the degree zero case.
For this discussion we will write $\Higgs_{0}$, $\Fix_{0}$, $X_{0}$,
etc.  By definition and by the stability comment in the previous
paragraph we have that a Higgs bundle $(E,\theta)$ belongs to $\Fix_{o}$
if and only if $\lim_{z \to 0} (E,z\theta)$ exists as a stable Higgs
bundle and is isomoprphic to
\[
\left( \kappa\oplus
  \kappa^{-1}, \begin{pmatrix} 0 & 0 \\ \gamma & 0 \end{pmatrix}\right),
  \ \text{with}
  \ \xymatrix@1@C-1pc@M+0.3pc{\gamma :  \kappa \ar[r]^-{\cong} &
    \kappa^{-1}\otimes \omega_{C}}.
\]
But if $(E,\theta)$ is a semistable Higgs bundle for which $E$ is not
semistable as bundle, then we can find a saturated line sub bundle
$\kappa \subset E$, such that $\deg \kappa > \deg E/2$, and
$\theta(\kappa) \not\subset \kappa\otimes \omega_{C}$. Since $E$ has
trivial determinant, this means that $E$ fits in a short exact
sequence
\begin{equation} \label{eq:ext.us}
0 \to \kappa \to E \to \kappa^{-1} \to 0,
\end{equation}
$\deg \kappa > 0$, and the composite map
\begin{equation} \label{eq:composite.us}
\xymatrix@1@C-0.5pc@M+0.5pc{ \kappa \ar@{^{(}->}[r] & E \ar[r]^-{\theta} &
  E\otimes \omega_{C} \ar[r] & \kappa^{-1} \otimes \omega_{C} }
\end{equation}
is non-zero. Thus we must have $0 < \deg \kappa \leq g(C) - 1 = 1$, Since
\eqref{eq:composite.us} is a non-zero map between two line bundles of
degree $1$, it must be an isomorphism and  $\kappa$ must be a theta
characteristic on $C$.  In particular  this implies that
$\theta(\kappa)\otimes \omega_{C}^{-1} \subset E$ is a line sub bundle
which projects isomorphically onto $\kappa^{-1}$. Thus the short exact
sequence \eqref{eq:ext.us} is split and we have an isomorphism 
\begin{equation} \label{eq:splitE}
(E,\theta) \  \cong \ \left( \kappa\oplus\kappa^{-1}, \begin{pmatrix}
    \alpha & \beta \\ \gamma  & -\alpha \end{pmatrix} \right),
\end{equation}
where $\alpha \in H^{0}(C,\omega_{C})$, $\beta \in
H^{0}(C,\omega_{C}^{\otimes 2})$, and $\xymatrix@1@C-1pc@M+0.3pc{\gamma :
  \kappa \ar[r]^-{\cong} & \kappa^{-1}\otimes \omega_{C}}$.  The shape
of the Higgs field in \eqref{eq:splitE} can be rigidified somewhat by
choosing the identification $E \cong \kappa\oplus \kappa^{-1}$ more
carefully.  Indeed, the automorphisms of the bundle $\kappa\oplus \kappa^{-1}$
that act trivially on the determinant line bundle are given by
matrices of the form
\[
\begin{pmatrix} a & b \\ 0 & a^{-1} \end{pmatrix}, \quad a\in \cstar,
\ b \in H^{0}(C,\omega_{C}).
\]
Conjugating by the authomorphism $\begin{pmatrix} u & -u^{-1}\alpha
  \\ 0 & u^{-1} \end{pmatrix}$ with $u^{2} = \gamma$, we
get an isomorphism
\[
(E,\theta) \ \cong \ \left( \kappa\oplus \kappa^{-1},
\begin{pmatrix} 0 & -\mathsf{b} \\ 1 & 0 \end{pmatrix} \right)
\]
where $\mathsf{b} = -\alpha^{2} - \gamma\beta = \hit(E,\theta) \in \Bb$ is
the quadratic differential defining the spectral cover of
$(E,\theta)$. This shows that the locus $\Higgs_{0}^{\text{unss}}
\subset \Higgs_{0}$ of all trivial determinant Higgs bundles whose
underlying vector bundle is not semi-stable is the disjoint union
$\Higgs_{0}^{\text{unss}} = \sqcup_{\kappa}  \Higgs_{0,\kappa}^{\text{unss}}$
of the $16$
Hitchin sections
\[
\Higgs_{0,\kappa}^{\text{unss}} =
\left\{\left. \, \left(\kappa\oplus \kappa^{-1},
\begin{pmatrix} 0 & -\mathsf{b} \\ 1 & 0 \end{pmatrix}\right)  \quad
\right| \quad \mathsf{b} \in H^{0}(C,\omega_{C}^{\otimes 2} \, \right\} 
\]
of $\hit : \Higgs_{0} \to \Bb = \mathbb{A}^{3}$
labeled by the theta characteristics $\kappa \in \mathsf{Spin}(C)$.
Also, for any
\[
(E,\theta) \ \cong \ \left(\kappa\oplus \kappa^{-1},
\begin{pmatrix} 0 &  -\mathsf{b} \\ 1 & 0 \end{pmatrix}\right) \in
\Higgs_{0}^{\text{unss}}
\]
we have that whenever  $z \neq 0$ the Higgs bundle $(E,z\theta)$ is
isomorphic to 
\[
\left(\kappa\oplus \kappa^{-1}, 
\begin{pmatrix} 0 & -z^{2}\mathsf{b}
  \\  1 & 0 \end{pmatrix}\right)
\]
where the isomorphism is given by conjugation by $\begin{pmatrix} u &
  0 \\ 0 & u^{-1} \end{pmatrix}$, with $u^{2} = z$.
Thus
\[
\lim_{z \to 0} (E,z\theta) = \left(\kappa\oplus \kappa^{-1}, 
\begin{pmatrix} 0 & 0
  \\ 1 & 0 \end{pmatrix}\right)
\]
and so $\Higgs^{\text{unss}}_{0} \subset \Fix_{0}$. Since both the
un(semi)stable locus and the incoming variety are conical Lagrangian
submanifolds which intersect at the $16$ non-unitary $\cstar$-fixed
points, we get that $\Higgs^{\text{unss}}_{0} = \Fix_{0}$. To shorten
the notation we will write $\Fix_{0,\kappa} =
\Higgs^{\text{unss}}_{0,\kappa}$ for the connected components of this
locus.

Note that each of the Hitchin sections
$\Higgs^{\text{unss}}_{0,\kappa}$ has a trivial normal bundle inside
$\Higgs_{0}$. Indeed, the Hitchin sections intersect each Hitchin
fiber at smooth points, and by the Lagrangian property we have that
the normal bundle of $\Higgs^{\text{unss}}_{0,\kappa}$ inside
$\Higgs_{0}$ is isomorphic to the cotangent bundle of
$\Higgs^{\text{unss}}_{0,\kappa}$. Since
$\Higgs^{\text{unss}}_{0,\kappa} \cong \Bb =
H^{0}(C,\omega_{C}^{\otimes 2})$ we conclude that
$N_{\Higgs^{\text{unss}}_{0,\kappa}/\Higgs_{0}} \cong
\mathcal{O}_{\Higgs^{\text{unss}}_{0,\kappa}} \otimes
H^{0}(C,\omega_{C}^{\otimes 2})^{\vee}$. This implies that when we
form the blow-up $\widetilde{\Higgs}_{0} =
\mathsf{Bl}_{\Fix_{0}}\Higgs_{0}$ the exceptional divisor
$\widetilde{\Fix}_{0}$ is a product
\[
\widetilde{\Fix}_{0} = \Fix_{0} \times \mathbb{P}\left(
H^{0}(C,\omega_{C}^{\otimes 2})^{\vee}\right) \cong  \bigsqcup_{\kappa \in
  \mathsf{Spin}(C)} \Fix_{0,\kappa} \times
\mathbb{P}^{2}.
\]
The question of whether the rational map $\Higgs_{0} \dashrightarrow
X_{0}$ extends to a morphism \linebreak $\widetilde{\Higgs}_{0} \to
X_{0}$ is local near $\Fix_{0} \subset \Higgs_{0}$. Therefore to check
that this happens it suffices to check that if $S$ is a scheme and
$({}^{S}E,{}^{S}\theta)$ is a relative semistable Higgs bundle on
$S\times C$ that gives an \'{e}tale map $S \to \Higgs_{0}$, then the
composite rational map $S \to \Higgs_{0} \dashrightarrow X_{0}$
extends to a morphism $\widetilde{S} = S\times_{\Higgs_{0}}
\widetilde{\Higgs}_{0} \to X_{0}$. This is equivalent to checking that
the vector bundle ${}^{S}E_{|(S - \Fix_{0})\times C}$ has a
canonical (unique up to unique isomorphism) extension to a vector
bundle ${}^{\widetilde{S}}E$ on $\widetilde{S}\times C$ which is
semistable of trivial determinant on all geometric fibers over
$\widetilde{S}$.

As a first approximmation to ${}^{\widetilde{S}}E$ one can take the
pullback
\[
  {}^{\widetilde{S}}E^{\text{naive}} := (\widetilde{S}\times C
  \to S\times C)^{*}\left({}^{S}E\right).
  \]
 This is a rank two vector bundle on $\widetilde{S}\times C$ which has
 trivial determinant on all geometric fibers over $\widetilde{S}$ and
 is semistable on all geometric fibers over points in $\widetilde{S} -
 \widetilde{\Fix}_{0} = S - \Fix_{0}$. However by construction the restriction
 of ${}^{\widetilde{S}}E^{\text{naive}}$ to $\{\tilde{y}\}\times C$
 for any closed point $\tilde{y} \in
 \Fix_{0,\kappa}\times \mathbb{P}^{2} \subset
 \widetilde{\Fix}$ is isomorphic to $\kappa\oplus \kappa^{-1}$ which
 is unstable.

 To construct the actual bundle ${}^{\widetilde{S}}E$ let
 ${}^{S}\Fix_{0} = S \times_{\Higgs_{0}}\Fix_{0}$ and let
 ${}^{\widetilde{S}}\Fix_{0}$ be the Cartier divisor
 $\widetilde{S}\times_{S} ({}^{S}\Fix_{0})$, i.e. the exceptional
 divisor of the blow-up $\widetilde{S} \to S$. We will also write
 ${}^{S}\Fix_{0,\kappa} = \times_{\Higgs_{0}}\Fix_{0,\kappa}$ and
 ${}^{\widetilde{S}}\Fix_{0,\kappa} = \widetilde{S}\times_{S}
 ({}^{S}\Fix_{0,\kappa})$ for the connected components of these loci.

Shrinking
$S$ if necessary we get that the restriction of the bundle
${}^{\widetilde{S}}E^{\text{naive}}$ to
$({}^{\widetilde{S}}\Fix_{0,\kappa})\times C$ has a natural surjective
homomorphism to $\op{pr}_{C}^{*} \kappa^{-1}$ which is unique up to
multiplication by an invertible function on ${}^{S}\Fix_{0,\kappa}$. Define
${}^{\widetilde{S}}E$ to be the Hecke transform of
  ${}^{\widetilde{S}}E^{\text{naive}}$ centered at this homomorphism,
  i.e.
\[
\xymatrix@1@M+0.5pc{
0 \ar[r] & 
{}^{\widetilde{S}}E \ar[r] & {}^{\widetilde{S}}E^{\text{naive}} \ar[r] &
\bigoplus_{\kappa \in \mathsf{Spin}(C)}
\left({}^{\widetilde{S}}\Fix_{0,\kappa}\times C \to
\widetilde{S}\times C\right)_{*}\op{pr}_{C}^{*}\kappa^{-1} \ar[r] & 0.
}
\]
Note that scaling the last map in the sequence by an invertible
function on ${}^{S}\Fix_{0}$ does not change the kernel sheaf
${}^{\widetilde{S}}E$ and so it is canonically defined as a subsheaf
  in ${}^{\widetilde{S}}E^{\text{naive}}$.

The fact that ${}^{\widetilde{S}}E$ is a relative semistable bundle
    on $\widetilde{S}\times C$ now comes  from the following
    
    \

\begin{lemma} \label{lem:rel.ss.deg0}
  Suppose $\tilde{y} = (y,\mathfrak{e}) \in
  {}^{\widetilde{S}}\Fix_{0,\kappa} = ({}^{S}\Fix_{0,\kappa})\times
  \mathbb{P}(H^{0}(\omega_{C}^{\otimes 2})^{\vee})$.  Then the vector
  bundle ${}^{\widetilde{S}}E_{|\{\tilde{y}\}\times C}$ is isomorphic
  to the semi-stable rank two bundle $E_{\mathfrak{e}}$ which is the
  extension
  \[
0 \to \kappa^{-1} \to E_{\mathfrak{e}} \to \kappa \to 0
\]
corresponding to the extension class $\mathfrak{e} \in
\mathbb{P}(H^{0}(\omega_{C}^{\otimes 2})^{\vee}) =
\mathbb{P}(H^{1}(\omega_{C}^{-1})) = \mathbb{P}(H^{1}(\kappa^{\otimes
  -2}))$.
\end{lemma}
\begin{proof}
Since $\mathfrak{e} \in \mathbb{P}(H^{1}(\kappa^{\otimes -2}))$ the
extension $E_{\mathfrak{e}}$ is non-split. But if $L \subset
E_{\mathfrak{e}}$ is a destibilazing line sub bundle, then $\deg L > 0$
  and hence as subsheaves in $E_{\mathfrak{e}}$ we must have $L\cap
  \kappa^{-1} = 0$ this shows that $L$ maps injectively into $\kappa$,
  and hence we will have $\deg L = 1$ and $L \cong \kappa$. This will
  split the extension which is contradiction. This proves that
  $E_{\mathfrak{e}}$ is semi-stable.

The identification of ${}^{\widetilde{S}}E_{|\{\tilde{y}\}\times C}$
and $E_{\mathfrak{e}}$ is a general fact in deformation theory which
we explain next.  Suppose $S$ is the spectrum of a DVR with
uniformizer $\unif$ and closed point $o \in S$. Let $\mathcal{E} \to
S\times C$ be an algebraic vector bundle and let $E =
\mathcal{E}_{|\{o\}\times C}$. Assume that $E$ fits in a short exact sequence
\begin{equation} \label{eq:Ehecke}
\xymatrix@1@M+0.5pc{ 0 \ar[r] & \mathsf{A} \ar[r]^-{\imath} & E
  \ar[r]^-{\jmath} & \mathsf{B} \ar[r] & 0 }
\end{equation}
of vector bundles on $C$.  Consider the Hecke transform
  \[
0 \to \mathcal{E}' \to \mathcal{E} \to \left(\{o\}\times C \to S\times
C\right)_{*} \mathsf{B} \to 0.
\]
Then the restriction $E' = \mathcal{E}'_{|\{o\}\times C}$ is a vector
bundle that fits in a short exact sequence
\begin{equation} \label{eq:Eprime}
  \xymatrix@1@M+0.5pc@R-0.5pc{
    0 \ar[r] & \mathsf{B}\otimes \mathcal{O}_{S\times C}
    (-\{o\}\times C)_{|\{o\}\times C}  \ar[r] \ar[d]_-{\cong}
    & E' \ar[r] & \mathsf{A} \ar[r] & 0 \\
    & \mathsf{B} & &  
  }
\end{equation}
Write 
\[
\mathsf{ks}^{\mathcal{E}} : H^{0}(S,T_{S}) \to H^{1}(C,\op{End}_{0}(E)) \subset
\op{Hom}(E,E[1])
\]
for the Kodaira-Spencer map of $\mathcal{E}$. A direct
calculation with either square zero extenions or \v{C}ech cocycles
gives now that the extension class of the short exact sequence
\eqref{eq:Eprime} is
\[
    \jmath\circ \mathsf{ks}^{\mathcal{E}}(\partial_{\unif}) \circ
    \imath \in \op{Hom}(\mathsf{A},\mathsf{B}[1]) =
    \op{Ext}^{1}(\mathsf{A},\mathsf{B}).
\]
If in addition the sequence \eqref{eq:Ehecke} is split, i.e. $E =
\mathsf{A}\oplus \mathsf{B}$, then we can write the Kodaira-Spencer class
$\mathsf{ks}^{\mathcal{E}}(\partial_{\unif})$ as a block matrix
\[
\mathsf{ks}^{\mathcal{E}}(\partial_{\unif}) = \begin{pmatrix}
  \eta_{\mathsf{A}\mathsf{A}} & \eta_{\mathsf{A}\mathsf{B}}
  \\ \eta_{\mathsf{B}\mathsf{A}} & \eta_{\mathsf{B}\mathsf{B}}
 \end{pmatrix}
\]
where $\eta_{\mathsf{A}\mathsf{A}} : \mathsf{A} \to \mathsf{A}[1]$,
$\eta_{\mathsf{A}\mathsf{B}} : \mathsf{B} \to \mathsf{A}[1]$,
$\eta_{\mathsf{B}\mathsf{A}} : \mathsf{A} \to \mathsf{B}[1]$, and
$\eta_{\mathsf{B}\mathsf{B}} : \mathsf{B} \to \mathsf{B}[1]$.  In
particular, the extension class of \eqref{eq:Eprime} will be
$\eta_{\mathsf{B}\mathsf{A}} : \mathsf{A} \to \mathsf{B}[1]$.

\

Suppose next that $S$ is the spectrum of a DVR which maps to
$\Higgs_{0}$ so that the closed point $o \in S$ maps to a closed point
$y \in \Fix_{0,\kappa}$ and the differential of the map $S \to
\Higgs_{0}$ maps $\partial_{\unif}$ to a tangent vector in
$T_{y}\Higgs_{0}$ which projects to the normal line $\mathfrak{e}
\subset H^{1}(\kappa^{\otimes -2})$. In particular the map $S \to
\Higgs_{0}$ is given by a relative Higgs bundle
$({}^{S}E,{}^{S}\theta)$ on $S \times C$ and the component
$\eta_{\kappa,\kappa^{-1}}$ of
$\mathsf{ks}^{({}^{S}E)}(\partial_{\unif})$ is just equal to
$\mathfrak{e}$. This concludes the proof of the lemma.
\end{proof}

\

These considerations prove assertion (a) in the statement of
Theorem~\ref{qresolves} for moduli of Higgs bundles of determinant
$\mathbf{d} = \mathcal{O}_{C}$.  In fact, the above discussion also
proves assertion (b) in the trivial determinant case. Indeed, when
$\mathbf{d} = \Oo_{C}$ the discussion above shows that the connected
component $\Fix_{0,\kappa}$ of the center of teh blow-up
$\widetilde{\Higgs} \to \Higgs$ intesects the Hitchin fiber
$\hit^{-1}(\mathsf{b}) = \Prym_{2} = \{ L \in \op{Jac}^{2}(\Ctilde)
\ | \ \mathsf{Nm}_{\pi}(L) \cong \omega_{C} \}$ at a single point,
namely the point $\pi^{*}\kappa \in \Prym_{2}$. By the universal
property of the blow-up it follows that $Y$ is the blow-up
of $\Prym_{2}$ at the $16$ distinct points $\pi^{*}\kappa$, $\kappa
\in \mathsf{Spin}(C)$. But the morphism $\Prym_{2} -
\{\pi^{*}\kappa\}_{\kappa \in \mathsf{Spin}(C)} \to X$ is
quasi-finite, and also by the construction in
Lemma~\ref{lem:rel.ss.deg0} the map $f = \fullf_{|Y} : Y \to X$ maps
the exceptional divisor $\ExY_{0,\kappa} \subset Y$ corresponding to
the point $\pi^{*}\kappa \in \Prym_{2}$ isomorphically to the plane in
$X = X_{0}$ parametrizing all semistable bundles which are non-trivial
extensions of $\kappa$ by $\kappa^{-1}$. Thus $f : Y \to X$ is
everywhere quasi-finite and since it is proper, it is
finite. Furthermore note that the plane parametrizing non-trivial
extensions of $\kappa$ by $\kappa^{-1}$ is precisely the trope
component $\trope_{\kappa}$ of $\Wob_{0}$ which completes the proof of
assertion (b) in this case

\

\noindent
Finally we analyze the degree one case. Again we will write
$\Higgs_{1}$, $\Fix_{1}$, $X_{1}$, etc.

\

\begin{lemma} \label{lem:unstable.deg1} Suppose $(E,\theta) \in
\Higgs_{1}$ is a semistable Higgs bundle, such that 
$E$ is not semistable as a bundle. Then 
\[
(E,\theta) \ \cong \ \left( A^{-1}(\pw)\oplus A, \begin{pmatrix}
\alpha & \beta  \\ \gamma & -\alpha  \end{pmatrix} \right)
\]
where $A \in \op{Jac}^{0}(C)$ satisfying $A^{\otimes 2}(\pw) =
\mathcal{O}_{C}(t)$ for some $t \in C$, $\alpha \in
H^{0}(C,\omega_{C})$, $\beta \in H^{0}(C,\omega_{C}(t'))$, $\gamma \neq 0
\in H^{0}(C,\mathcal{O}_{C}(t))$, and $t'$ is the image of $t$ under the
hyperelliptic involution on $C$.
\end{lemma}
\begin{proof}
Let $L \subset E$ be a saturated destabilizing line bundle, i.e. a
saturated line subbundle s.t. $\deg L \geq \deg E/2 = 1/2$. Since
$(E,\theta)$ is semistable as a Higgs bundle, $L$ can not be
$\theta$-invariant. In other words, we must have $\theta(L) \neq 0
\subset E\otimes \omega_{C}$ which is equivalent to the composite map
$\theta : L \to E\otimes \omega_{C} \to (E/L)\otimes \omega_{C}$ being
non-zero. Thus we must have $1 \leq \deg L \leq 1 - \deg L + \deg
\omega_{C} = 3 - \deg L$, i.e. $\deg L = 1$.  Writing $L =
A^{-1}(\pw)$ with $\deg A = 0$ we get a non-zero (hence injective) map
of locally free sheaves $A^{-1}(\pw) \to A(2\pw)$ and so we must have
$A^{\otimes 2}(\pw) = \mathcal{O}_{C}(t)$ for some point $t \in C$,
i.e. we must have that $(A,t) \in \Cbar$.

Also  $E$ fits in a short exact sequence
\begin{equation} \label{eq:deg1unstableseq}
0 \to A^{-1}(\pw) \to E \to A \to 0,
\end{equation}
and viewing $\theta$ as a map $E\otimes \omega_{C}^{-1} \to $ we get a
subsheaf $\theta\left(A^{-1}(\pw)\otimes \omega_{C}^{-1}\right)
\subset E$ which maps injectively into $A$ via the map $E \to A$. This
implies that the extension \eqref{eq:deg1unstableseq} splits over the
subsheaf $\theta : A^{-1}(\pw)\otimes \omega_{C}^{-1} \hookrightarrow
A$. Since $(A,t) \in \Cbar$ we have that $A^{-1}(\pw)\otimes
\omega_{C}^{-1} \cong A(-t)$ and under this isomorphism the map
$\theta$ gets identified with the natural inclusion $A(-t) \subset A$.

By the same token we have
\[
\op{Ext}^{1}(A,A^{-1}(\pw)) =
H^{1}(C,A^{\otimes -2})(\pw)) = H^{1}(C,\mathcal{O}_{C}(2\pw - t)) =
H^{1}(C,\mathcal{O}_{C}(t')),
\]
while
\[
\op{Ext}^{1}(A(-t),A^{-1}(\pw)) = H^{1}(C,A^{\otimes -2}(t + \pw)) =
H^{1}(C,\mathcal{O}_{C}(2\pw)) = H^{1}(C,\mathcal{O}_{C}(t'+t)),
\]
and under these identifications the map
\[
\op{Ext}^{1}(A,A^{-1}(\pw)) \to \op{Ext}^{1}(A(-t),A^{-1}(\pw))
\]
gets identified with the map on cohomology
\begin{equation} \label{eq:maponextindeg1}
H^{1}(C,\mathcal{O}_{C}(t')) \to H^{1}(C,\mathcal{O}_{C}(t'+t))
\end{equation}
induced from the natural inclusion $\mathcal{O}_{C}(t') \subset
\mathcal{O}_{C}(t'+t)$.

From the long exact sequence in cohomology associated to the short
exact sequence of sheaves
\[
0 \to \mathcal{O}_{C}(t') \to \mathcal{O}_{C}(t'+t) \to \mathcal{O}_{t} \to 0
\]
we see that the map \eqref{eq:maponextindeg1} is surjective. on the
other hand both $H^{1}(C,\mathcal{O}_{C}(t'))$ and \linebreak
$H^{1}(C,\mathcal{O}_{C}(t'+t))$ are one dimensional and so the map
\eqref{eq:maponextindeg1} is an isomorphism. This implies that if an
extension class in $\op{Ext}^{1}(A,A^{-1}(\pw))$ maps to zero in
$\op{Ext}^{1}(A(-t),A^{-1}(\pw))$, then this class is zero to begin
with. Hence the extension \eqref{eq:deg1unstableseq} is split as
claimed. The statement about the matrix entries of $\theta$ under the
decomposition $E \cong A^{-1}(\pw)\oplus A$ now follows
tautologically. Note that the entry $\gamma$ is precisely the composite map
\[
\xymatrix@1@M+0.5pc{
  A^{-1}(\pw) \ar@{^{(}->}[r] & E \ar[r]^-{\theta} &
  E\otimes \omega_{C} \ar[r] & A\otimes \omega_{C} = A(2\pw)
  }
\]
and hence is non-zero. This finishes the proof of the lemmma.
\end{proof}

\

\noindent
Again, the shape of the Higgs field in the pair $(E,\theta)$ can be
simplified by choosing the isomorphism $E \cong A^{-1}(\pw)\oplus A$
more carefully. The global automorphisms of the bundle $A^{-1}(\pw)\oplus A$
are all of the form
\[
\begin{pmatrix}
a & b \\ 0 & a^{-1}
\end{pmatrix}
\]
where $a \in \cstar$, and $b \in H^{0}(C,\mathcal{O}_{C}(t'))$.
Conjugating
$z\cdot \begin{pmatrix} \alpha & \beta \\ \gamma & -\alpha \end{pmatrix}$
by $\begin{pmatrix}
u & 0 \\ 0 & u^{-1}
\end{pmatrix}$ with $u^{2} = z$ we get
\[
\begin{pmatrix}
  z\alpha & z^{2}\beta \\
  \gamma & -z\alpha
\end{pmatrix}.
\]
Hence
\[
\lim_{z \to 0} (E,z\theta) = \left( \, A^{-1}(\pw)\oplus
A, \begin{pmatrix} 0 & 0 \\ \gamma & 0 \end{pmatrix}\, \right),
\]
and so the locus $\Higgs_{1}^{\text{unss}}$ parametrizing stable Higgs
bundles with unstable underlying vector bundles is contained in the
incoming variety $\Fix_{1}$. By the explicit description both of these
are four dimensional irreducible subvarieties in $\Higgs_{1}$ and
therefore $\Higgs_{1}^{\text{unss}} = \Fix_{1}$.

Similarly to the degree zero case the identifictaion
$\Higgs_{1}^{\text{unss}} = \Fix_{1}$ is compatible both with the
Hitchin map $\hit : \Higgs_{1} \to \Bb$ and with the flow limit map
$\Fix_{1} \to \Higgs_{1}^{\cstar,\text{nu}}$. Indeed,if
\[
\left( \, A^{-1}(\pw)\oplus
A, \begin{pmatrix} \alpha & \beta \\ \gamma & -\alpha \end{pmatrix}\, \right) \
\in \
\Higgs_{1}^{\text{unss}},
\]
then the determinant of the Higgs field is $\mathsf{b} = -\alpha^{2} -
\gamma\beta$. This shows that $\Higgs_{1}^{\text{unss}} = \Fix_{1}$
surjects onto the Hitchin base $\Bb$.  In Section~\ref{ssec:moduli} we
also checked that the fiber of $\Higgs_{1}^{\text{unss}}$ over a
general point $\mathsf{b} \in \Bb$ of the Hitchin base is a copy of
$\Chat_{\mathsf{b}} = \Cbar\times_{C} \Ctilde_{\mathsf{b}}$ embedded
in the corresponding Hitchin fiber $\Prym_{3,\mathsf{b}} =
\hit^{-1}(\mathsf{b})$. Thus we can view the smooth variety $\Fix_{1}$
as a family of curves over $\Bb$.

At the same time, the above caluclation of $\lim_{z \to 0} zy$ for $y
\in \Fix_{1}$ identifies the limiting flow map $\flow_{1} : \Fix_{1} \to
\Cbar$ with the natural projection of ``forgetting the Higgs
field'', i.e. for every
\[
y = \left(A^{-1}(\pw)\oplus A, \,
\begin{pmatrix} \alpha & \beta \\ \gamma & -\alpha
 \end{pmatrix} \right) \in \Fix_{1}
\]
we have $\flow_{1}(y) = A \in \Cbar$.

With this observation at hand we can now compute the normal bundle of
smooth subvariety $\Fix_{1} \subset \Higgs_{1}$ and the exceptional
divisor of the blow up $\widetilde{\Higgs}_{1} = \op{Bl}_{\Fix_{1}}
\Higgs_{1} \to \Higgs_{1}$. Specifically we have

\

\begin{lemma} \label{lem:normal.deg1}
  \begin{itemize}
  \item[(a)] The normal bundle of $\Fix_{1} \subset \Higgs_{1}$ is given by
    \[
N_{\Fix_{1}/\Higgs_{1}} = H^{0}(C,\omega_{C})^{\vee}\otimes
\flow_{1}^{*}\sq^{*}\omega_{C}.
\]
Thus the exceptional divisor $\widetilde{\Fix}_{1}$ of the blow up
$\widetilde{\Higgs}_{1} \to \Higgs_{1}$ is
\[
\widetilde{\Fix}_{1} = \mathbb{P}(N_{\Fix_{1}/\Higgs_{1}}) =
\Fix_{1}\times \mathbb{P}(H^{0}(C,\omega_{C})^{\vee}) =
\Fix_{1}\times \mathbb{P}^{1}.
\]
\item[(b)] Suppose $\pi: \Ctilde \to C$ is a smooth spectral
  curve. Let $\Prym_{3}$ denote the corresponding Hitchin fiber, and
  let $\Chat \subset \Prym_{3}$ be the correponding fiber of the map
  $\hit_{|\Fix_{1}} : \Fix_{1} \to \Bb$.  Then
  \[
  N_{\Chat/\Prym_{3}} = H^{0}(C,\omega_{C})^{\vee}\otimes
  \sqhat^{*}\pi^{*}\omega_{C},
  \]
 and the exceptional divisor $\ExY_{1} \subset Y_{1}$ of the blow up
 $\blo_{1} : Y_{1} \to \Prym_{3}$ is given by
 \[
\ExY_{1} = \Chat \times \mathbb{P}(H^{0}(C,\omega_{C})^{\vee}) =
\Chat\times \mathbb{P}^{1}.
 \]
  \end{itemize}
\end{lemma}
\begin{proof}
  To prove (a) note that since $\Fix_{1}$ is $\cstar$-stable, we have
  that the normal bundle $N_{\Fix_{1}/\Higgs_{1}}$ is
  $\cstar$-equivariant. Furthermore, as we explained in the previous
  section, in our situation there are no broken orbits in $\Fix_1$ and
  therefore the normal bundle will be the pullback of the bundle of
  outgoing directions to the non-unitary fix point set
  $\Higgs_{1}^{\cstar,\text{nu}} = \Cbar$, i.e.
  \[
  N_{\Fix_{1}/\Higgs_{1}} = \flow_{1}^{*} \left(T_{\Higgs_{1}}\right)^{(-1)}_{|\Cbar}
  \]
is the piece in the tangent bundle to $\Higgs_{1}$ of $\cstar$-weight
$-1$ along the fixed point curve $\Cbar$. Note that here $\Cbar$ is
embedded in $\Higgs_{1}$ by the map sending $A \in \Cbar$ to the
isomorphism class of the Higgs bundle
  \[
    (E_{A},\theta_{A}) :=
    \left( A^{-1}(\pw)\oplus A,\, \begin{pmatrix} 0 & 0
    \\ \gamma & 0 \end{pmatrix}\right)
  \]
with $\gamma : A^{-1}(\pw) \to A\otimes \omega_{C}$ being the unique
(up to scale) non-zero map.

From the cohomological description of the tangent space of
$\Higgs_{1}$ at a point we see that the piece of $\cstar$-weight $-1$
of the tangent space at a point $A \in \Cbar \subset \Higgs_{1}$ is
precisely the image of the composition map
\[
T_{\Higgs_{1},(E_{A},\theta_{A})} \to
\mathbb{H}^{1}\left(C,\op{End}_{0} (E_{A})
\stackrel{\left[\theta_{A},-\right]}{\longrightarrow}
\op{End}_{0} (E_{A})\otimes
\omega_{C}\right) \to H^{1}(C,\op{End}_{0} (E_{A})).
\]
This image is precisely the matrix
coefficient piece $H^{1}(C,A^{\otimes 2}(-\pw)) \subset
H^{1}(C,\op{End}_{0} (E_{A}))$ which gives a canonical identification 
\[
T^{(-1)}_{\Higgs_{1},(E_{A},\theta_{A})} = H^{1}(C,A^{\otimes 2}(-\pw)).
\]
Next recall, that $A$ is a point in $\Cbar$ precisely when $A^{\otimes
  2} \cong \mathcal{O}_{C}(t - \pw)$ for some point $t \in C$. Thus
$A^{\otimes 2}(-\pw)$ is isomorphic to $\mathcal{O}_{C}(-t')$ where
$t'$ is the hyperelliptic conjugate of $t$. But the inclusion
$\mathcal{O}_{C}(-t') \subset \mathcal{O}_{C}$ induces an isomorphism
on $H^{1}$'s and so we conclude that we get a natural, unique up to
scale, isomorphism
\begin{equation} \label{eq:T-1iso}
T^{(-1)}_{\Higgs_{1},(E_{A},\theta_{A})} \cong H^{1}(C,\mathcal{O}_{C}).
\end{equation}
More precisely, the isomorphism \eqref{eq:T-1iso} is the map induced
on $H^{1}$'s from the map of locally free rank one sheaves $A^{\otimes
  2}(-\pw) \cong \Oo_{C}(-t') \hookrightarrow \Oo_{C}$. Since both the
isomorphism $A^{\otimes 2}(-\pw) \cong \Oo_{C}(-t')$ and the inclusion
$\Oo_{C}(-t') \hookrightarrow \Oo_{C}$ are uniquely defined up to
scale it follows that isomorphism \eqref{eq:T-1iso} is also unique up
to scale.

This implies that we have a line bundle $L$ on $\Cbar$ so that
\[
(T^{(-1)}_{\Higgs_{1}})_{|\Cbar} \cong H^{1}(C,\mathcal{O}_{C})\otimes
L = H^{0}(C,\omega_{C})^{\vee}\otimes L.
\]
To finish the proof of (a) we need to compute the line bundle $L$. But
note that the restriction of the normal bundle of $\Fix_{1} \subset
\Higgs_{1}$ to a general fiber $\Chat$ of $\hit_{|\Fix_{1}}$ is the
normal bundle of $\Chat \subset \Prym_{3}$. Since the map $\pihat :
\Chat \to \Cbar$ induces an injective map on Picard varieties, we only
need to compute $\pihat^{*}L$. Thus (a) will follow immediately once we
prove (b).

\

\noindent
To prove (b) fix a smooth spectral cover $\pi : \Ctilde \to C$ whose
branch locus does not include any Weirstrass point on $C$. Let
$\Prym_{3} = \op{Prym}^{3}(\Ctilde,C) \subset \Higgs_{1}$ denote the
corresponding Hitchin fiber and let $\Chat = \Fix_{1}\cap \Prym_{3}$
be the unstable locus in $\Prym^3$. Recall that $\Chat =
\Cbar\times_{C}\Ctilde$ and the embedding $\Chat \hookrightarrow
\Prym_{3}$ is given by $(A,\tilde{t})
\mapsto\pi^{*}(A^{-1}(\pw))(\tilde{t})$. The tangent bundle to
$\Prym_{3}$ is trivial with fiber identified with the anti-invariants
$H^{1}(\Ctilde,\mathcal{O}_{\Ctilde})^{-} \subset
H^{1}(\Ctilde,\mathcal{O}_{\Ctilde})$ for the covering involution of
$\pi : \Ctilde \to C$.  By Serre duality we can identify
$H^{1}(\Ctilde,\mathcal{O}_{\Ctilde})^{-}$ with the anti-invariants in
$H^{0}(\Ctilde,\omega_{\Ctilde})^{\vee}$. As $\Ctilde \to C$ is a
spectral cover we have
\[
\begin{aligned}
  H^{0}(\Ctilde,\omega_{\Ctilde}) &
  = H^{0}(\Ctilde,\pi^{*}\omega_{C}^{\otimes 2}) \\
  & = H^{0}(C,\pi_{*}\pi^{*}\omega_{C}^{\otimes 2}) \\
  & = H^{0}(C,\omega_{C}^{\otimes 2}\oplus \omega_{C}) \\
  & = H^{0}(C,\omega_{C}^{\otimes 2})\oplus H^{0}(C,\omega_{C}),
  \end{aligned}
\]
and so the anti-invariants in $H^{0}(\Ctilde,\omega_{\Ctilde})^{\vee}$ are
identified intirinsically with
\[
H^{0}(C,\omega^{\otimes 2}_{C}).
\]
With this identification the normal sequence for the embedding $\Chat
\subset \Prym_{3}$ can be written as
\[
0 \to T_{\Chat} \to H^{0}(C,\omega^{\otimes 2}_{C})^{\vee}\otimes
\mathcal{O}_{\Chat} \to N_{\Chat/\Prym_{3}} \to 0,
\]
where the first map is the Kodaira-Spencer map for the varying family
of line bundles $\pi^{*}(A^{-1}(\pw))(\tilde{t})$ on
$\Ctilde$. Dually, the conormal sequence reads
\begin{equation} \label{eq:conormal.in.prym}
0 \to N_{\Chat/\Prym_{3}}^{\vee} \to H^{0}(C,\omega^{\otimes 2}_{C})\otimes
\mathcal{O}_{\Chat} \to \omega_{\Ctilde} \to 0,
\end{equation}
where the last map is the composition
\[
\xymatrix@1@M+0.5pc{ H^{0}(C,\omega^{\otimes 2}_{C})\otimes
  \mathcal{O}_{\Chat} \ar@{^{(}->}[r]^-{\pi^{*}} &
  H^{0}(\Ctilde,\omega_{\Ctilde}))\otimes \mathcal{O}_{\Ctilde}
  \ar[r]^-{\op{ev}} & \omega_{\Ctilde}.  }
\]
But $\omega_{\Ctilde} = \pi^{*}\omega_{C}^{\otimes 2}$ is a pull-back
of $\mathcal{O}(2)$ from the hyperelliptic $\mathbb{P}^{1}$ and
$H^{0}(C,\omega_{C}^{\otimes 2})$ is the pullback of
$H^{0}(\mathbb{P}^{1},\mathcal{O}(2))$. Thus the conormal sequence
\eqref{eq:conormal.in.prym} is a pullback of the canonical evaluation
sequence of vector bundles on the hyperelliptic $\mathbb{P}^{1}$
\[
\xymatrix@1@M+0.4pc{ 0 \ar[r] & \mathcal{K}  \ar[r] &
  H^{0}(\mathbb{P}^{1},\mathcal{O}(2))\otimes
  \mathcal{O}_{\mathbb{P}^{1}} \ar[r] & \mathcal{O}(2)
  \ar[r] & 0.}
\]
But this is a standard Koszul sequence and so we have a canonical
identification
\[
\mathcal{K} = H^{0}(\mathbb{P}^{1},\mathcal{O}(1))\otimes \mathcal{O}(-1).
\]
Equivalently we can note that by definition the bundle $\mathcal{K}$
has no cohomology and so is isomorphic to a two dimensional vector
space tensored with $\mathcal{O}(-1)$. Since the whole sequence is
$SL(2,\mathbb{C})$ equivariant, this two dimensional vector space has
to be the fundamental representation as claimed.

Altogether we get that
\[
N^{\vee}_{\Chat/\Prym_{3}} =
H^{0}(\mathbb{P}^{1},\mathcal{O}(1))\otimes
\sqhat^{*}\pi^{*}\hyp_{C}^{*}\mathcal{O}(-1) =
H^{0}(C,\omega_{C})\otimes \sqhat^{*}\pi^{*}\omega_{C}^{-1}
\]
which finishes the proof of part (b) and the lemma.
\end{proof}

\

Again, the question of whether the rational map $\Higgs_1
\dashrightarrow X_{1}$ is resolved by blowing up $\Fix_1 \subset
\Higgs_{1}$ is local near $\Fix_1$. Similarly to the degree zero case
it suffices to check that if $S$ is a scheme and
$({}^{S}E,{}^{S}\theta)$ is a relative semistable Higgs bundle on
$S\times C$, which corresponds to an \'{e}tale map $S \to \Higgs_{1}$,
then the composite rational map $S \to \Higgs_{1} \dashrightarrow
X_{1}$ extends to a morphism $\widetilde{S} = S\times_{\Higgs_{1}}
\widetilde{\Higgs}_{1} \longrightarrow X_{1}$. Thus we must show that
the vector bundle ${}^{S}E_{|(S - \Fix_{1})\times C}$ has a canonical
  extension to a vector bundle ${}^{\widetilde{S}}E$ on
  $\widetilde{S}\times C$ which is semistable and of determinant
  $\mathcal{O}_{C}(\pw)$ on all geometric fibers over $\widetilde{S}$.

Let ${}^{S}\Fix_{1} = S\times_{\Higgs_{1}} \Fix_{1}$ and let
${}^{\widetilde S}\Fix_{1} =
\widetilde{S}\times_{\widetilde{\Higgs}_{1}} \Fix_{1}$ denote the
exceptional divisor for the blow-up $\widetilde{S} \to S$.  Shrinking
$S$ if necessary we can assume without loss of generality that the
restriction of ${}^{S}E$ to ${}^{S}\Fix_{1}\times C$ is a direct sum
  \[
{}^{S}E =
\left({}^{S}\!\!\mathcal{A}\right)^{-1}(\op{pr}_{C}^{*}\pw)\oplus
     \left({}^{S}\!\!\mathcal{A}\right)
\]
where ${}^{S}\!\!\mathcal{A}$  denotes the pullback
\[
  {}^{S}\!\!\mathcal{A} = \left(
  \xymatrix@1@M+0.3pc@C+0.5pc{
    {}^{S}\Fix_{1}\times C \ar[r] &
\Fix_{1}\times C \ar[r]^-{\flow_{1}\times \op{id}} & 
\Cbar\times C}\right)^{*} \mathcal{A}
\]
of the normalized Poincar\'{e} line bundle $\mathcal{A} \to
\Cbar\times C$, which is uniquely characterized by the conditions
$\mathcal{A}_{\{A\}\times C} \cong A$ for all $A \in \Cbar$, and
$\mathcal{A}_{\Cbar\times \{\pw\}} \cong \mathcal{O}_{\Cbar}$.

Writing ${}^{\widetilde{S}}\!\!\mathcal{A}$ for the pullback line bundle
\[
{}^{\widetilde{S}}\!\!\mathcal{A} =
\left({}^{\widetilde{S}}\Fix_{1}\times C \to {}^{S}\Fix_{1}\times
C\right)^{*} {}^{S}\!\!\mathcal{A}
\]
we can repeat the reasoning  we used in the degree zero case we define
  ${}^{\widetilde{S}}E$ to be the Hecke transform
\[
\xymatrix@1@M+0.3pc{ 0 \ar[r] & {}^{\widetilde{S}}E \ar[r] &
  \left(\widetilde{S}\times C \to S\times C\right)^{*} \left({}^{S}E\right)
  \ar[r]
  & \left({}^{\widetilde{S}}\Fix_{1}\times C \to \widetilde{S}\times
  C\right)_{*} \left({}^{\widetilde{S}}\!\!\mathcal{A}\right) \ar[r] & 0.
}
\]
By the calculation in the previous lemma we see that the fiber of the
normal bundle to $\Fix_{1}$ at a point $(E,\theta)$ with $E =
A^{-1}(\pw)\oplus A$ are precisely the parametrizing the extension
space
\[
\op{Ext}^{1}(A^{-1}(\pw),A) = H^{1}(C,A^{\otimes 2}(-\pw)) =
H^{1}(C,\mathcal{O}_{C}(-t')) = H^{1}(C,\mathcal{O}_{C}).
\]
Suppose now that $S$ is the spectrum of a DVR which maps to
$\Higgs_{1}$ so that the closed point $o \in S$ maps to $y =
(E,\theta) \in \Fix_{1}$ with $E = A^{-1}(\pw)\oplus A$, and such that
the differential of $S \to \Higgs_{1}$ maps the tangent vector
$\partial_{\unif}$ to a normal vector $\mathfrak{e} \neq 0 \in
H^{1}(C,\mathcal{O}_{C}(-t')) = H^{1}(C,\mathcal{O}_{C})$. Again the
map $S \to \Higgs_{1}$ is given by a relative Higgs bundle
$({}^{S}E,{}^{S}\theta)$ on $S\times S$ for which the lower left
corner of the matrix representing the Kodairra-Spencer class
$\mathsf{ks}^{\left({}^{S}E\right)}(\partial_{\unif})$ is equal to
$\mathfrak{e}$. The deformation theory calculation in
Lemma~\ref{lem:rel.ss.deg0} now implies that the bundle
${}^{\widetilde{S}}E_{|\{[\mathfrak{e}]\}\times C}$ is the extension
of $A^{-1}(\pw)$ by $A$ given by the class $\mathfrak{e}$. This showas
that for every $\mathfrak{e}$ the bundle
${}^{\widetilde{S}}E_{|\{[\mathfrak{e}]\}\times C}$ is stable and
completes the proof of the Theorem~\ref{qresolves}(a) in the case of
Higgs bundles with determinant $\mathbf{d} = \Oo_{C}(\pw)$. We have
also proven Theorem~\ref{qresolves}(b) in this case. Indeed, in
Lemma~\ref{lem:normal.deg1}(b) we saw that $Y = Y_{1}$ is the blow-up
of $\Prym_{3}$ centered at a copy of $\Chat \subset \Prym_{3}$. Thus
$Y$ is a smooth projective threefold.  Furthermore,
Lemma~\ref{lem:normal.deg1}(b) and the discussion that follows
immediately the proof of Lemma~\ref{lem:normal.deg1} we see that the
map $f = \fullf_{|X} : Y \to X$ sends the exceptional divisor
$\ExY_{1}$ onto the wobbly divisor $\Wob_{1} \subset X_{1}$ and the map
$f_{1} : \ExY_{1} \to \Wob_{1}$ factors as a composition $\ExY_{1} =
\Chat\times \pp^{1} \to \Cbar\times \pp^{1} \to \Wob_{1}$, where
$\Chat\times \pp^{1} \to \Cbar\times \pp^{1}$ is the natural double
cover, and $\Cbar\times\pp^{1} \to \Wob_{1}$ is the normalization
map. In other words, in this case $\ExY_{1}$ and $\Wob_{1}$ are
irreducible and $f_{|\ExY_{1}} : \ExY_{1} \to \Wob_{1}$ is a double
cover. Similarly to the degree zero case this implies that $f : Y \to
X$ is quasi-finite everywhere and hence is finite since it is a proper
map. The theorem is proven.
\end{proof}

\

\

\subsection{The main construction} \label{subsec:main.construction}

\

\noindent
The description of the resolution in Theorem~\ref{qresolves}  shows in
particular that the exceptional divisor of the blowup maps onto the
wobbly locus in the moduli of bundles. This allows us to understand
the relationship between the moduli of Higgs bundles and the cotangent
bundle of the moduli of vector bundles not only over the very stable
locus but also in codimension one in the stable locus, i.e. over the
generic point of the wobbly divisor.

Let $\widetilde{\Fix}\subset \widetilde{\Higgs}$ denote the
exceptional divisor, that is  - the inverse image of $\Fix$. 
Recall that we have a natural morphism
$$
\psi : T^{\vee}X \rightarrow
\Higgs
$$
inducing an isomorphism on dense open subsets. 

Let $\Wob^{\rm sing}$ denote the singular locus of the
wobbly divisor $\Wob$, and set
$$
X^{\circ}:= X - \Wob^{\rm sing}, \;\;\;\;
\Wob^{\circ}:= \Wob\cap X^{\circ} = \Wob-\Wob^{\rm sing}.
$$
We note that $\Wob^{\circ}\subset X^{\circ}$ is now a smooth divisor. 

Let $\widetilde{\Higgs^{\circ}}:= \fullf^{-1}(X^{\circ})\subset
\widetilde{\Higgs}$ and let $\widetilde{\Fix}^{\circ}$ be the
intersection of $\widetilde{\Fix}$ with $\widetilde{\Higgs^{\circ}}$.
With this notation we now have

\

\begin{proposition} \label{prop:tot.space}
Let $T^{\vee}X^{\circ}(\log \Wob^{\circ})$ denote the total 
space of the locally free sheaf $\Omega ^1_{X^{\circ}}(\log \Wob^{\circ})$.
There is a commutative diagram
$$
\xymatrix{
  T^{\vee}X^{\circ}(\log \Wob^{\circ}) \ar[rr]^{\psi^{\log}} \ar[dr] & &
  \widetilde{\Higgs^{\circ}} \ar[dl]^{\fullf}  \\
& X^{\circ} &}
$$ such that the map $\psi^{\rm log}$ coincides with $\psi$ over a
dense open subset, and such that $\psi^{\rm log}$ maps the
zero-section of $T^{\vee}X^{\circ}(\log \Wob^{\circ})$ isomorphically
to $\widetilde{\Fix}^{\circ}$.
\end{proposition}
\begin{proof}
This follows tautologically from the Hecke description of the families
of Higgs bundles describing the map $\fullf$ over the exceptional
$\widetilde{\Fix}$ that we gave in the proof of
Theorem~\ref{qresolves}.
\end{proof}

\

\

\noindent
{\bf The construction:} Again we write $\Higgs$ for the moduli space
od Higgs bundles of determinant $\mathbf{d} = \Oo_{C}(k\cdot \pw)$
where $k$ isfixed to be either $0$ or $1$. Let $\hit : \Higgs \to \Bb$
denote the corresponding Hitchin map and let $\mathsf{b} \in \Bb$ be a
general point in the Hitchin base. Concretely, choose $\mathsf{b}$ so
that the correspondng spectral cover $\pi : \Ctilde =
\Ctilde_{\mathsf{b}} \to C$ is smooth and unramified above the
Weierstrass points of $C$.  Write $\Prym = \hit^{-1}(\mathsf{b})$ for
the corresponding fiber of the Hitchin map. Then $\Prym$ is smooth,
and isomorphic to the Prym variety of the cover $\Ctilde \to C$. Let
$Y$ denote the pullback of $\Prym$ in the blown up moduli space
$\widetilde{\Higgs}$. Then $Y$ is a smooth projective threefold, which
by Theorem~\ref{qresolves} is obtained by blowing up $\Prym$ in the
smooth subvariety $\Prym\cap \Fix$ which in the case of degree zero
consists of $16$ points - the $16$ line bundles
$\{\pi^{*}\kappa\}_{\kappa \in \mathsf{Spin}(C)}$ and in the case of
degree one is a copy of the curve $\Chat$. Furthermore from the proof
of Theorem~\ref{qresolves}(b) we see that the restriction $f :=
\fullf_{|Y} : Y \to X$ is a finite morphism which maps the exceptional
divisor in $Y$ onto a union of components of the wobbly divisor in
$X$. In the degree zero case the exceptional divisor obtained by
blowing up the point $\pi^{*}\kappa \in \Prym$ maps birationally onto the
trope plane $\trope_{\kappa} \subset \Wob_{0} \subset X_{0}$ labeled
by the theta characteristic $\kappa$. In the degree one case the
exceptional divisor obtained by blowing up $\Chat \subset \Prym$ maps
two-to-one onto the wobbly divisor $\Wob_{1} \subset X_{1}$. Recall
also that as explained in the introduction, the finite morphism $f
: Y \to X$ has degree $2^{3g(C) - 3} = 8$.

\

\begin{definition}
Let $\Omega^{1}_{X}(\log \Wob)^{+}$ denote the unique reflexive sheaf
on $X$ whose restriction to $X^{\circ}$ is the 
vector bundle $\Omega^{1}_{X^{\circ}}(\log \Wob^{\circ})$.
\end{definition}

\

If we let $j^{\circ}:X^{\circ}\hookrightarrow X$, then we can set more
precisely
$$
\Omega^{1}_{X}(\log \Wob)^{+}:= \left( 
j^{\circ}_{*}\Omega^{1}_{X^{\circ}}(\log \Wob^{\circ}) \right) ^{\vee\vee}.
$$
Note that $\Omega^{1}_{X}(\log \Wob)^{+}$ is locally free over $X^{\circ}$
and there it is equal to $\Omega^{1}_{X^{\circ}}(\log \Wob^{\circ})$.

\

Let $\Ll$ be a line bundle on $Y$.  Then set $\Vv := f_{*}(\Ll )$. It
is a rank $8$ vector bundle on $X$. It has a meromorphic Higgs field
coming from the fact that a dense open subset of $Y$ may be identified
with a subvariety of $T^{\vee}X$ via the inverse to $\psi$.

\begin{corollary}
\label{loghiggs}
Under the main construction 
$\Vv = f_{*}(\Ll )$ comes equipped with a meromorphic Higgs field 
which comes from a morphism
$$
\Phi : \Vv \rightarrow  \Vv \otimes 
\Omega^{1}_{X}(\log \Wob)^{+} .
$$
Viewed as a sheaf with operators, $\Vv$ has no nontrivial
$\Phi$-invariant subsheaves. 
\end{corollary}
\begin{proof}
Since $\Omega^{1}_{X}(\log \Wob)^{+}$ is defined as a reflexive
hull, it suffices to see this over $X^{\circ}$, and there it
comes from the commutative diagram of Proposition~\ref{prop:tot.space}. 

The spectral variety of $\Phi$ is birational to 
the covering $f : Y\rightarrow X$, and this is 
irreducible, so there can be no $\Phi$-invariant subsheaves
of $\Vv$, indeed such a subsheaf would correspond to a nontrivial
decomposition of $Y$ as a union of closed subsets. 
\end{proof}

\

\subsection{Nonabelian Hodge outside codimension \texorpdfstring{$2$}{2}}
\label{subsec-nah}

In Corollary \ref{loghiggs}, we are not able to say that we construct
a ``logarithmic Higgs field'' because the divisor $\Wob$ does not have
normal crossings. The objective of this section is to investigate how
the data of a logarithmic Higgs bundle defined outside of codimension
$2$ determines a local system by the nonabelian Hodge
correspondence. We will also include parabolic structures.

In this section we consider a general variety $X$, a divisor $W$, and
an open subset $X^{\circ}\subset X$ complement of a closed subset of
codimension $\geq 2$, such that $W^{\circ}:= W\cap X^{\circ}$ is
smooth.

Suppose we are given a reflexive sheaf $\Vv$ over $X$ that is locally
free over $X^{\circ}$, and a parabolic structure for $\Vv$ along
$W^{\circ}$, denoted $\Vv _{\bullet}$.  Suppose $\Phi: \Vv \rightarrow
\Vv \otimes \Omega^{1}_{X}(\log W)^{+}$ is an operator with values in
the sheaf of differentials that is the reflexive extension of
$\Omega^{1}_{X^{\circ}}(\log W^{\circ})$ as defined previously. We
suppose that the residues of $\Phi$ along components of $W^{\circ}$,
acting on the parabolic associated graded pieces, are nilpotent.

Consider a projective morphism $g: Z \rightarrow X$ and let
$Z^{\circ}:= g^{-1}(X^{\circ})$. Suppose $D\subset Z$ is a simple
normal crossings divisor such that $D^{\circ} := D\cap Z^{\circ}$
contains $g^{-1}(W^{\circ})$. The following statement is a recasting
of Mochizuki's extension theorem \cite{Mochizuki-D2} to our setting.

\begin{theorem}
\label{outsidecod2}
Suppose that $g:(Z,D)\rightarrow (X,W)$ exists as above,
and suppose that there exists a semistable logarithmic
parabolic Higgs sheaf $\Ff$ on $(Z,D)$ whose restriction
to $Z^{\circ}$ is isomorphic to $g^{*}((\Vv _{\bullet},\Phi )|_{X^{\circ}})$.
Suppose furthermore that $c_1^{\rm par}(\Ff ) = 0$ and 
$c_2^{\rm par}(\Ff )=0$.  Then there is a mixed twistor $\srD$-module
${\bf E}$ over $X$, with singularities along $W$, such that
the restriction of ${\bf E}$ to $X^{\circ}$ corresponds to  
$(\Vv _{\bullet},\Phi )|_{X^{\circ}}$. Furthermore, 
${\bf E}$ is unique up to isomorphism, and the pullback 
$g^{*}({\bf E})$ corresponds to $\Ff$. 

If $S\subset X^{\circ}$ is any projective curve then the
restriction of ${\bf E}$ to $S$ is the mixed twistor
$\srD$-module that corresponds to 
the logarithmic parabolic Higgs bundle
on a curve $(\Vv _{\bullet}, \Phi )|_S$. 
\end{theorem}

\

\noindent
The mixed twistor $\srD$-module ${\bf E}$ has Betti realization
that restricts
to a local system on $X-W$, and this local system has 
quasi-unitary monodromy transformations around the
components of $W^{\circ}$ such that the arguments of
eigenvalues of the monodromy are the parabolic weights
and the unipotent part of the monodromy for each eigenvalue
has the same Jordan form as the residue of $\Phi$ on the
corresponding graded piece. This comes from the restriction
to curves $S\subset X^{\circ}$. 

Stated more compactly, the conclusion of the above theorem is that
in order to construct such a local system on $X-W$
it suffices to know the parabolic structure and logarithmic
Higgs field over an open subset $X^{\circ}\subset X$ whose
complement has codimension two, provided we can find some
covering $Z/X$ where the divisor has normal crossings and
where there exists an extension with vanishing parabolic
Chern classes. 

This is what we will be doing in our concrete situations.

\subsection{Parabolic structures}
\label{subsec-parabolic}

In this subsection we discuss some aspects of parabolic structures. We
refer to the numerous available references for the general theory of
parabolic bundles and their Chern classes, including for example
\cite{BiswasOrbifold, BiswasChern, Boden,BodenHu, Borne, BorneVistoli,
  Konno, MaruyamaYokogawa, IyerSimpson, IyerSimpson-dr, Taher1,
  Taher2}.  The present discussion will be tailored to our specific
needs. Suppose $(X,D)$ is a pair consisting of a smooth variety and a
reduced divisor. Let $X^{\circ}\subset X$ be the complement of a subset of
codimension $\geq 2$, and let $D^{\circ}:= D\cap X^{\circ}$ (with this
notation extended in a similar way to other objects). We assume that
$D^{\circ}$ is smooth.  That is notably attained by taking $X^{\circ}$
to be the complement of the non-smooth locus of $D$, but one might
want to throw out other subsets of codimension $\geq 2$ as well.

Define the notion of \emph{\bfseries quasi-parabolic bundle} on
$(X,D)$ to consist of a reflexive sheaf $\Vv$ on $X$ together with a
filtration by strict subsheaves $0=F_0 \subset F_1 \subset \cdots
\subset F_k = \Vv |_{D^{\circ}}$ of the restriction of $\Vv $ to
$D^{\circ}$. A \emph{\bfseries parabolic bundle} is a quasi-parabolic
bundle together with an assignment of real parabolic weights to the
subquotients of the filtration over connected components of
$D^{\circ}$.

If $D$ has simple normal crossings, we obtain by extension filtrations
over each smooth irreducible component of $D$. Assuming $D$ has simple
normal crossings, a \emph{\bfseries locally abelian parabolic bundle}
on $(X,D)$ is a parabolic bundle such that $\Vv$ is locally free, the
filtrations on divisor components have locally free subquotients, and
locally at any point there exists a frame adapted to the filtrations
of all divisor components passing through that point. If, in addition,
the parabolic weights are rational, then these objects correspond to
vector bundles on a root stack of $(X,D)$
\cite{BiswasOrbifold,Boden,Borne,BorneVistoli}.

For the purposes of Theorem \ref{outsidecod2}, we are interested in
parabolic Chern numbers of the form $c_1^{\rm par} \cdot [\omega
]^{n-1}$ and $c_2^{\rm par} \cdot [\omega ]^{n-2}$, and for those it
suffices to know the parabolic structure outside of a subset of
codimension $3$ in $X$. If $D$ has normal crossings outside of
codimension $3$, any parabolic structure satisfies the locally abelian
condition outside of codimension $3$.

Our particular application is to the case where $D=\Wob$ is the wobbly
divisor in $X=X_1$, and furthermore the parabolic weights are $-1/2$
and $0$. We will be using a technique of pullback to a ramified covering
to compensate for the fact that $\Wob$ has non-normal crossings
singularities, namely cuspidal singularities as well as nodes up to
$X$-codimension $2$. The fact that there are only two distinct
parabolic weights means that the filtration involves a single subsheaf.

In this setting, the notion of parabolic structure can be simplified
further.  Given $(X,D)$ a \emph{\bfseries  crude quasi-parabolic structure}
consists of a reflexive sheaf $\Vv$ and a subsheaf $\Vv' \subset \Vv$
such that $\Vv$ is locally free on $X^{\circ}$, $\Vv'$ is also
reflexive, and $\Uu :=\Vv / \Vv'$ is a locally free sheaf over
$W^{\circ}$. We let $r_\Uu:= \text{rank}(\Uu )$ be the rank of this
quotient.

A \emph{\bfseries crude parabolic structure} consists of a crude
quasi-parabolic structure plus a \emph{\bfseries parabolic weight} which is a
single number $\alpha \in (0,1]$. To these data we associate a
parabolic vector bundle on $(X^{\circ},W^{\circ})$ whose weights are
$\alpha$ for the associated-graded piece $\Uu$, and $0$ for other the
associated-graded piece $\Uu':= \Vv'/\Vv(-W)$.

The filtered sheaf $\Ee ^{\circ}_{\cdot}$  on $X^{\circ}$ is given by 
$$
\Ee ^{\circ}_a = \Vv ', \;\; 0\leq a < \alpha, 
\;\;\;\;
\Ee ^{\circ}_a = \Vv , \;\; \alpha \leq a < 1.
$$
Let $\Omega ^1_X(\log D )^+$ be the reflexive extension to $X$ of
$\Omega ^1_{X^{\circ}}(\log D^{\circ} )$ which is well-defined since
we are assuming that $D^{\circ}$ is smooth. This isn't really the
correct ``sheaf of logarithmic differentials'' as we will confront when
dealing with the tacnodes of the wobbly locus in the degree $0$ moduli
space, but let us leave that discussion for later.

Say that a map $\Phi : \Vv \rightarrow \Vv \otimes \Omega ^1_X(\log D
)^+$ is a \emph{\bfseries pre-logarithmic Higgs field} if it satisfies
$\Phi \wedge \Phi = 0$ generically. It is called \emph{\bfseries
nil-compatible} with the parabolic structure if $\Phi (\Vv ') \subset
\Vv ' \otimes \Omega ^1_X(\log W)^+$ and if the residue of $\Phi$
along $D^{\circ}$ is nilpotent.  This nilpotence is equivalent to
demanding that the two associated-graded pieces of the residue,
operating on $\Uu$ and $\Uu '$, are nilpotent.

Given such a logarithmic Higgs field, then on the
parabolic structure associated to $(\Vv , \Vv ', \alpha )$
we obtain a parabolic logarithmic Higgs field with 
nilpotent residues.

The following construction shows how to get a crude
logarithmic parabolic structure from a spectral
covering. 

\begin{theorem}
\label{crudeconstruction}
Suppose $f : Y\rightarrow X$ is a finite morphism,
and suppose we are given a factorization through
an inclusion $\psi : Y^{\circ}\rightarrow T^{\vee}X^{\circ}(\log D^{\circ})$
over an open subset $X^{\circ}$ complement
of a subset of codimension $\geq 2$. Set $Y^{\circ}:= f ^ {-1}(X^{\circ})$
and so forth. 
Let $B\subset Y$
be a divisor such that $B^{\circ}=B\cap Y^{\circ}$ is smooth and $f(B)\subset D$. 
Assume that
$\psi^{\circ}$ maps the reduced inverse image
$f ^{-1}(D^{\circ})^{\rm red}$ to the zero-section
of $T^{\vee}X^{\circ}(\log D^{\circ})$. 
Let $\Ll$ be a line bundle on $Y$, and let $\Ll ':= \Ll (-B)$.
Then $\Vv := f _{*}(\Ll )$ has a subsheaf $\Vv ':= 
f _{*}(\Ll )$. For any choice of parabolic weight $\alpha$,
and letting $\Phi$ be the pre-logarithmic Higgs field determined by
$\psi$, we obtain a crude parabolic pre-logarithmic Higgs bundle
$(\Vv , \Vv ', \alpha, \Phi )$ on $(X,D)$ with nilpotent residues. 
\end{theorem}
\begin{proof}
Indeed, the quotient $\Uu := \Vv / \Vv '$ is isomorphic to
$f_{*}(\Ll |_B)$ and over $D^{\circ}$ this is a locally free
quotient of $\Vv |_D$.  The fact that $\psi$ maps $Y^{\circ}$ to
$T^{\vee}X^{\circ}(\log D^{\circ})$ implies that $\Phi$ is logarithmic
over $X^{\circ}$, for both $\Vv$ and its subsheaf $\Vv'$. On the
reflexive extensions this yields a crude parabolic pre-logarithmic
Higgs bundle. The condition on the zero-section insures that the
residues of $\Phi$ are nilpotent.
\end{proof}

\noindent
{\bfseries Terminology:} In the situation of the theorem, we will call
$\Ll '$ the \emph{\bfseries spectral line bundle} on $Y$, because in
terms of the filtered sheaf we have
$$
\Ee _0 = f_{*}(\Ll ').
$$

There are two main difficulties in the theoretical situation, both
occasioned by the fact that the wobbly divisor $D=\Wob$ generally
does not have normal-crossings singularities:

\begin{itemize}
\item  How to calculate effectively
the parabolic $c_2$ or the parabolic Bogomolov $\Delta$-invariant? 

\item How to recognize if the Higgs field $\Phi$ will correspond
to a genuinely logarithmic Higgs field on a normal crossings
resolution?

\end{itemize}

Over the degree $0$ moduli space $X_0$, the crude parabolic structure
is going to be trivial. The divisor $\Wob$ then represents the
location of singularities of the Higgs field. In that case, the
strategy will be to follow the prescription in Theorem
\ref{outsidecod2} to make a resolution of singularities of the divisor
$\Wob$ into a normal crossings divisor, and look for a parabolic
extension whose parabolic Chern class vanishes. In this case, the
parabolic structures will have more than two jumps so one needs to
consider filtrations of bigger length.

For the degree $1$ moduli space $X_1$, we will be using a crude
parabolic structure with parabolic weight $\alpha = 1/2$. This was
found partly by computation and partly by guessing.  Our strategy for
the degree $1$ case consists of passing to a ramified covering in
order to remove the singularities of $D=\Wob$ or transform them to
normal crossings.  Although the pullback of a parabolic bundle with
arbitrary weight under a ramified covering may be difficult to
compute, in the case $\alpha = 1/2$ we can give a description.

\begin{proposition}
\label{pullbackcrude}
Suppose $\underline{\Vv}=(\Vv , \Vv ', \alpha , \Phi )$ is a crude
parabolic pre-logarithmic Higgs bundle on $(X,D)$ with parabolic
weight $\alpha = 1/2$.  Suppose $g:Z\rightarrow X$ is a Galois
covering from a smooth variety, such that $g$ is ramified of order $2$
along $B^1= g^{-1}(D^{\circ})^{\rm red}$ over $D^{\circ}$.  Then the
pullback $g^{*}(\underline{\Vv})$ has trivial parabolic structure,
so it corresponds to a reflexive sheaf $\Vv _Z$ with logarithmic Higgs
field that is described in the following way. Let $\Uu := \Vv / \Vv
'$, let $\Uu _Z:= g^{*}(\Uu )$, and use a superscript $()^{\circ}$
to denote the restriction over $X^{\circ}$. Let $\Uu _{B^{\circ}}$ be
the restriction of $\Uu ^{\circ}$ to $B^{\circ}$ seen as a coherent
sheaf on $Z^{\circ}$, also equal to $\Uu _Z^{\circ} /
\mathcal{I}_{B^{\circ} \subset Z^{\circ}}\Uu _Z^{\circ}$.  We have an
exact sequence
$$
0\rightarrow \Vv ^{\circ}_Z \rightarrow g^{*}\Vv ^{\circ} \rightarrow
\Uu _{B^{\circ}}\rightarrow 0
$$
and $\Vv _Z$ is the reflexive extension of $\Vv ^{\circ}_Z$
from $Z^{\circ}$ to $Z$. The logarithmic Higgs field 
$\Phi$ is induced from that of $g^{*}(\Vv )$. 
\end{proposition}
\begin{proof}
Follows immediately from the definition of a pullback of parabolic
structures.
\end{proof}

\subsection{Pushforward statement}
\label{pushforward-statements}

In order to calculate the Hecke transforms of the local systems that
we are going to construct, we need a Dolbeault method for calculating
higher direct images. This is provided by the theory of \cite{dirim},
which is based in turn on the theory of twistor $\srD$-modules of
Sabbah and Mochizuki \cite{Mochizuki-D1,Mochizuki-D2,Sabbah}.

We will need to extend the discussion of \cite{dirim} in order to
apply it to the specific situations encountered for the Hecke
transform. Luckily, the fibers of the Hecke correspondence are
$1$-dimensional, since we are dealing with moduli spaces of rank $2$
bundles in this paper. In order to identify the higher direct image
local systems, it suffices to restrict to curves in the target space,
and to make things easy we can even use lines. Thus, the direct image
calculations are for maps from a surface to a curve, the main setup of
\cite{dirim}. The present situation differs in that the horizontal
divisor typically has simple ramification points, whereas the
hypothesis of \cite{dirim} was to have an etale horizontal
divisor. Thus, some work needs to be done in the direction of looking
at a ramified cover of the base.

The full discussion will be deferred to Chapter
\ref{chapter-pushforward}.  In this subsection we summarize
the basic knowledge to come out of that, as needed to treat the Hecke
transforms. The setup we discuss here is therefore closely tailored to
the situations that are encountered in the applications.

Suppose $f:X\rightarrow S$ is a projective morphism from a smooth
surface to a smooth curve. Suppose $D\subset X$ is a reduced divisor,
each of whose components dominates $S$.  Classify the points of
$(X,D)$ among the following kinds:

\

\

\punkt  {\bfseries Types of points of $(X,D)$}
  \label{def:typesofpts}

\begin{enumerate}[label=\upshape(\alph*),ref=\thesubsubsection (\alph*)]
\item\label{type1}
Points in $X-D$ where $f$ is smooth;

\item\label{type2}
Points in $X-D$ where the fiber of $f$ has a simple normal crossing;

\item\label{type3}
Points on $D$ where $D$ is etale over $S$; 

\item\label{type4}
Points on $D$ where $f$ is smooth and $D$ is smooth with simple
ramification over $S$;

\item\label{type5}
Points on $D$ where $f$ is smooth and $D$ has a normal crossing such
that both branches are etale over $S$;

\item\label{type6}
Other points. 

\end{enumerate}

\

Suppose we are given a normal variety $\Sigma$ that we will call
(under a small abuse of notation) the \emph{\bfseries spectral
variety}, with a map $\pi_{\Sigma} : \Sigma \rightarrow X$, together with a
line bundle $\Ll$ that we will call the \emph{\bfseries spectral line
bundle} and a $1$-form $\alpha$ on the smooth locus of $\Sigma$ that
we will call the \emph{\bfseries tautological form}.
 
Set $E_0:= \pi_{\Sigma *}(\Ll)$. We assume that $\alpha$ leads to a map
$\Sigma \rightarrow T^{\vee}(X,\log D)$, implying that the tautological form
provides $E_0$ with a Higgs field $\theta : E_0 \rightarrow E_0\otimes
\Omega ^1_X(\log D)$.

Notice here that the image of $\Sigma$ in the logarithmic cotangent
bundle might not be normal, but we look at the normalization $\Sigma$
and call that the spectral variety, whence the abuse of notation
mentioned above.

We will consider two cases: the \emph{\bfseries nilpotent case}, where
the parabolic structure is trivial and the Higgs field has a nilpotent
residue along $D$ coming from the ramification of $\Sigma$; and the
\emph{\bfseries parabolic case} where the parabolic weights are $0,
1/2$. In this case, we suppose $E_{1/2} = \pi_{\Sigma *}(\Ll(R))$ where
$R\subset \Sigma$ is the part of the ramification divisor of $\Sigma
\to X$ that sits over $D$.  Note that $\Sigma$ will in general have
other ramification that is not included in $R$, the hypothesis is that
all the ramification over $D$ contributes to the parabolic structure.

The Higgs field is given by the spectral $1$-form $\alpha \in
H^0(\Sigma, \Omega ^1_{\Sigma})$ We assume that $\alpha$ is nonzero on
normal vectors to $R$ at general points. This implies, in the
nilpotent case, that the nilpotent residues of the Higgs field are
nonzero, giving a size $2$ Jordan block at each point of $R$ over a
general point of $D$.

Define the \emph{\bfseries upper critical locus} 
$$
\widetilde{{\rm Crit}}(X/S, E_{\bullet}, \theta )  \subset \Sigma
$$
to be the closure of the zero-locus of the projection of $\alpha$
to a section of relative differentials on $X/S$, from the open subset
where $\Sigma$ is smooth.  The \emph{\bfseries lower critical locus}
${\rm Crit}(X/S,E_{\bullet}, \theta )$ is its image in $X$.

Let $G$ denote the normalization of $\widetilde{{\rm Crit}}(X/S,
E_{\bullet}, \theta )$, and let $g: G \rightarrow S$ be the induced
morphism. Let $j:G \rightarrow X$ denote the normalization of the inclusion map. 

\

\begin{hypothesis}
\label{pushforward-hyp}
We make the following genericity hypotheses, in addition to the basic assumptions above. 

\begin{itemize}

\item
The lower critical locus ${\rm Crit}(X/S,E_{\bullet}, \theta )$
does not contain any points of type (\ref{type6});

\item
The spectral variety $\Sigma$ is smooth, except that each point above
a point of type (\ref{type5}) is an ordinary double point;

\item 
If there are points of type (\ref{type5}) then we are in the nilpotent
case;

\item
The upper critical locus is smooth except over points of type
(\ref{type5}), in particular its normalization $G$ is a smooth curve.

\item 
The restriction of the spectral $1$-form $\alpha$ to the vertical
direction in $\Sigma$ over a point of type (\ref{type4}) is nonzero.

\item
At each ordinary double point over a point of type (\ref{type5}), the
curve $G$ has two smooth branches whose tangent vectors are distinct.

\end{itemize}
\end{hypothesis}

\

Let $Q:= (G\cap R)$ be the divisor in $G$ that is the intersection
with $R$. The hypothesis that $\alpha$ is nonzero in the vertical
direction at points of type (\ref{type4}) implies that $G$ is transverse
to $R$ at those points, so $Q$ is reduced at such points (see Lemma
\ref{alpha-nonzero-vertical}).

In practice, $G$ will only meet $R$ over points of type (\ref{type4})
or (\ref{type5}): in the nilpotent case that comes from constancy of
the Jordan form of the residue of the Higgs field along $D$, while in
the parabolic case it is a property of our setup, see Lemma
\ref{onlymeets}. We can however state the pushforward property without
including such an hypothesis.

\

The basic philosophy of the pushforward statement,
whose proof will be given in Chapter \ref{chapter-pushforward}, 
is that the cohomology along the
fibers is calculated by localizing at the zeros of the relative Higgs field, namely 
the upper critical locus. Its normalization $G$ is a smooth curve destined to be the 
spectral curve of the higher direct image Higgs bundle. 

The contribution at a local point will come from the spectral 
line bundle, tensored with the relative differentials of $X/S$. A correction is needed at
the points where $D$ is ramified over the base, so the expressions for the higher direct image spectral
line bundle will include correction
terms at the divisor $Q$ in addition to $\Ll |_G \otimes j^{*}\omega _{X/S}$. 

The points of type
(\ref{type5}) introduce a difficulty since the spectral variety is not smooth and the
upper critical locus has normal crossings. Our statement identifying the higher direct image Higgs
bundle will require, in this case, a global degree calculation---that does work out in our cases of interest.

\begin{theorem}
\label{pushforward-thm}
Assume the basic hypotheses explained at the start, assume that we are
either in the parabolic or nilpotent cases as described above, and
assume the hypotheses \ref{pushforward-hyp}.  Then, the parabolic
Higgs bundle $(F_{\bullet}, \Phi )$ on $S$ corresponding to the middle direct
image of the perverse sheaf associated to $(E_{\bullet},\theta )$ has the
following descriptions depending on the case.

\noindent
---In the parabolic case, $F_{\bullet}$ has trivial parabolic structure and 
$$
F_0 = g_{*} \left( 
\Ll |_G \otimes j^{*}\omega _{X/S} \otimes \Oo _G(Q)
\right) .
$$

\noindent
---In the nilpotent case, $F_{\bullet}$ has a parabolic structure with
weights $0,1/2$ and contains a parabolic subsheaf
$F'_{\bullet}\hookrightarrow F_{\bullet}$ defined in the following way:
$$
F'_0 = g_{*} \left( 
\Ll |_G \otimes j^{*}\omega _{X/S} 
\right) 
$$
and 
$$
F'_{-1/2} = g_{*} \left( 
\Ll |_G \otimes j^{*}\omega _{X/S} \otimes \Oo _G(-Q)
\right) .
$$ 
We have $F'_{\bullet}= F_{\bullet}$ 
away from the images of points of type (\ref{type5}), in particular if there are no such points
then this gives the expression for $F_{\bullet}$. 

If the parabolic sheaf $F'_{\bullet}$ defined as above
has global parabolic degree $0$ then
$F'_{\bullet} = F_{\bullet}$ and these give expressions for
$F_{\bullet}$.  

In both the nilpotent and parabolic cases, the Higgs field on $F_0$ comes from the
differential form on $G$ obtained by restricting the tautological
differential $\alpha$ from $\Sigma$.
\end{theorem}

\subsection{Chern classes}
\label{sec-Chern}

For reference below, we recall here the required facts about Chern classes. 
We'll almost always be interested in what happens up to codimension
two, so by convention---unless otherwise specified---our formulas will be 
truncated at
codimension two. Thus for example we write
$$
{\rm ch}(E) = r(E) + c_1(E) + (c_1(E)^2/2 -c_2(E).
$$
In the other direction,
$$
c_2=c_1^2/2 - {\rm ch}_2.
$$
An important invariant is Bogomolov's discriminant
\cite{bogomolov-tensors,bogomolov}
$$
\Delta = \frac{1}{2r}c_1^2 - ch_2 = c_2- \frac{r-1}{2r}c_1^2 ,
$$
having the property that $\Delta (E\otimes L)=\Delta (E)$
for a line bundle $E$. 

The same hold for parabolic Chern classes. 

Suppose $\pi : Y\rightarrow X$ is a map. Then the
Grothendieck-Riemann-Roch formula says that for a coherent sheaf $L$
on $Y$ we have
$$
{\rm ch}(\pi _{*}(L))
= \pi _{*}({\rm td}(Y/X) \cdot {\rm ch}(L)),
$$
where the relative Todd class is
$$
{\rm td}(Y/X) = {\rm td}(TY) \pi ^{*}{\rm td}(TX)^{-1}
$$
with
$$
{\rm td}(TY)= 1 + \frac{c_1(TY)}{2}  
+ \frac{c_1(TY)^2 + c_2(TY)}{12} 
$$
and similarly for ${\rm td}(TX)$. 

We will specialize to the main cases that arise in our study: projective
space, the intersection of two quadrics in $\pp^5$, and blow-ups
of abelian varieties. 

\

\noindent
{\bfseries Projective spaces:} \ Let us consider the case of 
projective space. Let $H$ denote the
hyperplane class on $\pp^n$, and as stated above by
convention we truncate to codimension $2$.  The Euler exact sequence
$$
0\rightarrow \Oo _{\pp^n}\rightarrow \Oo _{\pp^n}(1)^{n+1}
\rightarrow T\pp ^n\rightarrow 0
$$
gives
$$
{\rm ch}(T\pp^n) = n + (n+1)H + \frac{n+1}{2}H^2
$$
so 
$$
c_1(T\pp^n) = (n+1)H , \;\;\;\; c_2(T\pp^n) = \frac{n^2+n}{2}H^2.
$$
Thus
$$
\begin{aligned}
{\rm td}(T\pp^n) & = 1 + c_1/2 + (c_1^2+c_2)/12  \\[+0.5pc]
& = 1 + \frac{n+1}{2}H + \frac{3n^2 + 5n+2}{24}H^2.
\end{aligned}
$$

For $X=\pp^3$ this gives 
$$
{\rm td}(TX) = 1 + 2H + \frac{11}{6}H^2.
$$

\

\noindent
{\bfseries Intersection of quadrics:} \ Suppose $X\subset \pp^5$ a
complete intersection of two quadrics. Let $H$ also denote the
restriction of the hyperplane class of $\pp^{5}$ to $X$.  The normal
bundle of $X$ is $N_{X/\pp^5}\cong \Oo _X(2)^2$ so the exact sequence
$$
0\rightarrow TX \rightarrow T\pp^5|_X \rightarrow N_{X/\pp^5} 
\rightarrow 0
$$
gives 
$$
\begin{aligned}
{\rm ch}(TX) & = {\rm ch}(T\pp^5) - 2 {\rm exp}(2H) \\
& = 5 + 6H + 3H^2 - 2 - 4H -4H^2 \\
& = 3 + 2H -H^2.
\end{aligned}
$$
We get
$$
c_1(TX)= 2H,\;\;\;\; c_2(TX) = 3H^2. 
$$
Thus,
$$
\begin{aligned}
{\rm td}(TX) & = 1 + c_1/2 + (c_1^2+c_2)/12 \\
& = 1 + H + 7H^2 / 12.
\end{aligned}
$$

\

\noindent
{\bfseries Blow-ups of Pryms:} \ Suppose
$\Prym $ is a $3$-dimensional abelian variety, so its tangent
sheaf has trivial Chern classes, and we let $\blo :Y\rightarrow \Prym$ be
the blow-up along a smooth subvariety $A\subset \Prym$.  Let $E\subset
Y$ denote the exceptional divisor.  There are two cases of interest.

Suppose first $A$ is a disjoint collection of $a$ points. 
Then $E$ is a disjoint collection of $a$ planes $\pp^2$. 
A point $y\in E$ corresponds to a tangent direction
at the corresponding image point $\blo(y)\in  \Prym$, and the
tangent vectors on $Y$ at $y$ are vectors that are constrained
to map into this tangent direction in $T_{\blo(y)} \Prym$. 
Over $E$ there is the universal subbundle 
$\Oo _E(-1)=\Oo _E(E)$,
fitting into an exact sequence
$$
0\rightarrow \Oo _E(E)\rightarrow \blo^{*}T \Prym|_E \rightarrow 
R \rightarrow 0
$$
where the middle is a trivial bundle. 
Then, letting $R_Y$ be $R$ viewed 
as a coherent sheaf on $Y$ supported 
on $E$, we then have an exact sequence
$$
0\rightarrow TY \rightarrow \blo^{*}T \Prym \rightarrow R_Y \rightarrow 0.
$$
We have 
$$
{\rm ch}(R_Y) = {\rm ch}(\blo^{*}T \Prym ) {\rm ch}(\Oo _E) 
- {\rm ch}(\Oo _Y(E)){\rm ch}(\Oo _E).
$$
Now ${\rm ch}(\blo^{*}T\Prym )=3$ since the tangent sheaf of $\Prym $
is trivial of rank $3$. The exact sequence
$$
0\rightarrow \Oo (-E)\rightarrow \Oo \rightarrow \Oo _E\rightarrow 0
$$
gives
$$
{\rm ch}(\Oo _E) = E -E^2/2,
$$
so we get 
$$
\begin{aligned}
{\rm ch}(R_Y) & = (3 - (1+E + E^2/2)) (E-E^2/2) \\
& = 2E - 2E^2.
\end{aligned}
$$
Putting this into the other exact sequence we get
$$
{\rm ch}(TY) = 3 - 2E + 2E^2.
$$
We note that the shape of 
this formula is independent of the number $a$ of
points that were blown up (but of course the
full exceptional divisor $E$ contains that information). 

Turn now to the other case: when $A$ is a smooth irreducible
curve. The normal bundle $N_{A/\Prym }$ is a rank $2$ vector bundle
over $A$. Let $\delta$ denote the degree of the normal bundle,
i.e. the degree of its determinant line bundle. A point $y\in E$
corresponds to a normal direction to $A$ at $b(y)$, so as before the
tangent vectors to $Y$ at $y$ are constrained to have image in
$T_{b(y)}\Prym $ that maps into this normal direction under the map
$T\Prym |_A\rightarrow N_{A/\Prym }$.  Let $\Oo _E(-1)$ denote the
tautological sub-bundle over $E$. It fits into an exact sequence
$$
0\rightarrow \Oo _E(-1)\rightarrow \blo^{*}N_{A/\Prym } \rightarrow
R\rightarrow 0
$$
and as before, if $R_Y$ denotes the corresponding coherent
sheaf on $Y$ supported over $E$, then we have
the exact sequence
$$
0\rightarrow TY\rightarrow \blo^{*}T\Prym  \rightarrow R_Y
\rightarrow 0.
$$

The tautological subbundle is also the normal bundle
to $E$ in $Y$, thus $\Oo _E(-1)=\Oo _E(E)$ as in the
previous case. 

The rank $2$ bundle $N_{A/\Prym }$ over the curve 
$A$ is determined, up to rational numerical equivalence,
by its degree that we'll call $\delta$. 

Suppose $G$ is an ample divisor class on $\Prym $, then
$G\cdot A$ is the degree of $G$ restricted to $A$. Working
numerically, we can write
$$
{\rm det}N_{A/\Prym } \sim \frac{\delta}{G\cdot A} G|_A
$$
so $\blo^{*}N_{A/\Prym }$, considered as a
coherent sheaf on $Y$, is rationally 
numerically equivalent to the expression
(with a fractional divisor)
$(\blo^{*}\Oo _Y(\frac{\delta}{G\cdot A}G) \oplus \Oo _Y) \otimes 
\Oo _E$. This gives
$$
{\rm ch}(\blo^{*}N_{A/\Prym }) = 
\left(2 + \frac{\delta}{G\cdot A} G \right){\rm ch}(\Oo _E) .
$$
Here we can truncate the expression on the left at codimension
$1$ because it is multiplied by 
${\rm ch}(\Oo _E)$ of rank $0$. 

Recall as above that ${\rm ch}(\Oo _E)=E-E^2/2$. 
Similarly ${\rm ch}(\Oo _E(E))= E+E^2/2$ and from
the exact sequence,
$$
\begin{aligned}
{\rm ch}(R_Y) & = \left( 2 + \frac{\delta}{G\cdot A} G \right)\cdot (E-E^2/2) 
- E -E^2/2 \\
& = E - \frac{3}{2}E^2 + \frac{\delta}{G\cdot A} G\cdot E.
\end{aligned}
$$
From the other exact sequence and recalling that $\blo^{*}T\Prym $
is trivial of rank $3$ we get
$$
{\rm ch}(TY) = 3 - E + \frac{3}{2}E^2 - \frac{\delta}{G\cdot A} G\cdot E .
$$
The last term could perhaps more easily be understood 
as a sum of fibers: let $\mathsf{fib}$ denote the numerical
class of a fiber of the projection $E\rightarrow A$, viewed
as a codimension $2$ class on $Y$. We have
(writing also $G$ for its pullback to $Y$)
$$
\mathsf{fib} = \frac{1}{G\cdot A} G\cdot E
$$
so we can also write
$$
{\rm ch}(TY) = 3 - E + \frac{3}{2}E^2 - \delta\cdot \mathsf{fib}. 
$$
we get
$$
\begin{aligned}
c_1(TY) & = -E, \\
c_2(TY) & = \delta\cdot \mathsf{fib}-E^2,
\end{aligned}
$$
and 
$$
\begin{aligned}
{\rm td}(TY) & = 1 + c_1/2 + (c_1^2+c_2)/12 \\
& = 1 -E/2 + (\delta/12)\cdot {\rm fib}.
\end{aligned}
$$

\section{The degree one moduli space}
\label{chapter-d1}

For odd degree we have fixed the line bundle $\Oo _C(\pw )$ to be the
determinant, for the chosen Weierstrass point $\pw\in C$.  For
purposes of the present section, $X:= X_1$ is the moduli space of
stable bundles $E$ with ${\rm det}(E) \cong \Oo _C(\pw )$. Similarly,
$\Wob := \Wob _1$ denotes the wobbly locus.  Recall that $\Higgs _1$
denotes the moduli space of Higgs bundles with determinant $\Oo _C(\pw
)$ and trace of the Higgs field equal to zero.

Recall from Proposition \ref{prop:hitchin.map} that the $2:1$ spectral
covering $\pi : \Ctilde \rightarrow C$ has four branch points
$x,x',y,y'\in C$ such that $x$ and $x'$ are opposite under the
hyperelliptic involution, and $y$ and $y'$ are opposite under the
hyperelliptic involution. For a general point in the Hitchin base, the
resulting two points of the hyperelliptic $\pp^1$ are general with
respect to the choice of $C$.

In this chapter, various notations such as $X,Y,\Wob$ etc. are used
without the subscript $1$ since we are talking about the degree $1$
moduli space. For insertion into the rest of the paper one should read
$X_1,Y_1,\Wob _1$ and so on.

\subsection{Geometry of the wobbly locus in degree one}

We use the notations of subsections \ref{subsec-curves} and
\ref{logprop}.  In particular $\Higgs_{1}^{\cstar,{\rm nu}}$ denotes
the fixed point locus of the $\cstar$-action on $\Higgs_{1}$ and
$\Prym^{\rm unss}$ denotes the unstable locus in $\Prym$ which by  Theorem~\ref{qresolves} coincides with the intersection
of $\Fix\cap \Prym$ of $\Prym$ with the incoming variety $\Fix$ defined
in subsection \ref{logprop}. The geometric description of these loci
was already given in subsections \ref{subsec-curves} and \ref{logprop}
but we record it again in the following Lemma for ease of reference.

\begin{lemma}
  \label{geomw1}
  \begin{itemize}
\item $\Higgs_{1}^{\cstar,{\rm nu}}$ is equal to the curve $\Cbar
  \subset \Higgs_1$, the $16$-sheeted etale covering of $C$ defined in
  subsection \ref{subsec-curves}.
\item If $\Prym$ is a generic Hitchin fiber corresponding to spectral curve
$\Ctilde$, then $\Prym^{\rm unss} = \Fix\cap\Prym   \subset \Prym$ is
a smooth curve $\Chat$ that can be expressed as a fiber product
$$
\Chat = \Ctilde \times _C\Cbar .
$$
In particular, $\Chat $ is a $16$-sheeted \'{e}tale
covering of $\Ctilde$.
  \end{itemize}
These curves have genera respectively
$$
g_C=2, \;\;\; g_{\Ctilde }= 5, 
\;\;\; g_{\Cbar }=17, \;\;\;
g_{\Chat }= 65. 
$$
\end{lemma} 
\begin{proof}
The spectral curve $\Ctilde$ has genus $5$ and is a smooth double
cover $\Ctilde \rightarrow C$ ramified over $4$ points in two
conjugate pairs. This follows from our running genericity assumption
requiring that $\Ctilde$ is not ramified over any of the Weierstrass
points of $C$.  The Hitchin fiber $\Prym \subset \Higgs_{1}$
corresponding to such $\Ctilde$ is the Prym variety $\Prym = \Prym_3$ of line
bundles on $\Ctilde $ whose direct image down to $C$ is a rank $2$
bundle of determinant $\Oo_C(\pw )$.

The statement that $\Higgs_{1}^{\cstar,{\rm nu}} = \Cbar$ was already
proven in Proposition~\ref{descrip-fixed} but we recall the argument
here for completeness. A fixed point that is in
$\Higgs_{1}^{\cstar,{\rm nu}}$ (i.e. such that the underlying vector
bundle is not stable) is a Higgs bundle of the form $\left(A \oplus
A^{\vee}(\pw ),\theta\right)$, such that $\theta : A^{\vee}(\pw )
\rightarrow A\otimes \omega _C$.  Stability of this Higgs bundle means
$\deg(A) \leq 0$ and existence of the map implies that $\deg(A)=0$.
As $\omega_C=\Oo _C(2\pw )$ this means that $A^{\otimes 2}\cong \Oo _C(t
- \pw)$ for the point $t\in C$ where $\theta$ vanishes. This
description provides an isomorphism $\Higgs_{1}^{\cstar,{\rm nu}}\cong
\Cbar$.

The incoming variety $\Fix\subset \Higgs _1$ consists of Higgs bundles
$(E,\theta )$ such that \linebreak $\lim _{t\rightarrow 0}(E,t\theta) \in
\Higgs_{1}^{\cstar,{\rm nu}} \cong\Cbar $. In subsection~\ref{orbits}
we saw that the limiting map $\Fix\rightarrow \Cbar$ is well defined
and has fibers isomorphic to $\aaaa^3$. We claim that 
$$
\Prym^{\rm unss} = \Fix\cap \Prym 
$$
has the structure
$$
\Chat = \Ctilde \times _C \Cbar .
$$
In particular, it is seen to be a genus $65$ curve. 

For a general Hitchin fiber the intersection $\Prym^{\rm unss}$ is
smooth of dimension $1$.  A point of $\Prym^{\rm unss}$ is by
definition a line bundle $U$ of degree $3$ on $\Ctilde$ such that $\pi
_{*}U$ has determinant $\Oo_C(\pw )$ but is unstable. Instability
means that there is a sub-line bundle of degree $1$ that we will
denote $A^{\vee}(\pw )\hookrightarrow \pi _{*}U$.  By adjunction we
have $\pi^{*}(A^{\vee}(\pw )) \hookrightarrow U$, and
$\pi^{*}(A^{\vee}(\pw))$ has degree $2$. Thus, there is a unique point
$\tilde{t}\in \Ctilde$ such that $U =
\pi^{*}(A^{\vee}(\pw))(\tilde{t})$. If $t$ denotes the image of
$\tilde{t}$ in $C$, the determinant condition says that
$$
A^{\otimes -2}(2\pw + t) \otimes \omega_C^{-1}  \cong \Oo _C(\pw ).
$$
Thus, $A$ solves $A^{\otimes 2}\cong \Oo _C(t- \pw )$ and $(A,t)$
is the corresponding point of $\Cbar$.  This identifies $\Prym^{\rm
  unss}$ with the image of $\Chat = \Ctilde \times _C \Cbar$ under the
map \eqref{eq:ChatinY1}. But by Lemma~\ref{lemma:embedChat} this map
is a closed embedding, and so we have $\Prym^{\rm unss} = \Chat$.
\end{proof}

\

\begin{remark}
One can notice that the genera of the curves involved
here are of the form $4^k+1$. We do not know if that is
significant or not.
\end{remark}

\

\begin{lemma}
\label{overwob1}
Let $D$ denote the 
$\pp^1$-bundle over $\Cbar$ projectivization of the
bundle of outgoing directions. Then the map sending 
a point on $\Cbar$ together with an outgoing direction,
to the limit in $X$ of the resulting $\cstar$-orbit,
provides a map $D\rightarrow
\Wob$ that is finite and generically injective. Thus, $D$ is the
normalization of $\Wob$. 
\end{lemma}
\begin{proof}
Given a $\cstar$-fixed Higgs bundle $(A\oplus A^{\vee}(\pw ),\theta)$
corresponding to $(A,t)\in \Cbar$, the space of outgoing directions is
$\pp {\rm Ext}^1 (A^{\vee}(\pw ),A)$ which is by definition the fiber
of $D$ over $(A,t)$. As we saw in
Lemma~\ref{lem:normal.deg1}{\bfseries (b)} and we will see again by a
different argument in Proposition \ref{structureD} below, $D\cong
\Cbar \times \pp^1$. Each $\pp^1$ fiber of $D \to \Cbar$ maps to the
space of extensions which is the general form of a line in $X =
X_1$. Both the lines and the horizontal $\Cbar$'s map to positive
degree curves, and an effective divisor on the product $D$ is a
positive sum of vertical and horizontal ones, so any effective divisor
on $D$ maps to a positive degree curve in $X$, in particular it is not
contracted. For generic injectivity, we note that a general wobbly
bundle $E$ has a one dimensional space of nilpotent Higgs fields and
if we $\theta$ is any non-zero nilpotent Higgs field, we get a well
defined line subbundle $A = \ker \theta E$ and a realization of $E$ as
an extension of $A^{\vee}(\pw)$ by $A$. This determines a unque point
in $D$ which shows that the map $D \to \Wob$ has a rational
section. Since $D$ is irreducible this implies that $D \to \Wob$ is a
birational morphism. Alternatively, we can see the generic injectivity
from the geometry of the quadric line complex elucidated in
section~\ref{sec-linecomplex}. Recall that at a general point $w$ of
$\Wob$ the four lines in $X$ through that point include one with
multiplicity $2$.  That line corresponds to the point $A \in \Cbar
\subset \op{Jac}^{0}(C)$, and the position of $w$ on it corresponds to
the extension class, so we can recover the data of a point of $D$. It
follows that $D$ is the normalization of $\Wob$.
\end{proof}

\begin{lemma}
\label{wobbly8h}
The divisor class of the wobbly locus is $[\Wob ]=8H$. 
\end{lemma}
\begin{proof}
We will consider  lines $\ell \subset X$ and show that
$\ell$ intersects $W$ in $8$ points.  In the moduli point of view on
$X$ a line corresponds to fixing a line bundle $L$ of degree $0$ and
looking at the set of bundles $E$ fitting into an extension
$$
0\rightarrow L \rightarrow E \rightarrow L^{\vee}(\pw) \rightarrow 0.
$$ Such an extension is wobbly if there is $A\rightarrow E$ and $t\in
C$ with $A^{\otimes 2} = \Oo (t-\pw)$. Assuming that $L$ does not
satisfy this condition, it means that $A = L^{\vee}(\pw -y)$, so
$L^{\otimes 2} (2y-2\pw ) = \Oo (\pw -t)$ in other words $L^{\otimes
  2} = \Oo (3\pw -t-2y)$, and the pullback of the extension to $A$
splits. For a choice of $(y,t)$ there will be a unique extension up to
scalars, i.e. a unique point on the line, such that the pullback
extension splits, so we need to look at the number of solutions of the
equation $L^{\otimes -2}(3\pw ) = \Oo (2y + t)$. This is the number of
branch points of a trigonal map  $C\rightarrow \pp^1$, which is
$8$. Thus $\ell \cap \Wob$ has $8$ points for a general line
$\ell$. We will discuss the trigonal covers more in subsection
\ref{trigonal} below.
\end{proof}

\

The basic structure of $\Wob $ is captured by the following
proposition. 

\

\begin{proposition}
\label{structureD}
We have $D\cong \Cbar  \times \pp^1$ and the
factor $\pp^1$ is naturally identified with the base 
of the hyperelliptic covering $\pp^1=C/\hi_{C} $. The map
$D\rightarrow \Wob$ identifies together pairs of points on
$6$ curves of the form $\Cbar \times \{ \brx_i\}$,
yielding a nodal locus in $\Wob $. Furthermore, on the
graph of the natural projection $\Cbar \rightarrow \pp^1$
the map $D\rightarrow X$ is non-immersive, yielding a
cuspidal locus in $\Wob $. 
\end{proposition}
\begin{proof}
Recall that $\brx_1,\ldots , \brx_6\in \pp^1$ denoted the branch points
for the hyperelliptic map $\hyp_{C} : C \to \pp^{1}$, i.e. the images
of the Weierstrass points $\pw_1,\ldots , \pw_6\in C$; also we chose
$\pw :=\pw_1$ to fix the determinant $\det(E)\cong \Oo _C(\pw )$ of
our stable vector bundles.

The statement that the $\pp^{1}$-fiber of $D$ is naturally identified
with the hyperelliptic $\pp^{1}$ was already checked in
Lemma~\ref{lem:normal.deg1}{\bfseries (b)} but we 
recall the argument here for convenience. A point in $\Cbar $ is a
line bundle $L$ such that $L^{\otimes 2}(\pw)$ is effective, or
equivalently such that $L^{\otimes 2} \cong \Oo _C(t-\pw )$ for some
point $t \in C$. Thus $L$ defines uniquely the point $t$.

A point in $D$ is a pair $(L,\eta )$ where 
$L\in \Cbar $ and $\eta \in \pp \op{Ext}^1(L^{\vee}(\pw ),L)$,
defining the extension
$$
0\rightarrow L\rightarrow E\rightarrow L^{\vee}(\pw )
\rightarrow 0.
$$
We note that
$$
\op{Ext}^1(L^{\vee}(\pw ),L) = H^1(L^{\otimes 2}(-\pw ))
\cong H^1(\Oo _C(t-2\pw ))\cong H^1(\Oo _C(-t'))
$$
where $t'$ is the image of $t$ under the 
hyperelliptic involution. 
The map $\Oo _C(-t')\rightarrow \Oo_{C}$
induces an isomorphism on $H^1$ so we can write
$$
\op{Ext}^1(L^{\vee}(\pw ),L)\cong H^1(\Oo _C).
$$
The scalar multiplying this isomorphism is not canonical,
it depends on something and indeed that leads (as we
have discussed) to the degree $128$ of the normal
bundle of $\Chat $ in $\Prym $. However, it
does give a canonical isomorphism 
$$
\pp \op{Ext}^1(L^{\vee}(\pw ),L)\cong \pp H^1(\Oo _C) =\pp^1.
$$
We note that if $y\in C$ then the map
$$
H^1(\Oo _C)\rightarrow H^1(\Oo _C(y))
$$
is a quotient of rank $1$ corresponding to the inage $\hyp_{C}(y) \in 
\pp^1$  of $y$. This is the identification between
the $\pp^1$ factor of $D$, and the base of the
hyperelliptic projection of $C$. 

We would like to understand when two points of $D$ correspond
to the same bundle. Suppose given $(L,\eta )\in D$,
corresponding to an extension as above. We look for a
new line subbundle $L_1\subset E$ such that $L_1\in 
\Cbar $. We note first of all that if $L \neq L_{1}$ as subbundles in $E$,
then
there is an inclusion
$$
L_1\hookrightarrow L^{\vee}(\pw ),
$$
so we may write
$$
L_1 = L^{\vee}(\pw -y).
$$
We have 
$$
L_1^{\otimes 2} = L^{\otimes -2}(2\pw -2y) = \Oo _C(\pw - t)(2\pw -2y)
=\Oo _C(3\pw-t-2y)
$$
so the condition $L_1\in \Cbar$
becomes 
$$
\Oo _C(3\pw -t-2y)\cong \Oo _C(t_1-\pw )
$$
for some point $t_{1} \in C$.  In other words, 
$$
4\pw \sim t + t_1 + 2y.
$$
One may notice that for our hyperelliptic curve,
any time that $a+b+c+d \sim 2\omega_C = 4\pw $,
we have $a+b+c+d = x+x'+y+y'$ such that $x,x'$ (respectively $y,y'$)
are paired under the hyperelliptic involution $\hi_{C}$. 

Given the form $t+t_1+2y$, there are only two possibilities:
\begin{itemize}

\item[(i)]  either $t = t_1$, and $t$ is $\hi_{C}$-paired with $y$, 

\item[(ii)] or $(t,t_1)$ are $\hi_{C}$-paired, and $y=\pw_i$ is a Weierstrass
point on $C$. 

\end{itemize}

\

\noindent
{\bfseries Case (i):} \ If we are in the first case we get
$$
L_1=L^{\vee}(\pw -y) = L^{\vee}(\pw -(2\pw - t)) = L^{\vee}(t -\pw )
$$ and so we must have $L_1\cong L$, that is $L$ and $L_{1}$
correspond to two different embedings of the same line bundle in
$E$. However, for $L_1 \cong L$ the space $\op{Ext}^{1}(L_{1},L)$ of
extensions of $L_1$ by $L$ has dimension $2$ (it is $H^1(\Oo
_C)$. Furthermore, the map on extension sapces 
\[
\op{Ext}^1(L^{\vee}(\pw),L) \to
\op{Ext}^{1}(L_{1},L)
\]
induced by the inclusion $L_{1} \hookrightarrow L^{\vee}(\pw)$,
is identified with the map on cohomology 
\[
H^{1}(L^{\otimes 2}(-\pw)) = H^{1}(\mathcal{O}_{C}(t - 2\pw))
= H^{1}(\mathcal{O}_{C}(-t')) \to H^{1}(\mathcal{O}_{C})
\]
induced by the inclusion $\mathcal{O}_{C}(-t') \hookrightarrow
\mathcal{O}_{C}$, which is an isomorphism.  Hence if $L_{1} \cong L$,
then no choice of bundle $E$ admits a lifting from $L_1$ into
$E$. Therefore, case (i)  does not lead to any
identification.

As we will see below, it will however lead to an identification
infinitesimally, yielding the cuspidal locus. We ignore case (i)  for now.

\

\noindent
{\bfseries Case (ii):} \ In this case, we have $y=p_i$ and $t_1=t'$. 
This gives a linear equivalence $2\pw \sim 2y$ and we get 
$$
L_1= L^{\vee}(\pw -p_i) = L^{\vee}\otimes \mathfrak{a}_i
$$
where $\mathfrak{a}_i=\Oo_C(\pw -\pw_i)$ is a $2$-torsion point. There are
six possibilities in this case. Let us look at what
is the extension $\eta$ such that $L_1$ lifts into $E$. 
The extension $\eta$ should be in the kernel of 
the map 
$$
\op{Ext}^1(L^{\vee}(\pw ),L)\rightarrow \op{Ext}^1(L_1, L)
$$
i.e. the  map
$$
H^1(\Oo (-t'))\rightarrow H^1(\Oo (-t'+\pw_i)).
$$
We note that this map is isomorphic to the map
$$
H^1(\Oo _C)\rightarrow H^1(\Oo _C(\pw_i))
$$
so it means that $\eta$ should be the extension 
corresponding to the branch point $\brx_i\in \pp^1$. 

We obtain $6$ glueings, one for each point $\brx_i$, with the point
$(L,\brx_i)$ being glued to $(L\otimes\mathfrak{a},\brx_{i})$. In
other words the surface $D$ is glued to itself along each of the $6$
curves $\Cbar\times \{\brx_{i}\}$, $i = 1, \ldots, 6$. On each curve
$\Cbar\times \{\brx_{i}\}$, the glueing map is the automorphism of $\Cbar$
given by tensoring with the $2$-torsion point $\mathfrak{a}_{i}$. This
is the first part of the statement in the proposition.

The analysis of the cuspidal locus, to conclude the proof of the
proposition, will be done later in subsection \ref{app-proof-cusp}.
\end{proof}

\subsection{Computations in degree one}

By convention all answers to Chern calculations are truncated to the
terms of degree $\leq 2$.

By Lemma~\ref{lemma:embedChat}) and Theorem~\ref{qresolves} the curve
$\Chat$ embeds in $\Prym$ and the blow-up $\blo : Y\rightarrow \Prym$
of $\Prym$ along $\Chat$ resolves the natural rational map $\pi_{*}(-)
: \Prym \dashrightarrow X$ producing a finite morphism $f : Y \to X$.
Let $\ExY \subset Y$ denote the exceptional locus. It is a
$\pp^1$-bundle over $\Chat$, namely $\ExY \cong \pp (N)$ (this being
the projective bundle of subspaces), where $N:= N_{\Chat /\Prym }$. In
fact by Lemma~\ref{lem:normal.deg1}{\bfseries (b)}
we have $N = H^{1}(\mathcal{O}_{C})\otimes
\sqhat^{*}\pi^{*}\omega_{C}$ and so $\ExY = \Chat\times \pp^{1}$ where
the second factor is the hyperelliptic $\pp^{1}$ for $C$.

\begin{remark}
\label{line}
The images in $X$ of the $\pp^{1}$'s that are fibers of $b_{\ExY}:\ExY
\rightarrow \Chat $ are lines in $X$.  Indeed, these fibers are
isomorphic to the projectivized bundle of outgoing directions along
$\Cbar$ which map to lines in $\Wob$, see Lemma \ref{overwob1}.
\end{remark}

\

Let $H$ denote the hyperplane class on $X$. Since $X\subset \pp^5$
is a degree $4$ subvariety of dimension $3$, we get $H^3=4$.

Let $\FxY =f ^{*}H$ be the inverse image of the class $H$ on $Y$.
We know that for a generic Hitchin fiber the Neron-Severi group of
$\Prym $ is $\zz$, so the Neron-Severi group of $Y$ is 
generated rationally by the classes $\ExY $ and $\FxY $.

\begin{proposition}
\label{todd1}
The Todd classes for the degree $1$ moduli spaces are: 
$$
{\rm td}(TX) = 1 + H + 7H^2 / 12, 
$$
$$
\pi ^{*}{\rm td}(TX)^{-1} = 1-\FxY   + 5\FxY ^2/ 12, 
$$
$$
{\rm td}(TY) = 1 -\ExY /2 + (\ExY ^2+\ExY \FxY )/9, 
$$
$$
{\rm td}(Y/X) = (1-\FxY   + 5\FxY ^2/ 12)(1 -\ExY /2 + (\ExY ^2+\ExY \FxY )/9). 
$$
\end{proposition}
\begin{proof}

Recall from Section \ref{sec-Chern}, the Chern class calculations for
the intersection of two quadrics $X$ say
$$
c_1(TX)= 2H,\;\;\;\; c_2(TX) = 3H^2, 
$$
and
$$
{\rm td}(TX) = 1 + c_1/2 + (c_1^2+c_2)/12 
$$
$$
= 1 + H + 7H^2 / 12.
$$
A computation shows that 
$$
(1 + H + 7H^2 / 12) (1-H  + 5H^2/ 12) = 1, 
$$
thus
$$
{\rm td}(TX)^{-1} = 1-H  + 5H^2/ 12
$$
which pulls back to the same formula with $\FxY $ on $Y$. 

For the blow-up $Y$ of an abelian threefold $\Prym $ along
the curve $\Chat $ we need to know 
the degree of the normal bundle. This degree is 
$2g_{\Chat }-2$ for the embedding of a curve
in an abelian  variety, thus by either the explicit description of $N$ in
Lemma~\ref{lem:normal.deg1}{\bfseries (b)} or by the genus calculation
of Lemma \ref{geomw1} we have
$$
\delta = {\rm deg}(N) = 128. 
$$
This could also be seen from the intersection number
calculations above, one can see that it is the same as $-\ExY ^3$. 

Again from the calculations 
of Section \ref{sec-Chern} we have
$$
c_1(TY) = -\ExY ,
$$
$$
c_2(TY) = \delta \mathsf{fib}-\ExY ^2,
$$
and 
$$
{\rm td}(TY) = 1 + c_1/2 + (c_1^2+c_2)/12
$$
$$
= 1 -\ExY /2 + \delta \mathsf{fib}/12.
$$
We note that $(\ExY +\FxY)|_{\ExY }$ is a divisor on $\ExY $ whose
intersection with a fiber is zero, indeed a fiber intersect $\ExY $ is
$-1$ and it intersects $\FxY $ in $1$ point since the image of a fiber
is a line contained in $X$ (Remark \ref{line}).  Therefore $(\ExY
+\FxY )\cdot \ExY $ is a sum of fibers. The number may be calculated
as $(\ExY +\FxY )\cdot\ExY\cdot \FxY $ again from Remark \ref{line},
this gives
$$
(\ExY +\FxY )\cdot \ExY  = 96\mathsf{fib}.
$$
Therefore $\delta {\rm fib} = (128/96)(\ExY +\FxY)\cdot
\ExY  = 4(\ExY ^2+\ExY \FxY )/3$.
Our formulas become
$$
c_2(TY) = 4(\ExY ^2+\ExY \FxY )/3-\ExY ^2 = (\ExY ^2 + 4\ExY \FxY )/3,
$$
and 
$$
{\rm td}(TY) =  1 -\ExY /2 + (\ExY ^2+\ExY \FxY )/9 
$$
which completes the calculation.
\end{proof}

\

\noindent
Next we compute the triple intersections of divisor classes on $Y$.

\

\begin{proposition}
\label{intersections1}
The triple intersections of divisor classes on the degree $1$ modular
spectral
covering $Y$ are:
$$
\FxY ^3 = 32, \;\;\;\; 
\ExY \FxY ^2 = 64, \;\;\;\; 
\ExY ^2\FxY  = 32, \;\;\;\; 
\ExY ^3 = -128.
$$
\end{proposition}

\

\noindent
Recall that $H$ denotes the divisor class of the hyperplane section 
on $X$, and $\FxY := f ^{*}(H)$. 
We have 
$$
\begin{aligned}
H^{1,1}(Y)\cap H^2(Y,\qq ) & = \langle \ExY ,\FxY \rangle, \text{ and } \\
H^{2,2}(Y)\cap H^4(Y,\qq ) & = \langle \ExY \FxY ,\FxY ^2\rangle.
\end{aligned}
$$
The corresponding groups on $X$ are generated 
(over $\qq$) by
$H$ and $H^2$ respectively.

\begin{lemma}
\label{deg8c1}
The map $f : Y\rightarrow X$ has degree $8$.
\end{lemma}
\begin{proof}
Since over the very stable locus the map $\pi_{*}(-) : \Prym
\dashrightarrow X$ is a proper morphism (see
e.g. \cite{PalPauly,Zelaci,PeonNieto} it suffices to compute
the number of preimages of a very stable point in $X$ under $f : Y \to
X$.  Choose a general hence very stable point $E\in X$. The fiber of
the cotangent bundle $T^{\vee}_EX$ is the space
$H^0(\op{End}_{0}(E)\otimes \omega_C)$ of $\omega_C$-twisted
endomorphisms. The Hitchin base is the $3$-dimensional space
$H^{0}(C,\omega_{C}^{\otimes 2}) \cong \aaaa^3$ of quadratic
differentials on $C$. The map $T^{\vee}_EX \rightarrow \aaaa^3$ is the
restriction of the Hitchin map and is thus given by three quadratic
forms, so the inverse image of a general point $\mathsf{b}\in \aaaa^3$
is the intersection of three quadrics: it has $8$ points.  As $\Prym
=\hit^{-1}(\mathsf{b})$, and the exceptional divisor $\ExY$ maps to
the wobbly locus, we have that the intersection of $Y\cap T^{\vee}_EX$
is equal to $\Prym\cap T^{\vee}_{E}X$ and so is this set of $8$
points. Thus for a very stable $E$ and a generic $\Ctilde$ the inverse
image $f^{-1}(E)\subset Y$ consists of $8$ points.
\end{proof}

\

\begin{remark} \label{rem:deg8X0} Note that the proof
of Lemma~\ref{deg8c1} repeats verbatim to show that in the degree
  $0$ case the map $f_{0} : Y_{0} \to X_{0}$ has degree $8$ as well.
\end{remark}

\

\noindent
We have $H^3=4$ since $X\subset \pp^5$ is a degree $4$
subvariety. Hence,
$$
\FxY ^3 = 32. 
$$
We also  checked that 
$\Oo _Y(\ExY +\FxY )|_{\ExY }$ is a sum
of $96$ fibers of the map $\ExY \to \Chat$. In particular,
the divisor $(\ExY +\FxY )$ restricted to
$\ExY $, has trivial self-intersection. We get the formula
$$
(\ExY +\FxY )\cdot \ExY  \cdot (\ExY +\FxY ) = 0,
$$
or 
$$
\ExY ^3 + 2 \ExY ^2 \FxY  + \ExY \FxY ^2 = 0.
$$

\

\begin{lemma}
The wobbly locus is the image  
$\Wob =f (\ExY )\subset X$. The intersection of $\Wob $ with a line
$\ell$ in $X$ (in other words a subvariety
that is a line in $\pp^5$) has $8$ points. 
The map $f_{|\ExY }:\ExY \rightarrow \Wob $ is generically
$2$ to $1$, and if $D$ denotes the normalization of $\Wob $
(see Lemma \ref{overwob1} and Proposition \ref{structureD} above) 
this map factors through a $2:1$ map
$\ExY \rightarrow D$. 
\end{lemma}
\begin{proof}
Recall from Lemma \ref{overwob1} that $D$ can be described as the
space of downward or outgoing directions to the curve $\Cbar \subset
\Higgs _1$. As we saw in subsections \ref{orbits} and \ref{logprop}
the incoming flow gives a $2:1$ map $\Chat \rightarrow \Cbar$ and a
normal direction to a general point of $\Chat$ maps to a downward
direction at the image point in $\Cbar$. This gives the map $\ExY
\rightarrow D$. In the present case, downward flow lines are broken
only once and the broken flow line depends only on the first normal
derivative to the point of $\Chat$ in the Hitchin fiber, so the
closure of the downward flow map restricts on $\ExY $ to the
composition $\ExY \rightarrow D \rightarrow \Wob \subset X_1$.
\end{proof}

\

\begin{proof}[Proof of Proposition \ref{intersections1}]
We saw above that $\FxY ^3 = 32$. 

A line $\ell$ lying on $X$ satisfies $\ell 
\cap H= 1$, but also $(H^2)\cap H = 4$; hence in the group $H^{2,2}(Y)\cap
H^4(Y,\qq )$ we have
$$
H^2 = 4 \ell . 
$$
This gives $H^2 \cap \Wob = 32$. 

From the $2:1$ covering property we get
$$
\FxY ^2 \cap \ExY = H^2 \cap 2\Wob 
$$
so 
$$
\ExY \FxY ^2 = 64.
$$
Let us now calculate $\ExY ^3$. It is the self-intersection of the
divisor class $c_{1}(\mathcal{O}_{\ExY}(\ExY))$ on the surface $\ExY$.
By the usual blowup picture the line bundle $\mathcal{O}_{\ExY}(\ExY)$
is the universal subbundle for the projectivization $\mathbb{P}(N) = \ExY$ 
so the universal quotient bundle is
$$
UQ = b_{\ExY }^{*}(N) / \Oo _{\ExY }(\ExY ).
$$
In terms of divisor classes on $\ExY$ this gives 
$$
c_{1}(UQ) = b_{\ExY}^{*}c_{1}(N) - c_{1}(\mathcal{O}_{\ExY}(\ExY).
$$
The relative tangent bundle $T(\ExY/\Chat)$ of the map $b_{\ExY} :
\ExY \to \Chat$ is given by
$$
T(\ExY /\Chat )= {\rm Hom}(\Oo _{\ExY }(\ExY ), UQ),
$$
which gives
$$ c_{1}(T(\ExY /\Chat)) = c_{1}(UQ) - c_{1}(\Oo _{\ExY }(\ExY)) =
b_{\ExY}^{*}c_{1}(N) - 2 c_{1}(\Oo _{\ExY }(\ExY))
$$
on the level of divisor classes.

Since $\ExY \cong \Chat\times \pp^{1}$, the Neron-Severi group of
$\ExY$ is freely generated by two classes $\hat{\mathfrak{c}}$ and
$\mathsf{fib}$, where $\hat{\mathfrak{c}}$ is the class of
$\Chat\times \op{pt}$ and as before $\mathsf{fib}$ is the class of
$\op{pt}\times \pp^{1}$. In terms of these classes we have
$$
b_{\ExY}^{*}c_{1}(N) = 128\mathsf{fib}
\qquad \text{and} \qquad  c_{1}(T(\ExY /\Chat)) = 2\hat{\mathfrak{c}}.
$$
Therefore we get
$$
c_{1}(\mathcal{O}_{\ExY}(\ExY)) =
\frac{1}{2}\left( b_{\ExY}^{*}c_{1}(N) - c_{1}(T(\ExY /\Chat))\right)
= 64\mathsf{fib} - \hat{\mathfrak{c}},
$$
and so taking into account that $\hat{\mathfrak{c}}^{2} =
\mathsf{fib}^{2} = 0$ and $\hat{\mathfrak{c}}\cdot\mathsf{fib} = 1$ we
get
$$
c_{1}(\mathcal{O}_{\ExY}(\ExY))\cdot c_{1}(\mathcal{O}_{\ExY}(\ExY))
= (64\mathsf{fib} - \hat{\mathfrak{c}})^{2} = -128.
$$
That was the self-intersection on $\ExY $, so it gives
on $Y$, $\ExY ^3=-128$.

Now from the formula $\ExY ^3 + 2 \ExY ^2 \FxY  + \ExY \FxY ^2 = 0$ 
we get 
$$
(-128) + 2\ExY ^2\FxY + (64) = 0, \mbox{ so } \ExY ^2\FxY =
(128-64)/2 = 32.
$$
This completes our list of intersection numbers on $Y$.
\end{proof}

\

\subsection{Main construction}
Let $\Ll_0$ be the pullback to $Y$ of a degree zero line bundle on
$\Prym$, and set
$$ \Ll _{a,b} := \Ll _0(a\FxY + (b+1)\ExY ) \ \text{and} \ \Vv
_{a,b}:= f_{*}(\Ll _{a,b}).
$$
Following the notation we adopted in the definition of crude
parabolic structures in section~\ref{subsec-parabolic} we will also
set
$$
\Ll ' _{a,b}:= \Ll _0(a\FxY  + b\ExY ) \ \subset \ \Ll _{a,b} \ \text{and} \
\Vv'_{a,b}= f_{*}(\Ll'_{a,b}) \ \subset \ \Vv _{a,b}.
$$ The notations are chosen so that the spectral line bundle defined
in section~\ref{subsec-parabolic} is precisely $\Ll' = \Ll_0(a\FxY
+b\ExY)$. Notice that $\Vv'_{a,b} = \Vv_{a,b-1}$.

\

\noindent
The next step is to look at the calculations of the Chern characters
of the direct image \linebreak $\Vv_{a,b}=f_{*}(\Ll_{a,b})$.  The
Chern character of $\Ll_{a,b}$ is
$$
{\rm ch}(\Ll_{a,b}) =1+(a\FxY +
(b+1)\ExY ) + (a\FxY + (b+1)\ExY )^2/2.
$$
Putting the Todd class formula of Proposition \ref{todd1} into the
Grothendieck-Riemann-Roch formula gives
$$
\begin{aligned}
{\rm ch}(\Vv _{a,b}) & = {\rm ch}(f _{*}(\Ll _{a,b}))  \\
& = f _{*} (1+(a\FxY  + (b+1)\ExY )  \\
& \hspace{0.5in}
+ (a\FxY  + (b+1)\ExY )^2/2)(1-\FxY   + 5\FxY ^2/ 12)
(1 -\ExY /2 + (\ExY ^2+\ExY \FxY )/9) \\
& = f _{*}\left[ 
1+((a-1)\FxY  + (b+1/2)\ExY )  + (a\FxY  + (b+1)\ExY )^2/2 \right. \\
& \hspace{0.5in} \left.
+ 5\FxY ^2/ 12 +(\ExY ^2+\ExY \FxY )/9 + \ExY \FxY /2 - 
(\FxY + \ExY /2)(a\FxY  + (b+1)\ExY )
\right].
\end{aligned}
$$
The rational cohomology of $X_1$ is generated by $H$ in degree $2$
and $H^2$ in degree $4$, so if $a_1 \in H^{2}(Y)$ and $a_2 \in
H^{4}(Y)$ are classes on $Y$ we can write
$$
f _{*}(1 + a_1 + a_2)= 8 +  b_1H + b_2H^2,
$$
for some rational numbers $b_{1}, b_{2} \in \mathbb{Q}$.
Then, using $H^3=4$ we get 
$4b_1 = H^2\cdot f _{*}(a_1) = \FxY ^2a_1$ and 
$4b_2 = H\cdot f_{*}(a_2) = \FxY a_2$. 

We get, using the calculations of Proposition \ref{intersections1}:
$$
H^2\cdot {\rm ch}_1(\Vv _{a,b}) 
=
\FxY ^2 ((a-1)\FxY  + (b+1/2)\ExY )  =32 (a-1) + 64(b+1/2) = 32a + 64 b
$$
so 
$$
{\rm ch}_1(\Vv _{a,b}) = (8a + 16 b)H.
$$
And
$$
\begin{aligned}
H\cdot {\rm ch}_2(\Vv _{a,b}) & =
\FxY \left[
  (a\FxY  + (b+1)\ExY )^2/2 \right.
\\
& \ \hspace{0.5in} \left.  
+ 5\FxY ^2/ 12 +(\ExY ^2+\ExY \FxY )/9 + \ExY \FxY /2 - 
(\FxY + \ExY /2)(a\FxY  + (b+1)\ExY )
\right] \\
& = (a^2/2 + 5/12 - a)\FxY ^3
+ (a(b+1) + 1/9 + 1/2 - (b+1 + a/2))\ExY \FxY ^2 \\
& \ \hspace{0.5in} + ((b+1)^2/2 + 1/9 - (b+1)/2) \ExY ^2 \FxY \\
& = 16 a^2 - 32 a + 32 a + 64 ab + 16 b^2 + - 64 b + 16b 
+ 40/3 + 64/9 - 32 + 32/9 \\
& = 16 a^2 + 16 b^2 + 64 ab -48b - 8.
\end{aligned}
$$
We get 
$$
{\rm ch}_2(\Vv _{a,b}) = 4a^2 + 16 ab + 4 b^2 - 12 b - 2.
$$
Together these prove the following:

\begin{proposition}
\label{chv}
For $\Ll _{a,b}$ in the numerical class of $\Oo _Y(a\FxY  + (b+1)\ExY )$
the direct image $\Vv _{a,b}= \pi _{*}(\Ll _{a,b})$
has Chern character (truncated to codimension $2$ as usual)
$$
{\rm ch}(\Vv _{a,b}) = 8 +(8a + 16 b)H + 
(4a^2 + 16 ab + 4 b^2 - 12 b - 2)H^2.
$$
\end{proposition}

\

\noindent
Later we will consider a parabolic modification of $\Ll_{a,b}$ over
the divisor $\ExY $, giving a new smooth parabolic structure defined
locally upstairs over the modular spectral cover or somewhat
equivalently over the projectivization of the bundle.

\subsection{The cusp locus of the wobbly divisor}
\label{app-proof-cusp}

In this subsection we prove the part of Proposition \ref{structureD}
about the cuspidal locus. Specifically we compute the locus $D_{\rm
  cusp} \subset D$ on which the map $\fD : D \to X$ is not
immersive. Recall that this map is induced from the map $f : Y \to
X$. More precisely, the restriction $f_{|\ExY} : \ExY \to \Wob \subset
X$ factors through the double cover $\ExY = \Chat\times \pp^{1} \to
\Cbar\times \pp^{1} = D$ and $\fD : D \to X$ is just the induced
map.

Let
$(L,\eta) \in D = \Cbar \times \mathbb{P}^{1}$ be a point, and
let $E = \fD(L,\eta) \in X$ be the corresponding rank two bundle. Our
goal is to understand when the differential
  \[
  d\fD_{(L,\eta)} : T_{D, (L,\eta)} \to T_{X,E}
  \]
has a non-trivial kernel.

Recall that $L \in \Cbar \subset \op{Jac}^{0}(C)$ means that
$L^{\otimes 2}(\pw )$ is effective, i.e. that $L^{\otimes 2} \cong
\mathcal{O}_{C}(t-\pw )$ for some point $t \in C$. Thus $L^{\otimes
  2}(-\pw ) \cong \mathcal{O}_{C}(t - 2\pw ) = \mathcal{O}_{C}(-t')$,
where $t' \in C$ is the point corresponding to $t$ under the
hyperelliptic involution.

Hence we have 
\[
\op{Ext}^{1}(L^{\vee}(\pw ),L) = H^{1}(C,L^{\otimes 2}(-\pw )) =
H^{1}(C,\mathcal{O}_{C}(-t')). 
\]
From the long exact sequence in cohomology associated with the short
exact sequence \linebreak $0 \to \mathcal{O}_{C}(-t') \to \mathcal{O}_{C} \to
\mathcal{O}_{t'} \to 0$  it follows that the natural map
\[
\xymatrix@1{
  \iota_{t'} :  \hspace{-0.5pc}
  & H^{1}(C,\mathcal{O}_{C}(-t')) \ar[r]^-{\cong} & H^{1}(C,\mathcal{O}_{C})}
\]
is an isomorphism, 
and hence we can view $\eta$ as an element in $\iota_{t'}^{-1}(\eta)
\in \text{Ext}^{1}(L^{\vee}(\pw ),L)$. The rank two vector bundle $E =
f(L,\eta)$ is the extension
\[
0 \to L \to E \to L^{\vee}(\pw ) \to 0 
\]
given by the extension class  $\iota_{t'}^{-1}(\eta)$.

We have $T_{D,(L,\eta)} = T_{\Cbar ,L}\oplus
T_{\mathbb{P}^{1},\eta}$ and so to understand the map $df_{(L,\eta)}$
it is enough to understand its restrictions to the two coordinate
lines $T_{\Cbar ,L}\oplus \{ 0 \}$ and $\{ 0 \} \oplus 
T_{\mathbb{P}^{1},\eta}$ in $T_{D,(L,\eta)}$.
These restrictions which will abbreviate as 
\[
\xymatrix@R-2pc{
d\fD_{(L,\eta)|T_{\Cbar ,L}} : \hspace{-0.5pc} &
T_{\Cbar ,L}\ar[r] & T_{X,E}\\
d\fD_{(L,\eta)|T_{\mathbb{P}^{1},\eta}} : \hspace{-0.5pc} & 
T_{\mathbb{P}^{1},\eta} \ar[r] & T_{X,E} }
  \]
have a natural modular interpretation.

\subsubsection{Interpretation of
  \texorpdfstring{$d\fD_{(L,\eta)|T_{\mathbb{P}^{1},\eta}}$}{dfDTP}}
\label{ssec-resP1}

Fix $L \in \Cbar $ and hence the points $t, t' \in C$.  The
extensions of $L^{\vee}(\pw )$ by $L$ corresponding to a varying point
$\eta \in \mathbb{P}(H^{1}(C,\mathcal{O}_{C}))$ fit in a natural
universal family parametrized by $\mathbb{P}^{1}$.  Indeed, consider
the surface $C\times \mathbb{P}^{1}$. By K\"{u}nneth we have
\[
\begin{aligned}
  \text{Ext}^{1}_{C\times
    \mathbb{P}^{1}}(p_{C}^{*}L^{\vee}(\pw ),p_{C}^{*}L\otimes
  p_{\mathbb{P}^{1}}^{*}\mathcal{O}(1)) & = H^{1}(C\times
  \mathbb{P}^{1},p_{C}^{*} \mathcal{O}(-t')\otimes
  p_{\mathbb{P}^{1}}^{*}\mathcal{O}(1)) \\
  & =
  H^{1}(C,\mathcal{O}(-t'))\otimes H^{0}(\mathbb{P}^{1},\mathcal{O}(1)) \\
  & \cong H^{1}(C,\mathcal{O})\otimes H^{1}(C,\mathcal{O})^{\vee}.
\end{aligned}
\]
Here in the last step we used $\iota_{t'}$ to identify
$H^{1}(C,\mathcal{O}(-t'))$ with $H^{1}(C,\mathcal{O})$ and we used
$\mathbb{P}^{1} = \mathbb{P}(H^{1}(C,\mathcal{O}))$ to identify
$H^{0}(\mathbb{P}^{1},\mathcal{O}(1))$ with $H^{1}(C,\mathcal{O})^{\vee}$.

The extension
\[
0 \to p_{C}^{*}L\otimes p_{\mathbb{P}^{1}}^{*}\mathcal{O}(1) \to
\mycal{E} \to p_{C}^{*} L^{\vee}(\pw ) \to 0
\]
corresponding to the identity element in $H^{1}(C,\mathcal{O})\otimes
H^{1}(C,\mathcal{O})^{\vee}$ is a rank two bundle $\mycal{E}$ on $C\times
\mathbb{P}^{1}$ whose restrictiion to $C\times \{\xi\}$ is the extension
\[
0 \to L \to \mycal{E}_{\xi} \to L^{\vee}(\pw ) \to 0 
\]
given by the class $\iota_{t'}^{-1}(\xi)$. In particular we have
$\mycal{E}_{\xi} = E$.

Therefore $\fD_{|\{L\}\times \mathbb{P}^{1}} : \mathbb{P}^{1} \to X$ can
be viewed as the classifying map for the bundle \linebreak $\mycal{E}
\to C\times \mathbb{P}^{1}$. Let $\text{End}_{0}(E)$ denote the bundle
of traceless endomorphisms of $E$. By deformation theory we have
$T_{X,E} = H^{1}(C,\text{End}_{0}(E))$ and we can identify the
restricted map
\[
d\fD_{(L,\eta)|T_{\mathbb{P}^{1},\eta}} : T_{\mathbb{P}^{1},\eta} \to
H^{1}(C,\text{End}_{0}(E))
\]
with the Kodaira-Spencer class of the family $\mycal{E}\to C\times
\mathbb{P}^{1}$ at $\eta$.  But tensoring a family with a fixed line
bundle does not change the Kodaira-Spencer class and so equivalently
we can compute the Kodaira-Spencer class for the family $\mycal{F} =
\mycal{E}\otimes p_{C}^{*}L(-\pw ) \to C\times\mathbb{P}^{1}$. By
definition $\mycal{F}$ is the extension
\begin{equation} \label{eq-detq}
0 \to \mathcal{O}(-t')\boxtimes
\mathcal{O}(1) \to \mycal{F} \to
\mathcal{O}_{C\times \mathbb{P}^{1}} \to 0
\end{equation}
corresponding to the element
\[
\iota_{t'}^{-1} \in
H^{1}(C,\mathcal{O}(-t'))\otimes H^{1}( C,\mathcal{O}))^{\vee} =
\text{Ext}^{1}_{C\times \mathbb{P}^{1}}(\mathcal{O},
\mathcal{O}(-t')\boxtimes \mathcal{O}(1)
),
\]
and so is the universal family of extensions
\[
0 \to \mathcal{O}_{C}(-t') \to \mycal{F}_{\xi} \to \mathcal{O}_{C} \to 0,
\]
corresponding to points $\xi \in \mathbb{P}^{1}$. To simplify notation
we will denote $\mycal{F}_{\eta}$ by $F$.

Now from the exact sequence \eqref{eq-detq} it follows that
$\text{End}_{0}(\mycal{F}) =
\text{End}_{0}(\mycal{E})$ maps naturally onto
$\text{Hom}(\mathcal{O}(-t')\boxtimes
\mathcal{O}(1),\mycal{F})$ and so we have an exact sequence
\[
0 \to \mathcal{O}(-t')\boxtimes \mathcal{O}(1) \to \text{End}_{0}(\mycal{F}) \to
\mycal{F}\otimes (\mathcal{O}(t')\boxtimes \mathcal{O}(-1) ) \to 0
\]
which for each $\xi \in \mathbb{P}^{1}$ specializes to a short exact sequence
\begin{equation} \label{eq-end0}
0 \to \mathcal{O}_{C}(-t') \to \text{End}_{0}(\mycal{E}_{\xi})
   \to  \mycal{F}_{\xi}(t') \to 0.
\end{equation}
Since the family $\mycal{E}$ encodes the variation of the bundle
$\mycal{E}_{\xi}$ as the extension class $\xi \in
H^{1}(C,\mathcal{O}))$ varies, the Kodaira-Spencer class
$\mathsf{ks}^{\mycal{E}}_{\xi} : T_{\mathbb{P}^{1},\xi} \to
H^{1}(C,\text{End}_{0}(\mycal{E}_{\xi}))$ is captured by the maps
between cohomology of the terms in \eqref{eq-end0}.  Indeed, taking
into account that $H^{0}(C,\mathcal{O}(-t')) = 0$ and
$H^{0}(C,\text{End}_{0}(\mycal{E}_{\xi})) = 0$, the long exact
sequence of cohomology associated with \eqref{eq-end0} reads
\begin{equation} \label{eq-lesend0}
\lesone{0}{0}{H^{0}(C,\mycal{F}(t'))}
       {H^{1}(C,\mathcal{O}(-t'))}
       {H^{1}(C,\text{End}_{0}(\mycal{E}_{\xi}))}
       {H^{1}(C,\mycal{F}_{\xi}(t'))}
\end{equation}
Also,  in the long exact sequence associated with
\[
0 \to \mathcal{O}_{C} \to \mycal{F}_{\xi}(t') \to \mathcal{O}_{C}(t') \to 0 
\]
the first edge homomorphism $H^{0}(C,\mathcal{O}(t')) \to
H^{1}(C,\mathcal{O})$ is given simply by cup product with
$\iota_{t'}^{-1}(\xi) \in H^{1}(C,\mathcal{O}(-t'))$. Since $\xi$ is
not zero and $H^{0}(C,\mathcal{O}(t')) \cong \mathbb{C}$ this implies
that this edge homomorphism is injective and so the natural map
$H^{0}(C,\mathcal{O}) \to H^{0}(C,\mycal{F}_{\xi}(t'))$ is an
isomorphism.  If we identify $H^{0}(C,\mycal{F}_{\xi}(t'))$ with
$H^{0}(C,\mathcal{O})$ via this isomorphism, then the first edge
homomorphism in \eqref{eq-lesend0} becomes the cup product
\[
\iota_{q'}^{-1}(\xi)\cup \ : H^{0}(C,\mathcal{O}) \to H^{1}(C,\mathcal{O}(-t')).
\]
Using $\iota_{t'}$ to identify $H^{1}(C,\mathcal{O}(-t'))$ with
$H^{1}(C,\mathcal{O})$ we then see that the long exact
sequence \eqref{eq-lesend0}
reduces to a short exact sequence
\[
0 \to H^{1}(C,\mathcal{O})/\mathbb{C}\cdot\xi \to
{H^{1}(C,\text{End}_{0}(\mycal{E}_{\xi}))} \to 
       {H^{1}(C,\mycal{F}_{\xi}(t'))} \to 0.
\]
But $H^{1}(C,\mathcal{O})/\mathbb{C}\cdot\xi = T_{\mathbb{P}^{1},\xi}$
and the first map in this sequence by construction sends an
infinitesimal deformation of the extension $\xi$ to an infinitesimal
deformation of $\mycal{E}_{\xi}$.  Therefore this map is the
Kodaira-Spencer map and we can rewrite the previous sequence as
\[
\xymatrix@1{ 0 \ar[r] & T_{\mathbb{P}^{1},\xi}
  \ar[r]^-{\mathsf{ks}^{\mycal{E}}_{\xi}} &
     {H^{1}(C,\text{End}_{0}(\mycal{E}_{\xi}))} \ar[r] &
     {H^{1}(C,\mycal{F}_{\xi}(t'))} \ar[r] & 0.  }
\]
or in the special case $\xi =\eta$ as
\[
\xymatrix@1@C+2pc{ 0 \ar[r] & T_{\mathbb{P}^{1},\eta}
  \ar[r]^-{d\fD_{(L,\eta)|T_{\mathbb{P}^{1},\eta}}} &
     {H^{1}(C,\text{End}_{0}(E))} \ar[r] &
     {H^{1}(C,F(t'))} \ar[r] & 0.  }
\]
This gives the desired modular interpretation of
$d\fD_{(L,\eta)|T_{\mathbb{P}^{1},\eta}}$, and confirms that $\fD$ restricted to
any ruling $\{L\}\times \mathbb{P}^{1}$ is an immersion.

\

\bigskip

\subsubsection{Interpretation of
  \texorpdfstring{$d\fD_{(L,\eta)|T_{\Cbar ,L}}$}{dfDTC}}
\label{ssec-resCbar}

Fix $\eta \in \mathbb{P}^{1} = \mathbb{P}(H^{1}(C,\mathcal{O}))$. The
extensions of $L^{\vee}(\pw )$ by $L$ corresponding to the fixed
extension class $\eta$ and a variable $L \in \Cbar $ fit in a natural
family parametrized by $\Cbar $. Indeed, consider the surface $C\times
C$ with its diagonal divisor $\Delta \subset C\times C$. Pushing down
the short exact sequence
\[
0 \to \mathcal{O}(-\Delta) \to \mathcal{O} \to \mathcal{O}_{\Delta} \to 0 
\]
via the projection $p_{2} : C\times C \to C$ onto the second factor
gives a long exact sequence of direct images:
\[
\lesone{p_{2*}\mathcal{O}(-\Delta)}
       {p_{2*}\mathcal{O}}{p_{2*}\mathcal{O}_{\Delta}}
       {R^{1}p_{2*}\mathcal{O}(-\Delta)}
       {R^{1}p_{2*}\mathcal{O}}{R^{1}p_{2*}\mathcal{O}_{\Delta}}
\]
We have $p_{2*}\mathcal{O}(-\Delta) = 0$ since $\mathcal{O}(-\Delta)$
  has negative degree along the fibers of $p_{2}$. Since $p_{2} :
  \Delta \to C$ is an isomorphism we also have
  $p_{2*}\mathcal{O}_{\Delta} = \mathcal{O}_{C}$ and
  $R^{1}p_{2*}\mathcal{O}_{\Delta} = 0$.  Finally by K\"{u}nneth we
  have $p_{2*}\mathcal{O} = \mathcal{O}_{C}$ and
  $R^{1}p_{2*}\mathcal{O} = H^{1}(C,\mathcal{O})\otimes
  \mathcal{O}_{C}$. Substituting these in the long exact sequence
  gives
 \[
\lesone{0}
       {\mathcal{O}_{C}}{\mathcal{O}_{C}}
       {R^{1}p_{2*}\mathcal{O}(-\Delta)}
       {H^{1}(C,\mathcal{O})\otimes
  \mathcal{O}_{C}}{0}
\]
Therefore we have
\[
p_{2*}\mathcal{O}(-\Delta) = 0, \qquad\qquad
  R^{1}p_{2*}\mathcal{O}(-\Delta) = H^{1}(C,\mathcal{O})\otimes
  \mathcal{O}_{C}
\]
With this in mind we compute
\[
\begin{aligned}
  \text{Ext}^{1}_{C\times C}(\mathcal{O},\mathcal{O}(-\Delta)) 
& = H^{1}(C\times C, \mathcal{O}(-\Delta)) \\
  & = H^{0}(C,R^{1}p_{1*}\mathcal{O}(-\Delta)) \\
  & = H^{1}(C,\mathcal{O})\otimes
  H^{0}(C,\mathcal{O})
\end{aligned}
\]
where at the second step we used the Leray spectral sequence for
$p_{2}$.

Therefore the class
$\eta\otimes 1 \in H^{1}(C,\mathcal{O})\otimes H^{0}(C,\mathcal{O})$
corresponds to a rank two bundle $\mathbb{F}$ on $C\times C$ which is
an extension
\[
0 \to \mathcal{O}_{C\times C}(-\Delta) \to \mathbb{F} \to
\mathcal{O}_{C\times C} \to 0.
\]
By construction, for every point $x \in C$, the restriction
$\mathbb{F}_{x} = \mathbb{F}_{|C\times \{x\}}$ is the extension 
\[
0 \to \mathcal{O}_{C}(-x) \to \mathbb{F}_{x}\to \mathcal{O}_{C} \to 0
\]
given by the class $\iota_{x}^{-1}(\eta) \in
H^{1}(C,\mathcal{O}(-x))$.  In particular we have $\mathbb{F}_{t'} =
F$.

Now consider the degree $16$ map $C\times \Cbar  \to C\times C$
which is the identity on the first factor and on the second factor is
the composition of the natural map $\Cbar  \to C$ with the
hyperelliptic involution\footnote{The canonical map from
  $\Cbar $ to $C$ assigns to each $L$ the unique point $t \in C$
such that $L^{\otimes 2} = \mathcal{O}_{C}(t)$. The map we consider
here sends $L$ not to $t$ but to its image under the hyperelliptic
involution $t'$}. Let
$\mycal{F} \to C\times \Cbar $ be the pull back of
$\mathbb{F}$ by this map, and let $\mycal{L} \to C\times \Cbar $
be the Poincar\'{e} line bundle. Then the bundle
$\mycal{E} = \mycal{F}\otimes \mycal{L}^{\vee}(\{\pw\}\times
\Cbar )$ is an extension
\[
0 \to \mycal{L} \to \mycal{E} \to \mycal{L}^{\vee}(\{\pw\}\times
\Cbar ) \to 0.
\]
For each $M \in \Cbar $ restricts to the bundle 
$\mycal{E}_{M} := \mycal{E}_{|C\times \{M\}}$ which is the
extension
\[
0 \to M \to \mycal{E}_{M} \to M^{\vee}(\pw ) \to 0
\]
corresponding to the class $\eta$. In particular $\mycal{E}_{L} = E$
and the map $\fD_{|\Cbar \times \{\eta\}} : \Cbar  \to X$ is
the classifying map for the family $\mycal{E} \to
C\times\Cbar $. 

\

\noindent
This identifies the restricted differential $d\fD_{(L,\eta)|T_{\Cbar ,L}}$
with the Kodaira-Spencer class  
\[
\mathsf{ks}^{\mycal{E}}_{L} :
T_{\Cbar ,L} \to H^{1}(C,\text{End}_{0}(E))
\]
of the family $\mycal{E}$ computed $L \in \Cbar $. 

The family $\mycal{E}$ differs from the family of bundles
$\mycal{F}$ by a tensoring with a family of line bundles. So, even
though $\mycal{F}$ is a family of bundles with varying determinant,
the traceless parts of the Kodaira-Spencer classes for $\mycal{E}$
and $\mycal{F}$ are equal. Given a family of rank two bundles
$\mathcal{F}$ on $C$ let us write
$\kappa^{\mathcal{F}} = \mathsf{ks}^{\mathcal{F}} -
(\text{tr}(\mathsf{ks}^{\mathcal{F}})/2)\cdot \text{id}$ for the
traceless part of the Kodaira-Spencer class of $\mathcal{F}$. With
this notation we then have
\[
d\fD_{(L,\eta)|T_{\Cbar ,L}} = \mathsf{ks}^{\mycal{E}}_{L} =
\kappa^{\mycal{E}}_{L} = \kappa^{\mycal{F}}_{L}.
\]
As a last simplification, note that the map $\Cbar  \to C$,
$L \mapsto t'$ is \'{e}tale and so its differential is an isomorphism
$T_{\Cbar ,L} \stackrel{\cong}{\to} T_{C,t'}$ of tangent spaces. Taking
into account that the family $\mycal{F}$ is the pullback of the family
$\mathbb{F}$ we conclude that modulo the isomorphism
$T_{\Cbar ,L} \stackrel{\cong}{\to} T_{C,t'}$ the traceless
Kodaira-Spencer class $\kappa^{\mycal{F}}_{L}$ equals the traceless
Kodaira-Spencer class $\kappa^{\mathbb{F}}_{t'}$.

\

\medskip

\subsubsection{The kernel of the differential}
After these preliminaries we are now ready to understand where the map
$\fD : D \to X$ is not immersive.

\begin{proposition} \label{prop-kernel-df} Let $(L,\eta) \in D =
  \Cbar \times \mathbb{P}^{1}$ be a point, and let $E =
  \fD(L,\eta) \in X$ be the corresponding rank two bundle.
  Then the differential
  \[
  d\fD_{(L,\eta)} : T_{D, (L,\eta)} \to T_{X,E}
  \]
has a non-trivial kernel if and only if $L \in \Cbar $ maps to
$\eta \in \mathbb{P}^{1} = \mathbb{P}(H^{1}(C,\mathcal{O}_{C}))$ under
the natural degree $32$ map  $\Cbar  \to C \to \mathbb{P}^{1}$. 
\end{proposition}

\noindent {\bfseries Proof.} \ Since $d\fD_{(L,\eta)}$ is
determined by its restricions to $T_{\Cbar ,L}$ and
$T_{\mathbb{P}^{1},\eta}$ and since we saw that
$d\fD_{(L,\eta)|T_{\mathbb{P}^{1},\eta}}$ is always injective, it
suffices to characterize all points $(L,\eta) \in D$ for which the
image $d\fD_{(L,\eta)}(T_{\Cbar ,L})$ is contained in the line
$d\fD_{(L,\eta)}(T_{\mathbb{P}^{1},\eta}) \subset T_{X,E} =
H^{1}(C,\text{End}_{0}(E)) = H^{1}(C,\text{End}_{0}(F))$.

In view of the modular interpretations of the restricted differentials
in sections \ref{ssec-resP1} and \ref{ssec-resCbar} this question is
equivalent to the problem of characterizing all points $(t',\eta) \in
C\times \mathbb{P}^{1}$ such that 
\begin{equation} \label{eq-condition-cusp}
\kappa^{\mathbb{F}}_{t'}(T_{C,t'}) \subset
\ker\left[H^{1}(C,\text{End}_{0}(F)) \to H^{1}(C,F(t'))\right].
\end{equation}
Here
 as before $\mathbb{F}$ is the extension
\begin{equation} \label{eq-extension_on_CxC}
\xymatrix@1{
  0 \ar[r] &  \mathcal{O}_{C\times C}(-\Delta) \ar[r]^-{a} &
  \mathbb{F} \ar[r]^-{b} &   \mathcal{O}_{C\times C} \to
  0,
  }
\end{equation}
corresponding to $\eta$ and $F = \mathbb{F}_{t'}
:= \mathbb{F}_{|C\times \{t'\}}$.

More invariantly the condition \eqref{eq-condition-cusp} can be
rewritten as follows. The short exact sequence
\eqref{eq-extension_on_CxC} induces a short exact sequence
\begin{equation} \label{eq-filt.end}
\xymatrix@R-1.5pc@C+0.7pc{
  0 \ar[r] &  \text{Hom}(\mathcal{O},\mathcal{O}(-\Delta))
  \ar@{=}[d] \ar[r]^-{a\circ-\circ b}  & \text{End}_{0}(\mathbb{F})
   \ar[r]^-{-\circ a}  &
  \text{Hom}(\mathcal{O}(-\Delta),\mathbb{F}) \ar@{=}[d] \ar[r] & 0   \\
&  \mathcal{O}(-\Delta) & & \mathbb{F}(\Delta) &}
\end{equation}
Pushing this down the map $-\circ a : \text{End}_{0}(\mathbb{F}) \to
\mathbb{F}(\Delta)$ via the second projection $p_{2} : C\times C \to
C$ gives us a map of coherent sheaves
\[
R^{1}p_{2*}(-\circ a) : R^{1}p_{2*}\text{End}_{0}(\mathbb{F}) \to
R^{1}p_{2*}\mathbb{F}(\Delta) .
\]
We also have the universal traceless Kodaira-Spencer class of $\mathbb{F}$
\[
\kappa^{\mathbb{F}} : T_{C} \to R^{1}p_{2*}\text{End}_{0}(\mathbb{F}).
\]
Composing these two maps we get a map
\[
R^{1}p_{2*}(-\circ a) \circ
\kappa^{\mathbb{F}} : T_{C} \to R^{1}p_{2*}\mathbb{F}(\Delta),
\]
or equivalently a global section
\begin{equation} \label{eq-global.section}
R^{1}p_{2*}(-\circ a) \circ \kappa^{\mathbb{F}} \in
H^{0}\left(C,\Omega^{1}_{C}\otimes R^{1}p_{2*}\mathbb{F}(\Delta)\right).
\end{equation}
Note that the bundle $\mathbb{F}$ and the global section
\eqref{eq-global.section} depend only on the class $\eta \in
\mathbb{P}(H^{1}(C,\mathcal{O}))$. The condition
\eqref{eq-condition-cusp} is simply the condition that the point $t'$
belongs to the zero locus of this section. Therefore we need to
describe all pairs $(t',\eta) \in C\times \mathbb{P}^{1}$ for which
the section \eqref{eq-global.section} defined by $\eta$
vanishes at the point $t'$.

To understand these pairs observe that $\mathcal{O}(-\Delta) \subset
\mathbb{F}$ can be viewed as the family of filtered bundles
$\mathcal{O}(-x) \subset \mathbb{F}_{x}$ parametrized by $x \in C$. In
particular, the traceless Kodaira-Spencer class $\kappa^{\mathbb{F}}$
will preserve the filtration. In other words if ${}^{\rm
  filt}\text{End}_{0}(\mathbb{F})$ denotes the bundle of traceless
endomorphisms of $\mathbb{F}$ preserving the subbundle
$\mathcal{O}(-\Delta)$, then the global traceless Kodaira-Spencer
class $\kappa^{\mathbb{F}}$ is induced from the traceless part of the
Kodaira-Spencer class
  \[
{}^{\rm filt}\kappa^{\mathbb{F}} : T_{C} \to R^{1}p_{2*} {}^{\rm
  filt}\text{End}_{0}(\mathbb{F})
  \]
 classifying filtered deformations. This means that
 $\kappa^{\mathbb{F}}$ factors as
\[
\xymatrix@1{ T_{C} \ar[r]^-{{}^{\rm filt}\kappa^{\mathbb{F}}}
  \ar@/_1.5pc/[rr]_-{\kappa^{\mathbb{F}}} & R^{1}p_{2*} {}^{\rm
    filt}\text{End}_{0}(\mathbb{F}) \ar[r] &
  R^{1}p_{2*}\text{End}_{0}(\mathbb{F}).  }
\]
By definition
\[
  {}^{\rm filt}\text{End}_{0}(\mathbb{F}) = \left\{ \varphi \in \text{End}(F) \, \left| \,
b\circ \varphi \circ a = 0 \ \text{and} \ \text{tr}(\varphi)  = 0
  \right. \right\},
  \]
  and hence
\[
  {}^{\rm filt}\text{End}_{0}(\mathbb{F}) =
  \ker\left[\xymatrix@1{\text{End}_{0}(\mathbb{F})\ar[r]^-{b\circ-\circ a} &
      \mathcal{O}(\Delta)}
      \right].
\]
This implies that $R^{1}p_{2*}(-\circ a)\circ\kappa^{\mathbb{F}}$
maps $T_{C}$ into the subsheaf
\[
\ker\left[ R^{1}p_{2*}\mathbb{F}(\Delta) \to
R^{1}p_{2*}\mathcal{O}(\Delta)\right] \subset
  R^{1}p_{2*}\text{End}_{0}(\mathbb{F}.
\]
But from the long exact sequence of $p_{2}$ direct images of the short exact sequence
\[
0 \to \mathcal{O} \to \mathbb{F}(\Delta) \to \mathcal{O}(\Delta) \to 0
\]
it follows immediately that
\[
 \ker\left[\xymatrix@1{\text{End}_{0}(\mathbb{F})\ar[r]^-{b\circ-\circ a} &
      \mathcal{O}(\Delta)}
      \right] = (H^{1}(C,\mathcal{O})/\eta)\otimes \mathcal{O}_{C}. 
\]
This imlplies that the section $R^{1}p_{2*}(-\circ a)\circ\kappa^{\mathbb{F}}$ in
$\omega_{C} \otimes R^{1}p_{2*}\mathbb{F}(\Delta)$ is in fact the tautological
section in the natural subbundle
\[
\omega_{C}\otimes \left( (H^{1}(C,\mathcal{O})/\eta)\otimes
\mathcal{O}_{C}) \right) = (H^{1}(C,\mathcal{O})/\eta)\otimes \omega_{C},
\]
i.e. is the element in
\[
(H^{1}(C,\mathcal{O})/\eta)\otimes H^{0}(C,\omega_{C})
= (H^{1}(C,\mathcal{O})/\eta)\otimes H^{1}(C,\mathcal{O})^{\vee}
\]
corresponding to the graph of the quotient map $H^{1}(C,\mathcal{O})
\to H^{1}(C,\mathcal{O})/\eta$. This is precisely the section of
$\omega_{C}$ that vanishes at the preimage of $\eta$ under the
hyperelliptic map $\hyp_{C} : C \to \mathbb{P}^{1}$. Therefore the set of pairs
$(t',\eta)$ we were seeking is just the graph of the hyperelliptic map
and so the set of points $(L,\eta)$ at which $\fD$ is not immersive is
the graph of the projection $\Cbar  \to C \to
\mathbb{P}^{1}$. This completes the proof of the proposition. \ \hfill
$\Box$

\subsection{The trigonal cover of a general line}
\label{trigonal}

Recall from \cite[Theorem~2]{Newstead}, \cite[Theorem~5]{NR} that the
variety of lines on $X$ is isomorphic to the Jacobian
$\op{Jac}^{0}(C)$ of $C$.  Taking up the theme used in the proof of
Lemma \ref{wobbly8h} in more detail, consider a general line $\ell
\subset X$. It corresponds to a degree $0$ line bundle $A\in
\op{Jac}^0(C)$, and in terms of vector bundles the line is the set of
isomorphism classes of nonzero extensions of the form
$$
0\rightarrow A \rightarrow V \rightarrow A^{\vee}(\pw ) \rightarrow 0.
$$
A point corresponding to a given extension, will lie on a 
different line $\ell'$ corresponding to $A'\in  \op{Jac}^0(C)$
distinct from $A$,
if and only if there is a map $A'\rightarrow A^{\vee}(\pw)$ such that the
pullback extension splits. Such a map corresponds to a point $t\in C$
with $A'=A^{\vee}(\pw -t)$. 

The splitting of the extension may be analyzed as follows. We have
$$
{\rm Ext}^1(A^{\vee}(\pw ),A)\cong H^1(C,A^{\otimes 2} (-\pw ))
$$ which is dual to $H^0(C,A^{\otimes -2}(3\pw ))$. These spaces have
dimension $2$, so a nonzero element of $H^1(C,A^{\otimes 2} (-\pw ))$
modulo scalars also corresponds to a nonzero element of
$H^0(C,A^{\otimes -2}(3\pw ))$ modulo scalars, the unique element that
pairs with it to zero by the duality.

The extension splits when pulled back to $A^{\vee}(\pw -t)$
if and only if the class maps to zero under the map 
$$
H^1(C,A^{\otimes 2} (-\pw ))\rightarrow H^1(C,A^{\otimes 2} (t-\pw )),
$$
or equivalently if the dual element is in the image of the map 
$$
H^0(C,A^{\otimes -2}(3\pw -t))\rightarrow H^0(C,A^{\otimes -2}(3\pw )).
$$
In turn, this is equivalent to saying that our dual element, viewed as a 
section of the line bundle $A^{\otimes -2}(3\pw )$, vanishes at $t$. 

In conclusion, each extension corresponds to a class (up to scalar multiples)
of nonzero sections of $A^{\otimes -2}(3\pw )$, and the other lines 
$\ell'$ passing through the point determined by the extension correspond
to the points where this section vanishes. 

We may view this as saying that the complete linear system associated to
the line bundle $A^{\otimes -2}(3\pw )$ of degree $3$ determines a trigonal map 
$C\rightarrow \ell \cong \pp^1$, and the other three lines through a point on
$\ell$ correspond to the fibers. 

\begin{proposition}
\label{ellW}
The trigonal cover ramifies over $8$ points in $\ell$. These are the
$8$ intersection points of $\ell$ with the wobbly locus. If $x$ is
such an intersection point, and if we write the divisor in $C$ as $2u+v$
then $u$ corresponds to the line through that point in the wobbly locus. 
\end{proposition}
\begin{proof}
If $A$ is general then the line bundle $A^{\otimes -2}(3\pw )$ is
general and hence has no base points. For such an $A$ the trigonal map
will be a morphism. The Hurwitz formula implies that there are $8$
ramification points. A general fiber of the trigonal cover is a
section of $A^{\otimes -2}(3\pw )$ that vanishes at three points
$u,v,w\in C$.  A simple ramification is a point $u$ such that there is a
section of $A^{\otimes -2}(3\pw )$ that vanishes at $2u+v$. We get
$A^{\otimes -2}(3\pw ) \cong \Oo _C(2u+v)$. As discussed above, the
extension corresponding to this point of $\ell$ splits on
$$
A^{-1}(\pw -u) \cong A(u+v-2\pw ) \cong A(u-v')
$$
where $v'$ is the
hyperelliptic conjugate of $v$, and we get a diagram
$$
\xymatrix@M+0.5pc{
0 \ar[r] & A \ar[r] & E \ar[r] & A^{-1}(\pw) \ar[r] & 0 \\
& & & A(u - v') \ar[u] \ar[ul] &
}  
$$
and the bundle fits in an extension  
$$
0 \rightarrow A (u-v') \rightarrow E \rightarrow A^{-1}(3\pw -u-v)\cong A(u) \rightarrow 0.
$$
Now $A(u-v') \otimes \omega _C = A(2\pw +u-v') = A(u+v)$ so we get a nonzero map 
$$
A(u) \rightarrow A(u+v) \cong A(u-v') \otimes \omega _C
$$
yielding a nonzero nilpotent Higgs field. Therefore, the bundle $E$
associated to the given extension is in the wobbly locus.
\end{proof}

\

\noindent
\emph{\bfseries Remark:} \ We already argued synthetically in
section~\ref{sssec:lpatterns} that the line $\ell$ goes through the
cuspidal locus of the wobbly locus when the divisor in $C$ has the
form $3u$.

\

Choose a general line $\ell \subset X=X_1$. Over this line we have
the trigonal cover $C\rightarrow \ell $ associating to each point 
$x\in \ell$ the set of three lines passing through $x$ that are distinct from 
$\ell$. Let $u,v,w$ be the points of $C$ over $x\in \ell$. 

Let $Y_{\ell}:= Y\times _X \ell$. It is a degree $8$ ramified cover of 
$\ell$. 

\begin{lemma}
\label{three-eight}
Assume $x$ general in $\ell$. 
The fiber of $Y_{\ell}\rightarrow \ell$ over the point $x$ is
naturally isomorphic to the set of $8$ points obtained by choosing 
$\tilde{u},\tilde{v},\tilde{w} \in \Ctilde $ over $u,v,w\in C$ respectively. 
\end{lemma}
\begin{proof}
Let $L$ be the line bundle of degree $0$ on $C$ corresponding to the line
$\ell$, and let  $E$ be the rank $2$ vector bundle corresponding to $x$.
It fits into 
an exact sequence 
$$
0\rightarrow L \rightarrow E \rightarrow L^{-1}(\pw )\rightarrow 0
$$
whose extension class corresponds to the point in $\ell$.  The
fiber of $Y$ over $x$ is the set of line bundles $U$ on $\Ctilde$, of
degree $3$, such that $\pi _{*}(U)\cong E$. The map $L\rightarrow \pi
_{*}(U)$ corresponds by adjunction to a map $\pi ^{*}(L) \rightarrow
U$. We get
$$
U \cong 
\pi ^{*}(L) \otimes \Oo _{\Ctilde}(\tilde{u} + \tilde{v} + \tilde{w}). 
$$
Let $u,v,w$ be the images of 
$\tilde{u}$, $\tilde{v}$, $\tilde{w}$ in $C$. The determinant of $\pi _{*}(U)$
is $L^{\otimes 2} \otimes \Oo _C(u+v+w)\otimes \omega _C^{-1}$ so we get
$$
L^{\otimes 2}(u+v+w) \cong \Oo _C(3\pw ).
$$
It means that $u+v+w$ is in the linear system $| L^{\otimes
  -2}(3\pw )|$, in other words it is one of the fibers of the trigonal
cover of $\ell$. To show that it is the fiber over the point
corresponding to the given extension defining $E$, we need to see that
the extension splits when restricted to $L^{-1}(\pw -u)$ (and the same
for $v,w$). For this, in turn, it suffices to see that there is a
nonzero map
$$
\pi ^{*}(L^{-1}(\pw -u)) \rightarrow U.
$$
Using the expression 
$$
L^{-1}(\pw -u) = L \otimes  (L^{\otimes -2}(3\pw ) ) \otimes \Oo _C(-2\pw  - u) 
= L \otimes \Oo _C (-2\pw  - u + u + v + w) = L (v + w - 2\pw ),
$$
we are looking for a map
\begin{equation} \label{eq:taumap}
\pi ^{*}(L )(\tilde{v} + \tau \tilde{v} + \tilde{w} + \tau \tilde{w})
\rightarrow \pi ^{*}(L \otimes \omega _C ) \otimes \Oo
_{\Ctilde}(\tilde{u} + \tilde{v} + \tilde{w}).
\end{equation}
Here, as usual, $\tau$ denotes the covering involution of the
spectral cover $\pi: \Ctilde \to C$.

Having a map \eqref{eq:taumap}  is equivalent to asking for a section of
$$
\pi ^{*}(\omega _C ) \otimes \Oo _{\Ctilde}(\tilde{u} - \tau
\tilde{v} - \tau \tilde{w}).
$$
The ramified cover $\pi : \Ctilde \rightarrow C$ has $4$
ramification points $a,a',b,b'\in \Ctilde$ such that $\pi (a),\pi
(a')$ and $\pi (b),\pi (b')$ are opposite pairs under the
hyperelliptic involution of $C$.

Recall now that $\Ctilde$ itself is hyperelliptic of genus $5$, with
the map $\he_{\Ctilde} : \Ctilde \rightarrow \pp^1$ branched over the
$12$ inverse images in $\Ctilde$ of the $6$ Weierstrass points of
$C$. We have $\omega _{\Ctilde} = \he_{\Ctilde}^{*} \Oo _{\pp ^1}(4)$.  The
hyperelliptic involution $\sigma$ of $\Ctilde$ is one of the lifts of
the hyperelliptic involution of $C$, so $a,a'$ and $b,b'$ are also
opposite pairs under the hyperelliptic involution of $\Ctilde$. It
follows that $\Oo _{\Ctilde}(a+a') \cong \Oo _{\Ctilde}(b+b') \cong
\he_{\Ctilde}^{*}\Oo _{\pp ^1}(1)$. We get
$$
\pi ^{*}(\omega _C ) \cong \omega _{\Ctilde} (-a-a'-b-b') =
\he_{\Ctilde} ^{*} \Oo _{\pp^1}(2).
$$
In particular, $\pi ^{*}(\omega _C )$ admits a section vanishing at
the two points
$\tau \tilde{v}, \tau \tilde{w}$. This gives the required section of 
$\pi ^{*}(\omega _C ) (\tilde{u} - \tau \tilde{v} - \tau \tilde{w})$. 

This (together with the analogous arguments for $v$ and $w$) completes
the proof that the divisor $u+v+w$ is the fiber of the trigonal map
over $x\in \ell$. Thus, our line bundle $U$ corresponds to one of the
choices of $8$ liftings as in the lemma.

Running this argument backwards shows that each of the $8$ liftings
corresponds to a choice of $U$. This gives the isomorphism claimed in
the lemma.
\end{proof}

\

\noindent
Let $Y_{\ell}:= Y\times _X \ell$. It is a degree $8$ ramified cover of 
$\ell$ whose general fiber is described by the lemma. From this, we can
describe the ramification of $Y_{\ell} / \ell$:

\begin{corollary}
\label{describe-branch}
Suppose $\ell$ is general.  The branch locus of $Y_{\ell} / \ell$
consists two disjoint groups of points: the $x\in \ell $ of the branch
locus of $C/\ell$, plus the points that are images of the four points
$y\in C$ of the branch locus of $\Ctilde /C$.  The former class of
points is the intersection of $\ell$ with the wobbly locus. As $\ell$
is general, the second group consists of four distinct points.
\end{corollary}
\begin{proof}
The line $\ell$ gives a trigonal map $C\rightarrow \pp^1$, and by
Lemma \ref{three-eight}, the fiber of $Y_{\ell}$ over a point $x\in \ell$
consists of the $8$ liftings to $\Ctilde$ of the three points $u,v,w$ in
the trigonal fiber over $x$. 

When $x$ is a branch point of $C/\ell $, as we have seen in
Proposition \ref{ellW} it means 
$x\in \ell\, \cap \Wob$, and two of the three points
come together, let's say $u=v$. The monodromy around such a point 
interchanges $u$ and $v$. We may assume that $u=v$ is general in $C$,
in particular it isn't a ramification point of $\Ctilde / C$ so we may
identify the two sheets of $\Ctilde$ labeled with $0,1$ 
near this point, similarly near $w$. 
With these identifications we may write the $8$ points as $(a,b,c)$ where
$a,b,c \in \{ 0, 1\}$. The monodromy sends $(a,b,c)$ to $(b,a,c)$. 
This has two transpositions $(1,0,c)\leftrightarrow (0,1,c)$ for $c=0,1$. 

Suppose $x$ is the image of one of the branch points $z\in C$ of
$\Ctilde / C$. For $\ell$ general we may assume that $x$ isn't
wobbly. The trigonal fiber is $u,v,w$ with, say, $u$ the branch
point. As we move around such a point, the set of choices of lifting
of $u$ undergoes a transposition; there are four such transpositions
corresponding to the various liftings of $v,w$.

The last statement, that was also mentioned in the previous paragraph,
is that two different branch points of $\Ctilde / C$ are not contained
in the same trigonal fiber.  Indeed, the branch points consist of two
general pairs of opposite points under the hyperelliptic
involution. If the trigonal map identified two such points, they would
have to be two opposite points, but that would mean that the trigonal
map identifies hyperelliptically opposite points, hence that it
factors through the hyperelliptic map, which isn't the case.
\end{proof}

\begin{lemma}
\label{sum-origin}
Suppose $V$ is a vector bundle corresponding to a general point in
$X$. Let $A_1,\ldots , A_4$ be the points in $ \op{Jac}^0(C)$
corresponding to the four lines in $X$ through $V$. Then the sum of
these points in $\op{Jac}^{0}(C)$ is the origin of $ \op{Jac}^0(C)$.
\end{lemma}
\begin{proof}
Let $A=A_1$ be one of the line bundles, so we have 
$$
0\rightarrow A \rightarrow V \rightarrow A^{\vee}(\pw ) \rightarrow 0.
$$
As we have seen above, the other $A_i$ are of the form 
$A_i = A^{\vee}(\pw -t_i)$ where $t_2+t_3+t_4$ is a divisor
in the linear system $A^{\otimes -2}(3\pw )$. This yields the equation in 
$ \op{Jac}(C)$
$$
t_2+t_3+t_4 = [A^{\otimes -2}(3\pw )] = 3\pw  - 2[A].
$$
Now, the sum of all the points is the sum in $ \op{Jac}(C)$
$$
{\rm sum} = 
[A ]+ [A^{\vee}(\pw -q_2)] + [A^{\vee}(\pw -t_3)] + [A^{\vee}(\pw -t_4)]
$$
$$
= [A] + 3 [A^{\vee}] + 3\pw  - (t_2 + t_3 + t_4)
$$
$$
= [A] - 3 [A] + 3\pw  - (3\pw  - 2[A]) = 0. 
$$
\end{proof}

\

\noindent
Given two of the four lines corresponding to line bundles $A_1$ and
$A_2$, we obtain a point $q\in C$ such that $A_2 = (A_1)^{\vee}(\pw
-q)$, in other words $A_1 \otimes A_2 = \Oo _C(\pw -q)$.  In view of
the lemma, if $A_3$ and $A_4$ are the other two points, we have
$A_1\otimes A_2 \otimes A_3 \otimes A_4 = \Oo _C$ so $A_3\otimes A_4 =
\Oo _C(q-\pw ) = \Oo _C(\pw -q')$. Thus, the point of $C$
corresponding to $A_3,A_4$ is $q'$ the image of $q$ by the
hyperelliptic involution.

\begin{remark}
\label{maptoP3}
In particular, we obtain the conclusion that the point $V\in X_1$
yields a well-defined triple of points in $\pp^1$. This gives a map
$X_1 \rightarrow Sym ^3 (\pp ^1) \cong \pp^3$ that will be useful in
Section \ref{chapter-tensor}.
\end{remark}

\subsection{Local structure of the spectral cover along the cusp locus}

Let $R^{\rm mov} \subset Y$ be the movable locus of the ramification of the
map $Y\rightarrow X$. 

Recall that $c_1(\Omega ^1_X) = 2H$ and $c_1(\Omega ^1_Y) = \ExY $
since $Y$ is obtained from the abelian variety $\Prym$ by blowing up
a curve of codimension $2$. Comparing these, we conclude that
the ramification locus has class $\ExY +2\FxY $. 

\begin{lemma}
The order of ramification of $Y\rightarrow X$ along $\ExY $ is two
(i.e. it has simple ramification at a general point of $\ExY $). 
\end{lemma}
\begin{proof}
If we assume that the map has ramification of order $a$ along $\ExY $,
we get 
$$
2\ExY +\FxY  = (a-1)\ExY  + [R^{\rm mov}].
$$
On the other hand, $R^{\rm mov}$ is effective, so $2\FxY -(a-2)\ExY
$ is effective.  If $a>2$ this says that $2\FxY -\ExY $ is
effective. We note that $2\FxY -\ExY $ is not zero, since $\ExY $ is
clearly not ample whereas $\FxY $ is ample (it is the pullback of an
ample divisor by a finite map).

Since $\FxY $ is ample, if $2\FxY -\ExY $ is effective then $(2\FxY
-\ExY )\cdot \FxY ^2 > 0$. However, from Proposition
\ref{intersections1} saying $\FxY ^3 = 32$ and $\ExY \FxY ^2 = 64$ it
gives $(2\FxY -\ExY )\cdot\FxY ^2=0$, contradicting the effectivity.
\end{proof}

\begin{proposition}
For generic choice of point in the Hitchin base, the movable ramification
locus $R^{\rm mov}$ intersects $\ExY $ in a union of two constant sections
of the fibration $\ExY \rightarrow \Chat $  and 
$64$ fibers. The constant sections correspond to the two points in 
$\pp^1$ whose inverse images in $C$ make up the $4$ branch points
of $\Ctilde / C$.
\end{proposition}
\begin{proof}
From above, the full ramification locus has class $\ExY +2\FxY $ and the
fixed part has multiplicity $1$ along $\ExY $. Thus, the class of 
$R^{\rm mov}$ is $2\FxY $. 

We use the description of Corollary \ref{describe-branch}.
We are interested here in points where the movable ramification
meets the intersection with the wobbly locus. 

We can describe more precisely, over points where $\ell$ meets the
wobbly locus, the fiber of $Y_{\ell}$. Suppose it is a general point
of the wobbly locus, so the degree $3$ divisor in $C$ over $x\in \ell$ 
has the form 
$2u+v$. Suppose $u$ and $v$ are not in the branch locus of $\Ctilde / C$.
To describe the liftings to points of  $Y_{\ell}$ we need to consider liftings
of nearby divisors into $\Ctilde $: the piece $2u$ splits into $u_1+u_2$
and this has $4$ liftings. The monodromy action permutes two of those
liftings and leaves fixed the two other ones, depending on whether
$\tilde{u}_1$ and $\tilde{u}_2$ are in different or the same sheets 
of $\Ctilde / C$. Then each of those configurations becomes doubled
depending on the lifting of $v$. We obtain $2$ branch points and $4$
unbranched points of $Y_{\ell} /\ell$ over $x\in \ell$. This is as expected.

The points of $\ExY _{\ell}:= \ExY \cap Y_{\ell}$ are the two branch points.

On the other hand, the points of the movable ramification  locus 
correspond to divisors $u_1+u_2+v$ such that one of those points 
is a branch point of $\Ctilde / C$. When that branch point is $v$
it does not correspond to a point of $\ExY _{\ell}$. We would like to 
describe the cases where $u$ is a branch point of $\Ctilde / C$.

Suppose given a divisor $2u+v$ in the linear system associated to
$A^{\otimes -2}(3\pw )$. Then as we saw above, the other line bundle
of degree $0$ associated to $u$, that we'll now denote $L:= A'$, is of the form 
$$
L= A^{\vee}(\pw -u).
$$
Our bundle $V$ is an up-Hecke of $A\oplus L$ at the point $u$. 
We can write
$$
0\rightarrow L \rightarrow V \rightarrow L^{\vee}(\pw ) \rightarrow 0
$$
and $A=L^{\vee}(\pw -u)$. The extension is the one that vanishes when
pulled back to $A$, or equivalently vanishing in the extension group
$H^1(L^{\otimes 2}(u-\pw ))$. 

We have
$$
L^{\otimes2}(\pw ) = A^{\otimes -2}(3\pw -2u)
$$
but 
$$
A^{\otimes -2}(3\pw )\cong \Oo _C(2u+v)
$$
so
$$
L^{\otimes2}(\pw ) = \Oo _C(v).
$$
This is the equation saying that the line associated to $L$ is in the wobbly
locus. Furthermore, the set of extension classes corresponding
to points of this line is 
$$
H^1(L^{\otimes 2}(-\pw ))= H^1(\Oo _C(v-2\pw )) = H^1(\Oo _C(-v'))
$$
which maps isomorphically to $H^1(\Oo _C)$. Here $v'$ is the image of $v$
by the hyperelliptic involution. 

An extension vanishes at the point $u$ if and only if it 
maps to zero in $H^1(\Oo _C(u))$, which is equivalent to saying that
it corresponds to the
point in $\pp^1= \pp H^1(\Oo _C)$ that is the image of $u$ under the
hyperelliptic double cover. 

We conclude that under the birational map from the wobbly locus to 
$\pp^1$, the image of a point corresponding to the degree $3$ divisor
$2u+v$ with respect to the trigonal covering corresponding 
to the original line $\ell$, it the image of $u$ in $\pp^1$ corresponding to
the hyperelliptic map. 

Now, we are interested in cases where $u$ coincides with  one of the branch
points of $\Ctilde /C$. We know from Proposition \ref{prop:hitchin.map} that there
are $4$ branch points in $C$ that are the inverse images of two points in 
$\pp^1$. Therefore, the points in question are ones where $u$ maps
to one of these two points in $\pp^1$. 

We conclude from this discussion that the points of $\ExY $ over the
wobbly locus that map to one of these two points in $\pp^1$ are in
$R^{\rm mov}\cap \ExY $.

We have made some genericity assumptions in the above discussion. These
hold for the points we have just described, since the two points in 
$\pp^1$ are general because we chose a general spectral cover 
$\Ctilde /C$. This says that $R^{\rm mov}\cap \ExY $
contains at least two constant sections of the fibration $\ExY \rightarrow 
\Chat $. 

We now rely on the intersection number calculations to
conclude. Indeed, we have seen that $\FxY \cap \ExY $ has the class of
$\hat{\mathfrak{c}} + 32 \mathsf{fib}$ where $\hat{\mathfrak{c}}$ is a
constant section. This may be re-verified as follows: this class
intersects the fiber in one point, since the fiber maps to a line in
$X$ that meets the hyperplane class in one point. Our class is thus of
the form $\hat{\mathfrak{c}}+ a\cdot\mathsf{fib}$. On the other hand,
the self-intersection of this class on $\ExY \cong \pp^1\times \Chat $
is equal to the intersection number $\FxY ^2\ExY $ on $Y$; that is
$64$.  The self intersection of $\hat{\mathfrak{c}} +
a\cdot\mathsf{fib}$ is $2a$ giving $a=32$.

This shows that $2\FxY \cap \ExY $ is a divisor in the class
$2\hat{\mathfrak{c}} + 64\mathsf{fib}$.  Since we have already
exhibited a subset of this divisor consisting of two sections, we
conclude that it has to be of the form two sections plus $64$
fibers. This completes the proof of the proposition.
\end{proof}

For a point $\mathsf{b}$ in the Hitchin base, let ${\rm Branch}^{\rm
  mov}(\mathsf{b})\subset X$ be the image in $X$ of the ramification
locus $R^{\rm mov}\subset Y$ for the point $\mathsf{b}$.  These form a
family of closed subvarieties of $X$.

We define the {\em base locus} of this family to consist of the set of
points $x\in X$ that lie in ${\rm Branch}^{\rm mov}(\mathsf{b})$ for
$\mathsf{b}$ general with respect to $x$, or equivalently for all
$\mathsf{b}$.

\begin{corollary}
If the base locus contains an irreducible subvariety of codimension
$\leq 2$ then it is one of the lines in the wobbly locus.
\end{corollary}
\begin{proof}
If $x\in X$ is not on the wobbly locus then it is a very stable
bundle, and the space of Higgs fields on that bundle maps properly to
the Hitchin base. Thus, the intersection with a general Hitchin fiber
over $\mathsf{b}$ is transverse, so $x$ is not in the branch locus of
the map $Y\rightarrow X$ for general $\mathsf{b}$.

Suppose $x$ is in the wobbly locus and let $S$ be the space of Higgs
fields on the corresponding bundle. It maps to the Hitchin base, by a
map that is not proper.  However, for a general point $\mathsf{b}$ in
the Hitchin base, the intersection of the fiber over $\mathsf{b}$ with
$S$ is transverse. Therefore, for such a general $\mathsf{b}$, if
$y\in Y$ is a point over $x$ that is not in $\ExY $, then it is in
$S\cap \hit^{-1}(\mathsf{b})$ hence not in $R^{\rm mov}$.  This shows
that the points of $R^{\rm mov}$ lying over $x$ are in $R^{\rm
  mov}\cap \ExY $.

We have seen in the proposition that for $\mathsf{b}$ general, the
intersection of $R^{\rm mov}$ with $\ExY $ consists of two constant
sections of the map $\ExY = \Chat\times \pp^{1} \to \Chat$, that move
as a function of $\mathsf{b}$, and of $64$ fibers of the map $\ExY =
\Chat\times \pp^{1} \to \Chat$. The images of these in the wobbly
locus have the same description.

So, suppose $x$ is a point of the wobbly locus contained in the image
of $R^{\rm mov}$ for two different general values of
$\mathsf{b}$. Then it is in the intersection of two subsets of the
form $2$ sections plus images of $64$ fibers; but the sections move,
so $x$ has to be contained in at most one of the $64$ lines. This
completes the proof of the corollary.
\end{proof}

\

\begin{corollary}
A general point on the cuspidal locus of the wobbly locus is not contained
in the movable part of the branch locus.
\end{corollary}
\begin{proof}
This follows immediately from the previous corollary since the cuspidal
locus has codimension $2$ and is not a union of fibers.  
\end{proof}

\begin{corollary}
\label{reallymovable}
The ``movable part'' of the ramification locus
does really move, in the sense that the base locus does not contain any
divisors.
\end{corollary}
\begin{proof}
The above corollary implies that the codimension of the base locus is
$\leq 2$.
\end{proof}

\subsection{The step in the parabolic filtration}

Recall that $\Vv _{a,b} = f _{*}(\Ll _{a,b})$ and consider the sheaf
$$
\Vv _{a,b} / \Vv _{a,b-1}
$$
which is supported on $\Wob \subset X$ by definition. 
Over $X^{\circ}$ this is the quotient that will be used to define
the crude parabolic structure. For the computations, it will be
necessary to have a refined definition that holds over all of
$X$. Recall that $\fD : D\rightarrow X$ denotes the map from the
normalization $D$ of $\Wob $.

\

\begin{lemma}
There is a quotient defined over $D$,  
$$
\fD^{*}(\Vv _{a,b})\rightarrow \Uu_{a,b}
$$
such that $\Uu_{a,b}$ is a rank $2$ torsion-free sheaf on $D$,
hence a bundle
outside of a finite collection of points. 
Over $X^{\circ}$ (this restriction being denoted by 
a superscript as usual) we have \linebreak
$\left(\Vv_{a,b}/\Vv_{a,b-1}\right)|_{X^{\circ}} = \fD_{*}\Uu^{\circ}_{a,b}$
as quotients of $\Vv ^{\circ}_{a,b}$ and this characterizes
$\Uu_{a,b}$.  
\end{lemma}
\begin{proof}
Start by noting that the property of the lemma uniquely characterizes
$\Uu_{a,b}$.  Indeed, the condition at the end states that the kernel
of the map $\fD^{*}(\Vv _{a,b})\rightarrow \Uu_{a,b}$ is, over
$X^{\circ}$ (where $D^{\circ}=\Wob^{\circ}\subset X$) the subsheaf of
sections generated by $\Vv _{a,b - 1}$. The condition that the
quotient $\Uu_{a,b}$ is a torsion-free sheaf means that the kernel is
a saturated subsheaf, so it is defined by its restriction to
$D^{\circ}$.

This discussion tells us how to construct $\Uu_{a,b}$, namely set
$$
\Uu^{\circ}_{a,b}:=
\fD^{*}\left(\left(\Vv_{a,b}/\Vv_{a,b-1}\right)|_{X^{\circ}}\right) =
\left(\Vv_{a,b}/\Vv_{a,b-1}\right)|_{X^{\circ}}
$$
and then let
$\Uu_{a,b}$ be the unique torsion-free quotient of $\fD^{*}\Vv _{a,b}$
that restricts to $\Uu^{\circ}_{a,b}$ over $D^{\circ}$. It is
constructed by taking the saturated subsheaf extending the kernel and
taking the quotient.

For this discussion, one should note that $\Vv_{a,b}$
is a reflexive sheaf on $X$, indeed it is the dirct image
of a line bundle under a map from a normal (and
indeed smooth) variety to $X$. In particuler,
$\Vv_{a,b}$ is a vector bundle outside of codimension $3$, 
therefore $\fD^{*}(\Vv _{a,b})$
is a vector bundle outside of codimension $2$ on $D$. 

The fact that $\Uu_{a,b}$ has rank $2$ will become apparent
from the next lemma. 
\end{proof}

\begin{lemma}
Let $f _{\ExY /D}: \ExY \rightarrow D$ be the morphism induced by
$f$. That  is $f_{|\ExY} : \ExY \to X$ factors as
$$
\xymatrix@1{\ExY \ar[r]^-{f_{\ExY/D}} & D \ar[r]^-{\fD} & X.}
  $$
Then we have
$$
\Uu_{a,b} = f _{\ExY /D , *} (\Ll _{a,b} |_{\ExY}  ).
$$
\end{lemma}
\begin{proof}
Over $X^{\circ}$, we have $D^{\circ}= \Wob^{\circ}\subset X^{\circ}$
and $\fD_{*}\Uu ^{\circ}_{a,b} = \Uu ^{\circ}_{a,b} = f _{*}(\Ll
^{\circ}_{a,b}) / f_{*}(\Ll^{\circ} _{a,b}(-\ExY ^{\circ}))$.  This
quotient is equal, in turn, to $f_{*}\left(i_{\ExY^{\circ},*}\left(\Ll
^{\circ}_{a,b} |_{\ExY ^{\circ}}\right)\right)$, where
$i_{\ExY^{\circ}} : \ExY^{\circ} \to Y$ is the natural inclusion.
Over these smooth points $\fD_{*}$ induces an equivalence between
sheaves on $D^{\circ}$ and sheaves on $X^{\circ}$ supported on
$D^{\circ}$.  Also, $\fD\circ f _{\ExY /D} = f \circ i_{\ExY} $ so
$$
\fD_{*}f_{\ExY ^{\circ}/D^{\circ} , *} (\Ll ^{\circ}_{a,b} |_{\ExY ^{\circ}} ) = 
f_{*} (i_{\ExY ^{\circ},*}\Ll^{\circ} _{a,b} |_{\ExY ^{\circ}} ).
$$
From the equivalence we conclude that
$$ f_{\ExY ^{\circ}/D^{\circ} , *} (\Ll ^{\circ}_{a,b} |_{\ExY
  ^{\circ}} ) = \Uu ^{\circ}_{a,b}.
$$
This gives the claimed statement over $D^{\circ}$. 

Next, consider  the commutative diagram
$$
\xymatrix@R-0.5pc@C-0.5pc@M+0.5pc{
\ExY \ar[r] \ar[rd] & Y\times_{X} D \ar[d] \ar[r] & Y \ar[d] \\
& D \ar[r] & X
}
$$
we see that $\Vv _{a,b}=f _{*}(\Ll _{a,b})$
restricts over $D$ to the direct image under the 
middle downward map, and this maps to 
$f _{\ExY /D,*}(\Ll _{a,b}|_{\ExY} )$. 
In other words, we have a map expressing the 
relationship in question. This map is an isomorphism 
over $D^{\circ}$ by the preceding discussion. Since 
$\ExY $ is smooth $f _{\ExY /D,*}(\Ll _{a,b}|_{\ExY} )$
is torsion-free---in fact it is a bundle, by 
``miracle flatness''. Also, the map
$\ExY \rightarrow Y\times _XD$ is a closed embedding so
$$
f _{*}(\Ll _{a,b})\rightarrow 
f _{\ExY /D,*}(\Ll _{a,b}|_{\ExY} )
$$
is surjective. It follows that the right hand side is the unique
torsion-free quotient sheaf extending the quotient $\Uu^{\circ}_{a,b}$
already known over $D^{\circ}$.  This completes the proof that the
right hand side is $\Uu_{a,b}$. Notice that we obtain $\Uu_{a,b}$ as a
torsion free sheaf of rank $2$ because $\ExY \rightarrow D$ is a $2:1$
covering.
\end{proof}

\

\noindent
The next step is to use this
characterization to calculate the Chern character of $\Uu_{a,b}$.
This calculation may be truncated at codimension $1$ on $D$ 
since that corresponds to codimension $2$ on $X$. 

\begin{proposition}
\label{chu}
The Chern character of $\Uu_{a,b}$ truncated to 
codimension $1$ on $D$ is 
$$
{\rm ch}(\Uu_{a,b})= 2 +  (2a + b)H_D  - (3b+2) H_D^{\perp}
$$
where $H_D^{\perp}= [\Cbar ] - 16[\pp^1]$ is 
a class whose direct image to $X$ vanishes. 
\end{proposition}
\begin{proof}
Write 
$$
H_D:=\fD^{*}(H) = \alpha [\Cbar ] + \beta [\pp^1].
$$
Noting that $H_D\cdot [\pp^1]=1$ gives $\alpha =1$,
then noting that 
$$
H_D\cdot H_D= H^2 \cdot \Wob = 32
$$
(see Lemma \ref{wobbly8h}),
we get $2\alpha \beta =32$ so $\beta =16$. Thus
$$
H_D = \fD^{*}(H) = [\Cbar ] + 16[\pp^1].
$$
We have 
$$
\fD_{*}[\Cbar ]\cdot H = [\Cbar]\cdot H_D=16
$$
so (recalling $H^3=4$) we get $\fD_{*}[\Cbar]=4H^2$. 
We have that $\fD_{*}[\pp^1]$ is the class of
a line in $X$, this is $H^2/4$. 
We get
$$
\fD_{*}(H_D)= 8H^2.
$$
The class $H_D^{\perp}:= [\Cbar ] - 16[\pp^1]$
has the property that $\fD_{*}(H_D^{\perp})=0$.

Turn to the question of the proposition.
Consider the diagram
$$
\xymatrix@M+0.5pc{
\ExY   \ar[r]^-{i_{\ExY}}  \ar[d]_-{f_{\ExY /D}} & Y \ar[d]^-{f} \\
D \ar[r]_-{\fD}  & X.
}
$$
We have 
$$
i_{\ExY ,*}(\Ll _{a,b} |_{\ExY} ) = \Ll _{a,b}/ \Ll _{a,b - 1}.
$$
Applying $f _{*}$ to this gives $\Vv_{a,b}/\Vv _{a,b - 1}$.

Recall that $\Uu_{a,b}:= f_{\ExY /D,*}(\Ll _{a,b} |_{\ExY} )$.
From the diagram,
$$
\fD_{*}(\Uu_{a,b}) \cong \Vv_{a,b}/\Vv _{a,b - 1}.
$$
Therefore, in order to calculate ${\rm ch}(\Uu_{a,b})$ we would like
to know the relationship between it and 
${\rm ch}(\fD_{*}(\Uu_{a,b})$. Some further information will
be needed due to the fact that the codimension $1$ cycle
class group of $D\cong \Cbar \times \pp^1$
does not inject into that of $X$. 

First apply the Grothendieck-Riemann-Roch  formula for $\fD$. Recall that
$$
{\rm td}(TX)^{-1} = 1-H  + 5H^2/ 12.
$$
Truncated to codimension $1$ we get 
$\fD^{*}{\rm td}(TX)^{-1}=1-H_D$. 
On the other hand, $2g_{\Cbar }-2 = 32$ so 
$c_1(TD) = 2[\Cbar ] - 32 [\pp^1]$, and
$$
{\rm td}(TD) = 1 + [\Cbar ] - 16 [\pp^1].
$$
The Grothendieck-Riemann-Roch formula says
$$
\begin{aligned}
{\rm ch}(\fD_{*}(\Uu_{a,b}) & =
\fD_{*}(2 + c_1(\Uu_{a,b}))(1 + [\Cbar ] - 16 [\pp^1])
(1-H_D) \\
& = \fD_{*}\left(
2 + c_1(\Uu_{a,b}) + 2[\Cbar ] - 32 [\pp^1]
- 2[\Cbar ] - 32[\pp^1] \right).
\end{aligned}
$$
Noting that $\fD_{*}(2[D])=2[\Wob ]$ and recalling
from Lemma \ref{wobbly8h} that $[\Wob ]=8H$, then also using
$\fD_{*}[\pp^1]=H^2/4$, we get
$$
{\rm ch}(\fD_{*}(\Uu_{a,b})
= 16H + i_{*}c_1(\Uu_{a,b}) - 16H^2. 
$$
Denote by $\Delta _b$ the operation on a function $f$,
$\Delta _b(f):= f(a,b) - f(a,b-1)$. 
Thus for example  $\Delta _b(b)= 1$ and $\Delta _b(b^2) = 2b-1$
With this notation, 
$$
\begin{aligned}
{\rm ch}(\Vv _{a,b} /\Vv _{a,b -1}) & = 
\Delta _b \left[8+  (8a + 16 b)H + 
(4a^2 + 16 ab + 4 b^2 - 12 b - 2)H^2\right] \\
& =  16 H + (16a + 8 b - 16)H^2.
\end{aligned}
$$
The comparison ${\rm ch}(i_{*}(\Uu_{a,b})
=  {\rm ch}(\Vv _{a,b} /\Vv _{a,b -1})$ yields 
$$
\fD_{*}c_1(\Uu_{a,b}) = (16a+8b)H^2. 
$$
In particular, it follows that
$$
c_1(\Uu_{a,b}) = (2a + b)H_D + \xi H_D^{\perp} 
$$
for some $\xi$. 

On the other hand, we can note that the inverse
image in $\ExY$ of a line $\pp^1\subset D$ consists of two fibers in $\ExY$.
A fiber maps to a line in $X$ so its intersection with $\FxY $ has one point. 
On the other hand, $\ExY $ intersected with a fiber is $-1$. Thus, 
the bundle $\Ll _{a,b}$ restricted to a fiber
has class $\Oo _{\rm fib}(a-b-1 )$ so its direct image
to the  line $\pp^1\subset D$  
has $c_1=2 c_1[\Oo _{\pp^1}(a-b-1)]$. This tells us
that
$$
c_1(\Uu_{a,b}|_{\pp^1}) = 2a - 2 b - 2 .
$$
Thus
$$
( (2a + b)H_D +\xi H_D^{\perp}).[\pp ^1] = 2a - 2 b - 2 .
$$
This gives
$$
(2a+b) + \xi = 2a - 2 b - 2
$$
so $\xi = -2-3b$ and we get 
$$
c_1(\Uu_{a,b}) = (2a + b)H_D  - (3b+2) H_D^{\perp}.
$$
\end{proof}

\

\noindent
From this calculation we immediately get the following

\

\begin{corollary}
\label{udegzero}
Suppose $S\subset D$ is a curve in the class of $kD_H$,
and take $2a+b=0$. Then $\Uu_{a,b}|_S$ is a rank $2$ vector
bundle of degree $0$ on $S$. 
\end{corollary}

\subsection{Hyperplane section}

Let us now intersect the wobbly locus $\Wob \subset X$ with 
a general hyperplane section $X_H$, this denoting a smooth
divisor in the class of $H$ on $X$.

The class of the divisor $H$ on $D$ is $H_D$, 
described above as 
$$
H_D = [\Cbar]  + 16[\pp^1].
$$
One has $(H_D)^2 = 32$. 
We note also that $H_D$ is the trace of $X_{H}$ on $D$,
in particular its image in $X_{H}$ is $\Wob \cap X_{H}$. 
That is a curve $\Wob _H$ 
in $X_{H}$ whose normalization is 
isomorphic to $\Cbar $ since the 
projection $H_D\rightarrow \Cbar $
is an isomorphism. 

We would now like to count the double points of the curve $\Wob _H$. 
The intersection of the curve with the $6$ divisors 
$$
\Cbar \times \{\pw_i\}
$$
are the points in the normalization $\Cbar $
that get  identified in pairs to form the double points of
the image curve $\Wob _H=\Wob \cap X_H$. The number of points in
the intersection of $H_D$ with these $6$ horizontal
divisors is $6\cdot 16 = 96$, so there are $96$ points
lying over the double points. We conclude:

\begin{itemize}

\item There are $48$ double points on the curve $\Wob _H$.

\end{itemize}

We now consider the genus of the normalization of $\Wob _H$. 
The class of $\Wob _H$ is in $8H$ on $X_H$, so
$$
(\Wob _H)^2 = (8H)^2_{X_H} = 64H^3 = 64\cdot 4 = 256. 
$$
The genus of a smooth curve $S$  in this class is given
as follows: we recall that $K_{X_H} = -H$ so
$$
K_S = K_{X_H}\cdot S + \Oo _S (S) = -H\cdot (8H)_{X_H} +
(8H)\cdot (8H)_{X_H}
= 56H^3 = 224.
$$
This gives $2g_S-2 = 224$ so $g_S=113$. 

If there were no singularities except the $48$ double 
points, we would get the genus of the normalization of
$\Wob _H$ to be $113-48 = 65$. However, the normalization
is $\Cbar $ that has genus $17$. Hence, there
are additional singularities accounting for a
drop in the genus of $65-17= 48$. We first guessed then proved:

\begin{itemize}

\item The curve $\Wob _H$ has its $48$ nodes, plus $48$ cusps. 

\end{itemize}

The count of cusps comes from the structural result on the cuspidal
locus in Proposition~\ref{structureD}, proven in
Subsection~\ref{app-proof-cusp}. The number of cusps of is obtained as
follows. The structural result says that the cuspidal locus of $\Wob$
is the graph of a map $\Cbar \rightarrow \pp^1$ of degree $32$.  It
intersects $[\pp^1]$ in $1$ point and $[\Cbar]$ in $32$ points so it
intersects $H_D = [\Cbar] + 16 [\pp^1]$ in $48$ points.

\subsection{Summary of properties of our given sheaves}
\label{summary}

We are defining a vector bundle 
$\Vv_{a,b}$ on $X$ by taking the direct image of 
the line bundle $\Ll _0 (a\FxY  + (b+1)\ExY )$ for a line bundle 
$\Ll_0$ of degree $0$. We have calculated
its Chern character in Proposition \ref{chv},
$$
{\rm ch}(\Vv _{a,b}) = 8 +(8a + 16 b)H + 
(4a^2 + 16 ab + 4 b^2 - 12 b - 2)H^2.
$$

It is instructive to calculate the Bogomolov
$\Delta$-invariant. Notice that modifying $a$ corresponds to tensoring
with a line bundle pulled back from $X$, so we may assume $a=0$. Then
$$
\begin{aligned}
\Delta (\Vv _{0,b}) & = \frac{1}{2r}c_1^2 - {\rm ch}_2 \\
& = (16 b^2 -(4 b ^2 -12 b -2))H^2  \\
& = (12 b ^2 + 12b + 2)H^2.
\end{aligned}
$$
For integer values of $b$, the minimum of $2H^2$ 
is attained at both values $b=-1$
and $b=0$. Symmetry considerations suggest that
the desired parabolic weight should be $1/2$. 
with the parabolic structure creating an 
interpolation from $\Vv _{0,-1}$ to $\Vv _{0,0}$. 

One can see, by the way, that it is not correct to just
imagine a fractional value of $b$, indeed the
above formula at $b = -1/2$ would give $\Delta = -H^2$
contradicting the Bogomolov-Gieseker inequality. 
Nonetheless, the symmetry consideration leads us to try using
a parabolic weight of $1/2$ which will turn out to work, once the
second Chern class is adjusted appropriately at the singularities. 

Over $D$ we have a morphism 
$$
\fD^{*}(\Vv_{a,b}) \rightarrow \Uu_{a,b} \rightarrow 0
$$
and the Chern character of $\Uu_{a,b}$ on $D$ 
is 
$$
{\rm ch}(\Uu_{a,b})= 2 +  (2a + b)H_D  - (3b+2) H_D^{\perp}
$$
by Proposition \ref{chu}.

We would like to use this to define a parabolic 
structure along $\Wob $. It will be easier now to
pass to the hyperplane section $X_H$ and consider
the parabolic structure along $\Wob _H$. Let $D_H$ denote
the normalization of $\Wob _H$, and $H_{D_H}$ the
restriction of $H$ to here. 

Note that $H_D^{\perp}$ restricts to zero on the
hyperplane section $D_H$, so 
$$
{\rm ch}(\Uu_{a,b})|_{D_H} = 2 + (2a+b) H_{D_H}.
$$
It is not at all easy to understand how to calculate the
contributions of cusps to the parabolic $c_2$.  Taher's papers
\cite{Taher1,Taher2} consider a multiple point, rather than a cusp,
and in any case the procedure---which could theoretically be
implemented also for the cusp, or even imported by going to a covering
and transforming the cusp to a triple point---is rather complicated.

Therefore, we shall adopt a shortcut: we will just look at the
calculation for the parabolic weight $\alpha = 1/2$, and we will do
that by going to a cover that is ramified with ramification index of
$2$ over the curve $\Wob _H$.  Then, on the cover rather than defining
a parabolic structure we will use the quotient $\Uu$ to make a
modification of the bundle.  The construction follows the technique
described in Proposition \ref{pullbackcrude}.

\subsection{First parabolic Chern class}

The first parabolic Chern class may be described by a simplified integral
formula as just the average of ${\rm ch}_1(\Ee _a)$ over $a\in [0,1]$. 
In our situation of a crude parabolic structure with sheaves $\Vv = \Vv _{a,b}$ and
$\Vv '=\Vv _{a,b-1}$ and parabolic weight $\alpha = 1/2$ the formula becomes
$$
\begin{aligned}
{\rm ch}_1^{\rm par}(\Ee _{\cdot}) & = 
\int _{a=0}^{1} {\rm ch}_1(\Ee _a) da \\
& = \frac{{\rm ch}_1(\Vv ') + {\rm ch}_1(\Vv )}{2} \\
& = \frac{(8a + 16 (b-1))H + (8a + 16 b)H}{2} \\
& =( 8a + 16 b - 8)H.
\end{aligned}
$$
We should choose $a=1$ and $b=0$ to give ${\rm ch}_1^{\rm par}=0$. 
This is motivated by the symmetry consideration described above and the
fact that we are looking for a parabolic structure with vanishing
first Chern class. 

However, recall that modifying $a$ amounts to tensoring with a line
bundle on $X = X_1$. We can do that later. For simplicity we will
therefore assume $a=0,b=0$.

\subsection{Going to a cover of \texorpdfstring{$X_H$}{XH}}

In order to put this strategy into play, let us suppose
we are given a finite ramified covering $g:Z\rightarrow 
X_H$ together with a ramification  locus $R\subset Z$, such that
$g^{-1}(\Wob _H)=2R$, in other words $g$ is ramified with
ramification of order $2$ along $R$ over $\Wob _H$. 

We recall that the cusp is the downstairs branch locus
of the quotient 
$$
\cc ^2 \rightarrow \cc^2 /S_3 = \cc^2.
$$
Upstairs, the branch locus consists of $3$ crossed lines. 
Therefore, we may assume that this provides the local
model for our covering $Z\rightarrow X_H$ over cusps
of $\Wob _H$. This local model is a covering 
of degree $6$. 

On the other hand, we will suppose that the local model for
the covering over a double point, is just the product of
two ramified coverings with ramification two. This
local model is a covering of degree $4$. 

The covering $Z$ may be constructed globally by the Kawamata covering
trick \cite{Kawamata-cover}. There may be ramification elsewhere but
that does not need to worry us.

Let $d$ denote the degree of the covering $Z \to X_H$. 

\begin{itemize}
\item Over the $48$ cusps of $\Wob _H$ we will therefore have
$48d/6 = 8d$ triple points of $R$. 
\item Over the $48$ double points of $\Wob _H$ we will in turn have
$48 d/4 = 12d$ double points of $R$. 
\end{itemize}

\

Let $\widetilde{R}$ denote the normalization of $R$. 
Let $\Vv_Z$ denote the restriction of $\Vv=\Vv_{a,b}$ to $Z$. 
We have a map, defined over $\widetilde{R}$:
$$
\Vv_Z|_{\widetilde{R}} \rightarrow \Uu_{\widetilde{R}}.
$$
Here $\Uu_{\widetilde{R}}$ is the pullback of $\Uu_{a,b}$
to $\widetilde{R}$. Let 
$$
i_{\widetilde{R}} : \widetilde{R}\rightarrow Z
$$
denote the birational inclusion. We can define the ``normal bundle''
$N_{\widetilde{R}/Z}$ by
$$
i_{\widetilde{R}}^{*} K_Z  \otimes N_{\widetilde{R}/Z}
\cong K_{\widetilde{R}}.
$$
Let us calculate this in the following way. From the
picture, a double point leads to two extra self-intersection
points, whereas a triple point leads to $6$ extra
self-intersection points. Therefore, we get 
$$
N_{\widetilde{B}/Z} + 2(\mbox{double pts}) + 6 (\mbox{triple pts})
= R^2.
$$
Now $R$  has class $4H$ so $R^2 = (4H)^2 = 16 H^2$, which
pulled back to $Z$ together with $H^2=4$ on $X_H$, gives
$$
R^2 = 64d.
$$
Recall that we have $12d$ double points and $8d$ triple
points, so our formula gives
$$
\begin{aligned}
\deg N_{\widetilde{R}/Z} & = 64d - 2\cdot 12d - 6\cdot 8d \\
& = 64d - 24d - 48d \\
& = -8d.
\end{aligned}
$$
On the other hand, let us recall the Grothendieck-Riemann-Roch formula
(ignoring higher order terms)
$$
{\rm ch}(i_{\widetilde{R}*}\Uu_{\widetilde{R}}) 
= i_{\widetilde{R}*}{\rm ch}(\Uu_{\widetilde{R}}) - r_{\Uu_{\widetilde{R}}}  
i_{\widetilde{R}*} N_{\widetilde{R}/Z} /2 .
$$
Recall that $\Uu$ has rank $2$.  Thus the term
$-r_{\Uu_{\widetilde{R}}} i_{\widetilde{R}*} N_{\widetilde{R}/Z} /2$
is equal to $-(-8d)=+8d$ and so we get
$$
{\rm ch}(i_{\widetilde{R}*}\Uu_{\widetilde{R}}) 
= 2R + \deg \Uu_{\widetilde{R}} +8d.
$$
Recall that $2R = g^{*}(\Wob _H)$ has class $8H$. 

Let us now set $(a,b)=(0,0)$. As explained in the previous
subsections, this differs from the value we guess that it will be good
to look at, by a modification of $a$ that amounts to tensoring with a
line bundle on $X$, and we will recover that later.

With this choice, $\Uu$ has degree $0$
by Corollary \ref{udegzero}, so we get
$$
{\rm ch}(i_{\widetilde{R}*}\Uu_{\widetilde{R}}) 
= 8H  +8d.
$$
We now note that the map 
$$
V_Z\rightarrow i_{\widetilde{R}*}\Uu_{\widetilde{R}}
$$
is not going to be surjective. 

Here is a guess as to what the image is going to be.
The proof will be a consequence of the
Bogomolov-Gieseker inequality. 

A little like the principle that the two quotients
are independent at the double points, it looks likely
that the quotient over the double point is something
that is ``moving'', corresponding to a smooth curve
in the Grassmanian of the total bundle above the cusp. 
The guess is that the image of the above map
corresponds then, not to the structure sheaf of 
the embedded plane curve $R$, but rather to the
version of $R$ that has a triple point embedded as
a space curve. Let us call this curve 
$R^{\sharp}$, it is obtained by not glueing together
the points lying over double points, but by
glueing together the points lying over a triple point
in as small a way as possible. We get 
$$
0\rightarrow 
\Oo _{R^{\sharp}} \rightarrow 
\Oo _{\widetilde{R}} \rightarrow {\mathcal S} \rightarrow 0
$$
where ${\mathcal S}$ is a skyscraper sheaf that has length
$2$ over each triple point. (To get to the plane curve $R$
one would use a skyscraper sheaf of length $3$). 

\begin{lemma}
\label{sharp}
The image of the map 
$$
V_Z\rightarrow i_{\widetilde{R}*}\Uu_{\widetilde{R}}
$$
is the subsheaf 
$$
\Uu_{R^{\sharp}}\subset
i_{\widetilde{R}*}\Uu_{\widetilde{R}} .
$$
\end{lemma}

The proof will be deferred to below, since it is  motivated
by the following computations. 

With this hypothesis, we have 
$$
0\rightarrow \Uu_{R^{\sharp}}\rightarrow 
i_{\widetilde{R}*}\Uu_{\widetilde{R}} \rightarrow {\mathcal S}_U \rightarrow 
0
$$
where ${\mathcal S}_U$ has length $4$ over each triple point
($4$ because $U$ has rank $2$). 

Then, the hypothesis is that we have a surjection 
$$
V_Z\rightarrow \Uu_{R^{\sharp}}\rightarrow 0. 
$$
Let $V'_Z$ be the kernel, this will be our vector bundle over
the covering $Z$ corresponding to the ``pullback of the
parabolic bundle to $Z$''. 

Using the fact that there are $8d$ double points, we get
$$
\begin{aligned}
{\rm ch}(\Uu_{R^{\sharp}}) &
= {\rm ch}(i_{\widetilde{R}*}\Uu_{\widetilde{R}}) - 8d\cdot 4 \\
& = 8H +8d -32d = 8H -24d.
\end{aligned}
$$
Let us put this together with the Chern character of
$V_Z$. Recall that we are setting $(a,b)=(0,0)$ 
so ${\rm ch}_2(\Vv_{0,0}) = -2H^2$ on $X_H$.
Recall that $H^2$ has $4$ points on $X_H$, and then
we pullback by the covering $g$ of degree 
$d$, so we get
$$
{\rm ch}_2(V_Z) = -8d
$$
(here measuring in terms of numbers of points) on $Z$. 

Recall also that at $(a,b)=(0,0)$ we have ${\rm ch}_1(\Vv_{0,0})=0$. 

Altogether: 
$$
\begin{aligned}
{\rm ch}(V'_Z) & = {\rm ch}(V_Z) - 
{\rm ch}(\Uu_{R^{\sharp}})  \\
& = 8 - 8H -8d + 24d \\
& = 8 -8H + 16d \\
& = 8(1-H + 2d).
\end{aligned}
$$
We now note that $H^2/2 = 2d$ points on $Z$ so,
the term $2d$ above that is counted in terms of points on 
$Z$, may be written as $H^2/2$. We get the formula: 
$$
{\rm ch}(V'_Z)= 
8(1-H + H^2/2) .
$$
This clearly expresses the bundle $V'_Z$ as equivalent to 
$8$ times a line bundle since $1-H + H^2/2 = e^{-H}$. 
Therefore, just by inspection, its $\Delta$-invariant
vanishes and we have a projectively flat bundle. 

If we then tensor with $\Oo _{X}(H)$ we obtain a  Higgs bundle
with vanishing Chern classes.

\subsubsection{Proof of the lemma}

Following through what would happen in the above
calculations if the image were different, we 
can now give the proof of the lemma. 

\begin{proof}[Proof of Lemma \ref{sharp}]
We note that the map $V_Z\rightarrow
i_{\widetilde{R}*}\Uu_{\widetilde{R}}$ sends sections of $V_Z$ to
sections of $\Uu_{\widetilde{R}}$ that agree over the point where the
three branches come together. Thus, the morphism factors through a map
$V_Z\rightarrow \Uu_{R^{\sharp}}$. The question is to show that it is
surjective. Suppose not. Let $\Uu '$ denote the image, and let $\ell$
denote the total length of the quotient sheaf $\Uu '/\Uu_{R^{\sharp}}$
(take the sum over all the cuspidal points).

In the above discussion, let $V'_Z$ still denote the
kernel of the map $V_Z\rightarrow \Uu '$, so we have
a left exact sequence
$$
0\rightarrow V'_Z\rightarrow V_Z \rightarrow \Uu '\rightarrow 0.
$$
Now ${\rm ch}(\Uu ') = {\rm ch}(\Uu_{R^{\sharp}})  - \ell = 
8H -24d-\ell$, and following the same calculation as before
but with the extra term $\ell$ gives
$$
\begin{aligned}
{\rm ch}(V'_Z) & = 
{\rm ch}(V_Z) - 
{\rm ch}(\Uu ')  \\
& = 8 - 8H -8d + 24d +\ell \\ & = 8(1-H + 2d) + \ell  \\ & = 
8(1-H + H^2/2) + \ell .
\end{aligned}
$$
Recall that the Bogomolov-Gieseker inequality says that
$c_2$ can only go up from the flat case, hence
${\rm ch}_2$ can only go down. The bundle $V'_Z$ would
be a stable Higgs bundle on $Z$, and if $\ell >0$ this
would contradict Bogomolov-Gieseker. Thus, we conclude
$\ell = 0$ which completes the proof of the lemma. 
\end{proof}

\subsection{Degree one case---conclusion}

\begin{theorem}
\label{degonethm}
Suppose $\Ll _0 \in \op{Pic}(Y)$ is the pullback of a flat line bundle on the abelian variety $\Prym$, and set $\Ll := \Ll _0\otimes \Oo _Y(\ExY +\FxY )$. Let
$$
\Vv := f _{*}(\Ll ) = \Vv_{0,0}\otimes \Oo _X(H)
$$
provided with its natural crude parabolic structure with weight
$\alpha = 1/2$ (tensoring the one defined above with quotient $\Uu_{0,0}$,
by $\Oo _X(H)$). Then this is a crude parabolic logarithmic Higgs
bundle over $(X^{\circ},\Wob ^{\circ})$ that admits an extension to a
purely imaginary twistor $\srD$-module as in Theorem
\ref{outsidecod2}. The associated local system on $X-\Wob$ has rank
$8$ and its monodromy around $\Wob$ is semisimple with eigenvalues $1$
of multiplicity $6$ and $-1$ of multiplicity $2$.

The spectral line bundle on $Y$ is $\Ll ' =  \Ll _0\otimes \Oo _Y(\FxY )$. 
\end{theorem}
\begin{proof}
We saw in the previous subsections that choosing $\Vv _{a,b}$ for $(a,b)=(0,0)$ 
yields a projectively flat solution whose  (truncated) parabolic Chern character
is $8(1-H + H^2/2)$. 

Tensoring this by $\Oo _{X_1}(H)$ corresponds to choosing
$(a,b)=(1,0)$, taking the crude parabolic structure $\Vv _{1,0}$ with
subsheaf $\Vv ' _{1,0}= \Vv _{1,-1}$, and this yields a parabolic
Higgs bundle with vanishing parabolic Chern classes.  The rank is $8$
and the parabolic structure along $\Wob$ has parabolic weight $1/2$
with multiplicity equal to the rank of $\Uu$ which is $2$, completed
by parabolic weight $0$ with multiplicity $6$. This gives the stated
eigenvalues of the monodromy.

We note that the spectral line bundle, by definition the one whose
direct image from $Y$ to $X$ gives $\Ee _0$, is $\Ll '_{a,b} = \Ll
_{a,b-1}= \Ll _0\otimes \Oo _Y(\FxY )$.
\end{proof}

\section{The degree zero moduli space}
\label{chapter-d0}

In this section $X$ denotes the degree $0$ moduli space and $Y$
denotes its modular spectral covering of degree $8$ that is the
blow-up of the Prym variety $\Prym = \Prym_2$ in $16$ points.
Similarly $\Wob$ means $\Wob _0$ and so forth.

Narasimhan and Ramanan show  that $X\cong \pp^3$ \cite{NR}. There is a
\emph{\bfseries Kummer surface} of degree $4$ with $16$ nodes that we
will denote by $\Kum \subset X$. This is the moduli space of bundles
that are semistable but not stable, up to $S$-equivalence.  The
$S$-equivalence class corresponding to a point of $\Kum$ is
represented by a polystable bundle of the form $L\oplus L^{\vee}$ where
$L\in \op{Jac}^0(C)$ is a line bundle of degree $0$.

We have
$$
\Kum \cong \op{Jac}^0(C) / \pm 1 ,
$$
a point on $\Kum $ corresponding to a pair $(L,L^{-1})$ of
degree zero line bundles on $X$, up to interchangeing
$L$ and $L^{-1}$. The corresponding polystable bundle
is $E=L\oplus L^{-1}$. When $L=L^{-1}$, that is to say
at the $16$ points of order two on $ \op{Jac}^0(C)$,
we get a singular point of $\Kum$.  

Let $\Higgs _0$ be the corresponding Hitchin moduli space of Higgs bundles with
determinant $\Oo _C$, the Higgs field having zero trace.

\subsection{Geometry of the wobbly locus in degree zero}

\noindent
The Higgs bundle moduli space has a locus analogous to the Kummer surface. 

\begin{lemma}
\label{kumhiggs}
The Hitchin moduli space is singular along 
a codimension two subvariety $\Higgs _0 ^{\rm sing}$, and this
subvariety consists of the set of polystable but not 
stable Higgs bundles. We have
$$
\Higgs _0 ^{\rm sing}\cong \frac{  \op{Jac}^0(C)\times H^0(C,K_C)}{ \pm 1},
$$
a point here parametrizing
a Higgs bundle of the form $(L,\phi ) \oplus (L^{-1},-\phi )$. 
The set of $\cstar$ fixed points of $\Higgs _0 ^{\rm sing}$
is equal to the Kummer surface $\Kum \subset X$. Outside
of the $16$ singular points of $\Kum$, the  transverse local structure
of $\Higgs _0 $ along $\Higgs _0 ^{\rm sing}$ is that of a simple
double point. 
\end{lemma}
\begin{proof}
A stable Higgs bundle $\mathbb{E}$ has only scalar endomorphisms, and
by duality the trace-free $H^2$, that is to say $H^2_{\rm
  Dol}(End_0(\mathbb{E}))$, vanishes. This is the obstruction space to
deforming $\mathbb{E}$. Thus, a stable Higgs bundle is a smooth point
of $\Higgs_{0}$ (this is well-known, of course). Thus, a singular
point must be a strictly semistable point represented by a polystable
Higgs bundle. Since we are in rank $2$, it is a direct sum of line
bundles, and the condition that the determinant is trivial means that
it must have the stated form $(L,\phi ) \oplus (L^{-1},-\phi )$. The
structure is that of the moduli space of $(L,\phi )$ modulo the
involution $(L,\phi )\mapsto (L^{-1},-\phi )$.

A Higgs bundle of this form is fixed by $\cstar$ if and only if $\phi
= 0$, the fixed points are therefore inside the moduli space of
bundles $X\subset \Higgs_0$ and are in the Kummer surface.

Suppose we are at a point $\mathbb{L} \in \Kum$ corresponding to
$(L,\phi = 0)$ with $L\not\cong L^{-1}$. In other words $\mathbb{L} =
(L,0)\oplus (L^{-1},0)$, $H^2_{\rm Dol}(End^0(\mathbb{L}))=\cc$, and
$$
H^1_{\rm Dol}(End_0(\mathbb{L})) = H^1_{{\rm Dol}}(\Oo ) \oplus
H^1_{{\rm Dol}}(L^{\otimes 2} ) \oplus H^1_{{\rm Dol}}(L^{\otimes -2}
)
$$
with the Kuranishi map being given by the trace of the cup-product
(recall that the Goldman-Millson deformation theory is formal
\cite{GoldmanMillson,SimpsonHiggs}).  The cup product vanishes on
$H^1_{{\rm Dol}}(\Oo )$, which is the unobstructed deformation space
of $(L,0)$.  On $H^1_{{\rm Dol}}(L^{\otimes 2} ) \oplus
H^1_{\Dol}(L^{\otimes -2} )$ it is the same as the Poincaré duality
form. Thus, as a quadratic form it defines a simple double point in
the direction transverse to the deformation space of $(L,0)$ i.e. the
tangent space of $\Higgs _0 ^{\rm sing}$.  We see in particular that
$\Higgs _0$ is indeed singular along this locus.
\end{proof}

\

\noindent
In Proposition~\ref{descrip-fixed}{\bfseries (a)} we showed that the
set $\Higgs_{0}^{\cstar,\text{nu}}$ of $\cstar$-fixed points on
$\Higgs_{0}$ that are not on $X$, is a disjoint union of $16$ points.
These points cprrespond to the $16$ Higgs bundles
\[
\mathbb{E}_{\kappa} = (E_{\kappa},\theta_{\kappa}) = \left(
\kappa\oplus \kappa^{-1}, \, \begin{pmatrix} 0& 0 \\ 1 &
  0 \end{pmatrix}\right), \quad \kappa \in \mathsf{Spin}(C).
\]
From this picture, we can now describe the wobbly locus in $X$
completely.

\

\begin{proposition}
\label{wobblyc0}
The wobbly locus $\Wob \subset X$ decomposes as
$$
\Wob = \Kum \cup \left( \bigcup_{\kappa \in \mathsf{Spin}(C)}
\trope_{\kappa}\right)
$$ where each $\trope_{\kappa} \cong \pp^2\subset X=\pp^3$ is a plane
corresponding to the fixed point $\mathbb{E}_{\kappa}$.  Outside of
its $16$ nodes $\Kum$ is simply tangent to each $\trope_{\kappa}$ 
along a smooth conic $\conic_{\kappa} \subset \trope_{\kappa} \cong
\pp^2$, and any two trope planes intersect transversally.
\end{proposition}
\begin{proof}
We already described the components of the wobbly locus in the proof
of Theorem~\ref{qresolves} but for completeness we recall the relevant
arguments here.  A wobbly bundle $E$ is either strictly semistable, in
which case it corresponds to a point of $\Kum$, or it is stable but
has a nonzero nilpotent Higgs field. In the second case, the upward
limit of the $\cstar$ orbit corresponding to that Higgs field is a
fixed point not in $X$, so by
Proposition~\ref{descrip-fixed}{\bfseries (a)} the limiting Higgs
bundle has the form $\kappa\oplus \kappa^{-1}$ where $\kappa$ is a
square-root of $\omega_C$. In terms of the bundle $E$, it means that
there is an exact sequence
$$
0\rightarrow \kappa^{-1} \rightarrow E \rightarrow \kappa \rightarrow 0.
$$
Thus $E$ corresponds to an extension class $\xi \in
H^1(\kappa^{-2}) = H^1(\omega _C^{-1}) \cong H^0(\omega _C ^{\otimes
  2} ) ^{\vee}$.  The set of bundles corresponding to extension
classes like that is a $\pp^2$. Narasimhan-Ramanan's discussion
\cite[Proposition 6.1, Theorem 2]{NR} shows that the map from the
projectivized space of extension classes into $X=\pp^3$ is linear,
thus we obtain a linear $\trope_{\kappa} \cong \pp^2 \subset \pp^3$.
These define the $16$ divisors $\trope_{\kappa}$. They are different
since the fixed points are different.  The intersections
$\trope_{\kappa}\cap \trope_{\kappa'}$ are therefore transverse.  Note
also that, as we will see in Corollary \ref{nothree} below, three
distinct trope planes can not meet along a line, so the intersections
between trope planes are normal crossings in $X$-codimension $2$.

We consider next the condition that an extension of the above form lies in
the Kummer surface.  It means that there is a line bundle $U$ of
degree $0$ with a map $U\rightarrow \kappa$ such that the extension
splits. We may write $U = \kappa(-t)$ for a point $t\in C$. Splitting the
extension means that the image of $\xi$ under the map
$$
H^1(\omega _C^{-1}) \rightarrow H^1(\omega _C^{-1}(t)) \cong H^1(\Oo _C(-t'))
$$
should vanish,  where 
$t'$ is the hyperelliptic conjugate of $t$.
This determines an extension $\xi$ uniquely up to scalars. 
Furthermore, 
the extension determined by $t'$ yields an $S$-equivalent bundle since
$$
\kappa(-t) \otimes \kappa(-t') = \omega _C (-t-t') = \Oo _C.
$$
Thus, the map from $C$ to the space of extensions factors through a
map from $\pp^1$.  Notice that the map from $C$ to $\Kum$ comes by
sending $t$ to the line bundle $\kappa(-t) \in \op{Jac}^{0}(C)$ and
then projecting from the Jacobian to the Kummer.

Let us consider a different plane given by the space of extensions of the form
$$
0\rightarrow V^{-1} \rightarrow E \rightarrow V \rightarrow 0
$$
for $V$ some general line bundle of degree $1$.  This plane is the
space of bundles such that $H^0(V\otimes E) \neq 0$. If that holds for
$E$ then it holds for the polystable bundle in the $S$-equivalence
class of $E$, so it means that we look for the condition $H^0(V\otimes
\kappa(-t))\neq 0$.  We note that the degree $2$ bundle $V\otimes
\kappa$ is general, so it has a single section with two zeros. If $t$
is one of these zeros, then we get a solution. The case of a
hyperelliptically conjugate point is when the other piece
$H^0(V\otimes \kappa(-t'))\neq 0$. There are two points in $\pp^1$ that are
the images of the zeros of the section of $V\otimes \kappa$, corresponding
to the intersection with this plane. We conclude that the rational
curve image of $\pp ^1$ is a conic.

Thus, $\Kum \cap \trope_{\kappa}$
has reduced scheme equal to a conic, so it is twice a
conic in the
plane $\trope_{\kappa} \cong \pp^2$. 
\end{proof}

\

\begin{remark}
\label{sixteensix}
In the classical studies of the geometry of the Kummer surface (see
\cite{Klein,GH,NR,BeauvilleSurfaces,Keum,Dolgachev200,Hudson}), these
$16$ planes were known as the {\em trope planes}. Each trope plane
passes through $6$ of the $16$ singular points of the Kummer surface,
and each singular point is contained in $6$ of the $16$ trope
planes. This configuration was known as the \emph{\bfseries Kummer
$16_6$ configuration}.  The conics $\conic_{\kappa}$ are called
\emph{\bfseries trope conics}.  We have adopted the `trope'
terminology.
\end{remark}

\

\begin{proposition}
Inside one of the trope plane $\trope_{\kappa}$, the trope conic may
be identified with the hyperelliptic $\pp^1$. The six branch points in
the hyperelliptic $\pp^1$ correspond to the $6$ nodes of $\Kum$
contained in that trope plane. The $15$ lines of intersection with the
other trope planes, consist of all lines passing through pairs of
branch points (i.e. nodes of $\Kum$ on the trope conic.
\end{proposition}
\begin{proof}
In modular terms, a trope plane corresponds to a choice of a
square-root $\kappa$ of $\omega _C$, with the corresponding
$\trope_{\kappa} \subset X\cong \pp^3$ consisting of all the bundles
that are extensions of the form
$$
0 \rightarrow \kappa^{-1} \rightarrow E \rightarrow \kappa \rightarrow 0.
$$
Such a bundle $E$ belogs to $\Kum$ if it contains a degree $0$ line
bundle; such would be of the form $L = \kappa(-t)$ with the dual
$L^{-1} = \kappa(-t')$ where $t'$ is the hyperelliptic conjugate of
$t$. The bundle determines in this way, and is determined by, a point
on the hyperelliptic $\pp^1$. The nodes of the Kummer surface occur
when $L^{-1}\cong L$, thus $t'=t$ in other words these correspond to
branch points in the hyperelliptic $\pp^1$.  Now, given two different
branch points $\brx_1$ and $\brx_2$ we get $(\kappa')^{-1} :=
\kappa(-\brx_1-\brx_2)$. This is the dual of another square-root
$\kappa'$ of $\omega _C$. The set of bundles $E$ that admit non-zero
maps from $(\kappa')^{-1}$ is the intersection of $\trope_{\kappa}$
and $\trope_{\kappa'}$. This set contains in particular the bundles
containing $\kappa(-\brx_1)$ respectively $\kappa(-\brx_2)$, i.e. this
intersection line passes through the two branch points $\brx_1$ and
$\brx_2$ on the trope conic. We have identified a set of $15$ lines
that is on the one hand the set of intersection lines between
$\trope_{\kappa}$ and all the other planes, and on the other hand is
the set of lines through pairs of the $6$ branch points.
\end{proof}

\

\begin{corollary}
\label{nothree}
No three trope planes pass through the same line.
\end{corollary}
\begin{proof}
Recall that we are assuming genericity of the set of $6$ branch points
of the hyperelliptic curve. If there were three trope planes passing
through the same line, then another trope plane would contain a
configuration of three lines of the form $\mathsf{L}_{\kappa\kappa'} =
\trope_{\kappa}\cap \trope_{\kappa'}$ as in the proof of the
proposition, that pass through the same point. However, noting that
all smooth plane conics are projectively equivalent to a given fixed
one, our set of $6$ branch points is general on the conic. If three of
the lines $L_{\kappa\kappa'}$ were to pass through the same point, we
could jiggle one of the branch points to move one of the three lines
away from that intersection point, so this can not happen under our
genericity hypothesis.
\end{proof}

\

Let $\Prym \subset \Higgs_{0}$ be the fiber of the Hitchin fibration
corresponding to a spectral curve $\pi : \Ctilde \rightarrow C$. Notice
that, although the Hitchin fibrations are different for degree $0$ and
$1$, the base is still the space of quadratic differentials and the
spectral curves are the same. Thus, Proposition \ref{prop:hitchin.map}
applies, and $\Ctilde / C$ is branched over two pairs of
hyperelliptically conjugate points $a,a',b,b'$.

Then, as we saw in Theorem~\ref{qresolves} we need to blow up $\Prym =
\Prym_{2}$ at the $16$ points $\{ \pi^{*}\kappa \}_{\kappa \in
  \mathsf{Spin}(C)}$ to get $Y$ with a map $f : Y \rightarrow X$.
Let $H$ be the hyperplane class on $X$, let $\FxY$ be its pullback to
$Y$, and let $\ExY = \sum_{\kappa \in \mathsf{Spin}(C)} \ExY_{\kappa}$
be the exceptional divisor (union of $16$ disjoint $\pp^2$'s).

\

\begin{lemma}
\label{tropes}
Under the finite morphism $f : Y \to X$ each plane $\ExY_{\kappa}$
maps isomorphically to the corresponding plane $\trope_{\kappa}$ in
$X$, the map from $Y$ to $X$ being ramified with index $2$ there, and
this accounts for the ramification over general points of
$\trope_{\kappa}$.
\end{lemma}
\begin{proof}
We have seen in the proof of Theorem~\ref{qresolves} and in
Proposition \ref{wobblyc0} that the plane $\trope_{\kappa}$ is the
projectivized set of downward directions at the fixed point, which is
isomorphic to the projectivized normal bundle of the upward subset
emanating from the fixed point. This in turn is isomorphic to the
projectivized normal bundle of the intersection point of the upward
set with $\Prym$, which is the exceptional divisor $\ExY_{\kappa}$.

Let us recall why the upward subset intersects the Prym in a single
point (this explains why there are $16$ points to blow up, as stated
above). Our fixed point corresponds to a nilpotent Higgs bundle with
underlying bundle of the form $\kappa\oplus \kappa^{-1}$ where
$\kappa$ is a square-root of $\omega_C$. The upward subset is the set
of Higgs bundles whose underlying vector bundle is unstable with
detabilizing subbundle $L$. To intersect with the Prym, we look for
line bundles $U$ on $\Ctilde$ of degree $2$ whose norm is $\omega _C$,
such that there is a map $\kappa\rightarrow \pi _{*}(U)$ or
equivalently $\pi ^{*}\kappa \rightarrow U$. But $\pi ^{*}\kappa$ also
has degree $2$ on $\Ctilde$ so $U= \pi ^{*}\kappa$ is the unique
solution.

There is an involution of $\Prym$ given by applying the hyperelliptic
involution on $\Ctilde$ (which covers the hyperelliptic involution on
$C$). Up to translation by $\omega _C$, this amounts to doing $-1$ on
the abelian variety, and the $16$ points being blown up are the fixed
points. Normal directions to these points are also fixed, so the
$\ExY_{\kappa}$ are fixed points of the induced involution of $Y$.

The pullback of a semistable degree $0$ bundle by the hyperelliptic
involution on $C$ is $S$-equivalent to itself. This is known for local
systems in the work of Goldman and Heu-Loray \cite{Goldman,
  HeuLoray-flat}, in particular it applies to unitary local systems
and hence to polystable bundles.

The map $f : Y\rightarrow X$ therefore factors through the quotient by
the involution.  But the involution acts by $-1$ on the normal
directions, so it acts by $-1$ on the normal bundle of
$\ExY_{\kappa}$, giving a simple ramification divisor of the map from
$Y$ to its quotient by this involution along $\ExY_{\kappa}$.

We claim that this accounts for the ramification of the map $f :
Y\rightarrow X$ along $\ExY_{\kappa}$.  Suppose that a point of
$\ExY_{\kappa}$ counts for $m$ points in the fiber of $f :
Y\rightarrow X$ that has degree $8$ (Lemma \ref{deg8c1}). We are
claiming that $m=2$, and have shown in the preceding paragraph that
$m\geq 2$. A monodromy argument (when moving everything around, the
planes $\ExY_{\kappa}$ get interchanged transitively) tells us that
this coefficient $m$ is the same for all $\ExY_{\kappa}$. Now,
$\trope_{\kappa}$ is a plane inside $X=\pp^3$, so a different
component $\trope_{\kappa'}$ intersects it in a line. Consider yet a
third component $\trope_{\kappa''}$ that intersects this line in at
least one point. This is not a smooth point of the Kummer surface,
since the trope planes are tangent planes to the Kummer surface at
points on the trope conics that are different from the $16$ singular
points, and the tangent plane is unique so in that case the three
planes would be the same. The point might be one of the $16$ singular
points of the Kummer surface. However, that does not always
happen. Indeed, the other trope planes intersect $\trope_{\kappa}$ in
$15$ lines, and the various intersection points of these lines can not
be limited to only the singular points of the Kummer surface (of which
there are $6$ corresponding to the Weierstrass points inside the
$\pp^1$ of the trope conic). So we can choose $\trope_{\kappa'}$ and
$\trope_{\kappa''}$ to correspond to lines that intersect somewhere in
$\trope_{\kappa}$ other than a singular point.

Now, the exceptional divisors upstairs in $Y$ are three disjoint
planes $\ExY_{\kappa},\ExY_{\kappa'},\ExY_{\kappa''}$.  In particular,
lying over the intersection point inside $X$, there are at least three
different points in the
$\ExY_{\kappa},\ExY_{\kappa'},\ExY_{\kappa''}$. This would give
$3m\leq 8$, from which we conclude $m=2$.

We now know that in a fiber over a general point of $\trope_{\kappa}$,
there is a ramification point in $\ExY_{\kappa}$ counting for $2$
points in the degree $8$ fiber, plus $6$ other points on $Y$ that are
not in the exceptional divisor, in other words they are in the locus
of $\Higgs_0$ whose underlying bundle is stable. An argument similar
to the proof of Lemma \ref{ram-mov} below shows that for a general
point in the Hitchin base, and over a general point of
$\trope_{\kappa}$, the other $6$ points are unramified.
\end{proof}

\

\subsection{Computations in degree zero}

Recall that by convention we drop terms of degree $\geq 3$ in all expressions.

\begin{proposition}
\label{todd0}
The Todd classes for the degree $0$ moduli spaces are: 
$$
\begin{aligned}
{\rm td}(TX) & = 1 + 2H + 11H^2/6, \\ 
f^{*}{\rm td}(TX)^{-1} & = 1 - 2\FxY  + 13\FxY ^2/6, \\
{\rm td}(TY) & = 1 - \ExY  +  \ExY ^2/3, \\ 
{\rm td}(Y/X) & = (1 - 2\FxY  + 13\FxY ^2/6)(1 - \ExY  +  \ExY ^2/3).
\end{aligned}
$$ 
\end{proposition}
\begin{proof}
For $Y$ which is the blow-up of an abelian variety at $16$ points we saw
$$
{\rm ch}(TY) = 3 - 2\ExY  + 2\ExY ^2.
$$
In particular $c_2(TY) = 0$. 
Thus
$$
\begin{aligned}
{\rm td}(TY) & = 1 + c_1/2 + (c_1^2+c_2)/12 \\
& = 1 - \ExY  +  \ExY ^2/3 .
\end{aligned}
$$
Next,
$$
\begin{aligned}
{\rm td}(T\pp^n) & = 1 + c_1/2 + (c_1^2+c_2)/12 \\
& = 1 + \frac{n+1}{2}H + \frac{3n^2 + 5n+2}{24}H^2.
\end{aligned}
$$
For $X=\pp^3$ this gives 
$$
{\rm td}(TX) = 1 + 2H + 11H^2/6 ,
$$
hence 
$$
{\rm td}(TX)^{-1} = 1 - 2H + 13H^2/6 .
$$
Thus
$$
(f^{*} {\rm td}(X))^{-1} = 1 - 2\FxY  + 13\FxY ^2/6. 
$$
\end{proof}

\

\begin{proposition}
\label{intersections0}
The triple intersections of divisor classes on the degree $0$ modular
spectral covering $Y$ are:
$$
\FxY ^3 = 8, \qquad
\ExY \FxY ^2 = 16, \qquad
\ExY ^2\FxY  = -16, \qquad
\ExY ^3 = 16.
$$
\end{proposition}
\begin{proof}
The degree of $f : Y\rightarrow X$ is $8$ and we have $H^3=1$ in
$X=\pp^3$, so $\FxY ^3 = 8$. Next, $\FxY ^2$ is the pullback of a line
that intersects the $16$ trope planes in $16$ points. As the map $\pi$
identifies the trope planes with the corresponding components of $\ExY
$ we have $\ExY \FxY ^2 = 16$.

Next, for each $\ExY_{\kappa}$ the self-intersection $\ExY_{\kappa}^2$
is the class on $\ExY_{\kappa}$ of the normal bundle.  The normal
bundle of $\ExY_{\kappa}$ in $Y$ is $\Oo_{\ExY_{\kappa}}(-1)$ so it is
minus the class of a line in $\ExY_{\kappa}$.  We get that
$\ExY_{\kappa}^3$ is the intersection of two of these together, so it
is $1$. Adding $16$ of these together gives $\ExY^3=16$. The
intersection of the line with $\FxY$ is the same as the intersection
of its image in $X$, which is again a line in the trope plane, with
$H$. This is $1$ so with the minus sign and $16$ planes we get $\ExY
^2\FxY =-16$.
\end{proof}

\

\begin{proposition}
\label{piLzero}
Suppose $\Ll = \Oo _Y(a\FxY  + b\ExY)$ is a line bundle on $Y$ and let 
$\Ee := f_{*}\Ll$. Then 
$$
{\rm ch}_1(\Ee ) = (8(a-2) + 16 (b-1))H
$$
and 
$$
{\rm ch}_2(\Ee ) = 4(a^2 + 4ab - 2b^2 -8a - 4b + 11)H^2 .
$$
\end{proposition}
\begin{proof}
The Grothendieck-Riemann-Roch formula says
$$
\begin{aligned}
{\rm ch}(\Ee ) &  = 
f_{*} \left( {\rm td}(Y/X)  e^{\Ll}\right)  \\
& = f_{*} \left[
(1 - \ExY  + \ExY ^2 /3) (1 - 2\FxY  + 13\FxY ^2/6 ) (1 + 
(a\FxY  + b\ExY ) + (a\FxY  + b\ExY )^2 / 2) 
  \right] .
\end{aligned}
$$
Thus 
$$ H^2\cdot {\rm ch}_1(\Ee ) = \FxY ^2 ((a-2)\FxY + (b-1)\ExY ) =
8(a-2) + 16 (b-1).
$$
And 
$$
\begin{aligned}
H\cdot {\rm ch}_1(\Ee ) & = \FxY (\ExY ^2 /3 + 13\FxY ^2/6 +  
(a\FxY  + b\ExY )^2 / 2 + 2\ExY \FxY  - (a\FxY +b\ExY )(\ExY + 2\FxY )) \\
& =\ExY ^2 \FxY  (1/3 + b^2 / 2 - b) 
+ \ExY  \FxY ^2 (ab + 2 - a - 2b) 
+ \FxY ^3 (13/6 + a^2 / 2 - 2a)
\end{aligned}
$$
and using Proposition \ref{intersections0} this becomes
$$
\begin{aligned}
4 a^2 & + 16 ab - 8 b^2 - 32 a -16b + (52 / 3 + 32 - 16/3)  \\
& = 4(a^2 + 4ab - 2b^2 -8a - 4b + 11).
\end{aligned}
$$
Noting that $H^3=1$ on $\pp^3$ this gives the required formula. 
\end{proof}

\

\begin{corollary}
\label{piLzerocor}
If we impose the condition ${\rm ch}_1(\Ee ) =0$ by 
setting 
$$
\Ll = \Oo _Y (\ExY +2\FxY  + m(\ExY -2\FxY )) \;\;\mbox{ with }m\in \zz
$$
then
$$
{\rm ch}_2(\Ee ) =(-24m^2 + 4)H^2 .
$$
\end{corollary}
\begin{proof}
Indeed we then have $a=2(1-m)$ and $b= (m+1)$ so the formula
becomes 
$$
\begin{aligned}
4(4(1-m)^2 & + 8(1+m)(1-m) - 2(m+1)^2  - 16(1-m) -4(m+1) + 11)  \\
& =
4(4m^2 -8m + 4 + 8 -8m^2 -2m^2 - 4m -2 -16 + 16m -4m -4  + 11) \\
& =4(-6m^2 + 1) = (-24m^2 + 4).
\end{aligned}
$$
\end{proof}

\

\begin{corollary}
\label{extremalzerocor}
The extremal value for ${\rm ch}_{2}(\Ee)$ is at $m=0$ in other words
$\Ll = \Oo _Y(\ExY + 2\FxY )$ and then
$$
{\rm ch}_2(\Ee ) = {\rm ch}_2(\pi _{*}(\Oo _Y(\ExY +2\FxY ) ) =4H^2 .
$$
\end{corollary}

\

\noindent
Note that the extremal value for $\Ll$ is the same as $-{\rm
  td}_1$. This is a general phenomenon. At this value, ${\rm ch}_2(\Ee
)=4H^2$ contradicts the Bogomolov-Gieseker inequality. That indicates,
on the one hand, that $\Ee$ is not a stable bundle. On the other hand,
it is stable as a bundle with meromorphic Higgs field.  The Higgs
field is logarithmic along the smooth points of the wobbly divisor,
but the contradiction to the Bogomolov-Gieseker inequality indicates
that if we pull back to a resolution making the wobbly divisor into a
normal crossings divisor, the resulting Higgs field there is not going
to be logarithmic.  We will see that in more detail in the next
sections, in which we investigate what parabolic structure can be
put over the resolution of singularities in order to get a logarithmic
Higgs field with maximal value of ${\rm ch}_2$. We will see that the
maximal value is then $0$ yielding a flat bundle.

We may note here a general principle: the corrections to ${\rm ch}_2$
coming from the required parabolic structures at the singularities,
over a planar slice (so the singularities are tacnodes), are local and
do not depend on $\Ll$. Therefore, whatever they are, the extremal
value has to be obtained when the ${\rm ch}_2$ calculated above is
extremal. This means that $\Ll$ has to be numerically equivalent to
$\Oo _Y(\ExY +2\FxY )$.

\subsection{Ramification and Riemann-Roch}

Recall that $f : Y \rightarrow X$ is a finite covering of degree
$8$. We would like to understand the ramification. First is the
ramification over the wobbly locus.

\begin{proposition}
\label{descrip-ram}
The ramification of the map $f : Y \rightarrow X$ over general points
of the components of the wobbly locus is as follows. Over a general
point of the Kummer surface $\Kum \subset X$, the covering $f$ is
fully ramified, breaking into four pieces with simple
ramification. Over each of the $16$ trope planes, the covering $f$
breaks generically into a simply ramified map along the corresponding
component of $\ExY $ plus a degree $6$ etale cover.
\end{proposition}
\begin{proof}
The part about ramification over the trope planes was shown in Lemma
\ref{tropes}. We need to describe the ramification over a general
point of the Kummer surface.

Let us describe the points of $Y$ lying over a general point of
$\Kum$. If $E \in \Kum$ is such a point, there is a line bundle $L$ of
degree $0$ on $C$ with an inclusion $L\hookrightarrow E$. If $E=\pi
_{*}(U)$ for a degree $2$ line bundle $U$ on $\Ctilde$, by adjunction
it means that there is a nonzero map $\pi ^{*}(L) \rightarrow U$,
hence
$$
U = \pi ^{*}(L) \otimes \Oo _{\Ctilde} (\tilde{u} + \tilde{v})
$$
for two points $\tilde{u}, \tilde{v}\in \Ctilde$. Let $u,v$ denote
their images in $C$. The determinant of $\pi _{*}(U)$ is $L^{\otimes
  2} \otimes \omega _C^{-1}(u+v)$ so the trivial determinant condition
$\det E = \mathcal{O}_{C}$ says
$$
\Oo _C(u+v) = \omega _C \otimes L^{\otimes -2}  .
$$
This determines the points $u,v$ as the zeros of the unique section of
$\omega _C \otimes L^{\otimes -2}$ (for generic $L$). Then,
$\tilde{u}$ and $\tilde{v}$ are liftings of these points to
$\Ctilde$. There are $4$ combinations of liftings. Conversely, each
one determines a bundle $U$ such that there exists a map
$L\hookrightarrow \pi _{*}(U)$.

We claim that $\pi _{*}(U)$ is polystable. This may be seen by noting
that as we move around in $\Kum$, the line bundles $L$ and $L^{-1}$
interchange, so if there is an inclusion from $L$ at a general point
there has to be an inclusion from $L^{-1}$ too. It may also be seen as
follows.  Recall that $\Ctilde$ is hyperelliptic too, and let $\sigma$
denote its hyperelliptic involution, that covers the hyperelliptic
involution $\hi_{C}$ of $C$.  Thus, $\sigma ^{*} \pi ^{*}(L) = \pi
^{*}(L^{-1})$ since $L^{-1}$ is the pullback of $L$ by the
hyperelliptic involution on $C$. We obtain a map $L^{-1} \rightarrow
\pi _{*}(\sigma ^{*}(U)) = \hi_{C} ^{*}(E)$. As noted in the proof of
Lemma \ref{tropes}, $\hi_{C} ^{*}(E) \cong E$ so we get a map $L^{-1}
\hookrightarrow E$.  For $L$ general, this gives an isomorphism
$L\oplus L^{-1} \rightarrow E$.

Using the fact that $\pi _{*}(U)$ are polystable, we can calculate
the fiber using either one of the degree $0$ line bundles contained in
$E$. The previous calculation shows that there are $4$ points.

As $L\oplus L^{-1}$ moves around in $\Kum$, the points $u,v$ move
around in $C$ and this induces a transitive action on the set of $4$
liftings. Thus, the ramification degree of $Y/X$ at the different
points must be the same, and since $\deg(Y/X)=8$ we conclude that
the $4$ points each have ramification degree $2$, in other words
simple ramification. This completes the proof.
\end{proof}

\

\noindent
Next we calculate the class of the ramification divisor of the map $f
: Y\rightarrow X$.  We have $\omega _X = \Oo _X (-4H)$ so $f^{*}
\omega _X ^{-1} = \Oo _Y (4\FxY )$.  From the formula for the
canonical class of the blow-up, we have $\omega _Y= \Oo _Y(2\ExY
)$. This gives
$$
\Oo _Y(R) = \omega _Y \otimes f^{*} \omega_X^{-1} = \Oo _Y( 2\ExY +
4\FxY ).
$$
The ramification divisor includes a copy of the reduced
divisor of the inverse image of the
Kummer surface ${\rm Kum}$ since $\pi$ is fully ramified over ${\rm Kum}$.
Let us call this part $R^{\rm Kum}$. As $\Oo _X({\rm Kum}) = \Oo _X(4H)$,
but the pullback of ${\rm Kum}$ is twice  $R^{\rm Kum}$ we get
$$
R^{\rm Kum} = 2\FxY .
$$ Let $R^{\rm Tro} = \ExY $ be the part of the ramification over the
trope planes.  Let $R^{\rm Mov}$ denote the novable part $R^{\rm Mov}
:= R - R^{\rm Kum} - R^{\rm Tro}$ of the ramification divisor. We therefore have
$$
R^{\rm Mov} \in  \left|\Oo_{Y}\left(2\FxY  + \ExY\right)\right|.
$$

\

\begin{lemma}
\label{ram-mov}
The ramification $R^{\rm Mov}$ is movable in the following sense: all
components of the image of this divisor in $X$ move as a function of
the point in the Hitchin base.
\end{lemma}
\begin{proof}
Let $B^{\rm Mov}\subset X$ be the image of $R^{\rm Mov}$. This does
not contain any component of the wobbly locus, indeed we have seen in
Proposition \ref{descrip-ram} that the only ramifications over the
various components of the wobbly locus are those given by $R^{\rm
  Kum}$ and $R^{\rm Tro}$.  Suppose that some component was
fixed. That means that we would have a divisor $B'\subset B^{\rm Mov}$
remaining fixed as our point in the Hitchin base moves.  In
particular, for a general point $x\in B'$, which is very stable, the
well-defined fiber $\aaaa^3$ of the projection $\Higgs _0 \rightarrow
X$ over $x$ would have a ramification point inside every Hitchin fiber
that it intersects. This would give a $3$-dimensional family of
ramification points, which would have to be all of $\aaaa^3$. In that
case, the projection map from $\Higgs _0$ to $X$ (over the very stable
open subset) would be non-smooth there, but the general Hitchin fiber
is smooth so that can not be the case. This contradiction shows that
there are no non-movable components in $B^{\rm Mov}$.
\end{proof}

\

\noindent
Suppose $\ExY_{\kappa}$ is one of the components of $\ExY$.  It
follows from the lemma and the previous proposition, that $R^{\rm
  Mov}$ does not include $\ExY_{\kappa}$. Therefore $R^{\rm Mov} \cap
\ExY _{\kappa}$ is a transverse intersection. Recall that
$\ExY_{\kappa} \cong \pp^2$.  The map $\ExY_{\kappa} \rightarrow X$ is
an isomorphism to a trope plane in $X\cong \pp^3$, so $\Oo
_{\ExY_{\kappa}}(\FxY \cap \ExY_{\kappa}) = \Oo
_{\ExY_{\kappa}}(1)$. On the other hand, $\Oo _{\ExY_{\kappa}}(\ExY )
\cong \Oo _{\ExY_{\kappa}}(-1)$ since $\ExY_{\kappa}$ is the
exceptional divisor from blowing up $\Prym $ at a point. We conclude
that
$$
\Oo _{\ExY_{\kappa}}(R^{\rm Mov} \cap \ExY_{\kappa}) \cong \Oo _{\ExY_{\kappa}}(1).
$$
This proves the following

\

\begin{corollary}
The divisor $R^{\rm Mov} \cap \ExY_{\kappa} \subset \ExY_{\kappa}$ is a line
in the plane $\ExY_{\kappa}$. In particular, it does not contain the trope
conic which is the inverse image in $\ExY_{\kappa}$ of the intersection of
the trope plane with the Kummer surface.
\end{corollary}

\

\noindent
Consider a general point $y\in R^{\rm Kum} \cap R^{\rm Tro}$. From the
previous corollary, the ramification of the map $f$ in a
neighborhood of $y$ is only the ramification due to the two pieces
$R^{\rm Kum}$ and $R^{\rm Tro}$.

\begin{corollary}
\label{coords}
In the neighborhood of such a point $y$, the map $f$ has a piece of
degree $4$ given in local coordinates $(u,v,z)$ of $Y$ and $(x,y,z)$
of $X$ by
$$
x=u^2 + v, \qquad  y = v^2.
$$
There remain two pieces of degree $2$ each that are simply ramified
over the Kummer surface only.
\end{corollary}
\begin{proof}
The picture is local and occurs at a general point of the trope conic
consisting of a one-dimensional family of tacnodes, so we can use the
third coordinate $z$ in both coordinate systems. The trope conic is
given by $x=y=0$ in this picture.  Transverse to the $z$ direction the
Kummer surface is given by $y=0$ and the trope plane is given by
$y=x^2$. Our coordinates are, of course, not those of $X=\pp^3$.

The local fundamental group of this singularity is generated by loops
around these two pieces, and from the ramification picture we know
that they act on the covering by a single transposition (for the loop
around $y=x^2$) and a product of four distinct transpositions (for the
loop around $y=0$). The fact that the moving part of the ramification
does not meet our general point of the trope conic means that there is
no other ramification and $Y$ is a covering of $X$ given by the action
of these two elements on a set of $8$ elements.

There are only two possibilities: either the single transposition is
equal to one of the four, or else it connects two of them. In the
first case, the resulting covering is singular, but we know that $Y$
is smooth, so we must be in the second case. In particular, the
covering breaks into two pieces consisting of ramification along $y=0$
only (for the two out of four transpositions that do not touch the
single one) and a piece of degree $4$ with two transpositions over
$y=0$ and one transposition over $y=x^2$.  Furthermore, any local
model that has this type of ramification and with smooth total space
must be isomorphic locally to our covering.

We can construct such a model in the following way: first take a
degree $2$ cover ramified along $y=0$ given by $y=v^2$, and look at
the inverse image in here of $x^2-x$. It breaks into two irreducible
components as $x^2-v^2= (x+v)(x-v)$. Choose one, let's say $(v-x)$,
and take a covering of degree $2$ ramified a long there, so $u^2 =
x-v$. We have our coordinate system $(u,v)$ with the equations as
stated. Let us verify the Jacobian:
$$
\begin{vmatrix}
\frac{\partial x}{\partial u} & \frac{\partial x}{\partial v}\\
\frac{\partial y}{\partial u} & \frac{\partial y}{\partial v}
\end{vmatrix} \ 
= \
\begin{vmatrix}
2u & 1\\
0 & 2v
\end{vmatrix}
= 4 uv
$$
so the ramification locus upstairs is the union of $u=0$ and $v=0$
and these project to $y=x^2$ and $y=0$ respectively.
\end{proof}

\

\subsection{Blowing up the tacnodes}
\label{bltac}

Cut with a general plane in $X_0=\pp^3$. Locally using the coordinates
above, we may assume that our plane is $z=0$ so it has local
coordinates $(x,y)$ and the degree $4$ piece of the covering is given
in coordinates as described in Corollary \ref{coords}. We are going to
use the description as a composition of two covers of degree $2$ to
describe this piece of $Y_0$.

In what follows, we'll change temporarily the meaning of the notations
$X,Y$: assume that we have cut $X_0$ down to the general plane that
will be denoted by $X$, and $Y_0$ is cut down to the corresponding
covering, denoted by $Y$ for simplicity of notations. We will also
focus on a single tacnode, noting that altogether there are $32$ of
them in the plane.

Blow up twice to resolve the tacnode. Let $\widehat{X}$ denote the
blown-up variety and let $\widehat{Y}$ denote the local degree $4$
piece of the blown-up covering, normalized, with map $\widehat{f} :
\widehat{Y}\rightarrow \widehat{X}$.  Let $D^{\rm Tro}$ and $D^{\rm
  Kum}$ be the trope and Kummer divisors in $X$, locally near our
point.

When we blow up the first time, the strict transforms of $D^{\rm Tro}$
and $D^{\rm Kum}$ form, together with the exceptional divisor, an
ordinary triple point. Then blow that up again to get
$\widehat{X}$. Let $A$ be the strict transform in $\widehat{X}$ of the
first exceptional divisor, and let $B$ be the exceptional divisor of
the second blow-up. Let $T$ and $K$ denote the strict transforms of
the original divisors $D^{\rm Tro}$ and $D^{\rm Kum}$
respectively. Let $\alpha : \widehat{X}\rightarrow X$ be the map. Then
$$
\alpha ^{*}(D^{\rm Tro} ) = T + A + 2B
$$
and 
$$
\alpha ^{*}(D^{\rm Kum} ) = K + A + 2B.
$$
We can factorize our covering $\widehat{f}$ as a composition
$$
\widehat{Y} \stackrel{\mu}{\longrightarrow} \widehat{Z}  
\stackrel{\nu}{\longrightarrow} \widehat{X}
$$
where $\nu$ is a normalized double cover ramified over $\alpha
^{*}(D^{\rm Kum} )$ and $\mu$ is a double cover ramified over one half
of the pullback divisor $\nu ^{*} \alpha ^{*}(D^{\rm Tro} )$ as will
be described below.

One description of $\widehat{Z}$ is to say that we take the smooth
double cover $Z$ of $X$ ramified over $D^{\rm Kum}$, then
$\widehat{Z}$ is the normalization of $\widehat{X}\times _X Z$.

Another viewpoint is that
$\nu$ is a double cover ramified over $K$ and $A$, 
because normalization removes the double cover over $2B$. 
Since $K$ and $A$ are disjoint, the covering space 
$\widehat{Z}$ is smooth. 

The inverse image of $B$ in $\widehat{Z}$ is a double covering 
$B_Z \rightarrow B$ ramified over $A\cap B$ and $K\cap B$. 

The divisor $\alpha ^{*}(D^{\rm Tro} )$ pulls back in $\widehat{Z}$ to
a divisor that has multiplicity $2$ along $B_Z$, has two disjoint
pieces that compose $\nu ^{*}( T )$, and includes $2A_Z$ where $A_Z$
is the reduced inverse image of $A$ in $\widehat{Z}$. Taking half of
this will be one of the two pieces over $T$ plus $A_Z$ plus $B_Z$. The
double covering $\mu$ therefore has simple ramification over $A_Z$,
$B_Z$, and the one piece of $\nu ^{*}( T )$ that we have chosen.

Let $A'$, $B'$ and $T'$ denote the reduced inverse images of these
pieces in $\widehat{Y}$.  Let $K'$ be the reduced inverse image of $K$
in $\widehat{Y}$. It is a double covering of the ramification locus
$K_Z$ of $\widehat{Z}/\widehat{X}$ that lies over $K$.

We see that $\widehat{Y}$ has ordinary double points 
at $A'\cap B'$ and $T'\cap B'$. 

The map $\widehat{f} : \widehat{Y}\rightarrow \widehat{X}$ has the
following types of ramification: a cyclic covering of order $4$ along
$A'$; simple ramification along $B'$ that maps to $B$ by a generically
$2$-sheeted covering; simple ramification along $T'$ that maps
isomorphically to $T$, with the other part $T''$ of the inverse image
of $T$ being a generically $2$-sheeted covering; and simple
ramification along $K'$ that maps by a generically $2$-sheeted
covering to $K$. We may write altogether
$$
\begin{aligned}
  \widehat{f}^{*} (A) & = 4A' \\
\widehat{f}^{*} (B) & = 2B' \\
\widehat{f}^{*} (T) & = 2T' + T'' \\
\widehat{f}^{*} (K) & = 2K'.
\end{aligned}
$$
We can construct a parabolic Higgs bundle over $\widehat{X}$ in the following
way: choose a line bundle $\widehat{\Ll}$ on $\widehat{Y}$, take its
direct image, then use the natural filtrations over $A,B,T,K$ to put a
parabolic structure with some parabolic weights. The pullback of the
tautological $1$-form on $Y$ is a $1$-form on the smooth locus of
$\widehat{Y}$ and this will lead to a logarithmic Higgs field on
$\widehat{f}_{*}(\widehat{\Ll})$ away from $A\cap B$ and $T \cap B$.  As
$\widehat{f}_{*}(\widehat{\Ll} )$ is a bundle over $Y$, this logarithmic
Higgs field extends to a logarithmic Higgs field over all of
$Y$. Furthermore, it will respect the filtrations if we choose them to
be compatible with the covering.

On the other hand, if we let $\Ll$ denote the corresponding line
bundle on $Y$ (reflexive extension of the line bundle restricted to
the complement of the points over the tacnodes) then $f_{*}(\Ll )$ is
a bundle on $X$, and it has a meromorphic Higgs field that is
logarithmic along the smooth points of the divisor $D= D^{\rm Kum} +
D^{\rm Tro}$.  As we will see below, the Bogomolov-Gieseker inequality
indicates that this bundle with Higgs field cannot be considered as
being logarithmic over the tacnodes, since its pullback to
$\widehat{X}$ will violate the Bogomolov-Gieseker inequality. However,
it is the bundle $f_{*}(\Ll )$ for which we can obtain a calculation
of the Chern classes using the Grothendieck-Riemann-Roch formula.

We would therefore like to compare the Chern classes of our parabolic
structure on $\widehat{f}_{*}(\widehat{\Ll} )$ with those of
$f_{*}(\Ll )$ or rather the pullback $\alpha ^{*}(f_{*}(\Ll ))$.

The comparison between these two things is a local question. To see
this, let us assume that the parabolic weights are rational, as will
be the case for the structure to be used.  Then, there is a root stack
$\rho : \widetilde{X}\rightarrow \widehat{X}$ such that the parabolic
bundle may be viewed as a vector bundle on the root stack.  Let us call
this bundle $\Vv$. The root stack structure occurs over $A+B$. On the
other hand, let
$$
V:= \rho ^{*} \left( \alpha ^{*}(f_{*}(\Ll )) \right) 
$$ be the bundle pulled back from a vector bundle on $X$. Let $\jmath
: X^{\circ}\hookrightarrow X$ be the complement of the tacnode, and we
have $\widetilde{\jmath} : X^{\circ} \hookrightarrow
\widetilde{X}$. Let $V^{\circ}:= \widetilde{\jmath}^{*}(V)$, so we are
given a natural isomorphism $V^{\circ} \cong
\widetilde{\jmath}^{*}(\Vv )$.  Thus $f_{*}(\Ll ) =
\jmath_{*}(V^{\circ})$. In particular there are inclusions of quasicoherent
sheaves on $\widetilde{X}$
$$
V \hookrightarrow \widetilde{\jmath}_{*}(V^{\circ}) \hookleftarrow \Vv .
$$
These are isomorphisms away from (the inverse images in
$\widetilde{X}$ of) the divisors $A$ and $B$.  The images are both
contained in a coherent sheaf of the form, say, $V(nA + nB)$ so we may
use these to compare the Chern characters, namely we have coherent
sheaves $V(nA + nB) / V$ and $V(nA + nB) / \Vv$ supported on $A+B$ and
$$
{\rm ch}(\Vv) - {\rm ch}(V ) =
{\rm ch}(V(nA + nB) / V ) - {\rm ch}(V(nA + nB) / \Vv ) .
$$
The Chern characters of $V(nA + nB) / V$ and $V(nA + nB) / \Vv$
only depend on the local picture at the tacnode. 

Because of this locality, we may do the calculation assuming that $\Ll
= \Oo _Y$ is the trivial line bundle, then multiply by $32$ since
there are $32$ tacnodes, to get the global difference of Chern
characters, and add this to the Chern character of $f_{*}(\Ll )$
for the chosen global line bundle $\Ll$ on $Y$, to get the Chern
character of the parabolic logarithmic extension.

To further simplify the exposition, we are now going to just state
what are the good parabolic weights to use. These were found by doing
some computations and then solving the optimization problem in a crude
way using a computer, and in fact that was done in a couple of stages:
the first time, we noticed that the parabolic weights would involve
multiples of $1/4$. When we picked up the question some time later,
instead of re-doing this optimization we just did a grid search over
all the possible weights that were multiples of $1/4$. At the present
time, instead of exposing this (the details of the computation would
be difficult to explain, and the technique was not optimal) we will
just state how to get the resulting parabolic structure and then check
that it satisfies the parabolic Chern class vanishing conditions. This
is made easier by the description of the covering in two stages.

First calculate the vector bundle $\Uu := \widehat{f}_{*} (\Oo
_{\widehat{Y}})$.  The map $\mu : \widehat{Y}\rightarrow \widehat{Z}$
is a double covering in particular it is Galois, so
$$
\mu _{*} (\Oo _{\widehat{Y}}) = \Ww ^+ \oplus \Ww ^-
$$
where $\Ww ^+$ and $\Ww ^-$ are two line bundles on $\widehat{Z}$ where the
involution of the covering $\mu$ acts by $+1$ and $-1$. 
Notice that $\Ww ^+ = \Oo _{\widehat{Z}}$. 
Taking the direct image by $\nu$ of this decomposition gives
$$
\Uu = \Uu ^+ \oplus \Uu ^-
$$
where $\Uu ^+ = \nu _{*}(\Ww ^+)$ and  $\Uu ^- = \nu _{*}(\Ww ^-)$.
These are two bundles of rank $2$.

The variety $\widehat{Z}$ is obtained from $Z$ (the double cover of
$X$ ramified along $D^{\rm Kum}$) by blowing up twice, with the first
time generating an exceptional divisor whose strict transform is $B_Z$
and the second time generating an exceptional divisor $A_Z$; the other
strict transforms are $K_Z$ and $T_Z$. Note that $B_Z = \nu ^{*}(B)$
whereas $\nu ^{*}(A) = 2A_Z$ and $\nu ^{*}(K) = 2K_Z$ while $\nu
^{*}(T)$ is $T_Z$ plus another piece. The self intersections are
$A_Z^2 = -1$ and $B_Z^2=-2$.

The decomposition of $\hat{\nu}^{*} (T)$ into two pieces comes from a
decomposition of the inverse image of $T$ in the double cover
$Z\rightarrow X$, with one of the two pieces being the image in $Z$ of
our chosen $T_Z\subset \widehat{Z}$.  In particular, this divisor is
principal in the neighborhood in $Z$. Its inverse image in the blow-up
$\widehat{Z}$ is $A_Z + B_Z + T_Z$, which is therefore also principal.

In taking the double cover ramified over a principal divisor, it means
in our situation to use the trivial bundle as a square root. Thus, the
covering $\widehat{Y}$ is a double cover of $\widehat{Z}$ defined by
using the trivial bundle as the square root of the ramification
divisor $A_Z + B_Z + T_Z$. This divisor gives a line bundle whose
restriction to $B_Z$ is trivial as is its restriction to $A_Z$. It
follows that
$$
\Ww ^- \cong \Oo _{\widehat{Z}} 
$$
too. This gives to $\Ww ^-$ an equivariant structure for the double
covering below.

Now $\widehat{Z}$ is a double cover of $\widehat{X}$ branched over the
divisor $A+T$. It came from normalizing the inverse image of the
double cover branched over $D^{\rm Tro}$ and the pullback of that
divisor in $\widehat{X}$ is $A+2B + T$.  This is principal and its
square-root used for the covering is the trivial bundle. Taking the
normalization has the effect of declaring that the square-root of the
$2B$ term is $B$.  We have
$$
\Oo _{\widehat{X}} = \sqrt{\Oo _{\widehat{X}}(A+2B+T)} =
\Oo_{\widehat{X}}(B) \otimes \sqrt{\Oo _{\widehat{X}}(A+T)}
$$
so $\sqrt{\Oo _{\widehat{X}}(A+T)} = \Oo _{\widehat{X}}(-B)$. That is
to say, the covering $\widehat{Z}$ is defined by using $\Oo
_{\widetilde{X}}(-B)$ as square-root of $\Oo
_{\widetilde{X}}(A+T)$. Checking, $A$ has self-intersection $-2$ so
$(-B).A = (A+T).A / 2$.  Thus, the formula for the structure sheaf of
the double cover says
$$
\nu _{*}(\Oo _{\widehat{Z}}) \cong \Oo _{\widehat{X}} \oplus \Oo
_{\widehat{X}}(B).
$$
This gives for both pieces 
$$
\nu _{*} (\Ww ^+) = 
\nu _{*} (\Oo _{\widehat{Z}}) \cong \Oo _{\widehat{X}} \oplus \Oo _{\widehat{X}}(B)
$$
and
$$
\nu _{*}(\Ww ^-) \cong \nu _{*} (\Oo _{\widehat{Z}}) \cong 
\Oo _{\widehat{X}}\oplus \Oo _{\widehat{X}}(B).
$$
The subsheaves generated by global sections are given by the subsheaf
$\Oo _{\widehat{X}} \subset \Oo _{\widehat{X}}(B)$ in the second
factor in each case.

We can therefore write 
$$
V = V_1 \oplus V_2 \oplus V_3 \oplus V_4
$$
as a direct sum of trivial line bundles, with 
$$
\Uu ^+ = V_1 \oplus V_3(B), \qquad
\Uu ^- = V_2 \oplus V_4(B), 
$$ 
so
$$
\nu _{*} (\Oo _{\widehat{Y}})
= V_1 \oplus V_2 \oplus V_3(B) \oplus V_4(B).
$$

We now consider the constraints on the parabolic structure. Over $B$,
the $\Uu ^-$ piece constitutes the natural subspace of the
filtration. The parabolic weights of any piece of $\Uu ^-$ should be
$\leq$ the weights of any piece of $\Uu ^+$.

Near $A$, let $b$ be a local coordinate defining $A$ in $X$ (this is in
keeping with the notations we'll have later: $A$ is the $a$-axis defined by
$b=0$). The local coordinate
defining $A_Z$ is $b^{1/2}$, and the local coordinate defining $A'$ is 
$b^{1/4}$. This is the square-root of the coordinate defining $A_Z$. 
The functions $\Oo _{\widehat{Z}}$ include $1,b^{1/2}, \ldots $. Now
$\Ww ^+$ is generated over $\Oo _{\widehat{Z}}$ by $1$ so it 
includes the functions $1$ and $b^{1/2}$. On the other hand,
$\Ww ^-$ is generated over $\Oo _{\widehat{Z}}$ by $b^{1/4}$ so it 
includes the functions $b^{1/4}$ and $b^{3/4}$.

The subspaces $\Oo _{\widehat{X}}(B)$ correspond to the functions that
are multiples of $b^{1/2}$.

If we say that $V_1\oplus V_3(B)$ corresponds to $\Uu ^+$ and
$V_2\oplus V_4(B)$ corresponds to $\Uu ^-$ then $V_1$ is generated
over $\Oo _{\widehat{X}}$ by $1$, $V_3(B)$ is generated by $b^{1/2}$,
$V_2$ is generated by $b^{1/4}$ and $V_4(B)$ is generated by
$b^{3/4}$.

The parabolic filtration should therefore be adapted to our decomposition, in 
the order $1,2,3,4$ along $A$. 

We will now declare a precise collection of parabolic levels.  Our
parabolic structure is now a direct sum of four parabolic line
bundles. Recall that when we have a bundle $\Oo$ with filtration level
that is placed at parabolic level $-c\in (-1,0]$ at a divisor ${\rm
    div}$ it corresponds to a line bundle that is written as $\Oo
  (c\cdot {\rm div})$.  Let's denote the parabolic line bundles as
  $P_1,P_2,P_3,P_4$.

Define
$$
\begin{aligned}
P_1 & = \Oo _{\widehat{X}}(- \frac{1}{4}A -  \frac{1}{2}B ) \\
P_2 & = \Oo _{\widehat{X}} \\
P_3 & = \Oo _{\widehat{X}} \\
P_4 & = \Oo _{\widehat{X}}(\frac{1}{4}A + \frac{1}{2}B).
\end{aligned}
$$
Over $A$ the parts that
have been added are in the correct order. 

Over $B$, this is obtained from $V_1 \oplus V_3 (B)$ and $V_2 \oplus
V_4(B)$ by adding $-B/2$ to $V_1$, adding $-B$ to $V_3(B)$ (which
really means to make an elementary transformation), and adding $-B/2$
to $V_4(B)$. Thus on the $\Uu ^+$ pieces we added $-B/2$ and $-B$
whereas on the $\Uu ^-$ pieces we added $0$ and $-B/2$. This satisfies
the required criterion. The decomposition into a sum of line bundles
is compatible with the Higgs field in the tangential directions,
because of the blowing up: tangential vector fields along $B$ map to
zero in the tangent bundle of $Y$. This is the short reason why this
parabolic structure is allowable.  In subsection \ref{coordcalcs}
below, we will do a calculation in local coordinates to make sure that
the Higgs field is logarithmic with respect to this structure.

Technically the definition of the local contribution to ${\rm ch}_2$
is as was described previously, but it may now be expressed as
$$
\begin{aligned}
{\rm ch}_2(P_1& \oplus P_2 \oplus P_3 \oplus P_4)   = \\[+0.5pc]
& = 
\frac{
\left(
-\frac{1}{4}A - \frac{1}{2}B
\right) ^2}{2} 
+ 
\frac{
\left(
\frac{1}{4}A + \frac{1}{2}B
\right) ^2}{2} \\[+0.5pc]
& = \frac{1}{16} 
(A + 2B)^2   \\[+0.5pc]
& =
\frac{1}{16}  (A^2 + 4 AB +  4 B^2).
\end{aligned}
$$
We have $A^2 = -2$ and $B^2=-1$ with $AB=1$. So our contribution is
$$
(-2 + 4 - 4)/ 16 = -1/8. 
$$
The local contribution is therefore $-1/8$. When we multiply by the
$32$ tacnodes in a plane section, we obtain a global adjustment of
$-4$ to ${\rm ch}_2$. In view of the calculations in Corollary
\ref{extremalzerocor}, this adjustment leads to a parabolic Higgs
bundle with vanishing first and second Chern characters.

We need to check that the local parabolic structure we have been
considering in this section is one for which the Higgs field
becomes logarithmic.

\subsection{Calculations in coordinates}
\label{coordcalcs}

The easiest way to make sure that the parabolic structure we are
defining will induce a logarithmic property of the Higgs field, is to
write things out in local coordinates.

Recall that $X$ denotes a slice by a plane in $X_0=\pp^3$, and
localize at a tacnode point of the wobbly locus. Then, we will look near
a point in the covering $Y/X$ that is the center of the piece of
degree $4$ over the neighborhood in $X$. Fix coordinates $x,y$ for a
small ball around the point in $X$, such that the Kummer is given by
$y=0$ and the trope is given by $x^2-y=0$.  Let $\cc \{ x, y \}$ be
the coordinate ring of convergent series in $x$ and $y$.

Express the neighborhood in $Y$ as a composition of two coverings of
degree $2$.  The first has coordinates $x,v$ with $y=v^2$, so it has
coordinate ring $\cc \{ x,v\}$.  The inverse image of the tacnode
divisor in here has equation $x^2 - v^2=0$ that decomposes as
$(x+v)(x-v)$ so the pullback divisor decomposes into two irreducible
components. Choose one of these, say $x-v=0$, as the branch locus for
the second covering.  Introduce the coordinate $u$ with
$u^2=x-v$. This gives the system of coordinates $u,v$ for the
neighborhood in $Y$, with ring $\cc \{ u,v\}$ and equations for the
degree $4$ map $f$ to the neighborhood in $X$ are:
$$
x=u^2 + v, \qquad y=v^2.
$$
Calculation of the Jacobian matrix showed that the branch locus
downstairs is the union of the components $y=0$ and $x^2-y=0$. Set
$\Delta := x^2-y$.

The direct image $f_{*} \Oo $ is a rank $4$ vector bundle $V$ over
the neighborhood in $X$, which may be viewed as given by the module
$\cc \{ u,v\}$ considered as a module of rank $4$ over $\cc \{
x,y\}$. Write the decomposition into a direct sum of line bundles
$V=V_1 \oplus V_2 \oplus V_3 \oplus V_4$ corresponding to the module
decomposition
$$
\cc \{ u,v\} = 1\cdot \cc \{ x,y\} \oplus u\cdot \cc \{ x,y\} 
\oplus v \cdot \cc \{ x,y\} \oplus uv \cdot \cc \{ x,y\} .
$$
Suppose our tautological form on $Y$ is written as $\lambda = \varphi
(u,v) du + \psi (u,v) dv$. Note that it satisfies a constraint, due to
the fact that the trope plane is the exceptional divisor for a blow-up
of the Prym, and the form comes from a form on the Prym. This says
that $\psi$ is a multiple of $u$, although that condition does not seem
to be needed later since the $dv$ term does not pose a problem.

The equation $v^2=y$ tells us that
$$
dv = vdy / 2y .
$$
The equation $u^2 = x-v$ tells us that
$$
du = u^{-1}(dx -dv)/2
$$
where $dv$ may be written in terms of $dy$. 

The actions of multiplication by $u$ and $v$ my be expressed in the
form of $4\times 4$ matrices acting on the direct sum
decomposition. Use as basis vectors $1,u,v,uv$.  We have
$$
u = 
\left( \begin{array}{cccc}
0 & x  & 0 & -y \\
1 & 0  & 0 & 0  \\
0 & -1 & 0 & x  \\
0 & 0  & 1 & 0  
\end{array} \right), \qquad \qquad 
v = 
\left( \begin{array}{cccc}
0 & 0  & y & 0 \\
0 & 0  & 0 & y  \\
1 & 0  & 0 & 0  \\
0 & 1  & 0 & 0  
\end{array} \right) .
$$
One calculates:
$$
u ^{-1} = 
\frac{1}{\Delta} 
\left( \begin{array}{cccc}
0 & \Delta  & 0  & 0         \\
x & 0         & y  & 0         \\
0 & 0         & 0  & \Delta  \\
1 & 0         & x  & 0  
\end{array} \right) 
$$
as may be verified by multiplying together with the matrix for $u$. 
Also, $v^{-1} = v/y$ so
$$
v ^{-1}= 
\left( \begin{array}{cccc}
0    & 0     & 1 & 0 \\
0    & 0     & 0 & 1  \\
1/y & 0     & 0 & 0  \\
0    & 1/y  & 0 & 0  
\end{array} \right) .
$$
The differentials $du$ and $dv$ may now be expressed as matrices 
with entries in terms of $dx$ and $dy$: 
$$
dv = (v/2y)dy = 
\left( \begin{array}{cccc}
0       & 0        & dy / 2 & 0        \\
0       & 0        & 0        & dy / 2  \\
dy/2y & 0        & 0       & 0         \\
0       & dy/2y  & 0       & 0  
\end{array} \right) ,
$$
and
$$
du = u^{-1}(dx -dv)/2
=
\frac{1}{2\Delta} 
\left( \begin{array}{cccc}
0 & \Delta  & 0  & 0         \\
x & 0         & y  & 0         \\
0 & 0         & 0  & \Delta  \\
1 & 0         & x  & 0  
\end{array} \right) 
\cdot 
\left( \begin{array}{cccc}
dx        & 0          & -dy/ 2  & 0        \\
0          & dx        & 0         & -dy/ 2  \\
-dy/2y  & 0          & dx       & 0         \\
0          & -dy/2y  & 0         & dx  
\end{array} \right) 
$$
$$
=
\frac{1}{4\Delta} 
\left( \begin{array}{c:c:c:c}
0                   &  2 \Delta dx    & 0                & -\Delta dy    \\   \hdashline
2xdx -  dy      &  0                  & -xdy + 2ydx &    0             \\   \hdashline
0                   & -\Delta dy / y &  0                &  2\Delta dx   \\    \hdashline
2dx -x dy / y  & 0                   & 2xdx - dy     &  0        
\end{array} \right) .
$$
The terms $2xdx -  dy$ are equal to $d\Delta$. We also note that
$y = x^2 - \Delta$ so 
$ dy = 2x dx - d\Delta$ and 
$$
-xdy + 2y dx = (2y -2x^2)  dx + xd\Delta = xd\Delta - 2\Delta dx,
$$
hence
$$
2dx -x dy / y = (x/y)d\Delta - (2/y)\Delta dx .
$$
We may therefore write
$$
du = 
\frac{1}{4} 
\left( \begin{array}{c:c:c:c}
0 &  0  &  0 &   0   \\   \hdashline
0 &  0  &  0 &   0   \\   \hdashline
0 & 0   &  0 &   0   \\    \hdashline
(x/y)\frac{d\Delta}{\Delta} -2dx/ y  & 0   &  0 &   0        
\end{array} \right) 
$$
$$
+
\frac{1}{4} 
\left( \begin{array}{c:c:c:c}
0                                 &  0  & 0                                   & 0    \\   \hdashline
\frac{d\Delta}{\Delta} &  0   & x\frac{d\Delta}{\Delta} &  0    \\   \hdashline
0                                 & 0    &  0                                 &  0    \\    \hdashline
0                                 & 0    & \frac{d\Delta}{\Delta}   &  0        
\end{array} \right) 
+
\frac{1}{4} 
\left( \begin{array}{c:c:c:c}
0   &  2 dx    & 0       & - dy    \\   \hdashline
0   &  0         & -2dx  &    0             \\   \hdashline
0   & - dy / y &  0      &  2 dx   \\    \hdashline
0   & 0          & 0       &  0        
\end{array} \right) .
$$
The second and third matrices have terms that are either holomorphic
forms, or multiples of $d\Delta / \Delta$ and $dy/y$. When we pull
back to a blow-up, these will remain as (at worst) logarithmic
forms. The first matrix, with a single non-zero coefficient in the
lower left, is not logarithmic. This term is however logarithmic on
points of the divisors $y=0$ or $\Delta = 0$ away from the origin,
indeed that is clear from the expression here when $y\neq 0$ and it is
clear from the original expression $(1/4\Delta )(2dx -x dy /y)$ when
$\Delta \neq 0$.  This term will lead to a non-logarithmic term when
we blow up twice the tacnode, so the elementary transformations and
parabolic structures will need to take that into account.

Consider the filtration of $V$ that has three steps, with quotient
$V_1$, subquotient $V_2\oplus V_3$ in the middle, and subbundle
$V_4$. This filtration is preserved by the operators of multiplication
by $u$ or $v$, modulo the maximal ideal $(x,y)$.  Furthermore, it is
preserved by the residues, although not strictly because of the terms
$xd\Delta / \Delta$ in position $(2,3)$ and $-dy/y$ in position
$(3,2)$.
 
We can therefore use this filtration to put parabolic structures after
blowing up, as will be done in the next subsection.

\subsection{Pulling back to the blow-up and a root cover}

Blowing up the origin introduces the coordinate $a=y/x$ so $y=ax$ and
the coordinate chart on the blow up has coordinates $(x,a)$.

The main term in the lower left corner  of $du$ becomes
$$
(x/y)\frac{d\Delta}{\Delta} -2dx/ y =
\frac{1}{a} d\log \left( 
\frac{x-a}{x} \right) = \frac{-1}{x-a} d\log (ax) .
$$
Then blow up again using $b=x/a$ so $x=ab$ and $y=a^2b$. Our lower
left corner term becomes
$$
\frac{1}{a(b-1)}\left(
2\frac{da}{a} + \frac{db}{b} \right) .
$$
This is logarithmic along the divisor $A$ which is $b=0$ (the
$a$-axis) and has a pole of order $2$ along $B$ (the $b$-axis given by
$a=0$).

We therefore need to make an elementary transformation in order to get
a logarithmic pole along $B$. This corresponds to the normalization of
the spectral covering in the previous discussion.

In all, we are going to use the same decomposition into $4$ line
bundles, but putting parabolic levels on these pieces, to get a
decomposition into four parabolic line bundles of the form
$$
\Oo (-A/4 - B/2) \oplus \Oo \oplus \Oo \oplus \Oo (A/4 + B/2).
$$
We would like to calculate that this indeed gives a parabolic
logarithmic Higgs bundle when we use a Higgs field that is a
combination of $du$ and $dv$.

To do this calculation, let us pull back to a root covering. Since we
want the divisors $A/4$ and $B/2$, with $A$ given by $b=0$ and $B$
given by $a=0$, let's introduce the root coordinates
$$
\alpha = a^{1/2}, \;\;\; \beta = b^{1/4}.
$$
Thus 
$$
a=\alpha ^2, \;\; b= \beta ^4, \;\; x = \alpha ^2 \beta ^4, \;\;
y = \alpha ^4 \beta ^4,
$$
and we have
$$
dx = 2\alpha \beta ^4 d\alpha + 4 \alpha ^2 \beta ^3 d\beta , 
\;\;\; 
dy = 4\alpha ^3\beta ^4 d\alpha + 4 \alpha ^4 \beta ^3 d\beta 
$$
and
$$
\Delta = \alpha ^4 \beta ^8 - \alpha ^4 \beta ^4 = \alpha ^4
\beta^4 (\beta ^4 - 1).
$$
The lower left corner term is
$$
(* ) = 
\frac{-1}{x-a} d\log (y) = 
\frac{-4}{\alpha ^2 (\beta ^4 - 1)} 
\left( 
\frac{d\alpha}{\alpha} + \frac{d\beta}{\beta}
\right) .
$$
The full matrix for $du$ becomes:
$$
\frac{1}{4} 
\left( \begin{array}{c:c:c:c}
0   &  2 dx    & 0       & - dy    \\   \hdashline
d\Delta / \Delta    &  0         & -2dx + xd\Delta / \Delta  &    0             \\   \hdashline
0   & - dy / y &  0      &  2 dx   \\    \hdashline
(* )   & 0          & d\Delta / \Delta       &  0        
\end{array} \right)
$$
$$
= \frac{1}{4} 
\left( \begin{array}{c:c:c:c}
0   &  4\alpha \beta ^4 d\alpha + 8 \alpha ^2 \beta ^3 d\beta     
& 0       & - 4\alpha ^3\beta ^4 d\alpha + 4 \alpha ^4 \beta ^3 d\beta    \\   \hdashline
d\Delta / \Delta    &  0         & -4\alpha \beta ^4 d\alpha - 8 \alpha ^2 \beta ^3 d\beta   + 
\alpha ^2 \beta ^4d\Delta / \Delta  &    0             \\    \hdashline
0   & - d\log (\alpha ^4 \beta ^4)  &  0      &  
4\alpha \beta ^4 d\alpha + 8 \alpha ^2 \beta ^3 d\beta   \\    \hdashline
(* )   & 0          & d\Delta / \Delta       &  0        
\end{array} \right) .
$$
Instead of the basis $e_1,e_2,e_3,e_4$ (which was originally
$1,u,v,uv$), the new frame for the new bundle over the root covering
is $f_1 = \alpha ^{-1}\beta ^{-1} e_1, f_2 = e_2, f_3 = e_3,
f_4=\alpha \beta e_4$.  In this frame, the new matrix is obtained by
multiplying the top row and last column by $\alpha ^{-1}\beta^{-1}$,
and multiplying the last row and first column by $\alpha \beta$. This
yields the new matrix:
$$
du _{\bf f} = \frac{1}{4} 
\left( \begin{array}{c:c:c:c}
0   &  4\beta ^3 d\alpha + 8 \alpha  \beta ^2 d\beta     
& 0       & - 4\alpha \beta ^2 d\alpha + 4 \alpha ^2 \beta  d\beta    \\   \hdashline
\alpha \beta d\Delta / \Delta    &  0         & -4\alpha \beta ^4 d\alpha - 8 \alpha ^2 \beta ^3 d\beta   + 
\alpha ^2 \beta ^4d\Delta / \Delta  &    0             \\   \hdashline
0   & - d\log (\alpha ^4 \beta ^4)  &  0      &  
4\beta ^3 d\alpha + 8 \alpha  \beta ^2 d\beta   \\    \hdashline
\alpha ^2 \beta ^2 \cdot  (* )   & 0          & \alpha \beta  d\Delta / \Delta       &  0        
\end{array} \right) .
$$
Now,
$$
\alpha ^2 \beta ^2 \cdot  (* )  = 
\frac{-4\beta ^2}{(\beta ^4 - 1)} 
\left( 
\frac{d\alpha}{\alpha} + \frac{d\beta}{\beta}
\right) .
$$
This is logarithmic at general points of $\alpha = 0$ and $\beta = 0$.
Note that the locus $\beta ^4 = 1$ corresponds to a strict transform
of one of our original divisors and we have verified the logarithmic
property over those. The remaining terms of $du _{\bf f}$ are also
logarithmic or better.

Our Higgs field is obtained as a combination of $du$ and $dv$ times
$1,u,v,uv$ times functions of $x,y$. To complete the verification, we
need to express the matrices for $u$, $v$ and $dv$ in terms of the
coordinates $\alpha$ and $\beta$ and in the new frame. In the frame
$e_1,e_2,e_3,e_4$ we had
$$
u = 
\left( \begin{array}{cccc}
0 & x  & 0 & -y \\
1 & 0  & 0 & 0  \\
0 & -1 & 0 & x  \\
0 & 0  & 1 & 0  
\end{array} \right) 
=
\left( \begin{array}{cccc}
0 & \alpha ^2 \beta ^4  & 0 & -\alpha ^4 \beta ^4 \\
1 & 0  & 0 & 0  \\
0 & -1 & 0 & \alpha ^2 \beta ^4 \\
0 & 0  & 1 & 0  
\end{array} \right) 
$$
which becomes, in the new frame (multiplying the top and bottom rows and the
first and last columns as previously): 
$$
u _{\bf f}=
\left( \begin{array}{cccc}
0 & \alpha  \beta ^3  & 0 & -\alpha ^2 \beta ^2 \\
\alpha \beta & 0  & 0 & 0  \\
0 & -1 & 0 & \alpha  \beta ^3 \\
0 & 0  & \alpha \beta & 0  
\end{array} \right) .
$$
Similarly, in the original frame
$$
v = 
\left( \begin{array}{cccc}
0 & 0  & y & 0 \\
0 & 0  & 0 & y  \\
1 & 0  & 0 & 0  \\
0 & 1  & 0 & 0  
\end{array} \right) 
= 
\left( \begin{array}{cccc}
0 & 0  & \alpha ^4 \beta ^4 & 0 \\
0 & 0  & 0 & \alpha ^4 \beta ^4  \\
1 & 0  & 0 & 0  \\
0 & 1  & 0 & 0  
\end{array} \right) 
$$
so in the new frame 
$$
v_{\bf f} = 
\left( \begin{array}{cccc}
         0  & 0  & \alpha ^3 \beta ^3 & 0 \\
0 & 0  & 0 & \alpha ^3 \beta ^3  \\
\alpha \beta & 0  & 0 & 0  \\
0 & \alpha \beta  & 0 & 0  
\end{array} \right) .
$$
Also, in the previous frame 
$dv = (dy / 2y)\cdot v$
so in the new frame, 
$$
dv_{\bf f} = 2 \left( \frac{d\alpha}{\alpha} + \frac{d\beta}{\beta} \right)
\left( \begin{array}{cccc}
0 & 0  & \alpha ^3 \beta ^3 & 0 \\
0 & 0  & 0 & \alpha ^3 \beta ^3  \\
\alpha \beta & 0  & 0 & 0  \\
0 & \alpha \beta  & 0 & 0  
\end{array} \right) 
$$
$$
= 
\left( \begin{array}{c:c:c:c}
0 & 0  & \alpha ^2 \beta ^3 d\alpha + \alpha ^3 \beta ^2 d\beta   & 0
         \\ \hdashline
0 & 0  & 0 & \alpha ^2 \beta ^3 d\alpha + \alpha ^3 \beta ^2 d\beta
         \\ \hdashline
         \beta  d\alpha + \alpha d\beta & 0  & 0 & 0  \\ \hdashline
0 & \beta  d\alpha + \alpha d\beta   & 0 & 0  
\end{array} \right) .
$$
We see that any products of $1,u,v,uv$ times $du$ or $dv$ in the new
frame are logarithmic, so the Higgs field will induce a logarithmic
Higgs field on this bundle as desired. This completes the proof of
theorem to be stated in the concluding subsection below.

\subsection{Degree zero case---conclusion}

\begin{theorem}
\label{degzerothm}
Suppose $\Ll_0$ is a flat line bundle on the abelian
variety $\Prym $, and define the spectral line bundle
$\Ll = \blo_{0}^{*}(\Ll _0)\otimes 
\Oo _Y(\ExY +2\FxY )$ on $Y$. 
Put
$$
\Ee := f _{*}(\Ll )
$$
as a meromorphic Higgs bundle on $X=\pp^3$, then blow up twice at
the tacnodes in a general planar section, and put the parabolic
structure we have defined above so that the Higgs field becomes
logarithmic.  For this parabolic structure, ${\rm ch}_i^{\rm
  par} = 0$ for $i=1,2$. Therefore, this
parabolic Higgs bundle admits an extension to a purely imaginary
twistor $\srD$-module as in Theorem \ref{outsidecod2}. The associated
local system on $X_0 -\Wob_0$ has rank $8$. Its monodromy around the
trope planes in $\Wob_0$ is unipotent with a single Jordan block of
size $2$. Its monodromy around the Kummer surface in $\Wob _0$ is
unipotent consisting of a direct sum of $4$ Jordan blocks of size $2$.
\end{theorem}
\begin{proof}
In keeping with the result of Corollary \ref{extremalzerocor}, choose
the line bundle $\Ll$ to be anything numerically equivalent to $\Oo
_Y(\ExY +2\FxY )$, that is to say anything of the form \linebreak
$\Ll =
\blo_{0}^{*}(\Ll _0)\otimes \Oo _Y(\ExY +2\FxY )$ for $\Ll_0$ a flat line
bundle on the Hitchin fiber $\Prym$.

The tautological $1$-form on $\Prym$ pulls back to a $1$-form on $Y$,
which we may view as a meromorphic section of $f^{*}(\Omega
^1_X)$ on $Y$, yielding a meromorphic Higgs field $\Phi : \Ee
\rightarrow \Ee \otimes \Omega ^1_X$ satisfying the commutativity
condition $\Phi \wedge \Phi = 0$.  A local calculation at the
ramification points of $f: Y\rightarrow  X$ shows that
$\Phi$ has logarithmic singularities along smooth points of the branch
divisor of $f$. However, the fact that $\Prym$ may be viewed as a
subvariety of the cotangent bundle of $X$ over the very stable points,
shows that $\Phi$ has no poles along smooth points of the branch
divisor not on the wobbly locus. Thus, we obtain a logarithmic Higgs
bundle $(\Ee ^{\circ}, \Phi ^{\circ})$ over $X^{\circ}$ defined to be
the complement of the singular locus of $\Wob$.

Let $\zeta : \widehat{X}^+ \rightarrow X^+$ be obtained by blowing up
twice the trope conics, on the open subset $X^+\subset X$ complement
of the $16$ singular points of $\Kum$.

Our discussion in the transverse plane sections to the trope conics
shows how to create a parabolic bundle $\Ee ^+ _{\cdot}$ on
$(\widehat{X}^+, D^+)$ where $D^+$ is the reduced inverse image of the
wobbly divisor in $\widehat{X}^+$. The nontrivial parabolic structure
is concentrated on the parts of $D^+$ lying over the trope conics. The
parabolic structure extends the given bundle $\Ee ^{\circ}$ from the
open subset $X^{\circ} \subset \widehat{X}^+$. We have seen, in the
transverse plane sections, that the Higgs field $\Phi ^+$ (unique
extension of $\Phi ^{\circ}$) becomes logarithmic for this parabolic
structure. This verification was completed in the previous subsection.

Let $\Ee ^+ _{\rm raw}:= \zeta ^{*}\Ee |_{X^+}$. If $X_H\subset X$ is
a plane (general so that it misses the singular points of $\Kum$) it
intersects the trope conics in a total of $32$ points (two points for
each of the $16$ conics).  Its inverse image $\widehat{X}_H$ is
obtained by combining $32$ times the local picture we have seen above.
We took care to make sure that
$$
{\rm ch}_1^{\rm par} \left( 
\Ee ^+ _{\cdot}|_{\widehat{X}_H} \right) = 0.
$$
As pointed out in Subsection \ref{bltac}, the difference in second
Chern characters between $\Ee ^+ _{\rm raw}|_{\widehat{X}_H}$ and $\Ee
^+ _{\cdot}|_{\widehat{X}_H}$ is local, so it is $32$ times the
quantity $(-1/8)$ calculated in Subsection \ref{bltac}.  In other
words,
$$
{\rm ch}_2^{\rm par} \left( 
\Ee ^+ _{\cdot}|_{\widehat{X}_H} \right)
=
{\rm ch}_2\left( 
\Ee ^+ _{\rm raw}|_{\widehat{X}_H} \right)  + 32 \cdot (-1/8).
$$
On the other hand, from Corollary \ref{extremalzerocor}, for our
choice of line bundle $\Ll$ we have
$$
{\rm ch}_2\left( 
\Ee ^+ _{\rm raw}|_{\widehat{X}_H} \right) = 4.
$$
We conclude that ${\rm ch}_2^{\rm par} \left( 
\Ee ^+ _{\cdot}|_{\widehat{X}_H} \right) = 0$. 

We now have a parabolic logarithmic Higgs bundle $(\Ee ^+ _{\cdot} ,
\Phi ^+ )$ on $(\widehat{X}^+ , D^+)$ such that the first and second
parabolic Chern characters vanish on a plane section
$\widehat{X}_H$. The divisor $D^+$ has normal crossings.  Mochizuki's
theory \cite{Mochizuki-kh1, Mochizuki-kh2} implies that there exists a
tame purely imaginary harmonic bundle on $\widehat{X}^+ - D^+ \cong X
- \Wob$ whose corresponding parabolic Higgs bundle is this one. The
rank is $8$. The monodromy along smooth points of $\Wob$ is given by
the residue of the Higgs field, which in turn corresponds to the
ramification of $Y/X$ over the different components of the wobbly
locus, giving the stated properties.
\end{proof}

\begin{remark}
\label{localmon}
Our technique of construction yields more precise information about
the monodromy of the local system along the exceptional divisors in
$\widehat{X}$ lying over the trope conics. For example, the
eigenvalues around the $A$ divisors are $4$-th roots of unity and the
eigenvalues around the $B$-divisors are $\pm 1$. There is also a
unipotent piece around the $A$ divisor. It is left to the reader to
make a more explicit statement.
\end{remark}

\section{Hecke operators}
\label{chapter-hecke}

\subsection{Hecke transformations in terms of bundles}
\label{bundle-view}

In Chapter \ref{synthetic} we encountered the 
geometric picture of the Hecke correspondences in the context of
pencils of quadrics in $\pp^5$. Let us review the Hecke transformations on bundles. 
 
If $E$ is a rank $2$ vector bundle and $t\in C$ is a point, the
\emph{\bfseries Hecke line} of $E$ at $t$ is the projectivization $\pp
E_t$ parametrizing rank $1$ quotients $E_t \twoheadrightarrow \cc$. A
vector space quotient corresponds to a surjection of coherent sheaves
$E\twoheadrightarrow \cc _t$ where $\cc _t$ denotes the skyscraper
sheaf at $t$.  Let $E'$ denote the kernel of this map, usually known
as the \emph{\bfseries down-Hecke transform} of $E$ centered at $E_t
\twoheadrightarrow \cc$.  It is a torsion-free coherent sheaf, hence
locally free and therefore it is itself a bundle. We have
$$
\det(E') = \det(E)\otimes \Oo _C(-t).
$$
In order to get a map between points of our moduli spaces we need
to correct the determinant by tensoring with a line bundle. Thus,
consider a point $(A,t)\in \Cbar$, meaning that $A$ is a line bundle
with $A^{\otimes 2} = \Oo _C(t-\pw )$. Now, if $E$ is a vector bundle
with determinant $\det(E)\cong \Oo _C(\pw )$, and $E_{t}
\twoheadrightarrow \mathbb{C}$, we can use the associated down-Hecke 
$E' = \ker(E \to \cc_t)$ with determinant $\Oo_{C}(\pw -t)$
to form the transformed bundle
$$
E' \otimes A.
$$
Note that by definition
$$
\det(E'\otimes A) = \det(E') \otimes A^{\otimes 2}
= \det(E) \otimes \Oo _C(-t) \otimes \Oo _C(t-\pw ) = \Oo _C.
$$
Thus, if $E$ was stable, i.e. representing a point of the moduli
space $X_1$, then the Hecke transform $E' \otimes A$ is a bundle (that
one may verify is semistable) with trivial determinant so it
corresponds to a point of $X_0$. We obtain the Hecke $\pp^{1}$
parametrizing the trivial determinant down-Hecke transforms of $E$ at $t$.
It is always a line in $X_0\cong
\pp^3$ as we will see in Theorem \ref{linefibers}.

In the other direction, suppose $E$ is a stable bundle with trivial
determinant. We need to transform the ``down-Hecke'' into an
``up-Hecke'' to get a point of $X_1$ with determinant $\Oo _C(\pw )$.
Again, this will depend on the choice of a point $(A,t) \in \Cbar$ and
a choice of a quotient $E_{t} \twoheadrightarrow \cc$. Once these choices
are made we define the transform to be
$$
E' \otimes A (\pw ),
$$
where again $E' = \ker(E \to \cc_t)$.
By definition $\det(E') = \Oo_{C}(-t)$ and so 
$$
\det(E'\otimes A(\pw )) = \det(E') \otimes A^{\otimes 2}
(2\pw )
=
\det(E) \otimes \Oo _C(-t) \otimes \Oo _C(t-\pw )\otimes \Oo _C(2\pw )
= \Oo _C(\pw ).
$$
It is straightforward again to verify that $E' \otimes A (\pw )$ is
also stable so it is a point of $X_1$. The image of the Hecke line is
a smooth rational curve in $X_1$, which in fact is a smooth plane
conic as we will see in Theorem \ref{conicfibers}.  These are the
Hecke curves over points of $X_0-\Kum$.

If we start with a point of the Kummer surface, there are several
choices of semistable bundle corresponding to that $S$-equivalence
class. The above formulas, applied to these different bundles, yield
the Hecke fiber. For smooth points of $\Kum$ it is the union of two
lines corresponding to the two semistable bundles; for nodes of $\Kum$
it is a single line counted twice.  See Theorem \ref{conicfibers}.

In the next section we will put the Hecke correspondences together
in a family over the base $\Cbar$. The fiber of $\bigHeckebar$ over  $a=(A,t)\in \Cbar$ will be
$$
\bigHeckebar(a) := \left\{ (E,\ekk) \;\; \left| \;
\begin{minipage}[c]{3.5in}
  $E\in X_1$ and $\ekk : E_t
  \twoheadrightarrow \eKK$
  is a one dimensional quotient of the fiber of $E$ at $t$.
\end{minipage}
\right.\right\} .
$$
We have maps
$$
X_1 \stackrel{p}{\longleftarrow} \bigHeckebar(a) \stackrel{q}{\longrightarrow} X_0
$$
where $p$ is just the projection. It is a $\pp^1$-bundle whose fiber over $E\in
X_1$ is the projective line of rank $1$ quotients of $E_t$.
Existence of the moduli space and the fibration property of the
projection are due to the fact that all points of $X_1$ are stable.

For a point $(E,\ekk)\in \bigHeckebar (a)$ the corresponding
point $q(E')\in X_0$ is the $S$-equivalence class of 
the down-Hecke of $E$ at the point $t$, defined by the quotient $\ekk$ and normalized using 
$A$ as described above.

\subsection{The big Hecke correspondences}
\label{sec-bigHecke}

Putting together these fibers yields the big Hecke
correspondences
\begin{equation}  
\quad   \xymatrix@M+0.5pc@-0.5pc{
    & \bigHeckebar \ar[dl]_-{p} \ar[dr]^-{q} & \\
    X_{1} & & X_{0}\times \Cbar
} \qquad
\xymatrix@M+0.5pc@-0.5pc{
    & \bigHeckebar \ar[dl]_-{\pzo} \ar[dr]^-{\qzo} & \\
    X_{0} & & X_{1}\times \Cbar
    }
\end{equation}
whose notation was introduced in Subsection \ref{subsec-Hecke}. 

These Hecke correspondences are the ones that most conveniently encode
the eigensheaf property for the de Rham data ($\srD$-modules) or
the Dolbeault data (parabolic Higgs complexes) on the moduli of
$\mathbb{P}SL_{2}(\mathbb{C})$-bundles on $C$. These are the de Rham
and Dolbeault objects that under the Langlands correspondence should
correspond to flat $SL_{2}(\mathbb{C})$-bundles on $C$ or to
a semistable Higgs $SL_{2}(\mathbb{C})$-bundle on $C$ respectively.

To spell this out, note that the moduli of
$\mathbb{P}SL_{2}(\mathbb{C})$-bundles on $C$ is a disconnected
Deligne-Mumford stack
\[
\mathsf{M}_{0} \sqcup \mathsf{M}_{1}, \quad \mathsf{M}_{i} =
\left[X_{i}/\mathsf{J}[2]\right], \ \text{for} \ i =0,1.
\]
Here $\mathsf{J}[2]$ is the group of $2$-torsion points in the abelian
surface $\mathsf{J} := \op{Jac}^{0}(C)$, and $\mathfrak{a} \in
\mathsf{J}[2]$ acts on the moduli space $X_{0}\sqcup X_{1}$ of bundles
with fixed determinant by $E \mapsto E\otimes\mathfrak{a}$.

The big $\mathsf{M}_{1}$ to $\mathsf{M}_{0}$ Hecke correspondence is a
correspondence
\[
\xymatrix@M+0.5pc@-0.5pc{
  & \mathbf{Hecke} \ar[dl] \ar[dr] & \\
\mathsf{M}_{1} & & \mathsf{M}_{0}\times C
}
\]
It is described explicitly in terms of a correspondence
\[
\xymatrix@M+0.5pc@-0.5pc{
  & \bigHecke \ar[dl] \ar[dr] & \\
X_{1} & & \mycal{X}_{0}},
\]
Here
\begin{itemize}
\item $\mycal{X}_{0}$ is the moduli space of semistable rank two
  bundles with determinant of the form $\mathcal{O}_{C}(\pw-t)$ for
  some $t \in C$. It fibers $\mycal{X}_{0} \to C$ over $C$, the fiber
  over a given $t \in C$ being the moduli $\mycal{X}_{0}(t)$ is the
  moduli of bundles with determinant equal to the specific line bundle
  $\mathcal{O}_{C}(\pw-t)$.
\item $\bigHecke$ is the moduli space of triples
  \[
  \bigHecke = \left\{ (E,E',\beta) \ \left| \
  \text{
    \begin{minipage}[b]{2.2in}
      \[
      \begin{aligned}
        E & \in X_{1}, \
        E'  \in \mycal{X}_{0}(t) \\
        \beta : & E' \hookrightarrow E, \ 
        \op{supp}(\op{coker}(\beta)) = t
      \end{aligned} 
      \]
    \end{minipage}
  }
  \ \right.\right\}.
  \]
\item the South-West map is given by by $(E,E',\beta) \mapsto E$ and
  the South-East map is given by $(E,E',\beta) \mapsto E'$.
\end{itemize}

\

\noindent
In fact $\mycal{X}_{0} \to C$ is an algebraic  $\mathbb{P}^{3}$
bundle, equipped with a flat (Heisenberg-like) connection with monodromy
$\mathsf{J}[2]$. The group $\mathsf{J}[2]$ also acts on
$\mycal{X}_{0}$ by tensorization, and  the action preserves the flat
connection. The quotient
\[
[\mycal{X}_{0}/\mathsf{J}[2]] \to C
\]
inherits a flat structure with trivial monoromy and taking $\pw \in C$
as a base point we get an algebraic isomorphism
\[
\xymatrix@M+0.5pc@-0.5pc{
  [\mycal{X}_{0}/\mathsf{J}[2]] \ar[r] \ar[d]_-{\cong} & C \ar@{=}[d]  \\
  \mathsf{M}_{0}\times C \ar[r]_-{\pr_{C}} & C.
}
\]
The group $\mathsf{J}[2]$ also acts on $\bigHecke$ by tensoring with
$2$-torsion line bundles:
\[
  \mathsf{J}[2]\times \bigHecke \to \bigHecke, \quad (\mathfrak{a},
  (E,E',\beta)) \mapsto
  (E\otimes\mathfrak{a},E'\otimes\mathfrak{a},
  \beta\otimes\op{id}_{\mathfrak{a}}),
\]
and passing to quotients we get
\[
  \left[
    \begin{minipage}[c]{2.9in}
 $
\xymatrix@M+0.5pc@R-0.5pc@C-1pc{
  & [\bigHecke/\mathsf{J}[2]] \ar[dl] \ar[dr] & \\
[X_{1}/\mathsf{J}[2]] & & [\mycal{X}_{0}/\mathsf{J}[2]]}
$
\end{minipage}
\right]
\ \cong \
\left[
   \begin{minipage}[c]{2.3in}
 $ 
\xymatrix@M+0.5pc@R-0.5pc@C-1pc{
  & \mathbf{Hecke} \ar[dl] \ar[dr] & \\
\mathsf{M}_{1} & & \mathsf{M}_{0}\times C
}
$
\end{minipage}
\right].
\]
Similar comments apply to the action of the big Hecke correspondence
in the other direction. Thus we can recast the problem of finding a
Hecke eigensheaf on $\mathsf{M}_{0}\sqcup \mathsf{M}_{1}$ as the
equivalent problem of finding a $\mathsf{J}[2]$-equivariant Hecke
eigensheaf on $X_{1}\sqcup \mycal{X}_{0}$.

We can refine this further by observing that $\mycal{X}_{0}$
trivializes on a finite cover of $C$. Indeed, recall the
curve $\sq : \Cbar \to C$ which is an \'{e}tale $\mathsf{J}[2]$-Galois
cover of $C$ defined as the fiber product
\[
  \xymatrix@M+0.5pc@-0.5pc{
\Cbar \ar[r]^-{\emb_{\Cbar}} \ar[d]_-{\sq} & \op{Jac}^{0}(C)
\ar[d]^-{\mathsf{mult}_{2}} \\
    C \ar[r]_-{\mathsf{AJ}_{\pw}} & \op{Jac}^{0}(C)
}
\]
where the bottom horizontal arrow is the $\pw$-based Abel-Jacobi map
$\mathsf{AJ}_{\pw} : C \to \op{Jac}^{0}(C)$, $t \mapsto
\mathcal{O}_{C}(t - \pw)$. In other words we have
$\Cbar = \{ (A,t) \in \op{Jac}^{0}(C)\times C \, | \, A^{\otimes
  2}(\pw) = \mathcal{O}_{C}(t) \, \}$.

Because the monodromy of $\mycal{X}_{0}$ over $C$ is $\mathsf{J}[2]$,
the pullback of $\mycal{X}_{0}$ trivializes canonically if we use
$(\mathcal{O}_{C},\pw)$ as the base point on $\Cbar$. That is, we have
a fiber square
\[
  \xymatrix@M+0.5pc@-0.5pc{
X_{0}\times \Cbar  \ar[r] \ar[d]_-{\pr_{\Cbar}} & \mycal{X}_{0} 
\ar[d] \\
    \Cbar \ar[r]_-{\sq} & C
}
\]
where the top horizontal map is given by $(E,(A,t)) \mapsto E\otimes
A^{-1}$. Thus we get a base changed Hecke diagram
\[
  \xymatrix@M+0.5pc@-0.5pc{
    & \bigHeckebar \ar[dl]_-{p} \ar[dr]^-{q} & \\
    X_{1} & & X_{0}\times \Cbar
    }
\]
where the base changed big Hecke correspondence is  the moduli
\[
\bigHeckebar = \left\{ \left( (E,E',\beta), (A,t) \right) \left|
 \text{
    \begin{minipage}[b]{2.7in}
      \[
      \begin{aligned}
        E & \in X_{1}, \
        E'  \in X_{0}, \ (A,t) \in \Cbar \\
        \beta : & E'\otimes A^{-1} \hookrightarrow E, \ 
        \op{supp}(\op{coker}(\beta)) = t
      \end{aligned} 
      \]
    \end{minipage}
  }
  \ 
\right.\right\}
\]
and the maps $p$ and $q$ are defined by
\[
p((E,E',\beta),(A,t)) := E, \quad \text{and} \quad
q((E,E',\beta),(A,t)) := (E',(A,t)).
\]
For future reference, note that $\bigHeckebar$ can also be viewed as a
correspondence 
\[
  \xymatrix@M+0.5pc@-0.5pc{
    & \bigHeckebar \ar[dl]_-{\pzo} \ar[dr]^-{\qzo} & \\
    X_{0} & & X_{1}\times \Cbar
    }
\]
where $\pzo = \op{pr}_{X_{0}}\circ q$ and $\qzo = p\times
(\op{pr}_{\Cbar}\circ q)$.

Thus the Hecke eigensheaf problem on $\mathsf{M}_{0}\sqcup
\mathsf{M}_{1}$ can be reformulated as the problem of finding a
$\mathsf{J}[2]\times \mathsf{J}[2]$-equivariant
$\bigHeckebar$-eigensheaf on $X_{0}\sqcup X_{1}$.
Here an element $(\mathfrak{a}_{1},\mathfrak{a}_{2}) \in
\mathsf{J}[2]\times \mathsf{J}[2]$ acts by
\[
\begin{aligned}
(\mathfrak{a}_{1},\mathfrak{a}_{2})\cdot & E  = E\otimes
  \mathfrak{a}_{1}\otimes \mathfrak{a}_{2}^{-1}, \ \text{for} E\in
  X_{1}, \\ (\mathfrak{a}_{1},\mathfrak{a}_{2})\cdot & (E',(A,t))  =
  (E'\otimes \mathfrak{a}_{1}, (A\otimes \mathfrak{a}_{2}, t)),
  \ \text{for} \ (E',(A,t)) \in X_{0}\times \Cbar, \\
  (\mathfrak{a}_{1},\mathfrak{a}_{2})\cdot & \beta  =
  \beta\otimes \op{id}_{\mathfrak{a}_{1}\otimes \mathfrak{a}^{-1}_{2}}.
\end{aligned}
\]
With this in place we can formulate the Dolbeault version of the
$\bigHeckebar$-Hecke eigensheaf problem as follows.

\

\begin{problem}[{\bfseries Dolbeault $\bigHeckebar$-Hecke eigensheaf
      problem }] \label{prob:Hbareigen} Fix a general semistable
  $SL_{2}(\mathbb{C})$-Higgs bundle $(E,\theta)$ on $C$. Construct
  $\mathsf{J}[2]\times\mathsf{J}[2]$-equivariant tame parabolic Higgs
  bundles   $(\mycal{F}_{0,\bullet},\Phi_{0})$ on $X_{0}$ and 
  $(\mycal{F}_{1,\bullet},\Phi_{1})$ on   $X_{1}$ so that
\begin{itemize}
\item $\mycal{F}_{0,\bullet}$ and $\mycal{F}_{1,\bullet}$ have rank $8$.
\item The parabolic structure on $\mycal{F}_{i,\bullet}$ and the poles
  of $\Phi_{i}$ are along the wobbly divisors in $X_{i}$ for $i = 0,1$.
\item the first and second parabolic Chern classes of each
  $\mycal{F}_{i,\bullet}$ are trivial and $(\mycal{F}_{i,\bullet},\Phi_{i})$
  are stable for $i = 0,1$.
\item $\mycal{F}_{i,\bullet}$ satisfy the $\bigHeckebar$-eigensheaf
  property with eigenvalue $(E,\theta)$. In other words we have 
  \begin{description}
\item[$(X_{1} \text{ to }
  X_{0})$] \qquad ${\displaystyle q_{*}p^{*}\left(\mycal{F}_{1,\bullet},
  \Phi_{1}\right)  =
\left(\mycal{F}_{0,\bullet},\Phi_{0}\right)\boxtimes
\sq^{*}(E,\theta)}$,
\item[$(X_{0} \text{ to }
  X_{1})$]
  \qquad
${\displaystyle \qzo_{*}\pzo^{*}\left(\mycal{F}_{0,\bullet},\Phi_{0}\right)
    =
\left(\mycal{F}_{1,\bullet},\Phi_{1}\right)\boxtimes
\sq^{*}(E,\theta)}$.
\end{description}
where all pullbacks, pushforwards, and tensoring are induced from the
corresponding operations on polarized twistor $\srD$-modules.
\end{itemize}
\end{problem}

\subsection{Comparison with the synthetic approach}
\label{synth-comparison}

In the modular direction used for most of the present paper, we start with the
curve $C$ having a chosen Weierstrass point $\pw \in C$, define the
moduli spaces $X_0$ and $X_1$, and observe that $X_0$ is isomorphic to $\pp^3$  
and $X_1$ is the base locus of a pencil 
of quadrics in $\pp^5$  \cite{NR}. In the following discussion, the Narasimhan-Ramanan 
projective space will be  denoted by 
$\pp^5_{NR}$. 

Choose $a=(A,t)$ in the covering curve $\Cbar$. The Hecke correspondence
$$
\bigHeckebar(a) \subset X_0 \times X_1
$$
has Hecke fibers that are lines in $X_0$, over every point of $X_1$. This yields a 
map $X_1\rightarrow \mathsf{Grass}(2,4)$ to the Grassmanian of lines
in $X_0$. If we write $X_0 = \pp (V)$ for a four-dimensional vector space
$V$ then the Grassmanian embedds by the Plücker embedding 
$\mathsf{Grass}(2,4) \hookrightarrow \pp (\bigwedge ^2V)$ with image a quadric hypersurface. 
We get the composed map
$$
X_1 \rightarrow \mathsf{Grass}(2,4) \hookrightarrow \pp (\bigwedge ^2V).
$$

\begin{proposition}[cf Theorem 4, \cite{NR}]
\label{synthetic-modular-Hecke}
The restriction $\Oo _{\pp (\bigwedge ^2V)}(1)|_{X_1}$ is isomorphic to the line bundle $\Oo _{X_1}(1)$
corresponding to the embedding of $X_1$ as an intersection of two quadrics in $\pp^5_{NR}$.
There is an isomorphism $\pp (\bigwedge ^2V)\cong \pp^5_{NR}$ such that 
$\mathsf{Grass}(2,4)$ becomes one of the quadrics in the pencil containing $X_1$, and the above
composed map identifies with the embedding $X_1 \hookrightarrow \pp^5_{NR}$. Under these
identifications, the modular Hecke correspondence $\bigHeckebar(a)$
is equal to the synthetic incidence Hecke correspondence of Section \ref{ssec:synth.corr}. 
\end{proposition}
\begin{proof}
Consider a line $\ell \subset X_1$ corresponding to a family of bundles fitting into an exact
sequence 
$$
0 \rightarrow L \rightarrow E \rightarrow L^{\vee} (\pw ) \rightarrow 0.
$$
The family of Hecke lines over bundles $E$ in this family consists of the lines in 
$X_0$ contained in a plane and passing through a point. The point is
the $S$-equivalence class of $(L\otimes A) \oplus (L^{\vee} \otimes A(\pw - t))$, while
the plane is the subset of $X_0$ consisting of bundles that have a nontrivlal
map from $L\otimes A(-t)$ (as pointed out in the proof of Proposition \ref{wobblyc0}, 
the discussion of \cite[Proposition 6.1, Theorem 2]{NR} implies that this is a plane). 
 
The Schubert cycle $\sigma _1$ of codimension $1$ in $\mathsf{Grass}(2,4)$, 
the divisor of $\Oo _{\mathsf{Grass}(2,4)}(1)$, is the set of lines passing through
a given general line $V\subset X_0$. Then $V$ intersects the plane in a general point, and there
is exactly one line in the family that passes through that point. Thus, the pullback of the 
Schubert cycle to the line $\ell$ via the embedding $\ell \subset X_1 \rightarrow \mathsf{Grass}(2,4)$
is a single point. This shows that the pullback of $\Oo _{\mathsf{Grass}(2,4)}(1)$
is $\Oo _{X_1}(1)$ since the Picard group of $X_1$ is $\zz$. 

The $6$-dimensional space $\bigwedge^2V$ of sections that embedd the Grassmanian 
pulls back to a 
$6$-dimensional space of sections of $\Oo _{X_1}(1)$, and this is also the $6$-dimensional
space of sections of $\Oo _{\pp ^5_{NR}}(1)$. We get the required identification between
$\pp ^5_{NR}$ and $\pp (\bigwedge ^2V)$. The Grassmanian is
a quadric containing $X_1$ so it is one of the members of the pencil. 

Now $X_0 = \pp (V)$ is recovered as one of the rulings of the quadric Grassmanian,
namely the ruling of Schubert cycles $\sigma _2$: a point $x\in X_0$ corresponds to the Schubert cycle 
$\sigma_2(x)$ of lines through $x$, which is a plane in $\mathsf{Grass}(2,4)$. The dual
projective space parametrizes the Schubert cycles $\sigma_{1,1}(h)$ of lines in a given
plane $h \subset X_0$.

The synthetic Hecke correspondence of Section \ref{ssec:synth.corr} is the incidence correspondence
saying when a point of $X_1$ is in an element of the ruling. Now, the ruling is isomorphic to 
$X_0$ and a point $E$ of $X_1$ is in the element $\sigma _2(x)$ of the ruling corresponding to $x\in X_0$,
if and only if the Hecke line corresponding to $E$ contains $x$. This is the same as
the Hecke correspondence $\bigHeckebar(a)$. 
This equality is indeed tautological, because we already used
$\bigHeckebar(a)$ to get the map $X_1 \rightarrow \mathsf{Grass}(2,4)$ that
led to the identification $\pp (\bigwedge ^2 V) \cong \pp^5_{NR}$. 
\end{proof}

\subsection{Description of Hecke curves}
\label{descrip-Hecke-curves}

The previous proposition shows that
the Hecke correspondence $\bigHeckebar(a)$ may be
viewed, in the synthetic picture of Section~\ref{ssec:synth.corr}, as the incidence correspondence
between points $x$ of $X_1$, corresponding to lines $\ell_x \subset
\pp^3$ and hence also to points of the Grassmanian quadric $G$,
thought of as points in $\pp ^5$ (that happen to be on the other quadric
$G'$ defining the pencil too), and points $y$ of $X_0$, corresponding
to planes $\Pi_y \subset \pp^5$ (that happen to be contained in the
Grassmanian quadric $G$ and to belong to ruling $R$):
\[
\begin{aligned}
\bigHeckebar(a) & = \left\{ (x, y) \in X_1 \times X_0 \; \left|\; x
\in \Pi_y \subset \pp^5 \, \right.\right\} \\
 & = \left\{ (x, y) \in X_1 \times X_0 \; \left|\; y
\in \ell_x \subset \pp^3 \, \right.\right\}.
\end{aligned}
\]
More precisely, if $a=(A,t)$ we have that $t \in C$
corresponds to the ruling $R$ of the quadric $G$ which is identiﬁed as
the Grassmannian.

We now get a clear description of the Hecke curves:
\begin{itemize}
\item Given $x \in X_1$, the Hecke curve $p^{-1}(x)$ is the
  line $\ell _x$ itself. 
\item Given $y \in X_0$,  the Hecke curve $q^{-1}(y)$ is the conic 
$\Pi_y \cap X_1 = \Pi_y \cap G'$. 
    \end{itemize}

\

\noindent
As explained in section~\ref{ssec:synth.corr} we can interpret this
globally as describing a subvariety $\bigHeckebar$ of $\Cbar \times
X_0 \times X_1$.  Recall that we have three copies of $(\zz /2)^4$,
acting respectively on $\Cbar$, $X_0$, and $X_1$ with quotients $C$,
$\mathsf{M}_0$, and $\mathsf{M}_1$. The quotient $\mycal{R}$ of $\Cbar
\times X_0$ by the diagonal action of $(\zz /2)^4$ fibers $\rul :
\mycal{R} \to C$ over $C$ with fibers  non-canonically isomorphic
to $X_0$. In section~\ref{synthetic} we
described  $\mycal{R}$ as the universal ruling for the
pencil of quadrics: it parametrizes the family of planes $\Pi$ contained in
\emph{\bfseries some} member of our pencil.  If we choose
a Weierstrass point $\pw \in C$, we can also identify $\mycal{R}$ with
the family $\mycal{X}_{0}$ of pairs $(t,V)$ with $t \in C$ and $V$ a
rank $2$ bundle on $C$ with determinant $\Oo _C(\pw-t)$.  A point $a
\in \Cbar$ above $t \in C$ determines a square root $A$ of the line
bundle $\Oo_C(t-\pw)$, and tensoring with $A$ converts $V$ to an
$SL(2)$ bundle, i.e. kills its determinant. So the pullback of
$\mycal{X}_{0}$ to $\Cbar$ is identified with $\Cbar \times X_0$,
proving the identification of $\mycal{R}$ with $\mycal{X}_{0}$.
    
Most naturally, the incidence description above gives a Hecke
subvariety $X_{1} \times \mycal{R}$:
\begin{equation} \label{eq:synth.bigHeckebar}
\left\{ (x,y) \in X_1 \times \mycal{R} \, \left| \, x
\subset \Pi_y \subset \pp ^5. \right.\right\},
\end{equation}
which is identified with $\bigHecke \subset X_{1}\subset
\mycal{X}_{0}$ via the isomorphism $\mycal{R} \cong \mycal{X}_{0}$.
As noted above, the Hecke actual $\pp GL(2)$ Hecke correspondence that
we need is a substack $\mathbf{Hecke} \subset C \times \mathsf{M}_0
\times \mathsf{M}_1$.  The relationship is that the subvariety
\eqref{eq:synth.bigHeckebar} is the inverse image of $\mathbf{Hecke}$
under the $(\zz /2)^8$-quotient map $\mycal{R} \times X_{1} \to C
\times \mathsf{M}_0 \times \mathsf{M}_1$. Alternatively, what we are
doing mainly in this paper, is to pull back further to $\bigHeckebar
\subset \Cbar \times X_0 \times X_1$ where the Hecke eigensheaf
calculations may be viewed as most straightforward.

The set of points $p(q^{-1}(y))$ is the set of points in $X_1$ that
admit a Hecke transform equal to $y$. In the correspondence with the
quadric line complex, it is the set of points whose associated line
$\ell_x$ passes through $y$. That was called $X_p$ in \cite{GH}.  The
Hecke curves $p^{-1}(x)$ are projective lines, so they never
degenerate.  The following theorem uses the results of
section~\ref{synthetic} to restate the results of \cite{NR} in the
synthetic language.

\

\smallskip

\begin{theorem}
\label{linefibers}
The fiber $p^{-1}(x)$ over any point $x\in X_1$ is mapped by $q$
isomorphically to a line in $X_0 \cong \pp^3$. This provides a map
$X_1 \rightarrow G=\mathsf{Grass}(2,4)$ to the Grassmanian of lines in
$\pp^3$. Furthermore, the embedding $X_1\subset \pp^5$ extends
uniquely to an embedding of $\mathsf{Grass}(2,4)$ in $\pp^5$
identifying it with one of the quadrics in the pencil that cuts out
$X_1$. The family of projective planes in $\pp^{3}$ consisting of all
the lines through a given point, is a ruling of the quadric $G$, and
this ruling is identified with the point $t\in C$ (image of $a\in
\Cbar$) via the identification between $C$ and the set of pairs of a
quadric in the pencil and a ruling of that quadric.
\end{theorem}

\

\noindent
On the other hand, the Hecke curves $q^{-1}(y)$ are conics, embedded
in planes in $G \subset \pp^{5}$ so they can degenerate into  pairs of
lines. This happens when two distinct lines $\ell, \emm \subset
X_{1}$ intersect, and $\Pi_y$ is their span in $\pp^{5}$. Our choice of
a Weierstrass point $\pw \in C$ allows us to identify the variety of
lines in $X_{1}$ with the degree zero Jacobian $\op{Jac}^{0}(C)$ of
$C$. In particular a line $\ell \subset X_1$ corresponds to a line bundle
$L \in \op{Jac}^{0}(C)$ in such a way that distinct lines $\ell, \emm$
intersect if and only if for the corresponding points $L, M \in
\op{Jac}^{0}(C)$ the degree $1$ line bundle $L \otimes M \otimes
\Oo_C(\pw)$ is effective on $C$, i.e. it is $\Oo_C(t)$ for some $t \in
C$.  The plane $\Pi_y$ spanned by such a pair of lines $\ell, \emm$ is then in
the ruling $\mycal{R}_t$ indexed by this $t$.  The locus of such
$\Pi_{y}$'s in $R=\mycal{R}_{t}$, for a fixed $t \in C$, is isomorphic
to the Kummer surface $\op{Jac}^{0}(C)/({\pm} 1) \subset R$.

It is also possible for the Hecke curve $q^{-1}(y)$ to degenerate into
a double line.  This happens when $\ell = \emm$, which occurs for $16$ lines
$\ell$ if we fix $t$. If we vary $t$, these lines $\ell$ are parametrized by
the cover $\Cbar \to C$, and their union in $X_{1}$ is the wobbly locus of
$X_{1}$. However, in $X_0=R$, their image is a curve (=quotient of $\Cbar
$ by the action of the hyperelliptic involution of $C$) which is
contained in the Kummer.

\

\noindent
This is all analyzed and proved synthetically in
section~\ref{synthetic}. For ease of reference we summarize the
above conclusions as follows:

\

\begin{theorem} \label{conicfibers}
The Hecke fibers over points of the trivial determinant moduli space
$X_{0}$ fall into three categories:
\begin{itemize}
\item Over $X_0 - \Kum$, the map $q$ is smooth with fibers that are
  identified with their images in $X_1$. These images are conics
  inside planes in $\pp^5$. If $y \in X_0 - \Kum$, the corresponding
  conic $q^{-1}(y)$, isomorphic to $\pp^1$, is the intersection of the
  plane $\Pi_{y} \subset \pp^{5}$ with $X_{1}$. The conic $q^{-1}(y)$ is
  identified with the space of rank $1$ quotients of the fiber at
  $t\in C$ of the stable rank $2$ bundle with trivial determinant
  corresponding to $y$. This gives the viewpoint of a Hecke
  correspondence going back from $X_{0}$ to $X_{1}$.
\item Over smooth points of $\Kum$ the fibers of $q$
  degenerate into
  reducible conics composed of two distinct lines meeting at a single point.
\item Over the $16$ singular
  points of $\Kum$, the fiber degenerates further into a double line.
\end{itemize}
\end{theorem}
\begin{proof}
The trichotomy of possibilities was discussed above in Corollary \ref{Hecke-curve-trichotomy},
see \cite[pp 762-763]{GH}. The synthetic Hecke correspondence used there is
identified with the modular Hecke correspondence in question here, by 
Proposition \ref{synthetic-modular-Hecke}. 
\end{proof}

\

Let $\KumKthree (a) \subset \bigHeckebar (a)$ denote the closure of
the subset of points lying over $\Kum$ that are the intersection
points of the two lines in the fibers of $q$ over general points of
$\Kum$.

\

\begin{remark}
The subvariety $\KumKthree (a)$ is the K3 surface obtained by blowing
up the $16$ nodes of $\Kum$.  The double lines over the nodes are the
exceptional divisors.  Its image $p(\KumKthree (a))\subset X_1$ is the
K3 surface denoted by $\Sigma$ in \cite{GH}. This embedding depends on
$a\in \Cbar$.
\end{remark}

\begin{remark}
In fact the surface $\KumKthree(a)$ only depends on the choice of point
$t\in C$ rather than $a\in \Cbar$, due to the fact that the Kummer
surface is invariant by the action of $(\zz / 2\zz )^4$ and in another
viewpoint, the Hecke correspondence between $X_1$ and a moduli space
$\srX _0(t)$ depends only on the choice of $t$; the image of the
Kummer K3 surface back into $X_1$ should therefore only depend on
$t$. We do not need this fact so a full proof is not given. 
\end{remark}

\

\subsection{The Heisenberg group} 
\label{sec-heisenberg}

For our given curve $C$ the group $\mathsf{J}[2]$ of points of order
$2$ on $\op{Jac}^0(C)$ is isomorphic to $(\zz / 2\zz )^4$ and we will adopt
informally the latter notation for this group. It may be viewed as the
group of bundles on $C$ with structure group $\zz / 2\zz$, a group
that we can view in turn as the center of $SL_2(\cc )$. 

There is a natural \emph{\bfseries Heisenberg central extension}
\begin{equation} \label{eq:HeisenbergGm}
1 \rightarrow \Gm \rightarrow \Heisen \rightarrow (\zz / 2\zz )^4
\rightarrow 1
\end{equation}
whose extension class is given by the natural $\zz / 2\zz$-valued
symplectic bilinear form on $(\zz / 2\zz )^4 = H^{1}(C,\zz/2\zz)$
given by the intersection pairing. Indeed, this $\zz / 2\zz$-valued
symplectic bilinear form on $(\zz / 2\zz )^4 = H^{1}(C,\zz/2\zz)$ defines a
finite Heisenberg central extension
$$
1 \rightarrow \zz/2 \rightarrow \Heisen^{\op{fin}} \rightarrow (\zz / 2\zz )^4
\rightarrow 1
$$
which under the natural inclusion $\zz/2 \hookrightarrow \Gm$
induces the extension \eqref{eq:HeisenbergGm}. In particular, by
construction, the finite group $\Heisen^{\op{fin}}$ is a normal
subgroup in $\Heisen$.

Consider  the theta divisor
\[
\Theta =
\left\{ L \in \op{Jac}^{0}(C) \, \left| \, h^{0}(C,L(\pw)) \geq 1
\right.\right\}
\]
corresponding to the theta characteristic $\mathcal{O}_{C}(\pw)$ on
$C$. From the work of Mumford
\cite{mumford-equations1,mumford-theta1,Mumford-Abelian,mumford-theta3,bl}
it is known that the group $\Heisen$ can be identified with the
\emph{\bfseries theta group} of authomorphisms of the total space of
the line bundle $\Oo (2\Theta) \in
\op{Pic}(\op{Jac}^{0}(C))$ that lift the translation action of
$\mathsf{J}[2]$ on $\op{Jac}^{0}(C)$.

In particular $\Heisen$ acts on the vector space $H^0 (\op{Jac}^0(C),
\Oo (2\Theta)) \cong \cc^4$ and this action can be
identified \cite{mumford-equations1,Mumford-Abelian} with the
Schr\"{o}dinger representation, i.e. with  the unique irreducible
representation of $\Heisen$ with a tautological central character. The
projectivisation of $H^0 (\op{Jac}^0(C), \Oo (2\Theta))$
is the Narasimhan-Ramanan model of the moduli space of rank two
bundles with trivial determinant giving an identification $X_0 \cong
\pp^3$ \cite{NR}. The center of $\Heisen$ acts by scalars, so the action
descends to an action of $(\zz / 2\zz )^4$ on $X_0$.

The resulting action of $\Heisen$ on $\cc^6 = \bigwedge ^2(\cc^4)$ is
also an action for which the center of $\Heisen$ acts by scalars and
thus again provides an action of $(\zz / 2\zz )^4$ on $\pp^5$,
preserving and hence acting on $X_1$. If we view $(\zz / 2\zz )^4$ as
being the group of bundles on $C$ with structure group the center of
$SL_2(\cc )$ then these actions are the actions obtained by tensoring
semistable bundles with finite order line bundles that we described in
section~\ref{sec-bigHecke}. It is worth noting that the action of
$(\zz / 2\zz )^4$ on $X_{1}$ can be naturally linearized on the
hyperplane bundle $\Oo_{\pp^{5}}(1)_{|X_{1}}$ while the action of
$(\zz / 2\zz )^4$ on $X_{0}$ is only projective, i.e. it does not
linearize on the hyperplane bundle $\Oo_{X_{0}}(1) =
\Oo_{\pp^{3}}(1)$. Indeed, the group $\Heisen$ and hence its subgroup 
$\Heisen^{\op{fin}} \subset \Heisen$ both act linearly on $\cc^{4}$
and on $\cc^{6} = \bigwedge^{2}\cc^{4}$.  Since the center of
$\Heisen^{\op{fin}}$ acts by multiplication by $\pm 1$ on $\cc^{4}$ it
follows that this center acts trivially on $\cc^{6} =
\bigwedge^{2}\cc^{4}$. In particular the action of $\Heisen^{\op{fin}}$
factrors through a linear action of $(\zz/2\zz)^4$ on $\cc^{6}$. In
contrast, since $\cc^4$ is the Schr\"{o}dinger representation of
$\Heisen^{\op{fin}}$, and this representation determines the
non-trivial extension class defining $\Heisen^{\op{fin}}$, it follows
that the projective action of $(\zz/2\zz)^{4}$ on $\pp^{3}$ can not be
lifted to a linear action on $\cc^{4}$.

The natural family of moduli spaces $\srX _1 \rightarrow C$ whose
fiber over $t\in C$ is the moduli space $\srX _1(t)$ of stable bundles
of determinant $\Oo _C(t)$, is \'{e}tale locally trivial but not
globally trivial. It is obtained by dividing $X_{1}\times \Cbar$ by
the action of $(\zz / 2\zz )^4$, where a $2$-torsion line bundle
$\mathfrak{a} \in (\zz / 2\zz )^4 = \mathsf{J}[2]$ acts by sending a
piar $(E,(A,t)) \in X_{1}\times \Cbar$ to the pair $(E\otimes
\mathfrak{a}, (A\otimes\mathfrak{a},t))$.

This is similar to the family $\srX _0\rightarrow C$ of moduli spaces,
whose fiber over $t\in C$ we described before, where the fiber over $t
\in C$ is the space $\srX_0(t)$ of bundles determinant $\Oo _C(t-\pw
)$. Note that while $\srX_{1}$ depends only on the curve $C$, the
space $\srX_{0}$ depends on having fixed the Weierstrass point $\pw$.

As mentioned above, pulling back to the covering $\Cbar \rightarrow C$
gives a trivialization
$$
\srX _1 \times _C \Cbar \cong X_1 \times \Cbar .
$$
This yields the isomorphism between $\srX _1(t)$ and $X_1$, through
which we pass to obtain the Hecke correspondence between $X_0$ and
$X_1$ depending on the point $a = (A,t) \in \Cbar$.  We similarly have a
trivialization
$$
\srX _0 \times _C \Cbar \cong X_0 \times \Cbar .
$$
If one wants to work with structure group $G=\pp GL(2)$, the moduli of
$G$-bundles on $C$ is a disjoint union $\mathsf{M}= \mathsf{M}_0
\sqcup \mathsf{M}_1$. Each $\mathsf{M}_i$ is a quotient of $X_i$ by
the action of $ (\zz / 2\zz )^4$.  The action on $X_1$ flips the sign
of an even subset of the 6 coordinates. As the point $t$ varies,
$\mathsf{M}_1$ is constant, while $X_1$ is fixed only up to this
subgroup of its finite group of symmetries.  The action on $X_0$ is
the Heisenberg action on $\mathbb{P}^3=\mathbb{P}H^0(\op{Jac}^0(C),
\mathcal{O}(2\Theta))$. Neither $X_0$ nor $\mathsf{M}_0$
depend on $t \in C$. In this viewpoint, the Hecke correspondence is a
$5$-dimensional subvariety of $C \times \mathsf{M}_0 \times
\mathsf{M}_1$. Fixing $t \in C$ and lifting from $\mathsf{M}_i$ to
$X_i$ gives the Hecke correspondence between $X_0$ and $\srX _1(t)$,
and then fixing a lifting of $t$ to $a\in \Cbar$ yields the Hecke
correspondence between $X_0$ and $X_1$.

\subsection{Hecke fibers over the nodes}
\label{overnodes}

The Hecke fibers over nodes of $\Kum$ are lines in $X_1$ counted with
multiplicity two.  In this subsection we indicate their locations, in
particular they will be lines in the wobbly locus. For comparison, we
look also at the special lines on the wobbly locus that correspond to
the trope planes.

For the latter question, suppose given a line in the wobbly locus
parametrizing  all non-split extensions
\[
0 \to L \to E \to L^{\vee}(\pw) \to 0,
\]
for some line bundle $L$ for which $L^{\otimes 2} = \Oo(q-\pw )$.  If
$a = (A,t) \in \Cbar$ is given, then the corresponding Hecke
transforms of such an $E$ are of the form $E'\otimes A$ where $E'$ is
the kernel of some map $E \twoheadrightarrow \cc_{t}$. Now, for all
Hecke transforms except a special one that we will ignore here, the
bundle $E'$ will be an extension
\[
0 \to L(-t) \to E' \to L^{\vee}(\pw) \to 0.
\]
Thus the Hecke transform contains as subbundle $U = L\otimes A(-t)$.
The equation saying that this collection is equal to a trope plane is
the equation $U^{\otimes 2} = \Oo (-2\pw )$. For $U=L\otimes A(-t)$ we have
$$
U^{\otimes 2} = L^{\otimes 2}\otimes A^{\otimes 2}(-2t) = \Oo
(q-\pw + t - \pw - 2t) = \Oo (q-t-2\pw ).
$$
The condition that we are on a trope plane is that this bundle is $\Oo
(-2\pw )$, i.e.  the equation $\Oo (q-t) = \Oo$ in other words, $q=t$.
Thus, we conclude that the $16$ lines on the wobbly locus that
correspond to bundles whose Hecke transforms are trope planes, are the
$16$ lines corresponding to solutions of $L^{\otimes 2} = \Oo (q-\pw
)$ with $q=t$.

Now to the main question: suppose given a node of $\Kum \subset X_0$,
and let us do the Hecke transformation. The Hecke transformation is
unstable if we take the polystable representative to start with, so
let us look at a bundle $E$ that is an extension of $V$ by itself with
$V^{\otimes 2} = \Oo _C$. Then the subsheaf $E'$ has as line subbundle
$V(-t)$, so the Hecke transform $E'\otimes A(\pw )$ has as line
subbundle $L= V\otimes A (\pw -t)$. Let us verify that this is on the
wobbly locus: we have
$$
L^{\otimes 2} = V^{\otimes 2} \otimes A^{\otimes 2} (2\pw - 2t)
= \Oo (t - \pw ) \otimes \Oo (2\pw - 2t) = \Oo (\pw - t).
$$
We may also set $q=t'$ to be the conjugate point of $t$, and note that
$\Oo (t+q) = \Oo (2\pw )$ so we can write
$$
L^{\otimes 2} =\Oo (t - \pw ) \otimes \Oo (2\pw - 2t) =\Oo (t +
q - \pw - t) = \Oo (q - \pw ).
$$
Thus the bundle $E'\otimes A(\pw )$ being an extension containing $L$
as subbundle, is on the wobbly locus. It corresponds to the line over
the point $(L,q)$ of $\Cbar$ with $L=V\otimes A (\pw -t)$ and $q=t'$
is the conjugate point of $t$.

We now have a description of the $16$ lines, that are on the wobbly
locus, that are Hecke transforms of the nodes. They are the $16$ lines
corresponding to solutions of $L^{\otimes 2} = \Oo (q-\pw )$ with
$q=t'$ being the conjugate point of $t$.

Putting together these two collections of $16$ lines, we get $32$
lines on the wobbly locus that correspond to the $32$ points in
$\Cbar$ whose images in $C$ are either $t$ or $t'$.

It looks like these should be the $32$ lines under discussion on
\cite[pp 775-77]{GH}.  In particular, the $32$ lines of \cite{GH}
should have the property that they are `special' in the sense of the
definition of page 792 (i.e. being lines of the wobbly locus). This is
certainly known in the classical theory but does not seem to have been
mentioned in \cite{GH}.

One may also ask to describe the planes in $X_0$ that correspond to
the second collection of $16$ lines i.e. the Hecke fibers over
nodes. We think that the plane corresponding to the Hecke fiber over a
node will contain that node, and will have tangent cone that is the
tangent line to the conic (tangent cone of $\Kum$) at the point of the
hyperelliptic line given by the image of $t$. We do not have a proof of
that; it should follow from a closer look at the theory of the
relationship between the Kummer surface $\Kum$ and its dual
\cite{Keum,GH}, but that goes beyond our present scope.

\subsection{Pullbacks of wobbly divisors and ramification}
\label{pull-wobbly-ramification}

Let 
$$
\Wbar _0 := \pzo ^{-1}(\Wob _0)\subset \bigHeckebar ,\;\;\; 
\Wbar _1 := p^{-1}(\Wob _1) \subset \bigHeckebar
$$
be the inverse images of the wobbly divisors in the big Hecke correspondence. 

Fix a point $a = (A,t) \in \Cbar$.  Let
$\Wbar_{0}(a)$ and $\Wbar_{1}(a)$ denote the fibers of
$\Wbar_{0} \to \Cbar$ and $\Wbar_{1} \to \Cbar$ over $a \in
\Cbar$. We have  a diagram
$$
\xymatrix@M+0.5pc{
\Wbar _{1}(a) \ar@{^{(}->}[r] \ar[dr] & \bigHeckebar (a) \ar[d] \\
& X_0 
}
$$
Recall that the map $p \bigHeckebar(a) \rightarrow X_1$ is a
$\pp^1$-bundle. Therefore, the map $\Wbar_1(a) \rightarrow \Wob _1$ is
smooth.  It follows that the singularities of $\Wbar_1(a)$ are the same
as those of $\Wob_1$ pulled back. Namely, we have a locus of cusps and
a locus of nodes in codimension $1$.

Let $\Wob _1^{\mathsf{n}}\rightarrow \Wob_1$ denote the normalization,
so $\Wob _1^{\mathsf{n}} = \Cbar \times \pp^1$. Similarly denote by
$\Wbar_1(a)^{\mathsf{n}}$ the normalization of $\Wbar_1(a)$, which
maps by a smooth $\pp^1$-fibration to $\Wob_1^{\mathsf{n}} $.

\

We give here some statements about the wobbly from $X_1$ pulled back
and ramifying over $X_0$.

\

\punkt \ The map $\Wbar _{1}(a)\rightarrow X_0$ is a proper
morphism which is a finite $16$-sheeted cover away from a codimension
$2$ subset of $X_0$.

\

\noindent
{\bfseries Proof:}  \ We have computed that the class of the
divisor $\Wob_{1} \in X_{1}$ is $8H$, where $H$ is the hyperplane
class in $\pp^{5}$. This means that $\Wob_{1}$ intersects each line in
$X_{1}$ at $8$ points and each conic in $X_{1}$ at $16$ points. Since
the fibers of $\bigHeckebar(a) \to X_{0}$ are the conics in $X_{1}$ this
shows that  $\Wbar _{1}(a)\rightarrow X_0$ is a map of degree
$16$. \ \hfill $\Box$

\

\punkt  \ Let $J =  \op{Jac}^{0}(C))$ denote the Jacobian of $C$ viewed as
the moduli of lines $\ell \subset X_{1}$. Let
\[
  \Gamma = \{ (x,\ell) \in  X_{1}\times J \ | \ x \in \ell \}
\]
be the incidence correspondence. $\Gamma$ is a $\pp^{1}$-bundle over
$J$. In the interpretation of $X_{1}$ as the moduli of vector bundles
of determinant $\Oo_{C}(\pw)$ a point of the Jacobian $A \in J$
corresponds to the line $\ell_{A} \subset X_{1}$ parametrizing all
non-trivial extensions
\[
0 \to A \to E \to A^{\vee}(\pw) \to 0.
\]
Thus the fiber $\ell_{A}$ of $\Gamma \to J$ over $A$ is identified
canonically with the line $\pp(H^{1}(C,A^{\otimes 2}(-\pw)))$. In
particular if $A \in \Cbar \subset J$, then $A^{\otimes 2} = \Oo_{C}(t
- \pw)$ for some point $t \in C$ and so we have $H^{1}(C,A^{\otimes
  2}(-\pw)) = H^{1}(C,\Oo_{C}(-t'))$ where $t'$ is the hyperelliptic
conjugate of $t$. This space in turn is equal to $H^{1}(C,\Oo_{C})$
  under the natural embedding $\Oo_{C}(-t') \subset \Oo_{C}$. In other
  words the restriction of the $\pp^1$-bundle $\Gamma \to J$ to the
  curve $\Cbar \subset J$ is trivial, i.e.  $\Gamma_{|\Cbar} \cong
  \pp^{1}\times \Cbar$ - the product of the hyperelliptic
  $\pp^1$ and $\Cbar$.
  The projection map $\Gamma \to X_{1}$ is known \cite{Newstead}
  to be  a finite morphism
  of degree $4$.

\

\begin{subclaim}
  \begin{itemize}
\item[(a)] For every line $\ell \subset X_{1}$ there is natural line
  bundle $\gamma(\ell)$ of degree three giving a rational map $C
  \dashrightarrow \ell$. The map has a base point if and only if $\ell
  \subset \Wob_{1}$ ($\ell$ ``wobbles'').
\item[(b)] The branch divisor of the map $\Gamma \to X_{1}$ contains
  $\Wob_{1}$, and the ramification divisor above $\Wob_{1}$ is equal
  to
 \[
   \Gamma_{|\Cbar} \cong \pp^{1}\times \Cbar \to \pp^{1}\times (\Cbar
   \times \Cbar) = (\pp^{1}\times \Cbar)\times \Cbar \to
   \Wob_{1}\times J \subset X_{1}\times J,
\]
where the maps between products are the natural maps on the
components. Also, \linebreak
$\Gamma \to X_{1}$ is simply ramified at the general
point of $\pp^{1}\times \Cbar \subset \Gamma$.
\item[(c)] If $\ell \subset \Wob_{1}$ does not wobble, then the
  intersection points $\ell \cap \Wob_{1}$ are the $8$ branch points
  for the cover $C \to \pp^{1}$. 
\end{itemize}
\end{subclaim}

\begin{proof}
We identify the curve $C$ with the family of rulings of the pencil of
quadrics. A point $t \in C$ Corresponds to a ruling $R_t$ , which is a
family (whose parameter space is isomorphic to $\pp^3$) of planes
$\pp^2$ contained in quadric $Q_{\he(t)}$, where $\he : C \to \pp^1$
is the hyperelliptic map. We get a morphism $i: C \times J \to J$: if
$\ell$ is a line in $X_1$, it is contained in a unique plane of ruling
$R_t$. The intersection of this plane with $X_1$ equals its
intersection with the quadric $Q_{\lambda}$ for any $\lambda \neq
\he(t)$, so it consists of $\ell$ plus another line $ i(t,\ell)$. If
we use the Weierstrass point $\pw$ to embed $C$ in $J =
\op{Jac}^0(C)$, then this morphism becomes $i(t,\ell_{A}) =
\ell_{A^{\vee}(t-\pw)}$. For each $t \in C$, the restriction $i_t: J
\to J$ is an involution, so it has 16 fixed points, namely $\{\ba \ |
\ 2\ba =\Oo(t-\pw)\}$. We see that $\Cbar$ is the union of these fixed loci
as the point $t$ varies over $C$. On the other hand, for each $\ell
\in X_1$, the restriction $i_{\ell}$ identifies $C$ with the family of
lines in $X_1$ intersecting $\ell$. So we get a morphism $j: \Gamma
\times_{X_1} \Gamma \to \Gamma \times C.$

The fiber product $\Gamma \times_{X_1} \Gamma$ is reducible,
consisting of the diagonal plus a 3-sheeted cover $\Gamma' \to
\Gamma$. The image $j(\Gamma') \subset \Gamma \times C$ gives a
morphism $\Gamma \to \Sym^3 C$ which induces the desired $\gamma : J
\to \op{Jac}^3 C.$ For a given $\ell$, the sections of $\gamma(\ell)$
map $C$ to $\ell$ itself. A point $x \in C$ is a base point of
$\gamma(\ell)$ iff $\ell$ is a fixed point of the involution $i_x$. So
as we noted above, $\gamma(\ell)$ has a base point iff $\ell$ is a
line contained in $\Wob_1$. More generally, for a line $\ell \subset
X_{1}$, the points where $\ell$ meets $\Wob_1$ are the points above
which two of the three points of the corresponding divisor in
$|\gamma(\ell)|$ come together .
\end{proof}

\

\begin{subcorollary}
\label{subcor-kummer}
The general line $\ell \subset X_{1}$ intersects $\Wob_{1}$ at $8$
distinct points and therefore $\Wbar_{1}(a) \to X_{0}$ is
unramified over the general point of the Kummer surface $\Kum$. 
\end{subcorollary}

\

\noindent
The non-smooth divisor of the map $\Wbar _{1}(a)\rightarrow X_0$
decomposes into three pieces that we will denote as
$$
\Wbar _{1}(a)^{\rm ramif}, \;\;\; 
\Wbar _{1}(a)^{\rm node}, \;\;\; 
\Wbar _{1}(a)^{\rm cusp}
$$
where $\Wbar _{1}(a)^{\rm node}$ is the nodal locus of $\Wbar _{1}(a)$,
$\Wbar _{1}(a)^{\rm cusp}$ is the cuspidal locus, 
and $\Wbar _{1}(a)^{\rm ramif}$ is the remainder of the ramification locus.

We note that the ramification locus of the map from the normalization
$\Wbar _{1}(a)^{\mathsf{n}} \rightarrow X_0$ will consist of the
pieces mapping to $\Wbar _{1}(a)^{\rm ramif}$ and $\Wbar _{1}(a)^{\rm
  cusp}$, the latter because a cusp gives ramification. Furthermore,
the three pieces are not disjoint since they intersect in codimension
$2$ of $\Wbar _{1}(a)$ but we ignore this aspect.

As $\Wbar _{1}(a)$ is a $3$-dimensional variety these pieces are
two-dimensional, so they map to divisors in $X_0$. The previous
Corollary~\ref{subcor-kummer} tells us that none of these images meet
the Kummer surface.

\begin{subclaim}
The images of $\Wbar _{1}(a)^{\rm node}$ 
and $\Wbar _{1}(a)^{\rm cusp}$ do not contain any trope planes. 
\end{subclaim}
\begin{proof}
The locus $\Wbar _{1}(a)^{\rm cusp}$ is irreducible and the locus
$\Wbar _{1}(a)^{\rm node}$ has $6$ pieces. On the other hand, when all
the data moves around, the trope planes are permuted in an orbit of
size $16$, so a monodromy argument in terms of our general parameters
implies that these divisors can not contain trope planes.
\end{proof}

\begin{proposition} Fix a point  $a = (A,t) \in \Cbar$. 
Away from the inverse image of a codimension $2$ subset of $X_0$, 
$$
\Wbar _{1}(a)^{\rm ramif} =
\bigcup _{\mathfrak{a} \in \mathsf{J}[2]} \Wbar _{1}(a)^{\rm ramif}_{\mathfrak{a}}
$$
is a union of $16$ divisor pieces
$\Wbar _{1}(a)^{\rm ramif}_{\mathfrak{a}} \subset
\Wbar_{1}(a)$ each of which is a ruled surface $\mathbb{F}_1$
mapping to the corresponding trope plane
$\trope_{\mathfrak{a}(\pw)} \cong \pp^2$.

The map back to $\Wob _1$ decomposes as a sum of maps (ruled surfaces) 
$$
\Wbar _{1}(a)^{\rm ramif}_{\mathfrak{a}} 
\rightarrow ({\rm line})_{\mathfrak{a}}(a) \subset \Wob _1
$$
where $ ({\rm line})_{\mathfrak{a}}(a)$ is one of the $\pp^1$'s
in the expression
$\Wob _1^{\mathsf{n}} = \Cbar \times \pp^1\rightarrow \Wob _1$
and it corresponds to the point $b=(A\otimes\xi,t)\in \Cbar$. 
\end{proposition}
\begin{proof}
If $b=(B,u)\in \Cbar$ then the line $\{ b\} \times \pp^1\subset \Wob
_1$ pulls back to an ${\mathbb F} _1$-surface in $\bigHeckebar
(a)$. Such a surface then contracts to a plane that we will denote by
$\pp^2 (b)\subset X_0$.  The contraction has an exceptional $\pp^1$
that maps to a point we will denote by $x(b) \subset X_0$.

The line corresponds to bundles fitting into an extension 
$$
0 \rightarrow B \rightarrow E \rightarrow B^{-1}(\pw ) \rightarrow 0
$$
and the Hecke transformed bundle will be of the form $E'\otimes A$
where $E'$ is the kernel of a length $1$ quotient of $E$ supported on
$t \in C$. In particular, the Hecke transformed bundle contains the
line bundle $B\otimes A(-t)$, and this determines the plane $\pp^2
(b)$ as being the plane of semistable bundles in $X_0$ that contain
the degree $-1$ line bundle $B\otimes A(-t)$.

Let $\Wbar _1(a)^{\mathsf{p}}$ denote the resulting family of $\pp^2$'s. 
It is a $\pp^2$-bundle over $\Cbar$ with a dominant map to $X_0 = \pp^3$. 

The map $b\mapsto x(b)$ has image a curve being the subset of
codimension $2$ appearing in the statement of the proposition.
Outside of this curve, $\Wbar _1(a)^{\mathsf{p}}$ and $\Wbar
_1(a)^{\mathsf{n}}$ coincide.  Therefore it is enough to locate the
ramification of $\Wbar _{1}(a)^{\mathsf{p}} \rightarrow X_0$.

The map $b\mapsto \pp^2(b)$ may be viewed as a map to the dual projective
space that we will denote by $X_0^{\vee}$. It factors as
$$
\Cbar \rightarrow \op{Jac}^{-1}(C) \rightarrow \Kum ^{\vee} 
\rightarrow X_0^{\vee}
$$
where $\Kum ^{\vee}$ is the \emph{\bfseries dual Kummer surface}
appearing in the classical synthetic theory \cite{Keum,GH}. It is the
dual surface of $\Kum$.  The birational map sending a point to its
tangent plane blows up the nodes of $\Kum$ and blows down the trope
conics to get the dual Kummer surface $\Kum^{\vee}$. We remark that an
important element of the classical theory says that $\Kum^{\vee}$ is
projectively isomorphic to $\Kum$ although not by this birational map.

The map $ \op{Jac}^{-1}(C) \rightarrow \Kum ^{\vee}\subset X_0^{\vee}$
is not an immersion, exactly over the nodes of the dual surface, and
there the tangent map vanishes. Any translation of the map $\Cbar
\rightarrow \op{Jac}^{-1}(C)$ is an immersion.  Thus, the map $\Cbar
\rightarrow X_0^{\vee}$ is non-immersive exactly at points that go to
nodes of $\Kum ^{\vee}$.

This implies that the family of $\pp^2(b)$'s
becomes stationary exactly whenever the plane  in question is a trope plane. 
In particular, 
the map $\Wbar _{1}(a)^{\mathsf{p}} 
\rightarrow X_0$ is ramified along the planes in $\Wbar _{1}(a)^{\mathsf{p}}$
that map to trope planes. 

We claim that there are $16$ of these, and then we will see by a
characteristic class calculation that this accounts for all of the
ramification.

The condition that the plane $\pp^2(b)$ of vector bundles containing
the line bundle $B\otimes A(-t)$ should be a trope plane, is
equivalent to saying that $B\otimes A(-t)$ is the dual of a square
root of the canonical bundle. We calculate
$$
\left[ B\otimes A(-t)\right] ^{\otimes 2}
= B^{\otimes 2} \otimes A^{\otimes 2} \otimes \Oo _C(-2t) 
=\Oo _C(u - \pw + t - \pw - 2t) =
\Oo _C(u-t-2\pw ).
$$
This is isomorphic to the dual of the canonical bundle exactly when
$u=t$. Thus the set of choices of $b$ making $\pp^2(b)$ into a trope
plane is the set of choices of square-root $B$ such that $B^{\otimes
  2} = \Oo _C(t-\pw )$. There are $16$ of these.

To complete the proof, we need to see that we have accounted for all
of the ramification. Notice that there is another reason for
ramification, namely the fact that the map \linebreak $\Wbar
_1(a)^{\mathsf{n}}\rightarrow X_0$ factors through $\Wbar
_1(a)\rightarrow X_0$ that has a cuspidal locus. This provides a
divisor of ramification that is the inverse image of the cuspidal
locus in $\Wob _1$ and consists of a line in each of the planes
$\pp^2(b)$.

Let $\omega$ denote here the canonical class of $\Wbar
_1(a)^{\mathsf{p}}$.  Its restriction to each fiber $\pp^2(b)$ is $\Oo
_{\pp^2}(-3)$. Therefore, the relative canonical class of $\Wbar
_1(a)^{\mathsf{p}}$ over $X_0$, which is $\Oo _{X_0}(4)\otimes
\omega$, restricts to $\Oo _{\pp _2}(1)$ on each plane $\pp^2(b)$.

We already know a divisor, the ramification coming from the cuspidal
locus, that is a line in each $\pp^2(b)$. This therefore accounts for
any ramification that restricts to a divisor in the $\pp^2(b)$. The
leftover possibility is of ramification consisting of a union of
fibers. But, if the map were ramified along a fiber then the family of
planes $\pp^2(b)\subset X_0$ will be stationary at that fiber. Above
we saw that the cases where the planes are stationary correspond
exactly to the cases where they are trope planes. We conclude
that the ramification locus of the map $\Wbar _{1}(a)^{\mathsf{p}}
\rightarrow X_0$ consists of the $16$ planes corresponding to lines in
$\Wob _1$ as are described in the statement of the proposition, plus
the divisor that maps to $\Wbar _{1}(a)^{\rm cusp}$.
\end{proof}

\

\noindent
Similar analysis proves the following

\

\begin{subcorollary}
The ramification along the planes given in the proposition is simple. 
\end{subcorollary}

\section{Abelianized Hecke}
\label{chapter-abelianized}

\subsection{Setup for abelianization of the Hecke property}
\label{ssec:setup.abelianze}

Let as before $C$ be our smooth genus $2$ curve, and let $\pw \in C$
be our fixed Weierstrass point. Let $(E,\theta)$ be an eigenvalue
$SL_{2}(\mathbb{C})$-Higgs bundle for our Hecke eigensheaf
problem. That is
\begin{itemize}
\item $E$ is a rank two algebraic vector bundle on $C$
  with trivial determinant.
\item $\theta : E \to E\otimes \omega_{C}$ is a traceless Higgs field
  satisfyig the genericity condition that the spectral curve
 \[
\Ctilde \, : \, \det(\lambda\cdot\op{id} - \boldsymbol{\pi}^{*}\theta) = 0
\ \subset \ T^{\vee}C
\]
is smooth.
\end{itemize}

Note that this genericity condition automatically implies that the
Higgs bundle $(E,\theta)$ is stable, i.e. corresponds to a flat bundle
and a polarized twistor $\srD$-module on $C$ by the non-abelian Hodge
correspondence \cite{SimpsonHiggs}.

Let $\specN \in \op{Pic}(\Ctilde)$ be the spectral line bundle
corresponding to $(E,\theta)$. That is $\specN$ is a line bundle on
$\Ctilde$ such that
\[
E = \pi_{*}\specN, \quad \theta = \pi_{*}(\lambda\otimes (-)).
\]
The spectral data $(\Ctilde\subset T^{\vee}C,\specN)$ is the
{\em\bfseries abelianization} of the Higgs bundle $(E,\theta)$. We
used this spectral data to construct modular spectral data defininig a
tame parabolic rank $8$ Higgs bundle on the moduli of semistable
$\mathbb{P}SL_{2}(\mathbb{C})$-bundles on $C$. We already checked that
this rank $8$ Higgs bundle is stabe and has vanishing first and second
parabolic Chern class, hence by Mochizuki's theorem
\cite{Mochizuki-D1,Mochizuki-D2} it gives rise to a polarized twistor
$\srD$-module and a tame parabolic flat bundle. Now we will use the
modular spectral data to rewrite the Hecke eigensheaf property on this
rank $8$ Higgs bundle in terms of abelianized information and the
abelianized Hecke correspondence. We will then explain how to use the
abelianized picture to prove that the rank $8$ Higgs bundle is indeed
a Hecke eigensheaf satisfying the conditions in
Problem~\ref{prob:Hbareigen}.

\

\noindent
Our strategy in constructing $(\mycal{F}_{i,\bullet},\Phi_{i})$ and
checking the eigensheaf property is to use abelianization. To
abelianize the problem we look at the $\Ctilde$-Hitchin fibers inside
the moduli of Higgs bundles $\Higgs_{0}$ and $\Higgs_{1}$
respectively. These are the degree $2$ and degree $3$ Prym varieties
$\Prym_{2}$ and $\Prym_{3}$ of the cover $\pi : \Ctilde \to C$:
\[
\begin{aligned}
\Prym_{2} & = \left\{ \left. M \in \op{Jac}^{2}(\Ctilde) \ \right|
\ \mathsf{Nm}_{\pi}(M) = \omega_{C} \, \right\}, \\
\Prym_{3} & = \left\{ \left. M \in \op{Jac}^{3}(\Ctilde) \ \right|
\ \mathsf{Nm}_{\pi}(M) = \omega_{C}(\pw) \, \right\}.
\end{aligned}
\]
We have natural rational maps
\begin{equation} \label{eq:ratPrym}
  \xymatrix@M+0.5pc@R-2.5pc{
    \Prym_{2} \ar@{-->}[r] & X_{0}, \\
      M \ar@{|->}[r] & \pi_{*}M}
      \quad  \text{and} \quad
      \xymatrix@M+0.5pc@R-2.5pc{\Prym_{3} \ar@{-->}[r] & X_{1}, \\
       M \ar@{|->}[r] & \pi_{*}M}  
\end{equation}
which are surjective, quasifinite, and finite of degree $8$ over the
very stable loci in $X_{0}$ and $X_{1}$. These maps fit in  the following
commutative diagrams of correspondences 
\begin{equation} \label{eq:basicPrymdiag}
\begin{minipage}[c]{2in}
\[
\xymatrix@M+0.3pc@R-0.8pc{ & \Prym_{2}\times \Chat \ar@{-->}[dd]
  \ar[ld]_-{\summ} \ar[rd]^-{\mathsf{id}} & \\
  \Prym_{3} \ar@{-->}[dd]_-(0.3){\pi_{*}}  & &
  \Prym_{2}\times \Chat \ar@{-->}[dd]^-(0.3){\pi_{*}\times \pihat} \\
  & \bigHeckebar
  \ar[ld]_-{p} \ar[rd]^-{q} & \\ X_{1} & & X_{0}\times\Cbar}
\]
\end{minipage} \qquad  \qquad \qquad 
\begin{minipage}[c]{3in}
\[
\xymatrix@M+0.3pc@R-0.8pc{ & \Prym_{3}\times \Chat \ar@{-->}[dd]
  \ar[ld]_-{\diff} \ar[rd]^-{\mathsf{id}} & \\
  \Prym_{2} \ar@{-->}[dd]_-(0.3){\pi_{*}}  & &
  \Prym_{3}\times \Chat \ar@{-->}[dd]^-(0.3){\pi_{*}\times \pihat} \\
  & \bigHeckebar
  \ar[ld]_-{\pzo} \ar[rd]^-{\qzo} & \\ X_{0} & & X_{1}\times\Cbar}
\]
\end{minipage}
\end{equation}
where $\Chat$ is the curve defined by the fiber product
\[
 \xymatrix@M+0.5pc@R-0.5pc@C-0.5pc{
    \Chat \ar[r]^-{\pihat} \ar[d]_-{\sqhat} & \Cbar
    \ar[d]^-{\sq} \\ \Ctilde \ar[r]_{\pi} & C},
 \]
the maps $\summ/\diff$ are the sum/difference maps
\[
\xymatrix@M+0.5pc@R-2.5pc{
\summ : \hspace{-3.5pc} & \Prym_{2}\times \Chat \ar[r] &  \Prym_{3}, \\
& (L,(A,\tilde{t})) \ar@{|->}[r] & L\otimes \pi^{*}A^{-1}(\tilde{t})} \qquad
\xymatrix@M+0.5pc@R-2.5pc{
\diff : \hspace{-3.5pc} & \Prym_{3}\times \Chat \ar[r] &  \Prym_{2}, \\
& (L,(A,\tilde{t})) \ar@{|->}[r] & L\otimes \pi^{*}A(-\tilde{t})}
\]
while the rational maps
$\Prym_{2}\times \Chat \dashrightarrow \bigHeckebar$ and
$\Prym_{3}\times \Chat \dashrightarrow \bigHeckebar$
are defined by
\[
\begin{aligned}
\Prym_{2}\times \Chat & \dashrightarrow \bigHeckebar, \qquad
(L,(A,\tilde{t})) \mapsto ((E,E',\beta),(A,t)),  \\
\Prym_{3}\times \Chat & \dashrightarrow \bigHeckebar, \qquad
(M,(A,\tilde{t})) \mapsto ((V,V',\gamma),(A,t))
\end{aligned}
\]
where
\begin{equation} \label{eq:P2andP3toHbar}
\begin{minipage}[b]{3in}
  \[
\left| \
\begin{aligned}
  & E  = \pi_{*}\left((L\otimes \pi^{*}A^{-1})(\tilde{t})\right), \\[-0.2pc]
  & E' = \pi_{*}L, \\[-0.2pc]
  & \beta  = \pi_{*}\left[L\otimes \pi^{*}A^{-1}
  \hookrightarrow (L\otimes \pi^{*}A^{-1})(\tilde{t})\right], \\[-0.2pc]
  & t  = \pi(\tilde{t}).
\end{aligned}
\right.
\]
\end{minipage}
\quad \text{and} \quad
\begin{minipage}[b]{2in}
  \[
\left| \
\begin{aligned}
  & V  = \pi_{*}M, \\[-0.2pc]
  & V' = \pi_{*}\left((M\otimes \pi^{*}A(-\tilde{t})\right) , \\[-0.2pc]
  & \gamma  = \pi_{*}\left[M(-\tilde{t}) \hookrightarrow M\right], \\[-0.2pc]
  & t  = \pi(\tilde{t}).
\end{aligned}
\right.
\]
\end{minipage}
\end{equation}

\

\noindent
The modular spectral covers
$Y_{0}$ and $Y_{1}$ corresponding to $\Ctilde$ are minimal resolutions
\[
\xymatrix@M+0.5pc{
Y_{0} \ar[d]_-{\blo_{o}} \ar[r]^-{f_{0}} & X_{0} \\
\Prym_{2} \ar@{-->}[ur] &
} \quad \text{and} \quad
  \xymatrix@M+0.5pc{Y_{1} \ar[d]_-{\blo_{1}} \ar[r]^-{f_{1}} & X_{1} \\
\Prym_{3} \ar@{-->}[ur] &
}
  \]
of the rational maps in \eqref{eq:ratPrym}.
We will recall the explicit construction of these resolutions below.
The resulting covering maps $f_{i} : Y_{i} \to X_{i}$ are finite of
degree $8$ and we will use them to construct the Hecke eigensheaf.

\

\noindent
Let $(\Ctilde \subset T^{\vee}C,\specN)$ be the spectral data
abelianizing the Higgs bundle $(E,\theta)$ on $C$.  Our approach to
solving Problem~\ref{prob:Hbareigen} is to use the Fourier-Mukai
transform on the Jacobian of $\Ctilde$ to convert the skyscraper sheaf
$\mathcal{O}_{\specN}$ to line bundles on $\Prym_{2}$ and $\Prym_{3}$
respectively.  After that we pull back these line bundles to $Y_{0}$,
$Y_{1}$, modify them appropriately along the exceptional divisors of
the blow-up maps $\blo_{0} : Y_{0} \to \Prym_{2}$ and $\blo_{1} :
Y_{1} \to \Prym_{3}$ and push them forward via $f_{0}$ and $f_{1}$ to
get the eigensheaf: a pair of parabolic Higgs bundles
$(\mycal{F}_{0,\bullet},\Phi_{0})$ $(\mycal{F}_{1,\bullet},\Phi_{1})$
on $X_{)}$ and $X_{1}$ respectively, which are stable and satisfy the
vanishing Chern class conditions from Problem~\ref{prob:Hbareigen}.

Once the parabolic Higgs bundles $(\mycal{F}_{0,\bullet},\Phi_{0})$
and $(\mycal{F}_{1,\bullet},\Phi_{1})$ are constructed via this
procedure, we can use the diagram \eqref{eq:basicPrymdiag} to rewrite
the eigensheaf equation from Problem~\ref{prob:Hbareigen} as an
equation of line bundles. Specifically, the $(X_{1} \text{ to }
X_{0})$ part of the eigensheaf condition becomes an equation of line
bundles on $Y_{0}\times \Chat$ while the $(X_{0} \text{ to } X_{1})$
part of eigensheaf condition becomes an equation of line bundles on
$Y_{1}\times \Chat$.  To carry this out and to show that our
Fourier-Mukai constructions satisfy the equations on $Y_{0}\times
\Chat$ and $Y_{1} \times \Chat$ we will need the abelianized Hecke
correspondence $\bigHeckehat^{ab}$ which can be described either as
the minimal blow-up of $\Prym_{2}\times \Chat$ resolving the rational
map $\Prym_{2}\times \Chat \dashrightarrow \bigHeckebar$ or as the
minimal blow-up of $\Prym_{3}\times \Chat$ resolving the rational map
$\Prym_{3}\times \Chat \dashrightarrow \bigHeckebar$. In the next
section we will construct and analyze this minimal resolution in
detail.

\subsection{The abelianized Hecke correspondence in context}
\label{ssec:abelinized.in.context}

Recall that the (big) abelianized Hecke space $\bigHeckehat^{\rm ab}$
depends on the choice of a fixed spectral curve $\pi : \Ctilde \to C$
and can be viewed either as a correspondence between $Y_{1}$ and
$Y_{0}\times \Chat$ or as a correspondence between $Y_{0}$ and
$Y_{1}\times \Chat$:
\begin{equation} \label{eq:abHecke}
  \begin{minipage}[c]{2.5in}
 \[
\xymatrix@M+0.25pc{
  & \bigHeckehat^{\rm ab} \ar[dl]_-{p^{\rm ab}}  \ar[dr]^-{q^{\rm ab}}  & \\
 Y_{1} & & Y_{0}\times \Chat
} 
\]
  \end{minipage} \qquad
\begin{minipage}[c]{2.5in}
 \[
\xymatrix@M+0.25pc{
  & \bigHeckehat^{\rm ab} \ar[dl]_-{\pzo^{\rm ab}}  \ar[dr]^-{\qzo^{\rm ab}}  & \\
 Y_{0} & & Y_{1}\times \Chat
} 
\]
\end{minipage} 
\end{equation}
Here
\begin{itemize}
\item $Y_{0}$ is the blow-up of the Prym variety
\[
\Prym_{2} = \left\{ \left. M \in \op{Jac}^{2}(\Ctilde) \ \right|
\ \mathsf{Nm}_{\pi}(M) = \omega_{C} \, \right\},
\]
at the $16$ points $\{
  \pi^{*}\kappa\}_{\kappa \in \mathsf{Spin}(C)} \subset \Prym_{2}$,
  where  $\mathsf{Spin}(C) = \{ \kappa \in \op{Jac}^{1}(C)
  \, | \, \kappa^{\otimes 2} = \omega_{C} \}$ denotes the set of theta
  characteristics on $C$. We will write $\blo_{0} : Y_{0} \to
  \Prym_{2}$ for the blow up morphism   and $\ExY_{0} = \sqcup_{\kappa
    \in \mathsf{Spin}(C)}
  \ExY_{0,\kappa}$  for the corresponding exceptional divisor. 
\item $Y_{1}$ is the blow-up of the Prym variety
\[
\Prym_{3} = \left\{ \left. M \in \op{Jac}^{3}(\Ctilde) \ \right|
\ \mathsf{Nm}_{\pi}(M) = \omega_{C}(\pw) \, \right\},
\]
at the image of the map \eqref{eq:ChatinY1} which is explicitly given by
\[
\emb_{\Chat} : \Chat \to \Prym_{3}, \quad (A,\tilde{t}) \mapsto
\pi^{*}(A^{-1}(\pw))\otimes \mathcal{O}_{\Ctilde}(\tilde{t}).
\]
We will write $\blo_{1} : Y_{1} \to \Prym_{3}$ for the blow up
morphism and $\ExY_{1}$ for the corresponding exceptional divisor.
\end{itemize}
Note that $\ExY_{0}$ is disconnected with $16$ connected components
$\ExY_{0,\kappa}$ each isomorphic to $\mathbb{P}^{2}$. In contrast the
divisor $\ExY_{1}$ is irreducible.  In fact, as we saw in
Lemma~\ref{lemma:embedChat} and Lemma~\ref{lem:normal.deg1}{\bfseries
  (b)} it is smooth and isomorphic to $\Chat\times \mathbb{P}^{1}$.

\

\noindent
We are now ready to describe the space $\bigHeckehat^{\rm ab}$ as a
common blow-up of $Y_{0}\times \Chat$ and $Y_{1}\times \Chat$.  Namely
$\bigHeckehat^{\rm ab}$ is the blow up of $Y_{0}\times \Chat$
  centered at the strict transform of the surface $\Chat\times \Chat
  \subset \Prym_{2} \times \Chat$, where the inclusion $\Chat\times
  \Chat \subset \Prym_{2}\times \Chat$ is given by
\begin{equation} \label{eq:Hab_centerY_0}
  \xymatrix@M+0.5pc@R-2pc{
\emb_{\Chat\times\Chat} : \hspace{-3.5pc} & \Chat\times \Chat \ar@{^{(}->}[r] &
  \Prym_{2}\times \Chat
  \\
& ((A_{1},\tilde{t}_{1}),(A_{2},\tilde{t}_{2})) \ar[r] &
  \left(\pi^{*}\left(A_{2}\otimes A_{1}^{\vee}(\pw)\right)(\tilde{t}_{1}
  - \tilde{t}_{2}), (A_{2},\tilde{t}_{2})\right).  }
\end{equation}
The fact that \eqref{eq:Hab_centerY_0} is a closed embedding follows from
Lemma~\ref{lemma:embedChat} since the restriction of  the map
\eqref{eq:Hab_centerY_0} to $\Chat\times \{(A_{2},\tilde{t}_{2})\}$
is a translate of the map \eqref{eq:ChatinY1}.

Note next that for each theta characteristic $\kappa \in \mathsf{Spin}(C)$,
the curve $\{\pi^{*}\kappa\} \times \Chat \subset \Prym_{2}\times \Chat$
is contained in the image of \eqref{eq:Hab_centerY_0}. Indeed, if
$(A,\tilde{t}) \in \Chat$, then $(A\otimes \kappa^{\vee}(\pw),
\tilde{t})$ is also in $\Chat$, and the image of the pair $((A\otimes
\kappa^{\vee}(\pw),\tilde{t}), (A,\tilde{t}))$ under the map
\eqref{eq:Hab_centerY_0} is precisely $(\pi^{*}\kappa,(A,\tilde{t}))
\in \Prym_{2}\times \Chat$. This implies that the strict transform of
$\Chat\times\Chat$ in 
\[
Y_{0}\times \Chat  = \mathsf{Bl}_{\hspace{-0.7pc}
  \underset{\kappa \in \mathsf{Spin}(C)}{\bigsqcup}
\hspace{-0.3pc}   \{\pi^{*}\kappa\}\times
\Chat} \, \left(\Prym_{2}\times \Chat\right)
\]
is still equal to $\Chat\times\Chat$. Hence we can describe
the abelianized Hecke correspondence $\bigHeckehat^{\rm ab}$
as the iterated blow up
\[
\bigHeckehat^{\rm ab} =
\mathsf{Bl}_{\Chat\times \Chat} \, \left(Y_{0}\times \Chat\right) =
\mathsf{Bl}_{\Chat\times \Chat} \, \left( 
\mathsf{Bl}_{\hspace{-0.7pc} \underset{\kappa \in \mathsf{Spin}(C)}{\bigsqcup}
\hspace{-0.3pc}   \{\pi^{*}\kappa\}\times
\Chat} \, \left(\Prym_{2}\times \Chat\right)\right)
\stackrel{\blo}{\longrightarrow} \Prym_{2}\times \Chat.
\]
Furthermore, since the centers of these succesive blowups
are smooth and nested in each other
\[
\underset{\kappa \in
  \mathsf{Spin}(C)}{\bigsqcup}
\hspace{-0.3pc} \{\pi^{*}\kappa\}\times \Chat \subset \Chat\times \Chat
\]
we can also describe
$\bigHeckehat^{\rm ab}$ by performing the blow ups in the opposite
order, i.e. first blow up $\Prym_{2}\times \Chat$ along $\Chat\times
\Chat$, and then blow up  along the strict transforms of the curves
$\pi^{*}\kappa\times \Chat$. In other words, we also have
\[
\bigHeckehat^{\rm ab} = 
\mathsf{Bl}_{\hspace{-0.7pc} \underset{\kappa \in \mathsf{Spin}(C)}{\bigsqcup}
\hspace{-0.3pc}   \{\pi^{*}\kappa\}\times
  \Chat} \, \mathsf{Bl}_{\Chat\times \Chat} \, \Prym_{2}\times \Chat.
\]
For ease of reference we introduce notation for the exceptional
divisors of these blow ups. We let
$\bigExc_{1} \subset \bigHeckehat^{\rm ab}$ denote the exceptional
divisor of the blow up map
$\bigHeckehat^{\rm ab} \to Y_{0}\times \Chat$, and let
$\bigExc_{0} = \sqcup_{\kappa} \bigExc_{0,\kappa}\subset
\bigHeckehat^{\rm ab}$ denote the strict transform of the divisor
$\ExY_{0}\times \Chat \subset Y_{0}\times \Chat$. Equivalently
$\bigExc_{0}$ is the exceptional divisor of the blow up map
$\bigHeckehat^{\rm ab} \to \mathsf{Bl}_{\Chat\times \Chat} \,
\Prym_{2}\times \Chat$, and $\bigExc_{1}$ is the strict transform of
the exceptional divisor of the map
$\mathsf{Bl}_{\Chat\times \Chat} \, \Prym_{2}\times \Chat \to
\Prym_{2}\times \Chat$.

The map $q^{\rm ab} : \bigHeckehat^{\rm ab} \to Y_{0}\times \Chat$ is just
the blow up morphism, while the map $p^{\rm ab} : \bigHeckehat^{\rm ab} \to
Y_{1}$ lifts the Abel-Jacobi sum map
\[
\summ : \Prym_{2}\times \Chat \to \Prym_{3}, \quad
\summ(L,(A,\tilde{t})) = L\otimes \pi^{*}A^{-1}(\tilde{t})
\]
to the blow ups of the source and target. Recall that  $Y_{1} =
\mathsf{Bl}_{\Chat} \, \Prym_{3}$ is the blow up of $\Prym_{3}$
centered at the curve $\Chat \hookrightarrow \Prym_{3}$,
$(A,\tilde{t}) \to \pi^{*}(A^{-1}(\pw))(\tilde{t})$.  

The preimage
$\summ^{-1}(\Chat)$ of this curve under $\summ$ is exactly
the surface $\Chat\times \Chat \subset \Prym_{2}\times \Chat$ given by
the image of \eqref{eq:Hab_centerY_0}. So by the universal property of
blow ups we see that the map $\summ$ lifts to a morphism
\[
\mathsf{Bl}_{\Chat\times \Chat} \, \Prym_{2}\times \Chat \to
\mathsf{Bl}_{\Chat} \Prym_{3} = Y_{1}
\]
which we can further precompose with the second blow up map
\[
\bigHeckehat^{\rm ab} = 
\mathsf{Bl}_{\hspace{-0.7pc} \underset{\kappa \in \mathsf{Spin}(C)}{\bigsqcup}
\hspace{-0.3pc} \{\pi^{*}\kappa\} \times \Chat} \, \mathsf{Bl}_{\Chat\times
  \Chat} \, \Prym_{2}\times \Chat \longrightarrow \mathsf{Bl}_{\Chat\times
  \Chat} \, \Prym_{2}\times \Chat
\]
to obtain the morphism $p^{\rm ab} : \bigHeckehat^{\rm ab} \to Y_{1}$.
This exhibits  $\bigHeckehat^{\rm ab}$ as a blow-up of $Y_{0}\times
\Chat$ and gives the first correspondence diagram in \eqref{eq:abHecke}.

To realize $\bigHeckehat^{\rm ab}$ as a blow-up of $Y_{1}\times \Chat$
and describe the second correspondence diagram in \eqref{eq:abHecke}
note that the products $\Prym_{2}\times \Chat$ and $\Prym_{3}\times
\Chat$ are naturally isomorphic via the pair of mutually inverse maps
\[
\xymatrix@M+0.5pc@R-2.5pc{
  (\diff, \op{id}) : \hspace{-3.5pc} &
  \Prym_{3}\times \Chat \ar[r] & \Prym_{2}\times \Chat, \\
  & (M,(A,\tilde{t})) \ar@{|->}[r] &
  (M\otimes \pi^{*}A(-\tilde{t}),(A,\tilde{t}))
}
\]
and
\[
\xymatrix@M+0.5pc@R-2.5pc{
  (\summ,\op{id}) : \hspace{-3.5pc} &
  \Prym_{2}\times \Chat \ar[r] & \Prym_{3}\times \Chat. \\
  & (L,(A,\tilde{t})) \ar@{|->}[r] &
  (L\otimes \pi^{*}A^{-1}(\tilde{t}),(A,\tilde{t}))
}
\]
But postcomposing the inclusion \eqref{eq:Hab_centerY_0} with the map
$\summ\times \op{id}$  gives the map
\[
\xymatrix@M+0.5pc@R-2.5pc{
\emb_{\Chat}\times \op{id} : \hspace{-3.5pc} &
  \Chat \times \Chat \ar[r] & \Prym_{3}\times \Chat, \\
  & ((A_{1},\tilde{t}_{1}),(A_{2},\tilde{t}_{2})) \ar@{|->}[r] &
  (\pi^{*}(A_{1}^{-1}(\pw))(\tilde{t}),(A_{2},\tilde{t}_{2}))
}
\]
gives the product of $\Chat \subset \Prym_{3}$  with $\Chat$.
Therefore the composition
\[
\xymatrix@1@C+2pc@M+0.5pc{
\op{Bl}_{\Chat\times\Chat} \, (\Prym_{2}\times \Chat) \ar[r] & 
\Prym_{2}\times \Chat \ar[r]^-{(\summ,\op{id})} &
\Prym_{3}\times \Chat
}
\]
is simply the blow up of $\Prym_{3}\times \Chat$ along the product  surface
$\emb_{\Chat}(\Chat)\times \Chat$. Thus
\[
\op{Bl}_{\Chat\times\Chat} \, (\Prym_{2}\times \Chat) = Y_{1}\times \Chat,
\]
which in turn identifies $\bigHeckehat^{\rm ab}$ with a blow up of
$Y_{1}\times \Chat$.

The morphism $\qzo^{\rm ab} : \bigHeckehat^{\rm ab} \to Y_{1}\times
\Chat$ is again the blow-up morphism, while the morphism $\pzo^{\rm
  ab} : \bigHeckehat^{\rm ab} \to Y_{0}$ now lifts the difference map
$\diff : \Prym_{3}\times \Chat \to \Prym_{2}$ to the blow-ups of the
source and target. This describes the second correspondence in diagram
\eqref{eq:abHecke}.

Observe also that the two maps $(\diff,\op{id})$ and
$(\summ,\op{id})$ not only identify $\Prym_{3}\times \Chat$ and
$\Prym_{2}\times \Chat$ but also identify the two diagrams of
correspondences in \eqref{eq:basicPrymdiag}
when the maps between all of the other
corresponding nodes in the two diagrams are taken to be identities.

Recall next that the blown up Prym varieties $Y_{0}$ and $Y_{1}$ are
equipped with finite degree $8$ morphisms $f_{0} : Y_{0} \to X_{0}$
and $f_{1} : Y_{1} \to X_{1}$ which are the minimal resolutions of the
rational maps $\Prym_{2} \dashrightarrow X_{0}$ and $\Prym_{3}
\dashrightarrow X_{1}$ both given by $L \mapsto \pi_{*}L$, for $L$ in
either $\Prym_{2}$ or $\Prym_{3}$. These maps realize $Y_{0}$ and
$Y_{1}$ as the modular spectral covers of $X_{9}$ and $X_{1}$ which
are used to define the putative Hecke eigen Higgs bundle on
$X_{0}\sqcup X_{1}$.

Following our general strategy, we will check the Hecke eigensheaf
property for this prarbolic Higgs bundle by reducing the question to
checking an abelianized Hecke eigensheaf property for the
corresponding modular spectral data. This is facilitated by the
observation that the Pryms $\Prym_{2}$, $\Prym_{3}$, the moduli
$X_{0}$, $X_{1}$ of bundles, their modular spectral covers $Y_{0}$,
$Y_{1}$, the classical Hecke correspondence $\bigHeckebar$, and the
abelianized Hecke correspondence $\bigHeckehat^{\rm ab}$ can all be
organized in a single geometric context, which, for the $(X_{1} \text{
  to } X_{0})$ direction of the Hecke property, is most compactly
recorded in the following commutative diagram
\begin{equation} \label{eq:abHecke.context} \qquad
  \begin{minipage}[c]{5in}
\[
\begin{gathered} 
\xymatrix@M+0.3pc@R-0.8pc{ & \Prym_{2}\times \Chat \ar@{-->}[dd]
  \ar[ld]_-{\summ} \ar[rd]^-{\mathsf{id}} & & & & \bigHeckehat^{\rm ab}
  \ar[ld]_-{p^{\rm ab}} \ar[rd]^-{q^{\rm ab}} \ar[dd]_-{g}
  \ar@[thicker][llll]_-{\blo} & \\ \Prym_{3} \ar@{-->}[dd]_-(0.3){\pi_{*}}  & &
  \Prym_{2}\times \Chat \ar@{-->}[dd]^-(0.3){\pi_{*}\times \pihat}  & & Y_{1}
  \ar@/_1.5pc/@[thicker][llll]_-(0.3){\blo_{1}} \ar[dd]_-(0.3){f_{1}}
  & & Y_{0}\times \Chat \ar[dd]^-(0.3){f_{0}\times \pihat}
  \ar@/_1.5pc/@[thicker][llll]_-(0.3){\blo_{0}\times \op{id}} & \\
  & \bigHeckebar
  \ar[ld]_-{p} \ar[rd]^-{q} & & & & \bigHeckebar \ar[ld]_-{p} \ar[rd]^-{q}
  \ar@{=}[llll] & \\ X_{1} & & X_{0}\times\Cbar  & & X_{1} \ar@{=}@/^1.5pc/[llll]
  & & X_{0}\times \Cbar \ar@{=}@/^1.5pc/[llll] & }
\end{gathered}
\]
\end{minipage}
\end{equation}

\

\noindent
We used $f_{0}$, $f_{1}$, and $g$ to denote the resolutions of the
rational maps from the Pryms to the moduli of bundles and the Hecke
correspondence to maps from the modular spectral covers, and the
abelianized Hecke correspondence.  Note that the map $g :
\bigHeckehat^{\rm ab} \to \bigHeckebar$ is indeed a morphism because
the the maps $f_{1} : Y_{1} \to X_{1}$ and $f_{0}\times \pihat :
Y_{0}\times \Chat \to X_{0}\times \Cbar$ are morphisms, and also the
map $p\times q : \bigHeckebar \to X_{1}\times (X_{0}\times \Cbar)$ is
a closed embedding.

Finally note that we also have a companion diagram which compactly
records all spaces and maps needed to abelianize and check the
$(X_{0} \text{
  to } X_{1})$ direction of the Hecke property, namely

\begin{equation} \label{eq:companion.context} \qquad
  \begin{minipage}[c]{5in}
\[
\begin{gathered} 
\xymatrix@M+0.3pc@R-0.8pc{ & \Prym_{3}\times \Chat \ar@{-->}[dd]
  \ar[ld]_-{\diff} \ar[rd]^-{\mathsf{id}} & & & & \bigHeckehat^{\rm
    ab} \ar[ld]_-{\pzo^{\rm ab}} \ar[rd]^-{\qzo^{\rm ab}} \ar[dd]_-{g}
  \ar@[thicker][llll]_-{(\summ,\op{id})\circ \blo} & \\ \Prym_{2}
  \ar@{-->}[dd]_-(0.3){\pi_{*}} & & \Prym_{3}\times \Chat
  \ar@{-->}[dd]^-(0.3){\pi_{*}\times \pihat} & & Y_{0}
  \ar@/_1.5pc/@[thicker][llll]_-(0.3){\blo_{0}} \ar[dd]_-(0.3){f_{0}}
  & & Y_{1}\times \Chat \ar[dd]^-(0.3){f_{1}\times \pihat}
  \ar@/_1.5pc/@[thicker][llll]_-(0.3){\blo_{1}\times \op{id}} & \\
  & \bigHeckebar
  \ar[ld]_-{\pzo} \ar[rd]^-{\qzo} & & & & \bigHeckebar \ar[ld]_-{\pzo}
  \ar[rd]^-{\qzo} \ar@{=}[llll] & \\ X_{0} & & X_{1}\times\Cbar & & X_{0}
  \ar@{=}@/^1.5pc/[llll] & & X_{1}\times \Cbar \ar@{=}@/^1.5pc/[llll]
  & }
\end{gathered}
\]
\end{minipage}
\end{equation}

\

\subsection{The pullback of the relative dualizing sheaf}
\label{sssec:pull.rel}

For the verification of the Hecke property we need to understand the
line bundles $g^{*} \omega_{q}$ and $g^{*}\omega_{\qzo}$ on
$\bigHeckehat^{\rm ab}$, where $\omega_{q} =
\omega_{\bigHeckebar}\otimes q^{*} \omega_{X_{0}\times \Cbar}^{-1}$
and $\omega_{\qzo} = \omega_{\bigHeckebar}\otimes \qzo^{*}
\omega_{X_{1}\times \Cbar}^{-1}$ are the relative dualizing sheaf of
the maps $q$ and $\qzo$ respectively. The key thing is to understand
the pullbacks of $\omega_{\bigHeckebar}$, $\omega_{q}$, and
$\omega_{\qzo}$ under the rational map $\Prym_{2} \times \Chat
\dashrightarrow \bigHeckebar$.  Note that the base locus of this map
has codimension at least two, and so the pullback of line bundles is
well defined.

First we have the following

\begin{proposition} \label{prop:omegaHbar}
  The Picard group of $\bigHeckebar$ is
\[
\op{Pic}(\bigHeckebar) =  \op{Pic}(X_{1})\times \op{Pic}(X_{0})\times
\op{Pic}(\Cbar) \cong \mathbb{Z}\times \mathbb{Z}\times  \op{Pic}(\Cbar).   
\]
In particular there is a unique line bundle $M \in \op{Pic}(\Cbar)$ so
that
  \[
  \omega_{\bigHeckebar} = p^{*}\mathcal{O}_{X_{1}}(-1)\otimes
  q^{*}(\mathcal{O}_{X_{0}}(-2)\boxtimes M).
  \]
\end{proposition}
\begin{proof} Recall that the  map $p\times (p_{\Cbar}\circ q) :
\bigHeckebar \to X_{1}\times \Cbar$ is a $\mathbb{P}^{1}$ bundle. In
moduli terms this $\mathbb{P}^{1}$-bundle is described as follows.

The moduli space $X_{1}$ has a universal bundle $\mycal{E} \to
X_{1}\times C$, satisfying $\mycal{E}_{|\{E\}\times C} \cong E$ for
all $E \in X_{1}$ \cite{Newstead,NR}. The universal bundle $\mycal{E}$
is not quite unique\footnote{In fact in our setting there the
universal sheaf can be normalized by further requiring that
$\det(\mycal{E})_{|X_{1}\times \{\pw\}} \cong
\mathcal{O}_{X_{1}}(1)$. One can show that such a normalized universal
bundle is unique but we will not need this fact.} but is well defined
up to tensoring with a line bundle of the form
$p_{X_{1}}^{*}\mathcal{O}_{X_{1}}(k)$ for some $k$.  We also have the
normalized Poincar\'{e} line bundle $\mathsf{Poinc} \to
\op{Jac}^{0}(C)\times C$, i.e. the unique line bundle satisfying
$\mathsf{Poinc}_{|\{A\}\times C} \cong A$ for all $A$, and
$\mathsf{Poinc}_{|\op{Jac}^{0}(C)\times \{p\}} \cong
\mathcal{O}_{|\op{Jac}^{0}(C)}$. Since the curve
\[
\Cbar = \{ (A,t) \in \op{Jac}^{0}(C)\times C \, | \, A^{\otimes 2}(\pw)
\cong \mathcal{O}_{C}(t)
\}
\]
is embedded in $\op{Jac}^{0}(C)\times C$, we can consider the pulled
back line bundle $\mycal{A} := \mathsf{Poinc}_{|\Cbar}$. Then we have
\begin{equation} \label{eq:P1bundle}
\bigHeckebar = \mathbb{P}\left((\op{id}\times \sq)^{*}\mycal{E}\otimes
  p_{\Cbar}^{*}\mycal{A}\right) \longrightarrow  X_{1}\times \Cbar.
\end{equation} 
Also, recall that the fibers of the $\mathbb{P}^{1}$-bundle
$\bigHeckebar \to X_{1}\times \Cbar$ map to straight lines in $X_{0}$
under the projection $\pr_{X_{0}} = p_{X_{0}}\circ q : \bigHeckebar
\to X_{0}$, and so the line bundle
$\pr_{X_{0}}^{*}\mathcal{O}_{X_{0}}(1)$ has degree one on the fibers
of $\bigHeckebar \to X_{1}\times \Cbar$.

Suppose $\xi$ is a line bundle $\bigHeckebar$. If $\xi$ has degree $k$
on the fibers of $\bigHeckebar \to X_{1}\times \Cbar$, then
$\xi\otimes \pr_{X_{0}}^{*}\mathcal{O}_{X_{0}}(-k)$ is trivial on the
fibers of $\bigHeckebar \to X_{1}\times \Cbar$, and so by cohomology
and base change
\begin{equation} \label{eq:firstpush}
(\bigHeckebar \to X_{1}\times \Cbar)_{*}(\xi\otimes
  \pr_{X_{0}}^{*}\mathcal{O}_{X_{0}}(-k))
\end{equation}
is a line bundle on $X_{1}\times \Cbar$.  Since $\op{Pic}(X_{1}) =
\mathbb{Z}$ and is generated by $\mathcal{O}_{X_{1}}(1)$ it follows
that the restriction of the line bundle \eqref{eq:firstpush} to any
slice $X_{1}\times \{\bar{t}\}$, $\bar{t} \in \Cbar$ is isomorphic to
$\mathcal{O}_{X_{1}}(r)$ for some fixed integer $r \in \mathbb{Z}$,
independent of $\bar{t}$. Again by cohomology and base change we
conclude that the further pushforward of $(\bigHeckebar \to X_{1}\times
\Cbar)_{*}(\xi\otimes \pr_{X_{0}}^{*}\mathcal{O}_{X_{0}}(-k))\otimes
p^{*}\mathcal{O}_{X_{1}}(-r)$ under the map $p_{\Cbar} : X_{1}\times
\Cbar \to \Cbar$ is some line bundle $\alpha$ on $\Cbar$. But now the
projection formula implies that
\[
\xi = \pr_{X_{1}}^{*}\mathcal{O}_{X_{1}}(r)\otimes
\pr_{X_{0}}^{*}\mathcal{O}_{X_{0}}(k)\otimes \pr_{\Cbar}^{*}\alpha,
\]
where $\pr_{X_{1}} = p$, $\pr_{X_{0}} = p_{X_{0}}\circ q$, and
$\pr_{\Cbar} = p_{\Cbar}\circ q$.

We can now apply this reasoning to $\omega_{\bigHeckebar}$.
We will have
\[
\omega_{\bigHeckebar} = \pr_{X_{1}}^{*}\mathcal{O}_{X_{1}}(r)\otimes
\pr_{X_{0}}^{*}\mathcal{O}_{X_{0}}(k)\otimes \pr_{\Cbar}^{*}M.
\]
Since the restriction of the dualizing sheaf to a smooth fiber of a
map is the dualizing sheaf of that fiber, and since the fibers of
$\bigHeckebar \to X_{1}\times \Cbar$ map to lines in $X_{0}$ we
conclude that $k = -2$. Also,  a smooth fiber of the map
$\bigHeckebar \to X_{0}\times \Cbar$ maps to a conic inside $X_{1}
\subset \mathbb{P}^{5}$, and so we conclude that $r = -1$. Therefore 
\[
\omega_{\bigHeckebar} = \pr_{X_{1}}^{*}\mathcal{O}_{X_{1}}(-1)\otimes
\pr_{X_{0}}^{*}\mathcal{O}_{X_{0}}(-2)\otimes \pr_{\Cbar}^{*}M,
\]
for some line bundle $M \in \op{Pic}(\Cbar)$. This completes the proof
of the proposition
\end{proof}

\

\noindent
Let now
$\omega_{q} = \omega_{\bigHeckebar}\otimes q^{*}\omega_{X_{0}\times
  \Cbar}^{-1}$ denote the relative dualizing sheaf for the map $q$.
Then the previous proposition gives
\[
\omega_{q}  = p^{*} \mathcal{O}_{X_{1}}(-1) \otimes q^{*}\left(
\mathcal{O}_{X_{0}}(2)\boxtimes
\left( M\otimes \omega_{\Cbar}^{-1}\right)\right), 
\]
and therefore
\begin{equation} \label{eq:g*omegaqprelim}
  \begin{minipage}[c]{5in}
\[
\begin{aligned}
\left( \Prym_{2}\times \Chat \dashrightarrow
  \bigHeckebar\right)^{*}\omega_{q} & =  \\
& \hspace{-4pc} \left(\summ^{*}\left(\Prym_{3}
  \dashrightarrow X_{1}\right)^{*} \mathcal{O}_{X_{1}}(-1)\right)
\otimes
\left( \left(\Prym_{2}
  \dashrightarrow X_{0}\right)^{*} \mathcal{O}_{X_{0}}(2) \boxtimes
\pihat^{*}\left( M\otimes\omega_{\Cbar}^{-1} \right) \right).
\end{aligned}
\]
\end{minipage}
\end{equation}
Similarly, if $\omega_{\qzo} = \omega_{\bigHeckebar}\otimes
\qzo^{*}\omega_{X_{1}\times \Cbar}^{-1}$ denotes the relative dualizing
sheaf for the map $\qzo$, then the formula for $\omega_{\bigHeckebar}$ implies
\[
\omega_{\qzo}  = \pzo^{*} \mathcal{O}_{X_{0}}(-2) \otimes \qzo^{*}\left(
\mathcal{O}_{X_{1}}(1)\boxtimes
\left( M\otimes \omega_{\Cbar}^{-1}\right)\right), 
\]
and therefore
\begin{equation} \label{eq:g*omegaqprelim0to1}
  \begin{minipage}[c]{5in}
\[
\begin{aligned}
\left( \Prym_{3}\times \Chat \dashrightarrow
  \bigHeckebar\right)^{*}\omega_{\qzo} & =  \\
& \hspace{-4pc} \left(\diff^{*}\left(\Prym_{2}
  \dashrightarrow X_{0}\right)^{*} \mathcal{O}_{X_{0}}(-2)\right)
\otimes
\left( \left(\Prym_{3}
  \dashrightarrow X_{1}\right)^{*} \mathcal{O}_{X_{1}}(1) \boxtimes
\pihat^{*}\left( M\otimes\omega_{\Cbar}^{-1} \right) \right).
\end{aligned}
\]
\end{minipage}
\end{equation}

\

\bigskip

\noindent
Next we will compute the line bundle $M \in \op{Pic}(\Cbar)$ appearing
in Proposition~\ref{prop:omegaHbar} and in the formulas
\eqref{eq:g*omegaqprelim} and \eqref{eq:g*omegaqprelim0to1}.

\

\begin{lemma} \label{lem:M}
$M = \sq^{*}\omega_{C}(\pw)$.
\end{lemma}

\begin{proof} Let $E \in X_{1}$ be a general point. Consider the
  surface  $P_{E} := p^{-1}(E) \subset \bigHeckebar$. Then $P_{E}$ is
  a smooth geometrically ruled surface over $\Cbar$ and in fact from the
  $\mathbb{P}^{1}$-bundle description \eqref{eq:P1bundle} of $\bigHeckebar$ we get
\[
  P_{E} = \mathbb{P}\left((\op{id}\times
    \sq)^{*}\mycal{E}_{|\{E\}\times \Cbar}\otimes \mycal{A}\right) =
  \mathbb{P}(\sq^{*}E\otimes \mycal{A}) 
  \cong \mathbb{P}(\sq^{*}E).
\]
The map
\[
  \xymatrix@-1pc{
    P_{E} \ar[rr]^-{q_{|P_{E}}} \ar[rd] & & X_{0}\times \Cbar \ar[ld] \\
    & \Cbar &
  }
\]
embeds $P_{E}$ into $X_{0}\times \Cbar$ by embedding each fiber of
$P_{E} \to \Cbar$ as a straight line in $X_{0} \cong
\mathbb{P}^{3}$. Also, since $P_{E} = p^{-1}(E)$ is a fiber of a map, 
we have  $\omega_{\bigHeckebar|P_{E}} \cong \omega_{P_{E}}$ and so
by Proposition~\ref{prop:omegaHbar} we have
\[
\omega_{P_{E}} \cong \omega_{\bigHeckebar|P_{E}}  = \left(
  q_{|P_{E}}\right)^{*} \left( \mathcal{O}_{X_{0}}(-2)\boxtimes M\right).
\]
So we can compute $M$ by computing the restrictions of
$\pr_{X_{0}}^{*} \mathcal{O}_{X_{0}}(-2)$ and $\omega_{P_{E}}$ to some
section  $\zeta : \Cbar \to P_{E}$ of the ruled surface $P_{E} \to
\Cbar$.

Any line subundle $L \subset E$ of degree zero will give rise to such
a section. There are four such sububundles  for a generic $E$ and we
can choose any one of them.  Let $L \subset E$ be a line subbundle of
degree zero. Then we have a short exact sequence
\[
\xymatrix@1@M+0.5pc{ 0 \ar[r] & L \ar[r] & E \ar[r] & L^{\vee}(\pw) \ar[r]
  & 0,}
\]
and $L$ gives a section of $\mathbb{P}(E) \to C$ which pulls back to a
section $\zeta : \Cbar \to P_{E}$ of $P_{E} \to \Cbar$. Note that for
any point $(A,t) \in \Cbar$ the fiber of $P_{E} \to \Cbar$ over
$(A,t)$ is canonically $\mathbb{P}(E_{t}\otimes A_{t})$ and the value
of the section $\zeta$ at $(A,t)$ is given by the line $L_{t}\otimes
A_{t}$, i.e.
\[
\zeta(A,t) = \left[ L_{t}\otimes A_{t}\right] \in \mathbb{P}(E_{t}\otimes A_{t}).
\]
The bundle $E\otimes A$ fits in a short exact sequence
\[
\xymatrix@M=0.5pc{
0 \ar[r] & L\otimes A \ar[r] & E\otimes A \ar[r] & L^{\vee}\otimes A
(\pw) \ar[r] \ar@{=}[d] 
& 0. \\
& & & (L\otimes A)^{\vee}(t) & 
}
\]
By definition the map $P_{E} \subset \bigHeckebar
\stackrel{\pr_{X_{0}}}{\longrightarrow} X_{0}$ sends the point
$\zeta(A,t)$ to the down Hecke transform of $E\otimes A$ centered at
the line $(L\otimes A)_{t} \subset (E\otimes A)_{t}$, i.e. to the
locally free sheaf $\widetilde{E}$ defined by the commutative diagram
\[
 \xymatrix@M+0.5pc@R-1pc@C-1pc{
    & & 0 & 0 & \\
& & ((L\otimes
 A)^{\vee}(t))_{t} \ar@{=}[r] \ar[u]  &  ((L\otimes
 A)^{\vee}(t))_{t} \ar[u] &  \\
 0 \ar[r] & L\otimes A \ar[r] & E\otimes A \ar[r] \ar[u] &  (L\otimes
 A)^{\vee}(t) \ar[r] \ar[u] & 0 \\
0 \ar[r] & L\otimes A \ar[r] \ar@{=}[u] & \widetilde{E}  \ar[r]
\ar[u] &  (L\otimes
A)^{\vee}\ar[r] \ar[u] & 0 \\
& & 0 \ar[u] & 0 \ar[u] &
}
\]
The third row of this diagram exhibits $\widetilde{E}$ as an extension
of  the degree zero line bundle $(L\otimes A)^{\vee}$ by $L\otimes A$
and so the bundle $\widetilde{E}$  goes to the  point
corresponding to the $S$-equivalence class $\left[ (L\otimes A)\oplus
  (L\otimes A)^{\vee}\right] \in \Kum \subset X_{0}$. 

Hence the map $\pr_{X_{0}}\circ \zeta : \Cbar \to X_{0}$  factors as
\begin{equation} \label{eq:CbarsecX0}
\xymatrix@M+0.5pc@-1pc{
\Cbar \ar[rr]^-{\pr_{X_{0}}\circ \zeta} \ar@{^{(}->}[d] & & X_{0} \\
\op{Jac}^{0}(C) \ar[rr] & & \Kum \ar@{^{(}->}[u]
}
\end{equation}
Here the left vertical map $ \Cbar \hookrightarrow \op{Jac}^{0}(C)$
is given by  $(A,t) \mapsto A\otimes L$, i.e. is the translation by
$L$ of the defining embedding of $\Cbar$ in $\op{Jac}^{0}(C)$.   The
bottom horizontal map $\op{Jac}^{0}(C) \to \Kum$ is the quotient of
$\op{Jac}^{0}(C)$ by $(-1)$, and  the composition  $\op{Jac}^{0}(C) \to
\Kum \to X_{0} \cong  \mathbb{P}^{3}$ is given by the linear
system
$|2\Theta_{\pw}|$, with $\Theta_{\pw}$ being the theta divisor
\[
\Theta_{\pw} = \left\{\left. \alpha \in \op{Jac}^{0}(C) \ \right| \
  h^{0}(C,\alpha(\pw)) \geq 1 \, \right\}.
\]

\begin{remark} \label{rem:self.duality} Strictly speaking the moduli
  space $X_{0}$ is isomorphic to
  $\mathbb{P}(H^{0}(\op{Jac}^{0}(C),\mathcal{O}(2\Theta_{\pw})))$,
  rather than
  $\mathbb{P}(H^{0}(\op{Jac}^{0}(C),\mathcal{O}(2\Theta_{\pw}))^{\vee})$. Indeed,
  by \cite{NR} the identification
  \[
    X_{0} =
    \mathbb{P}(H^{0}(\op{Jac}^{0}(C),\mathcal{O}(2\Theta_{\pw})))
  \]
  is
  given by the map $V \mapsto D_{V} \in |2\Theta_{\pw}|$, where
  $D_{V} = \{ \alpha \in \op{Jac}^{0}(C) \ | \ h^{0}(C,V\otimes
  \alpha(\pw)) \, \}$. To conclude that the map
  $\op{Jac}^{0}(C) \to X_{0}$ is indeed given by the linear system
  $|2\Theta_{\pw}|$ we use the classical fact 
  \cite{Hudson,Dolgachev200,Keum} 
  that the embedding of the Kummer
  surface in $\mathbb{P}^{3}$ is projectively self-dual.  
\end{remark}

\

\noindent
Suppose $x \in C$ is a fixed  point, and
  $\mathsf{AJ}_{x} : C \to \op{Jac}^{0}(C)$,
  $t \mapsto \mathcal{O}(t-x)$ is the $x$-based \linebreak Abel-Jacobi
  map. Suppose $\kappa \in \mathsf{Spin}(C)$ is a theta characteristic
  and let \linebreak  $\Theta_{\kappa} = \{ \alpha \in \op{Jac}^{0}(C) \, | \,
  h^{0}(C,\alpha\otimes \kappa) \geq 1\, \}$ be the associated theta
  divisor on $\op{Jac}^{0}(C)$. By Rieman's theorem \cite{GH} we have
  that  $\mathsf{AJ}_{x}^{*}\mathcal{O}(\Theta_{\kappa}) = \kappa(x)$.
In particular, since $\mathcal{O}_{C}(\pw)$ is a theta
characteristic we get that 
\[
 \mathsf{AJ}_{\pw}^{*}\mathcal{O}(\Theta_{\pw})
 \cong \mathcal{O}_{C}(2\pw) = \omega_{C}.
\]
If we denote the natural embedding of $\Cbar$ in $\op{Jac}^{0}(C)$ by
$\emb_{\Cbar} : \Cbar \hookrightarrow \op{Jac}^{0}(C)$,
$(A,t) \mapsto A$, then we can write the left vertical map in the
diagram \eqref{eq:CbarsecX0} as the composition \linebreak
$\trans_{L}\circ \emb_{\Cbar} : \Cbar \hookrightarrow
\op{Jac}^{0}(C)$, where for a line bundle $L$ on $C$ we use
\[
\trans_{L} := L\otimes (-) : \op{Jac}^{k}(C) \to \op{Jac}^{k + \deg
  L}(C)
\]
to denote the map of tensoring by $L$. With this notation we
now have a commutative diagram
\begin{equation} \label{eq:Cbarpullback}
\xymatrix@M+0.5pc@-0.5pc{
\Cbar \ar[rr]^-{\sq} \ar@{^{(}->}[d]_-{\trans_{L}\circ \emb_{\Cbar}} &
& C \ar@{^{(}->}[d]^-{\trans_{L^{\otimes 2}}\circ \mathsf{AJ}_{\pw}}  \\
\op{Jac}^{0}(C) \ar[rr]_-{\mathsf{mult}_{2}}  & & \op{Jac}^{0}(C) 
}
\end{equation}
But $\mathsf{mult}_{2}^{*} \mathcal{O}(\Theta_{\pw}) =
\mathcal{O}(4\Theta_{p})$, and so we get
\[
  \begin{aligned}
(\trans_{L}\circ \emb_{\Cbar})^{*} \mathcal{O}(4\Theta_{p}) & =
\sq^{*}\circ  \mathsf{AJ}_{\pw}^{*}\circ \trans_{L^{\otimes 2}}^{*}
\mathcal{O}(\Theta_{\pw}) \\
& = \sq^{*}\left( L^{\otimes -2}\otimes
  \mathsf{AJ}_{\pw}^{*}\mathcal{O}(\Theta_{\pw}) \right) \\
& = \sq^{*}\left( L^{\otimes -2}(2\pw)\right).
\end{aligned}
\]
But by \eqref{eq:CbarsecX0} we have
\[
(\pr_{X_{0}}\circ \zeta)^{*}\mathcal{O}_{X_{0}}(-2) = (\trans_{L}\circ
\emb_{\Cbar})^{*} \mathcal{O}(- 4\Theta_{p}),
\]
and hence
\begin{equation} \label{eq:pullfromX0}
(\pr_{X_{0}}\circ \zeta)^{*}\mathcal{O}_{X_{0}}(-2) = \sq^{*}\left(
  L^{\otimes  2}(-2\pw)\right).
\end{equation}
Next we need to compute $\zeta^{*}\omega_{P_{E}}$. By adjunction we have
\[
  \zeta^{*}\omega_{P_{E}} = \omega_{\Cbar}\otimes
  \zeta^{*}\mathcal{O}_{P_{E}}(-\zeta).
\]
Using the fact that $\zeta : \Cbar \to P_{E} = \mathbb{P}(\sq^{*}E)$ corresponds to the
subbundle $\sq^{*}L$  we  compute
\[
\begin{aligned}
  \zeta^{*}\mathcal{O}_{P_{E}}(\zeta) & = N_{\zeta/P_{E}} \\
  & = \underline{Hom}_{\mathcal{O}_{\Cbar}}\left( \sq^{*}L,
    \sq^{*}(E/L)\right) \\
  & = \sq^{*}(L^{\vee}\otimes L^{\vee}(\pw)) \\
  & = \sq^{*}\left(L^{\otimes -2}(\pw)\right).
\end{aligned}
\]
Hence
\begin{equation} \label{eq:pullomegaPE}
  \zeta^{*}\omega_{P_{E}} = \omega_{\Cbar}\otimes
  \sq^{*}\left(L^{\otimes 2}(-\pw)\right) =
  \left(\sq^{*}\mathcal{O}_{C}(2\pw)\right) \otimes
  \sq^{*}\left(L^{\otimes 2}(-\pw)\right) = \sq^{*}\left(L^{\otimes 2}(\pw)\right).
\end{equation}
Substituting \eqref{eq:pullfromX0} and  \eqref{eq:pullomegaPE} in the
the identity
\[
\zeta^{*}\omega_{P_{E}} = (\pr_{X_{0}}\circ
\zeta)^{*}\mathcal{O}_{X_{0}}(-2) \otimes M
\]
we get
\[
\sq^{*}\left(L^{\otimes 2}(\pw)\right) = \sq^{*}\left(L^{\otimes
    2}(-\pw)\right) \otimes M,
\]
that is $M = \sq^{*}\mathcal{O}_{C}(3\pw) = \sq^{*}\omega_{C}(\pw)$. 
\end{proof}

\

\noindent
Since $\omega_{\Cbar} = \sq^{*}\omega_{C}$ we therefore get
that
\[
\omega_{q} = p^{*}\mathcal{O}_{X_{1}}(-1) \otimes q^{*}\left(
  \mathcal{O}_{X_{0}}(2)\boxtimes \pihat^{*}\sq^{*}\mathcal{O}_{C}(\pw)\right)
\]
and so we get a slightly simpler version of \eqref{eq:g*omegaqprelim}, namely:
\begin{equation} \label{eq:g*omegaqprelimsiple}
  \begin{minipage}[c]{5in}
\[
\begin{aligned}
\left( \Prym_{2}\times \Chat \dashrightarrow
  \bigHeckebar\right)^{*}\omega_{q} & =  \\
& \hspace{-4pc} \left(\summ^{*}\left(\Prym_{3}
  \dashrightarrow X_{1}\right)^{*} \mathcal{O}_{X_{1}}(-1)\right)
\otimes
\left( \left(\Prym_{2}
  \dashrightarrow X_{0}\right)^{*} \mathcal{O}_{X_{0}}(2) \boxtimes
\pihat^{*}\sq^{*}\mathcal{O}_{C}(\pw)\right).
\end{aligned}
\]
\end{minipage}
\end{equation}
Finally, pulling \eqref{eq:g*omegaqprelimsiple} by $\blo :
\bigHeckehat^{\rm ab}  \to \Prym_{2}\times \Chat$ and using the commutativity of
the main diagram \eqref{eq:abHecke.context} we get
\begin{equation} \label{eq:g*omegafinal}
\boxed{
g^{*}\omega_{q} = \left(p^{\rm ab}\right)^{*}f_{1}^{*}\mathcal{O}_{X_{1}}(-1)
\otimes \left(q^{\rm ab}\right)^{*}\left(f_{0}^{*}\mathcal{O}_{X_{0}}(2) \boxtimes
\pihat^{*}\sq^{*}\mathcal{O}_{C}(\pw)\right)}
\end{equation}
This is exactly what we will need to abelianize and check the $(X_{1}
\text{ to } X_{0})$ direction of the Hecke property.  For
understanding the abelianization of the $(X_{0} \text{ to } X_{1})$
direction we also need to rewrite \eqref{eq:g*omegafinal} in terms of
the maps in the second correspondence diagram in \eqref{eq:abHecke}. This is
straighforward. Substituting $M = \sq^{*}\omega_{C}(\pw)$ in 
\eqref{eq:g*omegaqprelim0to1} gives the simplified identity 
\begin{equation} \label{eq:g*omegaqprelimsiple0to1}
  \begin{minipage}[c]{5in}
\[
\begin{aligned}
\left( \Prym_{3}\times \Chat \dashrightarrow
  \bigHeckebar\right)^{*}\omega_{\qzo} & =  \\
& \hspace{-4pc} \left(\diff^{*}\left(\Prym_{2}
  \dashrightarrow X_{0}\right)^{*} \mathcal{O}_{X_{0}}(-2)\right)
\otimes
\left( \left(\Prym_{3}
  \dashrightarrow X_{1}\right)^{*} \mathcal{O}_{X_{1}}(1) \boxtimes
\pihat^{*}\sq^{*}\mathcal{O}_{C}(\pw)\right).
\end{aligned}
\]
\end{minipage}
\end{equation}
Finally, pulling back \eqref{eq:g*omegaqprelimsiple0to1} via the map
$(\summ,\op{id})\circ \blo : \bigHeckehat \to \Prym_{3}\times \Chat$ 
and using the commutativity of
the companion diagram \eqref{eq:companion.context} we get
\begin{equation} \label{eq:g*omegafinal0to1}
\boxed{ g^{*}\omega_{\qzo} = \left(\pzo^{\rm
    ab}\right)^{*}f_{0}^{*}\mathcal{O}_{X_{0}}(-2) \otimes
  \left(\qzo^{\rm ab}\right)^{*}\left(f_{1}^{*}\mathcal{O}_{X_{1}}(1)
  \boxtimes \pihat^{*}\sq^{*}\mathcal{O}_{C}(\pw)\right)}
\end{equation}

\subsection{The Hecke property via abelianization} \label{sssec:abelianizeHecke}

Throughout this section we will fix a base point $\pwtilde \in
\Ctilde$ such that $\pi(\pwtilde) = \pw$. Now suppose $(\Ctilde
\subset T^{\vee}C,\specN)$ is the spectral data for $(E,\theta)$,
i.e.  $(E,\theta) = (\pi_{*}\specN,\pi_{*}(\lambda\otimes
(-)))$. Then $\specN \in \Prym_{2}$ and it determines two natural
line bundles
\[
  \Lprym_{0} \in \op{Pic}^{0}(\Prym_{2}), \quad \text{and} \quad
  \Lprym_{1} \in \op{Pic}^{0}(\Prym_{3}),
\]
where $\Lprym_{1} = \trans_{-\pwtilde}^{*}\Lprym_{0}$, and
$\Lprym_{0}$ is the appropriately defined Fourier-Mukai transform
of the skyscraper sheaf $\mathcal{O}_{\specN}$ on
$\Prym_{2}$. Explicitly, consider the abelian subvariety $\Prym
\subset \op{Jac}^{0}(\Ctilde)$ defined by
\[
\Prym = \left\{ \, L \in \op{Jac}^{0}(\Ctilde) \, \left| \,
    \mathsf{Nm}_{\pi}(L) = \mathcal{O}_{C} \, \right.\right\}.
\]
This abelian subvariety comes with a natural polarization $\thetaprym \to
\Prym$ defined by pulling back the canonical theta line bundle
on
$\op{Jac}^{4}(\Ctilde)$ via the natural map
\[
\Prym \subset  \op{Jac}^{0}(\Ctilde) \to  \op{Jac}^{4}(\Ctilde),
\quad L \mapsto L\otimes \mathcal{O}_{\Ctilde}(2\pi^{*}\pw).
\]
Explicitly $\thetaprym$ is defined as
$\thetaprym = \mathcal{O}_{\Prym}\left(\Xi_{2\pi^{*}\pw}\right)$,
where $\Xi_{2\pi^{*}\pw}$ is the Prym theta divisor \linebreak 
$\Xi_{2\pi^{*}\pw} := \left\{ L \in \Prym \ \left|
\ h^{0}(\Ctilde, L(2\pi^{*}\pw)) \geq 1  \,\right.\right\}$.
The line bundle $\thetaprym$ has polarization type $(1,2,2)$ and the
kernel of the associated polarization homomorphism is isomorphic to
$\mathsf{J}[2]$ \cite{BNR}.  Concretely, the polarization homomorphism
$\phi_{\thetaprym}$ to the dual abelain variety
$\Prym^{\vee} = \op{Pic}^{0}(\Prym)$ is defined by
\[
\phi_{\thetaprym} : \Prym \to \Prym^{\vee}, \quad L \mapsto
\trans_{L}^{*}\thetaprym\otimes \thetaprym^{-1}.
\]
It is a surjective homomorphism of abelian varieties and $\ker
\phi_{\thetaprym} = \pi^{*}\mathsf{J}[2]$.

Since $\Prym$ and $\Prym^{\vee}$ are dual abelian varieties we have a
canonical normalized Poincar\'{e} line bundle
\[
  {}^{\Prym}\mathsf{Poinc}  \longrightarrow \Prym\times \Prym^{\vee}
\]
characterized by the conditions
\[
\begin{aligned}
  {}^{\Prym}\mathsf{Poinc}_{|P\times \{\mycal{L}\}} & \cong \mycal{L},
  \ \text{for all } \mycal{L} \in \Prym^{\vee}, \ \text{and}
  \\ {}^{\Prym}\mathsf{Poinc}_{|\{\mathcal{O}_{\Ctilde}\}\times
    \Prym^{\vee}} & \cong \mathcal{O}_{\Prym^{\vee}}.
\end{aligned}
\]
Let $\mathsf{FM} : D^{b}(\Prym) \stackrel{\cong}{\longrightarrow}
D^{b}(\Prym^{\vee})$, $F \mapsto p_{2*}\left(p_{1}^{*}F\otimes
{}^{\Prym}\mathsf{Poinc}\right)$ denote the Fourier-Mukai transform
with kernel ${}^{\Prym}\mathsf{Poinc}$. For any $L \in \Prym$ we
have
\[
\mathsf{FM}(\mathcal{O}_{L}) =
{}^{\Prym}\mathsf{Poinc}_{|\{L\}\times \Prym^{\vee}}
\]
which is a line bundle of degree zero on $\Prym^{\vee}$.

With this notation we can now define the line bundles on $\Prym_{2}$,
and $\Prym_{3}$ which become part of the modular spectral data for the
Hecke eigensheaf.

\

\begin{definition} \label{def:modular.spectral.lb}
  For any $\specN \in \Prym_{2}$ define degree zero line bundles on
  $\Prym$, $\Prym_{2}$, and $\Prym_{3}$ by setting
  \[
   \begin{aligned}
\Lprym & := \phi_{\thetaprym}^{*} \mathsf{FM}\left(
  \mathcal{O}_{\specN(-\pi^{*}\pw)}\right) \  \in \
  \op{Pic}^{0}(\Prym), \\
\Lprym_{0} & :=
\trans_{\mathcal{O}_{\Ctilde}(-\pi^{*}\pw)}^{*}  \Lprym
\  \in \
  \op{Pic}^{0}(\Prym_{2}), \\
 \Lprym_{1} & :=
 \trans_{\mathcal{O}_{\Ctilde}(-\pwtilde -\pi^{*}\pw)}^{*}
 \Lprym
\  \in \
\op{Pic}^{0}(\Prym_{3}).
\end{aligned}
   \]
\end{definition}

\

\noindent
{\em\bfseries Note:} By construction The line bundles $\Lprym$,
$\Lprym_{0}$, and $\Lprym_{1}$ all depend on the Higgs bundle
$(E,\theta)$ or more precisely on the corresponding spectral data
$(\Ctilde \subset T^{\vee}C,\specN)$. We suppress the explicit
dependence on $\specN$ in the labeling of these line bundles to
avoid cluttering the notation.

We can describe these line bundles explicitly in terms of $\thetaprym$.

\begin{lemma} \label{lemma:LN}
$\Lprym =
\trans_{\specN(-\pi^{*}(\pw))}^{*}\thetaprym\otimes \thetaprym^{-1}$.
\end{lemma}
\begin{proof}
  The Poincar\'{e} line bundle ${}^{\Prym}\mathsf{Poinc}$ is built explicitly 
  starting from  the bi-extension line bundle
  \[
\mathcal{Q} \to \Prym\times \Prym, \qquad \mathcal{Q} =
\boldsymbol{m}^{*}\thetaprym\otimes p_{1}^{*}\thetaprym^{-1}\otimes
p_{2}^{*}\thetaprym^{-1},
  \]
where $\boldsymbol{m} : \Prym\times \Prym \to \Prym$ is the group
operation on the abelian variety $\Prym$.

\

\noindent
Indeed, the translation action of $\mathsf{J}[2]$ on the second copy
of $\Prym$ in $\Prym\times\Prym$ lifts to an action on
$\mathcal{Q}$. Indeed, since $\thetaprym$ is invariant under translation
by elements in $\mathsf{J}[2]$ , it follows that $\mathcal{Q}$ is
invariant under the translation action by elements of
$\mathsf{J}[2]\times \mathsf{J}[2] \subset \Prym\times\Prym$.  But the
extension class defining the theta group of $\mathcal{Q}$ is given by
\cite{Mumford-Abelian,bl} a pairing $e^{\mathcal{Q}}$
  which is the exponentiation of the first Chern class of
  $\mathcal{Q}$. By definition the first Chern class of $\mathcal{Q}$
  restricts to zero on each of the two summands
  $H_{1}(\Prym,\mathbb{Z})\oplus H_{1}(\Prym,\mathbb{Z})$ in the
  K\"{u}nneth decomposition of the first homology of $\Prym\times
  \Prym$. Hence $e^{\mathcal{Q}}$ restricts to the trivial character on
    each of the two factors in $\mathsf{J}[2]\times
    \mathsf{J}[2]$. Thus the theta group of $\mathcal{Q}$ splits over
    $\mathsf{J}[2]$ which implies that $\mathcal{Q}$ is equivariant
    under the the translation action of
    $\{\mathcal{O}_{\Ctilde}\}\times \mathsf{J}[2]$, as claimed.  Note
    that by definition the restriction of $\mathcal{Q}$ to
    $\{\mathcal{O}_{\Ctilde}\} \times \Prym$ is canonically trivial
    and we can use this trivializatio to normalize the
    $\{\mathcal{O}_{\Ctilde}\}\times \mathsf{J}[2]$-action on
    $\mathcal{Q}$ by requiring that $\mathsf{J}[2]$ acts
    tautologically on $\mathcal{O}_{\Prym} =
    \mathcal{Q}_{|\{\mathcal{O}_{\Ctilde}\} \times \Prym}$.

    \

    \noindent
    Using this normalized equivariant structure on $\mathcal{Q}$ we can descend $\mathcal{Q}$ to a biextension line bundle
    \[
\underline{\mathcal{Q}} \to \Prym\times \Prym^{\vee}
    \]
on the quotient $\Prym\times (\Prym/\mathsf{J}[2]) \cong \Prym\times
\Prym^{\vee}$. If $\mycal{L} \in \Prym^{\vee}$ we have that
\[
\underline{\mathcal{Q}}_{|\Prym\times \{\mycal{L}\}} \cong \mathcal{Q}_{|\Prym\times \{\alpha\}},
\]
where $\alpha \in \Prym$ is any point, s.t. $\phi_{\thetaprym}(\alpha)
= \mycal{L}$.  But then by the definition of $\mathcal{Q}$  we have
\[
\mathcal{Q}_{\Prym\times \{\alpha\}} =
\trans_{\alpha}^{*}\thetaprym\otimes \thetaprym^{-1} =
\phi_{\thetaprym}(\alpha) = \mycal{L}.
\]
Hence
\[
\underline{\mathcal{Q}}_{|\Prym\times \{\mycal{L}\}} \cong \mycal{L},
\quad \text{for all } \mycal{L} \in \Prym^{\vee}.
\]
Also, by our normalization of the equivariant structure,
$\underline{\mathcal{Q}}_{|\{\mathcal{O}_{Ctilde}\}\times
  \Prym^{\vee}}$ is the descent of $\mathcal{O}_{\Prym} =
\mathcal{Q}_{|\{\mathcal{O}_{\Ctilde}\}\times \Prym}$ with respect to
the tautological $\mathsf{J}[2]$-action. Hence
$\underline{\mathcal{Q}}_{|\{\mathcal{O}_{\Ctilde}\}\times
  \Prym^{\vee}} \cong \mathcal{O}_{\Prym^{\vee}}$. This shows that
\[
\underline{\mathcal{Q}} \cong {}^{\Prym}\mathsf{Poinc},
\]
and so
\[
\mathcal{Q} \cong \left(\op{id}\times
\phi_{\thetaprym}\right)^{*}\left({}^{\Prym}\mathsf{Poinc}\right).
\]
This implies that
\[
\begin{aligned}
  \Lprym & = \phi_{\thetaprym}^{*}\left({}^{\Prym}
  \mathsf{Poinc}_{|\{\specN(-\pi^{*}\pw)\}\times
    \Prym^{\vee}}\right) \\ & = \left(\left(\op{id}\times
  \phi_{\thetaprym}\right)^{*}\left({}^{\Prym}\mathsf{Poinc}\right)
  \right)_{|\{\specN(-\pi^{*}\pw)\}\times \Prym} \\ & =
  \mathcal{Q}_{|\{\specN(-\pi^{*}\pw)\}\times \Prym} \\ & =
  \trans^{*}_{\specN(-\pi^{*}\pw)}\thetaprym\times \thetaprym^{-1}.
\end{aligned}
\]
This completes the proof of the lemma.
\end{proof}

\

\noindent
Recall next that our candidate Hecke eigensheaf for the eigenvalue
$(E,\theta)$ was  a  tame parabolic Higgs bundle
$(\mycal{F}_{0,\bullet},\Phi_{0})\sqcup
(\mycal{F}_{1,\bullet},\Phi_{1})$ on $X_{0}\sqcup X_{1}$ which was
constructed from the modular spectral covers
\[
f_{0} : Y_{0} \to X_{0}, \quad \text{and} \quad f_{1} : Y_{1} \to X_{1}
\]
by setting
\[
\begin{aligned}
  \left(\mycal{F}_{0,0},\Phi_{0}\right) & =
  \left(f_{0*}\LY_{0},f_{0*}\left(\modalpha_{0}\otimes (-)\right)\right), \\
  \left(\mycal{F}_{1,0},\Phi_{1}\right) & = \left(f_{1*}\LY_{1},
  f_{1*}\left(\modalpha_{1}\otimes (-)\right)\right).
\end{aligned}
\]
Here
\begin{itemize}
\item $\modalpha_{0} : Y_{0} \to T^{\vee}_{X_{0}}(\log \Wob_{0})$ and
  $\modalpha_{1} : Y_{1} \to T^{\vee}_{X_{1}}(\log \Wob_{1})$ are the
  tautological maps (defined away from the preimage of codimension two
  loci in $X_{0}$ and $X_{1}$) from the modular spectral covers
  $Y_{0}$ and $Y_{1}$ to the logarithmic cotangent bundles of $X_{0}$
  and $X_{1}$ with poles along the (normal crossings part of the)
  wobbly divisors.
\item The modular spectral line bundles $\LY_{0}$ and $\LY_{1}$ are given by
  \[
 \begin{aligned}
  \LY_{0} & = (\blo_{0}^{*}\Lprym_{0})(\ExY_{0})\otimes
  f_{0}^{*}\mathcal{O}_{X_{0}}(2), \text{ and} \\
  \LY_{1} & = \blo_{1}^{*}\Lprym_{1}\otimes
  f_{1}^{*}\mathcal{O}_{X_{1}}(1)
 \end{aligned}
 \]
  and $\ExY_{0}$ and $\ExY_{1}$  denote the
  exceptional divisors of the blowup maps $\blo_{0} : Y_{0} \to
  \Prym_{2}$ and $\blo_{1} : Y_{1} \to \Prym_{3}$.
\item $\left(\mycal{F}_{0,0},\Phi_{0}\right)$ and
  $\left(\mycal{F}_{1,0},\Phi_{1}\right)$ are equipped with parabolic
  structures along $\Wob_{0}$ and $\Wob_{1}$ respectively.
\end{itemize}

\

\noindent
In  chapters~\ref{chapter-heckex1x0} and \ref{chapter-heckex0x1} we
used these constructions for $\left(\mycal{F}_{i,0},\Phi_{i}\right)$
together with the formula for the $L^{2}$ Dolbeault pushforward of
parabolic Higgs bundles from \cite{dirim} and its refinements proven
in Proposition~\ref{localcalc} to rewrite the Hecke condition
\begin{equation} \label{eq:X1toX0Hecke}
\begin{aligned}
q_{*}p^{*}\left(\mycal{F}_{1,\bullet},\Phi_{1}\right) & =
\left(\mycal{F}_{0,\bullet},\Phi_{0}\right)\boxtimes
\sq^{*}(E,\theta), \\
\qzo_{*}\pzo^{*}\left(\mycal{F}_{0,\bullet},\Phi_{0}\right) & =
\left(\mycal{F}_{1,\bullet},\Phi_{1}\right)\boxtimes
\sq^{*}(E,\theta)
\end{aligned}
\end{equation}
as a condition on the modular spectral data $(Y_{0},\LY_{0})$ and
$(Y_{1},\LY_{1})$. Namely, chasing the maps in the basic diagram
\eqref{eq:abHecke.context} we get that the $(X_{1} \text{ to } X_{0})$
direction of the Hecke eigensheaf property \eqref{eq:X1toX0Hecke} is
equivalent to the $(Y_{1} \text{ to } Y_{0})$ abelianized Hecke
condition
\begin{equation} \label{eq:Y1toY0abHecke}
  \boxed{
    \begin{aligned}
 (p^{\rm ab})^{*}\big(\blo_{1}^{*}\Lprym_{1}\otimes &
  f_{1}^{*}\mathcal{O}_{X_{1}}(1)\big) \otimes \left(
g^{*}\omega_{q}\right)(\bigExc_{0}+\bigExc_{1}) \\
& =
  (q^{\rm ab})^{*}\left[\left(\blo_{0}^{*}\Lprym_{0}(\ExY_{0})\otimes
  f_{0}^{*}\mathcal{O}_{X_{0}}(2)\right)\boxtimes
  \sqhat^{*}\specN\right],
\end{aligned}}
\end{equation}
where \eqref{eq:Y1toY0abHecke} is understood as an isomorphism of line
bundles on $\bigHeckehat^{\rm ab}$ which holds away from the $q^{\rm
  ab}$-pullback of any codimension two subvariety in $Y_{0}\times
\Chat$. In particular it suffices to prove that
\eqref{eq:Y1toY0abHecke} holds modulo multiples of the divisor
$\bigExc_{1}$.

Similarly, by Proposition~\ref{localcalc}, the analysis of the
contribution from the singularities of the horizontal divisor to be carried
out in 
section~\ref{ssec:degreeX0X1}, and  the companion diagram
\eqref{eq:companion.context} we get that the $(X_{0} \text{ to }
X_{1})$ direction of the Hecke eigensheaf property
\eqref{eq:X1toX0Hecke} is equivalent to the $(Y_{0} \text{ to }
Y_{1})$ abelianized Hecke condition
\begin{equation} \label{eq:Y0toY1abHecke}
  \boxed{
    \begin{aligned}
 (\pzo^{\rm ab})^{*}\big(\blo_{0}^{*}\Lprym_{0}(\ExY_{0})\otimes &
  f_{0}^{*}\mathcal{O}_{X_{0}}(2)\big) \otimes
g^{*}\omega_{\qzo} \\
& =
  (\qzo^{\rm ab})^{*}\left[\left(\blo_{1}^{*}\Lprym_{1}\otimes
  f_{1}^{*}\mathcal{O}_{X_{1}}(1)\right)\boxtimes
  \sqhat^{*}\specN\right].
\end{aligned}}
\end{equation}
Again \eqref{eq:Y0toY1abHecke} should be understood as an equality of
line bundles on $\bigHeckehat^{\rm ab}$ away from the $\qzo^{\rm
  ab}$-pullback of any codimension two subvariety in $Y_{1}\times
\Chat$. In particular, it suffices to check that
\eqref{eq:Y0toY1abHecke} holds modulo integral combinations of the
components of the divisor $\bigExc_{0}$.

Substituting the formula \eqref{eq:g*omegafinal} into
\eqref{eq:Y1toY0abHecke} and \eqref{eq:g*omegafinal0to1} we see that
the $(Y_{1} \text{ to } Y_{0})$ abelianized Hecke condition
\eqref{eq:Y1toY0abHecke} reduces to checking that 
\begin{equation} \label{eq:Y1toY0withexc}
\left((p^{\rm ab})^{*}\blo_{1}^{*}\Lprym_{1}\right)(\bigExc_{0}+\bigExc_{1}) =
  (q^{\rm ab})^{*}\left[\blo_{0}^{*}\Lprym_{0}(\ExY_{0})\boxtimes
  \sqhat^{*}(\specN(-\pi^{*}\pw))
\right], \quad \text{modulo } \mathbb{Z}\cdot \bigExc_{1}.
\end{equation}
Similarly  the $(Y_{0} \text{ to } Y_{1})$ abelianized Hecke condition
\eqref{eq:g*omegafinal0to1} reduces to showing that 
\begin{equation} \label{eq:Y0toY1withexc}
(\pzo^{\rm ab})^{*}\left(\blo_{0}^{*}\Lprym_{0}(\ExY_{0})\right) =
(\qzo^{\rm ab})^{*}\left[\left(\blo_{1}^{*}\Lprym_{1}\right)\boxtimes
  \sqhat^{*}(\specN(-\pi^{*}\pw))\right], \quad \text{modulo } \sum_{\kappa \in
  \mathsf{Spin}(C)} \mathbb{Z}\cdot \bigExc_{0,\kappa}.
\end{equation}
But in section~\ref{ssec:abelinized.in.context} we saw that
\[
(q^{\rm
  ab})^{*}\left[\mathcal{O}_{Y_{0}}(\ExY_{0})\boxtimes
  \mathcal{O}_{\Chat}\right] = (q^{\rm
  ab})^{*} \mathcal{O}_{\Prym_{2}\times \Chat}(\ExY_{0}\times \Chat) =
\mathcal{O}_{\bigHeckehat^{\rm ab}}(\bigExc_{0} + \bigExc_{1}),
\]
and that
\[
(\pzo^{\rm ab})^{*}\mathcal{O}_{Y_{0}}(\ExY_{0}) =
\mathcal{O}_{\bigHeckehat^{\rm ab}}(\bigExc_{0}). 
\]
This shows that full exceptional divisor content cancels in both equations
\eqref{eq:Y1toY0withexc} and \eqref{eq:Y0toY1withexc} and so
these become equations on line bundles on $\Prym_{2}\times \Chat$ and
$\Prym_{3}\times \Chat$ respectively. Concretely we are reduced to
checking that
\begin{subequations}
  \begin{equation} \label{eq:eqonP2xChat} \vspace{-1pc}
    \summ^{*} \Lprym_{1} = \Lprym_{0}\boxtimes \sqhat^{*}(\specN(-\pi^{*}\pw))
    \text{ in } \op{Pic}(\Prym_{2}\times
  \Chat).
  \end{equation} 
\begin{equation} \label{eq:eqonP3xChat} 
    \diff^{*} \Lprym_{0} = \Lprym_{1}\boxtimes \sqhat^{*}(\specN(-\pi^{*}\pw))
    \text{ in } \op{Pic}(\Prym_{3}\times
  \Chat).
  \end{equation}  
\end{subequations}
Since $(\summ,\op{id})$ and $(\diff,\op{id})$ are inverse isomorphisms,
\eqref{eq:eqonP2xChat} and \eqref{eq:eqonP3xChat} are clearly
equaivalent.  Therefore, to verify both the $(Y_{1} \text{ to }
Y_{0})$ and $(Y_{0} \text{ to } Y_{1})$ abelianized Hecke conditions
we only need the following

\

\bigskip

\begin{proposition} \label{prop:Y1toY0abHecke}  For
  any $\specN \in \Prym_{2}$ the corresponding modular spectral
  line bundles $\Lprym_{0}$ and $\Lprym_{1}$ satisfy the identity
  $\summ^{*}\Lprym_{1} = \Lprym_{0}\boxtimes
  \sqhat^{*}(\specN(-\pi^{*}\pw))$ in $\op{Pic}(\Prym_{2}\times
  \Chat)$.
\end{proposition}
\begin{proof}
Let us first compute $\summ^{*}\Lprym_{1}$. For this it will be
convenient to express $\summ^{*}\Lprym_{1}$ as a pullback of a
line bundle on $\Prym\times \Prym$.  We have a commutative diagram
\[
\xymatrix@M+0.5pc@-0.5pc{ \Prym_{2}\times\Chat \ar[r]^-{\summ}
  \ar[d]_-{\left(\trans_{-\pi^{*}\pw}\right)\times \embj} & \Prym_{3}
  \ar[d]^-{\trans_{-\pwtilde - \pi^{*}\pw}} \\
\Prym\times \Prym \ar[r]_-{\boldsymbol{m}} & \Prym
}
\]
where $\embj : \Chat \to \Prym$ is given by $(A,\tilde{t})
\mapsto \pi^{*}A^{-1}(\tilde{t} - \pwtilde)$, and $\boldsymbol{m} :
\Prym\times \Prym \to \Prym$ is the group law on the Prym, i.e.
$(L_{1},L_{2}) \mapsto L_{1}\otimes L_{2}$.

Since by definition $\Lprym_{1} = \trans_{-\pwtilde - \pi^{*}\pw}^{*}
\Lprym$ we get that
\[
\summ^{*} \Lprym_{1} =
\left(\left(\trans_{-\pi^{*}\pw}\right)\times \embj\right)^{*}
\boldsymbol{m}^{*}\Lprym.
\]
Also
\[
\Lprym_{0} = \trans_{-\pi^{*}\pw}^{*}\Lprym =
\left(\left(\trans_{-\pi^{*}\pw}\right)\times \embj\right)^{*}
p_{1}^{*}\Lprym,
\]
and so it suffices to understand the line bundle
\[
\boldsymbol{m}^{*}\Lprym \otimes p_{1}^{*}\Lprym^{-1}
\]
on the abelian variety $\Prym\times \Prym$. By Lemma~\ref{lemma:LN} we
have $\Lprym = \trans_{\specN(-\pi^{*}(\pw))}^{*}\thetaprym\otimes
\thetaprym^{-1}$. To simplify notation write $\boldsymbol{a} =
\specN(-\pi^{*}\pw) \in \Prym$. Thus we would like to understand
the line bundle
\[
\boldsymbol{m}^{*}(\trans_{\boldsymbol{a}}^{*}\thetaprym\otimes
\thetaprym^{-1})\otimes
p_{1}^{*}(\trans_{\boldsymbol{a}}^{*}\thetaprym\otimes
\thetaprym^{-1}).
\]
To simplify notation consider the $\boldsymbol{a}$-translated
biextension line bundle ${}^{\boldsymbol{a}}\mathcal{Q}$ on
$\Prym\times \Prym$ defined by
\[
  {}^{\boldsymbol{a}}\mathcal{Q} = \boldsymbol{m}^{*}\trans_{\boldsymbol{a}}^{*}
  \thetaprym \otimes p_{1}^{*}\trans_{\boldsymbol{a}}^{*}
  \thetaprym^{-1} \otimes p_{2}^{*}\trans_{\boldsymbol{a}}^{*}
  \thetaprym^{-1}.
\]
With this notation we now have
\[
\begin{aligned}
\boldsymbol{m}^{*}\trans_{\boldsymbol{a}}^{*}\thetaprym\otimes
p_{1}^{*}\trans_{\boldsymbol{a}}^{*}\thetaprym^{-1} 
& = {}^{\boldsymbol{a}}\mathcal{Q} \otimes
p_{2}^{*}\trans_{\boldsymbol{a}}^{*}\thetaprym, \text{ and} \\
\boldsymbol{m}^{*}\thetaprym^{-1} \otimes 
p_{1}^{*}\thetaprym & = \mathcal{Q}^{-1}\otimes p_{2}^{*}\thetaprym^{-1}.
\end{aligned}
\]
In particular we get
\[
\boldsymbol{m}^{*}\Lprym\otimes p_{1}^{*}\Lprym^{-1} =
\left({}^{\boldsymbol{a}}\mathcal{Q}\otimes \mathcal{Q}\right) \otimes
p_{2}^{*}(\trans_{\boldsymbol{a}}^{*}\thetaprym\otimes \theta^{-1}) =
\left({}^{\boldsymbol{a}}\mathcal{Q}\otimes \mathcal{Q}\right)
\otimes p_{2}^{*}\Lprym.
\]
On the other hand,  for any $\boldsymbol{b} \in \Prym$ we have
\[
{}^{\boldsymbol{a}}\mathcal{Q}_{|\{\boldsymbol{b}\}\times \Prym} =
\trans_{\boldsymbol{a}+\boldsymbol{b}}^{*} \thetaprym \otimes
\trans_{\boldsymbol{a}}^{*}\thetaprym^{-1} =
\trans_{\boldsymbol{a}}^{*}\left(
\trans_{\boldsymbol{b}}^{*}\thetaprym \otimes \thetaprym^{-1}\right).
\]
But $\trans_{\boldsymbol{b}}^{*}\thetaprym \otimes
\thetaprym^{-1} \in \op{Pic}^{0}(\Prym)$, and so is translation invariant.
Hence
\[
{}^{\boldsymbol{a}}\mathcal{Q}_{|\{\boldsymbol{b}\}\times \Prym} \cong
\trans_{\boldsymbol{b}}^{*}\thetaprym \otimes \thetaprym^{-1} =
\mathcal{Q}_{|\{\boldsymbol{b}\}\times \Prym}.
\]
Similarly we have that ${}^{\boldsymbol{a}}\mathcal{Q}_{|\Prym\times
    \{\boldsymbol{b}\}} \cong \mathcal{Q}_{|\Prym\times
    \{\boldsymbol{b}\}}$. So by the see-saw principle we have
  ${}^{\boldsymbol{a}}\mathcal{Q} \cong \mathcal{Q}$ and hence
  \[
  \boldsymbol{m}^{*}\Lprym\otimes p_{1}^{*}\Lprym^{-1} = p_{2}^{*}\Lprym.
  \]
  This implies that
\begin{equation} \label{eq:sL1=L0sqN}
\summ^{*}\Lprym_{1}\otimes p_{1}^{*}\Lprym_{0}^{-1} =
\left(\left(\trans_{-\pi^{*}\pw}\right)\times
\embj\right)^{*}p_{2}^{*}\Lprym =
\pr_{\Chat}^{*}\embj^{*}\Lprym
\end{equation}
and the proposition reduced to the following

\

\begin{lemma}
$\embj^{*}\Lprym = \sqhat^{*}\left(\specN(-\pi^{*}\pw)\right)$.
\end{lemma}
\begin{proof}
First observe that $\embj^{*}\Lprym$ is a pullback of a line bundle on
$\Ctilde$ via the map $\sqhat : \Chat \to \Ctilde$.  Indeed, by
definition
\[
\Lprym = \phi_{\thetaprym}^{*}\left(
  {}^{\Prym}\mathsf{Poinc}_{|\{\specN(-\pi^{*}pw)\}\times \Prym^{\vee}}\right).
\]
On the other hand, the map
\[
\phi_{\thetaprym}\circ \embj : \Chat \longrightarrow \Prym^{\vee} 
\]
is given by
\[
(\phi_{\thetaprym}\circ \embj)(A,\tilde{t}) =
\phi_{\thetaprym}(\pi^{*}A^{-1}(\tilde{t} - \pwtilde)).
\]
But  $\embj : \Chat  \to \Prym$ is $\mathsf{J}[2]$ equivariant for the
Galois action of $\mathsf{J}[2]$ on the source $\Chat$ and the
translation action of $\mathsf{J}[2] = \pi^{*}\mathsf{J}[2] = \ker
\phi_{\thetaprym} \subset \Prym$. Indeed we have
\[
\embj\left(\mathfrak{a}\cdot (A,\tilde{t})\right) =
\embj(A\otimes\mathfrak{a},\tilde{t}) =
\pi^{*}(A\otimes\mathfrak{a})^{-1}(\tilde{t} - \pwtilde) =
\pi^{*}A^{-1}\otimes\pi^{*}\mathfrak{a} (\tilde{t} - \pwtilde)  =
\pi^{*}\mathfrak{a}\otimes \embj((A,\tilde{t})).
\]
Hence
\[
(\phi_{\thetaprym}\circ \embj)\left(\mathfrak{a}\cdot
  (A,\tilde{t})\right) = (\phi_{\thetaprym}\circ \embj)(A,\tilde{t}),
\]
and so the map $\phi_{\thetaprym}\circ \embj$ factors through
$\Chat/\mathsf{J}[2] = \Ctilde$. In other words we have a comutative
triangle
\[
  \xymatrix@-0.5pc@M+0.5pc{
\Chat \ar[rr]^-{\phi_{\thetaprym}\circ \embj} \ar[dr]_-{\sqhat} & &
\Prym^{\vee} \\
& \Ctilde \ar[ur]_-{\psi} &
  }
\]
for a well defiined map $\psi : \Ctilde \to \Prym^{\vee}$.

In particular
\[
  \embj^{*}\Lprym =
  \embj^{*}\phi_{\thetaprym}^{*}\left({}^{\Prym}
    \mathsf{Poinc}_{\{\specN(-\pi^{*}\pw)\}\times
      \Prym^{\vee}}\right) = \sqhat^{*}\psi^{*}
  \left({}^{\Prym}\mathsf{Poinc}_{\{\specN(-\pi^{*}\pw)\}\times
      \Prym^{\vee}}\right).
\]
So the question reduces to computing
$\psi^{*}
  \left({}^{\Prym}\mathsf{Poinc}_{\{\specN(-\pi^{*}\pw)\}\times
      \Prym^{\vee}}\right)$. This calls for a better understanding of
  the map $\psi$ which is a version of the Abel-Prym map.

  We have an addition map
  \[
\add: \op{Jac}^{0}(C)\times  \Prym \longrightarrow \op{Jac}^{0}(\Ctilde),
\quad (A,L) \mapsto \pi^{*}A\otimes L,
\]
which is a surjective homomorphism of abelian varieties with kernel
isomorphic to $\mathsf{J}[2]$, embedded in
$\op{Jac}^{0}(C)\times \Prym$ by the map
$\mathfrak{a} \mapsto (\mathfrak{a},\pi^{*}\mathfrak{a})$. Next note
that we have a natural  map
\[
\pihat\times \embj : \Chat \longrightarrow, \quad (A,\tilde{t})
\mapsto \left(A,\pi^{*}A^{-1}(\tilde{t}-\pwtilde)\right).
\]
This map is clearly $\mathsf{J}[2]$-equivariant and so the composition
$\Chat \stackrel{\pihat\times \embj}{\longrightarrow}
\op{Jac}^{0}(C)\times \Prym \stackrel{\add}{\longrightarrow}
\op{Jac}^{0}(\Ctilde)$ will factor as
\[
\xymatrix@-0.5pc@M+0.5pc{
\Chat \ar[d]_-{\sqhat} \ar[r]^-{\pihat\times \embj} &
\op{Jac}^{0}(C)\times \Prym \ar[r]^-{\add} &
\op{Jac}^{0}(\Ctilde) \\
\Ctilde \ar[rru]_{\mathsf{AJ}_{\pwtilde}} &&
}
\]
where $\mathsf{AJ}_{\pwtilde}$ is the $\pwtilde$-based Abel-Jacobi map 
$\mathsf{AJ}_{\pwtilde} : \Ctilde \to \op{Jac}^{0}(\Ctilde)$, $\tilde{t}
\mapsto \mathcal{O}_{\Ctilde}(\tilde{t} - \pwtilde)$.

Consider the theta line bundle $\thetatilde$ on $\op{Jac}^{0}(\Ctilde)$
defined by
the theta characteristic
$\mathcal{O}_{\Ctilde}(2\pi^{*}\pw) \in \mathsf{Spin}(\Ctilde)$. In
other words
\[
\thetatilde =
\mathcal{O}_{\op{Jac}^{0}(\Ctilde)}(\Thetatilde_{2\pi^{*}\pw}), \quad
\text{where} \quad
\Thetatilde_{2\pi^{*}\pw} = \left\{ L \in \op{Jac}^{0}(\Ctilde) \ \left| \
h^{0}(\Ctilde,L(2\pi^{*}\pw)) \geq 1 \, \right.\right\}.
\]
Recall that we used $\thetatilde$ to define the line bundle
$\thetaprym$ on $\Prym$, i.e. we had $\thetaprym =
\thetatilde_{|\Prym}$. Also, $\thetatilde$ is a principal polarization
on $\op{Jac}^{0}(\Ctilde)$  which defines a surjective  group homomorphism
\[
\qfrak_{\thetatilde} : \op{Jac}^{0}(\Ctilde) \to \Prym^{\vee}, \quad L \mapsto
\left(\trans_{L}^{*}\thetatilde \otimes \thetatilde^{-1} \right)_{|\Prym},
\]
which fits in the commutative diagram
\[
\xymatrix@M+0.5pc{
  \Chat \ar[d]_-{\op{id}} \ar[r]^-{\pihat\times \embj} &
  \op{Jac}^{0}(C)\times \Prym \ar[d]_-{\pr_{\Prym}}
  \ar[r]^-{\mathsf{sum}} & \op{Jac}^{0}(\Ctilde) \ar[d]^-{\qfrak_{\thetatilde}} \\
\Chat \ar[r]_{\embj} & \Prym \ar[r]_-{\phi_{\thetaprym}} & \Prym^{\vee}
}
\]
Therefore
\[
\phi_{\thetaprym}\circ \embj = \qfrak_{\thetatilde}\circ \mathsf{sum}
\circ \left(\pihat\times \embj\right) = \qfrak_{\thetatilde}\circ
\mathsf{AJ}_{\pwtilde}\circ \sqhat.
\]
This implies that the map $\psi : \Ctilde \to \Prym^{\vee}$ factors as
\[
\xymatrix@-0.5pc@+0.5pc{
\Ctilde \ar[rr]^-{\psi} \ar[rd]_-{\mathsf{AJ}_{\pwtilde}} & & \Prym^{\vee} \\
& \op{Jac}^{0}(\Ctilde) \ar[ur]_-{\qfrak_{\thetatilde}} &
}
\]
and so we need to compute pullbacks by $\mathsf{AJ}_{\pwtilde}$ and
$\qfrak_{\thetatilde}$.

But for any $\alpha \in \Prym$ and any $\mycal{A} \in \Prym^{\vee}$ we have
that the fiber of the Poincare line bundle ${}^{\Prym}\mathsf{Poinc}$
at the point $(\alpha,\mycal{A}) \in \Prym\times \Prym^{\vee}$ is
canonically isomorphic to the fiber of the line bundle $\mycal{A} \in
\Prym^{\vee} = \op{Pic}^{0}(\Prym)$ at $\alpha \in \Prym$. Hence for any $L
\in \op{Jac}^{0}(\Ctilde)$ we have equality of fibers
\[
\begin{aligned}
\left(\qfrak^{*}_{\thetatilde}\left(
     {}^{\Prym}\mathsf{Poinc}_{|\{\alpha\}\times
       \Prym^{\vee}}\right)\right))_{L} & =
     {}^{\Prym}\mathsf{Poinc}_{(\alpha,\qfrak_{\thetatilde}(L))} \\
     & =
     \left(\qfrak_{\thetatilde}(L)\right)_{\alpha} \\
     & =
     \left(\trans_{L}^{*}\thetatilde\otimes
     \thetatilde^{-1}\right)_{\alpha}.
\end{aligned}
\]
By the see-saw principle this equality of fibers implies that
we have an isomorphism
\[
\qfrak^{*}_{\thetatilde}\left(
     {}^{\Prym}\mathsf{Poinc}_{|\{\alpha\}\times
       \Prym^{\vee}}\right) \cong
     \widetilde{\mathcal{Q}}_{|\{\alpha\}\times \op{Jac}^{0}(\Ctilde)}
\]
of line bundles on $\op{Jac}^{0}(\Ctilde)$, where
$\widetilde{\mathcal{Q}} \to \op{Jac}^{0}(\Ctilde)\times
\op{Jac}^{0}(\Ctilde)$ is the biextension line bundle given by
\[
\widetilde{\mathcal{Q}} = \boldsymbol{m}^{*}\thetatilde \otimes
p_{1}^{*}\thetatilde^{-1}\otimes p_{2}^{*}\thetatilde^{-1}.
\]
By definition, the pullback of
$\widetilde{\mathcal{Q}}_{|\{\alpha\}\times \op{Jac}^{0}(\Ctilde)}$ by
the Abel-Jacobi map $\mathsf{AJ}_{\pwtilde} : \Ctilde \to
\op{Jac}^{0}(\Ctilde)$ is the line bundle on $\Ctilde$ whose fiber at
a point $\tilde{t} \in \Ctilde$ is equal to the line
$\alpha_{\tilde{t}}\otimes
\alpha_{\pwtilde}^{-1}$. Thus
\[
\mathsf{AJ}_{\pwtilde}^{*}
\left(
\widetilde{\mathcal{Q}}_{|\{\alpha\}\times
  \op{Jac}^{0}(\Ctilde)}\right) = \alpha\otimes \alpha_{\pwtilde}^{-1}
\cong \alpha.
\]
All together this gives an isomorphism
\[
\begin{aligned}
\psi^{*}\left({}^{\Prym}\mathsf{Poinc}_{|\{\specN(-\pi^{*}\pw)\}\times
  \Prym^{\vee}}\right) & = \mathsf{AJ}_{\pwtilde}^{*}\qfrak_{\thetatilde}^{*} 
\left({}^{\Prym}\mathsf{Poinc}_{|\{\specN(-\pi^{*}\pw)\}\times
  \Prym^{\vee}}\right) \\ & = \mathsf{AJ}_{\pwtilde}^{*} \left(
\widetilde{\mathcal{Q}}_{|\{\specN(-\pi^{*}\pw)\}\times
  \op{Jac}^{0}(\Ctilde)}\right) \cong \specN(-\pi^{*}\pw),
\end{aligned}
\]
and hence $\embj^{*}\Lprym = \sqhat^{*}\specN(-\pi^{*}\pw)$
as claimed. This proves the lemma.
\end{proof}

\

\noindent
The previous lemma together with the identity \eqref{eq:sL1=L0sqN} now
implies that $\summ^{*}\Lprym_{1} = \Lprym_{0}\boxtimes
\sqhat^{*}\specN(-\pi^{*}\pw)$ which completes the proof of the
proposition.
\end{proof}

\subsection{Disjointness statements}
\label{sec:disjointness}

\begin{lemma}
\label{injectionH}
The projection map $\bigHeckebar \rightarrow X_0 \times X_1 \times
\Cbar$ is an isomorphism onto its image, which is smooth.
\end{lemma}
\begin{proof}
If $a=(A,t)\in \Cbar$ and a stable bundle $E_1\in X_1$ are fixed, the
fiber of $\bigHeckebar$ over $(a,E_1)$ is the Hecke line consisting of
all the $E_0$ obtained by Hecke transformations along $1$-dimensional
quotients of $(E_1\otimes A)(t)$. This is a $\pp^1$ mapping to a line
in $\pp^3$.  We get a morphism $X_1\times \Cbar \rightarrow {\rm
  Grass}(2,4)$ to the Grassmanian of lines in $\pp^3$. Since the
source is smooth, this map is a regular function.  The pullback of the
universal fibration over ${\rm Grass}(2,4)$ is a smooth $\pp^1$-bundle
over $X_1\times \Cbar$.  A point on one of the lines determines
uniquely the sheaf $E_0$ and the map $\zeta$, so $\bigHeckebar $ is
isomorphic to this $\pp^1$-bundle.
\end{proof}

\

\noindent
Recall that $\Prym_m$ denotes the Prym variety of degree $m$ line bundles
$L$ on $\Ctilde$ such that the norm $\mathsf{Nm}_{\Ctilde / C}(L)$ of
the associated divisor down to $C$ is linearly equivalent to $m\pw$. This
condition is equivalent to stating that $\det(\pi _{*}L) \cong
\Oo _C((m-2)\pw)$ since $\pi _{*}(\Oo _{\Ctilde}) = \omega _C ^{-1} =
\Oo _C(-2\pw)$.
In particular, $\Prym _2$ is the Prym variety of degree $2$ line
bundles $L$ on $\Ctilde$ such that ${\rm det}(\pi _{*} L) \cong \Oo
_C$, and $\Prym _3$ is the Prym variety of degree $3$ line bundles $L$
on $\Ctilde$ such that ${\rm det}(\pi _{*} L) \cong \Oo _C(\pw)$.

We have the maps $\blo_{0} : Y_0\rightarrow \Prym_2$ and $\blo_{1} :
Y_1 \rightarrow \Prym_3$, that are respectively blowing up of $16$
points or a curve $\Chat \subset \Prym_3$.
We also have $\Prym_1$ the Prym variety of degree $1$ line bundles on
$\Ctilde$ whose norm is $\Oo_{C}(\pw)$. Subtraction gives a map
$$
{\bf m}_{3,2} : \Prym_3 \times \Prym_2 \rightarrow \Prym_1.
$$
On the other hand, we have a shifted Abel-Jacobi map 
$$
{\bf j} : \Chat \rightarrow \Prym _1
$$ defined by ${\bf j}(A,\tilde{t}):= \Oo _{\Ctilde}(\tilde{t})
\otimes \pi ^{*}(A^{-1})$.  Note that if $L= \Oo
_{\Ctilde}(\tilde{t}) \otimes \pi ^{*}(A^{-1})$ and $\pi(\tilde{t}) = t$, then
$$
\mathsf{Nm}_{\Ctilde / C}(L) = A^{\otimes -2}(t) = \Oo_{C}(t -
(t-\pw)) = \Oo_{C}(\pw)
$$
as required. 

Using the product of the blowup maps  and ${\bf m}_{3,2}$
gives a composed map 
$$
\xymatrix@1@M+0.3pc@C+0.3pc{
  Y_1\times Y_0 \ar[r]^-{\blo_{1}\times \blo_{0}} & 
  \Prym _3 \times \Prym _2 \ar[r]^-{{\bf m}_{3,2}} & \Prym _1.
  }
$$
The abelianized Hecke maps to the fiber product 
$$
\bigHeckehat^{\rm ab} \rightarrow  (Y_1 \times Y_0)\times _{\Prym _1} \Chat .
$$
If we fix a point $\widehat{a} \in \Chat$, and consider the
restriction of $\bigHeckehat^{\rm ab}$ to the subvariety $(Y_{1}\times
Y_{0})\times_{\Prym_{1}} \{\widehat{a}\}$
we get a correspondence
\[
\xymatrix@M+0.25pc{
  & \bigHeckehat^{\rm ab}(\widehat{a})
  \ar[dl]_-{p^{\rm ab}}  \ar[dr]^-{q^{\rm ab}}  & \\
 Y_{1} & & Y_{0}
} 
\]
The map $\bigHeckehat^{\rm ab} (\widehat{a})\rightarrow Y_1$ is the
blow-up of $Y_1$ along $16$ lines contained in the exceptional divisor
$\ExY_{1} \cong \Chat\times \pp^{1}$ of the map $\blo_{1} : Y _{1} \to
\Prym_{3}$. Similarly the map $\bigHeckehat^{\rm ab}
(\widehat{a})\rightarrow Y_0$ is the blow-up of $Y_0$ along a curve
isomorphic to $\Chat$ which intersects transversally at a single point
each of the $16$ plane components 
$\ExY_{0,\kappa} \cong \pp^{2}$ of  the exceptional divisor $\ExY_{0} =
\sqcup_{\kappa \in \mathsf{Spin}(C)} \ExY_{0,\kappa}$.

Consider the maps $r_0 = \op{id}\times f_{0} : Y_1\times Y_0
\rightarrow Y_1\times X_0$ and $r_1 = f_{1}\times \op{id} : Y_1\times
Y_0 \rightarrow X_1\times Y_0$. With this notation we have

\

\begin{theorem}
\label{disjointness0}
Fix a point $\widehat{a} = (A,\tilde{t}) \in \Chat$ and let $\tau \widehat{a}
= (A,\tau(\tilde{t}))$ denote its conjugate under the covering
involution of the double cover $\pihat : \Chat \to \Cbar$. Assume that
$\widehat{a}$ is general in $\Chat$. Then the image in $X_0$ of the
intersection
$$
r_0\left( \bigHeckehat^{\rm ab} (\widehat{a}) \right) \cap r_0\left(
\bigHeckehat^{\rm ab} (\tau \widehat{a}) \right) \subset Y_1\times X_0
$$
has dimension at most $1$. 
\end{theorem}
\begin{proof}
  Suppose we are given two points in $\bigHeckehat^{\rm ab}$ that map
  to the points $(y_0,y_1,\widehat{a})$ and $(y'_0,y'_1,\tau \widehat{a})$
  in $(Y_{0}\times Y_{1})\times_{\Prym_{1}} \Chat$. We assume that
  these points further map to the same point $(x_0,x_1,a)$ of
  $\bigHeckebar$, and that $y_{1} = y_{1}'$. We would
like to conclude $x_0$ lies in a dimension $\leq 1$ subset of
$X_0$. Here this dimension is measured for a given fixed $\widehat{a}$.

Suppose first that the point $y_1=y'_1$ lies in the
exceptional divisor $\ExY_{1} \subset Y_1$. This part of the
exceptional locus in $\bigHeckehat^{\rm ab} (\widehat{a})$ projects to a
curve in $Y_0$ and so to a curve in $X_0$, so we may ignore this case.

We may therefore assume that $y_1=y'_1$ corresponds to a line bundle
$L$ on $\Ctilde$ whose direct image is a stable vector bundle $V=\pi
_{*}L$ of determinant $\det V = \mathcal{O}_{C}(\pw)$. The Hecke
transformations corresponding to the abelianized Hecke at $\widehat{a}$
and $\tau \widehat{a}$ are along the lines $L_{\tilde{t}}$ and
$L_{\tau(\tilde{t})}$ in $V_{t}$ corresponding to the two different
points $\tilde{t}$ and $\tau(\tilde{t})$ of the spectral
curve $\Ctilde$. Thus $L_{\tilde{t}}, L_{\tau(\tilde{t}} \subset V_{t}$
correspond to distinct points of the Hecke line $\pp^1 = \pp((V\otimes
A)_{t})$, so they are distinct points in $\bigHeckebar (a)$ and hence
distinct points in $X_0$. This contradicts the hypothesis, which
completes the proof starting with $y_1=y'_1$. 
\end{proof}

\

\begin{corollary}
\label{disjointness-cor0}
Suppose $\ell$ is a general line in $X_0$ and let $H_{\ell} =
q^{-1}(\ell) \subset \bigHeckebar(a)$ be the restriction of
$\bigHeckebar(a)$ to $\ell$. Then the restrictions $\bigHeckehat^{\rm
  ab} (\widehat{a})_{\ell}$ and $\bigHeckehat^{\rm ab} (\tau
\widehat{a})_{\ell}$ over $\ell \subset X_0$ are disjoint in the spectral
surface $\Sigma \rightarrow H_{\ell}$ that is the pullback
$$
\Sigma := H_{\ell} \times _{X_1}Y_1 .
$$
\end{corollary}
\begin{proof}
By Theorem \ref{disjointness0}, the image in $X_0$ of the intersection
of the two abelianized Hecke spaces inside $Y_1\times X_0$, is at most a
curve. Thus, a general line does not meet it, and we get the required
disjointness property.
\end{proof}

\

\noindent
For the Hecke transform in the direction from $X_0$ to $X_1$, the
following lemma is used.

\

\begin{lemma}
\label{onlyK3}
For a general point $x_0 \in \Kum$ of the Kummer in $X_0$, the image
$\bigHeckehat^{\rm ab}(\widehat{a}) \rightarrow \bigHeckebar (a)$
intersects the Hecke fiber over $x_0$ only at the point on the K3
surface where the two lines intersect, corresponding to the Hecke
transform of a polystable bundle.
\end{lemma}
\begin{proof}
Choose one of the four points $y_0\in Y_0$ over $x_0$, and let $L$ be
the corresponding line bundle on $\Ctilde$. Let $E:= \pi _{*}(L)$. We
claim that $E$ is polystable.  Since $x_{0} \in \Kum$ we know that $E$
is strictly semistable, i.e. there is a degree $0$ line bundle $U$ with
$U\hookrightarrow E$.  We would like to show that there is also a
nonzero map $U^{\vee} \rightarrow E$.

By adjunction the inclusion $U \subset E$ on $C$ corresponds to a map
$\pi^{*}U \hookrightarrow L$ of rank one locally free sheaves on
$\Ctilde$. The latter map will have cokernel which is a torsion sheaf
of length two since $L$ has degree $2$. We can therefore write
$$
L = \pi ^{*}U\otimes \Oo _{\Ctilde}(\tilde{t}_1 + \tilde{t}_2).
$$ 
Let $t_1,t_2\in C$ be the image points of $\tilde{t}_1, \tilde{t}_2$. 
Let $t'_1$ and $t'_2$
be their conjugates by the hyperelliptic involution $\hi_C$. 
The norm of $L$ down to $C$ is $\omega _C$, yielding the formula
$$
U ^{\otimes 2}\cong \omega _C(-t_1-t_2).
$$
In particular, 
$$
U^{\vee} \cong U \otimes \Oo _C(t_1 + t_2 - 2\pw).
$$
Thus
$$
\pi ^{*}(U^{\vee}) \cong \pi ^{*}(U) \otimes \Oo _{\Ctilde}
(\tilde{t}_1 + \tilde{t}_2 +\tau \tilde{t}_1 + \tau \tilde{t}_2 -
2\pwtilde - 2 \tau \pwtilde ).
$$
To get a map $\pi ^{*}(U^{\vee})\rightarrow L$ we therefore need a
section of $\Oo _{\Ctilde}(2\pwtilde + 2 \tau \pwtilde -\tau
\tilde{t}_1 - \tau \tilde{t}_2)$.

Recall that $\Ctilde$ is also a hyperelliptic curve, whose
hyperelliptic involution $\sigma$ is a lift of the hyperelliptic
involution $\hi_C$ of $C$.  The fixed points of $\sigma$ are the $12$
inverse image points of the Weierstrass points of $C$.  In particular,
$2\pwtilde$ and $2\tau \pwtilde $ are linearly equivalent divisors
coming from $\pp^1$ by pullback along the hyperelliptic map $\Ctilde
\rightarrow \pp^1$.  Therefore
$$
\Oo _{\Ctilde}(2\pwtilde - \tau \tilde{t}_i) \cong \Oo
_{\Ctilde}(\sigma \tau \tilde{t}_i) \cong \Oo _{\Ctilde}(2\tau
\pwtilde - \tau \tilde{t}_i).
$$
Thus 
$$
\Oo _{\Ctilde}(2\pwtilde + 2 \tau \pwtilde -\tau \tilde{t}_1 - \tau
\tilde{t}_2) \cong \Oo _{\Ctilde}(\sigma \tau \tilde{t}_1 + \sigma
\tau \tilde{t}_2)
$$
is effective. This gives the required nonzero section so there is a
nonzero map $\pi ^{*}(U^{\vee})\rightarrow L$ corresponding by adjunction
to $U^{\vee} \rightarrow E$.  We get a map $U\oplus U^{\vee}\rightarrow E$
that, for general $U$, has to be an isomorphism.  Thus $E$ is
polystable. Now, a Hecke transformation of $E$ that becomes a stable
bundle has to be by a diagonal line, corresponding to the stated point
of $\bigHeckebar (a)$.
\end{proof}

\

\begin{theorem}
\label{disjointness1}
Fix a point $\widehat{a} = (A,\tilde{t})\in \Chat$ and its conjugate $\tau
\widehat{a} = (A,\tau(\tilde{t}))$. Assume that $\widehat{a}$ is general in
$\Chat$. Then, up to a subset whose image in $X_1$ has dimension $\leq
1$, the intersection
$$
r_1\left( \bigHeckehat^{\rm ab} (\widehat{a}) \right) \cap 
r_1\left( \bigHeckehat^{\rm ab} (\tau \widehat{a}) \right) \subset X_1\times Y_0
$$
consists of the set of points of the form $(x_1,y_0)$ such that
$y_0$ is a point of $Y_0$ lying over the Kummer variety, and $x_1$ is
a point in $X_1$ that is in the image of the K3 surface in $\bigHeckebar
(a)$.
\end{theorem}
\begin{proof}
This will follow the same lines as the proof of Theorem
\ref{disjointness0} and we keep the same notations. However, instead
of supposing that $y_1 = y'_1$, we suppose now that $y_0=y'_0$. If
$y_0$ is a point of one of the $16$ exceptional divisors
$\ExY_{0,\kappa} \subset Y_0$ lying over the trope planes, these
project down to lines in the wobbly locus of $X_1$ so as before we can
ignore this case.

If $y_0$ corresponds to a line bundle $L$ on $\Ctilde$ whose direct
image $E=\pi _{*}L$ is a stable bundle, then by the same argument as
in the proof of Theorem~\ref{disjointness0}, the two points of the abelianized
Hecke correspond to Hecke transformations at different points of the
spectral curve, so they are distinct in the Hecke curve $\pp^1 =
\pp((E\otimes A)_{t})$ that maps to a conic in $X_1$. Thus, they map
to distinct points in $X_1$ and so this case can not happen.

We are left with the case that the point $x_0 = f_{0}(y_{0})$ is a
point $x_{0} \in \Kum$ of the Kummer variety. We may assume that $x_0$
is general in the Kummer variety, as points on a divisor of the Kummer
will lead to subsets of $X_1$ of dimension $\leq 1$.  By Lemma
\ref{onlyK3}, the point $(x_1,x_0)$ on $\bigHeckebar (a)$ lies on the
Kummer K3 surface, so $x_1$ is in the image of the K3 surface in
$X_1$. We get the stated divisor as the image of the intersection of
the two abelianized Hecke pieces.
\end{proof}

\

\begin{corollary}
\label{disjointness-cor1}
Suppose $\ell$ is a general line in $X_1$ and let $H_{\ell} =
p^{-1}(\ell) \subset \bigHeckebar(a)$ be the restriction of
$\bigHeckebar(a)$ to $\ell$.  Let
$$
\Sigma := H_{\ell} \times _{X_1}Y_1 
$$
be the spectral surface  of the pullback Higgs sheaf over $H_{\ell}$.
Then the restrictions 
$\bigHeckehat^{\rm ab} (\tilde{a})_{\ell}$ and 
$\bigHeckehat^{\rm ab} (\tau \tilde{a})_{\ell}$ over $\ell \subset X_1$
are curves in $\Sigma$ that
meet at each of the double points of $\Sigma$ lying over one of the
two nodes of the curve to be denoted $\bT\cup \bN$ in
Section \ref{chapter-heckex0x1}. 
\end{corollary}
\begin{proof}
Consider the image $\KumKthree(a)\subset X_1$ of the Kummer K3 surface inside
$\bigHeckebar (a)$.  The line $\ell$ intersects $\KumKthree(a)$ in two
points. This may be seen from the discussion of Section
\ref{chapter-heckex0x1} where the two points are the images in $\ell$
of $\bPP$ and $\bQQ$, or alternatively in Remark
\ref{thirdquadric} where $K(a)$ is viewed as the intersection of $X_1$
with an additional quadric that depends on $a$.

Suppose $(x_1,y_0)$ is an intersection point in $\Sigma$ of the two
curves, that means that $x_1\in \ell$ and this point is in the
intersection of the two abelianized Hecke pieces. We have seen that
this means that $y_0$ lies over a point of the Kummer surface, and
indeed that $y_0$ corresponds to a line bundle $L \in \Prym_{@}$ such
that $\pi_{*}L \cong U\oplus U^{\vee}$ for some line bundle $U \in
\op{Jac}^{0}(C)$ with
$$
L = \pi ^{*}U \otimes \Oo _{\Ctilde}(\tilde{t}_1 + \tilde{t}_2).
$$
But $\widehat{a} =(A, \tilde{t})$ with $\pi(\tilde{t}) = t$ and
$A^{\otimes 2} = \Oo _C(t-p)$, and so  
the point $x_1$ corresponds to a bundle $V$ that fits in an exact sequence
$$
0\rightarrow U\otimes A^{-1} \rightarrow V \rightarrow U^{\vee} \otimes 
A^{-1} \otimes \Oo _C(t) \rightarrow 0 .
$$
Assuming that $x_0$ is a general point of the Kummer surface, then we
are given two sub-line bundles of degree $0$ in the same Hecke
transformation, namely $U\otimes A^{-1}$ and $U^{-1}\otimes
A^{-1}$. It means that these two agree as sub-lines in the fiber
$V_{t}$.

Suppose we now vary the point in a family $x_1(z)$ and thus the bundle
$V(z)$.  Starting with a line that is not multiple as one of the four
lines through $x_1(0)$, and extend it to a family of lines locally at
$x_1(z)$; follow this line, making the Hecke transformation at the
resulting sub-line of $V(z)_{t}$. This gives a locally well-defined
section of the Hecke correspondence that maps into the Kummer surface
in $X_0$. We see in this way that if two distinct lines through $x_1$
(neither of which is doubled in the set of four) have the property
that they agree in $V_{t}$, then the family of intersections of the
Hecke lines with the Kummer surface has two distinct branches.

This applies to our previous situation. In Section~\ref{ssec:crit} we
will see that the family of intersections of the Hecke lines with the
Kummer surface is a reducible  curve $\bT\cup \bN$ inside $H_{\ell}$. We have
said that at our intersection point $(x_1,y_0)$ we are given two
distinct sub-lines of the bundle $V$ (the point $x_1$) such that the
Hecke transformations at $a$ agree.  Our argument then says that
$(x_1,x_0)$ is a point in $H_{\ell}$ that is a double point of the
curve $\bT\cup \bN$. This shows that it is one of the two points we have
identified, over which $\Sigma$ has four ordinary double points.

We notice that there are four points $y_0$ over $x_0$, and this will
give intersection points in the four ordinary double points of
$\Sigma$ lying over the given point of $H_{\ell}$.
\end{proof}

\begin{remark}
\label{effective-disjointness1}
In the Hecke transformation from $X_0$ to $X_1$, the restriction of the
resulting rank
$16$ Higgs bundle on a general line $\ell\subset X_1$ is the direct
sum of two rank $8$ Higgs bundles isomorphic to the restrictions of
our constructed Higgs bundle to $\ell$.  The discussion of Section
\ref{chapter-heckex0x1} will show that the contributions from the two
branches of the critical locus became equivalent to direct images of
line bundles on the the disjoint unions of the two branches, via the
mechanism of the blow-up of the double points of $\Sigma$ and the
correction due to a line bundle on the exceptional locus.
\end{remark}

\subsection{Pullbacks of the wobbly locus in the abelianized Hecke}
\label{WH-pullbacks}

We have two composed maps
$$
\bigHeckehat ^{\rm ab} \stackrel{p^{\rm ab}}{\longrightarrow}
Y_1 \stackrel{f_1}{\longrightarrow} X_1
$$
and
$$
\bigHeckehat ^{\rm ab} \stackrel{\pzo ^{\rm ab}}{\longrightarrow}
Y_0 \stackrel{f_0}{\longrightarrow} X_0.
$$
We would like to describe the inverse images of the wobbly loci 
$\Wob _1$, $\Wob _0$ under these.

For $f_1: Y_1 \to X_1$ recall that 
$\ExY_1 \subset Y_1$ is the exceptional divisor of the blowup
$\blo_{1} : Y_{1} \to \Prym _{3}$ (see 
section~\ref{ssec:moduli}). Then 
\[
(f_1)^{-1}(\Wob_1) = 2 \ExY_1 + \Residual_1.
\]
Here $\Residual_1$ is a 4-sheeted cover of $\Wob_1$, because $f_1$
has degree 8, while $\ExY_1$ is a double cover of $\Wob_1$.  By
definition $\Residual_1$ parametrizes line bundles in the Prym
that push forward to wobbly bundles.

\

\noindent
Look at the map $ p^{\rm ab}: 
\bigHeckehat ^{\rm ab}  \to Y_1$. Birationally
this is the Abel-Jacobi map: $\Prym_2 \times \Chat \to \Prym _3$, sending a
pair $(L,\hat{y} = (\tilde{y}, (A,y)))$ to $L \otimes {\pi}^*A^{-1}
(\tilde{y})$.  The generic fiber of $ p^{\rm ab}$ is $\Chat$. Over
the generic point of the exceptional divisor $\ExY_1 = \Chat \times
\pp^1$, the fiber becomes the union of $\Chat$ with $16$ surfaces, each
an $\mathbb{F}_1$.

We want to pull back this divisor by the Abel-Jacobi map. Recall from
the discussion in Section~\ref{ssec:abelinized.in.context} that the
map $\bigHeckehat^{\rm ab} \rightarrow Y_1\times \Chat$ is a
blowing-up along $16$ disjoint subvarieties $\pp^1\times \Chat$ in
$\ExY_1 \times \Chat$. In other words, the variety $\bigHeckehat ^{\rm
  ab}$ is the blowup of $Y_1 \times \Chat$ along the $16$ surfaces,
all contained in $\ExY_1 \times \Chat = \pp^1 \times \Chat \times
\Chat$. The $16$ surfaces are
$$ 
\pp^1 \times \Chat \times_{\Ctilde} \Chat =
   {\sqcup}_{\kappa \in \mathsf{Spin}(C)}
\pp^1
\times \Chat \subset \pp^1 \times \Chat \times \Chat.
$$

When we pull back our divisor $(f_1)^{-1}(\Wob_1) = 2 \ExY_1 +
\Residual_1$ by $ p^{\rm ab}$, we may first just take its product with
$\Chat$ and then do the blowing-up.  Not much happens over $\Residual
_1 \times \Chat$ (we get its strict transform $\Residualhat_1$), while
the pullback of $\ExY_1 \times \Chat$ is $\bigExc_1$ plus $16$ new
components $\bigExc_{0,\kappa}$. This proves the following:

\

\begin{lemma}
\label{double-pullback-w1}
The pullback of $\Wob _1$ by the composed map 
$f_1\circ p^{\rm ab}$ is 
$$
\bigExc_1 \cup (\bigcup_{\kappa \in \mathsf{Spin}(C)} \bigExc_{0,\kappa}) \cup
\Residualhat_1. 
$$
\end{lemma}

\

\smallskip

\noindent
Turn now to 
the pullback of $\Wob_{0}$ in $Y_{0}$ and $\bigHeckehat ^{\rm ab}$. 
For the map
$f_0: Y_0 \to X_0$ and the exceptional divisor 
$\ExY_0 \subset Y_0$ we get
\[
(f_0)^{-1}(\Wob_0) = 2 \ExY_0 + \Residual_0 + K
\]
Here $K$ is the full inverse image of the Kummer surface $\Kum \subset
Y_{0}$. It is irreducible, in fact it is dominated by the second
symmetric product of $\Chat$.  

Left over, $\Residual_0$ is a $6$-sheeted cover of $\cup_{\kappa}
\trope_{\kappa}$, because $f_0$ has degree $8$, while $\ExY_0$ is birational
to $\cup_{\kappa} \trope_{\kappa}$.  In particular,
$\Residual_0$ is a union of $16$ pieces corresponding to the $16$
trope planes (we are not saying here that these pieces are necessarily
irreducible but of course that is strongly suspected). 

Consider $ \pzo ^{\rm ab}: \bigHeckehat ^{\rm ab}  \to Y_0$. 
Here $\bigHeckehat ^{\rm ab} $ is
the blowup of $Y_0 \times \Chat$ along the surface $\Chat \times \Chat
\subset Y_0 \times \Chat$. The exceptional divisor for this blowup is
$\bigExc _1$. The pullback of $\Residual_0$ is the strict transform
$\Residualhat_0$ of the divisor $\Residual_0 \times \Chat$. For
$\ExY_0 = {\sqcup}_{\kappa} \ExY_{0,\kappa}$, the pullback of each
component consists of the strict transform of $\bigExc_{0,\kappa}$.

\

\begin{lemma}
\label{double-pullback-w0}
The pullback of $\Wob _0$ by the composed map 
$f_0\circ p^{\rm ab}$ is 
$$
\left(\bigcup_{\kappa} \bigExc_{0,\kappa}\right) \cup \bigExc_1
\cup \widehat{K} \cup \Residualhat_0. 
$$
\end{lemma}
\begin{proof}
We claim that the surface $\Chat \times \Chat$ is
contained in the extra component $K \times \Chat$. 
We have a commutative diagram
$$
    \begin{array}{ccc}
      \bigHeckehat^{\rm ab} & \rightarrow & Y_{0}\times \Chat \\
      \downarrow & & \downarrow \\
      \bigHeckebar & \rightarrow & X_{0}\times \Cbar
    \end{array} .
    $$
If we take a general point $(V,\hat{y}) \in \Chat\times \Chat \subset
Y_{0}\times \Chat$ then the inverse image of $(V,\hat{y})$ in
$\bigHeckehat^{\rm ab}$ is a $\mathbb{P}^{1}$. We want to say that
this $\pp^{1}$ maps to a positive dimensional subvariety in
$\bigHeckebar$. By construction the image of this $\pp^{1} \subset
\bigHeckehat^{\rm ab}$ to $Y_{1}$ is one of the rulings of the divisor
$\Chat\times \pp^{1} = \ExY_{1} \subset Y_{1}$. In particular if we
map the $\pp^{1} \subset \bigHeckehat^{\rm ab}$ all the way down to
$X_{1}$, then the image is one of the lines in $\Wob_1$, i.e. is
positive dimensional.  But the map $\bigHeckehat^{\rm ab} \to X_{1}$
factors through $\bigHeckebar$. Hence the image of $\pp^{1} \subset
\bigHeckehat^{\rm ab}$ in $\bigHeckebar$ will be a $\pp^{1}$ as well.
Note also that the image $\pp^1 \subset \bigHeckebar$ is contained in a
fiber of the map $q : \bigHeckebar \to X_{0}\times \Cbar$. But the
only components of fibers of $q$ that map to lines in $X_{1} \subset
\pp^{5}$ are components of fiber over points of the Kummer (note that
by assumption our $\pp^1$ maps to a general point in $\Cbar$ in
particular not a preimage of a Weierstrass point of $C$).

This shows that the locus to be blown up is contained in $K\times
\Chat$.  Therefore, the pullback of $K$ in $\bigHeckehat ^{\rm ab}$ is
the union of the strict transform $\widehat{K}$, and the exceptional
divisor $\bigExc _1$. The pullbacks of the other pieces of
$(f_0)^{-1}(\Wob_0)$ yield the other pieces of the stated
decomposition.
\end{proof}

\

\begin{theorem}
\label{whyarewedoingthis}
The intersection of the two pullbacks of the wobbly loci 
up to codimension $1$ in $\bigHeckehat ^{\rm ab}$ 
consists of just the
exceptional pieces
$$
(f_0\circ p^{\rm ab})^{-1} (\Wob _0) 
\cap 
(f_1\circ p^{\rm ab})^{-1} (\Wob _1) 
=
\bigExc_1 \cup (\bigcup_{\kappa} \bigExc_{0,\kappa}).
$$
\end{theorem}
\begin{proof}
Using the previous lemmas, it means we need to look at the intersection
$$
\left[
(\bigcup_{\kappa} \bigExc_{0,\kappa}) \cup \bigExc_1
\cup \widehat{K} \cup \Residualhat_0
\right] 
\; \cap \;
\left[
\bigExc_1 \cup (\bigcup_{\kappa} \bigExc_{0,\kappa}) \cup
\Residualhat_1
\right] .
$$
We need to show that $\Residualhat_1$ does not have any irreducible
components in common with $\widehat{K} \cup \Residualhat_0$.

Choose $a=(A,t)\in \Cbar$ with a lifting $\widehat{a} = (A,\tilde{t})$
and look at the fiber over $\widehat{a}$. Here, the subvarieties have
dimension $2$ and we want to show (for general $a$) that they do not
share any $2$-dimensional components.

In view of the decomposition of the fiber $\Residualhat_0(a)$ into
$16$ pieces that move around under the monodromy operation when we
move $a$, whereas on the other hand $\Residualhat_1(a)$ can't have
more than $4$ irreducible components since it is a $4$-sheeted
covering of the irreducible $\Wob _1$, we see that $\Residualhat _0$
and $\Residualhat _1$ can not share irreducible components. Thus, we are
reduced to the question of showing that $\Residualhat _1$ and
$\widehat{K}$ do not share any components.

Now, by Lemma \ref{onlyK3}, the image of $\bigHeckehat^{\rm ab}
(\tilde{a})$ in $\bigHeckebar (a)$ only intersects the Hecke fibers
over the Kummer $\Kum$, in the intersection points of the two lines.
But this collection of points maps to the Kummer K3 surface inside
$X_1$ that isn't the same as $\Wob _1$. Thus, the intersection
$\Residualhat _1 \cap \widehat{K}$ does not have any $2$-dimensional
pieces.  This completes the proof of the theorem.
\end{proof}

\

\begin{corollary}
\label{maybeforthis}
Inside the big Hecke correspondence 
$\bigHeckebar$, the intersection of the two pullbacks of the wobbly loci with
the abelianized Hecke is, as for its $3$-dimensional pieces: 
$$
p^{-1}(\Wob _1) \cap \pzo ^{-1} (\Wob _0) \cap 
g ( \bigHeckehat ^{\rm ab} )
= g (\bigExc _0 \cup \bigExc _1). 
$$
\end{corollary}
\begin{proof}
This is the image by $g: \bigHeckehat ^{\rm ab} \rightarrow \bigHeckebar$
of the statement in the previous theorem. 
\end{proof}

\

\noindent
Recall that for each $a\in \Chat$, the piece $g (\bigExc _0)(a)$ 
in the fiber $\bigHeckebar (a)$ contracts to a lower-dimensional
subvariety (in this case the union of $16$ lines in the wobbly locus) 
under the projection to $X_1$. And similarly,
the piece $g (\bigExc _1)(a)$ contracts to a lower-dimensional
subvariety (a birational copy of the curve $\Cbar$ in $\Kum$) 
under the projection to $X_0$. 

The piece $g (\bigExc _0)(a)$ mapping to $X_0$ consists of $16$ planes
over the trope planes, which are the ramification of the horizontal
divisor in the Hecke correspondence that contribute to the
singularities of the Hecke transformed system from $X_1$ to $X_0$. The
singularities over the Kummer surface come from the degenerations of
Hecke conics to pairs of lines.

On the other hand, the piece $g (\bigExc _1)(a)$ mapping to $X_1$ is the
ramification of the inverse image of the Kummer surface in the Hecke
correspondence, giving the singularities of the Hecke transformed
system from $X_0$ to $X_1$.

In the next chapters we will be testing the Hecke transforms by
restricting over lines in the target spaces, so the ramification
locations contributing to singularities will show up as intersected
with the inverse images of the lines in the Hecke correspondence.

\section{Hecke transformation from \texorpdfstring{$X_0$}{X0} to
  \texorpdfstring{$X_1$}{X1}}
\label{chapter-heckex0x1}

We recall the standard diagrams. The big Hecke correspondence fits
into a diagram of the form
\[
  \xymatrix@M+0.5pc@-0.5pc{
    & \bigHeckebar \ar[dl]_-{\pzo} \ar[dr]^-{\qzo} & \\
    X_{0} & & X_{1}\times \Cbar .
    }
\]
Fix a point $a =(A,t)\in \Cbar$ where $A^{\otimes 2} = \Oo
_C(t-\pw )$ and consider the Hecke correspondence $\bigHeckebar
(\overline{a})$ fitting into the diagram
\[
  \xymatrix@M+0.5pc@-0.5pc{
    & \bigHeckebar (a) \ar[dl]_-{\pzo} \ar[dr]^-{\qzo} & \\
    X_{0} & & X_{1} .
    }
\]
The objective is to pull-back the constructed Higgs bundle from $X_0$
and take the higher direct image along $\qzo$ to $X_1$.

In order to simplify the measurement of the result, fix a general line
$\ell \subset X_1$. Let $\bigHeckebar (a) _{\ell}$ be the
inverse image of $\ell$ in the Hecke variety $\bigHeckebar
(a) $. This is a Hirzebruch surface $\mathbb{F}_{1}$, which
is a $\pp^1$-bundle over $\ell \cong \pp^1$.

Let $P\subset X_0$ be the image of the map from
$\bigHeckebar(a)_{\ell}$ to $X_0$.  This is a plane $P \cong \pp^2
\subset \pp^3$. The map $\pzo :\bigHeckebar(a)_{\ell} \rightarrow P$
is blowing up a point $\bPP \in P$. This point is a general point
contained in the Kummer $\Kum \subset X_0$. There is also a point
$\bQQ \in P\cap \Kum $ where $P$ is tangent to $\Kum $.  In particular
$\bPP$ and $\bQQ$ are not on the trope planes.

Let $\bT\subset \bigHeckebar (a) _{\ell}$ be the strict transform of
$\Kum \cap P$, and let $\bN\subset \bigHeckebar(a)_{\ell}$ be the
exceptional divisor. It is the unique section of the $\mathbb{F}_1$
surface that has self-intersection $-1$.

We have that $\bT\cap \bN = \bPP'$ is a single point in
$\bigHeckebar(a)_{\ell}$. The point $\bQQ$ corresponds to a unique
point $\bQQ'\in \bigHeckebar(a)_{\ell}$ where the curve $\bT$ has a
node. The distinction in notation here is that $\bPP\in P$ is a point
of the plane inside $X_0$, while $\bPP'$ is one of the points in the
exceptional divisor above $\bPP$, namely the intersection point of
$\bT$ and $\bN$; it also corresponds to the tangent direction of the
Hecke line at $\bPP$ or equivalently the tangent direction of $P\cap
\Kum$ at $\bPP$. The points $\bQQ$ and $\bQQ'$ are basically identical
since there is no blowing-up, we just use the notation $\bQQ'$ rather
than $\bQQ$ for uniformity.

\

\begin{figure}[!ht]
\begin{center}
\psfrag{Hal}[c][c][1][0]{{$\bigHeckebar(a)_{\ell} \cong \mathbb{F}_{1}$}}
\psfrag{p'}[c][c][1][0]{{$\bPP'$}}
\psfrag{q'}[c][c][1][0]{{$\bQQ'$}}
\psfrag{T}[c][c][1][0]{{$\bT$}}
\psfrag{N}[c][c][1][0]{{$\bN$}}
\psfrag{P1}[c][c][1][0]{{$\ell \cong \pp^{1}$}}
\psfrag{f1}[c][c][1][0]{{$\qzo$}}
\epsfig{file=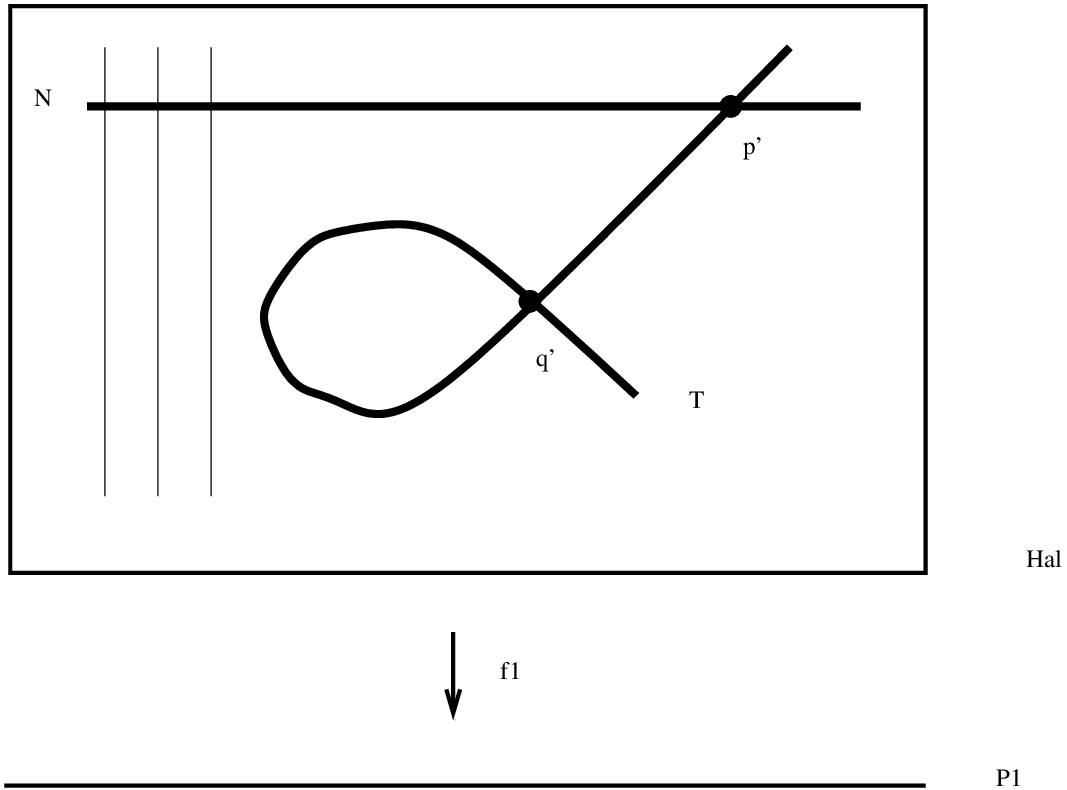,width=4in} 
\end{center}
\caption{The Hirzebruch surface $\bigHeckebar(a)_{\ell}$}\label{fig:Hal}  
\end{figure}

\

Let $Y_{1,\ell}\subset Y_1$ be the inverse image of $\ell$.  This is a
smooth curve in view of the moving properties of the family of lines
$\ell$, indeed, locally this family looks like the family of complete
intersections of two moving divisors in $X_1$.

\subsection{The abelianized Hecke as a critical locus} \label{ssec:crit}

The big abelianized Hecke fits into a diagram of the form 
 \[
\xymatrix@M+0.25pc{ & \bigHeckehat^{\rm ab} \ar[dl]_-{\pzo^{\rm ab}}
  \ar[dr]^-{\qzo^{\rm ab}} & \\ Y_{0} & & Y_{1}\times \Chat .  }
\]
The inverse image in $\Chat$ of $a\in \Cbar$ consists of two points
that we will note $\widehat{a}$ and $\tau \widehat{a}$. We obtain two
abelianized Hecke varieties $ \bigHeckehat^{\rm ab}(\widehat{a})$ and $
\bigHeckehat^{\rm ab}(\tau \widehat{a})$. Denote their disjoint union by
$$
\bigHeckehat^{\rm ab}(\widehat{a},\tau \widehat{a}):= 
\bigHeckehat^{\rm ab}(\widehat{a}) \sqcup 
\bigHeckehat^{\rm ab}(\tau \widehat{a}).
$$
The map $\bigHeckehat^{\rm ab}(\widehat{a}) \rightarrow Y_1$ is the
blow-up along $16$ lines contained in the wobbly locus, as was
discussed in Section \ref{ssec:abelinized.in.context}.

A general line $\ell \subset X_{1}$ does not meet those lines.  Thus,
if we let $\bigHeckehat^{\rm ab}(\widehat{a})_{\ell}$ denote the inverse
image of $\ell$ in $\bigHeckehat^{\rm ab}(\widehat{a})$, the projection
induces an isomorphism
$$
\bigHeckehat^{\rm ab}(\widehat{a})_{\ell}\stackrel{\cong}{\rightarrow} Y_{1,\ell}.
$$ The same holds for the other piece $\bigHeckehat^{\rm ab}(\tau
\widehat{a})_{\ell}$.

The curve $\bigHeckehat^{\rm ab}(\widehat{a})_{\ell}$ has degree $8$ over $\ell$. 
The full curve of degree $16$ over
$\ell$ is 
$$
\bigHeckehat^{\rm ab}(\widehat{a},\tau \widehat{a})_{\ell} =
\bigHeckehat^{\rm ab}(\widehat{a})_{\ell} 
\sqcup \bigHeckehat^{\rm ab}(\tau \widehat{a})_{\ell}.
$$
The $16$ exceptional divisors in $\bigHeckehat^{\rm
  ab}(\widehat{a})_{\ell}$ contract to the $16$ exceptional planes in
$Y_0$. It follows that the image of $\bigHeckehat^{\rm
  ab}(\widehat{a},\tau \widehat{a})_{\ell}$ in $Y_0$ does not meet those.

Let $Y_{0,P}:= Y_0\times _{X_0}P$, and set
$$
\Sigma := Y_{0,P} \times_P \bigHeckebar (a) _{\ell}.
$$
This is the spectral variety for the Higgs bundle
$(\mycal{F}_{0,\bullet},\Phi_{0})_{\bigHeckebar(a)_{\ell}}$ which is
$(\mycal{F}_{0,\bullet},\Phi_{0})$ pulled
back to $\bigHeckebar(a)_{\ell}$ via the map $d : \bigHeckebar(a)_{\ell} \to X_{0}$. In particular $\Sigma$  has degree $8$
over $\bigHeckebar(a)_{\ell}$. We have a lifting 
$$
\bigHeckehat^{\rm ab}(\widehat{a},\tau \widehat{a})_{\ell} \rightarrow \Sigma .
$$

\

\begin{proposition}
\label{abheckecritical}
The map $\bigHeckehat^{\rm ab}(\widehat{a},\tau \widehat{a})_{\ell}
\rightarrow \Sigma$ identifies the abelianized
Hecke as the upper critical locus (see Proposition \ref{relcrit})
$$
\bigHeckehat^{\rm ab}(\widehat{a},\tau \widehat{a})_{\ell} 
\cong 
\widetilde{{\rm Crit}}\left(\bigHeckebar(a)_{\ell}/\ell ,
(\mycal{F}_{0,\bullet},\Phi_{0}) _{\bigHeckebar (a) _{\ell}}\right)
$$
away from the points $\bPP'$ and $\bQQ'$. 
\end{proposition}
\begin{proof}
Suppose we are given a point of $\Sigma$. It corresponds to a point of
$Y_0$ together with a Hecke operation leading to a Hecke transformed
bundle which is a point of $\ell \subset X_1$.

If the point of $Y_0$ is on an exceptional divisor, it means that the
point of $\Sigma$ is a ramification point of $\Sigma$ over the
horizontal divisor in $\bigHeckebar(a)_{\ell}$. The residues of the
Higgs field $\Phi_{0}$ are nilpotent, so their Jordan types must be
constant along the parts of the horizontal divisor that are etale over
$\ell$, and this implies that ramification points of $\Sigma$ over
this subset of the horizontal divisor can not be in the relative
critical locus (we will see this argument again in the next
subsection).  The points $\bPP'$ and $\bQQ'$, as well as the
ramification points of the horizontal divisor, are contained in the
image of $\bigHeckehat^{\rm ab}(\widehat{a},\tau \widehat{a})_{\ell} $ and
also in the relative critical locus.

We claim that the upper critical locus is smooth at the ramification
points of the horizontal divisor, i.e. at points where $\bT$ ramifies
over $\ell$.  In local coordinates $x,t$ where the map to $\ell$ is
given by $t$, we can write $\bT$ as $x^2-t=0$, and $\Sigma$ is $x^2 - t
= w^2$. Thus $x,w$ are coordinates on $\Sigma$ and $t=x^2 -
w^2$. Write $\alpha$ as $adx + bdw$. The vertical direction over the
ramification point is the $w$-axis, and the fact that the Jordan form
of the residue of the Higgs field is constant along smooth points of
$\Wob _1$ (our ramification is such a smooth point) implies that
$\alpha$ is nonzero on the vertical direction, as was pointed out in
Lemma \ref{alphaatram01}.  This says that $b(0)\neq 0$. Now, the
relative differentials are obtained by setting $dt=0$ so $xdx =
wdw$. Thus $\alpha^{\rm rel}$ is written as
$$
\alpha^{\rm rel} = (a + (x/w) b)dx 
$$
and the condition $\alpha^{\rm rel}=0$ may be written as $wa + xb =
0$. As $b(0)\neq 0$, the linear term $wa(0) + xb(0)$ is nonzero, so
this defines a smooth curve which is the local branch of the relative
critical locus at this point of $\Sigma$.

This shows that the relative critical locus is smooth at ramification
points of $\bT/\ell$, so a set-theoretical identification between the
two subvarieties (neither of which has embedded points) is an
isomorphism.

Over the points $\bPP'$ and $\bQQ'$, we will see that the relative
critical locus has two branches in each of the four nodes of $\Sigma$
over one of these points, and the abelianized Hecke is isomorphic to
the normalization.

In view of the previous paragraphs, we can now assume that the point
of $Y_0$ is not on an exceptional divisor, so it corresponds to a line
bundle $L$ of degree $2$ on $\Ctilde$ whose norm to $C$ is $\omega
_C$. The corresponding vector bundle on $C$ is $E=\pi _{*}L$. Our
Hecke transformation is given by a pair $(A,x)$ where $x\in C$ and
$A^{\otimes 2} = \Oo _C(x-\pw)$. The Hecke transformation is centered at
a line $V\subset E_x$ and we let $E'$ be the kernel of the map
$E\rightarrow E_x / V$.  This is a bundle of degree $-1$ with
determinant $\Oo _C(-x)$. Then $E'\otimes A(\pw)$ has degree $1$ and
determinant $\Oo _C(\pw)$ so it is a point in $X_1$.

In turn we have a line $V'\subset E'(x)_x$ and the kernel of $E'
\rightarrow E'(x)_x / V'$ is $E$. The line $V'$ corresponds to a line
in $E'\otimes A(\pw)$ and the opposite Hecke transform involving
tensoring again with a square-root, gets us back to $E$. The Hecke
line through the point $[E]$ consists of all the bundles obtained as
kernels of $E'\rightarrow E'(x) / W$ where $W$ is another line
$W\subset E'(x)_x$. Letting $W$ be an infinitesimal deformation of
$V'$ we get a tangent vector to $X_0$ at the point $[E]$, represented
by a class $\eta \in H^1(\text{End}_{0}(E))$. We'll give an expression
for this class using an exact sequence, below.

The tautological $1$-form on $Y_{0}$ is a section of
$f_{0}^{*}T_{X_{0}}^{\vee}$ At the point $[L] \in Y_{0}$ corresponding
to the line bundle $L$ that projects to $[E] \in X_0$, $\alpha([L])
\in \left(f_{0}^{*}T^{\vee}_{X_{0}}\right)_{[L]} =
T^{\vee}_{X_{0},[E]}$ can be paired with a tangent vector of the form
$\eta$.  The condition that we are at a point of the upper critical
locus means that the value of the pairing is zero.

Our point on $Y_{0}$ actually represents a Higgs bundle over $C$,
because $Y_{0}$ was the blow-up of $\Prym_{2} \subset \Higgs _0$. The
underlying bundle is $E$ (since we are assuming that we are not over
the wobbly locus), and the Higgs field $\theta : E \to E\otimes
\omega_{C}$ is induced by the tautological $1$-form of
$\Ctilde$. In these terms, the pairing is just the pairing between
$\theta \in H^0(\text{End}_0(E)\otimes \omega _C)$ and $\eta \in
H^1(\text{End}_0(E))$. The upper critical locus condition is that this
pairing is $0$.

Let us now look more closely at the deformation class $\eta$. The
deformation of $V'$ is given by an element of $\text{Hom}(V',E'(x)_x /
V')$.  Tensor
$$
0 \rightarrow E \rightarrow E'(x) \rightarrow E'(x)_x / V' \rightarrow 0
$$
with $E^{\vee}$ to get
$$
0 \rightarrow \op{End}(E) \rightarrow E^{\vee} \otimes E'(x) \rightarrow 
E_x^{\vee}  \otimes E'(x)_x / V' \rightarrow 0.
$$ 
Note that $E_x \rightarrow V'$ so $(V')^{\vee} \rightarrow E_x^{\vee}$.
Compose with the connecting map for the previous exact sequence as
follows:
$$
(V')^{\vee} \otimes (E'(x)_x / V')
\rightarrow 
E_x^{\vee}  \otimes (E'(x)_x / V')
\rightarrow H^1(\op{End}(E)).
$$
The image of an element of $\op{Hom}(V', (E(x)_x / V'))$ is the
deformation class of $E$ generated by doing the Hecke operation back
using the infinitesimally close $W$.  The determinant of the new
bundle is the same as that of $E$, so the deformation class lies in
the trace-free part $H^1(\op{End}_0(E))$.

The exact sequence for the fiber of the Hecke transformation is
$$
0 \rightarrow E'_x / V'(-x) \rightarrow E_x \rightarrow V' \rightarrow 0.
$$
The exact sequence 
$$
0 \rightarrow \op{End}(E) \rightarrow E^{*} \otimes E'(x) \rightarrow 
E_x^{\vee}  \otimes (E'(x)_x / V') \rightarrow 0.
$$ 
has as dual,  taking $\op{Ext}^1(-,\omega _C)$:
$$
0\rightarrow E \otimes (E')^{\vee} (-x)\otimes \omega _C \rightarrow
\op{End}(E)\otimes \omega _C \rightarrow \op{Ext}^1(E_x^{\vee} \otimes (E'(x)_x /
V' ),\omega _C) \rightarrow 0.
$$
We have
$$
\op{Ext}^1(E_x^{\vee}  \otimes (E'(x)_x / V' ),\omega _C) 
= \op{Hom} (\Oo _C, E_x^{\vee}  \otimes (E'(x)_x / V' ) )^{\vee}
= \op{Hom}( (E'(x)_x / V' ), E_x). 
$$
This maps to $\op{Hom}( E'(x)_x / V' , V')$. 
Use the residue identification $\omega _C(x) _x \cong \cc$ 
and hence $V'\otimes \omega _C(x) \cong V'$, to say that 
$$
\op{Hom}( E'(x)_x / V' , V') \cong  
\op{Hom}( E'(x)_x / V' , V'\otimes \omega _C(x) )
\cong  \op{Hom}(E'_x/V'(-x),V'\otimes \omega _C).
$$
The previous map now becomes
$$
\op{Ext}^1(E_x^{\vee}  \otimes (E'(x)_x / V' ),\omega _C) 
\rightarrow 
\op{Hom}(E'_x/V'(-x),V'\otimes \omega _C).
$$
Altogether, the second map in the dual exact sequence becomes
$$
\op{End}(E)\otimes \omega _C \rightarrow
\op{Hom}(E'_x/V'(-x),V'\otimes \omega _C),
$$
and this is just the evaluation map evaluating a Higgs field on the 
subspace $E'_x/V'(-x)\subset E_x$ and projecting the answer to the quotent
$V'$ tensored with $\omega _C$. 

The Serre dual of the map $(V')^{\vee} \otimes (E'(x)_x / V') \rightarrow
H^1(\op{End}(E))$ is the action of the previous map on global sections:
$$
H^0(\op{End}(E)\otimes \omega _C ) \rightarrow
\op{Hom}(E'_x/V'(-x),V'\otimes \omega _C),
$$
To say that this vanishes for a Higgs field $\theta : E \rightarrow
E\otimes \omega _C$ is equivalent to saying that the filtration with
subspace $E'_x/V'(-x)\subset E_x$ and quotient $E_x\rightarrow V'$ is
respected by the Higgs field.

This is now the same thing as saying that our original Hecke
transformation was in the abelianized Hecke.
\end{proof}

\

\begin{lemma}
\label{notinttropeconics}
For a general $\ell$, $\bigHeckehat^{\rm ab}(\widehat{a},\tau
\widehat{a})_{\ell}$ does not intersect the inverse images in
$\bigHeckebar (a) _{\ell} $ of the $32$ points in $P$ given by the
intersection of $P$ with the trope conics.
\end{lemma}
\begin{proof}
Recall that $\bigHeckehat^{\rm ab}(\widehat{a})_{\ell} \rightarrow Y_0$ is
the blow-up along a curve isomorphic to $\Chat$, this curve
intersecting the planes of the exceptional locus in $Y_0$ in a finite
set. The trope conics are conics in these planes, so their inverse
image in $\bigHeckehat^{\rm ab}(\widehat{a})_{\ell}$ has dimension $1$.
This projects to a dimension $1$ subset of $X_1$, so a general line
$\ell$ does not intersect that. The same holds for the other piece. So,
$\bigHeckehat^{\rm ab}(\widehat{a},\tau \widehat{a})$ intersected with the
pullbacks of the trope conics, intersected with the pullback of
$\ell$, will be empty.
\end{proof}

\

\subsection{Description of the configuration}

We look more closely at the details of the configuration inside and
above $\bigHeckebar(a)_{\ell}$.  Inside $\Sigma$, there are $16$ lines
$L_{\kappa,\Sigma} \subset \Sigma$ lying over the (strict transforms
in $\bigHeckebar(a)_{\ell}$ of the) intersections $P\cap
\trope_{\kappa}$ of $P$ with the trope planes, on which
$\Sigma/\bigHeckebar(a) _{\ell}$ has a simple ramification. Lemma
\ref{notinttropeconics} says that the curve $\bigHeckehat^{\rm
  ab}(\widehat{a},\tau \widehat{a})_{\ell}$ does not meet these
$L_{\kappa,\Sigma}$.

Over general points of $\bT\cup \bN \subset \bigHeckebar(a) _{\ell}$
the ramification of $\Sigma$ consists of four sheets each of which is
a simple ramification. We also know that $\Sigma$ has four ordinary
double points over each of $\bPP'$ and $\bQQ'$.  There are probably
other singularities for example over the intersection of $P$ with the
trope conics, however the curve $\bigHeckehat^{\rm ab}(\widehat{a},\tau
\widehat{a})_{\ell}$ does not touch those.

There is a spectral line bundle $\LY_0$ on $Y_0$, which pulls back to
a line bundle $\LY_{\Sigma}$ on $\Sigma$, such that
$\Ee_{\bigHeckebar(a)_{\ell}}$ is the pushforward of $\LY_{\Sigma}$
from $\Sigma$ to $\bigHeckebar(a)_{\ell}$.  The Higgs field of
$\Ee_{\bigHeckebar(a)_{\ell}}$ is given by a holomorphic $1$-form
$\alpha$ on the smooth locus of $\Sigma$. It projects to a section of
$\Omega^1_{\bigHeckebar(a)_{\ell} / \ell}(\log D)|_{\Sigma}$ where the
divisor $D$ is the union of $\bT$, $\bN$ and the strict transforms in
$\bigHeckebar(a)_{\ell}$ of the $16$ trope lines $P\cap
\trope_{\kappa}$.

Along a point where $\bT\cup \bN$ is horizontal over $\ell$, the
residue of the Higgs field consists of a sum of four nonzero nilpotent
transformations, since the Jordan type of the residue has to stay
fixed along the divisor because it corresponds to
monodromy. Furthermore, over $\bN$ at least, the inverse image of
$\bN$ in $\Sigma$ is a union of four lines along which the map has
simple ramification. Thus, the form $\alpha$ is nonvanishing in the
transverse direction to these lines so as a section of
$\Omega^1_{\bigHeckebar(a) _{\ell}/\ell}(\log D)|_{\Sigma}$ it is
nonvanishing along these lines except at the points over $\bPP'$. On
the other hand, the curve $\bigHeckehat^{\rm ab}(\widehat{a},\tau
\widehat{a})_{\ell}$ is the vanishing locus of the section of
$\Omega^1_{\bigHeckebar(a) _{\ell}/\ell}(\log D)|_{\Sigma}$. It
follows that $\bigHeckehat^{\rm ab}(\widehat{a},\tau
\widehat{a})_{\ell}\subset \Sigma$ does not intersect the four lines over
$\bN$ except at points over $\bPP'$. In sum, the image of
$\bigHeckehat^{\rm ab}(\widehat{a},\tau \widehat{a})_{\ell}$ in
$\bigHeckebar(a)_{\ell}$ only intersects $\bN$ at $\bPP'$.  Similar
considerations show that this also holds for the horizontal part of
$\bT - \{ \bPP',\bQQ'\}$.

We would like to know how many branches of $\bigHeckehat^{\rm
  ab}(\widehat{a},\tau \widehat{a})_{\ell}$ pass through the point
$\bPP'$. We will see by the local calculations below that each double point
of $\Sigma$ over $\bPP'$ corresponds to two branches of the polar curve
of zeros of the relative Higgs field.  This says that the full curve
$\bigHeckehat^{\rm ab}(\widehat{a},\tau \widehat{a})_{\ell}$ has $2$
branches in each of the four double points, thus it has $8$ branches
over $\bPP'$. A monodromy argument says that these have to be evenly
distributed, so $\bigHeckehat^{\rm ab}(\widehat{a})_{\ell}$ has $4$
branches passing through $\bPP'$. The local argument will also show that
they intersect $\bN$ transversally. Thus, if we denote by
$[\bigHeckehat^{\rm ab}(\widehat{a})_{\ell}]$ the image of this curve in
$\bigHeckebar(a)_{\ell}$, we have
$$
[\bigHeckehat^{\rm ab}(\widehat{a})_{\ell}]\cdot \bN = 4.
$$
We will discuss later the distribution of these $4$ branches among the
$4$ double points of $\Sigma$. 

We next consider intersection numbers inside the $\mathbb{F}_1$
surface $\bigHeckebar(a)_{\ell}$.  Let $\mathsf{fib}$ denote a general
fiber of $\qzo : \bigHeckebar(a)_{\ell} \rightarrow \ell$. Then the
divisors $\mathsf{fib}$ and $\bN$ generate the Picard group of
$\bigHeckebar(a)_{\ell}$, with $\mathsf{fib}^2=0$, $\bN^2=-1$ and
$\mathsf{fib}\cdot \bN=1$.

The curve $\bigHeckehat^{\rm ab}(\widehat{a})_{\ell}$ has degree $8$ over $\ell$,
so
$$
[\bigHeckehat^{\rm ab}(\widehat{a})_{\ell}]\cdot \mathsf{fib} = 8.
$$
The intersection numbers with $\bN$ and $\mathsf{fib}$ uniquely determine
the class of $\bigHeckehat^{\rm ab}(\widehat{a})_{\ell}$ inside
$\bigHeckebar(a)_{\ell}$:
$$
[\bigHeckehat^{\rm ab}(\tilde{a})_{\ell}] \sim 8 \bN + 12 \mathsf{fib}. 
$$
On the other hand, we also have that $\bT$ has degree $3$ over
$\ell$ since the intersection with a fiber is all the points in that
line intersected the Kummer, except the point $\bPP'$ which
corresponds to the point of $\bN$ intersected that fiber.  We have
$\bT\cdot \mathsf{fib} = 3$ and $\bT\cdot \bN=1$ so
$$
\bT \sim 3\bN + 4\mathsf{fib}.
$$
Combining these we get
$$
[\bigHeckehat^{\rm ab}(\widehat{a})_{\ell}]\cdot \bT = (8\bN+12
\mathsf{fib})(3\bN+4\mathsf{fib}) = 36 + 32 - 24 = 44.
$$ We already know the intersections over the point $\bPP'$, there are
$4$ branches.  Also over $\bQQ'$ there will similarly be $4$ branches,
but $\bQQ'$ is a double point of $\bT$ (whereas the double point of
$D$ at $\bPP'$ involved both $\bT$ and $\bN$).  Thus, the intersection
counts for $8$ points over $\bQQ'$.  This leaves $32$ intersection points
over the ramification points of $\bT/\ell$. We note that there are $8$
such points, mapping to the $8$ points in $\ell$ where $\ell \subset
X_1$ intersects ${\rm Wob}_1$.  We get an intersection number of $4$
at each ramification point (a monodromy argument shows that they need
to be evenly distributed).  This will come from two tacnodes between
the curve $\bigHeckehat^{\rm ab}(\widehat{a})_{\ell}$ and the divisor
$\bT$, in accordance with our calculations for the local contribution of
the direct image at a simple ramification point of the horizontal
divisor.

The curve $\bigHeckehat^{\rm ab}(\widehat{a})_{\ell}$ has to meet the
ramified locus of $Y_0$ over $\Kum $ in two of the four sheets, as
each one will contribute a tacnode and from the above discussion there
should not be more than $2$ tacnodes. The curve $\bigHeckehat^{\rm
  ab}(\tau \widehat{a})_{\ell}$ meets the ramified locus in the other
two sheets. On the other hand, for the nodal points we'll see that
$\bigHeckehat^{\rm ab}(\widehat{a})_{\ell}$ has one branch in each sheet
and $\bigHeckehat^{\rm ab}(\tau \widehat{a})_{\ell}$ has the other
branch.

For each of the direct summands of the direct image Higgs bundle,
we will get two ramification points of the spectral variety over each
point of $\ell \cap \Wob _1$. This agrees with the ramification
pattern of $Y_1\rightarrow X_1$ over $\Wob _1$.

We need to verify the part of Hypothesis \ref{pushforward-hyp} about
points of type \ref{type4}, namely the ramification points of $\bT$ over
$\ell$.

\begin{lemma}
\label{alphaatram01}
For a general line $\ell$, the value of the spectral $1$-form $\alpha$
on the vertical direction in the tangent space of $Y_0$ at any of the
ramification points of $\bT/\ell$, is nonzero.
\end{lemma}
\begin{proof}
As $\ell$ varies, each of these ramification points varies and
constitutes a general point of $\Kum$. So it suffices to note that the
value of $\alpha$ on the vertical direction is nonzero.  This is
because the nilpotent residue of the Higgs field is nonzero along
$\Kum$, which in turn is a consequence of the Bogomolov-Gieseker
inequality as explained in the proof of Lemma \ref{nonzeronormal}
below.
\end{proof}

\subsection{Apparent singularities}

The discriminant of $(\bigHeckebar (a) _{\ell} , \bT\cup \bN)$
relative to $\ell$ consists of the following kinds of points:
\begin{itemize}
\item
images of simple ramification points of $\bT$ over $\ell$, these are the
points of $\ell \cap \Wob _1$ and we expect singularities of the local
system at these points;

\item
  images of points where the Hecke line goes through a trope conic---
  by Lemma \ref{notinttropeconics}
these points are not contained in the lower critical locus, so 
by Proposition \ref{nonisolated} they do not contribute singularities
of the higher direct image local system; and

\item
images of the points $\bPP'$ and $\bQQ'$ where $\bT\cup \bN$ has a
node---we'll see for topological reasons below that these do not
contribute singularities, and indeed the Dolbeault higher direct image
consideration will also show that.
\end{itemize}

\

The following lemma gives the topological proof for the third part. 

\

\begin{lemma}
\label{topproof}
Suppose given a point where the horizontal divisor has a simple normal
crossing with both branches \'{e}tale over the base, and suppose that the
monodromy transformations around each branch are direct sums of
identical size $2$ unipotent Jordan blocks. Then such a point does not
contribute to the monodromy of the higher direct image local system.
\end{lemma}
\begin{proof}
This description characterizes the nontrivial local rank $2$ pieces of
the local system upstairs.  However, we also obtain the same type of
situation if we start with a unipotent rank $2$ local system having
singularities along a horizontal divisor that is simply ramified over
the base, then pull back to a simply ramified double cover of the
base. We will see  in subsection \ref{nilpotentcase} that the
monodromy for a simply ramified horizontal divisor with unipotent
local system, is a transformation of order $2$. Therefore, its
pullback by a ramified double cover has trivial monodromy
transformation.
\end{proof}

\

Even though the topological monodromy transformation for such a point
is trivial, we still need to use the description of \cite{dirim} to
get the description of the spectral line bundle for the higher direct
image at such a point, since the map is singular.

\

\begin{corollary}
\label{singx0x1}
The singularities of the higher direct image local system over $\ell$
are located at the points of $\ell \cap \Wob_1$.
\end{corollary}

\

\begin{remark}
\label{othercaseeff}
In the situation of Lemma \ref{topproof}, if we had monodromy over the
nodal horizontal divisor decomposing into size $2$ blocks with a
reflection instead of a nilpotent monodromy transformations, then the
higher direct image would have the square of a nilpotent
transformation (subsection \ref{parcase}), so again a nontrivial
nilpotent transformation. Applied to the Hecke situation, this means
that if we do the Hecke transformation starting with a local system
with parabolic weights $1/2$ over the Kummer, we will get singularities
at the images of $\bPP'$ and $\bQQ'$ in $\ell$. These are probably the
intersection of $\ell$ with a singular K3 surface inside $X_1$, the
image of the K3 inside the Hecke variety, and which is the
intersection of $X_1$ with another quadric in $\pp^5$.
\end{remark}

\subsubsection{General position arguments}

Let us now look locally at the point $\bPP'$. The horizontal divisor 
$\bT\cup \bN$ has a node there, and the spectral variety $\Sigma$
decomposes as a union of four ordinary double points. 

The point $\bQQ'$ is a place where $\bT$ has a node, with both
branches being \'{e}tale over the base $\ell$. This corresponds to a
point where the plane $P$ is tangent to $\Kum$. Again, $\Sigma$
decomposes into $4$ pieces, and each piece restricted to the tangent
plane gives an ordinary double point.

In order to apply the construction of Theorems \ref{pushforward-thm}
and \ref{pushforward-appli-t5}, we need the following result.

\

\begin{theorem}
\label{pqtype5}
The points $\bPP'$ and $\bQQ'$ are points of type \ref{type5} in the
classification of Subsection \ref{pushforward-statements}.
\end{theorem}
\begin{proof}
The geometric picture shows that the inverse image of $\bPP'$ resp. $\bQQ'$
in $\Sigma$ decomposes into four ordinary double points.  Recall that
over a general point of $\Kum$, the covering $Y_0\rightarrow X_0$
decomposes into four local pieces each of which is a double cover
simply ramified over the local piece of $\Kum$.  Each of these four
pieces leads to a simple double point in the fiber of $\Sigma$, either
over $\bPP'$ or over $\bQQ'$.

Over $\bPP'$, we are restricting to a transverse plane where it is
again a smooth double cover ramified over a smooth curve, then blowing
up the origin. The point $\bPP'$ corresponds to the place where the
exceptional divisor $\bN$ meets the strict transform $\bT$ of $\Kum$, and
the double cover gives there an ordinary double point.

Over $\bQQ'$, we have a plane that is tangent to $\Kum$ with
nondegenerate second fundamental form, such that the two branches of
the intersection between the plane and $\Kum$ correspond to the two
branches of $\bT$ at $\bQQ'$. The restriction of the local double cover to
the tangent plane is again an ordinary double point.

To complete the verification of the conditions required for the
classification of our points in the category ``type \ref{type5}'', we
need to show that the relative critical locus has two branches at each
node of $\Sigma$, with the branches having distinct tangent vectors in
the tangent cone of $\Sigma$.

A degree and monodromy calculation shows that the relative critical
locus, that we have identified also with $\bigHeckehat^{\rm
  ab}(\widehat{a},\tau \widehat{a})_{\ell}$, has two branches at each node of
$\Sigma$. We need to show that they are smooth with distinct tangent
vectors.  This statement amounts to a general position argument
concerning the relationship between the local geometry at a general
point of $\Kum$, and the tautological $1$-form on the covering
$Y_0/X_0$.

The cases of $\bPP'$ and $\bQQ'$ are different although similar. The basic
idea is to find three (or more) directions in a plane that come from
the geometry of the moduli spaces, independent of the point in the
Hitchin base; then the tautological $1$-form is a fourth direction. If
the cross-ratio between these directions did not move, it means that
the data coming from the tautological $1$-form would be constant, and
we try to get a contradiction. That will show that the cross-ratio is
general for general points, which in turn gives the required general
position property.

The proofs for $\bPP'$ and $\bQQ'$ will be completed in detail in the
following subsections.
\end{proof}

Let $Y_{\Kum} \subset Y_0$ be the reduced inverse image of $\Kum$ in
$Y_0$.  We saw in Proposition \ref{descrip-ram} that
$Y_{\Kum}\rightarrow \Kum$ is a $4$-sheeted covering with an explicit
description.

Let $H_{\Kum}$ be the set of points in the Hecke correspondence
$\bigHeckebar (a)$ which are intersection points of two
lines in the Hecke fiber over points of $\Kum$. Thus, $H_{\Kum}$ is
the Kummer K3-surface obtained by blowing up the $16$ nodes of $\Kum$,
with a map $H_{\Kum}\rightarrow \Kum$ that is an isomorphism outside
of the nodes.  We obtain a map $H_{\Kum} \rightarrow X_1$ that we'll
view as a rational map $\nu : \Kum \dashrightarrow X_1$, whose image
is the Kummer K3 surface in $X_1$
\cite{BeauvilleSurfaces,Dolgachev200,GH,Hudson,Keum}.

\

\begin{lemma}
\label{tankum}
Suppose $x\in \Kum$, not a node. Let $z:= \nu (x) \in X_1$ be the
image of the point of $H_{\Kum}$ lying over $x$. Let $V_x \subset X_0$
be the Hecke line over $z$.  Then $V_x$ is tangent to $\Kum$ at $x$.
\end{lemma}
\begin{proof}
Let $\pzo ^{-1}(z)\subset \bigHeckebar (\widehat{a})$ be the inverse
image, which is a $\pp^1$ mapping isomorphically to the Hecke line
$V_x$. If $x'\in H_{\Kum}\subset \bigHeckebar (\widehat{a})$ is the
point over $x$, so that $z=\pzo (x')$, then $\pzo ^{-1}(z)$ passes
through $x'$.  The divisor $\qzo ^{-1}(\Kum ) \subset \bigHeckebar
(\widehat{a})$ is, locally near $x'$, a union of two branches each of
which is a part of the $\pp^1$-bundle over the Jacobian of $C$. The
curve $\pzo ^{-1}(z)$ meets each of the branches, so its intersection
with the divisor $\qzo ^{-1}(\Kum )$ has multiplicity $\geq 2$.  This
intersection number is the same as the local intersection number of
$V_x$ with $\Kum$ at $x$.  The fact that it is $\geq 2$ means that the
line $V_x$ is tangent to $\Kum$ at $x$.
\end{proof}

\

The inverse image of $\Kum$ in $Y_0$ is $2Y_{\Kum}$ as a divisor,
since the map $Y_0\rightarrow X_0$ is fully simply ramified along
$Y_{\Kum}$. At points where $Y_{\Kum} /\Kum$ is etale, there is a
specified tangent direction to $Y$, normal to $Y_{\Kum}$, namely the
directions that map to zero in the tangent space of $X_0$. Call this
subsheaf
$$
N^v \hookrightarrow T(Y_0) |_{Y_{\Kum}}.
$$

\begin{lemma}
\label{nonzeronormal}
The restriction of the spectral $1$-form $\alpha$ of $Y_0$ to $N^v$ is
nonzero.
\end{lemma}
\begin{proof}
If it were zero everywhere, this would imply that the Higgs field on
the Higgs bundle $\pi _{*}(\LY_0 )$ does not have singularities over
$\Kum$. Recall from Corollary \ref{extremalzerocor} that for an
appropriate choice of $\LY_0$, the second Chern character violates the
Bogomolov-Gieseker inequality. But, if there were no singularities of
the Higgs field along $\Kum$, then the Higgs field would be
logarithmic along just the trope part of the wobbly divisor. From
Corollary \ref{nothree}, the trope part of the wobbly divisor has
normal crossings up to codimension $2$. Since the spectral variety $Y$
is irreducible, the Higgs bundle is stable. This would contradict
Mochizuki's Bogomolov-Gieseker inequality
\cite{Mochizuki-kh1}. Therefore, the restriction of $\alpha$ to $N^v$
is nonzero.

In fact, it is nonzero at all the points where it is well-defined over
the smooth locus of $\Kum$ minus the trope conics, because the residue
of the Higgs field has to have constant Jordan type along the smooth
points of the wobbly divisor.
\end{proof}

\

\begin{lemma}
\label{kervaries}
Consider the subspace $\ker (\alpha ) \cap T(Y_{\Kum})$ varying as a
function of the point in $Y_{\Kum}$. Over a general point $x\in \Kum$,
then the projections of these subspaces to $T_x(\Kum )$ vary as a
function of the point in the Hitchin base.
\end{lemma}
\begin{proof}
If the subspaces depend only on $C$ and $x\in \Kum$ then, in
particular, their values on different sheets of $Y_{\Kum}$ would be
the same. We'll see that this is not the case.

Points of $Y_{\Kum}$ are represented by $(L,\tilde{u}+\tilde{v})$
where $\tilde{u},\tilde{v} \in \Ctilde$ are points lying over their
images denoted $u,v\in C$, and $L^{\otimes 2} = \Oo _C(u+v-2\pw
)$. The restriction of $\alpha$ to $Y_{\Kum}$ is given by adding the
values of the tautological form on $\Ctilde$ evaluated on $\tilde{u}$
and $\tilde{v}$. Specialize near a point where $\tilde{u}$ is not a
ramification of $\Ctilde / C$ but $\tilde{v}$ is a ramification
point. As $v$ moves around the branch point, $\tilde{v}$ changes
branches and the tautological form $\alpha$ changes sign when viewed
as a dual element of $T_v(C)$. Then, $\alpha$ is the sum of a fixed
part namely the tautological form at $\tilde{u}$, plus a part that
changes sign namely the value at $\tilde{v}$. Thus, the directions of
$\ker (\alpha )$ over these two points in different sheets of
$Y_{\Kum}$ are different.
\end{proof}

\subsubsection{Arguments for \texorpdfstring{$\bPP'$}{Pprime}}

Look first at the point $\bPP'$. View it as a general point $\bPP \in \Kum
\subset X_0$. The Hecke fiber over $x$ has two lines that meet in a
point, and the image of this point in $X_1$ is a point $z$ where the
two lines meet. The line $\ell$ in our picture is one of these
two. Going back, the line $\ell$ corresponds to a plane $P\subset
X_0$, transverse to $\Kum$ at $\bPP$ and tangent to it at $\bQQ$. On the
other hand, the Hecke fiber over $z$ is a line in $X_0$, and this line
is tangent to $P\cap \Kum$ at $\bPP$. Indeed, after blowing up the Hecke
line and the strict transform of $\Kum$ both meet the exceptional
divisor $\bN$ at the same point, meaning that the Hecke line was tangent
to $\Kum$ before blowing up.

The double cover $Y_P\rightarrow P$ ramified along $\Kum$ has a point
$y$ over $\bPP$, and its tangent space at $y$ is a $2$-dimensional
space containing the following subspaces: the vertical space $N^v(y)$
of the ramification; the tangent space $T_{y}(Y_{P,\Kum})$ of the
ramification divisor $Y_{P,\Kum} := Y_P\cap Y_{\Kum}$; and the two
directions of the pullback of the Hecke line $V_{\bPP}$ which pulls
back to a pair of crossed lines in the double cover since it is
tangent to $\Kum$ (Lemma \ref{tankum}). The cross-ratio of these four
points is fixed, because the $4$ directions are symmetric for the
involution of the double cover. Furthermore, Lemma \ref{nonzeronormal}
implies that the restriction of the spectral $1$-form to $Y_P$ is
nonzero. It therefore defines a $1$-dimensional subspace $A\subset
T_y(Y_P)$.

\

\smallskip

\begin{proposition}
\label{amoves}
For general global parameters $C,\Ctilde$, as the point $\bPP$ moves
around in $\Kum$, the subspace $A$ moves with respect to the framing
of $T_y(Y_P)$ given by the four previously discussed directions.
\end{proposition}
\begin{proof}
Proceed by contradiction: suppose that the subspace $A$ has a fixed
direction with respect to the framing. We first show that it must then
be the direction $T_{y}(Y_{P,\Kum}$. By Lemma \ref{nonzeronormal}, $A$ is
not the direction $N^v(y)$. It seems likely that the two directions
coming from $V_{\bPP}$ should interchange under a global symmetry or
monodromy operation, and if $A$ were in a fixed direction it would
have to be invariant under that operation, which would imply that it
is the direction $T_{y}(Y_{P,\Kum}$. However, we have not been able to
specify such an interchange operation. So, instead, let us calculate at
a special point.

Recall from Proposition \ref{descrip-ram} that the points of
$Y_{\Kum}$ are represented by $(L,\tilde{u}+\tilde{v})$ where
$\tilde{u},\tilde{v} \in \Ctilde$ are points lying over their images
denoted $u,v\in C$, and $L^{\otimes 2} = \Oo _C(u+v-2\pw
)$. Specialize near a point on $\Kum$ where both $u$ and $v$ are
(different) branch points of $\Ctilde / C$. This corresponds to a
place where two other pieces of the movable ramification locus of
$Y/X$ meet the ramification over $\Kum$.  Both $\tilde{u}$ and
$\tilde{v}$ then branch. We can write this in local coordinates
$y_1,y_2,y_3$ for $Y$ over coordinates $x_1,x_2,x_3$ for $X$, such
that $x_1=0$ is the equation of $\Kum$, $y_1=0$ is the equation of
$Y_{\Kum}$, and the map is given by $x_i=y_i^2$. The spectral $1$-form
along $Y_{\Kum}$ is obtained from the tautological $1$-form on
$\Ctilde$ by adding the values at $\tilde{u}$ and $\tilde{v}$.  The
tautological $1$-form on $\Ctilde$ looks locally like $z^2dz$ if $z$
is the coordinate on $\Ctilde$.  We can therefore write the leading
term as
$$
\alpha = a(y)d(y_1) + y_2^2 d(y_2) + y_3^2 d(y_3).
$$
Since $\alpha$ comes from a linear form on the abelian variety, it
does not vanish at any point, thus $a(0)$ has to be nonzero.

Consider a point $x = (x_{1},x_{2},x_{3})$ near to $\bPP = (0,0,0)$.  The
Hecke line $V_{x}$ becomes, in these coordinates, a curve tangent to
$(x_1=0)$. Let us write it (again looking at the highest order term) as
$$
(x_1,x_2,x_3) = (b_1 t^2, \epsilon _2 + b_2t, \epsilon _3 + b_3 t)
$$
for some coefficients $b_1\neq 0$ and $(b_2,b_3)\neq (0,0)$, and small
$(\epsilon _2,\epsilon _3)$.  Note that as $(\epsilon _2,\epsilon
_3)\rightarrow (0,0)$, the second order term $b_1t^2$ approaches a
nonzero limit ($(0,0,0)$ being itself a point of $\Kum$ that is general with
respect to $C$ although not necessarily with respect to $\Ctilde$).

Specifying a lifting of $(0,\epsilon _2,\epsilon _3)$ in $Y$
i.e. extracting $\epsilon _2^{1/2}$ and $\epsilon _3^{1/2}$, the curve
lifts into $Y$ in two branches corresponding to $\pm b_1^{1/2}$.  One
of these branches is
$$
(y_1,y_2,y_3) = \left( 
b_1^{1/2} t, \epsilon _2^{1/2} (1 + b_2 t / 2\epsilon _2 + \ldots ), 
\epsilon _3^{1/2} (1 + b_3 t / 2\epsilon _3 + \ldots ) 
\right)  .
$$
Its tangent vector is $(b_1^{1/2}, b_2 \epsilon _2^{-1/2} / 2, b_3
\epsilon _3^{-1/2} / 2 )$.

The vertical tangent vector is $(1,0,0)$.  We use the middle direction
to normalize the horizontal tangent vector with respect to the
vertical one, in other words the middle vector should be the sum of
the same multiple of the vertical and normalized horizontal
vectors. The multiple is $b_1^{1/2}$, so the normalized horizontal
vector is
$$
\left( 0, b_1^{-1/2}b_2 \epsilon _2^{-1/2} / 2, b_1^{-1/2}b_3 \epsilon _3^{-1/2} / 2 \right) .
$$
The values of $\alpha$ on the vertical tangent vector is $a(0)$. On
the normalized horizontal vector, evaluating at the point
$(y_1,y_2,y_3)=(0, \epsilon _2^{1/2},\epsilon _2^{1/2})$, it is
$$
\alpha (0, \epsilon _2^{1/2},\epsilon _2^{1/2}) \cdot \left( 0,
b_1^{-1/2}b_2 \epsilon _2^{-1/2} / 2, b_1^{-1/2}b_3 \epsilon _3^{-1/2}
/ 2 \right) = b_1^{-1/2}b_2 \epsilon _2^{1/2} / 2 + b_1^{-1/2}b_3
\epsilon _3^{1/2} / 2 .
$$

This approaches $0$ as $(\epsilon _2,\epsilon _3)\rightarrow
(0,0)$. This tells us that the kernel line $A$ of $\alpha $ on
$T(Y_P)$ approaches the horizontal direction as we approach the
special point.

This shows that if $A$ is some fixed direction with respect to the
framing, it must be the horizontal direction. Let's now take up that
possibility.  In that case, it means that the line $T_{y}(Y_{P,\Kum})$ is
always in the kernel of $\alpha$.  In particular, the restriction of
$\alpha$ to $Y_{\Kum}$ has this fixed foliation as a kernel.  But
Lemma \ref{kervaries} says that this does not happen, completing the
contradiction.
\end{proof}

\

\smallskip

\noindent
The local piece of $\Sigma$ corresponds to a locally defined rank $2$
Higgs bundle on $X_{0}$. It comes originally from a rank $2$ piece of
$\mycal{F}_0$ whose spectral variety is a double cover of $X_0$ simply
ramified at the point of $\Kum $. We'll use the Higgs bundle to
calculate the lower critical locus, knowing that it is the image of
the upper critical locus so if we separate the two branches of the
lower critical locus that will separate them upstairs too.

Use a coordinate system $(x,t)$ on $P\subset X_0$ with
$P\cap \Kum $ given by $t=0$. The double cover looks locally
like $t=w^2$. We can write the form $\alpha$ as
$$
\alpha = (a^+(x,t) + wa^-(x,t))dx 
+ (b^+(x,t) + wb^-(x,t))dw
$$
$$
= (a^+(x,t) + wa^-(x,t))dx +
(wb^+(x,t) + tb^-(x,t))(dt/t).
$$
Locally we assume that the line bundle is trivial $\LY _0 = \Oo
_{P}$.  Then the direct image down to $P$, which is the rank $2$ piece
of our Higgs bundle $\mycal{F}_0 |_P$, has basis $1,w$. In these terms, the
form $\alpha$ leads to a Higgs field in matrix form, with coefficients
being functions of $(x,t)$:
$$
\varphi = \left(
\begin{array}{cc}
a^+ dx  +  b^-dt &    ta^- dx + b^+ dt  \\
a^- dx + b^+ (dt/t) & a^+ dx + tb^- dt
\end{array}
\right) .
$$

\

\noindent
We next blow up the point $(0,0)$ with coordinates 
$(x,v)$ with $t=xv$. This gives $dt=vdx + xdv$. The Higgs field becomes
$$
\varphi = \left(
\begin{array}{cc}
(a^+ + vb^-)dx  +  xb^-dv &    (xva^- + vb^+)dx + xb^+ dv  \\
a^- dx + b^+ (dx/x + dv / v) & (a^+ + vb^-) dx + x^2vb^- dv
\end{array}
\right) .
$$
This is the formula for the Higgs field on
$\bigHeckebar(a)_{\ell}$. We may assume that the map $\bigHeckebar(a)
_{\ell} \rightarrow \ell$ is given by the function $u - v$, so for the
relative differentials it induces the relation $dx=dv$.  The relative
Higgs field becomes
$$
\varphi_{\bigHeckebar (a) _{\ell}/\ell }
= \left(
\begin{array}{cc}
(a^+ + vb^-  +  xb^-)dv &    (xva^- + vb^+ + xb^+ )dv  \\
(a^-  + b^+ /x + b^+ / v)dv  & (a^+ + vb^-  + x^2vb^- )dv
\end{array}
\right).
$$
The support of the cokernel of this matrix is given by its
determinant
$$
{\rm det} (\varphi _{\bigHeckebar (\overline{a}) _{\ell}/\ell } ) 
=
(a^+ + vb^-  +  xb^-)(a^+ + vb^- + x^2vb^- ) 
-
(xva^- + (x+v)b^+ )(a^-  + ((x+v)/xv)b^+ )
$$
$$
= (xv)^{-1} \left[
xv (a^+)^2 - (x+v)^2 (b^+)^2  + \ldots 
\right]
$$
where the next terms have degree $\geq 3$ in $x,v$. 
The vanishing locus, away from the axes, is therefore given to
first order by the equation
$$
(b^+)^2 x^2 + (2 (b^+)^2 - (a^+)^2) xv + (b^+)^2 v^2 = 0.
$$
Dividing by $(b^+)^2$ and setting $c:= (a^+ / b^+)^2$ this becomes
$$
x^2 + (2-c)xv + v^2.
$$
Its discriminant is $(2-c)^2 - 4 = c^2 - 4c$ so if $c\neq 0,4$, in
other words the quotient is not $0$, $-2$ or $2$, then there are two branches. 

\

\begin{corollary}
If $\bPP$ was a general point on $\Kum $ then the quotient $a^+(0,0) /
b^+(0,0)$ is general.
\end{corollary}
\begin{proof}
This restates the result of Proposition \ref{amoves}. 
\end{proof}

\

Therefore, at a general point on $\Kum$, the upper critical locus has
two branches meeting the exceptional divisor at two distinct
points. This completes the proof of Theorem \ref{pqtype5} for the
point $\bPP'$.

\

\subsubsection{Arguments for \texorpdfstring{$\bQQ'$}{Qprime}}

At the point $\bQQ = \bQQ'$ the plane $P$ is tangent to
$\Kum$. Introduce coordinates $x,y$ for the plane and $z$ in the
transverse direction so that $\Kum$ is given by $z=0$ and the plane is
given by $z=xy$. The covering $Y$ has coordinates $x,y,w$ with
$w^2=z$, and $Y_{\Kum}$ has equation $w=0$.

We may scale the coordinates such that the map from the plane to the
line $\ell$ is given by $t=x-y$. Notice that the fibers of the map are
not tangent to the two principal directions of $P\cap \Kum$, as we
have seen above that the two branches of $\bT$ at the point $\bQQ'$ are
etale over $\ell$.

Write the spectral differential form on $Y$ as 
$$
\alpha = a(x,y,w) dx + b(x,y,w) dy + c(x,y,w) dw.
$$
We have $c(0,0,0) \neq 0$ by Lemma \ref{nonzeronormal}. Also, $\alpha$
does not vanish identically on $Y_{\Kum}$, by Lemma \ref{kervaries},
and the kernel of $\alpha$ is not generically given by any
construction depending only on $C$ (such as the lines $x=0$, $y=0$ or
$x=y$ in our coordinate system).  Thus, at our general point we can
assume that $a(0), b(0), a(0) + b(0)$ are nonzero.

Proceed now with the calculation of the directions of the upper
critical locus. The covering $\Sigma $ has equation $w^2 = xy$ mapping
to the plane $P$ with coordinates $x,y$. Blow up the origin in the
$x,y$ plane, so we introduce a coordinate $v$ with $y=xv$.  The
exceptional divisor is $x=0$ and $v$ is the coordinate along it. Now
the form is $\alpha = adx + b(xdv + vdx) + cdw$ and $\Sigma$ has
equation $w^2 = x^2 v$. This normalizes to two branches $w=xv^{1/2}$
and $w=-xv^{1/2}$.  Let $u=v^{1/2}$ be the coordinate at a general
point of the covering of the exceptional divisor, so we have $w=xu$
(there are two branches as $u,-u$ go to the same point $v$).  We have
$v=u^2$ so $dv = 2udu$. We have $y=xu^2$.

This normalization is smooth over $u\neq 0$ (actually everywhere) and
the form becomes
$$
\alpha = adx + b(2xudu + u^2 dx) + c(x du+ udx) .
$$
The mapping function from the plane to the line is $t=x-y$ so
$dt=dx-dy$.  Relative differentials are obtained by working modulo
$dt$, identifying
$$
dx \sim dy=2xudu + u^2 dx.
$$
Modulo $dt$ we get 
$$
dx = \frac{2xu}{1-u^2} du.
$$
Thus 
$$
\frac{\alpha^{\rm rel}}{du}  =a \frac{2xu}{1-u^2} +
b\left( 2xu +\frac{2xu^3}{1-u^2}\right) +
c\left( x+ \frac{2xu^2}{1-u^2}\right) .
$$
The zeros of $\alpha ^{\rm rel}$ are given by the equation
$$
a\cdot 2xu + b\cdot (2xu - 2xu^3 + 2xu^3 ) + c\cdot (x - xu^2 + 2xu^2) 
$$
$$
= (a+b) (2xu) + c(x+xu^2).
$$
One can factor $x$ out of this expression. 

We notice that $a(0)+b(0)$ is the value of $\alpha$ on the tangent
direction of the Hecke line $t=0$, and $c(0)$ is the value of $\alpha$
in the vertical direction. In our current coordinate system, the unit
vectors of these directions are normalized to be related by the lines
described in the previous subsection: the Hecke line is $y=x,z=x^2$
and its two lifts are $y=x,w=\pm x$.

Therefore, curiously enough, we are now in the same general position
setting as in the previous subsection.  Proposition \ref{amoves}
implies that the ratio $(a(0) + b(0))/c(0)$ moves as a function of the
point, so we may assume that it is general. The factored expression
for the zeros of the relative $1$-form $\alpha ^{\rm rel} / xdu$ along
the exceptional curve becomes
$$
1 + \frac{a(0) + b(0)}{c(0)} u + u^2 =0.
$$
There are two distinct points $u$ in the exceptional divisor
corresponding to limits of zeros of $\alpha ^{\rm rel}$. This gives
two branches of the upper critical locus with distinct tangent vectors
in the nodal point of our local piece of $\Sigma$.

This completes the proof of the $\bQQ'$ part of Theorem \ref{pqtype5}. 

\

\subsection{Spectral line bundle for the higher direct image}

Recall the \emph{\bfseries spectral line bundle}, on $Y_0$ or $Y_1$,
is the line bundle $\mycal{L}_{0}$ or $\mycal{L}_{1}$ whose direct
image down to $X_0$ or $X_1$ respectively is the level $0$ piece of
the parabolic Higgs bundles $(\mycal{F}_{0,\bullet},\Phi_{0})$ or
$(\mycal{F}_{1,\bullet},\Phi_{1})$ respectively. For $X_0$ the
parabolic Higgs bundle has trivial parabolic structure (away from the
codimension $2$ tacnode points). On $X_1$ the parabolic structure has
levels $0$ and $-1/2$, so this definition of spectral line bundle
involves a choice.

We would like to calculate the spectral line bundle $\Uu$ on
$\bigHeckehat^{\rm ab}(\widehat{a},\tau \widehat{a})_{\ell}$ corresponding to
the higher direct image of the pullback Higgs bundle
$\mycal{F}_{0,\bigHeckebar(a_{\ell}}$ down to $\ell$, in terms of the
spectral line bundle $\LY_0$ over $Y_0$ and its restriction
$\LY_{\ell}$ to $\bigHeckehat^{\rm ab}(\widehat{a},\tau \widehat{a})_{\ell}$.

For this, we use the statement of Theorem \ref{pushforward-thm}, in
particular the part proven in Theorem \ref{pushforward-appli-t5}. The
previous discussion shows that the points $\bPP'$ and $\bQQ'$ are of
type \ref{type5}, and the remaining points are covered by the other
parts of the classification in Subsection
\ref{pushforward-statements}.  To use the theorem in the presence of
the points $\bPP'$, $\bQQ'$, we should take the answer
$$
\Uu = \LY_{\ell} \otimes \omega _{\bigHeckebar (\overline{a}) _{\ell}/\ell}. 
$$
The normalization of the relative critical locus separates the two
branches at each of the nodes. By Theorem \ref{disjointness1} and
Corollary \ref{disjointness-cor1}, the two components of the
abelianized Hecke are disjoint away from $\bPP', \bQQ'$, so the normalization
that was denoted by $G$ in Subsection \ref{pushforward-statements} is
the same as the disjoint union
$$
G = \bigHeckehat^{\rm ab}(\widehat{a},\tau \widehat{a})_{\ell} =
\bigHeckehat^{\rm ab}(\widehat{a})_{\ell} \sqcup 
\bigHeckehat^{\rm ab}(\tau \widehat{a})_{\ell} . 
$$
The prescription of Theorem \ref{pushforward-thm} gives a parabolic
sheaf $F'_{\cdot}$ such that 
$$
F'_0 = f_{*} (\Uu |_G)
$$
and $F'_{-1/2}$ is the standard subsheaf coming from the ramification
points of $f:G\rightarrow \ell $ over $\ell \cap \Wob _1$. Theorem
\ref{pushforward-thm} says that $F'_{\cdot}$ is a parabolic subsheaf
of the parabolic Higgs bundle $F_{\cdot}$ we are looking for. On the
other hand, we know from general principles that $F_{\cdot}$ has
parabolic degree $0$. Also, $F$ and $F'$ agree except possibly over
the points $f(\bPP'),f(\bQQ')\in \ell$ and these are different from the
points where there is a parabolic structure.

So, if we can show that the parabolic degree of $F'_{\cdot}$ is zero,
this will imply that $F_{\cdot} = F'_{\cdot}$ and we get the
computation of $F_{\cdot}$. The degree will be calculated in the next
subsection.

\subsection{Degree calculation}
\label{degree-calc-X0X1}

Recall our notation that inside $Y_0$ that $\ExY_0$ denotes the
exceptional divisor of the blow-up $\blo_{0} : Y_{0} \to \Prym_{2}$,
and $\FxY _0$ denotes the inverse image of the hyperplane class of
$X_0$.

\

\begin{lemma}
\label{reldiff1}
For the relative canonical class over the degree $1$ moduli space, we
have
$$
\omega _{{\bigHeckebar(a)}/X_1}= \Oo _{X_0}(-2) 
|_{\bigHeckebar (a)} \otimes \Oo _{X_1}(1)|_{\bigHeckebar(a)}.
$$ 
If $\Xi$ denotes the theta divisor on the Prym, or its pullbacks to
$Y_0$, $Y_1$ or $\bigHeckehat^{\rm ab}(\widehat{a},\tau \widehat{a})_{\ell}$,
we have $\FxY _{0} = \Xi - \ExY _{0}$ and $\FxY _{1} = 2\Xi - \ExY
_{1}$. We note that $\Xi^3 = 24$, concording with the calculations
in Propositions \ref{intersections1} and \ref{intersections0}.
\end{lemma}
\begin{proof}
The canonical bundle $\omega _{\bigHeckebar (a)}$ restricts to
$\Oo_{\pp^1}(-2)$ on the fibers of the projection to $X_1$ and to $\Oo
_{\pp^1}(2)$ on the fibers of the projection to $X_0$.  It means that
it restricts to the pullback of $\Oo _{X_0}(2)$ on the Hecke lines
over points of $X_1$, and to the pullback of $\Oo _{X_1}(1)$ on the
Hecke conics over points of $X_0$.  The Picard group of ${\bigHeckebar
  (a)}$ is generated by these two things, and the pullbacks restrict
to trivial bundles on their own Hecke fibers, so we conclude
$$
\omega_{\bigHeckebar (a)} = 
\Oo _{X_0}(-2) |_{\bigHeckebar (a)} \otimes \Oo _{X_1}(-1)|_{\bigHeckebar (a)}.
$$
On the other hand, as $X_1$ is an intersection of two quadrics in
$\pp^5$, its canonical class is $\Oo _{X_1}(-6 + 2 + 2) = \Oo
_{X_1}(-2)$. We conclude the stated formula.

For the last parts, we note that the linear system $|\Xi|$ that
produces the rational map $\Prym_{2} \dashrightarrow \pp^3$ has base
points on the $16$ points that we blow-up to get $Y_{0}$ so on the
resulting map on the blow-up, we subtract $\ExY_{0}$.  The
$6$-dimensional subsystem (anti-invariant part) of $|2\Xi|$ that
provides the rational map $\Prym_{3} \dashrightarrow \pp^5$ has $\Chat$ as
base locus so we subtract $\ExY_{1}$. For the verification, the
formulas of Propositions \ref{intersections0} and \ref{intersections1}
give
$$
(\FxY _{0} + \ExY _{0})^3  = 24 \;\;\; \mbox{ and } \;
(\FxY _{1} + \ExY _{1})^3  = 192
$$
compatible with $\Xi^3=24$. 
\end{proof}

\

\noindent
Recall that $\bigExc_0$ and $\bigExc_1$ are the exceptional divisors
in $\bigHeckehat^{\rm ab}$. We could denote their restrictions over
the point $\widehat{a} = (A,\tilde{t})\in \Chat$ by $\bigExc _0(\widehat{a})$
and $\bigExc _1(\widehat{a})$. Those are thus the divisors in
$\bigHeckehat^{\rm ab}(\widehat{a})_{\ell}$ that are strict transforms
of the divisors $\ExY_0$ and $\ExY_1$.

We get that
$$
\ExY_{0} |_{\bigHeckehat^{\rm ab}(\widehat{a})_{\ell}} = \bigExc _0(\widehat{a})
$$
whereas 
$$
\ExY_{1} |_{\bigHeckehat^{\rm ab}(\widehat{a})_{\ell}} = \bigExc _0(\widehat{a}) + 
\bigExc _1(\widehat{a})
$$
because starting  from $Y_1$ we blow up lines contained in $\ExY _1$ to get 
$\bigHeckehat^{\rm ab}(\widehat{a})_{\ell}$. 
This gives 
$$
\FxY_{0} |_{\bigHeckehat^{\rm ab}(\widehat{a})_{\ell}} = 2 \Xi -
\bigExc _0(\widehat{a})
$$
and 
$$
\FxY_{1} |_{\bigHeckehat^{\rm ab}(\widehat{a})_{\ell}} = 4 \Xi - 
\bigExc _0(\tilde{a}) - \bigExc _1(\tilde{a}) .
$$
This gives for the relative differentials
$$
\omega _{{\bigHeckebar (a)}/X_1}|_{\bigHeckehat^{\rm ab}(\widehat{a})_{\ell}} =
\Oo (-2 \FxY_{0} |_{\bigHeckehat^{\rm ab}(\widehat{a})_{\ell}}
+
\FxY_{1} |_{\bigHeckehat^{\rm ab}(\widehat{a})_{\ell}})
$$
$$
= \Oo ( \bigExc _0(\widehat{a}) - \bigExc _1(\widehat{a})).
$$
Recall that---up to tensoring with degree $0$ line bundles
$(\blo_{0}^{*}\Lprym_{0})$ resp. $(\blo_{1}^{*}\Lprym_{1})$---the
spectral line bundle on $Y_0$ is $\LY _{0}= \Oo _{Y_0}(2\FxY _{0}+
\ExY _{0})$ and the spectral line bundle on $Y_1$ is $\LY _{1}= \Oo
_{Y_1}(\FxY_{1})$.  It follows that up to numerical equivalence
$$
\LY _{0} |_{\bigHeckehat^{\rm ab}(\widehat{a})_{\ell}} =
\Oo (4\Xi - \bigExc _0(\widehat{a}))
$$
and
$$
\LY _{1} |_{\bigHeckehat^{\rm ab}(\widehat{a})_{\ell}} = \Oo (4\Xi - 
\bigExc _0(\widehat{a}) - \bigExc _1(\widehat{a}) ).
$$

\

\smallskip

\begin{proposition}
\label{spectraldirim}
Define the line bundle $\Uu$ over $\bigHeckehat^{\rm
  ab}(\widehat{a},\tau \widehat{a})_{\ell}$ as the pullback of the
spectral line bundle $\LY_0$ on $Y_0$, tensored with
$\omega_{\bigHeckebar (a) _{\ell}/\ell}$. Then $\Uu$ is the
spectral line bundle for the $L^2$ Dolbeault higher direct image Higgs
bundle on $\ell$.  Over each of the two pieces of $\bigHeckehat^{\rm
  ab}(\widehat{a},\tau \widehat{a})_{\ell}$, $\Uu$ is isomorphic to the
restriction of the spectral line bundle that we construct for $X_1$,
from $Y_1$ to $\bigHeckehat^{\rm ab}(\widehat{a})_{\ell}$ (respectively
$\bigHeckehat^{\rm ab}(\tau \widehat{a})_{\ell}$).
\end{proposition}
\begin{proof}
It will be convenient to prove the second part first.  Let us use the
expressions for the spectral line bundles that we are constructing on
$Y_0$ and $Y_1$, as calculated prior to the statement of the
proposition.  We get
$$
\begin{aligned}
\Uu |_{\bigHeckehat^{\rm ab}(\widehat{a})_{\ell}} & = 
\LY_0 \otimes \omega_{\bigHeckebar (\widehat{a}) _{\ell}/\ell} 
|_{\bigHeckehat^{\rm ab}(\widehat{a})_{\ell}} \\
& = \Oo _{\bigHeckehat^{\rm ab}(\widehat{a})_{\ell}}((4\Xi - \bigExc _0(\widehat{a})) + 
(\bigExc _0(\widehat{a}) - \bigExc _1(\widehat{a})  )) \\
& = 
\Oo _{\bigHeckehat^{\rm ab}(\widehat{a})_{\ell}}(4\Xi - \bigExc _1(\widehat{a})).
\end{aligned}
$$
We would like to compare this with the spectral line bundle coming
from $Y_1$ which is, as we have seen above,
$$
\LY _{1} |_{\bigHeckehat^{\rm ab}(\widehat{a})_{\ell}}  = 
\Oo _{\bigHeckehat^{\rm ab}(\widehat{a})_{\ell}}(4\Xi - 
\bigExc _0(\widehat{a}) -\bigExc _1(\widehat{a})).
$$ These two differ by $\bigExc _0(\widehat{a})$. However, remember that
we are looking at the restriction over a line $\ell \subset X_1$. The
divisor components of $\bigExc _0(\widehat{a})$ are obtained by blowing up
$16$ disjoint $\pp^{1}$'s in $Y_{1}$. These $\pp^{1}$'s are sixteen
fibers of $\ExY_{1} = \Chat\times \pp^{1}$, so they map to $16$ lines
in the wobbly locus $\Wob _1\subset X_1$. A general line $\ell$ will
miss these. Thus, as restricted over $\ell$, the two expressions are
the same.

The same discussion holds on $\bigHeckehat^{\rm ab}(\tau
\widehat{a})_{\ell}$. Thus, $\LY _{1}|_G\cong \Uu |_G$.

Let $F_{\cdot}$ denote the $L^2$ Dolbeault pushforward (on $\ell$)
with its parabolic structure, and let $F'_{\cdot}$ be the parabolic
bundle on $\ell$ whose spectral line bundle is $\Uu$ over $G$. The
parabolic degree of $F_{\cdot}$ is zero since it corresponds to a
harmonic bundle.

The spectral line bundle $\LY _{1}$ pushes forward to a parabolic
Higgs bundle with first parabolic Chern class equal to $0$, over $X_1$
and hence over a line $\ell \subset X_1$. The parabolic structure here
is also the standard one coming from the ramification over $\Wob _1$.
As we have identified the two line bundles $\LY _{1}|_G$ and $\Uu
|_G$, and since the parabolic structures are standard coming from
ramification points of $G$ over $\ell \cap \Wob _1$, it follows that
the parabolic degree of $F'_{\cdot}$ is zero. Thus, $F_{\cdot} =
F'_{\cdot}$. We conclude that $\Uu$ is the spectral line bundle for
$F_{\cdot}$.
\end{proof}

\

\subsection{Restriction to a line}
\label{sub-resline}

We have been considering the local system corresponding to the
parabolic Higgs bundle $\left(\mycal{F}_{1,\bullet},\Phi_{1}\right)$
via the non-abelian Hodge and Riemann-Hilbert correspondences. So far
we focused on the restriction of this local system to a general line
$\ell \subset X_1$.  We show in this subsection that knowing this
restriction suffices in order to identify the local system on
$X_1-\Wob_1$. Throughout the subsection, $\ell$ denotes a general line
in $X_1$. Let $\ell ^{\circ}:= \ell -(\ell \cap \Wob_1)$.

\

\smallskip

\begin{lemma}
\label{irredL}
The local system on $\ell ^{\circ}$ constructed in Section
\ref{chapter-d1}, is irreducible.
\end{lemma}
\begin{proof}
It suffices to show that the spectral covering is irreducible. Recall
from Lemma \ref{three-eight} that the spectral covering has the
following description. The line $\ell$ has a natural trigonal covering
$k:C\rightarrow \ell$ whose fiber over $x\in \ell$ is identified with
the set of three lines through $x$ that are different from $\ell$. The
$8$ branch points are the points of $\ell \cap \Wob_1$.

We have the spectral cover $\Ctilde / C$. For a point $x\in \ell$ we
can make the following set with $8$ elements: it is the set of
liftings of the subset $k^{-1}(x)$ to a subset of three elements of
$\Ctilde$. This family determines a covering of $\ell$ of degree $8$,
and that is the same as $Y_1\times _{X_1} \ell \rightarrow \ell$.

For $\ell$ general, the set of $4$ branch points of $\Ctilde / C$ maps
to a subset of $4$ distinct points in $\ell$. Indeed, the set of four
branch points is the inverse image of a general pair of points in
$\pp^1$ under the hyperelliptic map $\he_{C} : C \to \pp^{1}$, and a
general trigonal map $C\rightarrow \pp^1$ does not identify opposite
points under the hyperelliptic involution on $C$, so the images of the
four branch points are distinct.  Thus, as $x$ moves around in
$\ell^{\circ}$, we can change individually the parity of any one of
the liftings of the three points. This shows that the monodromy action
on the $8$ points is transitive, so $Y_1\times _{X_1} \ell$ is
irreducible. That was the spectral covering of the Higgs bundle on
$\ell$ corresponding to the restricted local system, so the local
system is irreducible.
\end{proof}

\

\smallskip

\begin{proposition}
\label{prop-compare}
Suppose $V$ and $V'$ are two local systems on $X_1-\Wob _1$, such that
$$
V|_{\ell^{\circ}} \cong V '|_{\ell^{\circ}} 
$$
and such that this local system on $\ell^{\circ} $ is
irreducible. Then $V\cong V'$.
\end{proposition}
\begin{proof}
Let $\rho : \pi_{1}(X_{1}-\Wob_{1},o) \to GL_{N}(\cc)$ and $\rho '
\pi_{1}(X_{1}-\Wob_{1},o) \to GL_{N}(\cc) $ denote the monodromy
representations of $V$ and $V'$ with respect to a basepoint $o \in
\ell^{\circ}$ chosen in $\ell ^{\circ}$. Let $\zeta :
\pi_{1}(\ell^{circ},o) \to GL_{N}(\cc)$ denote the monodromy
representation of $V|_{\ell^{\circ}}$ and $V'|_{\ell^{\circ}}$,
assuming we choose framings making these representations the same.

Suppose $a\in \pi _1(X_1-\Wob _1,o)$. Then $\rho (a)$ and $\rho '(a)$
are morphisms of representations of $\pi _1(\ell ^{\circ},o)$ between
$\zeta$ and its conjugate $\zeta ^a$. Since $\zeta$ is irreducible,
these two morphisms differ by a scalar. This gives a rank $1$
character $\chi : \pi_{1}(X_{1}-\Wob_{1},o) \to \cc^{\times}$ such
that $\rho\otimes \chi \cong \rho'$, and $\chi$ is trivial on $\pi
_1(\ell ^{\circ},o)$.

To finish the proof, we need to note that $\chi$ is trivial. This is
because the map
\begin{equation} \label{eq:loops.surject}
H_1(\ell^{\circ},\zz) \rightarrow H_1(X_1 - \Wob _1,\zz)
\end{equation}
is surjective. Indeed, this follows from the following

\

\begin{claim} \label{claim:loop}
The fist homology group $H_{1}(X_{1} - \Wob_{1}, \zz)$ is cyclic and
generated by the linking loop in $X_{1}$ going around some smooth point
of $\Wob_{1}$.
\end{claim}
\begin{proof} Note first that the statement of the claim
makes sense since any two linking loops at smooth points of $\Wob_{1}$
are conjugate in $\pi_{1}(X_{1} - \Wob_{1})$ and hence are
homologous. This follows immediately since $\Wob_{1}$ is irreducible
and hence the smooth locus of $\Wob_{1}$ is connected. Using the
tubular neighborhood theorem we can view the linking loops as two
different fibers of the circle bundle in the normal bundle
$N_{\Wob_{1}^{\op{smooth}}/X_{1}}$ and hence they are homotopy
equivalent up to a conjugation via a path connecting the two points in
$\Wob_{1}^{\op{smooth}}$ over which these circle fibers sit.

\

Next recall that if $M$ is a connected, oriented, not necessarily
compact, $C^{\infty}$-manifold and if $Z \subset M$ is a connected,
oriented $C^{\infty}$ submanifold of dimension $d$, which is closed in
the Eucledian topology, then $Z$ defines a $d$-dimensional cohomology
class $[Z] \in H^{d}(M;\mathbb{Z})$.  Indeed, let $Z \subset
\mathfrak{T} \subset M$ be a tubular neighborhood of $Z$.  Then for
every $k$ we have canonical identifications $H^{k}(M,M-Z;\zz) =
H^{k}(\mathfrak{T},\partial \mathfrak{T};\zz)$ (by excision) and also
$H^{k}(\mathfrak{T},\partial \mathfrak{T};\zz) = H^{k-d}(Z;\zz)$ (by
the Thom isomorphism theorem). In particular we get that
\[
H^{d}((M,M-Z;\zz) = H^{d}(\mathfrak{T},\partial \mathfrak{T};\zz) =
H^{0}(Z;\zz) = \zz,
\]
where the last equality holds since $Z$ is connected.  But from the
long exact cohomology sequence of the pair $(M,M-Z)$ we have a
canonical map
\begin{equation} \label{eq:pairmap}
H^{d}(M,M-Z;\zz) \to H^{d}(M;\zz)
\end{equation}
and we define $[Z] \in H^{d}(M;\zz)$ to be the image of $1 \in \zz =
H^{d}(M,M-Z;\zz)$ under the map \eqref{eq:pairmap}.  The cap product
\[
[Z]\cap (-)  : H_{k}(M;\zz) \to H_{k-d}(Z;\zz) \to H_{k-d}(M;\zz) 
\]
fits in the classical \emph{\bfseries tube exact sequence} in homology
\[
\xymatrix@1@M+0.2pc@C+0.6pc{\cdots \ar[r] & H_{k-d+1}(Z;\zz)
  \ar[r]^-{\mathsf{tube}_{Z/M}} & H_{k}(M-Z;\zz) \ar[r] & H_{k}(M;\zz)
  \ar[r]^-{[Z]\cap (-)} & H_{k-d}(Z;\zz) \ar[r] & \cdots,}
\]
where the middle map is induced from the inclusion $M-Z \subset M$ and
the tube map $\mathsf{tube}_{Z/M}$ sends a $(k-d+1)$-cycle $A$ in $Z$ to the
$k$-cycle in $M-Z$ which is the total space of the $S^{d-1}$-bundle
$\partial\mathfrak{T}_{|A} \subset \partial\mathfrak{T} \subset M-Z$.

Taking this into account, consider $Z \subset M$ to be the
$C^{\infty}$ manifolds underlying the smooth complex varietties
$(\Wob_{1} - \op{Sing}(\Wob_{1})) \subset (X_{1} -
\op{Sing}(\Wob_{1}))$. Thus $Z$ and $M$ are oriented, $Z$ is of real
codimension $2$, and since $X_{1}$ and $\Wob_{1}$ are irreducible,
both $Z$ and $M$ are connected. Furthermore $M - Z = X_{1} - \Wob_{1}$
as topoloical spaces and so the piece of the tube sequence 
corresponding to $k=1$ reads
\[
\xymatrix@1@M+0.2pc@C+0.6pc{\cdots \ar[r] & H_{0}(Z;\zz)
  \ar[r]^-{\mathsf{tube}_{Z/M}} & H_{1}(X_{1}-\Wob_{1};\zz) \ar[r] & H_{1}(M;\zz)
  \ar[r]^-{[Z]\cap (-)} & 0 \ar[r] & \cdots.}
\]
Since $\op{Sing}(\Wob_{1})$ is a compact subvariety of complex
codimension $2$ in $X_{1}$ we have that $H_{1}(M;\zz) =
H_{1}(X_{1}-\Wob_{1};\zz) = H_{1}(X_{1};\zz) = 0$. Therefoe
\[
\mathsf{tube}_{Z/M} : H_{0}(Z:\zz) \to H_{1}(X_{1}-\Wob_{1};\zz)
\]
is surjective. But the $H_{0}(Z;\zz) = \zz$ is generated by the class
of of a point in $Z$, and by the definition of $\mathsf{tube}_{Z/M}$
for any point $pt \in Z$, we have that $\mathsf{tube}_{Z/M}(pt)$ is
the class of the circle in $X_{1} - \Wob_{1}$ linking to $Z = \Wob_{1}
- \op{Sing}(\Wob_{1})$ this point. This proves the claim.
\end{proof}

\

\noindent
Finally, observe that a general line $\ell$ will intersect $\Wob_{1}$
transversally at a set of $8$ smooth points, and so any simple loop in
$\ell$ that goes once around a point $x \in \ell\cap \Wob_{1}$ will be
a a loop in $X_{1}$ that links to $\Wob_{1}$ at $x$. Thus the image of
$H_{1}(\ell^{\circ},\zz)$ in $H_{1}(X_{1} - \Wob_{1},\zz)$ contains a
linking loop and so the map \eqref{eq:loops.surject} is
surjective. This implies that $\chi$ is trivial, and hence $V\cong V'$
which completes the proof of the proposition.
\end{proof}

\

\noindent
Apply this now to the Hecke transform, a rank $16$ local
system on $X_1-\Wob_1$. 

\

\smallskip

\begin{lemma}
The rank $16$ Hecke transform local system on $X_1-\Wob_1$ decomposes
as a direct sum of two rank $8$ local systems which are isomorphic,
when restricted to a general line, to the local system constructed in
Section \ref{chapter-d1}.
\end{lemma}
\begin{proof}
The arguments of this section show that the restriction of the rank
$16$ local system to a general line decomposes as a direct sum of two
copies of the rank $8$ local system constructed in Section
\ref{chapter-d1}. Furthermore, such a decomposition can be obtained by
considering the two pieces $\bigHeckehat^{\rm ab}(\widehat{a})$ and
$\bigHeckehat^{\rm ab}(\tau \widehat{a})$. We obtain a family of
decompositions over the lines, with the property that when two lines
intersect the decompositions correspond. Now, a complete intersection
of $X_1$ with two general hyperplanes is an elliptic curve that can
degenerate into a cycle of $4$ lines $Z=\ell _1 \cup \ell _2 \cup \ell
_3 \cup \ell _4$.  We get the decomposition on each of these lines and
the decompositions coincide on the intersection points. This gives a
decomposition of the rank $16$ local system over the singular curve
$Z$, into two pieces of rank $8$.  The map $\pi _1(Z-Z\cap \Wob _1)
\rightarrow \pi _1(X_1-\Wob _1)$ is surjective, so we get a
decomposition of local systems on $X_1-\Wob _1$.
\end{proof}

\

\smallskip

\begin{corollary}
The rank $16$ Hecke transform local system on $X_1-\Wob_1$ decomposes
as a direct sum of two rank $8$ local systems which are isomorphic to
the local system constructed in Section \ref{chapter-d1}.
\end{corollary}
\begin{proof}
The lemma gives the decomposition, and the pieces restrict to a
general line to the given rank $8$ local system. By Lemma \ref{irredL}
and Proposition \ref{prop-compare}, the pieces of the decomposition
are globally isomorphic to the rank $8$ local system constructed in
Section \ref{chapter-d1}.
\end{proof}

\

\section{Hecke transformation from \texorpdfstring{$X_1$}{X1} to
  \texorpdfstring{$X_0$}{X0}}
\label{chapter-heckex1x0}

The standard diagram for the direction of the Hecke transformation
$(X_1\rightarrow X_0)$ shows the big Hecke correspondence fitting into
\[
  \xymatrix@M+0.5pc@-0.5pc{
    & \bigHeckebar \ar[dl]_-{p} \ar[dr]^-{q} & \\
    X_{1} & & X_{0}\times \Cbar .
    }
\]
Fix a point $a = (A,t)\in \Cbar$ where $A^{\otimes 2} = \Oo _C(t-\pw
)$ and consider the Hecke correspondence $\bigHeckebar(a)$ with its
diagram
\[
  \xymatrix@M+0.5pc@-0.5pc{
    & \bigHeckebar (a) \ar[dl]_-{p} \ar[dr]^-{q} & \\
    X_{1} & & X_{0} .
    }
\]
The objective in this section is to pull-back the constructed Higgs
bundle from $X_1$ to $\bigHeckebar (a)$ and take the higher
direct image along $q$ to $X_0$.

The big abelianized Hecke fits into a diagram of the form 
 \[
\xymatrix@M+0.25pc{
  & \bigHeckehat^{\rm ab} \ar[dl]_-{p^{\rm ab}}  \ar[dr]^-{q^{\rm ab}}  & \\
 Y_{1} & & Y_{0}\times \Chat .
} 
\]
The two points $\widehat{a}$ and $\tau\widehat{a}$ over $a\in
\Cbar$ give two abelianized Hecke varieties $ \bigHeckehat^{\rm
  ab}(\widehat{a})$ and $ \bigHeckehat^{\rm ab}(\tau \widehat{a})$. These
are the same as the varieties with the same notation in the previous
section, and recall that their disjoint union is denoted by
$$
\bigHeckehat^{\rm ab}(\widehat{a},\tau \widehat{a}):= 
\bigHeckehat^{\rm ab}(\widehat{a}) \sqcup 
\bigHeckehat^{\rm ab}(\tau \widehat{a}).
$$
The map $\bigHeckehat^{\rm ab}(\widehat{a}) \rightarrow Y_1$ is the
blow-up along $16$ lines contained in the wobbly locus, while the map
$\bigHeckehat^{\rm ab}(\widehat{a}) \rightarrow Y_0$ is the blow-up along
a curve isomorphic to $\Chat$ (and that will be called by the same
name) inside $Y_0$.

\

\subsection{Restriction to a line} \label{ssec:lineinX0}

As before, fix a general line $\ell \subset X_0$. Let $\bigHeckebar
(a) _{\ell}$ be the inverse image of $\ell$ in the Hecke variety
$\bigHeckebar(a)$.  As $\ell$ is general, it does not meet the image
in $X_0$ of the curve $\Chat \subset Y_0$, so if $\bigHeckehat^{\rm
  ab}(\widehat{a})_{\ell}$ denotes the inverse image of $\ell$ in
$\bigHeckehat^{\rm ab}(\widehat{a})$, the projection induces an
isomorphism
$$
\bigHeckehat^{\rm ab}(\widehat{a})_{\ell}\stackrel{\cong}{\rightarrow} Y_{0,\ell}.
$$
The same holds for the other piece $\bigHeckehat^{\rm ab}(\tau
\widehat{a})_{\ell}$.

Recall from Corollary~\ref{disjointness-cor0} that for $\ell$ general,
the images of $\bigHeckehat^{\rm ab}( \widehat{a})_{\ell}$ and
$\bigHeckehat^{\rm ab}(\tau \widehat{a})_{\ell}$ are disjoint in the
spectral variety
$$
\Sigma := Y_1 \times _{X_1} \bigHeckebar(a) _{\ell} \rightarrow 
\bigHeckebar (a) _{\ell} .
$$
Thus, we may treat each piece separately. We will look mainly at
$\bigHeckehat^{\rm ab}(\widehat{a})_{\ell}$ with the understanding that
the arguments for $\bigHeckehat^{\rm ab}(\tau \widehat{a})_{\ell}$ are
the same.

\

\begin{lemma}
\label{hlhyperplane}
The Hecke correspondence $\bigHeckebar(a) _{\ell}$ over
$\ell$ maps by a closed immersion onto a hyperplane section that we will
denote by $H_{\ell} \subset X_1$. The image of $\bigHeckehat^{\rm
  ab}(\widehat{a})_{\ell}$ is a curve of genus $25$. The map
$$
\phi : \bigHeckehat^{\rm ab}(\widehat{a})_{\ell}\rightarrow \ell
$$
has degree $8$. It has: $16$ ramification points over the
intersection points of $\ell$ with the trope planes, and $4$
ramification points over each of the $4$ intersection points in $\ell
\cap \Kum$.  It has $32$ other 'movable' ramification points not
mapping to points in the wobbly $\Wob_0$, so there are altogether $64$
branch points.
\end{lemma}
\begin{proof}
Recall that one of the quadrics in the pencil is identified with the
Grassmanian of lines in $\pp^3$, embedded in $\pp^5$ by the Plücker
coordinates. The line $\ell$ itself corresponds to a vector in $v(\ell
)\in \bigwedge ^2\cc^4$. The condition for another line $\emm \subset
\pp^{3}$ to meet $\ell$ is that its vector $v(\emm)\bigwedge ^2\cc^4$
satisfies $v(\emm)\wedge v(\ell ) = 0$ in $\bigwedge ^4\cc^4\cong
\cc$. This is a linear condition on $v(\emm)$ so it corresponds to a
hyperplane in $\pp^6$. The image of $\bigHeckebar(a)_{\ell}$ in $X_1$
corresponds to the subset of points in $X_1$ whose corresponding line
meets $\ell$, in other words it is this hyperplane intersected with
$X_1$. This yields the hyperplane section $H_{\ell}$.

Since $\ell$ is general, its inverse image in $\bigHeckehat^{\rm
  ab}(\widehat{a})$ is the same as its inverse image in $Y_0$ because
the abelianized Hecke is the blow-up of $Y_0$ on a subset that maps to
a curve in $X_0$ that will be missed by a general $\ell$. Now, we may
proceed to calculate the normal bundle of this curve in $Y_0$. It is
the pullback of the normal bundle of $\ell$ in $\pp^3$ which is to say
$\Oo(1)\oplus \Oo(1)$, so that has degree $2$ on $\ell$. Its pullback has
degree $16$.  Now, $Y_0$ is the blow-up of an abelian variety at $16$
points; the exceptional divisors map to the trope planes so our
general $\ell$ meets each of the exceptional divisors once.  The
canonical bundle of $Y_0$ is twice the exceptional divisor, so it has
degree $32$.  The canonical bundle of our curve is the canonical of
$Y_0$ restricted to the curve, plus a divisor of degree $16$, so it
has degree $32 + 16 = 48$. Therefore, the image of $\bigHeckehat^{\rm
  ab}(\widehat{a})_{\ell}$ in $H_{\ell}$ has genus $25$.

The ramification of this curve over $\ell$ is just the restriction to
$\ell$ of the ramification of $Y_0/X_0$. We know that has fixed pieces
including a simple ramification over each trope plane, plus four
simple ramifications over general points of $\Kum$. The remaining
ramification points are really movable, as was noted in Lemma
\ref{ram-mov} (see also Corollary \ref{reallymovable} for the
corresponding statement in the other direction).
\end{proof}

\

\

\noindent The ramification points of the horizontal divisor lie over
points of $\left(\bigcup_{\kappa} \trope_{\kappa}\right) \cap \ell$.
They come from $16$ lines inside $\Wob_1$, and correspond to the
first $16$ ramification points mentioned in Lemma
\ref{hlhyperplane}. The locations of these lines in $\Wob_1$ depend
on the choice of point $a = (A,t)$ used to make the Hecke
correspondence.

We note that the modular spectral covering $\Sigma$ has two
ramification points over each general point of the wobbly locus, so
these $16$ points correspond to $32$ points in the full abelianized
Hecke.  They are distrubuted as $16$ in each of the two pieces
$\bigHeckehat^{\rm ab}(\widehat{a})_{\ell}$ and $\bigHeckehat^{\rm
  ab}(\tau \widehat{a})_{\ell}$.

\

\

\begin{remark}
\label{no-nodes-cusps}
The horizontal divisor in $H_{\ell}$ is $\Wob_1\cap H_{\ell}$. It has
nodes and cusps coming from the nodes and cusps of $\Wob_1$. However,
for a general $\ell$ these do not meet the image of $\bigHeckehat^{\rm
  ab}(\widehat{a})_{\ell}$. In particular, they are not going to
contribute singularities to the higher direct image -- this is another
version of the ``apparent singularities'' encountered in the previous
chapter.
\end{remark}

\

\

\noindent
The proof of this statement is that the nodal and cuspidal loci of the
wobbly locus $\Wob_{1}$ are curves; they pull back to curves in
$Y_1$. The centers of the blow-up $\bigHeckehat ^{\rm ab}(\widehat{a})$
are the $16$ lines in $\Wob_1$ that get blown up to form $\bigExc
_{0,\kappa }$. These are transverse to both the nodal and cuspidal
loci.  Therefore, the pullbacks of the nodal and cuspidal loci in
$\bigHeckehat ^{\rm ab}(\widehat{a})$ are $1$-dimensional. Their
$1$-dimensional images in $X_0$ do not meet a general line $\ell$ so
$\bigHeckehat^{\rm ab}(\widehat{a})_{\ell}$ does not meet the nodal and
cuspidal loci of $\Wbar_1$ and we get the statement of the remark.

Among the hypotheses 
of Subsection \ref{pushforward-statements} is the following statement.

\

\

\begin{lemma}
\label{alphaatram}
For a general line $\ell$, the intersection of the plane $H_{\ell}$
with any of the $16$ lines in $\Wob_1$ that provide ramification of
the Hecke correspondence, is a general point on that line. At such a
point, the value of the spectral $1$-form $\alpha$ on the vertical
direction in the tangent space of $Y_1$, at either of the two
ramification points of $Y_1/X_1$ over this point, is nonzero.
\end{lemma}
\begin{proof}
For the first part, we note that one of these $16$ lines in $\Wob_1$
gets blown up to an $\mathbb{F}_{1}$-surface, that then maps to the
corresponding trope plane in $X_0=\pp^3$. The exceptional divisor
blows down to a point, that we will call the origin, in the trope plane
(corresponding to the intersection in $Y_0$ of the plane with the
$\Chat$ curve). The location of the intersection of $H_{\ell}$ with
the line corresponds to the slope of the line from this origin to the
intersection point of $\ell$ with the trope plane. For $\ell$ general,
this direction is a general point of the line.

Inside $\Prym_{2}$, the curve $\Chat$ passes through $16$ points that
are blown up to get $Y_0$.  The identification $\Prym_{2} \cong
\Prym_{3}$ depends on $a$.  Instead of blowing up the points, we blow
up the curve $\Chat \subset \Prym_{3}$ to get $Y_1$, and this
generates $16$ lines inside $\ExY_1\subset Y_1$.  Let $\mathfrak{v}$
denote one of these lines.  The normal bundle of $\mathfrak{v}$ in $Y_1$
is $\Oo \oplus \Oo (-1)$, with the $\Oo$ direction being the normal
bundle of $\mathfrak{v}$ in $\ExY _1$.

Let us look at how this maps to the normal bundle of the image line
$\lv \subset \Wob _1$. The normal bundle of a line in $\pp^5$ is
$\Oo(1)^{\oplus 4}$, and to get the normal bundle of the line in the
intersection of two quadrics we take the kernel of a map
$\Oo(1)^{\oplus 4} \rightarrow \Oo (2)^{\oplus 2}$. For a typical
line, this kernel will be $\Oo ^{\oplus 2}$.

However, for a line in the wobbly locus, that we recall counts twice
in the set of four lines through each point of $X_1$, we claim that
the kernel is $\Oo (-1) \oplus \Oo (1)$. This may be seen by recalling
that the lines in the wobbly locus are the tangent lines to the copy
of $\Cbar$ that forms the cuspical locus of $\Wob_1$, and when we move
to first order along this $\Cbar$ the tangent line undergoes a
deformation that vanishes at that point in the normal direction. Thus,
these first order deformations are sections of the normal bundle that
vanish. But, the bundle $\Oo ^{\oplus 2}$ does not have any nonzero
sections that vanish somewhere. The only other possibility is the
bundle $\Oo (-1) \oplus \Oo (1)$ since it has to be a subbundle of
$\Oo (1)^{\oplus 4}$. This proves the claim.

As we move in the curve $\Cbar$, the line moves in $\Wob _1$ so the
normal direction has sections, and it has to be the $\Oo (1)$
subbundle. The map $\Oo \rightarrow \Oo (1)$ has a zero at the points
where the line $\lv$ crosses the cuspidal locus of $\Wob_1$. However,
by looking at the local picture of the covering $Y_1\rightarrow X_1$
near such a cuspidal locus, we can see that the map from the full
tangent space of $Y_1$ into the tangent space of $X_1$ has image of
dimension $2$, so the map from the full normal bundle of $\mathfrak{v}
\subset Y_{1}$ to the $\Oo (1)$ piece in the normal bundle of $\lv
\subset X_{1}$ is surjective. We conclude that the bundle of vertical
directions, which is the kernel of this map
$$
\Oo \oplus \Oo (-1) \rightarrow \Oo (1),
$$ is $\Oo (-2)$ sitting in $\Oo\oplus \Oo(-1)$ as a saturated
subbundle. In particular it is not contained in any subbundle of the
form $\Oo (-1)$. The map from this bundle of vertical normal
directions, into the space of normal directions at the original point
of $\Prym_{3}$, therefore does not have image in a plane. Thus, the
tautological $1$-form, that is a nonzero linear form on the tangent
space of $\Prym_{3}$, does not vanish on the vertical normal direction
at a general point of $\mathfrak{v}$.  This proves the second part of
the lemma, in view of the generality statement of the first part.
\end{proof}

\

\

Recall that the fiber of $q : \bigHeckebar(a)_{\ell}\rightarrow \ell$
over a general point $\ell$ is a smooth conic, but over points of the
Kummer it degenerates to a union of two lines.

\

\begin{lemma}
\label{branchK} 
The branch points of $\bigHeckehat^{\rm ab}(\widehat{a})_{\ell}$ lying
over points of the Kummer, all map to points of the Kummer K3 surface
$H_{\Kum}\subset \bigHeckebar (\widehat{a})$, that is to say in
$\bigHeckebar(a)_{\ell}$ they map to points where the two
lines in the fiber meet.
\end{lemma}
\begin{proof}
The reasoning for this is as follows: the horizontal divisor in
$\bigHeckebar(a)_{\ell}$ is the intersection of $\bigHeckebar(a)
_{\ell}$ (considered as a hyperplane $H_{\ell}$ in $X_1$) with the
wobbly $\Wob_1$. But, the lines in the fibers over points of the
Kummer, represent general lines in $X_1$ since $\ell$ intersects the
Kummer in general points.  A general line in $X_1$ will intersect
$\Wob_1$ transversally.

We note that it is necessary to have four Jordan blocks of size two in
the monodromy of the Hecke-transformed local system (the rank $8$
piece corresponding to our chosen branch of $H^{\rm ab}$ out of two,
that we are hoping is our chosen flat bundle of rank $8$ on $X_0$) at
each point of the Kummer. Otherwise, the Hecke transform back in the
other direction from $X_0$ to $X_1$ will not have the right rank.

So, these all have to come from simple ramification points of
$\bigHeckehat^{\rm ab}(\widehat{a})_{\ell}$ over the Kummer
points. However, if those ramification points were to occur on smooth
points of the lines, that wouldn't contribute anything to the
monodromy in the direct image (i.e. the residue of the Higgs field)
since the map is not singular at those locations. This heuristic
argument suggests that $\bigHeckehat^{\rm ab}(\widehat{a})_{\ell}$
should have four points over each point where a fiber of $\phi$ breaks
into two lines.

This can be shown in terms of bundles.  Suppose $L \in \Prym_{2}$ is a
line bundle on $\Ctilde$ whose direct image $V=\pi_{*}L$ is stable
with trivial determinant, so it is a point in $X_0$, and suppose that
$V$ contains a line subbundle $U$ of degree $0$. Then we claim that
$V$ also contains $U^{-1}$.

To see this, we note that the group $\zz/2\times \zz/2$ acts on 
$\Ctilde$, as may be seen for example by expressing
$$
\Ctilde = C \times _{\pp^1} \pp^1
$$
as a fiber product with the double cover $\pp^1\rightarrow \pp^1$
ramified at two points (so the ramification of $\Ctilde/C$ is the
preimage in $C$ of those two points of $\pp^1$).  Explicitly the
non-trivial elements of $\zz/2\times \zz/2$ are given by the covering
involution $\tau : \Ctilde \to \Ctilde$ for thr map $\pi : \Ctilde \to
C$, the covering involution $\sigma : \Ctilde \to \Ctilde$ for the
hyperelliptic map $\hyp_{\Ctilde} : \Ctilde \to \pp^{1}$, and their
composition $\rho = \tau\circ\sigma$.

Now, note that by construction the group element $\rho \in \zz /
2\times \zz / 2$ acts trivially on the Prym of $\Ctilde / C$ but acts
by $-1$ on $ \op{Jac}(C)$.  Applying $\rho$ to the map $\pi ^{*}(U)
\rightarrow L$ we obtain a map $\pi ^{*}(U^{-1}) \rightarrow L$ and
hence by adjunction we get an injective map $U^{-1} \rightarrow V$.
This will give by degree considerations $V=U\oplus U^{-1}$. In the
Hecke correspondence this will correspond to a point at the
intersection of the two $\pp^{1}$ components of the fiber of $q$ over
$V$ (we recall that the affine parts of the two lines themselves were
bundles that were semistable but not polystable).

We can also note that $L\cong \pi ^{*}(U)\otimes \Oo _{\Ctilde}(a+b)$
for an effective degree $2$ divisor $a+b$ on $\Ctilde$, amd that the
image divisor $\pi (a) + \pi (b)$ is fixed by the determinant
condition for $V$. Thus, there are four choices of $(a,b)$ lifting
this divisor to $\Ctilde$.  We get the four claimed branches of
$\bigHeckehat^{\rm ab}(\widehat{a})_{\ell}$ going through a crossing
point in $\bigHeckebar(a)_{\ell}$.
\end{proof}

\

\

\begin{proposition}
\label{abheckecritical10}
The map $\bigHeckehat^{\rm ab}(\widehat{a},\tau \widehat{a})_{\ell}
\rightarrow \Sigma$ identifies the abelianized
Hecke as the upper critical locus (see Proposition \ref{relcrit})
$$ \bigHeckehat^{\rm ab}(\widehat{a},\tau \widehat{a})_{\ell} \cong
\widetilde{{\rm Crit}}\left(\bigHeckebar(a)_{\ell}/\ell,
\left(\mycal{F}_{0,\bigHeckebar(a)_{\ell}},\Phi_{0}\right)\right).
$$
These are smooth curves in $\Sigma$ that decompose into a disjoint
union of two pieces corresponding to $\widehat{a}$ and $\tau \widehat{a}$.
\end{proposition}
\begin{proof}
This is similar to the proof of Proposition \ref{abheckecritical},
using Lemma \ref{alphaatram} for smoothness of the upper critical
locus at ramification points of the horizontal divisor.

For other points, we note that the abelianized Hecke is identified
with $Y_0$ outside of $\bigExc_1(\widehat{a})$, but $\bigExc_1(\widehat{a})$
has image equal to a curve in $X_0$ that does not intersect a general
$\ell$ (as will be pointed out again in Lemma~\ref{onlymeets} below).
Thus, the inverse image of a general line is smooth in this open
subset of the abelianized Hecke. This shows that $\bigHeckehat^{\rm
  ab}(\widehat{a},\tau \widehat{a})_{\ell} $ is smooth away from the
ramification points of the horizontal divisor.

The pointwise identification between the abelianized Hecke and the
upper critical locus, whose proof is the same as for Proposition
\ref{abheckecritical}, therefore gives an identification of subschemes
(again, neither of them has embedded points). We have seen that one or
the other is smooth at every point, so the are both
smooth. Disjointness of the two pieces comes from Corollary
\ref{disjointness-cor0} .
\end{proof}

\

\

\noindent
We next turn to the calculation of the direct image. Suppose given a
spectral line bundle $\LY_1$ over $Y_1$ whose direct image to $X_1$
corresponds to our parabolic Higgs bundle.  Recall that this means,
more precisely, that $f_{1*}(\LY _1)$ is the parabolic level $0$
piece of the Higgs bundle, with parabolic level $1/2$ piece equal to
$f_{1*}(\LY_1(\ExY_1))$.

The direct image formula for the holomorphic Dolbeault complex leads
to a line bundle over $\bigHeckehat^{\rm ab}(\widehat{a})_{\ell}$. At
points where $\bigHeckehat^{\rm ab}(\widehat{a})_{\ell}$ goes through the
crossing points of the fibers over Kummer points, the contribution is
just $\LY_1 |_{\bigHeckehat^{\rm ab}(\widehat{a})_{\ell}} \otimes \omega
_{{\bigHeckebar(a)}_{\ell} / \ell}$, since the horizontal divisor does
not intervene. Recall that $\omega_{{\bigHeckebar(a)}_{\ell}/\ell}$
may also be viewed as the quotient of the forms with logarithmic
singularities along the fiber, modulo logarithmic forms from the base.

In the notations of Subsection \ref{pushforward-statements}, the
present situation is covered by the situation there, and there are no
points of type \ref{type5}. Indeed we are in the ``parabolic''
case where there is a parabolic structure with weights $0,1/2$ on the
source space, and such points are not allowed as the horizontal divisor
does not have nodes.

The direct image calculations were summarized in
Theorem~\ref{pushforward-thm} and proven in
Theorem~\ref{pushforward-appli-main} based on
Proposition~\ref{localcalc}.  Near a point $\mathbf{q}$
where the horizontal divisor (which is $\Wob_1\cap H_{\ell}$) has a
ramification over a point of a trope plane, the curve
$\bigHeckehat^{\rm ab}(\widehat{a})_{\ell}$ also passes through that
point, and the contribution for the direct image is
$$
\LY_1 |_{\bigHeckehat^{\rm ab}(\widehat{a})_{\ell}} \otimes 
\omega _{H_{\ell}/\ell}(\mathbf{q}).
$$
This formula would also hold near a point $\mathbf{q}$ where the upper critical
locus intersects the ramification divisor $\mathsf{Ram} \subset \Sigma$ over
points where the horizontal divisor is \'{e}tale over the base.  Indeed,
the proof of Theorem~\ref{pushforward-appli-main} included a
discussion of that possibility. The ramification divisor is the same
as the pullback of $\ExY_1$ to $\Sigma$. However, in fact, these
kinds of points do not occur:

\

\

\begin{lemma}
\label{onlymeets}
For general $\ell$, the curve $\bigHeckehat^{\rm
  ab}(\tilde{a})_{\ell}$ viewed inside $Y_1$, only meets the divisor
$\ExY_1$ at points $\mathbf{q}$ where the horizontal divisor has a
ramification over a point of a trope plane intersected with $\ell$.
\end{lemma}
\begin{proof}
The pullback of $\ExY_1$ to $\bigHeckehat^{\rm ab}(\widehat{a})$ is the
sum $\bigExc_0(\widehat{a}) \cup \bigExc_1(\widehat{a})$ of the strict
transform $\bigExc_1(\widehat{a})$ and the exceptional divisors over the
$16$ lines that form $\bigExc_0(\widehat{a})$. However,
$\bigExc_1(\widehat{a})$ is contracted and maps to a curve inside $Y_0$,
hence also in $X_0$. A general line does not meet this curve. Thus,
taking the inverse image of a general line $\ell$ inside
$\bigHeckehat^{\rm ab}(\widehat{a})$ and projecting back to $Y_1$, gives
a curve that only meets $\ExY_1$ at points of the $16$ lines that
generate ramification over a trope plane.
\end{proof}

\

\

\noindent
Either using this lemma, or in any case by the remark of the preceding
paragraph, the contribution to the higher direct image coming from
points near $\ExY_1$ is the restriction to the curve $G$ of
$$
\LY_1 (\ExY_1) |_{\bigHeckehat^{\rm ab}(\widehat{a})_{\ell}} \otimes 
\omega _{H_{\ell}/\ell} .
$$
This globalizes to other points of $\bigHeckehat^{\rm
  ab}(\widehat{a})_{\ell}$.  For the points lying over the Kummer, notice
that the relative dualizing sheaf is the same as the relative sheaf of
logarithmic differentials that enters into the higher direct image
calculations, so the spectral line bundle is obtained by taking
$\LY_1$ and tensoring with $\omega _{H_{\ell}/\ell} $ near these
points. Since these points are not on $\ExY_1$ (as follows from
Lemma~\ref{branchK}), the contributions for those points are given by
the same expression. Putting these all together, we obtain the
computation of the spectral line bundle:

\

\

\begin{proposition}
\label{spectrallinex1x0}
Taking the Hecke transform and restricting to the line $\ell$, the
spectral line bundle on $\bigHeckehat^{\rm ab}(\widehat{a})_{\ell}$ is
the line bundle
$$
\LY_1(\ExY_1) |_{\bigHeckehat^{\rm ab}(\widehat{a})_{\ell}} \otimes 
\omega _{H_{\ell}/\ell} .
$$
The same holds for the spectral line bundle on $\bigHeckehat^{\rm
  ab}(\tau \widehat{a})_{\ell}$, and these two are disjoint in $\Sigma$,
so this expression gives the spectral line bundle over the disjoint
union $\bigHeckehat^{\rm ab}(\widehat{a},\tau \widehat{a})_{\ell}$.
\end{proposition}

\

\subsection{Calculation of the pushforward}

Calculate in the same way as at the end of the previous chapter.
Recall that if we write $\Xi$ for the theta divisor on the Prym, or
its pullbacks to $Y_0$, $Y_1$ or $\bigHeckehat^{\rm ab}(\widehat{a},\tau
\widehat{a})_{\ell}$, we have 
$ \FxY_{0} = 2\Xi -  \ExY_{0}$ and 
$\FxY_{1} = 4\Xi -  \ExY_{1}$. 
Also
$$
\ExY_0 |_{\bigHeckehat^{\rm ab}(\widehat{a})} = \bigExc_0(\widehat{a})
$$
whereas 
$$
\ExY_1 |_{\bigHeckehat^{\rm ab}(\widehat{a})} = \bigExc_0(\widehat{a}) +
\bigExc_1(\widehat{a}),
$$
and 
$$
\FxY_0 |_{\bigHeckehat^{\rm ab}(\widehat{a})} = 2 \Xi - \bigExc_0(\widehat{a})
$$
$$
\FxY_1 |_{\bigHeckehat^{\rm ab}(\widehat{a})} = 4 \Xi - 
\bigExc_0(\widehat{a}) - \bigExc_1(\widehat{a}).
$$

\

\

\begin{lemma}
\label{reldiff0}
For the relative canonical class over the degree $0$ moduli space, we
have
$$
\omega_{\bigHeckebar (a)/X_0}= 
\Oo _{X_0}(2) |_{\bigHeckebar (a)} \otimes \Oo _{X_1}(-1)|_{\bigHeckebar (a)}
$$
and this pulls back to $\Oo _{\bigHeckehat^{\rm ab}(\widehat{a})}
(\bigExc_1(\widehat{a}) - \bigExc _0(\widehat{a}))$.
\end{lemma}
\begin{proof}
As in Lemma~\ref{reldiff1},
$$
\omega _{\bigHeckebar(a)} = 
\Oo_{X_0}(-2) |_{\bigHeckebar(a)} \otimes \Oo _{X_1}(-1)|_{\bigHeckebar (a)}.
$$
On the other hand, $\omega _{X_0} = \Oo_{X_0}(-4)$, giving the first formula. 
Then
$$
\omega _{\bigHeckebar (a)/X_0}|_{\bigHeckehat^{\rm ab}(\widehat{a})} =
\Oo (2  \FxY_0 |_{\bigHeckehat^{\rm ab}(\widehat{a})}
-  \FxY_1 |_{\bigHeckehat^{\rm ab}(\widehat{a})})
$$
$$
= \Oo (2(2 \Theta - \bigExc _0(\widehat{a})) -(4 \Xi - 
\bigExc _0(\widehat{a}) - \bigExc _1(\widehat{a}))  )=
\Oo (\bigExc _1(\widehat{a}) - \bigExc _0(\widehat{a})).
$$
\end{proof}

\

\

\noindent
Recall that the spectral line bundle on $Y_0$ is
$\LY_{0}= \Oo _{Y_0}(2 \FxY_0+  \ExY_0)$ and the spectral line bundle on 
$Y_1$ is $\LY _{1}= \Oo _{Y_1}( \FxY_1)$.
As before, we get
$$ \LY _{0}|_{\bigHeckehat^{\rm ab}(\widehat{a})} = \Oo (4\Xi -
\bigExc_0(\widehat{a}))
$$
and
$$
\LY _{1} |_{\bigHeckehat^{\rm ab}(\widehat{a})} = \Oo (4\Xi - 
\bigExc _0(\widehat{a}) - \bigExc_1(\widehat{a}) ) .
$$
From Proposition \ref{spectrallinex1x0}, the spectral line bundle of
the Hecke transform restricted over $\ell$ is the bundle
$$
\begin{aligned}
  \LY _{1} \otimes &
  \Oo _{Y_1}( \ExY_1) |_{\bigHeckehat^{\rm ab}(\widehat{a})_{\ell}} \otimes 
\omega _{H_{\ell} / \ell}  \\
& = \Oo _{\bigHeckehat^{\rm ab}(\widehat{a})_{\ell}} (4\Xi - \bigExc _0(\widehat{a})
- \bigExc _1(\widehat{a}) )
\otimes \Oo _{\bigHeckehat^{\rm ab}(\widehat{a})_{\ell}}
(\bigExc _0(\widehat{a}) + \bigExc _1(\widehat{a})) \\
& \hspace{1.5in}
\otimes 
\Oo _{H^{\rm ab}}
(\bigExc _1(\widehat{a}) - \bigExc _0(\widehat{a}))|_{\bigHeckehat^{\rm ab}(\widehat{a})_{\ell}} \\
& 
= \Oo _{\bigHeckehat^{\rm ab}(\widehat{a})_{\ell}} (4\Xi + 
\bigExc_1(\widehat{a}) - \bigExc_0(\widehat{a})).
\end{aligned}
$$
This compares with the spectral line bundle coming from $Y_0$ which is
$\Oo _{\bigHeckehat^{\rm ab}(\widehat{a})_{\ell}}(4\Xi -
\bigExc_0(\widehat{a}))$.

As before, the divisor $\ExY_1$ is the exceptional divisor of
blowing-up $Y_1$ along a copy of the curve $\Chat$. The image of the
curve in $\pp^3$ does not intersect a general line $\ell$ so
$\bigExc_1(\widehat{a})$ does not intersect $\bigHeckehat^{\rm
  ab}(\widehat{a})_{\ell}$. Thus, the spectral line bundle of the Hecke
transform coincides with the spectral line bundle coming from $Y_0$ on
$\bigHeckehat^{\rm ab}(\widehat{a})_{\ell}$.

\

\

\begin{proposition}
\label{spectraldirim10}
Define the line bundle $\Uu$ over $\bigHeckehat^{\rm ab}(\widehat{a},\tau
\widehat{a})_{\ell}$ as the pullback of $\LY_1(\ExY_1)$ on $Y_1$, tensored
with $\omega_{H_{\ell}/\ell}$. Then $\Uu$ is the spectral line bundle
for the $L^2$ Dolbeault higher direct image Higgs bundle on $\ell$.
Over each of the two pieces of $\bigHeckehat^{\rm ab}(\widehat{a},\tau
\widehat{a})_{\ell}$, $\Uu$ is isomorphic to the restriction of the
spectral line bundle that we construct for $X_0$, from $Y_0$ to
$\bigHeckehat^{\rm ab}(\widehat{a})_{\ell}$ (respectively  $\bigHeckehat^{\rm
  ab}(\tau \widehat{a})_{\ell}$).
\end{proposition}
\begin{proof}

The above calculations show this, it is also described in Section
\ref{chapter-abelianized}.  Notice that the proof in this direction is
significantly less complicated than in the previous section since we
do not need to deal with the points of type \ref{type5}.
\end{proof}

\section{The big Hecke correspondences}
\label{chapter-bighecke}

In this section we will consider the ``big'' Hecke correspondences
fitting into a diagram of the form, for the $(X_1\rightarrow X_0)$
direction:
\[
  \xymatrix@M+0.5pc@-0.5pc{
    & \bigHeckebar \ar[dl]_-{p} \ar[dr]^-{q} & \\
    X_{1} & & X_{0}\times \Cbar
    }
\]
or similarly in the opposite direction $(X_0\rightarrow X_1)$ pictured
below.  We would like to show that the local system obtained by
pulling back our local system from $(X_1,\Wob _1)$ then taking
$R^1q_{*}$, is an exterior tensor product of our local system on
$(X_0,\Wob _0)$ by the initially given rank two local system on $C$,
pulled back to $\Cbar$.

The main part of the proof on the spectral data was given in
Subsection \ref{sssec:abelianizeHecke}.  The objective of this section
is to prove some complementary statements designed to deal with
possible apparent singularities of the higher direct image operation.
For this, we introduce the notion of \emph{\bfseries effective discriminant
  divisor}, this is the part of the discriminant divisor on which the
higher direct image local system really does have singularities.

\subsection{From \texorpdfstring{$X_1$}{X1} to \texorpdfstring{$X_0$}{X0}} 

For the moment we will work in the direction from $X_1$ to $X_0$ as
pictured in the previous diagram.  Let $\Delta \subset X_0\times
\Cbar$ be the discriminant divisor of the map $q$ with respect to the
pair $(\bigHeckebar, p^{-1}\Wob _1 )$. This includes points in
$X_0\times \Cbar$ over which the map $q$ is not smooth, and points
over which the horizontal divisor $p^{-1}\Wob _1$ is not \'{e}tale. The
singularities of the $R^1q_{*}$ local system are \emph{\bfseries a priori}
contained in $\Delta$. Let $\Delta _{\rm eff}$ be the \emph{\bfseries effective}
singular divisor, namely the divisor over which the local system has
singularities. Thus $\Delta _{\rm eff} \subset \Delta$.

From the previous section we know the following statement: 

\begin{proposition}
If $a$ is a general point of $\Cbar$ then along the fiber $X_0 \times
\{ a \}$, the singular divisor $\Delta _{\rm eff}$ consists of just
$\Wob _0\subset X_0$.
\end{proposition}

\

\noindent
This proposition readily implies the following:

\

\begin{corollary}
The effective singular divisor consists of $\Wob _0 \times \Cbar$,
possibly union with a finite number of fibers of the form $X_0 \times
\{ a_i\}$ for points $a_i\in \Cbar$.
\end{corollary}
\begin{proof}
In general a divisor such as $\Delta_{\rm eff}$ in a product
$X_0\times \Cbar$ decomposes as
$$
\Delta _{\rm eff}= \Delta _{\rm eff}^{\rm vert} + 
\Delta _{\rm eff}^{\rm horiz} + \Delta _{\rm eff}^{\rm mov}
$$
where $\Delta _{\rm eff}^{\rm vert}$ is a sum of vertical
components $X_0 \times \{ a_i\}$, $\Delta _{\rm eff}^{\rm horiz}$ is a
sum of divisors of the form $D_i \times \Cbar$, and $\Delta _{\rm
  eff}^{\rm mov}$ is given by a moving family of divisors parametrized
by $\Cbar$.  Over a general point of $\Cbar$, the intersection of
these divisors with the fiber will be: for $\Delta _{\rm eff}^{\rm
  vert}$, empty; for $\Delta _{\rm eff}^{\rm horiz}$, the union of the
$D_i$; and for $\Delta _{\rm eff}^{\rm mov}$, a divisor that moves as
a function of the point. The proposition says that these all consist
of just $\Wob _0\subset X_0$, not moving as a function of the point of
$\Cbar$. It follows that $\Delta _{\rm eff}^{\rm mov} = \varnothing$ and
$\Delta _{\rm eff}^{\rm horiz} = \Wob _0 \times \Cbar$. There remains
the possibility of a nonempty $\Delta _{\rm eff}^{\rm vert}$.
\end{proof}

\

\

\noindent
In order to rule out the possibility of having a vertical piece in
$\Delta _{\rm eff}$ we'll just rule that out for the full discriminant
divisor.

\

\begin{proposition} 
\label{novert0}
The discriminant $\Delta$ of the map $q$ from $(\bigHeckebar,
p^{-1}\Wob _1 )$ to $X_0\times \Cbar$ does not contain any vertical
pieces of the form $X_0\times \{ a\}$.
\end{proposition}
\begin{proof}
We need to show that for any point $a\in \Cbar$, the full fiber
$X_0\times \{ a\}$ is not contained in the discriminant.

The point $a\in \Cbar$ corresponds to a pair $a=(A,t)$ where $A$ is a
line bundle of degree $0$ and $t\in C$ such that $A^{\otimes 2}(\pw) =
\Oo _C(t)$ (recall that $\pw$ is our fixed Weierstrass point). We need
to show that for a general point $\Ff \in X_0$, the Hecke curve
corresponding to $(\Ff , (A,t))$ is smooth and intersects $\Wob _1$
transversally.

A general $\Ff$ is stable. A bundle obtained by a Hecke transformation
is the kernel in
$$ 0\rightarrow \Ee \rightarrow \Ff \otimes A (\pw ) \rightarrow \cc
_t \rightarrow 0
$$
where the quotient corresponds to a rank $1$ quotient of the fiber
$(\Ff \otimes A(\pw ))_t$ over the point $t$. If $\Ff$ is stable of
degree $0$ then its maximal degree line subbundles have degree $-1$,
so any line subbundle of $\Ee$ must have degree $\leq 0$. Thus, $\Ee$
is stable. The Hecke curve is therefore isomorphic to the space of
such rank $1$ quotients, so it is $\pp^1$ and is hence smooth.

A little more precisely, let $\bigHeckebar^s\subset \bigHeckebar$ be
the moduli space of line bundles $B\in \Cbar$ paired with inclusions
$\Ee \subset \Ff \otimes A (\pw ) $ of colength $1$ such that $\Ee$
and $\Ff$ are stable. Let $X_0^s\subset X_0$ denote the open subset of
stable bundles.  By the previous paragraph, $\bigHeckebar^s$ is the
inverse image in $\bigHeckebar$ of $X_0^s\times \Cbar$.  The
projection from here to $X_0\times \Cbar$ is a $\pp^1$-bundle over the
open subset $X_0^s\times \Cbar$. This was also verified synthetically
in section~\ref{ssec:synth.corr}.

Thus, $\bigHeckebar$ is smooth over $X_0^s\times \Cbar$. This shows
that for any $a=(A,t)$ and for a general $\Ff$ the Hecke curve over
$(\Ff , a)$ is smooth. Next, we would like to understand its
intersection with $p^{-1}\Wob_1$. Recall that $\Ee \in \Wob _1$ if and
only if there exists a line subbundle $B\subset \Ee$ of degree $0$
such that $B^{\otimes 2}(\pw )$ is effective.  If $\Ee$ is a Hecke
transformation of $(\Ff , (A,t))$ and $B$ is such a line subbundle
then we get an injection $B\hookrightarrow \Ff \otimes A(\pw )$ or
equivalently
$$
L:= B\otimes A^{\vee}(-\pw ) \hookrightarrow \Ff .
$$
Here $L$ is a line bundle of degree $-1$. If $\Ff$ is stable of
degree $0$, that is the maximal degree of a locally free subsheaf, in
particular $L\subset \Ff$ is a saturated locally free subsheaf, i.e. a
strict subbundle.

In the other direction, given a degree $-1$ subbundle $L\subset \Ff$
there is a unique rank $1$ quotient over the point $t$ such that
$B=L\otimes A(\pw )$ maps into the kernel of $\Ff \otimes A(\pw
)\rightarrow \cc _t$. Note that 
$$
B^{\otimes 2}(\pw ) = L^{\otimes 2} \otimes A^{\otimes 2} (3\pw ) = 
L^{\otimes 2} (2\pw + t).
$$
This makes an isomorphism between the set of intersection points of
the Hecke curve with $\Wob _1$, and the set of solutions $L$ of the
pair of conditions
\begin{enumerate}
\item[(1)] $h^0(L^{\vee} \otimes \Ff ) > 0$
\item[(2)] $h^0(L^{\otimes 2} (2\pw + t)) > 0$. 
\end{enumerate}
We will look at the solutions as a subset of $\op{Jac}^{-1}(C)$ the
Jacobian of line bundles of degree $-1$ on $C$. Solutions of Condition
(1) form a divisor $D_{\Ff}$ in the linear system $|2\Theta |$.
Indeed, $D_{\Ff}$ is the Narasimhan-Ramanan point corresponding to
$\Ff$ \cite{NR}.  We have
$$
\xymatrix@M+0.5pc@R-2.5pc{
  \op{Jac}^{-1} \ar[r] & \pp^3 &
  \hspace{-3pc} = \pp\left(H^{0}(\op{Jac}^{-1}(C),\Oo(2\Theta)^{\vee}\right) \\
\cup & \cup & &  \\
D_{\Ff } \ar[r] &  H_{\Ff} & &
}
$$
where the $\pp^3 =
\pp\left(H^{0}(\op{Jac}^{-1}(C),\Oo(2\Theta)^{\vee}\right)$ here is
the dual of $X_0 =
\pp\left(H^{0}(\op{Jac}^{-1}(C),\Oo(2\Theta)\right)$ and $H_{\Ff}$ is
the hyperplane corresponding to the point $\Ff \in X_0$. The top map
is the mapping given by $|2\Theta|$ that is basepoint-free
\cite{mumford-theta1,bl}.

The second Condition (2) corresponds to the pullback of
$C\hookrightarrow \op{Jac}^1$ by the composed map
$$
\xymatrix@1@M+0.5pc@C+3pc{
\op{Jac}^{-1}(C) \ar[r]^-{(-)^{\otimes 2}} &
  \op{Jac}^{-2}(C) \ar[r]^{(-)\otimes \Oo(2\pw + t)} & \op{Jac}^{1}(C).
 }
$$
Let $\Cbar^{t}\subset \op{Jac}^{-1}$ denote this pullback curve. 

Solutions of both conditions together correspond to the intersection
of these two subspaces.  It corresponds to mapping $\Cbar^{t}
\rightarrow \pp^3$ by the composed map
$$
\Cbar^{t} \rightarrow \op{Jac}^{-1} \rightarrow \pp^3
$$
and then pulling back a hyperplane section $H_{\Ff}$. If $\Ff$ is a
general point of $X_0$ then this is a general hyperplane section.

The map $\Cbar^{t} \rightarrow \pp^3$ is nonconstant, so the pullback
of a general hyperplane section is reduced, consisting of a collection
of distinct points.  The number of points is the degree of the map,
which we claim is $16$. To prove that, note that our arguments below
will show that it can not be $>16$ otherwise the scheme-theoretical
intersection of the Hecke curve, a conic, with $\Wob_1$ would be too
big. The degree is the intersection number of the pushforward of
$2\Theta$ by the squaring map. Let us push forward the original
$\Theta$ divisor, which is just $C\subset \op{Jac}^1$. The square is the
set of divisors of the form $2x$ for $x\in C$, translated back to
$\op{Jac}^1$ as the set of divisors of the form $2x-t$. We want to know
when this is effective, i.e.  how many pairs $(x,y)$ solve
$2x-t=y$. This may be written as the equation $2x + y' = t + 2\pw$, so
it is the set of ramification points of the trigonal curve associated
to the linear system $|\Oo (t+2\pw )|$. We know that this is $8$, that
is the intersection of the pushforward of $\Theta$ with $C$ has $8$
points, so degree $\geq 8$. Thus, the intersection of the pushforward
of $2\Theta$ with $C$ has degree $\geq 16$.  This completes the proof
that the degree of the map $\Cbar ^y\rightarrow \pp^3$ is $16$.

Thus, there are $16$ distinct solutions $L$ of the two conditions.
Each point $L$ yields a line bundle $B=L\otimes A(\pw )$, with
$(B,s)\in \Cbar$ for some $s\in C$, and for distinct $L\in
\op{Jac}^{-1}$ the line bundles $B$ are distinct. They correspond
therefore to distinct points in $\Cbar$, hence to distinct points of
the normalization of $\Wob_1$ that fibers over $\Cbar$. We next note
that since our $16$ points were obtained by intersecting the image
$\Cbar \rightarrow \pp^3$ with a general hyperplane section, any pair
of two points are general with respect to each other.  Pairs of points
can be glued together under the normalization map for $\Wob_1$, but
this generality condition implies that our $16$ solutions $L$, leading
to $16$ line bundles $A$, can not contain pairs of points that are glued
together. This shows that they give $16$ distinct points in $\Wob_1$.

Now, the scheme theoretic intersection of a Hecke curve with $\Wob _1$
has length $16$, so if there are $16$ distinct points then they have
to be reduced points in the scheme theoretic intersection. This
implies that $p^{-1}\Wob _1$ is unramified, hence etale over a general
points of $X_0^s\times \Cbar$. This completes the proof of the
proposition.
\end{proof}

\

\

\begin{corollary}
\label{effwob0}
The effective singular locus $\Delta _{\rm eff}$ of the $R^1q_{*}$ local system 
on $X_0\times \Cbar$ is $\Wob _0 \times \Cbar$. 
\end{corollary}

\

\

\begin{corollary}
\label{extensor0}
The $R^1q_{*}$ local system of rank $16$ on $X_0\times \Cbar$
decomposes as an exterior tensor product of the rank $8$ local system
we have constructed on $X_0$ with a rank $2$ local system on $\Cbar$
whose spectral curve is $\Chat \rightarrow \Cbar$.
\end{corollary}
\begin{proof}
We have seen that on $X_0\times \{ a\}$ the Higgs bundle is a direct
sum of two copies of the rank $8$ Higgs bundle we construct over
$X_0$. Furthermore, the direct sum consists of two copies that are
preserved by the Higgs field in the $\Cbar$ direction, indeed the
Higgs field comes from the section of the sheaf of total differentials
on the upper relative critical locus which is a disjoint union.  The
two components vary in a covering given by $\Chat \rightarrow
\Cbar$. Also, the spectral $1$-form is the canonical $1$-form over
$\Chat$. This may be seen by restricting to horizontal copies of
$\Chat$ inside $Y_0\times \Cbar$ which map to translates of standard
copies of $\Chat$ in $Y_1$, on which the spectral $1$-form is the
canonical one for $\Chat$.

It follows that the Higgs bundle is not a direct sum of two copies of
a rank $8$ Higgs bundle on $X_0\times \Chat$. In view of the theorem
on irreducible representations of product groups, the only other
possibility is that it is an exterior tensor product. The spectral
curve of the rank $2$ local system on $\Cbar$ is $\Chat$ with
embedding given by the canonical $1$-form on $\Chat$, so this
identifies the spectral curve as a curve in $T^{\vee} \Cbar$.
\end{proof}

\

\

\subsection{From \texorpdfstring{$X_0$}{X0} to
  \texorpdfstring{$X_1$}{X1}} 
\label{ssec:degreeX0X1}

Consider next the diagram 
\[
  \xymatrix@M+0.5pc@-0.5pc{
    & \bigHeckebar \ar[dl]_-{\pzo} \ar[dr]^-{\qzo} & \\
    X_{0} & & X_{1}\times \Cbar  .
    }
\]
In this case, $\qzo$ is a fibration with fibers $\pp^1$, so there is
no discriminant for the map $\qzo$.  Let $\Delta \subset X_1 \times
\Cbar$ now denote the discriminant for the horizontal divisor $\pzo
^{-1}\Wob _0$ over $X_1\times \Cbar$. Here again, let $\Delta _{\rm
  eff}\subset \Delta$ denote the effective singularities of the local
system $R^1\qzo _{*}$ of the pullback by $\qzo$ of the local system we
construct over $(X_0,\Wob _0)$.

As before, the effective discriminant decomposes into potentially
nonempty pieces as
$$
\Delta _{\rm eff} = \Delta _{\rm eff}^{\rm vert} + 
\Delta _{\rm eff}^{\rm horiz}  + \Delta _{\rm eff}^{\rm mov} .
$$

\begin{lemma}
We have $ \Delta _{\rm eff}^{\rm mov}= \emptyset$, and
$\Delta _{\rm eff}^{\rm horiz} = \Wob _1 \times \Cbar$. 
\end{lemma}
\begin{proof}
We have seen in chapter \ref{chapter-heckex0x1} that the effective
discriminant in the fiber over a general point $a\in \Cbar$ is just
$\Wob _1 \times \{ a\} \subset X_1 \times \{ a\}$.
\end{proof}

We would like to rule out the possibility of a component $X_1\times \{ a\}$ in 
$\Delta _{\rm eff}^{\rm vert}$. As before, for this we show that
the full discriminant $\Delta$ itself 
does not contain any vertical components. We would like to show that
for any point $a\in \Cbar$, and for a general $\Ff \in X_0$, the Hecke line
associated to $(\Ff , a)$ intersects $\Wob_0$ transversally.

Since it is a line in $\pp^3$, it will be transverse to the trope
planes unless it is contained in one of them. Therefore, we would like
to show that the general Hecke line is not contained in a trope plane, and
that it intersects the Kummer in $4$ distinct points.

The bundle $\Ff$ is a stable bundle of determinant $\Oo _C(\pw )$.  If
$a=(A,t)$ with $A^{\otimes 2} (\pw ) = \Oo _C(t)$, then the bundles in
the Hecke line are kernels of the form
$$
0\rightarrow \Ee \rightarrow \Ff \otimes A \rightarrow \cc_{t}
\rightarrow 0
$$
noting that the determinant of $\Ee$ is $A^{\otimes 2}(\pw -t) = \Oo _C$. 

Notice that $\Ee$ is semistable, indeed if it had a subbundle of
degree $1$ that would give a degree $1$ subbundle of $\Ff$
contradicting stability of $\Ff$.

We have that $\Ee$ is in a trope plane if there is a line subbundle
$B\subset \Ee$ of degree $-1$ such that $B^{\otimes 2}(2\pw) = \Oo
_C$. The $16$ trope planes correspond to the $16$ solutions of this
equation. If this is the case, it gives
$$
L:= B \otimes A^{\vee} \rightarrow \Ff .
$$
The set of possible $L$'s is finite. 

If $h^0(L^{\vee} \otimes \Ff ) \geq 2$ then there is a morphism
$L^{\oplus 2} \hookrightarrow \Ff$, with cokernel a nontrivial
skyscraper sheaf. If a point $y \in C$ is in the support of this
skyscraper shead, then we will get a map $L(y)\rightarrow \Ff$,
showing that $\Ff \in \Kum$. Thus, for general $\Ff$ there is at most
one map up to scalars from $L$ to $\Ff$. For similar reasons, the
image is a saturated subsheaf of $\Ff$. Thus, for a given $L$ there is
at most one rank $1$ quotient of $(\Ff \otimes A)_{t}$ whose kernel
contains $B=L\otimes A$. As there are finitely many $L$, the Hecke
transformed bundles $\Ee$ can not all be in the wobbly locus. This shows
that, for a general $\Ff$, the Hecke line is not contained entirely in
a trope plane.

Let's now look at the intersection of the trope line with the
Kummer. The point of $X_0$ corresponding to the bundle $\Ee$ is a
point of the Kummer if $\Ee$ has a subbundle of degree $0$. So we need
to consider the possibility that there is a subbundle
$$
B\hookrightarrow \Ee 
$$
with $A$ of degree $0$. This gives
$$
L:= B\otimes A^{\vee} \hookrightarrow \Ff 
$$
with $L$ of degree $0$. In that case, note that the inclusion has
to be strict, since $\Ff$ can't contain a subbundle of degree
$1$. Then, there is a unique rank $1$ quotient of $\Ff_{t}$ such that
this subbundle corresponds to a subbundle of $\Ee$.

This reduces to our classical situation: such a subbundle $L\subset
\Ff$ corresponds to a line in $X_1$ through $\Ff$, and we know from
section~\ref{synthetic} that for a general $\Ff$ there are four
  distinct lines. So, for a general $\Ff$ there are four distinct
  subbundles $L\subset \Ff$ of degree $0$.

These in turn correspond to points $(L\otimes A) \oplus (L\otimes A)^{-1}$ 
of the Kummer.

\

\begin{lemma}
Suppose $a=(A,t)\in \Cbar$ is fixed. Choose $\Ff \in X_1$ general with
respect to $a$, and let $L_1,L_2,L_3,L_4 \subset \Ff$ be the four
subbundles of degree $0$ corresponding to the four lines through
$\Ff$. Then the four points $(L_i\otimes A)\oplus (L_i\otimes A)^{-1}$
of the Kummer are distinct.
\end{lemma}
\begin{proof}
First the $L_i$ are distinct, so $L_i\otimes A$ is not isomorphic to
$L_j\otimes A$ for $i\neq j$.  We need to show that $L_i\otimes A$ is
not isomorphic to $(L_j\otimes A)^{-1}$ for $i\neq j$.  If they were
isomorphic we would have
$$
L_i\otimes L_j \cong A^{\otimes -2} = \Oo _C(p-t).
$$
Look at a bundle $L_1$ and choose $\Ff$ general along the
corresponding line. This $L_1$ is the first subbundle of $\Ff$, the
other three being line bundles corresponding to points in one of the
fibers of the trigonal covering $C\rightarrow \pp^1$.  Here, more
precisely, if $y\in C$ is a point then $L_1^{-1}(\pw - y)$ is the
corresponding line, and this for the three points $y_1,y_2,y_3$ in the
fibers of the trigonal cover over the point $\Ff$ in the line
corresponding to $L_1$.

In particular, $L_1\otimes L_j$ is among a moving family of line
bundles, so for a general point on the line, $L_1\otimes L_j$ is not
equal to the fixed $\Oo _C(\pw -t)$.  We conclude that if $\Ff$ is
general, and contained in some line, then none of the other three
lines corresponds to the same point of the Kummer. This holds for all
the lines through the general point, concluding the proof of the
lemma.
\end{proof}

\

\

Fixing $a=(A,t)$ and for a general $\Ff \in X_1$ with respect to $a$,
then the Hecke line corresponding to $(\Ff , a)$ is transverse to the
trope planes and meets the Kummer in four distinct points.
We note that the points where the Hecke line meets the Kummer are not
on the trope conics. Indeed, for $\Ff$ general with respect to $a$,
the lines $L$ are general points of the Jacobian with respect to $A$
that is fixed, so the points of the Kummer are general.

This completes the proof of the following proposition. It gives the
same corollaries as in the degree $0$ case.

\

\begin{proposition}
\label{novert1}
The discriminant $\Delta$ does not contain any components of the form 
$X_1\times \{ a\}$ for $a\in \Cbar$. 
\end{proposition}

\

\

\begin{corollary}
\label{effwob1}
The effective singular locus $\Delta _{\rm eff}$ of the $R^1p_{*}$ local system 
on $X_1\times \Cbar$ is $\Wob _1 \times \Cbar$. 
\end{corollary}

\

\begin{corollary}
\label{extensor1}
The $R^1p_{*}$ local system of rank $16$
on $X_1\times \Cbar$ decomposes as an exterior tensor product of the
rank $8$ local system
we have constructed on $X_1$ with a rank $2$ local system on $\Cbar$
whose spectral curve is $\Chat \rightarrow \Cbar$. 
\end{corollary}

\

\subsection{Identification of the eigenvalues}

\begin{proposition}
The eigenvalue rank $2$ local systems on $\Cbar$ in Corollaries
\ref{extensor0} and \ref{extensor1} are the same as the original local
system $\Lambda$ associated to the Higgs bundle $(E,\theta )$, pulled
back to $\Cbar$.
\end{proposition}
\begin{proof}
We have seen that the spectral curve of the eigenvalue is $\Chat$. This is
the spectral curve of $(E,\theta )|_{\Cbar}$. 

The correspondence of Proposition \ref{prop:Y1toY0abHecke}, and the
analogous statement in the $(Y_0\rightarrow Y_1)$ direction, tell us
that the spectral line bundle on $\Chat$ is the same as the one for
$(E,\theta )|_{\Cbar}$, away from the ramification points of $\Chat /
\Cbar$. Our arguments did not apply to those ramification points,
because the full abelianized Hecke does not decompose into two pieces,
rather it is non-reduced over those points.

The spectral line bundle of the eigenvalue, and the spectral line
bundle of $(E,\theta )|_{\Cbar}$, are two line bundles of degree zero
that differ possibly by a divisor of degree $0$ supported on this
ramification set. But as we move around in the Hitchin base, and in
the moduli of genus $2$ curves, the set of ramification points is
permuted transitively. Thus, the coefficients of each point in the
divisor must be the same. As the divisor has degree $0$, this implies
that it vanishes, and we obtain the identification of spectral line
bundles.

It remains to identify the spectral $1$-forms. In the $(Y_0\rightarrow
Y_1)$ direction, the abelianized Hecke correspondence gives for each
point of $Y_0$ a translated map $\Chat \rightarrow Y_1$.  The spectral
$1$-form for the eigenvalue Higgs bundle is the pullback of the
spectral $1$-form on $Y_1$ to a $1$-form on $\Chat$. This is the same
as our original $1$-form, the spectral $1$-form of $(E,\theta )$. The
proof in the $(Y_1\rightarrow Y_0)$ direction is identical.

This completes the identification between spectral data for the
eigenvalue Higgs bundle and $(E,\theta )$, so the associated local
systems are isomorphic.
\end{proof}

\

\section{Construction of a third kind of local system}
\label{chapter-thirdconstruction}

In this section we take note of a different but similar construction
of a flat parabolic Higgs bundle over $(X_0, \Wob_0)$. Namely, we use
the same spectral covering $f_{0} : Y_0\rightarrow X_0$, still putting
trivial parabolic structures over the trope planes, but putting a
parabolic structure with levels $0,1/2$ over the Kummer. We observe
that this still leads to the construction of a rank $8$ local system
on $X_0-\Wob_0$, however it does not have the properties necessary to
yield a rank $8$ Hecke transform. We conjecture that this construction
gives the Hecke eigensheaf associated by geometric Langlands to an odd
degree $\pp GL(2)$ local system on $C$.

For the construction, more precisely, let $\RxY_{0} \subset Y$ be the
reduced inverse image of the Kummer surface $\Kum \subset X_{0}$. Thus
the inverse image of the Kummer is $2\RxY_{0}$ since \linebreak $f_{0} :
Y_{0}\rightarrow X_{0}$ has four simple ramification points over each
point of the Kummer.  As the Kummer has degree $4$ in $X_{0}$ we get
that $\RxY_{0} \sim 2\FxY_{0}$, where as usual $\FxY_{0}$ denotes the
pullback of the hyperplane class from $X_{0} \cong \pp^{3}$.

Let $\LY_{0}$ be a spectral line bundle over $Y_{0}$. Use this to define a
parabolic Higgs bundle $(\mycal{F}_{\bullet},\Phi )$ by
$$ \mycal{F}_a := \quad
\begin{cases} \ f _{*}(\LY_{0}), & \;\;\;\; 0
  \leq a < 1/2, \\ \ f _{*}(\LY_{0} (\RxY_{0})), & \;\;\;\; 1/2 \leq a <
  1.
\end{cases}
$$

\

\noindent
With these conventions we now have the following

\

\begin{proposition} \label{prop:X0other}
Let $\Lprym _0$ be a flat line bundle over $\Prym_2$.
If we choose \linebreak $\LY_{0} = (\blo _0^{*}\Lprym _0) \otimes
\Oo _Y(\ExY _0 +\FxY _0 )$ then 
$$
H ^2\cdot {\rm ch}^{\rm par}_1(\mycal{F}_{\bullet}) = 0
\;\; \mbox{ and }
H\cdot {\rm ch}^{\rm par}_2(\mycal{F}_{\bullet}) = 0.
$$
Furthermore, there is no need for a correction term at the
tacnodes, so we get a Higgs bundle corresponding to a local system on
$X_0 - \Wob_0$. 
\end{proposition}
\begin{proof}
The computations are left to the reader using the formulae of
Propositions \ref{todd1} and \ref{intersections1}.  We note that the
parabolic structure may be viewed as a bundle on the root stack
$X_0[\frac{1}{2}\Kum ]$. Over this root stack, a tacnode of a trope
plane with the Kummer pulls back to a normal crossings intersection.
Thus, over the root stack where the parabolic structure over the Kummer
disappears, we have a logarithmic Higgs field along a normal crossings
divisor, so no additional correction term is needed to ${\rm ch}^{\rm
  par}_2$.
\end{proof}

\

\begin{remark}
Let $X'_{0} \rightarrow X_{0}$ be the smooth double cover branched
over the Kummer. This may be viewed as the moduli space for parabolic
rank $2$ vector bundles on $\pp^1$ with parabolic structure over $6$
points (the $6$ branch points of $C/\pp^1$), via the correspondence of
Goldman and Heu-Loray \cite{HeuLoray-flat}. The parabolic Higgs bundle
$(\mycal{F}_{\bullet},\Phi )$ from Proposition~\ref{prop:X0other}
pulls back to a Higgs bundle with trivial parabolic structures over
$X'_{0}$. This is probably the Higgs bundle corresponding to a
Langlands local system over $X'_{0}$. It would go outside our current
scope to pursue this here.
\end{remark}

\

\noindent
Through the non-abelian Hodge correspondence, the parabolic Higgs
bundles $(\mycal{F}_{\bullet},\Phi )$ from
Proposition~\ref{prop:X0other} give rise to parabolic flat bundles on
$(X_{)},\Wob_{0})$ and hence to local systems on $X_{0} -
\Wob_{0}$. We expect that the middle perversity extensions of these
local systems are the Hecke eigensheaves on $Bun_{SL_{2}}$
corresponding to the $\pp GL_2$-local systems on $C$ of degree $1$ by
the geometric Langlands correspondence. Notice that since we expect to
get a perverse sheaf on $Bun_{SL_2}$,   \emph{\bfseries a priori}
there is no corresponding perverse sheaf on $X_1$.

Nonetheless, we can take the Hecke transform of this perverse sheaf
over to $X_1$. In this case, Lemma \ref{topproof} does not apply, see
Remark \ref{othercaseeff}.  The Hecke transformed local system will
therefore have singularities over $\Wob _1 \cup \KumKthree_1$ where
$\KumKthree_1$ is another subvariety of $X_1$. It is the
\emph{\bfseries Kummer K3 surface}
\cite{BeauvilleSurfaces,Dolgachev200,GH,Hudson,Keum} that, we recall,
may be described as follows: inside the Hecke space we have a K3
surface birational to the Kummer, it is the set of points where the
two lines meet in the Hecke fiber over a point of the Kummer. Then
project this to $X_1$ to get $\KumKthree_1$.

The intersection of $\KumKthree_1$ with a general line $\ell \subset X_1$
consists of the two points on $\ell$ that are the images of the points
$\bPP$ and $\bQQ$ appearing in the description of chapter
\ref{chapter-heckex0x1}. In particular $\KumKthree_1$ has degree $2$.

\

\begin{remark}
\label{thirdquadric}
Indeed, $\KumKthree_1$ is the intersection of $X_1$ with a third
quadric in $\pp^5$.  It is the K3 surface obtained by resolving the
$16$ singularities of the Kummer surface $\Kum \subset X_{0}$.  The
embedding $\KumKthree_1 \hookrightarrow X_1$ depends on the point of
$\Cbar$ over which we make the Hecke transformation, and the
intersection with the third quadric is identified in the synthetic
description with the K3 surface denoted by $\Sigma$ in \cite{GH}.
\end{remark}

\

\noindent
The Hecke transform of our rank $8$ local system on $X_0-\Wob_0$ is a
rank $16$ local system over an open subset of $X_1\times \Cbar$. A
direct cohomology indication shows that it reduces to a direct sum
of two copies of a rank $8$ local system on $X_1-\Wob_1 - \KumKthree_1$ over
each point of $\Cbar$.

Calculation of the rank indicates that the rank $8$ local system on
$X_1-\Wob_1 - \KumKthree_1$ will then Hecke transform back to a local
system of rank $24$ on an open subset of $X_0\times \Cbar$. This would
be supposed to correspond to the Hecke eigensheaf property expected of
the perverse sheaf corresponding to an odd degree $\pp GL_2$ local
system.

We close this topic for now, with the prospect of further discussion
elsewhere.

\

\section{Some pushforward calculations}
\label{chapter-pushforward}

The objective of this section is to arrive at a proof of Theorem
\ref{pushforward-thm}.  We will do that by applying the theory of
\cite{dirim} in some specific cases. The situation of Theorem
\ref{pushforward-thm} is fairly specific and tailored to our needs for
the Hecke transform calculations. We will go through some intermediate
steps of varying degrees of generality that might be of independent
interest as complements to the discussion of \cite{dirim}, for example
we extend the general theory of \cite{dirim} to the case of morphisms
with multiple fibers.

For the most part, the notations in this section will be general, not
related to our moduli spaces of stable bundles on the genus $2$ curve.

\

\subsection{The relative critical locus}

Suppose $f:X\rightarrow S$ is a projective morphism from a smooth
surface to a smooth curve. Suppose $D\subset X$ is a divisor, and let
$D_H\subset D$ be the union of components that map surjectively to
$S$. Assume that $D_H$ is reduced. Let $D^{\circ}$ be the smooth locus
of $D$ and let $X^{\circ}$ be the complement in $X$ of the singular points of
$D$.

Suppose $\Ee _{\bullet}$ is a parabolic bundle on
$(X^{\circ},D^{\circ})$ with a compatible Higgs field $\varphi$, and
let $\Sigma^{\circ} \rightarrow X^{\circ}$ be the spectral covering of
$\varphi$. Extend this to a finite covering $\Sigma \rightarrow X$. As
$\Sigma^{\circ}$ is the spectral covering of the logarithmic Higgs field
$\varphi$ acting on $\Ee_0$, there is a natural inclusion
$$
\Sigma^{\circ} \hookrightarrow T^{\vee}(X^{\circ}, \log D^{\circ}).
$$
The tautological section of $T^{\vee}(X^{\circ})|_{\Sigma ^{\circ}}$ restricts to
a section of $T^{\vee}(X^{\circ}/S)|_{\Sigma ^{\circ}}$ that we will call the 
relative tautological $1$-form, denoted $\alpha^{\rm rel}$. 

Define the \emph{\bfseries upper critical locus} 
$$
\widetilde{{\rm Crit}}(X/S,\Ee_{\bullet}, \varphi )^{\circ} \subset
\Sigma^{\circ} 
$$
to be the  zero-scheme of $\alpha^{\rm rel}$. Define the \emph{\bfseries
lower critical locus}
$$
{\rm Crit}(X/S,\Ee_{\bullet},\varphi )^{\circ}  \subset X^{\circ} 
$$
to be its image in $X^{\circ}$. Denote by 
$\widetilde{{\rm Crit}}(X/S,\Ee_{\bullet},\varphi)$ and 
${\rm Crit}(X/S,\Ee_{\bullet},\varphi)$
their closures in $\Sigma$ and $X$ respectively.

\

\noindent
We now have the following spectral description of the direct image
Higgs bundle:

\begin{proposition}
\label{relcrit}
Consider the relative $L^2$-Dolbeault complex
$$
\Dol(X^{\circ} / S, \Ee_{\bullet},\varphi) \;\;\; \mbox{ over } X^{\circ}
$$
from \cite{dirim}. Assume that the (upper or lower) critical locus
has dimension $1$. Then 
the cohomology sheaf of the $L^2$-Dolbeault complex in degree $0$ vanishes,
and the cohomology sheaf in degree $1$ (i.e. the cokernel) is supported on 
${\rm Crit}(X/S,\Ee_{\bullet},\varphi)^{\circ}$. 

Assume that $\Sigma$ is smooth in codimension $\leq 1$. There is a
dense Zariski open subset $U\subset S$ such that, the restricted map
$\widetilde{{\rm Crit}}(X/S,\Ee_{\bullet},\varphi)_{|U} \rightarrow U$
is provided with a section of the pullback of $T^{\vee}U$ making it
into the spectral variety for the restriction of the higher direct
image Higgs bundle
$\left(\mathbb{R}^{1}f_{*}
\Dol(X^{\circ}/S,\Ee_{\bullet},\varphi),f_{*}\mathfrak{d}(\varphi)\right)$
to $U$.
\end{proposition}
\begin{proof}
 See \cite{dirim}.  The differential in the $L^2$-Dolbeault complex
 $\Dol(X^{\circ} / S, \Ee_{\bullet},\varphi)$ is a morphism
 between locally free sheaves of the same rank.  Thus, the condition
 that the critical locus has dimension $1$ means that this map has
 maximal rank at the general point of $X$, in particular it is
 injective. This shows the first paragraph of the statement.

For the second part, we may assume that over $U$, the critical locus
is relatively $0$-dimensional, that it is contained in the smooth
points of $\Sigma$, that the map $X_{|U} \rightarrow U$ is smooth, and
that the horizontal divisor is \'{e}tale. We may also assume that the
critical locus does not meet the horizontal divisor.  The Higgs bundle
$\Ee$ is the direct image of a spectral line bundle on $\Sigma$, and
the cokernel of the Dolbeault complex is the spectral line bundle,
tensored with the relative differentials and then restricted to
$\widetilde{{\rm Crit}}(X/S,\Ee_{\bullet},\varphi)_{|U}$.  The higher
direct image of the $L^2$-Dolbeault complex is the usual direct image
of this line bundle on the critical locus, down to $U$. The spectral
embedding of $\Sigma$ gives a section $\Sigma \rightarrow T^{\vee}X$. In
view of the exact sequence
$$
0\rightarrow f^{*} T^{\vee} U  \rightarrow T^{\vee}X_{|U}
\rightarrow T^{\vee}(X_{|U}/U) \rightarrow 0,
$$ the critical locus being the zero set of the projection
$\Sigma_{|U} \rightarrow T^{\vee}(X_{|U}/U)$ is provided with a map
$\widetilde{{\rm Crit}}(X/S,\Ee_{\bullet},\varphi )_{|U} \rightarrow
f^{*} T^{\vee}U$.  This gives the tautological differential making the
critical locus into the spectral variety of the higher direct image
Higgs bundle \cite{dirim}.
Indeed, recalling  
that the Higgs field on the higher direct image is calculated through 
the Dolbeault analogue of the
classical connecting map for the Gauss-Manin connection, one verifies that
over $\widetilde{{\rm Crit}}(X/S,\Ee_{\bullet},\varphi)_{|U}$ the action of
tangent vectors to $U$ is given by the restriction of the tautological $1$-form, the same one
that gave the tautological embedding. 
\end{proof}

\

\begin{proposition}
\label{nonisolated}
Given a point where the map $f$ is not smooth, or where the horizontal
divisor is not \'{e}tale over the base, if the point is not in the
lower critical locus then it does not contribute a singularity of the
Dolbeault higher direct image Higgs bundle on $S$.
\end{proposition}
\begin{proof}
When we blow up to resolve the singularities of the map and apply the
technique of \cite{dirim}, such a point could result in an isolated
vertical component of the critical locus for the blown-up
map. However, this can't contribute anything to any of the parabolic
levels of the higher direct image of the $L^2$-Dolbeault complex,
since we know that the higher direct image parabolic Higgs bundle has
level pieces that are locally free over $S$, in particular they can not
have sections supported over a finite set in $S$.
\end{proof}

\

\subsection{Vertical divisors with multiple components}

We would like to understand an extended pushforward setup in which we
allow for non-reduced vertical divisors.  Let $f:X \rightarrow S$ be a
projective morphism, with $S$ a smooth curve and $X$ a smooth surface.
Assume that the discriminant of $f$ is a single point $o \in S$.
Suppose that $D_V = f^{-1}(o)$ is a simple normal crossings divisor,
possibly non-reduced.  Write the components as $D_1,\ldots, D_r$ with
multiplicities $m_1,\ldots , m_r$.

Suppose given $D_H$ a smooth divisor that is transversal to $f$ in
the following sense:
\begin{itemize}
\item $D_H$ meets reduced
  components of fibers of $f$ transversally at smooth points.
\item $D_H$ does not pass through the nodes of $D_{V}$.
\item At points where $D_H$
intersects a component of $D_V$ that has multiplicity $>1$, the map
$D_H \rightarrow S$ is ramified, with order of ramification equal
to the multiplicity of $D_V$ at that point.
\end{itemize}
These conditions imply that $D = D_{H} + D_{V}$ is a 
normal crossing divisor with multiplicities whose horizontal piece is \'{e}tale over $S -
\{o\}$. The reduced divisor $D^{\rm red} = D_H + D^{\rm red}_V$ is a strict
normal crossings divisor. 

\

Suppose we are given a logarithmic parabolic Higgs bundle
$(E_{\bullet},\varphi)$ on $X$ with respect to the reduced divisor
$D^{\rm red}$.  We assume the usual hypothesis on $\varphi$
(parabolic with nilpotent graded parts of the residue).

Let $E_0$ be the bundle obtained by assigning parabolic levels $0$
along all components.  The Higgs field is
$$
\varphi : E_0 \rightarrow E_0\otimes \Omega^1_X(\log D^{\rm red}).
$$
The residue of $\varphi$ along $D_H$ induces a weight filtration on
$E|_{D_H}$.  We denote by $W_{i}E$ the subsheaf of $E$ of sections that
restrict to sections of $E|_{D_H}$ that are in $W_i(E|_{D_H})$.

The weight filtration extends over the crossing points as a strict
filtration, in the same way as was discussed \cite{dirim} in the
original case of a reduced vertical divisor.  Locally at a point where
the horizontal and vertical divisors meet we can take a root of the
function $f$ and reduce this question to the original case.

Define 
$$
\Omega^1_{X/S}(\log D) := \Omega^1_X(\log D) / f^{*} \Omega^1_S(\log o ).
$$
This is a line bundle over $X$. Note that since $\Omega^{1}_{X}(\log D) =
\Omega^{1}_{X}(\log D^{\rm red})$, we have a surjective homorphism
$\Omega^{1}_{X}(\log D^{\rm red}) \to \Omega^1_{X/S}(\log D)$, and so
can project $\varphi$ to a Higgs field with coefficients in
$\Omega^1_{X/S}(\log D)$.

\

\

\noindent
{\em Caution:} The line bundle  $\Omega^1_{X/S}(\log D)$
    is not the same as the relative dualizing sheaf
$\omega _{X/S}(D_H)$ if any multiplicities of $D_V$ are $> 1$.  The
relation between these two is given by
$$
\omega _{X/S}(D_H) = \Omega ^1_{X/S}(\log D)  \otimes \Oo _X(D - D^{\rm red}).
$$
See for an example the footnote on page \pageref{footnote7}.

\

\

The projection of $\varphi$ to a Higgs field with coefficients in the
vertical cotangent bundle provides a map
$$
W_{i}E_0 \rightarrow W_{i-2} E_0 \otimes \Omega ^1_{X/S}(\log D)
$$
which we will still denote by $\varphi$.

\

\begin{theorem}
\label{multidol}
Define the relative Dolbeault complex,  in parabolic level $0$,  to be 
$$
\mathsf{Dol}(X/S,E)_{0} = \mathsf{Dol}(E_0,\varphi) := 
\left[
  W_0 E_0 \stackrel{\varphi}{\longrightarrow} W_{-2} E_0 \otimes
  \Omega ^1_{X/S}(\log D) 
\right] .
$$
This is quasi-isomorphic to the modified relative Dolbeault complex 
$$
\mathsf{Dol}'(E_0,\varphi )= 
\mathsf{Dol}'(X/S,E)_{0} = \left[
  W_1 E_0
  \stackrel{\varphi}{\longrightarrow}
  W_{-1} E_0 \otimes \Omega ^1_{X/S}(\log D)  
\right] .
$$
Define 
$$
F^i_0 := \mathbb{R}^i f_{*} \mathsf{Dol}(X/S,E)_0 \cong \mathbb{R}^i f_{*}
\mathsf{Dol}'(X/S,E)_0.
$$
Then $F^i_0$ is the parabolic level $0$ piece of the parabolic
Higgs bundle on $S$ corresponding to the $i$-th higher direct image
local system of the local system corresponding to
$(E_{\bullet},\varphi)$.
\end{theorem}

\

\noindent
In order to get the other parabolic level pieces of $F^{i}_{\bullet}$ we
proceed as follows.  Consider the parabolic line bundle $\Oo _S(a\cdot
o)$ that has a jump at parabolic level $-a$.

We get a pullback parabolic bundle $\Oo _X(a\cdot D_V):= f^{*} \Oo
_S(a\cdot o)$.  We can tensor this with $E_{\bullet}$ to get the
parabolic Higgs bundle
$$
(E(a\cdot D_V)_{\bullet},\varphi).
$$
One should be careful that the component sheaves along components
of $D_V$ of higher multiplicity need to be calculated using the
correct notions of pullback and tensor product of parabolic bundles.

Then we can also form the complexes 
$$
\mathsf{Dol}(E(a\cdot D_V)_0,\varphi ) \;\;\;\; \mbox{and} \;\;\;\; 
\mathsf{Dol}'(E(a\cdot D_V)_0,\varphi )
$$
in the same way as before: 
$$
\mathsf{Dol}(E(a\cdot D_V)_0,\varphi ) = 
\left[
  W_0 E(a\cdot D_V)_0
  \stackrel{\varphi}{\longrightarrow}
  W_{-2} E(a\cdot D_V)_0 \otimes \Omega ^1_{X/S}(\log D) 
\right] 
$$
and 
$$
\mathsf{Dol}'(E(a\cdot D_V)_0,\varphi ) = 
\left[
  W_1 E(a\cdot D_V)_0 \stackrel{\varphi}{\longrightarrow}
  W_{-1} E(a\cdot D_V)_0 \otimes \Omega ^1_{X/S}(\log D) 
\right] 
$$
and again these are quasi-isomorphic. 

\

Then 
$$
F^i_{a} = \mathbb{R}^i f_{*} \mathsf{Dol}(E(a\cdot D_V)_0,\varphi )
\cong \mathbb{R}^i f_{*} \mathsf{Dol}'(E(a\cdot D_V)_0,\varphi )
$$
is the level $a$ piece of the parabolic Higgs bundle on $S$
corresponding to the $i$-th higher direct image local system of the
local system corresponding to $(E_{\bullet},\varphi)$.

\

We will prove Theorem~\ref{multidol} in the next section. In
section~\ref{sec:applications} we will use Theorem~\ref{multidol} to
describe a Dolbeault complex which algebraically computes the $L^{2}$
pushforward of a parabolic Higgs bundle with parabolic structure and
poles on a normal crossing divisor whose horizontal part is allowed to
be ramified over the base.

\

\subsection{Proof of \texorpdfstring{Theorem~\ref{multidol}}{Theorem 12.3}}

For any point $p\in D_i\cap D_j$ let $k(p)$ be the $\mathsf{lcm}$  of the
multiplicities $m_i$ and $m_j$.  Let $t$ be a local coordinate on $S$
vanishing at $o$. Consider a cover $\widetilde{S}\rightarrow S$ given
by $t=s^n$ with $n$ chosen so that for any crossings point $p$ in
$D_{V}$ we have that $k(p)$ divides $n$. Let $\tilde{o}$ be the point of $\widetilde{S}$
over $o \in S$. Let $\widetilde{X}$ be the
normalization of $X\times _S\widetilde{S}$.

\

\begin{proposition}
\label{anprop}
Over any crossing point $p$, the space $\widetilde{X}$ has a
singularity of type $A_{m-1}$ where $m=n/k(p)$.  Let $\widehat{X}
\rightarrow \widetilde{X}$ be the minimal
resolution of all these points.  Then
$$
\widehat{f} : \widehat{X}\rightarrow \widetilde{S}
$$
is a map from a smooth surface to $\widetilde{S}$ such that
$\hat{f}^{-1}(\tilde{o})$ is a reduced normal crossings divisor.
\end{proposition}

\

This will be seen in subsection \ref{combocase} below.

\

\noindent
Let $G$ be the cyclic group of symmetries of $\widetilde{S}/S$.  Let
$\eta : \widehat{X} \rightarrow X$ denote the composition map
$\widehat{X} \to \widetilde{X} \to X$.  Suppose we are given a parabolic
Higgs bundle $(E_{\bullet},\varphi )$ on $X$ with respect to a divisor $D$
including $D_V=f^{-1}(o)$ and a horizontal part $D_H$ not in our
picture.

We have defined above the Dolbeault complex $\mathsf{Do}l(X/S,E)_0$
for parabolic level $0$ in the case of multiple components. On the
other hand, if $\eta^{*}E_{\bullet}$ is the pullback parabolic bundle
on $\widehat{X} $, we can also form the parabolic level $0$ Dolbeault
complex upstairs for $\mathsf{Dol}(\widehat{X}
/\widetilde{S},\eta^{*}E)_0$.  This is for a case of a reduced
divisor.

\

\begin{theorem}
\label{coverdol}
With these notations,  we have
$$
(\mathbb{R}\eta _{*}(\mathsf{Dol}(\widehat{X}/\widetilde{S},\eta^{*}E)_0 ))^G
= \mathsf{Dol}(X/S,E)_0.
$$
\end{theorem}

\

\noindent
The proof of this theorem will occupy subsections~\ref{nkpcase},
\ref{alreadyred}, and \ref{combocase}.  Theorem~\ref{coverdol} implies
Theorem \ref{multidol}, using Lemma \ref{parzerolemma} of the next
subsection.

\

\subsubsection{A lemma on  parabolic structures on
  \texorpdfstring{$S$}{S}} Suppose we have a parabolic bundle $F$ on
$S$, and consider a base change $\xi : \widetilde{S} \rightarrow S$
given by $t=s^n$.  Recall that $\tilde{o}$ is the point of $\widetilde{S}$
sitting over $o \in S$, and let $G$ be the cyclic group acting on
$\widetilde{S}/S$.

\

\begin{lemma}
\label{parzerolemma}
We have
$$
\xi_{*} ((\xi^{*}F)_0 )^G \cong F_0.
$$
\end{lemma}
\begin{proof}
We may assume that $F = \Oo_S (a\cdot o)$ is a parabolic line bundle. 

Then 
$$
\xi ^{*}F = \Oo _{ \widetilde{S}}(na\cdot \tilde{o}).
$$
Thus 
$$
(\xi ^{*}F)_0 = \Oo _{ \widetilde{S}}(\floor{na}\cdot \tilde{o}).
$$
Let us write things in terms of modules.  We have that 
$$
(\xi ^{*}F)_0 \leftrightarrow s^{-j}\mathbb{C}[s] , \;\;\;\; j = \floor{na} .
$$
We then think of this as a $\mathbb{C}[t]$ module. 

The $G$ fixed part is the sum of all the monomials $t^i=s^{ni}$ that
are contained in this module.  Being contained in here is equivalent to $ni
\geq - \floor{na}$, so in terms of monomials in $t$ we have the sum of
the $t^i$ whenever
$$
i \geq  - \floor{na} / n.
$$
On the other hand, $F_0 = \Oo _S(\floor{a}\cdot 0)$ and this is the
sum of monomials $t^i$ for $i\geq - \floor{a}$.  Thus, the statement
of our lemma is equivalent to saying, for integers $i$, that
$$
i \geq -  \floor{na} / n \Leftrightarrow i\geq -\floor{a}.
$$
Changing the sign of $i$, this is equivalent to the statement for all $i$
$$
i \leq   \floor{na} / n \Leftrightarrow i\leq \floor{a},
$$
which in turn is equivalent to
$$
\floor{  \floor{na} / n } = \floor{a}.
$$
Write $a= \floor{a} + (b/n) + c$ where $0\leq b < n$ is an integer
and $0\leq c < 1/n$ is a real number.  Then
$$
\floor{na} = n \floor{a} + b
$$
and 
$$
\floor{ \floor{na} / n } = \floor{ n\floor{a} / n + b/n} =
\floor{a} + \floor{b/n} = \floor{a}.
$$
\end{proof}

\

\subsubsection{Proof of \texorpdfstring{Theorem~\ref{coverdol}}{Theorem 12.5}:
  the local case when \texorpdfstring{$n=k(p)$}{neqkp}}
\label{nkpcase}

We begin the proof of Theorem~\ref{coverdol} by first examining the
case when $n = k(p)$. Consider the local picture at a crossings point
$p$ for a covering with $n=k(p)$.  Let $D_1$ and $D_2$ denote the two
components of $D_{V}$ that meet at $p$. The coordinates on $X$ are
$x,y$ in our local neighborhood, with $D_1$ defined by $x=0$ and $D_2$
defined by $y=0$. Let $m_i$ be the multiplicity of $D_i$ in
$f^{-1}(o)$, so we can assume that $f$ is given by $t=x^{m_1}
y^{m_2}$.

Write $n= m_1a = m_2b$ with $a$ and $b$ relatively prime. Let 
$\widetilde{X}\rightarrow X$ be the covering of order $ab$
given in coordinates by
$$
x=u^a, \;\;\;\; y = v^b.
$$
The map $\widetilde{X}\rightarrow S$ is given in coordinates by 
$$
t= u^n v^n
$$
so it factors through the map $\widetilde{X}\rightarrow \widetilde{S}$
given by $s=uv$. 

If $d$ is the $\mathsf{gcd}$ of $m_1$ and $m_2$, the normalization of
the fiber product
$$
X\times _S \widetilde{S}
$$
consists of the disjoint union of $d$ copies of the above covering
$\widetilde{X}$.

We can take the relative Dolbeault complex on the normalization of the
fiber product, then take the direct image by an \'{e}tale map down to
$\widetilde{X}$, then take the $\zz / d\zz$ invariants; this is the
same as the relative Dolbeault complex of $\widetilde{X} /
\widetilde{S}$.

So, to show the desired property in this case we can consider the
direct image from $\widetilde{X}$ to $X$ and take the $\zz / ab
\zz$-invariants.

The pullback of the logarithmic forms on $X/S$ is the same as the
logarithmic forms on $\widetilde{X}/\widetilde{S}$.

Lemma \ref{parzerolemma} applied in each coordinate yields the
statement that, under push-forward and taking invariants,
the parabolic level $0$ of the pullback parabolic
bundle goes  to the level $0$
piece on $X$. This yields the required statement for
$\widetilde{X}\rightarrow X$.

\subsubsection{Proof of
  \texorpdfstring{Theorem~\ref{coverdol}}{Theorem 12.5}:
  the local case when the divisor is already reduced}
\label{alreadyred}

Suppose now we are  given $f: X\rightarrow S$ that locally looks like
$(x,y)\mapsto t=xy$. Consider a covering $\widetilde{S}\rightarrow S$
given by $t=s^n$.  Then $X\times _S \widetilde{S}$ has a single
$A_{n-1}$ surface singularity over the node of $f^{-1}(0)$.  Indeed
the equation is $xy=s^n$.

Resolve this by the minimal resolution to get 
$$
\xymatrix{
  \widehat{X} \ar[r] \ar[d]  & X\times _S \widetilde{S} \ar[r] \ar[d]
  & X  \ar[d] \\
\widetilde{S} \ar@{=}[r]  & \widetilde{S} \ar[r] &  S
}
$$

\

\noindent
{\bfseries Claim 1}: the fiber of $\widehat{f}: \widehat{X} \rightarrow
\widetilde{S}$ over $s=0$ is reduced and of the form
$$
\widehat{X}_0 \cup E_1 \cup \cdots \cup E_{n-1}
$$
where $\widehat{X}_0$ is the strict transform of $f^{-1}(0)$ and 
$E_1,\ldots , E_{n-1}$ form an $A_{n-1}$ configuration of $(-2)$ curves.
This is a classical statement \cite{duVal}, see for example 
\cite{Brieskorn}.

\

\noindent
{\bfseries Claim 2}: Let $G$ be the cyclic group of symmetries of
$\widetilde{S}/S$.  Let $\eta : \widehat{X} \rightarrow X$ be the map.
Suppose given a parabolic Higgs bundle $(E_{\bullet},\varphi )$ on $X$ with
respect to a diviror $D$ including $D_V=f^{-1}(0)$ and a horizontal
part $D_H$ not in our picture.  Form the Dolbeault complex
$\mathsf{Dol}(X/S,E)_0$ for parabolic level $0$.

Similarly, if $(\eta^{*}E_{\bullet},\eta^{*}\varphi)$ is the pullback
parabolic bundle on $ \widehat{X}$ we can also form the parabolic
level $0$ Dolbeault complex
$\mathsf{Dol}(\widehat{X}/\widetilde{S},\eta ^{*}E)_0$ upstairs.

We claim that the natural map 
$$ \mathsf{Dol}(X/S,E)_0 \rightarrow (\mathbb{R} \eta _{*}
(\mathsf{Dol}(\widehat{X} /\widetilde{S},\eta ^{*}E)_0 ) )^G
$$
is a quasi-isomorphism.

We prove this in the following way. We show that the map induces an
isomorphism on the $G$-invariant $R\eta _{*}$ of each piece of the
Dolbeault complex.  The statement is local at the normal crossings
point, and in particular it does not concern the horizontal
divisor. Also, since the statement is made for each individual term of
the Dolbeault complex, it does not concern the Higgs field which
provides the differential between the individual pieces.  Thus, we may
decompose $E_{\bullet}$ into a direct sum of parabolic line bundles.
A parabolic line bundle is determined by the parabolic levels along
the two divisor components $D_1$ and $D_2$. Any parabolic levels can
occur, as may be seen by considering a toy example where the
fundamental group of the total space is $\zz \times \zz$.

We can apply the general theory for reduced normal crossings divisors,
to both $X\rightarrow S$ and to $\widehat{X} \rightarrow
\widetilde{S}$. We know in both cases that the higher direct image of
the Dolbeault complex calculates the parabolic bundle, respectively
$F_{\bullet}$ and $\widehat{F}_{\bullet}$, corresponding to the higher
direct image local systems.  We have seen in Lemma \ref{parzerolemma}
that taking the $G$-invariants of the direct image of $\widehat{F}_0$
from $\widetilde{S}$ down to $S$ yields $F_0$. Therefore, the natural map 
$$
\mathsf{Dol}(X/S,E)_0 \rightarrow 
(\mathbb{R}\eta _{*} (Dol(\widehat{X} /\widetilde{S},\eta ^{*}E)_0 ))^G
$$
induces a quasi-isomorphism between the higher direct images on $S$. 
Thus, if we denote the
two terms of  $\mathsf{Dol}(X/S,E)_0$
by $\mathsf{Dol}(X/S,E)_0^i$ for $i=0,1$, and the two terms of $\mathsf{Dol}
(\widehat{X} /\widetilde{S},\eta ^{*}E)_0$ by $\mathsf{Dol}
(\widehat{X} /\widetilde{S},\eta ^{*}E)_0^i$ for $i=0,1$,
we get that
$$
Rf_{*} \left(\mathsf{Dol}(X/S,E)_0^i\right) \rightarrow Rf_{*}\left(
\left(\mathbb{R} \eta _{*} (\mathsf{Dol}
(\widehat{X} /\widetilde{S},\eta ^{*}E)_0^i) \right)^G\right)
$$
is a quasi-isomorphism for $i=0,1$. 

Let $C^i$ be the cone 
$$
\mathsf{Dol}(X/S,E)^i_0 \rightarrow 
(R \eta _{*} (\mathsf{Dol}(\widehat{X} /\widetilde{S},\eta ^{*}E)^i_0 ) )^G
\rightarrow C^i .
$$
It is a complex of sheaves on $X$ supported at the crossings point
$p\in X$. From the above, we conclude that
$$
Rf_{*} (C^i) = 0.
$$
This implies that $C^i$ is quasi-isomorphic to $0$. Indeed, the
cohomology sheaves of $C^i$ are coherent sheaves on $X$ concentrated
at $p$ so their higher direct images in strictly positive degrees
vanish. The Leray spectral sequence going from the higher direct
images of the cohomology sheaves to the cohomology of the higher
direct image, therefore starts out with only a single line, so it
degenerates right away. If any of the cohomology sheaves were nonzero
this would give a nonzero $Rf_{*} (C^i) $ contradicting the previous
statement; thus, $C^i$ must be acyclic.

\

The cone being acyclic implies that the map 
$$
\mathsf{Dol}(X/S,E)^i_0 \rightarrow 
\left(  R\eta _{*}
\left(\mathsf{Dol}(\widehat{X} /\widetilde{S},\eta ^{*}E)^i_0 \right) \right)^G
$$
is a quasi-isomorphism. 
This holds for the case of a parabolic line bundle,
and hence for the case of any parabolic bundle.

This is a local statement at the point $p$.
In the presence of a
Higgs field,  the hypercohomology
spectral sequence
for $\mathbb{R}\eta_{*}$ of the Dolbeault complex implies that the map
$$
\mathsf{Dol}(X/S,E)_0 \rightarrow 
(\mathbb{R}\eta _{*} (Dol(\widehat{X} /\widetilde{S},\eta ^{*}E)_0 ) )^G
$$
is a quasi-isomorphism. This proves {\bfseries Claim 2}.

\

\subsubsection{Proof of
  \texorpdfstring{Theorem~\ref{coverdol}}{Theorem 12.5}:
  Combining the cases}
\label{combocase}

Suppose now that we are given a projective map $f: X\rightarrow S$
such that the vertical divisor $D_{V} =f^{-1}(o)$ is simple normal
crossings, with smooth components having various multiplicities
$m_i$. Let $n$ be a number divisible by all the $m_i$. Let
$\widetilde{S}\rightarrow S$ be a cyclic covering of degree $n$ fully
ramified over $o$.

Let $\widetilde{X}$ be the normalization of $X\times_S \widetilde{S}$.
For $p\in D_i \cap D_j$, let $k(p)$ denote the $\mathsf{lcm}$ of $m_i$
and $m_j$. Then $k(p)$ divides $n$ and we set $G^k:= \zz / k(p)\zz$
and $G':= \zz / (n/k(p))\zz$. Let $X'$ be the intermediate covering of
$X$ of degree $k(p)$, over $S'\rightarrow S$ of degree $k(p)$.  The
map $X'\rightarrow X$ is a finite covering with $X'$ smooth near the
preimages of $p$ and the fiber of $f':X'\rightarrow S'$ has reduced
components with normal crossings passing throught the preimages of
$p$, as we saw in subsection \ref{nkpcase}. Note that $p$ splits into
several points $p'$, the number is the $\mathsf{gcd}$ of the
multiplicities at $p$.

Then, the map  $\widetilde{X}\rightarrow X'$ is the covering 
considered in subsection \ref{alreadyred} of degree $n/k(p)$. Over each
point $p'$ we obtain a point $\widetilde{p}\in \widetilde{X}$
with an $A_{n/k(p)-1}$ singularity.  This shows Proposition \ref{anprop}. 

Let $\widehat{X}$ be obtained by taking, at each point $\tilde{p}$,
the minimal resolution of the $A_{n/k(p)-1}$ singularity. Then the map
$\widehat{f}: \widehat{X} \rightarrow \widetilde{S}$ has a reduced
normal crossings fiber over $\tilde{o} \widetilde{S}$.

We show the desired property locally at each point $p$. 
The covering $S'\rightarrow S$ and $X'\rightarrow X$ are covered by
the case of subsection \ref{nkpcase}. 
Thus, the map
$$
\mathsf{Dol}(X'/S')_0 \rightarrow 
(\mathbb{R}\eta '_{*} Dol(\widehat{X} / \widetilde{S})_0)^{G'}
$$
is a quasi-isomorphism, 

The map $X'\rightarrow S'$ has normal crossings fibers,
so the covering $\widetilde{S}\rightarrow S'$ and
upstairs $\widehat{X} \rightarrow \widetilde{X} \rightarrow X'$ are covered
by the already-reduced case of subsection \ref{alreadyred}. Therefore 
$$
\mathsf{Dol}(X/S)_0 \rightarrow 
(\mathbb{R}\eta ^k_{*}Dol(X'/S')_0)^{G^k}
$$
is a quasi-isomorphism. Putting these together we conclude the desired
statement that
$$
\mathsf{Dol}(X/S)_0 \rightarrow 
(\mathbb{R}\eta _{*} \mathsf{Dol}(\widehat{X} / \widetilde{S})_0)^{G}
$$
is a quasi-isomorphism.

\

\noindent
This completes the proof of Theorem~\ref{coverdol} and hence of
Theorem~\ref{multidol}.

\

\subsubsection{Application to the higher direct image}

Let $F_{\bullet}$ be the parabolic Higgs bundle on $S$ corresponding
to the higher direct image local system from $X$.  Then $\xi
^{*}F_{\bullet}$ is the parabolic Higgs bundle on $\widetilde{S}$
corresponding to the higher direct image local system coming from the
local system on $\widehat{X}$ that corresponds to the pullback Higgs
bundle $\eta ^{*}E_{\bullet}$.

On the other hand, since the fiber $\widehat{f}^{-1}(\widehat{o})$ is
reduced normal crossings and the horizontal divisor meets it
transversally, we know that this is the same as the parabolic bundle
on $\widetilde{S}$ corresponding to the higher direct image of the
Dolbeault complex $\mathsf{Dol}(\widehat{X} /\widetilde{S},\eta
^{*}E)_{\bullet}$.  In particular this gives, for the level $0$ pieces,
that
$$ (\xi ^{*}F)_0 = \mathbb{R} \widehat{f}_{*} \mathsf{Dol}(\widehat{X}
/\widetilde{S},\eta ^{*}E)_0.
$$
Lemma~\ref{parzerolemma} says that
$$
\begin{aligned}
F_0 & = \xi _{*} ((\xi ^{*}F)_0)^G \\
&  =  \left[ 
\xi _{*} 
\left( 
\mathbb{R} \widehat{f}_{*} \mathsf{Dol}(\widehat{X} /\widetilde{S},\eta ^{*}E)_0 \right) 
\right]  ^G.
\end{aligned}
$$
We can write
$$
\xi _{*} 
\left( 
\mathbb{R} \widehat{f}_{*}
\mathsf{Dol}(\widehat{X} /\widetilde{S},\eta ^{*}E)_0 \right) 
= 
\mathbb{R}f_{*} \left(
\mathbb{R}\eta _{*} \mathsf{Dol}(\widehat{X} /\widetilde{S},\eta ^{*}E)_0 
\right).
$$
Thus
$$
F_0 = \left[ \mathbb{R}f_{*} \left( \mathbb{R}\eta _{*}
  \mathsf{Dol}(\widehat{X} /\widetilde{S},\eta ^{*}E)_0 \right)
  \right] ^G.
$$
We can put the $G$-invariants inside before taking
$\mathbb{R}f_{*}$,  so this becomes
$$
F_0 =
\mathbb{R}f_{*} \left[ 
\left(
\mathbb{R}\eta _{*} Dol(\widehat{X} /\widetilde{S},\eta ^{*}E)_0 \right) ^G
\right]  .
$$
Now the theorem says
$$
\left(
\mathbb{R}\eta _{*} \mathsf{Dol}(\widehat{X} /\widetilde{S},
\eta ^{*}E)_0 \right) ^G
= \mathsf{Dol}(X/S,E)_0.
$$
Thus we get
$$
F_0 =
\mathbb{R}f_{*} 
\mathsf{Dol}(X/S,E)_0  .
$$
This is the statement we want.  The analogous
statement holds with $\mathsf{Dol}$ replaced by $\mathsf{Dol}'$.

\

\subsection{Applications} \label{sec:applications}

Consider the following situation: $f:X\rightarrow S$ is a map from a
surface to a curve, with a divisor $D\subset X$ such that in local
coordinates $(x,y)$ we have $f(x,y) = y$, with $y$ being a local
coordinate on $S$ vanishing at $o \in S$, and so that $D$ is given by
$y=x^2$. We will denote the origin point $(x,y) = (0,0)$ by $P \in X$.

Let $E_{\bullet}$ be a parabolic logarithmic Higgs bundle on $X$ with
singularities along $D$. Although in the end we need to apply the
theory to Higgs bundles of higher rank ($8$ to be exact), the local
picture involves a piece of rank $2$, so we will suppose here that
$E_{\bullet}$ has rank $2$.

We would like to identify a complex on $X$ that calculates the higher
direct image of the $L^2$ Dolbeault complex on a resolution. For this,
let $\beta : \widetilde{X}\rightarrow X$ be obtained by blowing up $X$
twice, first at the origin $P \in X$ and then at the resulting triple
intersection point of the exceptional line with the strict transforms
of the fiber $\mathsf{Fib}$  and $D$.

\begin{figure}[!ht]
\begin{center}
\psfrag{X}[c][c][1][0]{{$X$}}
\psfrag{Xt}[c][c][1][0]{{$\widetilde{X}$}}
\psfrag{Y}[c][c][1][0]{{$S$}}
\psfrag{p}[c][c][1][0]{{{\footnotesize $P$}}}
\psfrag{f}[c][c][1][0]{{$f$}}
\psfrag{ft}[c][c][1][0]{{$\widetilde{f}$}}
\psfrag{F}[c][c][1][0]{{\footnotesize $\mathsf{Fib}$}}
\psfrag{Ft}[c][c][1][0]{{\footnotesize $\widetilde{\mathsf{Fib}}$}}
\psfrag{beta}[c][c][1][0]{{$\beta$}}
\psfrag{D}[c][c][1][0]{{{\footnotesize $D$}}}
\psfrag{Dt}[c][c][1][0]{{{\footnotesize $\widetilde{D}$}}}
\psfrag{A}[c][c][1][0]{{{\footnotesize $A$}}}
\psfrag{B}[c][c][1][0]{{{\footnotesize $B$}}}
\psfrag{o}[c][c][1][0]{{$o$}}
\epsfig{file=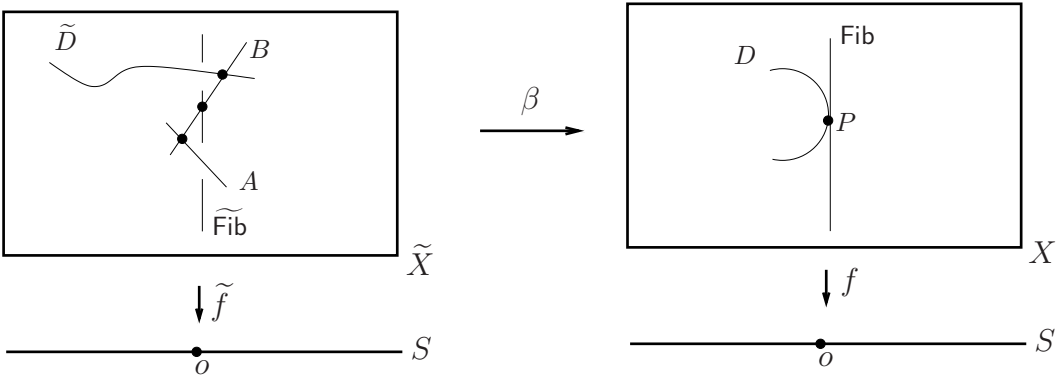,width=5.5in} 
\end{center}
\caption{Divisor configuration after two blowups}\label{fig:2blowups}  
\end{figure}

Notationally, call $\widetilde{D}$ the strict transform of $D$ in
$\widetilde{X}$ and $\widetilde{\mathsf{Fib}}$ the strict transform of
$\mathsf{Fib}$ in $\widetilde{X}$.  Let $A$ be the strict transform of
the first exceptional divisor, and $B$ the second exceptional divisor.
Thus, $\widetilde{D}$ meets $B$ transversally at a single point, $A$
meets $B$ transversally at a single point, and
$\widetilde{\mathsf{Fib}}$ meets $A$ transversally at a single point,
as depicted in Figure~\ref{fig:2blowups}.  The parabolic divisor for
$\beta^{*}E_{\bullet}$ is $A + 2B + \widetilde{D}$, and the fiber of
the the map $\widetilde{f} : \widetilde{X} \rightarrow S$ over the
origin is equal to $A + 2B+\widetilde{\mathsf{Fib}}$.

\

\

In other words, the fiber of $\widetilde{f} : \widetilde{X}
\rightarrow S$ over the origin $o \in S$
is a normal crossings divisor, but with
nontrivial mulitplicity $2$ on the middle component. The higher direct
image statement of Theorem \ref{multidol}, for maps whose fibers have
nontrivial muliplicity, yields the Dolbeault complex on
$\widetilde{X}$ that calculates the higher direct image to $Y$.  Then
take its higher direct image down to $X$, this is the Dolbeault
complex on $X$ to be identified.

We will consider two cases, both for the local situation when
$E_{\bullet}$ has rank $2$.  One case is for parabolic levels $0,-1/2$,
the other is for trivial parabolic structure but a nontrivial
nilpotent residue of the Higgs field along $D$.

One main idea for doing the calculations is to write the rank $2$
bundles as a direct sum of two rank $1$ pieces, even though such a
decomposition is not compatible with the Higgs field. This allows for
computation of the pieces in holomorphic Dolbeault complexes, which
are then put back together before inputting the Higgs field.

\

\subsubsection{Parabolic levels \texorpdfstring{$0,-1/2$}{0andhalf}}
\label{parcase}

In this case the level zero piece of $E_{\bullet}$
is a rank $2$ bundle $E = E_{0}$ on $X$,
provided with $U\subset E_{|D}$ of rank $1$.  Let $Q:= E_{|D}/U$ be
the quotient line bundle over $D$.  The parabolic structure is given
by $E_0=E$ and $E_t = \ker (E \rightarrow Q)$ for $-1/2 \leq t < 0$,
then $E_t = E(-D)$ for $-1\leq t < -1/2$.

Assume given a logarithmic Higgs field $\varphi$ that is strictly
parabolic for the filtration, so $\varphi : E \rightarrow E_{-1/2}
\otimes \Omega ^1_X(\log D)$.  In this case, the direct image
$F_{\bullet}$ on $Y$ is going to have trivial parabolic structure so
we just want to calculate $F_0$.

The $L^2$ Dolbeault complex is defined in the usual way of
\cite{dirim} outside of the origin in $X$. It has a unique extension
to a two-term complex of locally free sheaves, the \emph{\bfseries
locally free  extension of the Dolbeault complex} expressed as
$$
\mathsf{Dol}(E_{\bullet},X/S) _{\rm lf}
:=
\left[
E_0 \rightarrow E_{-1/2} \otimes \omega _{X/S}(D)
\right] .
$$
Use the notations for the blow-up
$\widetilde{X}\stackrel{\beta}{\rightarrow} X$ established above, with
exceptional divisors $A$ and $B$.

Let $\widetilde{E}$ be the pullback of $E=E_0$ to $\widetilde{X}$.
This is not quite the same as the level zero value of the pullback
parabolic bundle.  We note that $\widetilde{E}$ is constant on $A$ and
$B$, equal to the vector space $E_P$ with its quotient $Q_P$.  The
parabolic structure of $\beta^{*}E_{\bullet}$ along $A$ is just the
pullback one, with 
$$
\begin{aligned}
(\beta^{*}E_{\bullet})_{A,0} &  = \widetilde{E}, \ \text{and}  \\
  (\beta^{*}E_{\bullet})_{A, -1/2} &  = \ker \left((\beta^{*}E_{\bullet})_{A,0}
  \rightarrow Q_{P,A}\right),
\end{aligned}
$$
where $Q_{P,A}$ means the vector space $Q_P$ considered as a trivial
bundle on $A$.

Along $B$ we get extra sections in $(\beta^{*}E_{\bullet})_{B,0}$,
namely in a neighborhood of  a general point of $B$, we have
$$
(\beta^{*}E_{\bullet})_{B,0} =
\ker \left( \widetilde{E} (B) \rightarrow Q_{P,B}(B) \right) .
$$
The parabolic structure is trivial along $B$, in other words there
are no non-integral levels.  This is because the levels $1/2$ on both
divisors $A$ and $\widetilde{D}$ combine together to give a piece
without parabolic structure (it may be seen in greater detail by
looking at $E_{\bullet}$ as a direct sum of two parabolic line bundles).

Let $E'$ denote the level $0$ piece of the parabolic bundle
$\beta^{*}E_{\bullet}$ on $\widetilde{X}$.  It then has a weight
filtration, with $W_0E' = E'$ and $W_{-2} E' = \ker (E' \rightarrow
Q'_{\widetilde{D}} )$ where $Q'$ means the quotient adjusted to be a
quotient of $E'$.

Applying the general theory we have
$$
\widetilde{\mathsf{Dol}} =
\left[
E' \rightarrow 
W_{-2} E' \otimes \Omega_{\widetilde{X} / S} (\widetilde{D})
\right] .
$$

In our situation where $B$ has multiplicity two,  we have 
$$
\Omega_{\widetilde{X}/S}^{1} = \beta^{*} \Omega_{X/S}^{1}(A + B)
$$
since it should have degree $-1$ on $A$ and degree $0$ on $B$,
and should agree
with the relative dualizing sheaf 
\label{footnote7}
except\footnote{  
For comparison, the calculation of the relative dualizing sheaf is as follows:
We have 
$$
\omega_{\widetilde{X}/S}  = \beta ^{*} \omega _{X/S} (A + 2B).
$$
To see this, note that we divided by $\omega_S$ on both sides so we
can just do the calculation for the $\omega_X$ and
$\omega_{\widetilde{X}}$.

We will write $\beta^{A} : X^{A} \to X$ for the blowup of $X$ at the
origin $P \in X$ and let $A_{X^{A}} \subset X^{A}$ denote the
exceptional divisor for this map.  Similarly we will write $\beta^{B}
: \widetilde{X} \to X^{A}$ for the blowup of $X^{A}$ at the point in
$X^{A}$ where the exceptional divisor $A_{X^{A}}$ and the proper
transforms of $D$ and $\mathsf{Fib}$ intersect, and as before we will
write $B \subset \widetilde{X} $ for the exceptional divisor of
$\beta^{B}$ and $A$ for the proper transform of $A_{X^{A}}$. With this
notation we then have $\beta = \beta^{A}\circ \beta^{B}$.

Suppose we blow up the origin  in a surface $X$  with coordinates $x,y$,
yielding
new coordinates $u,v$ with $x=uv$ and $y=v$ for example.   Then 
$$
dx\wedge dy = (u dv + v du) \wedge dv = v du \wedge dv .
$$
This gives 
$$
\begin{aligned}
\omega _{X^A/S}  & = (\beta^{A})^{*} \omega _{X/S} (A_{X^A}), \quad \text{and} \\
\omega _{\widetilde{X}/S} & = (\beta^{B})^{*} \omega _{X^A/S} (B).
\end{aligned}
$$
However,  note that 
$$
\beta^{B,*} \Oo _{X^A} (A_{X^A}) = \Oo _{\widetilde{X}} (A + B).
$$
Thus we get 
$$
\omega _{\widetilde{X}/S} = \beta ^{*} \omega _{X/S} (A + 2B),
$$
and hence
$$
\omega _{\widetilde{X}/S} = \Omega_{\widetilde{X}/S}^{1}(B).
$$
} 
along $B$.  Note that $A^2=-2$ while $B^2=-1$.

Using the birational transformation $\beta : \widetilde{X} \rightarrow
X$, we would like to calculate the complex
$$
\mathbb{R}\beta _{*} \widetilde{\mathsf{Dol}} 
= \mathbb{R}\beta _{*} \left[
E' \rightarrow 
W_{-2} E' \otimes \Omega_{\widetilde{X}/S}^{1} (\widetilde{D})
\right] 
$$
on $X$.  It will be the same as the previous complex away from the
origin, and it will turn out also that the components are locally free
at the origin.

Let, as in the footnote, $\beta^A : X^A \rightarrow X$
be the first blowup,  and
$$
\beta^B : \widetilde{X} \rightarrow X^A
$$
denote the second blowup.
Start by looking at $\mathbb{R}\beta^B_{*} \widetilde{\mathsf{Dol}}$
as a complex on $X^A$. 

The first claim is that near points of $A$ that are not the
intersection point with $B$ (which may be identified via $\beta^{B}$
with the points of $A_{X^{A}}$ which are not the triple intersection
point with the strict transforms of $D$ and $\mathsf{F}$), this complex is the
same as the $\beta^{A}$-pullback of the locally free extension
$\mathsf{Dol}(E_{\cdot},X/S) _{\rm lf}$.

In degree $0$ we have $E'$ on the one hand, and the
$\beta^{A}$-pullback of $E$ on
the other hand.  These are the same along points of $A$.

In degree $1$, noting that our points of $A$ under consideration do
not touch $\widetilde{D}$, we have on the one hand
$E'\otimes \Omega_{X^A/S}^{1}$ and on the other hand,
the pullback of $E_{-1/2} \otimes
\Omega _{X/S}^{1} (D)$.  These are again the same, so that shows the first
claim.

The remaining question is what happens over $B$.  In order to
calculate the terms of the complex, let us suppose that $E$ has rank
$1$. The rank $2$ case with trivial Higgs field is a direct sum of two
parabolic line bundles. This will serve to calculate the two terms of
the complex. Then the Higgs field gives a map between these two and we
get the desired higher direct image. A standard hypercohomology
spectral sequence argument then shows the desired quasiisomorphism in the
presence of the Higgs field.

If the line bundle has trivial parabolic structure: then 
$E'$ is the same as $\widetilde{E}$. 
The degree $0$ piece is $\widetilde{E}$ and the degree $1$ piece is
$$
\widetilde{E} (-\widetilde{D})\otimes
\Omega_{\widetilde{X}/S}^{1}(\widetilde{D})
= \widetilde{E} \otimes \Omega_{\widetilde{X}/S}^{1} 
= \beta ^{*}(E\otimes \Omega_{X/S})(A+B).
$$
The higher direct image of $\widetilde{E}$ down to $X^A$ is $E^A$
(the $\beta^{A}$-pullback of $E$ to $X^{A}$), and then the higher direct
image down to $X$ is $E$ respectively   $E\otimes \Omega _{X/S}^{1}$.
For these,
note that on $\pp^1$, the $H^1$ of a trivial bundle vanishes; and
furthermore, along the blowing-down of an exceptional $\pp^1$, the
higher direct image of a bundle that is trivial along the exceptional
divisor is the same as the usual pushforward.

Now suppose we have a level $-1/2$ parabolic structure.  This means
$E_t = E$ for $-1/2 \leq t \leq 0$ and $E_t = E(-D)$ for $-1 \leq t <
-1/2$.  Here $W_0E=E$ and $W_{-2} E = E$.

From the parabolic structure we get
$$
E' = \widetilde{E} (B).
$$
Thus, the degree $0$ piece is $\widetilde{E}(B)$ and the degree $1$
piece is $\widetilde{E}(B) \otimes \Omega _{\widetilde{X}/S}$ (here as
above the $-\widetilde{D}$ and $\widetilde{D}$ cancel).

For the degree $0$ piece,
$$
\mathbb{R}\beta _{*}  \widetilde{E} (B) = E.
$$
This is because the bundle $\widetilde{E} (B)$ restricted to the
exceptional $B\cong \pp^1$ is of the form $\Oo _{\pp^1}(-1)$ that has
vanishing $H^1$ so the higher direct image is zero, and the direct
image is the space of sections of $E$ by Hartogs' theorem.

Also note that 
$$
\beta ^{*} \Oo _X (D) = (\beta^{B})^{*} \Oo _{X^A} (D_{X^A} + A_{X^A})
= \Oo _{\widetilde{X}}(\widetilde{D} + A + 2B).
$$
Thus,  for the degree $1$ piece we get
$$
W_{-2}E' \otimes \Omega _{\widetilde{X}/S}^{1}(\widetilde{D}) 
=
\beta^{*}(E \otimes  \Omega _{X/S}^{1}) (A + 2B + \widetilde{D})
= \beta ^{*}(E \otimes  \omega _{X/S}(D)).
$$
We first take the $\mathbb{R}\eta^B_{*}$ from $\widetilde{X}$ down to $X^A$. 
This yields 
$$
(\beta^{A})^{*}(E \otimes  \omega _{X/S})(D))
$$
and then taking the direct image down to $X$ gives 
$E \otimes  \omega _{X/S}(D)$.  

In conclusion, back to the case when $E$ has rank two so it is a
direct sum of line bundles for the two cases discussed above, if we
define the locally free extension Dolbeault complex
$\mathsf{Dol}(E_{\bullet},X/S) _{\rm lf}$ as at the start of this subsection,
then the higher direct image of $\mathsf{Dol}(E_{\bullet},X/S) _{\rm lf}$ from
$X$ to $S$ calculates $F_0$. For the line bundle $L$ with parabolic
level $0$, $E_0=L$ is the bundle and $E_{-1/2}=L(-D)$ so $E_{-1/2}
\otimes \omega _{X/S}(D) =L\otimes \omega _{X/S}$.  For the line
bundle $L'$ with parabolic level $-1/2$ we have $E_0=L'$ and $E_{-1/2}
= L'$ so $E_{-1/2} \otimes \omega _{X/S}(D) =L'\otimes \omega
_{X/S}(D)$.

\subsubsection{Alternative method}
\label{parcase-alt}

For comparison, we also do the calculations by going to a double cover
$\alpha : Z\rightarrow X$ ramified over $D$, with involution $\sigma :
Z \rightarrow Z$.  Let $R\subset Z$ be the upper ramification divisor
mapping isomorphically to $D$. We will write $f_{Z} f\circ \alpha : Z
\to S$ for the natural map to $S$ and will denote the fiber of $f_{Z}$
over $o \in S$ by $\mathsf{Fib}_{Z}$. 

The parabolic pullback of $E$ to $Z$ is a bundle that we denote $E_Z$,
given by
$$
E_Z = \ker \left(  \alpha ^{*} (E) (R) \rightarrow Q(R) \right)  .
$$
This has a Higgs field without poles.  It projects to a relative
Higgs field
$$
\varphi_Z : E_Z \rightarrow E_Z \otimes \omega_{Z/S} ,
$$
where we recall that $\omega_{Z/S} = \Omega^1_Z(\log \mathsf{Fib}_{Z})/f_Z^{*}
\Omega^1_S(\log o)$.

The upstairs Dolbeault complex is 
$$
\mathsf{Dol}_Z = \left[
E_Z \rightarrow  E_Z \otimes \omega _{Z/S}
\right] .
$$
This has an action of $\sigma$ covering the action on $Z$.  We can
take the direct image down to $X$, and consider the $\sigma$-invariant
part
$$
(\alpha_{*} \mathsf{Dol}_Z)^{+} .
$$
This is a direct summand in the complex of locally free sheaves
$\\alpja_{*} \mathsf{Dol}_Z$
so it is itself a complex of locally free sheaves.

The first term in the complex is just 
$$
(\alpha_{*} E_Z)^{+} = E.
$$
We also note that
$$
\omega_{Z/S} = \alpha^{*} \omega_{X/S}(R).
$$
Thus,  the second term in the complex is
$$
(\alpha_{*}((\alpha^{\mathrm{par},*}(E\otimes\omega_{X/S}))(R)).
$$
We have
$$
(\alpha^{\mathrm{par},*}( E) \otimes \alpha^{*}\omega_{X/S}))(R) =
\alpha^{*}\omega_{X/S}) \otimes 
\ker \left(
\alpha^{*}(E) (2R) \rightarrow Q(2R)
\right)
$$
and when we take the invariant part of the direct image back to $X$ we get
$$
\omega_{X/S} \otimes 
\ker ( E(D) \rightarrow Q(D) )
= E_{-1/2} \otimes \omega _{X/S} (D).
$$
This is the second sheaf in the Dolbeault complex, and it is indeed
the reflexive extension of the sheaf we would get away from the
ramification point of $D$ by taking $W_{-2} E \otimes \Omega
^1_{X/S}(\log D)$.

So, we conclude in this calculation that the Dolbeault complex on
$X/S$ is just the reflexive i.e. locally free extension of the one we
would get by the usual formula away from the ramification point.

\

\subsubsection{The nilpotent case}
\label{nilpotentcase}

Let us go back to the geometric setup from the beginning of
Section~\ref{sec:applications}. That is, we have the map $f : X \to
S$, the smooth divisor $D \subset X$ simply ramified over $o \in
S$ and the
double blowup modification
\[
\xymatrix@R=0.1pc{
  \widetilde{D} & & D \\ 
  \cap & & \cap \\ 
  \widetilde{X} \ar[rddd]_-{\widetilde{f}}
  \ar[rr]^-{\beta} & & X \ar[lddd]^-{f} \\
  \\
  \\
  & S & 
}
\]
where $\widetilde{D}$ denotes the strict transform of $D$ in
$\widetilde{X}$. 

We start with a logarithmic Higgs bundle on $(X,D)$ which has no
parabolic structure and has a nilpotent residue along $D$, and we want
to describe a $L^{2}$ Dolbeault complex on $X$ that computes the
pushforward of such Higgs bundle to $S$. In our case the Higgs bundle
will be of rank $8$ but locally near the origin $P \in X$ it will
decompose into a direct sum of rank two pieces so that for each piece
the Higgs field will have a residue which is a $2\times 2$ nilpotent
Jordan block.

Thus, to understand the requisite Dolbeault complex near $P \in X$ it
is enough to consider the case of a rank $2$ Higgs bundle
$(E,\varphi)$ on $X$ with nontrivial nilpotent residue and no
parabolic structure along $D$.  We still want to decompose $E$ into
line bundles locally near the origin $P \in X$ in a robust way. In
other words, we would like to pay attention to some additional
structure on $E$, and then work with a decomposition into a direct sum
of line bundles in which the line bundle summands carry the same type
of structure in a compatible way. Normally we use parabolic structures
for this, as we did for example in section~\ref{parcase}. In the
current setting the parabolic structure on $E$ is trivial, so instead
we will work with the weight filtration induced from $\varphi$.

For this, introduce the notion of weight-filtered bundle, that is a
bundle $E$ with weight filtrations (inspired by the monodromy weight
filtration from mixed Hodge theory) on the parabolic graded pieces. In
this case the parabolic level filtration is trivial so the parabolic
graded piece is just the restriction of the bundle $E$ to $D$.  The
weight filtration on $E$ is the filtration by locally free subsheaves
that are the pullbacks of the steps of the (monodromy)  weight filtration on
$E_{|D}$ given by the nilpotent endomorphism $\op{Res}_{D} \varphi :
E_{|D} \to E_{|D}$. The weight filtration $W_{\bullet}E$ determines
the pieces of the Dolbeault complex, and we can look at a direct sum
decomposition of $E$ into line bundles with weight filtartions.
In our case the ``mixed Hodge'' weights will be $1$ and $-1$, since
the residue of the Higgs field is a nonzero nilpotent $2\times 2$
matrix so its monodromy weight filtration has those weights.

Again we will write $E'_{\bullet}$ for the parabolic pullback of $E$
to $\widetilde{X}$.  In particular the bundle $E'_0$ is the same as
$\widetilde{E}=\beta^{*}E$ since $E$ had no parabolic structure to
begin with. Similarly if we split $E$ locally into a sum of two weight
filtered line bundles, the paraboolic level zero of the parabolic
pullbacks of each of the line bundles will just be equal to the
pullback of that line bundle.

Let as before $\mathsf{Fib} \subset X$ denote the fiber of $f$ over
the origin $o \in S$, and let $\widetilde{\mathsf{Fib}}$ be its strict
transform in $\widetilde{X}$.

Putting parabolic level $-1/2$ along the
fiber of $\widetilde{f} : \widetilde{X} \to S$ over $o \in S$ yields  
$$
E'_{-1/2} = \widetilde{E}(-A-B-\widetilde{\mathsf{Fib}})
$$
Note that $E'_{-1}$ is $\widetilde{E}
(-A-2B-\widetilde{\mathsf{Fib}})$, that is just
$\beta^{*}(E(-\mathsf{F}))$.

Recall that 
$$
\Oo_{\widetilde{X}}(-\widetilde{D}) = \beta^{*}(\Oo _X(-D)) (A + 2B).
$$
Also similarly 
$$
\Oo_{\widetilde{X}}(-\widetilde{\mathsf{Fib}}) =
\beta^{*}(\Oo_X(-\mathsf{Fib})) (A + 2B).
$$
Consider now a line bundle $E$ with weight $1$.  Then $W_0E =
E(-D)$ and $W_{-2}E = E(-D)$ (away from the origin $P \in X$).

Then, the degree $0$ piece in the Dolbeault complex of $E'_{\bullet}$
on $\widetilde{X}/S$ is
$$ W_0E'_0 =\beta^{*}(E(-D))(A+2B) =
(\beta^{B})^{*}(((\beta^{A})^{*}E(-D))(A_{X^A}))(B)
$$
since  $(\beta^{B})^{*}(\Oo_{X^A}(A_{X^A}) = \Oo _{\widetilde{X}}(A + B)$. 
So the direct image down to $X^A$ is locally free and then the direct
image down to $X$ is locally free, and we get
$$
R\beta _{*}  W_0E'_0 = E(-D).
$$
For level $-1/2$ we have
$$
\begin{aligned}
W_0E'_{-1/2} & =\beta^{*}(E(-D))(A + 2B - A - B -
\widetilde{\mathsf{Fib}}) \\
& = \beta^{*}(E(-D))(B - \widetilde{\mathsf{Fib}} ) \\
& = \beta^{*} (E(-D-\mathsf{Fib}))(A + 3B).
\end{aligned}
$$
In this case,  there is an $R^1\eta _{*}$ term.
Because of that we will do a modified
version of the calculation later.

Look now at the degree $1$ piece of the Dolbeault complex of
$E'_{\bullet}$ on $\widetilde{X}/S$.  We have
$$
W_{-2}E'_0 =\beta^{*}(E(-D))(A+2B) .
$$
Also recall
$$
\Omega_{\widetilde{X}/S}^{1}(\widetilde{D}) =
\beta^{*} \Omega _{X/S}^{1} (A + B + \widetilde{D})
=  \beta^{*} (\Omega_{X/S}^{1}(D))(-B).
$$
Putting these together, the degree $1$ piece of the relative
Dolbeault complex of $E'_{\bullet}$ is
$$
W_{-2}E'_0\otimes \Omega_{\widetilde{X}/S}^{1}(\widetilde{D})
=
\beta^{*}(E\otimes \Omega_{X/S}^{1}))(A+B) .
$$
When we push down to $X^A$ then to $X$ the result is locally free,
equal to
$$
R\beta _{*}(W_{-2}E'_0\otimes \Omega_{\widetilde{X}/S}^{1} (\widetilde{D}))
= E\otimes \Omega _{X/S}.
$$
For the parabolic level $-1/2$ piece we have
$$
\begin{aligned}
W_{-2}E'_{-1/2} & = \beta^{*}(E(-D))(A+2B)(-A-B-\widetilde{\mathsf{Fib}})  \\
& = \beta^{*}(E(-D))(B-\widetilde{\mathsf{Fib}}),
\end{aligned}
$$
and so 
$$
\begin{aligned}
W_{-2}E'_{-1/2}\otimes & \Omega _{\widetilde{X}/S} (\widetilde{D}) = \\
& = W_{-2}E'_{-1/2}\otimes \beta^{*} (\Omega_{X/S}^{1}(D))(-A-2B)(A+B) \\
& = \beta^{*}(E\otimes \Omega_{X/S}^{1})(- \widetilde{\mathsf{Fib}}) \\
& = \beta^{*}(E\otimes \Omega_{X/S}^{1}(-\mathsf{Fib}))(A + 2B).
\end{aligned}
$$
Done differently we note that 
$$
\begin{aligned}
W_{-2}E'_{-1/2}\otimes & \Omega_{\widetilde{X}/S}^{1} (\widetilde{D})  = \\
& = \widetilde{E} (-A-B- \widetilde{\mathsf{Fib}}  - \widetilde{D})\otimes
\Omega_{\widetilde{X}/S}^{1}(\widetilde{D}) \\
& = \widetilde{E} (-A-B-
\widetilde{\mathsf{Fib}})\otimes \Omega_{\widetilde{X}/S}^{1} \\
  & = \widetilde{E} (-A-B-  \widetilde{\mathsf{Fib}})\otimes
    \beta^{*}(\Omega_{X/S}^{1})(A + B) \\
& = \beta^{*}(E\otimes \Omega_{X/S}^{1}(-\mathsf{Fib}))(A + 2B).
\end{aligned}
$$
The higher direct image down to $X$ is locally free as seen by doing
in two stages, and
$$
R\beta _{*} W_{-2}E'_{-1/2}\otimes \Omega_{\widetilde{X}/S}^{1} (\widetilde{D})
= E\otimes \Omega _{X/S}(-\mathsf{Fib}).
$$

\

\noindent
Let us now  look at the case of weight $-1$.
In this case $W_{-2}E' = E'(-\widetilde{D})$ and 
$W_0E=E$.

\

\noindent
The degree zero piece in parabolic level zero yields just
$$
\beta^{*}(E)
$$
and the direct image down to $X$ is just $E$.

Consider the parabolic level $-1/2$ piece.   Upstairs on $\widetilde{X}$
this is
$$
\beta^{*}(E)(-A-B-\widetilde{\mathsf{Fib}}) = \beta^{*}(E(-\mathsf{Fib}))(B).
$$
This has direct image $E(-\mathsf{Fib})$ down on $X$.  

Look now at the degree $1$ piece.  We have as before
$$
W_{-2}E'_0 =\beta^{*}(E(-D))(A+2B) .
$$
Also recall
$$ \Omega_{\widetilde{X}/S}^{1} (\widetilde{D}) = \beta^{*} \Omega _{X/S}^{1}
(A + B + \widetilde{D}) = \beta^{*}(\Omega_{X/S}^{1}(D))(-B).
$$
Putting these together,  the degree $1$ piece is
$$
W_{-2}E'_0\otimes \Omega_{\widetilde{X}/S}^{1} (\widetilde{D})
=
\beta^{*}(E\otimes \Omega_{X/S}^{1}))(A+B) .
$$
As before when we push down to $X^A$ then to $X$ the result is
locally free, equal to
$$
R\beta_{*} W_{-2}E'_0\otimes \Omega_{\widetilde{X}/S}^{1}(\widetilde{D})
= E\otimes \Omega_{X/S}^{1}).
$$
For the parabolic level $-1/2$ piece we have
$$
W_{-2}E'_{-1/2} =\beta^{*}(E(-D))(A+2B)(-A-B-\widetilde{\mathsf{Fib}} ) 
= \eta ^{*}(E(-D))(B- \widetilde{\mathsf{Fib}}),
$$
so 
$$
\begin{aligned}
W_{-2}E'_{-1/2}\otimes \Omega_{\widetilde{X}/S}^{1}(\widetilde{D}) & = \\
& = \beta^{*}(E\otimes \Omega_{X/S}^{1})(- \widetilde{\mathsf{Fib}}) \\
& = \beta^{*}(E\otimes \Omega_{X/Y}^{1}(- \mathsf{Fib}))(A + 2B).
\end{aligned}
$$
Done differently we note that 
$$
\begin{aligned}
W_{-2}E'_{-1/2}\otimes & \Omega_{\widetilde{X}/S}^{1}(\widetilde{D})
 = \\
& = \widetilde{E} (-A-B- \widetilde{\mathsf{Fib}} - \widetilde{D})\otimes
\Omega_{\widetilde{X}/S}^{1}(\widetilde{D}) \\
& = \widetilde{E} (-A-B- \widetilde{\mathsf{Fib}} )\otimes
  \Omega_{\widetilde{X}/S}^{1} \\
  & = \widetilde{E} (-A-B- \widetilde{\mathsf{Fib}} )\otimes
    \beta^{*}(\Omega _{X/S}^{1})(A + B) \\
    & = \beta^{*}(E\otimes \Omega_{X/S}^{1}(- \mathsf{Fib}))(A + 2B).
\end{aligned}
$$
The higher direct image down to $X$ is locally free as seen by doing
in two stages, and
$$
R\beta _{*} W_{-2}E'_{-1/2}\otimes \Omega_{\widetilde{X}/S}^{1} (\widetilde{D})
= E\otimes \Omega _{X/S}(-\mathsf{Fib}).
$$

\

\subsubsection{Modified version}
\label{nilcase-modone}

Since we had an $R^1\beta _{*}$ term, let us do an alternative version
of the calculation.  We note that, in the case where the 
weighted bundles come from a $\varphi$, we could replace $W_0E$ and
$W_{-2}E$ by $W_1E$ and $W_{-1}E$. in degrees $0$ and $1$
respectively.  This observation originally due to Zucker \cite{Zucker}
for $L^2$ Dolbeault complexes of VHS was explained for the present
setting in \cite{dirim}.  So let us look at those.

Suppose $E$ has weight $1$.  Then (away from the origin $P \in X$)
$W_1E = E$ and $W_{-1}E = E(-D)$.  We did this calculation above (for
the mixed weight $-1$ case) and for the degree $0$ part we got:
$$
W_1E'_0  = \beta^{*}(E), \qquad \quad
R\beta_{*} W_1E'_0 = E
$$
and
$$
W_1E'_{-1/2} = \beta ^{*}(E)(-A-B-\widetilde{\mathsf{Fib}}) =
\beta ^{*}(E(-\mathsf{Fib}))(B)
$$
giving $R\beta _{*} W_1E'_{-1/2} = E(-\mathsf{Fib})$. 

\

For the degree $1$ part we got:
$$
W_{-1}E'_0\otimes \Omega_{\widetilde{X}/S}^{1} (\widetilde{D})
= \beta^{*}(E\otimes \Omega_{X/S}^{1}))(A+B) 
$$
so 
$$
R\beta_{*} W_{-1}E'_0\otimes \Omega_{\widetilde{X}/S}^{!} (\widetilde{D})
=
E\otimes \Omega _{X/S}^{1}.
$$
At parabolic level $-1/2$ we have
$$
W_{-1}E'_{-1/2} \otimes \Omega_{\widetilde{X}/S}^{1}(\widetilde{D})
= \beta^{*}(E\otimes \Omega_{X/S}^{1}(-\mathsf{Fib}))(A + 2B),
$$
so 
$$
R\beta_{*} W_{-1}E'_{-1/2} \otimes \Omega_{\widetilde{X}/S}^{1} (\widetilde{D})
= E\otimes \Omega_{X/S}(-\mathsf{Fib}).
$$

\

Turn now to the case of mixed weight $-1$.  In this case, 
$$
W_1E = E, \qquad W_{-1} E = E.
$$
For the degree $0$ piece the calculation is the same as above: 
$$
R\beta _{*} W_1E'_0 = E, \qquad \quad
R\beta _{*} W_1E'_{-1/2} = E(-C).
$$
For degree $1$ we have 
$$
\begin{aligned}
W_{-1}E'_0\otimes & \Omega_{\widetilde{X}/S}^{1} (\widetilde{D})
=  \\
& = \beta^{*}(E) \otimes \Omega_{\widetilde{X}/S}^{1} (\widetilde{D}) \\
& = \beta^{*}(E \otimes \Omega_{X/S}^{1}(D))(A + B -A -2B) \\
& = \beta^{*}(E \otimes \Omega_{X/S}^{1}(D))(-B).
\end{aligned}
$$
The higher direct image down to $X^A$ is 
$$
(\beta^{A})^{*}(E \otimes \Omega_{X/S}^{1}(D))\otimes \mathcal{I}_{P_A}
$$
where $P_A \in X^{A}$ is the point that is blown up the second
  time, and $\mathcal{I}_{P_{A}}$ denotes its ideal sheaf .  This fits
  in an exact sequence
$$
\xymatrix@1@M+0.4pc@C+1pc{
0\ar[r] & (\beta^{A})^{*} (E \otimes \Omega_{X/S}^{1}(D))\otimes
\mathcal{I}_{P_A}
\ar[r] &  (\beta^{A})^{*}(E \otimes \Omega _{X/S}^{1}(D))
\ar[r] & \cc_{P_A} \rightarrow 0.
}
$$
Taking the direct image down to $X$ we get a long exact sequence that
shows
$$
R^1(\beta^A)_{*}\left((\beta^{A})^{*}(E \otimes
\Omega_{X/S}^{1}(D))\otimes \mathcal{I}_{P_A}\right)
$$
as the cokernel of the map 
$$
R^0(\beta^A)_{*}\left((\beta^{A})^{*}(E \otimes\Omega_{X/S}^{1}(D)\right)
\longrightarrow
\cc_{P} 
$$
induced by the previous map (with $P \in X$ as before denoting the origin).
We have 
$$
R^0(\beta^A)_{*}\left((\beta^{A})^{*}(E \otimes \Omega_{X/S}^{1}(D)\right) =
E \otimes \Omega_{X/S}^{1}(D)
$$
and a local section not vanishing at the  origin corresponds  to a section of 
$(\beta^{A})^{*}(E \otimes \Omega_{X/S}^{1}(D))$ that does not vanish at $P_A$.
Therefore,  this map is surjective and the cokernel is $0$.  We get
$$
R^1(\beta^A)_{*}\left((\beta^{A})^{*}(E \otimes
\Omega_{X/S}^{1}(D))\otimes \mathcal{I}_{P_A}\right) = 0
$$
and 
$$
R^0(\beta^A)_{*}\left((\beta^{A})^{*}(E \otimes
\Omega_{X/S}^{1}(D))\otimes \mathcal{I}_{P_A}\right) = E \otimes
\Omega_{X/S}^{1}(D)\otimes \mathcal{I}_{P}.
$$
We conclude that in the weight $-1$ case,   for parabolic level $0$ and
in degree $1$,
$$
R\beta_{*} W_{-1}E'_0\otimes
\Omega_{\widetilde{X}/S}^{1}(\widetilde{D}) = E \otimes \Omega
_{X/S}(D)\otimes \mathcal{I}_{P}.
$$
Look at the case of parabolic level $-1/2$.  Now,  
$$
\begin{aligned}
W_{-1}E'_{-1/2}\otimes & \Omega_{\widetilde{X}/S}^{1}(\widetilde{D})
=  \\
& = \beta^{*}(E) \otimes \Omega_{\widetilde{X}/Y} (\widetilde{D})
(-A-B-\widetilde{\mathsf{Fib}}) \\
& = \beta ^{*} (E \otimes \Omega_{X/S}^{1}(D))(A + B -A -2B - A - B
- \widetilde{\mathsf{Fib}}) \\
& = \beta^{*} (E \otimes \Omega_{X/S}^{1}(D-\mathsf{Fib})).
\end{aligned}
$$
Then 
$$
R\beta_{*} W_{-1}E'_{-1/2}\otimes \Omega_{\widetilde{X}/S}^{1}(\widetilde{D})
= E \otimes \Omega_{X/S}^{1}(D-\mathsf{Fib}).
$$

\

\subsubsection{Conclusion for the Dolbeault complexes}
\label{nilcase-modtwo}

Let us now put this back into the situation of our rank $2$ bundle $E$
with no parabolic structure but a Higgs field whose residue is
nilpotent.  We will denote by $W_iE$ the weight filtered bundles given
their reflexive extensions as bundles across the origin.

We get,  for the modified Dolbeault complex $\mathsf{Dol}'$ using
$W_1$ and $W_{-1}$: 
$$
\mathsf{Dol}'_0 = 
\left[
W_1 E \longrightarrow 
W_{-1}E \otimes \Omega_{X/S}^{1}(D) 
\longrightarrow 
W_{-1}(E_D\otimes \ldots )_{P} 
\right]
$$
and
$$
\mathsf{Dol}'_{-1/2} = 
\left[
W_1 E(-\mathsf{Fib}) \longrightarrow 
W_{-1}E \otimes \Omega_{X/S}(D-\mathsf{Fib}) 
\right] .
$$
We could also write the parabolic level $0$ piece as follows.
Introduce notation for the kernel
$$
0\rightarrow W_{-1}^{-P} E \rightarrow W_{-1}E \rightarrow W_{-1}(E_D)_{P}
\rightarrow 0
$$
then 
$$
\mathsf{Dol}'_0 = 
\left[
W_1 E \longrightarrow 
W_{-1}^{-P}E \otimes \Omega_{X/S}^{1}(D) 
\right].
$$
We note that, since the divisor $D$ becomes vertical at the point $P$,
the nonzero residue map of the Higgs field at that point becomes zero
when projected into the relative differentials, so $\varphi$ does
induce a map
$$
W_1 E \stackrel{\varphi}{\longrightarrow}
W_{-1}^{-P}E \otimes \Omega_{X/S}^{1}(D) 
$$
to be used for the above complex.

\

\subsection{Smooth spectral variety}

We now specialize the previous calculations further, to the local
situation where our rank $2$ Higgs bundle\footnote{
Recall that we are working in a local neighborhood---in the
applications to Hecke transformations the spectral variety will have
higher degree, decomposing into a disjoint union of local pieces that
are either \'{e}tale  over $X$ or of the present degree $2$ form.}
is the direct image of a line bundle $\Ll$ on its spectral variety
$\pi_{\Sigma} :\Sigma \rightarrow X$, where $\pi_{\Sigma}$ has degree
$2$ with simple ramification over the divisor $D$.  Let $R\subset
\Sigma$ be the reduced preimage of $D$, so it is a smooth divisor in
$\Sigma$ mapping isomorphically to $D$.

We have 
$$
\Omega ^1_{\Sigma} (\log R) = \pi_{\Sigma}^{*} \Omega^1_X(\log D).
$$
thus we have a inclusion map
$$
\Omega^1_{\Sigma} \hookrightarrow \Omega^1_{\Sigma} (\log R) =
\pi_{\Sigma}^{*} \Omega ^1_X(\log D).
$$
We suppose that the Higgs field $\varphi : E \to E \otimes
\Omega^{1}_{X}(\log D)$ on $E=\pi_{\Sigma *}\Ll$ is given by
\linebreak multiplication $(-)\otimes \alpha : \Ll \to \Ll\otimes
p^{*} \Omega^1_X(\log D)$ by a spectral $1$-form $\alpha$ which is a
section $\alpha \in H^0(\Sigma , \Omega^1_{\Sigma}) \subset
H^{0}(\Sigma, \pi_{\Sigma}^{*} \Omega^1_X(\log D))$.

Choose local coordinates $x,y$ on $X$, such that the map $X\rightarrow S$
is given by $(x,y)\mapsto y$. Choose local coordinates $(x,z)$ on 
$\Sigma$ such that the map $\pi_{\Sigma} : \Sigma \rightarrow X$ is given by 
$$
(x,z) \mapsto (x, y) = (x, x^2 - z^2).
$$
The divisor $R$ is given by $z=0$ and $D$ is given by $y=x^2$. 

Write 
$$
\alpha = a(x,z) dx + b(x,z) dz.
$$
The projection $\alpha_{\Sigma /S}$ of $\alpha$ to a section of 
$\omega _{\Sigma /S}$ is given as follows. The relative differentials
are the reflexive closure of 
$$
\begin{aligned}
\Oo_{\Sigma } \langle dx,dz & \rangle / \langle d(x+z)(x-z) \rangle \ =  \\
& = \Oo _{\Sigma } \langle du/u , dv/v\rangle / \langle du/u + dv/v \rangle 
\end{aligned}
$$
where $u=x+z$ and $v=x-z$. 
Thus we may think of $du/u$ as a generator of $\omega_{\Sigma /S}$
subject to the relation $dv/v = -du/u$. 

We have 
$$
\alpha = \frac{a+b}{2} du + \frac{a-b}{2} dv
$$
and 
$$
\alpha_{\Sigma /S} = 
\left[
u\frac{a+b}{2} + v\frac{a-b}{2}
\right] \frac{du}{u}.
$$
We assume that $a$ and $b$ take generic nonzero values at the origin.
Then the curve $G$ defined by $(\alpha_{\Sigma /S} = 0)$ is smooth at
the origin. Also, under this genericity hypothesis the residue of the
resulting Higgs field is nontrivially nilpotent along $D$.

There are two cases: the parabolic case, as discussed in
section~\ref{parcase}, and the nilpotent case, as discussed in
section~\ref{nilpotentcase}. In the first case, we have a nontrivial
parabolic level $-1/2$, in which case the parabolic structure is given
by the quotient
$$
\xymatrix@M+0.5pc@C+1pc@R-3pc{
  \pi_{\Sigma *}\left(\Ll |_{2R}\right)  \ar[r] & \pi_{\Sigma *}
  \left(\Ll |_R\right) \\
  || & \\
  E |_{D} &
}
$$
over $D$. 
Thus we have $E_{-1/2} = \pi_{\Sigma *} \left(\Ll (-R)\right)$. 

In the second case, we have no parabolic structure on $E$ but the
weight filtration is given by $W_{-1}E = W_0E$ being the kernel of the
above map, and $W_1E/W_0E$ is the quotient $\Ll |_R$. We have $W_{-1}E
= \pi_{\Sigma *} \left(\Ll (-R)\right)$.

\subsubsection{The parabolic case}

We now plug these descriptions in the previous calculations. Consider
first the case with parabolic level $-1/2$. Then we saw in Subsections
\ref{parcase} and \ref{parcase-alt} that
$$
\mathsf{Dol}(E,X/S)_0 = 
\left[
E_0 \rightarrow E_{-1/2} \otimes \omega_{X/S}(D)
\right] .
$$
We have $E_0 = \pi_{\Sigma *}\Ll$ and $E_{-1/2} = \pi_{\Sigma *} (\Ll (-R))$,
whereas 
$$
\pi_{\Sigma}^{*} \omega_{X/S}(D) \cong \omega_{\Sigma /S}(R).
$$
We may therefore write
$$
\mathsf{Dol}(E,X/S)_0 = 
\ \pi_{\Sigma *} 
\left[\hspace{-0.4pc}\text{\begin{minipage}[c]{1.72in} 
$\xymatrix@M+0.5pc@C+1.5pc{\Ll \ar[r]^-{\alpha _{\Sigma / S}} &
\Ll \otimes \omega_{\Sigma/S}}$ \end{minipage}}\right]
$$
In turn, this complex is quasiisomorphic to the second 
bundle restricted over the curve $G$ of zeros of the differential here: 
$$
\mathsf{Dol}(E,X/S)_0 \ \stackrel{q.i.}{\sim}  \
\left(\Ll \otimes
\omega_{\Sigma / S}\right) |_G [-1].
$$
Taking the direct image down to $S$ we obtain $G$ as spectral variety
(see Proposition \ref{relcrit}), 
with spectral line bundle being the restriction of
$\Ll \otimes \omega _{\Sigma / S}$ to $G$. 

\

\

\noindent
{\em Caution:} This is a local calculation near the ramification point 
of the horizontal divisor, and
on branches of the curve $G$ that pass through the ramification point. 
The formula might be different at other points where $G$ passes
through points of $\Sigma$ on other branches lying over $D$.

\

\subsubsection{The nilpotent case}

Next consider the case where there is no parabolic structure but the
residue of the Higgs field is a nontrivial nilpotent transformation.
The weight filtration is given by the quotient.  In this case, there
are two pieces of the Dolbeault complex (modified as in Subsections
\ref{nilcase-modone} and \ref{nilcase-modtwo}) to consider.

First, we had
$$
\mathsf{Dol}(X/S,E)_0 =
\left[
W_1E \rightarrow W_{-1}E\otimes \Omega ^1_{X/S}(D)\otimes \mathcal{I}_{P}
\right]
$$
where the tensoring with $\mathcal{I}_{P}$ that we take the kernel
of the map to the fiber of the quotient sheaf over $P$.

We have $W_1E = \pi_{\Sigma *} \Ll$ and $W_{-1}E = \pi_{\Sigma *}(\Ll (-R))$.
Thus
$$
W_{-1}E\otimes \Omega^1_{X/S}(D)\otimes \mathcal{I}_{P}
=
\pi_{\Sigma *} (\Ll \otimes \omega_{\Sigma/S}\otimes \mathcal{I}_{Q})
$$
where here $\mathcal{I}_{Q} \subset \Oo_{\Sigma}$ is the ideal sheaf of the
point $Q \in \Sigma$ lying over the origin $P \in X$ and
$\Ll \otimes \omega_{\Sigma/S}\otimes \mathcal{I}_{Q}$ is the kernel
$$
0\rightarrow \Ll \otimes \omega _{\Sigma/S} \otimes \mathcal{I}_{Q}
\rightarrow \Ll \otimes \omega _{\Sigma / S} 
\rightarrow \cc _Q\rightarrow 0
$$
of the evaluation at the point $Q\in \Sigma$.
We get
$$
\mathsf{Dol}(X/S,E)_0 =
\pi_{\Sigma *}
\left[\hspace{-0.6pc}\text{\begin{minipage}[c]{2.16in} 
$\xymatrix@M+0.5pc@C+1.5pc{\Ll \ar[r]^-{\alpha _{\Sigma / S}} &
        \Ll \otimes \omega_{\Sigma/S}}\hspace{-0.6pc} \otimes \mathcal{I}_{Q}$
\end{minipage}}
\right].
$$
This is quasi-isomorphic to the second piece restricted over the curve 
$G$: 
$$
\mathsf{Dol}(X/S,E)_0 \
\stackrel{q.i.}{\sim}  \
(\Ll \otimes \omega _{\Sigma / S})|_G(-Q) [-1]
$$
We obtain the following statement: in the situation where the original
Higgs bundle had trivial parabolic structure and nilpotent residue,
the line bundle on $G$ 
yielding the level $0$ part of the parabolic bundle is
$$
F_0=g_{*}(\Ll \otimes \omega _{\Sigma / S}|_G(-Q))
$$
where $g: G\rightarrow S$ is the covering map.

We recall that the level $-1/2$ Dolbeault complex was the same, but
without the $\mathcal{I}_{P}$, then twisted by $-\mathsf{Fib}$ where
$\mathsf{Fib} = f^{-1}(o)\subset X$ was the fiber. When we restrict to
$G$ the fiber becomes $g^{-1}(o)=2Q$ so putting these together gives
$$
\mathsf{Dol}(X/S,E)_{-1/2} \ 
\stackrel{q.i.}{\sim} \
(\Ll \otimes \omega _{\Sigma / S}|_G(-2Q) [-1])
$$
and hence
$$
F_{-1/2} =
g_{*}(\Ll \otimes \omega _{\Sigma / S}|_G(-2Q)).
$$

\

\

\noindent
{\em Caution:} As before, this is only a local calculation near the
ramification point of the horizontal divisor, and on branches of the
curve $G$ that pass through the ramification point.  Again, the
formula might be different at other points where $G$ passes through
points of $\Sigma$ on other branches lying over $D$.

\

\subsubsection{Rephrasing}

Let us rephrase the above calculations in terms of the pullback to
$G\subset \Sigma$ of the relative differentials $\omega_{X/S}$. 
This is useful because, at other points, this is the relevant term. 

We have, as was used above, 
$$
\omega_{\Sigma / S}(R) = \pi_{\Sigma}^{*} (\omega _{X/S}(D)) 
= (\pi_{\Sigma}^{*} \omega _{X/S})(2R).
$$
This gives 
$$
\omega_{\Sigma/S}=(\pi_{\Sigma}^{*} \omega_{X/S})(R).
$$
On the curve $G$, the divisor $R$ is the same as $Q$. 
We conclude the following statement. 

\

\begin{proposition}
\label{localcalc}
Suppose given a Higgs bundle with smooth spectral variety $\Sigma$
having simple ramification over a horizontal divisor $D_H \subset X$
such that $D_H$ is smooth but such that the map $D_{H} \to S$ to the
base $S$ has a simple ramification over the point $o \in S$.  We have
considered two scenarios for producing a parabolic Higgs bundle with a
spectral line bundle $\Ll$. These result in a higher direct image
Higgs bundle $F_{\bullet}$ whose local expression near $o \in S$ is as
follows.

\smallskip

\noindent
\quad \emph{\bfseries (parabolic case)} \
In the case of parabolic levels $0,-1/2$ along $D_H$ for the Higgs
bundle on $X$, we get
$$
F_0 = g_{*} \left(\left(\Ll |_{G} \otimes j^{*} \omega _{X/S}\right)(Q)\right))
$$
where $j: G\rightarrow X$ denotes the composed map. Recall here that
$\Ll$ being the spectral line bundle means that its direct image to
$X$ is the level $0$ piece of the parabolic structure.

\smallskip

\noindent
\quad \emph{\bfseries (nilpotent case)} \
In the case where the parabolic levels along $D_H$ are trivial
but the residue of the Higgs field along $D_H$ is nilpotent, we get
$$
F_0 = g_{*} \left(\Ll |_{G} \otimes j^{*} \omega _{X/S}\right)
$$
and 
$$
F_{-1/2} = g_{*} \left(\left(\Ll |_{G} \otimes j^{*}
\omega_{X/S}\right)(-Q)\right).
$$
\end{proposition}

\

\

\begin{remark}
\label{fillingin}
The result of the proposition is local near a singular point. The
statement refers to the way of filling in the structure of the
spectral line bundle on the spectral covering $G$ over $S$, noting
that over a general point of $S$ we have a canonical identification
between the spectral line bundle on $G$, and the restriction of $\Ll$
to $G$ tensored with $\omega_{X/S}$.
\end{remark}

\

\

\subsection{Globalization}

Let us now look at a global situation. The notations and hypotheses of
Subsection \ref{pushforward-statements} are heretofore in
effect. Thus, $f:X\rightarrow S$ is a map from a smooth projective
surface to a smooth projective curve, and $D \subset X$ is a simple
normal crossings divisor, $T \subset S$ is a divisor consisting of a
finite set of points $t_1,\ldots , t_k$, and let $K\subset X$ be a
closed subset containing the ``other'' points of type \ref{type6}.

Write $D = D_H + D'_V$ such that the irreducible components of 
$D_H$ dominate $S$ and the irreducible components of $D'_V$ map
into $T$. We let $D_V$ denote the full inverse image of $T$. 
Assume that $f(K)$ is a finite subset of $S$.

\

\noindent
Suppose that for any point $x\in X$, one of the following holds: 
\begin{enumerate}
\item[(1)] $x\in X -D$ and $f$ is either smooth at $x$ 
(type \ref{type1})or has
a simple normal crossing (type \ref{type2}); 
\item[(2)] $x\in D_H$, $f$ is smooth at $x$ and $D_H$ is étale over
$S$ at $x$ (type \ref{type3}); 
\item[(3)] $x\in D_H$ , $D_H$ is smooth at $x$, 
$f$ is smooth at $x$, $f(x)=t_i\in T$, and
$f|_{D_H}$ has a simple ramification point at $x$ (type \ref{type4}); 
\item[(4)] or $x\in K$ (type \ref{type6}). 
\end{enumerate}

\

\

\noindent
For the moment, this supposes that there are no points of type
\ref{type5}; those will be treated in the next subsection.
Notice that $D'_V\subset K$.

\

Let $(E_{\bullet},\varphi)$ be a parabolic logarithmic Higgs bundle
over $(X,D_H)$. We assume that over $D_H$, either the parabolic
structure is trivial and $\varphi$ has nilpotent residue, or that the
parabolic structure has jumps at  $0,-1/2$ and $\varphi$ is strictly
parabolic.

Let $\pi_{\Sigma} : \Sigma \rightarrow X$
be the spectral covering of $(E_{0},\varphi )$
with spectral line bundle $\Ll$ over $\Sigma$
so that $E_0 = p_{*}(\Ll )$, and spectral 
$1$-form $\alpha$ inducing the Higgs field $\varphi$. 

We assume that away from the subset $K$, $\Sigma$ is smooth and the
covering map $\pi_{\Sigma}$ has at most simple ramification; let
$R\subset \Sigma$ be the reduced inverse image of $D_H$, assume $R$ is
smooth and $\pi_{\Sigma}$ has simple ramification along $R$ away from
points of $K$.  Let $\alpha _{\Sigma/S}$ be the relative spectral form
viewed as a section of
$\left(\pi_{\Sigma}^{*}\omega_{X/S}\right)(R)$. Let $G$ be the curve
of zeros of this form, defined away from $K$.

We work under the following genericity hypotheses:
\begin{itemize}
\item We assume that the closure in $X$ of the image of $G$
  does not meet $K$. In particular $G$ is proper.
\item We also assume that the restriction of the spectral $1$-form
  $\alpha$ to the vertical direction in $\Sigma$ over a point of type
  \ref{type4} is nonzero.
\end{itemize}

\

\

\noindent
Let $g: G\rightarrow S$ be the map, factoring as $g=f\circ j$ through
the morphism $j:G\rightarrow X$. Let $\Ll |_G$ be the restriction of
$\Ll$ to $G$.  Let $Q:= R\cap G$ be the trace of the reduced divisor
$R$ onto the curve $G$.

\

\noindent
With this notation we have the following

\

\begin{lemma}
\label{alpha-nonzero-vertical}
The hypothesis that $\alpha$ is nonzero in the vertical direction, at
points in $\Sigma$ lying over points of type \ref{type4},
implies that $G$ is transverse to $R$ at points of type \ref{type4},
so $Q$ is reduced at such points.
\end{lemma}
\begin{proof}
Let us calculate in coordinates $x,t$ on $X$, such that $t$ gives the
map $f : X\rightarrow S$.  Assume the horizontal divisor is $x^2 - t= 0$
and the covering $\Sigma$ is $w^2 = x^2-t$. Thus, $x$ and $w$ give a
system of coordinates on $\Sigma$. Write $\alpha = adx + bdw$ with
$a=a(x,w)$ and $b=b(x,w)$ holomorphic functions of $x,w$. The equation
$t=x^2 - w^2$ tells us that dividing out by $dt$ is equivalent to
setting $xdx = wdw$. The form $dx$ provides a frame for the sheaf of
relative differentials $\Omega ^1_{X/S}$ and in terms of this frame,
$$
\alpha_{\Sigma/S} = (a + bx/w) dx. 
$$
Thus, the equation for the upper critical locus $\alpha_{\Sigma/S}
=0$ becomes
$$
aw+bx=0.
$$
The hypothesis $b(0,0)\neq 0$ implies that the linear term of this
equation at the origin is nonzero. The ramification divisor $R$ is
given by $w=0$, that is to say it is the $x$-axis in this coordinate
system, and $b(0,0)\neq 0$ tells us that the above equation has a
simple zero along the $x$-axis, so the zero set $G$ is transverse to
$R$. This completes the proof of the lemma.
\end{proof}

\

\

\noindent
The following statement gives Theorem \ref{pushforward-thm} in the
case when there are no points of type \ref{type5}.

\

\begin{theorem}
\label{pushforward-appli-main}
With all the above hypotheses, in the case of parabolic levels
$0,-1/2$ we have
$$
F_0 = g_{*} \left( 
\Ll |_G \otimes j^{*}\omega_{X/S} \otimes \Oo _G(Q)
\right).
$$
In the nilpotent case, we have
$$
F_0 = g_{*} \left( 
\Ll |_G \otimes j^{*}\omega _{X/S} 
\right) 
$$
and 
$$
F_{-1/2} = g_{*} \left( 
\Ll |_G \otimes j^{*}\omega _{X/S} \otimes \Oo _G(-Q)
\right).
$$
\end{theorem}
\begin{proof}
In the nilpotent case, $G$ does not intersect $R$ at points distinct
from ramification points, because we know that the map $gr_1^W
\rightarrow gr_{-1}^W$ is an isomorphism at any point of $D_H$ \'{e}tale
over the base. We notice that, at points of $G$ mapping to points of
$D_H$ but in sheets of $\Sigma$ that are \'{e}tale, the formula is as
stated since the term of the Dolbeault complex has $W_{-1}E = E(-D_H)$
and $\Omega ^1_{X/S}(\log D_H) = \omega _{X/S}(D_H)$ locally on
$\Sigma$ at those points.

Also, the divisor consisting of ramification points
is the same as $Q=R\cap G$ in this case, in view
of the non-intersection of $G$ with other points of $R$. 
Thus, for the nilpotent case our previous
calculations give the required results. 

For the parabolic case, there might be points where $G$ intersects $R$
at points where $D_H$ is \'{e}tale over $S$. In this case, the
component of  $E_{-1/2}$ coming from the local neighborhood in $\Sigma$
is $\pi_{\Sigma *} (\Ll (-R))$ and 
$$
\pi_{\Sigma}^{*} \Omega^1_{X/S}(\log D_H)
= (\pi_{\Sigma}^{*} \omega_{X/S})(2R)
$$
so the required bundle in the degree $1$ term of the Dolbeault complex 
is
$$
\pi_{\Sigma *} \left(\Ll \otimes \pi_{\Sigma}^{*} \omega_{X/S}(R)\right).
$$
When we restrict to $G$, we get the line bundle 
$$
\Ll |_G \otimes j^{*}\omega _{X/S} \otimes \Oo _G(Q)
$$
at these points. 
At the ramification points of $D_H$, we notice that
$Q$ is the same as the divisor of ramification points considered above. 
Therefore, near the ramification points the required bundle is also
$$
\Ll |_G \otimes j^{*}\omega _{X/S} \otimes \Oo _G(Q).
$$
This yields the stated formula. 
\end{proof}

\

\

\subsection{Type \texorpdfstring{\ref{type5}}{type5} points}
\label{t5points}

The calculation needs to be extended to cover the points of type
\ref{type5} in the classification of Subsection
\ref{pushforward-statements}.  Suppose $x\in X$ is a point of type
\ref{type5}. Recall that this means that the horizontal divisor $D$
has a node at $x$, with both branches etale over $S$, and that over a
small neighborhood of $x \in X$ the spectral variety $\pi_{\Sigma} :
\Sigma \rightarrow X$ decomposes into a disjoint union of two sheeted
covers, each having an ordinary double point over $x$.

Choose such a neighborhood $x\in U\subset X$, and let ${}^{U}\Sigma
\subset \Sigma$ be the corresponding neighborhood of one of the double
points over $x$.  Consider the blow-up map $\blowb :
\widetilde{U}\rightarrow U\subset X$. Let $Z$ be the normalization of
${}^{U}\widetilde{\Sigma}:= {}^{U}\Sigma \times _U
\widetilde{U}$. Notice that $Z$ is the blow-up of ${}^{U}\Sigma$ at
the double point.

Let $D^+$ and $D^-$ be the two branches of the strict transform of $D$
in $\widetilde{U}$, and let $\mathsfit{F}$ be the strict transform of
the fiber of $f : U \to S$ over $o \in S$. Let $D_Z^+$ and $D_Z^-$
denote the reduced inverse images of the divisors $D^+$ and $D^-$ in
$Z$.  We will write $\mathsfit{B}\subset \widetilde{U}$ be the
exceptional divisor of $\blowb$ and $\mathsfit{B}_{Z}\subset Z$ for
its reduced inverse image in $Z$. Thus $\mathsfit{B}\cong \pp^1$ and
$\mathsfit{B}_{Z} \cong \pp^1$ are copies of the projective line with
the map $\mathsfit{B}_{Z}\rightarrow \mathsfit{B}$ a double cover
ramified over the points $d^+ := D^+\cap \mathsfit{B}$ and $d^- :=
D^-\cap \mathsfit{B}$.

Let $\Ll$ be the spectral line bundle on $\Sigma$, and let $\Ll_Z$
denote its pull-back to $Z$.  Write $\pi_{{}^{U}\Sigma} : {}^{U}\Sigma
\to U$ for the restriction of the spectral covering map $\pi_{\Sigma}
: \Sigma \to X$ to ${}^{U}\Sigma$ and let
\[
(E,\varphi) =
\left(\pi_{{}^{U}\Sigma *}\left(\Ll |_{{}^{U}\Sigma}\right),
\pi_{{}^{U}\Sigma *}\left( \alpha |_{{}^{U}\Sigma}\right)\right)
\]
be the corresponding Higgs bundle on $U\subset X$, with
$(\widetilde{E},\widetilde{\varphi}) :=
\left(\blowb^{*}E,\blowb^{*}\varphi\right)$ its inverse image on
$\widetilde{U}$. Since we are working with neighborhoods both
downstairs and upstairs, $E$ has rank $2$ here.

The assumption on the nilpotent residue of the Higgs field $\varphi$
means that there is a decomposition $E=E_1 \oplus E_{-1}$ into a
direct sum of two line bundles, such that the image of the residue of
$\varphi$ is $E_{-1}$ over both branches of the divisor $D$. The
pullback then decomposes as $\widetilde{E}=\widetilde{E}_1 \oplus
\widetilde{E}_{-1}$

Let $\mathsfit{p} : Z\rightarrow \widetilde{U}$ denote the natural
map. The bundle $\widetilde{E}$ is obtained from
$\mathsfit{p}_{*}(\Ll_Z(-\mathsfit{B}_{Z}))$ by glueing together the
line bundles on the two branches of $Z$ over general points of
$\mathsfit{B}$. Since $\Ll_Z$ is pulled back from a line bundle on $\Sigma$, the
restriction $\Ll_Z|_{\mathsfit{B}}$ is a trivial line bundle that we
will denote $\Ll _{\mathsfit{B}}$. The glueing is done using this
trivialization, so we have an exact sequence
$$
0\rightarrow \widetilde{E} \rightarrow \mathsfit{p}_{*} (\Ll _Z )
\rightarrow \mathsf{coker} \rightarrow 0
$$
where the cokernel sheaf $\mathsf{coker}$ is supported on
$\mathsfit{B}$ and is also the cokernel in the sequence
$$
0\rightarrow \Ll _{\mathsfit{B}}  \rightarrow \mathsfit{p}_{*} (\Ll _Z ) \rightarrow \mathsf{coker} \rightarrow 0
$$ 
over $B$. 
The exact sequence of the elementary transformation relating
$\mathsfit{p}_{*}(\Ll_Z(-\mathsfit{B}_{Z}))$ and $\widetilde{E}$  is 
$$
0\rightarrow \mathsfit{p}_{*} (\Ll _Z (-\mathsfit{B}_{Z}))
\rightarrow \widetilde{E} \rightarrow \Ll _{\mathsfit{B}} \rightarrow
0.
$$
Notice now that $\Ll _\mathsfit{B} \cong \widetilde{E}_{1}
|_{\mathsfit{B}} = (E_{1})_{x}\otimes \Oo_{\mathsfit{B}}$ viewed as quotient
  of $\widetilde{E}|_B$, as it is the trivial subbundle whose values
  over the two ramification points correspond to the unramified part.

On the subsheaf given by the weight filtration we get
$$
0\rightarrow \mathsfit{p}_{*} (\Ll _Z (-\mathsfit{B}_{Z}-D_Z^+
-D_Z^-)) \rightarrow W_{-1}\widetilde{E} \rightarrow \widetilde{E}_{1}
|_{\mathsfit{B}}( -d^+ -d^-)\rightarrow 0.
$$
The sheaf of relative logarithmic differentials is then
$$
\omega_{\widetilde{U}/S}(\log ) = (\mathsfit{b}^{*} \omega _{X/S})(D^+
+ D^- + \mathsfit{B}).
$$
Pulling back to $Z$ gives
$$
\mathsfit{p}^{*} \omega _{\widetilde{U}/S}(\log ) = (\mathsfit{b}
\circ \mathsfit{p})^{*} \omega _{X/S}
\otimes \Oo _Z (2D_Z^+ + 2 D_Z^-  +  \mathsfit{B}_{Z}).
$$
Also, $\omega _{\widetilde{U}/S}(\log )|_{\mathsfit{B}} \cong
\Oo_{\mathsfit{B}}(1)$.  We have
$$
\begin{aligned}
\mathsfit{p}_{*} (\Ll_Z & (-\mathsfit{B}_{Z}-D_Z^+ -D_Z^-))
\otimes \omega_{\widetilde{U}/S}(\log ) 
= \\
& = 
\mathsfit{p}_{*} (\Ll_Z (-\mathsfit{B}_{Z}
-D_Z^+ -D_Z^- + 2D_Z^+ + 2 D_Z^-  +  \mathsfit{B}_{Z}))
\otimes \mathsfit{b}^{*} \omega _{X/S}  \\
& =
\mathsfit{p}_{*} (\Ll_Z (D_Z^+ + D_Z^- )) \otimes
\mathsfit{b}^{*} \omega _{X/S}.
\end{aligned}
$$
Thus we have an exact sequence
$$
0\rightarrow  \mathsfit{p}_{*} (\Ll _Z (D_Z^+ + D_Z^- )) \otimes
\mathsfit{b}^{*} \omega _{X/S} \rightarrow 
W_{-1}\widetilde{E}\otimes  \omega _{\widetilde{U}/S}(\log )\rightarrow 
\widetilde{E}_{1} |_{\mathsfit{B}}( -d^+ -d^-)\otimes \Oo _{\mathsfit{B}}(1)
\rightarrow 0.
$$

\

\begin{lemma}
\label{t5subsheaf}
Suppose $G\subset \Sigma$ is a curve that decomposes, near the double
point, into a disjoint union of two smooth branches whose tangent
vectors at the double point are distinct.  Then the map on spaces of
sections on the neighborhood
$$
\Gamma (U, \mathsfit{p}_{*}
(\Ll _Z (D_Z^+ + D_Z^- )) \otimes \mathsfit{b}^{*} \omega _{X/S} )
\rightarrow  (G\cap {}^U \Sigma , \Ll |_G \otimes \omega _{X/S}|_G)
$$
is an isomorphism. 
\end{lemma}
\begin{proof}
The curve $G$ is isomorphic to its strict transform inside $Z$. The
two branches intersect $\mathsfit{B}_{Z}$ in different points. The
line bundle $\Ll _Z (D_Z^+ + D_Z^- )$ is $\Oo_{\mathsfit{B}_{Z}}(2)$
on $\mathsfit{B}_{Z}$, so its sections (which extend to local sections
in a neighborhood) span the space of sections over $G$.
\end{proof}

\

\

\begin{theorem}
\label{pushforward-appli-t5}
In the presence of points of type \ref{type5}, let $F'$ be the
parabolic bundle defined by
$$
F'_0 = g_{*} \left( 
\Ll |_G \otimes j^{*}\omega _{X/S} 
\right) 
$$
and 
$$
F'_{-1/2} = g_{*} \left( 
\Ll |_G \otimes j^{*}\omega _{X/S} \otimes \Oo _G(-Q)
\right) ,
$$ and let $F$ be the parabolic Higgs bundle corresponding to the
higher direct image local system. Then there is an injective morphism
$F'\hookrightarrow F$.
\end{theorem}
\begin{proof}
The expression of the higher direct image Higgs bundle comes from a
direct image of a sheaf supported on the relative critical locus, so
it is local over the relative critical locus. Near a point of type
\ref{type5}, it comes from blowing up once $\mathsfit{b} :
\widetilde{X}\rightarrow X$. There is a sheaf $\Gg$ supported on
$\widetilde{X}$ whose direct image down to $S$ is the local piece of
$F$, and $\Gg$ is supported on the image of $G$ union the exceptional
divisor $\mathsfit{B}$. In turn, the expression is local near each
exceptional curve in the blowup $\widetilde{\Sigma}$ of $\Sigma$ at
its ordinary double points, and the sheaf $\Gg$ is a quotient of
$W_{-1}\widetilde{E}\otimes \omega _{\widetilde{X}/S}(\log )$ (where
here $E$ is the piece of the full Higgs bundle, corresponding to the
local piece of $\Sigma$).  Thus, any sections of
$W_{-1}\widetilde{E}\otimes \omega_{\widetilde{X}/S}(\log )$ over a
neighborhood of $\mathsfit{B}$ will generate a subsheaf of $F$.

The exact sequences above give a map from sections of
$\mathsfit{p}_{*} (\Ll _Z (D_Z^+ + D_Z^- )) \otimes \mathsfit{b}^{*}
\omega _{X/S}$ to sections of $W_{-1}\widetilde{E}\otimes \omega
_{\widetilde{X}/S}(\log )$, and Lemma \ref{t5subsheaf} says that these
restrict on $G$ to sections of $\Ll |_G \otimes
\omega_{X/S}|_G)$. Those are the sections that appear in the
definitions of $F'_0$ and $F'_{-1/2}$ (note that there is no
difference in the two parabolic level spaces locally at a point of
type \ref{type5}).

Let $F'$ be the subsheaf of $F$ generated by such sections near the
points of type \ref{type5}, and equal to $F$ as calculated in Theorem
\ref{pushforward-appli-main} elsewhere. This is the stated subsheaf.
\end{proof}

\

\

\begin{remark}
We  apply this statement by calculating the degree of the
parabolic sheaf $F'$ given by Theorem \ref{pushforward-appli-t5}. If
it has parabolic degree $0$, then since we know that $F$ also has
parabolic degree $0$, the map $F'\hookrightarrow F$ is an isomorphism,
and this will yield the same calculation as in Theorem
\ref{pushforward-appli-main} for the case when there are points of
type \ref{type5}.
\end{remark}

\

\section{Drinfeld's construction}
\label{chapter-drinfeld}

Drinfeld's original construction of Hecke eigensheaves was done in
\cite{Drinfeld} and later consolidated and extended by Laumon
\cite{Laumon95}, Gaitsgory \cite{GaitsgoryThesis,Gaitsgory} and others.
In order to make a comparison with our spectral cover constructions,
we will give here a preliminary approach to the interpretation of
Drinfeld's construction in the setting of Higgs bundles.

A main ingredient of Drinfeld's construction is the following general
definition: given a local system $\Lambda$ on $C$, we obtain a local
system $\Lambda ^{( \otimes m)}$ on the symmetric power
$\op{Sym}^m(C)$ with singularities along the big diagonal.

In the case ${\rm rk}(\Lambda )=2$, Drinfeld starts with $\Lambda ^{(
  \otimes m)}$ as input and uses Radon transform to construct a local
system on an open subset of a projective space bundle over
$\mathbf{Bun}$. The main theorem of \cite{Drinfeld} says that it
descends from the projective space bundle down to a perverse sheaf on
the moduli stack $\mathbf{Bun}$. Going to the coarse moduli spaces of
stable bundles we get a local system on an open subset of the moduli
space of stable bundles. We would like to calculate the Higgs sheaf
associated to the Radon transform and use that to show that the Higgs
sheaves we constructed above are the same as the ones corresponding to
Drinfeld's perverse sheaf.

Fix a line bundle $M$ on $C$ of sufficiently high degree $m$. Denote by
$$
P:= |M| = \pp H^0(C,M)
$$
the linear system, i.e. the space of divisors $x=x_1+\ldots + x_m$
on $C$ such that $\Oo _C(x_1+\ldots + x_m)\cong M$. It is the fiber
over the point $M$ of the map on the right here:
$$
P\stackrel{\bii}{\hookrightarrow} \op{Sym}^m(C) \rightarrow \op{Pic}^m(C)
$$
where we took the opportunity to give the name $\bii$ to the inclusion. 

Serre duality says that 
$$
H^0(C,M)^{\vee} \cong H^1(C,M^{\vee}\otimes \omega _C).
$$
Set 
$$
Q:=  \pp H^1(C,M^{\vee}\otimes \omega _C),
$$
so $Q$ is naturally the dual projective space to $P$. In other words,
$Q$ is identified with the space of hyperplanes in $P$ and vice-versa.

In a certain sense dual to $\bii$ is a map $\bLL : Q \rightarrow
\mathbf{Bun}$ to the moduli stack, or with the same notation a
rational map to the moduli space of stable bundles. This map
associates to $\xi \in H^1(C,M^{\vee}\otimes \omega_C)$, considered
as an extension class i.e.  an element of $\op{Ext}^1(M,\omega_C)$, the
bundle in the middle of the corresponding extension
\[
\tag{$\xi$}
0\rightarrow \omega _C \rightarrow E \rightarrow M \rightarrow 0.
\]
It is a point in the moduli stack of rank $2$ bundles having
determinant $M\otimes \omega _C$.  The appearance of $\omega_C$ on
the left of the extensions we look at is just there to accommodate the
twist in Serre duality.

Let $X_{M\otimes \omega _C}$ denote the coarse moduli space of
semistable rank $2$ bundles up to $S$-equivalence, having determinant
$M\otimes \omega _C$.  Thus $\bLL$ may be interpreted (losing some
information) as a rational map $\bLL :Q \dashrightarrow X_{M\otimes
  \omega _C}$.

The duality between $P$ and $Q$ is reflected in the incidence correspondence
$$
I\subset P \times Q, \;\;\;\; I:= \left\{ (x,\xi ) \mbox{ s.t. } x\in
\xi \mbox{ i.e. } \langle \xi , x \rangle = 0 \, \right\} .
$$
Let $p:I\rightarrow P$ and $q:I\rightarrow Q$ be the
projections. Thus, $p$ induces an isomorphism between $q^{-1}(\xi )$
and the hyperplane in $P$ associated to $\xi $, whereas $q$ induces an
isomorphism between $p^{-1}(x)$ and the hyperplane in $Q$ associated
to $x$.

Let $\bii :P \hookrightarrow \op{Sym}^m(C)$ be the
inclusion. Drinfeld's basic \emph{\bfseries Radon transform
construction} starts with a rank $2$ local system $\Lambda$ on $C$ and
assigns to it the perverse sheaf
$$
\mathbf{Rad} := R q_{*} \left( p^{*} \bii^{*} \Lambda ^{( \otimes m)}\right) 
$$
on $Q$. 

Let $\Delta^{(m)}_P\subset P$ be the intersection of the big diagonal
in $\op{Sym}^m(C)$ with $P$.  This is the singular set of $\bii^{*}
\Lambda ^{( \otimes m)}$. Thus, the singular set of $p^{*} i^{*}
\Lambda ^{( \otimes m)}$ on $I$ is
$$
D_H:= p^{-1}(\Delta ^{(m)} ).
$$
Let $U\subset Q$ be the open set over which $D_H$ is a relative normal
crossings divisor, that is to say where all strata of $D_{H}$ are etale
under the projection $q$. Then $\mathbf{Rad}_U:= \mathbf{Rad} |_U$ is a
local system over $U$.

The first main part of the geometric Langlands correspondence for rank
$2$ bundles, according to Drinfeld and Laumon, may be formulated as
follows:

\begin{theorem}[Drinfeld-Laumon \cite{Drinfeld,Laumon95}] 
\label{drinfeldthm}
There is a Zariski open subset $U'\subset U$
such that the map $\bLL':U'\rightarrow X_{M\otimes \omega _C}$ 
is well-defined, and maps into the complement of
the wobbly locus $\Wob _{M\otimes \omega _C} \subset X_{M\otimes \omega _C}$. 
The restriction of $\mathbf{Rad} _U$ to $U'$ denoted 
$\mathbf{Rad}_{U'}$ is a local system on $U'$, constant on the fibers
of $\bLL'$. 
It is therefore isomorphic to the pullback of a local system $\mathbf{E}$
on $X_{M\otimes \omega _C} - \Wob _{M\otimes \omega _C}$. 
The local system $\mathbf{E}$ is a Hecke eigensheaf with elgenvalue $\Lambda$
(in an appropriate sense
taking into account the determinant $M\otimes \omega _C$). 
\end{theorem}

\

\noindent
Another proof was given by Gaitsgory in his thesis \cite{GaitsgoryThesis}. 

\

\

\subsection{A spectral variety}

The perverse sheaf $\mathbf{Rad}$ corresponds to a $\srD$-module, and
by Sabbah's and Mochizuki's theory
\cite{Mochizuki-D1,Mochizuki-D2,Sabbah}, it has a structure of purely
imaginary pure twistor $\srD$-module. Over the open subset $U'\subset
Q$ the fiber at $z = 0$ in the twistor line, is a Higgs
bundle. Let
$$
\Sigma_{\mathbf{Rad},U'}\subset T^{\vee} U'
$$
be the spectral variety of this Higgs bundle.  We would like to
express this in terms of the fiber of the Hitchin fibration.

Consider the moduli space $X_{M\otimes \omega _C}$ of semistable rank
$2$ bundles up to $S$-equivalence, with determinant $M\otimes \omega
_C$, with the smooth open subset $X^{\circ}_{M\otimes \omega _C}$ of
stable bundles. Let $\Higgs_{M\otimes \omega _C}$ denote the
moduli space of semistable Higgs bundles with determinant $M\otimes
\omega _C$, and let $\hit : \Higgs_{M\otimes \omega _C}
\rightarrow \aaaa^N$ be its Hitchin fibration.
A general point $\mathsf{b}\in \aaaa^N$ corresponds to a spectral
curve $\Ctilde \stackrel{\pi}{\rightarrow} C$ provided with a
tautological differential $\lambda \in H^{0}(\Ctilde,\pi^{*}\omega_{C}) \subset
H^0(\Ctilde, \omega_{\Ctilde})$.

Let $(E,\theta)$ be the Higgs bundle associated to $\Lambda$. Let
$\Ctilde$ be the spectral covering of $(E,\theta)$, $\lambda \in
H^0(\Ctilde, \omega_{\Ctilde})$ be the spectral $1$-form, and let
$\specN$ be the corresponding spectral line bundle on $\Ctilde$.
Denote by $(E^{(m)}, \theta^{(m)})$ the Higgs bundle associated to
$\Lambda ^{( \otimes m)}$, that we view as having a parabolic
structure in codimension $\leq 1$.

The spectral covering of $(E^{(m)}, \theta ^{(m)})$ is described as is
$\op{Sym}^m(\Ctilde )$, see Lemma \ref{spectralEm} below.  If $m$ is
big enough, this is a projective space bundle over $\op{Pic}^0(\Ctilde
)$.

We assume that the point $\mathsf{b} \in \aaaa^N$ corresponds to
$(E,\theta )$ and hence to $\Lambda$, in that the spectral covering
of $(E,\theta )$ is $\Ctilde$ and spectral $1$-form of $(E,\theta )$
is $\lambda$.

Let $\Prym _{M\otimes \omega_C}$ denote the Prym variety of line
bundles $V$ on $\Ctilde$ such that $\pi_{*}(V)$ has determinant
$M\otimes \omega _C$. Equivalently it means that the norm of the
divisor defining $V$ down to $C$ is the divisor of $M\otimes \omega
_C^{\otimes 2}$, in particular the line bundles $V$ have degree
$m+4g-4$.

The tautological one-form $\alpha_{\Prym}$ on $\Prym _{M\otimes \omega_C}$
leads to a rational map $\Prym_{M\otimes \omega_C} \dashrightarrow 
T^{\vee}X^{\circ}_{M\otimes \omega _C}$,
and this may be pulled back using the dominant
rational map $Q\dashrightarrow X^{\circ}_{M\otimes \omega _C}$ to get a map 
$$
\Prym _{M\otimes \omega_C} \times _{X_{M\otimes \omega _C}} Q
\dashrightarrow T^{*} Q.
$$

\

\noindent
With this setup we now have the following

\

\

\begin{theorem}
\label{rad}
The restriction $\Sigma _{\mathbf{Rad},U'}$ over $U'\subset Q$ is
isomorphic, as a variety mapping to $T^{\vee} U'$, with the pullback
$\Prym _{M\otimes \omega_C} \times _{X_{M\otimes \omega _C}} U'$.
\end{theorem}

\

\begin{corollary}
\label{comparison}
In the case when $C$ is a curve of genus $g=2$, the local system
constructed by Drinfeld-Laumon-Gaitsgory on $X_{M\otimes \omega _C} -
\Wob _{M\otimes \omega _C}$ is the same as the local system we have
constructed in the previous chapters of this paper.
\end{corollary}

\

\

\noindent
The next subsections are devoted to the proofs. Some parts work for
rank $2$ local systems and bundles on a curve $C$ of arbitrary genus
$g$. For the calculations of parabolic structures along the wobbly
divisor we will restrict to the case of curves of genus $g=2$ where we
understand well the geometry.

\

\
 
\subsection{The incidence correspondence}

We have fixed a line line bundle $M$ of degree $m$ on $C$, and defined 
$$
P:= \pp H^0(M), \;\;\; Q:= \pp H^1(M^{\vee} \otimes \omega _C).
$$
These projective spaces are dual by Serre duality. We assume $m\gg 0$,
so they both have dimension $m+1-g$. Let 
$$
(p,q) : I \hookrightarrow P\times Q
$$
be the incidence correspondence.  

There is a rational map $\bLL:Q \dashrightarrow X_{M\otimes \omega _C}$
from $Q$ to the moduli space $X_{M\otimes \omega _C}$ of rank $2$
vector bundles $\Ll$ with determinant $M\otimes\omega _C$, sending
a point represented by a nonzero class
$\xi \in H^1(M^{*}\otimes \omega _C)= \op{Ext}^1(M,\omega _C)$
to the isomorphism class of the bundle in the extension $\xi$
\[
\tag{$\xi$}
0 \rightarrow \omega _C \rightarrow \Ll \rightarrow M \rightarrow 0.
\]
On the other hand there is a map $\bii : P\rightarrow
\op{Sym}^m(C)$ sending a point represented by a nonzero section $f\in
H^0(C,M)$ to the divisor $z=z_1+\ldots + z_m$ of zeros of $f$. The
image is the set of points of $\op{Sym}^m(C)$ mapping to $[M]\in
Pic^m(C)$.

Let $(E,\theta )$ be a rank $2$ Higgs bundle corresponding to a local
system $\Lambda$.  We recall that there is a perverse sheaf denoted
$\Lambda^{(m)}$ on $Sym^m(C)$ obtained by descending
$\Lambda^{\boxtimes m}$ on $C^m$ via the action of the symmetric
group.

This perverse sheaf corresponds to a parabolic Higgs bundle
on$\op{Sym}^{m}(C)$ that we will denote by $(E^{(m)}, \theta ^{(m)})$.
Since the divisor in $\op{Sym}^m(C)$ has non-normal crossings
singularities, we do not exactly know what a parabolic structure
means, so we will instead say that we work with the pure twistor
$\srD$-module and look at the Higgs fiber. This has a parabolic structure
given by the $V$-filtration in codimension $\leq 1$, and we get a
parabolic Higgs bundle on an open subset $\op{Sym}^m(C)^{\circ}
\subset \op{Sym}^{m}(C)$ complement of a set of codimension $2$.

Let $\Ctilde \stackrel{\pi}{\rightarrow} C$ be the spectral covering
of $(E,\theta )$. We assume that $\Ctilde $ is smooth, and since the
degree is $2$ it automatically has simple ramification.  Let $\lambda
\in H^0(\Ctilde , \pi ^{*} \omega _C)$ be the tautological $1$-form,
and let $\specN$ be the spectral line bundle. Thus $E\cong \pi _{*}
\specN$, and $\theta = \pi_{*}(\lambda\otimes (-))$.

We have a covering 
$$
\op{Sym}^m(\Ctilde ) \rightarrow \op{Sym}^m(C)
$$ and a tautological $1$-form $\lambda^{(m)}$ given as the descent
from $\Ctilde ^m$ of the sum of the pullbacks of $\lambda$ from the
components. The covering has a line bundle $\specN^{(m)}$ descended from
the line bundle $\specN^{\boxtimes m}$ on $\Ctilde ^m$; it is characterized
by the condition that $\specN^{\boxtimes m}$ is the pullback of $\specN^{(m)}$.

Let $\ExY \subset \op{Sym}^m(\Ctilde )$ be the divisor consisting of
points $\tilde{t}_{1} + \cdots + \tilde{t}_{m}$ such that $t_i=t_j$
but $\tilde{t}_i \neq \tilde{t}_j$, where $t_i$ are the images in $C$
of $\tilde{t}_i$.

\

\

\begin{lemma}
\label{spectralEm}
The spectral data in codimension $\leq 1$ for the parabolic Higgs
bundle $(E^{(m)}, \theta ^{(m)})$ is given by the covering
$\op{Sym}^m(\Ctilde ) $ with its tautological $1$-form $\lambda^{(m)}$
and spectral line bundle $\specN^{(m)}$. The parabolic structure is
the standard one with weights $0,1/2$ using the divisor $\ExY$.

The inverse image of $P$ in $\op{Sym}^m(\Ctilde )$ is the subvariety
$\op{Sym}^m(\Ctilde )_M$ consisting of divisors whose associated line
bundle has norm down to $C$ equal to $M$. If $m$ is big enough, it is
smooth, being a projective space bundle over a translate of the Prym
variety that is the kernel of $\op{Pic}^m(\Ctilde ) \rightarrow
\op{Pic}^m(C)$.
\end{lemma}
\begin{proof}
The exterior tensor product Higgs bundle $E^{\boxtimes m}$ over $C^m$
has spectral variety $\Ctilde ^m$. On there, the spectral one-form for
the tensor product is the sum of the pullbacks of $\lambda$ to each of
the components, and the spectral line bundle on $\Ctilde ^m$ is
$\specN^{\boxtimes m}$.  These descend to the given spectral data outside
of codimension $2$. Notice here that we might also need to remove a
codimension $2$ subset bigger than just the singular locus of the
parabolic divisor, in case the image of $\op{Sym}^m(\Ctilde )$ in the
logarithmic cotangent bundle has singularities at the branch points of
$\Ctilde / C$ along the multidiagonal (we did not calculate if this
happens or not).

In order to understand the line bundle, we note the following
remark. Let $g: C^m \rightarrow \op{Sym}^m(C)$ be the projection. The
level $0$ piece of the parabolic structure, on the parabolic bundle
over $\op{Sym}^m(C)$ obtained by descending $E^{\boxtimes m}$ from $C^m$
down to $\op{Sym}^m(C)$, is equal to the subsheaf of sections of
$g_{*}(E{\boxtimes m})$ invariant by the symmetric group action. This
can be seen using the metric interpretation of parabolic structures.

Consider now the commutative (although
not cartesian) diagram 
$$
\xymatrix{
\Ctilde ^m \ar[r] \ar[d] &  C^m \ar[d] \\
\op{Sym}^m(\Ctilde ) \ar[r] & \op{Sym}^m(C) 
}
$$
The sheaf $E^{\boxtimes m}$ is thee direct image of $\specN^{\boxtimes
  m}$ from $\Ctilde^m$.  The permutation group action preserves the
diagram. We can take the direct image of $\specN^{\boxtimes m}$ by the
left vertical and then bottom arrows, and the invariant sections
therein are the direct image of the invariant sections on
$\op{Sym}^m(\Ctilde )$.  This says that our sheaf of invariant
sections is the direct image of $\specN^{(m)}$ from
$\op{Sym}^m(\Ctilde )$ to $\op{Sym}^m(C)$, which means in turn that
$\specN^{(m)}$ is the spectral line bundle. One might have expected
that there could be a correction term by some multiple of the
ramification divisor but this argument shows that that is not the
case.

The parabolic structure is the standard one coming from reflections in
the monodromy of the local system; this happens along the divisor
$\ExY$ in the spectral variety.

The last part follows from the standard properties of symmetric powers
of curves, noting that $P$ is the projective space of the linear
system $|M|$.
\end{proof}

\

\

\noindent
Let $F_B := Rq_{*}(p^{*} E_B^{(m)})$ be the perverse sheaf on $Q$
obtained by the Radon transform. Drinfeld's theorem \ref{drinfeldthm} says
that this is constant on the fibers of the map $\bLL$, so it descends
to a perverse sheaf on $X_{M\otimes \omega _C}$.

\

\begin{theorem}[Deligne cf \cite{LaumonDuke}]
\label{dcalc}
The rank of this sheaf at a general point is $2^{3g-3}$.
\end{theorem}

\

Laumon states in \cite[Remarque 5.5.2]{LaumonDuke} that this was
communicated by Deligne.

\

\

We would like to approximate the Dolbeault calculation of a Higgs
bundle $F$ on $Q$ corresponding to the perverse sheaf $F_B$. We know
that this exists by applying Mochizuki's theory, and we will 
calculate it over a Zariski open subset.

Let 
$$
\xymatrix{
\widetilde{I} \ar[r] \ar[d] &  \widetilde{P} \ar[r] \ar[d] & \op{Sym}^m(\Ctilde ) \ar[d] \\
I \ar[r] &  P \ar[r] & \op{Sym}^m(C)
}
$$
be the cartesian diagram of pullbacks of the spectral variety
of $E^{(m)}$. Thus, 
pulling back the statement from Lemma \ref{spectralEm},
$\widetilde{P}$ is the spectral variety of  $E^{(m)}|_P$ and 
$\widetilde{I}$ is the spectral variety of $p^{*} E^{(m)}$.

\

\

\begin{lemma} \cite{Mumford-Prym}
\label{wtpirr}
The spectral variety $\widetilde{P}$ is irreducible.
\end{lemma}
\begin{proof}
From Lemma \ref{spectralEm}, $\widetilde{P}$ is the subvariety
$\op{Sym}^m(\Ctilde )_M \subset \op{Sym}^m(\Ctilde )$ consisting of
divisors whose associated line bundle has norm down to $C$ equal to
$M$. It is a projective bundle over a translate of the kernel of
$\op{Pic}^m(\Ctilde ) \rightarrow \op{Pic}^m(C)$. We need to show that
this kernel is irreducible.

The map between real tori viewed in terms of the 
exponential exact sequence as
$$
\frac{H^1(\Ctilde , \Oo )}{H^1(\Ctilde , \zz )} 
\rightarrow 
\frac{H^1(C , \Oo )}{H^1(C , \zz )} 
$$
may be written, using Poincaré duality, in the form
$$
\frac{H_1(\Ctilde , \rr )}{H_1(\Ctilde , \zz )} 
\rightarrow 
\frac{H_1(C , \rr )}{H_1(C , \zz )} .
$$
The map $H_1(\Ctilde , \zz ) \rightarrow H_1(C , \zz )$ is surjective,
since $\Ctilde / C$ is a double cover and the ramification set is
nonempty. An element of the kernel is a point in $H_1(\Ctilde , \rr )$
that maps to $H_1(C , \zz )$, and by the surjectivity it can be
modified by an element of $H_1(\Ctilde , \zz )$ so that it maps to
zero in $H_1(C , \rr )$; thus we have a point in the kernel of
$H_1(\Ctilde , \rr )\rightarrow H_1(C , \rr )$ and that covers the
connected component of the identity in the kernel we are looking
at. This shows that our kernel is connected.
\end{proof}

\

\

\noindent
Let $\tilde{q}: \widetilde{I}\dashrightarrow Q$ be the composition of
the map $q$ with the covering.  The holomorphic $L^2$ Dolbeault
complex of $p^{*} E^{(m)}$ relative to the map $q$ has a cokernel
sheaf in top degree
$$
\Dol(I\stackrel{q}{\rightarrow} Q, p^{*} E^{(m)})\rightarrow \Gg .
$$
We can write what this is, away from the parabolic divisors. The spectral
line bundle on $\widetilde{I}$  is 
$$
\specN_{\widetilde{I}}:=\tilde{p}^{*} (\specN^{(m)}|_{\widetilde{P}}.
$$
There is a tautological $1$-form denoted $\alpha_{\widetilde{I}}$. 

The relative dimension of $I/Q$ is one less than the dimension of $P$ or $Q$
since $I$ is a family of hyperplanes; it is $m-g$. Consider the sequence 
$$
\specN_{\widetilde{I}} \otimes \Omega ^{m-g-1}_{\widetilde{I} / Q}
\stackrel{\wedge \alpha_{\widetilde{I}}}{\longrightarrow}
\specN_{\widetilde{I}} \otimes \Omega ^{m-g}_{\widetilde{I} / Q}
\rightarrow 
\widetilde{\Gg}
\rightarrow 0
$$
where $\widetilde{\Gg}$ is defined to be the cokernel. Let $\Gamma
\subset \widetilde{I}$ be the support of $\widetilde{\Gg}$.

This is the relative critical locus, and 
one may alternatively say that $\Gamma$ is the subset of points on which 
the projection $\alpha_{\widetilde{I}}^{\rm rel}$ of the tautological $1$-form
into a section of  $\Omega ^{1}_{\widetilde{I} / Q}$ vanishes.

\

\

\begin{lemma}
\label{projectionG}
Over an appropriate open subset, $\Gg$ is the direct image from
$\widetilde{I}$ down to $I$ of $\widetilde{\Gg}$.
\end{lemma}
\begin{proof}
The relative Dolbeault complex $\Dol(I\stackrel{q}{\rightarrow} Q,
p^{*} E^{(m)})$ is the pushforward of the complex $[\cdots \rightarrow
  \specN_{\widetilde{I}} \otimes \Omega ^{m-g-1}_{\widetilde{I} / Q}
  \rightarrow \specN_{\widetilde{I}} \otimes \Omega
  ^{m-g}_{\widetilde{I} / Q}]$ so over the open subset where
$\widetilde{I} / I$ is finite and flat, the cokernel $\Gg$ of the
pushforward is the pushforward of the cokernel.
\end{proof}

\

\noindent
The main result we need is the following, whose proof is deferred
until after the statements of some lemmas---the lemmas in turn being
proven later too.

\

\

\begin{theorem}
\label{main1}
The subset $\Gamma$ has a unique irreducible component $\Gamma ^{\rm main}$ that
surjects onto $P$. It contains
a dense Zariski open subset that maps isomorphically to an open subset
of $\widetilde{P}$.
This irreducible component $\Gamma ^{\rm main}$ 
maps to $Q$ by a generically finite map of 
degree $2^{3g-3}$. 
\end{theorem}

\

We are going to define a rational map $v:\widetilde{P}\dashrightarrow Q$. 
Suppose $\tilde{z} = \tilde{z}_{1} + \cdots + \tilde{z}_{m}
\in \widetilde{P}$. Use this to define a line bundle 
$$
V(\tilde{z}):= \pi ^{*} (\omega _C) \otimes \Oo _{\Ctilde
}(\tilde{z}_1 + \cdots + \tilde{z}_m)
$$
on $\Ctilde $, and set $\Ll _{\tilde{z}}:= \pi _{*} (V(\tilde{z}))$.

One calculates ${\rm det}(\Ll _{\tilde{z}}) = M\otimes \omega
_C$. This uses the condition that $\tilde{z} \in \widetilde{P}$ saying
in particular that the image $z\in \op{Sym}^m(C)$ is in the linear
system $|M| = P$, i.e. $\Oo _C(z_1+\ldots + z_m)\cong M$.

Therefore $V(\tilde{z})\in \Prym _{M\otimes \omega _C}$.  The moduli
point of $\Ll _{\tilde{z}}$ in $X_{M\otimes \omega _C}$ is the image
of $V(\tilde{z})$ under the natural projection from the Prym. To
define $v(\tilde{z})$ we need to specify an expression of
$\Ll_{\tilde{z}}$ in an extension.

The fact that $\Ctilde $ is a spectral curve implies that there is a natural isomorphism
\begin{equation}
\label{imomega}
\pi _{*} (\pi ^{*} (\omega _C)) \cong \omega _C \oplus \Oo _C,
\end{equation}
in particular we get a morphism $\omega _C \rightarrow
\Ll_{\tilde{z}}$.  For general $\tilde{z}$ (the condition being that
there aren't opposite pairs of points), this is a subbundle. Because
of the determinant calculation we obtain an exact sequence
$$
0\rightarrow \omega _C \rightarrow \Ll _{\tilde{z}} \rightarrow M
\rightarrow 0
$$
and hence an extension $v(\tilde{z}) \in H^1(M^{*} \otimes \omega _C)$
(well-defined up to scalars). This defines our map $v:
\widetilde{P}\dashrightarrow Q$.

\

\begin{remark}
\label{commutativity}
It follows from the above discussion that the diagram of rational maps
$$
\xymatrix{
\widetilde{P} \ar[r] \ar[d] &  Q \ar[d] \\
\Prym _{M\otimes \omega _C} \ar[r] &  X_{M\otimes \omega _C}
}
$$
is commutative, so we get a map $\widetilde{P}
\dashrightarrow \Prym _{M\otimes \omega _C} \times _ {X_{M\otimes \omega _C}} Q$.
\end{remark}

\

\

\begin{lemma}
\label{inhyperplane}
If $\tilde{z}\in \widetilde{P}$ maps to $z\in P$ and to $v(\tilde{z})$
by the above construction then the point $(z,v(\tilde{z}))$ is in
$I$. In particular, the point $(\tilde{z},v(\tilde{z}))$ defines a
point of $\tilde{I}$ so we get a (rationally defined) section $\sigma
: \widetilde{P} \dashrightarrow \widetilde{I}$.
\end{lemma}

\

\

\begin{lemma}
\label{degreecalc}
The degree of $v:\widetilde{P}\dashrightarrow Q$ is $2^{3g-3}$ and
indeed the fiber over a general point $\xi \in Q$ is naturally
identified with the inverse image of $\Ll _{\xi}\in X_{M\otimes \omega
  _C}$ in the Prym variety that is the fiber of the Hitchin map
corresponding to our given spectral curve $\Ctilde /C$.
\end{lemma}

\

\

\begin{lemma}
\label{siggam}
The section $\sigma$ maps $\widetilde{P}$ into $\Gamma$. 
\end{lemma}

\

\

\begin{proof}[Proof of Theorem \ref{main1}]
We claim that over a general point of $\widetilde{P}$ there is exactly
one point of $\Gamma$.  Suppose $\tilde{z}\in \widetilde{P}$ is a
general point mapping to $z\in P$.  Consider the subspace
$\widetilde{V}\subset T_{\tilde{z}}\widetilde{P}$ on which the
tautologial $1$-form $\alpha_{\widetilde{P}}$ vanishes. For
$\tilde{z}$ general, this maps to a codimension one subspace $V\subset
T_zP$.  There will be a unique hyperplane of the projective space $P$
containing $z$ and such that the tangent space of the hyperplane at
$z$ contains $V$.  This hyperplane represents a point of $I$, and its
pair with $\tilde{z}$ defines a point of $\widetilde{I}$. Assuming
that $\widetilde{P}/P$ is \'{e}tale at $\tilde{z}$, the condition on the
tangent spaces is equivalent to the vanishing of the tautological
vertical one-form $\alpha_{\widetilde{I}}^{\rm rel}$ at the
point. This shows the claim.

That implies that the section $\sigma$ provided by Lemmas
\ref{inhyperplane} and \ref{siggam} is an isomorphism of
$\widetilde{P}$ (which is irreducible by Lemma \ref{wtpirr}) onto the
irreducible component $\Gamma ^{\rm main} \subset \Gamma$. The general
point of this component is therefore finite over $Q$ since they have
the same dimension. It has the required degree by Lemma
\ref{degreecalc}.  This gives the statement of the theorem.
\end{proof}

\

\

For the proofs of the lemmas we begin with the following observation
using $m\gg 0$.

\

\begin{remark}
\label{ichar}
Suppose $\xi \in Q$ and $\Ll _{\xi}$ is stable. Then $H^1(\Ll
_{\xi})=0$, indeed an element of $H^0(\Ll _{\xi}^{\vee} \otimes
\omega_C)$ would be a map $\Ll_{\xi}\rightarrow \omega _C$
contradicting stability. There is an exact sequence
$$
0\rightarrow H^0(\omega _C) \rightarrow H^0(\Ll _{\xi}) 
\rightarrow H^0(M) \rightarrow H^1(\omega _C) \rightarrow 0.
$$
Thus the image of $H^0(\Ll _{\xi})$ is a codimension $1$ subspace of
$H^0(M)$ corresponding to a hyperplane in $P$. This hyperplane is
$I_{\xi} = p(q^{-1}\{ \xi \})= \xi ^{\perp}$.
\end{remark}

\

\begin{proof}[Proof of Lemma \ref{inhyperplane}]
Let $\xi = v(\tilde{z})$. We have $\Ll _{\xi} \cong
\pi_{+}(V(\tilde{z}))$.  Recall that $\pi ^{*} (\omega _C) \subset
V(\tilde{z})$, so the formula \eqref{imomega} gives a map
$$
\omega _C \oplus \Oo _C = \pi _{*}(\pi ^{*} (\omega _C) )
\hookrightarrow \Ll _{\xi}.
$$
This gives a map $\Oo _C \rightarrow \Ll _{\xi}$ not factoring
through $\omega _C$, so it projects to a nonzero section in $H^0(M)$
and hence gives a point of $P$.  The cokernel of $\omega _C \oplus
\Oo_C\rightarrow \Ll _{\xi}$ is the sheaf $\pi
_{*}(V(\tilde{z})_{\tilde{z}}$ which is supported on $z$. For the
proof it suffices to assume that $\tilde{z}$ is general, and in this
case at least, the zero scheme of the section $\Oo _C \rightarrow M$
is $z$. The characterization of Remark \ref{ichar} tells us that $z\in
\xi ^{\perp}$, so $(z,\xi )=(z,v(\tilde{z}))$ is in $I$. 
\end{proof}

\

\begin{proof}[Proof of Lemma \ref{degreecalc}]
In the space of Higgs bundles whose determinant is $M\otimes \omega
_C$, let $\Prym$ denote the fiber of the Hitchin fibration over the
point corresponding to the spectral curve $\Ctilde /C$. Let
$$
f : \Prym \dashrightarrow X_{M\otimes \omega _C}
$$
be the rationally defined map of forgetting the Higgs field. The points of 
$\Prym$ are line bundles $V$ on $\Ctilde $ such that the determinant
of $\pi _{*}(V)$ is $M\otimes \omega _C$, and $f(V) = \pi _{*}(V)$.

Suppose $\xi \in Q$ is a general point. As is well-known (cf eg Lemma
\ref{deg8c1}), $f^{-1}(\Ll _{\xi})$ is a finite set with
$2^{3g-3}$ elements.

On the other hand, $v^{-1}(\xi )\subset\widetilde{P} $ 
is the subset of points
$\tilde{z} \in \widetilde{P}$ such that the resulting extension
$v(\tilde{z})$, displaying as 
$$
0 \rightarrow \omega _C \rightarrow \Ll _{\tilde{z}} \rightarrow M
\rightarrow 0,
$$
is equal to the extension $\xi$. 

We are going to establish an isomorphism between these two sets, from
which it follows that $v^{-1}(\xi )\subset\widetilde{P} $ is a finite
set with $2^{3g-3}$ elements.

Suppose given $\tilde{z} \in v^{-1}(\xi )$. Recall that $\Ll _{\tilde{z}}
:=\pi _{*}(V(\tilde{z}))$ and the condition
$\tilde{z} \in v^{-1}(\xi )$ says that 
this bundle is isomorphic to $\Ll _{\xi}$. Therefore, 
$V(\tilde{z}) \in f^{-1}(\Ll _{\xi})$. 

Suppose given $V\in f^{-1}(\Ll _{\xi})$. Then the map
$\omega _C \rightarrow \Ll _{\xi}$
gives a nonzero map $\pi ^{*}(\omega _C) \rightarrow V$, so there is a
divisor $\tilde{z}$ on
$\Ctilde $ such that $V \cong \pi ^{*}(\omega _C)\otimes
\Oo _{\Ctilde }(\tilde{z})$,
and the extension $v(\tilde{z})$ is equal to $\xi$. This gives a point
$\tilde{z} \in v^{-1}(\xi )$.

Let us note that these two constructions are inverses.  In the
previous paragraph, recall that by definition $V(\tilde{z})$ is the
line bundle $ \pi ^{*}(\omega _C)\otimes \Oo _{\Ctilde }(\tilde{z})$
that is isomorphic to $V$, so the composition
$$
V\mapsto \tilde{z} \mapsto V(\tilde{z})
$$
is the identity. For the composition in the other direction, notice
that if given $\tilde{z}\in v^{-1}(\Ll _{\xi})$ then since $\Ll
_{\xi}$ is stable (that is the case for a general $\xi$) the
identification $\Ll _{\xi} \cong \pi _{*} (V(\tilde{z}))$ is unique
and gives rise to a uniquely defined (up to a scalar) map $\pi
^{*}(\omega _C) \rightarrow V(\tilde{z})$, which by definition is our
given one. When we make the construction of the previous paragraph we
get back to the given point $\tilde{z}$.  This gives the required
isomorphism.
\end{proof}

\

\

For the proof of Lemma \ref{siggam} we are going to look at a general
point $\tilde{z}\in \widetilde{P}$. In particular, the $\tilde{z}_i$
are distinct, not ramification points of $\Ctilde /C$, and there are
no opposite pairs under the covering involution $\tau$ of $\Ctilde /C$.

Fix $\xi = v(\tilde{z})$.  We would like to consider a path
$\tilde{z}(t)$ with the given point as $\tilde{z}(0)$, and such that
$z(t)=\pi (\tilde{z}(t))\in P$ stay in the hyperplane $I_{\xi}\subset
P$ corresponding to $\xi$. By Remark \ref{ichar}, we obtain general
such paths of points $x(t)$ by considering a path $\gamma(t)\in
H^0(\Ll_{\xi})$. If we suppose $\gamma(0)$ is the canonical section
corresponding to $\tilde{z}$, then as we deform $\gamma(0)$ we will get a
deformation of the points $\tilde{z}(t)$. This projects to a general
tangent vector of $I_{\xi}$ based at the original point, and because
of our genericity hypothesis $\widetilde{P}\rightarrow P$ is \'{e}tale
near our points, so we get in this way the required general tangent
vector to the fiber $q_{\widetilde{I}}^{-1}(\xi )$.

Let $\lambda$ be the tautological form on $\Ctilde $, and let
$\lambda_V$ be $\lambda$ viewed as a section of $\pi ^{*}(\omega _C)$
and then in turn viewed as a section of $V:= V(\tilde{z}) = \pi
^{*}(\omega _C)(\tilde{v})$.

We have $\Ll _{\xi}= \pi _{*} (V)$ and $\lambda_V$, viewed now as a section of
$\Ll _{\xi}$, is the same as the section $\gamma(0)$. 

For any connection $\nabla$ on $V$, holomorphic near the $\tilde{z}_i$
but possibly meromorphic elsewhere,we have
$$
(\nabla \lambda_V)_{\tilde{z}_i} = T_{\tilde{z}_i}(\Ctilde )
\stackrel{\cong}{\rightarrow}
V_{\tilde{z}_i}.
$$
Letting $\epsilon$ denote an infinitesimal value of $t$, we get
$$
\epsilon^{-1} \gamma(\epsilon ) \in V_{\tilde{z}_i}.
$$
The zero $\tilde{z}_i(\epsilon )$ of the section $\gamma(\epsilon )$,
infinitesimally near to 
$\tilde{z}_i=\tilde{z}_i(0)$, is the displacment of
$\tilde{z}_i$ by $\epsilon$ times the
tangent vector 
$$
\vec{w}_i := 
[(\nabla \lambda_V)_{\tilde{z}_i}]^{-1}(\epsilon ^{-1} \gamma(\epsilon )  )
\in T_{\tilde{z}_i}(\Ctilde ).
$$
If we now apply $\lambda$ considered as a $1$-form to this vector, we
get a number.

\

\begin{lemma}
\label{residue}
The number $\lambda(\vec{w}_i)$ is equal to the residue at $\tilde{z}_i$ of
$$
\gamma'(0):= \frac{\gamma(\epsilon ) - \gamma(0)}{\epsilon} \in H^0(\Ll _{\xi})
$$
where $\gamma'(0)$ is viewed as a section of
$V$ and hence as a differential form on $\Ctilde $
with poles at the points of $\tilde{z}$. 
\end{lemma}
\begin{proof}
This calculation may be done in local coordinates. Let $x$ be the
coordinate near the point $\tilde{z}_i$. The bundle $\pi
^{*}(\omega_C)$ is the same as $\omega _{\Ctilde }$ near our point, so
this bundle has frame $(dx)$.  The bundle $V$ locally is sections
having a simple pole at the origin, so it has frame $x^{-1}dx$. The
section $\gamma(0)$ is $a(x)dx$ and $\gamma'(0) = b(x)x^{-1}dx$.  In other
words,
$$
\gamma(\epsilon ) = a(x)dx + \epsilon b(x)x^{-1}dx  + o(\epsilon ).
$$
The zero of this section near the origin is given by 
$$
x(\tilde{z}_i (\epsilon )) = \epsilon b(0) / a(0) + o(\epsilon ).
$$
The derivative of this in $\epsilon$ is 
$$
\vec{w}_i = b(0)/a(0) \frac{\partial}{\partial x}
$$
viewed as a tangent vector at the origin using our coordinate $x$. The
tautological $1$-form is the same as the section $\gamma(0)$ but viewed as
a $1$-form instead of a section of $V$. In our notations it is still
called $a(x)dx$. When we evaluate this on the tangent vector we get
$$
(a(x)dx) (\vec{w}_i) = b(0).
$$
This is exactly the residue of $\gamma'(0)$ at the origin. 
\end{proof}

\

\

\begin{proof}[Proof of Lemma \ref{siggam}]
The tautological $1$-form on $\op{Sym}^m(\Ctilde )$ evaluated at a
tangent vector that is composed of tangent vectors at distinct points
of $\Ctilde $, is the sum of the values of $\lambda$ evaluated on those
vectors. In the situation of the lemma, our family $\gamma(t)$ yields
curves of points $\tilde{z}_i(t)$ based at the $\tilde{z}_i$, whose
first derivatives give tangent vectors $\vec{w}_i$ based at the
$\tilde{z}_i$. The evaluation of the tautological form on this tangent
vector to the fiber $q_{\widetilde{I}}^{-1}(\xi )$ is therefore the
sum of the $\lambda(\vec{w}_i)$. By Lemma \ref{residue}, this is the
sum of the residues of the section $\gamma'(0)$, viewed as a section of $V$
and hence as a meromorphic section of $\pi^{*}(\omega _C)$, at the
$\tilde{z}_i$. Since these comprise all the polar locus of the
section, the residue theorem says that the sum of the residues is
$0$. We have now shown that the evaluation of the tautological form on
a general tangent vector to the fiber of the map
$q_{\widetilde{I}}^{-1}(\xi )$, is zero. This is exactly the condition
for inclusion of our point $(\tilde{z},\xi )$ in $\Gamma$, completing
the proof of the lemma.
\end{proof}

\

\

\subsection{Dolbeault part of a twistor \texorpdfstring{$\srD$}{D}-module}

In order to use the result of Theorem \ref{main1} in conjunction with
Deligne's calculation \ref{dcalc}, we need to delve into some general
theory of twistor $\srD$-modules
\cite{Mochizuki-D1,Mochizuki-D2,Sabbah}.
Suppose $\Ee$ is a pure twistor $\srD$-module on a variety $X$. Then
there is an open subset $U\subset X$ over which $\Ee$ is smooth. We
get a vector bundle $\Ee _U$ on $\aaaa^1 \times U$ with relative
integrable $z$-connection $\nabla$.  For each $z \in
\aaaa^1$ this gives a vector bundle $\Ee^{z} _U$ with
$z$-connection $\nabla^{z}$.

We will call the fiber at $z = 0$ the \emph{\bfseries Dolbeault part}
of $\Ee$. This is a vector bundle $\Ee_{\op{Dol},U}$ with Higgs field
$\varphi := \nabla^0$. In particular, it has a spectral variety
$\Sigma \hookrightarrow T^{\vee} U$ finite and dominant over $U$.

\

\begin{remark} 
There is a notion of extension of $\Ee_{\op{Dol},U}$ to a parabolic
bundle in codimension $1$, that is to say over $X^{\leq 1}:= X -
D^{\rm sing}$ where $D:= X-U$ is the complementary divisor. The
parabolic structure is given by the $V$-filtration construction. We do
not discuss that here, as we are looking at generic constructions over
$U$.
\end{remark}

\

\

\noindent
The Dolbeault complex of $\Ee$ over $U$ is the Dolbeault complex of
the holomorphic Higgs bundle $\left(\Ee_{\op{Dol},U},\varphi\right)$:
$$
\mathsf{Dol}\left(U; \Ee _{Dol,U},\varphi\right)=
\left[
\Ee _{Dol,U} \stackrel{\wedge \varphi}{\longrightarrow} 
\cdots \stackrel{\wedge \varphi}{\longrightarrow} 
\Ee _{Dol,U}  \otimes \Omega ^n _U
\right] .
$$
Let $\Gg (U; \Ee _{Dol,U},\varphi)$ denote the cokernel of the last map.

\

\

\begin{proposition}
Assume that $\Gg(U; \Ee_{Dol,U},\varphi)$ has finite support, and
let
$$
d:= \dim H^0 (\Gg (U; \Ee _{Dol,U},\varphi ))
$$
be the total length of the support.  Let $H^n(X,\Ee)$ be the
cohomology of $X$ with coefficients in the twistor $\srD$-module.
This is a twistor $\srD$-module over a point, that is to say a vector
bundle over $\pp^1$.  If $d\geq \op{rk} H^n(X,\Ee)$, then equality
holds, and the Dolbeault fiber of $H^n(X,\Ee)$ (i.e. the fiber over
$z = 0$) is naturally isomorphic to
the vector space $\mathbb{H}^0 (X,\Gg (U; \Ee_{Dol,U},\varphi ))$.
\end{proposition}
\begin{proof}
We will not spell the details of this here.  The idea is to look at a
fibration of $X$ as a family of curves, and use that to calculate the
cohomology. We can then apply the calculation of \cite{dirim}, being
careful about the difference between the cohomology of a resolution
and the cohomology of the twistor $\srD$-module using the
decomposition theorem.
\end{proof}

\

\

We can now formulate a relative version.  Let $\Ee$ be a pure twistor
$\srD$-module on a smooth projective variety $X$. If $f : X\rightarrow
Y$ is a morphism to a projective variety $Y$, then choose an
open set $U_Y\subset Y$ and an open subset $U\in f^{-1}(U_Y)$ over
which the map is smooth of relative dimension $m$ and $\Ee$ is smooth.
We can form the relative Dolbeault
complex
$$
\mathsf{Dol}(U/U_Y; \Ee _{Dol,U},\varphi)=
\left[
\Ee _{Dol,U} \stackrel{\wedge \varphi _{X/Y}}{\longrightarrow} 
\cdots \stackrel{\wedge \varphi_{X/Y}}{\longrightarrow} 
\Ee _{Dol,U}  \otimes \Omega ^m _{U/U_Y}
\right],
$$
where $\varphi_{X/Y}$ denotes the projection of $\varphi$ to a
relative Higgs field along the fibers of $f$.  Let $\Gg (U/U_Y; \Ee
_{Dol,U},\varphi )$ be the cokernel of the last map.

\

\

\begin{proposition}
\label{spectral-gen-rel}
Suppose $\Gg (U/U_Y; \Ee _{Dol,U},\varphi )$ has support that is
finite over $U_Y$, and let
$$
d:= \op{rk}  f_{*} (\Gg (U/U_Y; \Ee _{Dol,U},\varphi ))
$$
be the relative length. 
Let $\Ff := R^mf_{*} (\Ee )$ be the higher direct image.  This is a
twistor $\srD$-module over $Y$. Assume that the open set $U_Y$ has
been chosen so that $\Ff$ is smooth on $U_Y$.  If $d\geq \op{rk}
(\Ff)$, then equality holds, and
$$
\Ff_{Dol,U_Y} \cong f_{*} (\Gg (U/U_Y; \Ee _{Dol,U},\varphi )).
$$
Furthermore, the decomposition of $\Ff _{Dol,U_Y}$ over $y\in U_Y$
into a direct sum of pieces indexed by the points in the support of
$\Gg (U/U_Y; \Ee _{Dol,U},\varphi )$ lying over $y$, is the spectral
decomposition of the Higgs bundle $(\Ff _{Dol,U_Y},\phi )$.
\end{proposition}

\

\noindent
The Higgs field on $\Ff _{Dol,U_Y}$ is determined in the same way as in 
\cite{dirim}.

\

\noindent
We can now apply this to the Drinfeld situation.

\

\begin{proof}[Proof of Theorem \ref{rad}]
Comparing Deligne's calculation in Theorem \ref{dcalc} with Theorem
\ref{main1}, the dimension of the cohomology on each geneal fiber of
the map $q$ is equal to the number of points in the support of the
upstairs cokernel sheaf $\widetilde{\Gg}$.  It follows from
Proposition \ref{spectral-gen-rel} that the length of
$\widetilde{\Gg}$ at each of these points is $1$, and then that the
spectral variety for the direct image Higgs bundle is the support of
$\widetilde{\Gg}$, which is to say $\Gamma^{\rm main}$.

To complete the proof of Theorem \ref{rad}, we need to show that the
map $\Gamma^{\rm main} \dashrightarrow T^{\vee} Q$ obtained by
interpreting $\Gamma^{\rm main}$ as the critical locus, is the same as
the map $\widetilde{P}\rightarrow T^{\vee} Q$ that is the pullback of
the map from the Hitchin fiber to the cotangent bundle of the moduli
space of bundles $X_{M\otimes \omega_C}$.

Both maps are given, over general points of $\widetilde{P}$, by a
spectral $1$-form. In the case of the map on the relative critical
locus, this is the same as the spectral $1$-form
$\alpha_{\widetilde{I}}$ on $\widetilde{I}$ restricted to $\Gamma^{\rm
  main}$ that is birational to $\widetilde{P}$.
We have a commutative diagram 
$$
\xymatrix{
\widetilde{P} \ar[r] \ar[d] &  \op{Pic}^{m+4g-4}(\Ctilde ) \ar[d] \\
\op{Sym}^m(\Ctilde ) \ar[r] & \op{Pic}^m(\Ctilde )
}
$$
where the right vertical arrow is the isomorphism given by tensoring
with $\pi^{*}(\omega_C)$.  The spectral $1$-form $\lambda$ on
$\Ctilde$ corresponds to a unique $1$-form
$\alpha_{\op{Pic}(\Ctilde)}$ whose pullback to $\op{Sym}^m(\Ctilde )$
is the $1$-form given by summing up the pullbacks of $\lambda$ on
$\Ctilde^m$.

In the case of the map from the Hitchin fiber, the spectral $1$-form
is given by restricting the $1$-form on $\op{Pic}^{m+4g-4}(\Ctilde )$
coming from the spectral $1$-form on $\Ctilde$, to the Prym
variety. The top map of the above diagram factors through the
projection $\widetilde{P} \rightarrow \Prym _{M\otimes \omega _C}$, so
the top and right pullback of $\alpha_{\op{Pic}(\Ctilde)}$ is the
$1$-form corresponding to the Hitchin fiber.

The left and bottom pullback of $\alpha_{\op{Pic}(\Ctilde)}$ is the
spectral form on $\widetilde{P}$ that was used to define the Higgs
bundle over $P$. By definition the spectral form
$\alpha_{\widetilde{I}}$ is the pullback of this form to
$\widetilde{I}$. In turn, that pulls back to the spectral $1$-form on
$\Gamma^{\rm main}$. In other words, the isomorphism
$$
\Gamma ^{\rm main}\rightarrow \widetilde{P}
$$
identifies the two forms. This implies that they provide the same map
to $T^{\vee}Q$ over a Zariski open subset, and therefore they give the
same map whenever they are defined.  This completes the verification that
the two isomorphic spectral varieties sit in the same way inside
$T^{\vee}Q$.
\end{proof}

\

\

\subsection{Uniqueness over the degree \texorpdfstring{$1$}{1} space}

In this subsection, we will prove a uniqueness result for Higgs
bundles over the degree one moduli space $X_1$ that have the blown-up
Prym $Y_1$ as spectral variety. Let $X$ denote $X_1$ (the intersection
of two quadrics in $\pp^5$) and let $Y$ denote $Y_1 =
\op{Bl}_{\Chat}(\Prym_{3})$. Let $\Wob := \Wob_1$ be the wobbly
locus in $X = X_{1}$.

Let $X^{\circ}$ be the complement of the singular locus of $\Wob$, let
$\Wob ^{\circ} = \Wob \cap X^{\circ}$ and let $Y^{\circ}$ be the
inverse image of $X^{\circ}$ in $Y$. Over $X^{\circ}$ the tautological
$1$-form on $Y^{\circ}$ yields a map $Y^{\circ} \rightarrow
T^{\vee}_{X^{\circ}}(\log \Wob )$.

\

\begin{lemma}
Away from a subset of codimension $2$, this map is an embedding. 
\end{lemma}
\begin{proof}
It is an embedding away from the wobbly locus, since $Y$ is isomorphic
there to a fiber of the Hitchin fibration, the total space of which is
the cotangent bundle of $X$. In fact, that holds true over the
complement $Y-\ExY$ of the exceptional divisor, since the other points
of $Y$, that sit over the wobbly locus but do not belong to $\ExY$,
map into $T^{\vee}(X)$. On the wobbly locus, the first question is
whether the map separates tangent directions. The tangent directions
along $\ExY$ map to tangent directions along $\Wob$. The normal
directions will map to nontrivial vectors in the logarithmic cotangent
bundle, as soon as we know that the tautological $1$-form $\alpha$ on
the Prym variety $\Prym$ is nonzero in a general normal direction to
$\Chat \subset \Prym$. This is true for a general normal direction
since $\alpha$ is a linear form on the abelian variety and there are
two normal directions at each point. This shows that the map
$Y^{\circ} \rightarrow T^{\vee}_{X^{\circ}}(\log \Wob )$ separates
tangent vectors away from a codimension $2$ subset.

To show that it is an embedding, recall that above every (general)
point of $\Wob$ there are $2$ points of $\ExY$. We need to show that
these do not get glued together by the map to
$T^{\vee}_{X^{\circ}}(\log \Wob )$. Notice that no other pairs of
points on $\ExY$, one of which is general, can be glued together since
the map from $\ExY$ to $\Wob$ factors through this double cover.

By abuse of notation we will denote the canonical lift of the covering
involution of $\Ctilde$ over $C$, or equivalently the covering
involution of $\Chat$ over $\Cbar$, again by $\tau$.  The involution
$\tau$ acts naturally on $\Prym$, preserving the curve $\Chat \subset
\Prym$ and thus acts on $Y$.  Two points of $\ExY \subset Y$ mapping
to the same point of $\Wob$ correspond to two points of $\Chat$ of the
form $\tilde{a}$ and $\tau \tilde{a}$; together with corresponding
normal directions $\ell \in \pp(N_{\Chat/\Prym,\tilde{a}})$ and $d\tau
(\ell) \in \pp(\pp(N_{\Chat/\Prym,\tau \tilde{a}})$ at these
points. We note that both points go to the zero-section of the residue
map
$$
T^{\vee}_{X^{\circ}}(\log \Wob)|_{\Wob}  \rightarrow \Oo_{\Wob}.
$$
The kernel of the residue map projects to a map to $T^{\vee}
\Wob$. To distinguish the points we would like to show that their
images in $T^{\vee} \Wob$ are different.  The tautological form is
going to send the points into the image of
$$
T^{\vee} \Cbar \times _{\Cbar} \Wob \rightarrow T^{\vee} \Wob .
$$
The condition that $\Ctilde / \Cbar$ is a smooth spectral curve
means that the two points $\tilde{a}$ and $\tau \tilde{a}$ map to
distinct (opposite) points of $T^{\vee} \Cbar$, showing that they map to
distinct points in $T^{\vee} \Wob$.
\end{proof}

\

\

\begin{lemma}
\label{parstruct}
Suppose $\Ee _{\bullet}$ is a vector bundle with a parabolic structure
over $\Wob^{\circ}$ together with a meromorphic Higgs field that is
logarithmic along $\Wob^{\circ}$, such that the spectral variety of
the Higgs field is the image of $Y^{\circ}$ in
$T^{\vee}_{X^{\circ}}(\log \Wob )$. Suppose that it corresponds to a
flat bundle such that the cohomology of the restriction to a generic
Hecke conic has dimension $16$.  Then the parabolic structure is
obtained from a parabolic level $0<\mu \leq 1$ as follows: there is a
line bundle $\Ll$ on $Y$ such that $\Ee_a = \pi_{*}(\Ll )$ for $0\leq
a < \mu$ and $\Ee_a = \pi _{*}(\Ll (\ExY))$ for $\mu \leq a < 1$.
\end{lemma}
\begin{proof}
If there is no parabolic structure along the exceptional divisor
$\ExY$, then each of the $2$ points of $\ExY$ over a point of $\Wob$
will contribute a unipotent block of size $2$ to the monodromy. Let us
count the contributions of these to the cohomology: over the conic
$\pp^1$ we have a local system of rank $8$, with $16$ points (the
intersections of the conic with $\Wob$) on which there are two
unipotent blocks of size $2$. This gives a total contribution of $-32$
to the Euler characteristic $\chi = h^0 - h^1$. On the other hand,
from the rank $8$ local system over $\pp^1$ the contribution of the
Euler characteristic is $8\cdot 2 = 16$. The sum is then $-16$.  This
means that $h^1 \geq 16$. The hypothesis that $h^1=16$ implies that
there can't be any further parabolic structure. This is the case
$\mu = 1$ (that is equivalent to $\mu = 0$ by an elementary
transformation; for simplicity below we prefer calling this $\mu =
1$).

Suppose now that there is a nontrivial parabolic structure defined
along $\ExY$. If it has a single parabolic level different from $0$, then
the calculation is the same as above giving $h^1=16$. If there were
two parabolic levels at each point of $\ExY$ (or a single level distinct
from $0$ with multiplicity $2$) then the $h^1$ would be too big
compared to the hypothesis.  Thus, there is only one parabolic level
$\mu$ different from $0$ or $1$, and one can see (looking at a
transverse section and thinking about parabolic structures on Higgs
bundles over a curve) that the only possibility to create the
parabolic Higgs bundle is the one described at the end of the
statement.
\end{proof}

\

\
\begin{corollary}
\label{parchoices}
The parabolic level in the previous lemma is either $\mu = 1/2$ or
$\mu = 1$.
\end{corollary}
\begin{proof}
The parabolic first Chern class on $Y$ pushes down to a
class on $X$, which has to vanish for a flat bundle.  The formula of
\cite{IyerSimpson} for this parabolic first Chern class on $X$
simplifies to a simple average over the interval $[0,1]$:
$$
\op{ch}_1^{\rm par}(\Ee_{\bullet}) = 
\int _{a=0}^1 \op{ch}_1(\Ee_{a}) .
$$
This should vanish. In view of the formulas for $\Ee _{a}$ in the
lemma, we get
$$
0 = \op{ch}_1^{\rm par}(\Ee_{\bullet}) = 
\mu \op{ch}_1(\pi _{*}(\Ll )) + (1-\mu ) \op{ch}_1(\pi _{*}(\Ll (\ExY) )).
$$
The GRR formula gives
$$
0 = \pi _{*} ({\rm td}_1(Y/X) + \Ll + (1-\mu) \ExY) .
$$
Equivalently, 
$$
\FxY^2\cdot (\op{td}_1(Y/X) + \Ll + (1-\mu) \ExY) =0.
$$
Recall from Proposition \ref{todd1} that
$$
\op{td}(Y/X) = (1- \FxY + 5\FxY^2/ 12)(1 - \ExY/2 +
(\ExY^2+\ExY\cdot \FxY)/9)
$$
so $\op{td}_1(Y/X) = -\FxY -\ExY/2$. 
If $\Ll = \Oo _Y(a\FxY + b\ExY)$ with $a,b\in \zz$ we get 
$$
0 = \FxY^2 ((-\FxY-\ExY/2) + a\FxY + b\ExY + (1-\mu)\ExY) =
(a-1)\FxY^3 + (b+1/2 - \mu ) \FxY^2\cdot \ExY.
$$
Using the calculations of Proposition \ref{intersections1} this gives 
$$
32 (a-1) + 64 (b+ 1/2 - \mu) = 0.
$$
It follows from the condition that $a,b\in \zz$ that $\mu = 1/2$ or
$\mu = 1$. 
\end{proof}

\

\

\begin{remark}
\label{unique-rmk}
The case $\mu = 1/2$ corresponds to the case we have been treating,
in which we found our flat bundle. We claim that the numerical class
of the line bundle $\Ll$ is uniquely determined by the condition of
$\op{ch}_2$ being extremal. 
\end{remark}

\

\

\noindent
Notice that there is a relation on the coefficients $a$ and $b$ of $\FxY$ and $\ExY$
respectively, 
that is fixed by requiring $\op{ch}_1=0$.  In the case $\mu = 1/2$ this
condition from the previous corollary becomes 
$a=1-2b$. This leaves a single
parameter, which we can view as being the coefficient $b$ of $\ExY$.

The uniqueness of the value for which $\op{ch}_2$ is extremal may be
seen by calculating the parabolic Chern class, up to some correction
terms of the kind we have seen in Chapter \ref{chapter-d1}. The
correction terms are local at the non-normal crossings points of the
wobbly divisor, and don't depend on the choice of $\Ll$. The resulting
function of the coefficient of $\ExY$ in the divisor of $\Ll$ is a
strictly concave quadratic function, whose maximum is at an integer
value; in case the reader is interested, the calculation is reproduced
below, but the interesting point is that this fact comes from the
factor $(1/2)\ExY$ in the relative Todd class combined with the parabolic
level $\mu =1/2$. By the Bogomolov-Gieseker inequality, we can't
choose an integral line bundle $\Ll$ such that $\op{ch}_2 >0$. It
follows that the integral value for which $\op{ch}_2=0$ is unique.

\bigskip

\

Here is an approach to the calculation referred to above. The
condition $\op{ch}_1=0$ tells us that $\Ee_{\alpha} $ may be written
as $\pi _{*} \Oo ((1-2b) \FxY + b\ExY)$ for $0\leq \alpha < 1/2$ and $\pi
_{*} \Oo ((1-2b) \FxY + (b+1)\ExY)$ for $1/2 \leq \alpha < 1$.  Call
these two bundles $\Ee (b)$ and $\Ee '(b)$ respectively. At all steps
below we will allow ourselves to ignore any terms that are constant as
functions of $b$. Set
$$
c(b):= H\cdot  \op{ch}_2(\Ee (b)), \;\;\;
\Delta c(b):= H\cdot (\op{ch}_2(\Ee ' (b)) - \op{ch}_2(\Ee (b))).
$$
Using the convention that we ignore terms that are constant in $b$,
the integral formula for the second Chern class becomes much easier in
that it no longer depends on $\op{ch}_1$:
$$
H\cdot \op{ch}_2^{\rm par} = c(b) + \Delta c(b)/2.
$$
From the GRR formula (and dropping terms not depending on $b$) we have: 
$$
\begin{aligned}
c(b) & = \;\; \mbox{ the degree 3 terms of }  \FxY\cdot ( 1-\FxY  + 5\FxY^2/ 12)
(1 -\ExY/2 + (\ExY^2+\ExY\FxY)/9) \\
& \qquad \qquad \qquad \qquad \qquad \qquad \cdot 
\left(1 + ((1-2b) \FxY + b\ExY)  +
((1-2b) \FxY + b\ExY)^2 / 2\right)  \\
& =
\FxY\cdot (((1-2b) \FxY + b\ExY)^2 / 2 +
((1-2b) \FxY + b\ExY)(-\FxY-\ExY/2)) \\
& =
(1-2b)^2 \FxY^3 / 2 + b(1-2b)\FxY^2 \ExY + b^2 \FxY\ExY^2 / 2 - (1-2b)\FxY^3 
- ((1/2)-b +b)\FxY^2 \ExY - b\FxY\ExY^2/2 \\
& =
2b^2 \FxY^3 + (b - 2b^2) \FxY^2 \ExY + (b^2 - b) \FxY\ExY^2 / 2 \\
& = 2b^2 \cdot 32 + (b-2b^2) \cdot 64 + (b^2 - b) \cdot 16 = 48(b-b^2). 
\end{aligned}
$$
Denote by $c'(b)$ the same formula for $\Ee '(b)$, calculated similarly to be
$$
c'(b) =  -48(b+b^2).
$$
One may alternatively check that $c'(b)$ and $c(b+1)$ are the same up to 
terms not depending on $b$---not a tautology, 
this depends on the calculations of Proposition \ref{intersections1}. 
We get
$$
\Delta c(b) = c'(b) -c(b) = -48(b+b^2) - 48(b-b^2) = -96 b.
$$
Now 
$$
\begin{aligned}
H\cdot \op{ch}_2^{\rm par} & = c(b) + \Delta c(b) / 2 \; + \mbox{ terms constant in }b \\
& = 48 (b-b^2) - 48 b \; + \mbox{ terms constant in }b \\
& = - 48 b^2 \; + \mbox{ terms constant in }b .
\end{aligned}
$$
This has its extremum at $b=0$ as claimed. One may note that this implies the spectral line bundle
is numerically $\Oo (\FxY )$, agreeing with the conclusion of Theorem \ref{degonethm}. This completes our
parenthetical calculation.

\bigskip

\

Moving on, in order to prove Corollary \ref{comparison} we need to
rule out the possibility that $\mu = 1$.  This case corresponds to
the situation that our bundle $\Ee $ has no parabolic structure along
$\Wob ^{\circ}$.

\

\begin{corollary}
\label{alphazero}
With trivial parabolic structure, for $\Ll = \Oo _Y(a\FxY + b\ExY)$
and $\Ee = \pi _{*}(\Ll )$ we have
$$
H^2\cdot \op{ch}_1(\Ee ) = 32(a-1) + 64(b-1/2).
$$
If this vanishes then $a=2-2b$. 
\end{corollary}
\begin{proof}
As in the proof of the previous corollary we get
$$
\begin{aligned}
  H^2\cdot \op{ch}_1(\Ee ) & = H^2 \pi _{*} (\op{td}_1(Y/X) + \Ll ) \\
  & =
  \FxY^2 ((a-1)\FxY + (b-1/2)\ExY) \\
  & = 32(a-1) + 64(b-1/2) = 32(a+2b - 2).
\end{aligned}
$$
This vanishes for $a=2-2b$. 
\end{proof}

\

\

Put $m:= b-1$ so $b=m+1$ and $a=-2m$. (Note that this $m$ is unrelated to the one used in the beginning of this section.)

\begin{lemma}
Let $\Ll = \Oo _Y(-2m\FxY+(m+1)\ExY)$ and  $\Ee = \pi _{*}(\Ll )$, then
$$
H\cdot \op{ch}_2(\Ee ) = -48m^2 - 48m - 8 .
$$
This value is $\leq -8$ if $m\in \zz$. 
\end{lemma}
\begin{proof}
We have that $\op{ch}_2(\Ee ) $ is the degree $2$ part of
$$
\begin{aligned}
\pi _{*} 
\left[
  (1-\FxY \right.  & + 5\FxY^2/ 12)(1 -\ExY/2  \\
  & \left.+ (\ExY^2+\ExY\FxY)/9)(1 -2m\FxY+(m+1)\ExY
  + (-2m\FxY+(m+1)\ExY)^2/2)
  \right]
\end{aligned}
$$
which is $\pi _{*}$ of
$$
5\FxY^2 / 12 +(\ExY^2+\ExY\FxY)/9  +
(-2m\FxY+(m+1)\ExY)^2/2 + \ExY\FxY / 2 -
(\FxY + \ExY/2)(-2m\FxY+(m+1)\ExY) .
$$
We get 
$$
\begin{aligned}
H\cdot \op{ch}_2(\Ee )  
& = 
\FxY\cdot 
\left[ 
  5\FxY^2 / 12 +(\ExY^2+\ExY\FxY)/9  +
  (-2m\FxY+(m+1)\ExY)^2/2  \right. \\
  & + \left.  \ExY\FxY / 2 - (\FxY + \ExY/2)(-2m\FxY+(m+1)\ExY)
  \right].
\end{aligned}
$$
This expands to:
$$
\begin{aligned}
\FxY^3 & (5/12 + 2m^2 + 2m) + \\
& + \ExY \FxY^2 (1/9 - 2m(m+1) + 1/2 + m -(m+1)) \\
& + \ExY^2 \FxY (1/9 + (m+1)^2/2 - (m+1)/2)
\end{aligned}
$$
which, in view of Proposition \ref{intersections1}, becomes
$$
\begin{aligned}
32(2m^2 & + 2m + 5/12) 
+ 64(-2m^2 -2m + 1/9-1/2)
+ 32( m^2 / 2 + m/2 + 1/9) \\
& = -48 m^2 - 48m  + 40/3 + 32/3 - 32 \\
& = -48 m^2 -48m - 8.
\end{aligned}
$$
This is
$$
-48m^2 - 48m - 8 = -48(m+1/2)^2 + 4.
$$

The extremal value $m=-1/2$ is not allowed since $m$ is supposed to be
an integer. The extremal values for integers $m$ are at $m=0$ and $m=-1$ and
there the values are $-8$. 
\end{proof}

\

\

\noindent
Let $\widetilde{X}\rightarrow X$ be a resolution of singularities of
the cusps of $\Wob$ in codimension $2$.

\

\

\begin{lemma}
\label{nopar}
Suppose $\Ee$ is a vector bundle over a surface $Z$. Let $b :
\widetilde{Z}\rightarrow Z$ be a birational map from another smooth
surface obtained by blowing up some points.  Suppose
$\widetilde{\Ee}_{\bullet}$ is a parabolic bundle on $\widetilde{Z}$
that is isomorphic to $\Ee$ over an open subset, complement of a
finite collection of points in $Z$, where $b$ is an
isomorphism. Suppose $\op{ch}_1^{\op{par}}(\widetilde{\Ee}_{\bullet})=0$.  Then
$$
\op{ch}_2^{\op{par}}(\widetilde{\Ee}_{\bullet}) \leq 0
$$
with equality only in the case where the parabolic structure is
trivial and $\widetilde{\Ee}_0=b^{*}(\Ee )$.
\end{lemma}
\begin{proof}
The corrections due to the parabolic structure are local. If any of
the corrections were $>0$ then one could fill in such structures an
arbitrary number of times to a stable vector bundle and contradict the
Bogomolov-Gieseker inequality.

If any of the corrections is $=0$ then one can fill that into an
irreducible flat unitary bundle. The parabolic case of the
Donaldson-Uhlenbeck-Yau theorem
\cite{SteerWren,Mochizuki-kh1,Mochizuki-kh2} would imply that the
resulting bundle is flat, but a flat bundle can not have a nontrivial
parabolic structure only on an exceptional divisor, so the parabolic
structure would have to be trivial.
\end{proof}

\

\noindent
The previous lemma immediately implies the following

\

\begin{corollary}
If $\Ee = \pi _{*}(\Ll )$ for $\Ll$ a line bundle on $Y$, then for any
parabolic extension $\widetilde{\Ee} _{\bullet}$ of $\Ee
|_{X^{\circ}}$ across the exceptional divisors in $\widetilde{X}$ such
that the parabolic first Chern class is $0$, then we have
$$
H\cdot \op{ch}_2^{\op{par}}(\widetilde{\Ee} _{\bullet} )\leq -8.
$$
In particular, no such Higgs bundle can correspond to a local system.
\end{corollary}

\

\

\begin{proof}[Proof of Corollary \ref{comparison}]
This corollary rules out the possibility of a parabolic level $\mu=0$,
so by Corollary \ref{parchoices} the level must be $\mu = 1/2$.
The parabolic Higgs bundle therefore has the same structure, by Lemma
\ref{parstruct}, as the parabolic Higgs bundle that we construct. This
proves Corollary \ref{comparison} for the degree $1$ moduli space, up
to the choice of line bundle of degree $0$ on $Y$.  The
Drinfeld-Laumon construction gives a Hecke eigensheaf with the
original rank $2$ local system as eigenvalue. We have also shown the
Hecke eigensheaf property.  On the one hand, this fixes the choice of
the line bundle of degree $0$ on $Y$, and it also shows that our
construction coincides with the Drinfeld-Laumon construction on the
degree $0$ moduli space.
\end{proof}

\

\

\subsection{Tensor description in the degree \texorpdfstring{$1$}{1} case}
\label{chapter-tensor}

The description of the spectral variety given in Theorem \ref{rad} is
used above in the case of $m\gg 0$. However, it turns out that in our
special case of the moduli space $X_1$ of bundles of rank $2$ and odd
degree on a curve of genus $2$, this construction almost leads
directly to a description of the Hecke eigensheaf.

As pointed out to us by Hitchin, see \cite{Hit22}, Atiyah showed in
1955 \cite{Ati55} that the moduli space of odd degree $\pp GL(2)$-bundles
on $C$ is a double covering of $\pp^3 = \op{Sym}^3(\pp^1)$. See
Proposition~\ref{prop:bottommap} and Theorem~\ref{theo:discrwobbly}
below.  Points of $X_1$ are thus in correspondence with unordered
triples of points in the hyperelliptic $\pp^1$, and up to a further
covering, with unordered triples of points of $C$.  This allows us to
use the symmetric exterior tensor product $\Lambda ^{(\boxtimes
  3)}$. Because of the coverings involved, we will first look at the
description over a general line.

Let $\Lambda$ denote the eigenvalue rank $2$ local system on $C$. This
leads to a rank $8$ local system $\Lambda ^{( \boxtimes 3)}$ on the
third symmetric power $\op{Sym}^3(C)$. If $(E,\theta )$ is the rank
$2$ Higgs bundle associated to $\Lambda$, then the Higgs bundle
associated to $\Lambda ^{( \boxtimes 3)}$ is $(E^{(3)}, \theta
^{(3)})$. By Lemma \ref{spectralEm}, its spectral variety is
$\op{Sym}^3(\Ctilde )$.

\

Suppose $\ell \subset X_1$ is a general line, consisting of bundles
$E$ fitting into an exact sequence
$$
0 \rightarrow A \rightarrow E \rightarrow A^{-1}(\pw) \rightarrow 0
$$ for a degree $0$ line bundle $A$. Set $M:= A^{\otimes - 2} (3\pw)$
be the resulting degree $3$ line bundle, and let $|M |$ be the linear
system of sections of $M$. We have $|M | \cong \pp ^1 \hookrightarrow
\op{Sym}^3(C)$. Points of $| M|$ are lines in $H^0(Hom
(A,A^{-1}(\pw)\otimes \omega _C))$ and this space is Serre dual to
$\op{Ext}^1(A^{-1}(\pw),A)$, so $| M|$ can be identified with the
projectivized space of extensions, which is $\ell$.  Recall from
Subsection \ref{trigonal} that the elements of $|M|$ are the fibers of
the trigonal map $C\rightarrow \pp ^1 \cong \ell$ associated to the
line $\ell$.

\

\

\begin{theorem}
\label{tensor-descrip-line}
The Hecke eigensheaf on $X_1$ with eigenvalue $\Lambda$, 
according to Drinfeld's construction or our construction in  
Section \ref{chapter-d1}, restricted to $\ell$ 
becomes isomorphic to $\Lambda ^{( \otimes 3)}|_{\ell}$.
\end{theorem}
\begin{proof}
By Corollary \ref{comparison}, the Hecke eigensheaf on $X_1$ we
constructed in Section \ref{chapter-d1} is the same as that of
Drinfeld's construction. By construction, the spectral variety of the
associated parabolic Higgs bundle $(\mycal{F}_{1,\bullet},\Phi_{1})$
is the covering $Y_1\rightarrow X_1$. The description in Lemma
\ref{three-eight} of the covering $Y_1$ over a line yields the
isomorphism
$$
Y_1 \times _{X_1} \ell \cong \op{Sym}^3(\Ctilde ) \times _{\op{Sym}^3(C)} \ell .
$$
These spaces are isomorphic to the subvariety $\op{Sym}^3_M(\Ctilde )$
of divisors on $\Ctilde$ whose norm to $C$ is in the linear system
$|M|$.

The parabolic structure has levels $0,1/2$ and one verifies that the
parabolic structure for $E^{(3)}$ is the standard one coming from the
ramification points of the spectral covering, as is the case for
$\mycal{F}_{1,\bullet} |_{\ell}$.

The spectral line bundle for $E^{(3)}$ is pulled back from the line
bundle $\specN ^{(3)}$ on $\op{Sym}^3(\Ctilde )$ obtained using the
spectral line bundle $\specN$ on $\Ctilde$.

Let us check first that this has the right degree to get a degree $0$
parabolic bundle of rank $8$ (this will provide a check of Lemma
\ref{spectralEm}). The spectral line bundle $\specN$ has degree $2$ on
$\Ctilde$. Let $\Ctilde^3_{M}$ denote the subvariety in
$\Ctilde^{\times 3}$ consisting of ordered triples
$(\tilde{t}_{1},\tilde{t}_{2},\tilde{t}_{3})$ of points in $\Ctilde$
for which the divisor $t_{1} + t_{2} + t_{3} = \pi(\tilde{t}_{1}) +
\pi(\tilde{t}_{2}) + \pi(\tilde{t}_{3})$ is in the linear system
$|M|$. The map $\Ctilde ^3_M\rightarrow \Ctilde$ given by, say,  the first
projection has degree $8$. Indeed once we fix one of the points
$\tilde{t}_1$ with image $t_1$ in $C$, the remaining divisor $t_2+t_3$
in $C$ is fixed, leading to two choices of ordered triple of points;
there are $4$ lifts of each pair $t_2,t_3$ to
$\tilde{t}_2,\tilde{t}_3$ so the fiber over $\tilde{t}_1$ has $8$
points in all.  Thus for each $i =1,2,3$ the
pullback $\op{pr}_i^{*}(L)$ therefore has degree $16$
on $\Ctilde ^3_M$, and the tensor product of three of these has degree
$48$.

The spectral line bundle $\specN^{(3)}$ on $\op{Sym}^3_M(\Ctilde )$ 
pulls back to this bundle of degree $48$ via the $6:1$ map 
$\Ctilde ^3_M\rightarrow \op{Sym}^3_M(\Ctilde )$, so $\specN ^{(3)}$ has 
degree $8$.

The map $\op{Sym}^3_M(\Ctilde )\rightarrow \ell $ ramifies at two
types of points: there are the points of the divisor
$\op{Sym}^3_M(\Ctilde )\cap \ExY_1$ that are of the form
$\tilde{t}_{1} + \tilde{t}_{2} + \tilde{t}_{3} \in
\op{Sym}^3_M(\Ctilde)$, where $t_i=t_j$ but $\tilde{t}_{i}\neq
\tilde{t}_{j}$ for some $i \neq j$. There are two such points over
each point of $\ell \cap \Wob_1$. Then, there are ramification points
$\tilde{t}_{1} + \tilde{t}_{2} + \tilde{t}_{3}$ such that one of the
$\tilde{t}_i$ is a ramification point of $\Ctilde / C$. For each of
the $4$ ramification points of $\Ctilde / C$, the other two points in
the trigonal fiber are specified and there are $4$ ways of lifting
these to pairs of points in $\Ctilde$, so we get $16$ such
ramification points.  These constitute the movable ramification locus.

The parabolic degree of a bundle created using the standard parabolic
structure at some simple ramification points, and not at others, is
calculated by the same formula as the degree of a usual direct image,
but not counting the ramification points that are used for the
parabolic structure. Thus, we should count $16$ ramification points
instead of $32$, and the required parabolic degree is one-half of this
number, that is to say $16/2=8$. This is indeed the degree of the
bundle $\specN ^{(3)}$ as is to be expected.

The pullback of the spectral line bundle from $Y_1$ under the map
$\op{Sym}^3_M(\Ctilde )\rightarrow Y_1$ also has degree $8$, since the
spectral line bundle is chosen to create a parabolic structure over
$X_1$ that has vanishing $\op{ch}_1$.

To identify the two spectral line bundles, we need to show that these
two line bundles of degree $8$ are the same, when the spectral line
bundle on $Y_1$ is chosen as a function of $\specN$ in the specified
way that will be described next.

We recall from Subsection \ref{sssec:abelianizeHecke} that the
spectral line bundle $\LY _1$ on $Y_1$ for the Hecke eigensheaf is
related to the spectral line bundle $\specN$ on $\Ctilde$ for the
original eigenvalue Higgs bundle in the following way.

First, one has from Lemma \ref{lemma:LN}  the line bundle 
$$
\Lprym = \trans_{\specN(-\pi^{*}(\pw))}^{*}\thetaprym\otimes \thetaprym^{-1}
$$
on the Prym variety $\Prym$ of degree $0$ line bundles on $\Ctilde$
with trivial norm down to $C$.  This gives by translation the line
bundle
$$
\Lprym_{1} =
\trans_{\mathcal{O}_{\Ctilde}(-\pwtilde -\pi^{*}\pw)}^{*}
\Lprym
$$
on the Prym variety $\Prym _3$ of degree $3$ line bundles with norm
$\Oo _C(3\pw )$.  Then
$$
\LY_{1}  = \blo_{1}^{*}\Lprym_{1}\otimes
f_{1}^{*}\mathcal{O}_{X_{1}}(1)
$$
is the spectral line bundle for the Hecke eigensheaf $\mycal{F}_{1,\bullet}$. 

The following lemma, whose proof will be given below, shows the comparison
between the two spectral line bundles on $\op{Sym}^3_M(\Ctilde )$.

\

\

\begin{lemma}
\label{comparison38}
The pullback of $\LY _1$ by the map $\op{Sym}^3_M(\Ctilde ) \rightarrow
Y_1$ is isomorphic to $\specN^{(3)}$.
\end{lemma}

\

To finish the proof of the theorem, we leave it to the reader to check
that the tautological $1$-form for $E^{(3)}$ given in Lemma
\ref{spectralEm} is the same as the restriction of the tautological
form on $Y_1$. Therefore, the parabolic Higgs bundles are isomorphic,
giving the desired isomorphism of local systems.
\end{proof}

\

\

\begin{proof}[Proof of Lemma \ref{comparison38}:]
Express the degree $0$ line bundle as $A=\Oo _C(a+b-2\pw )$. Thus
$M=\Oo _C(2a+2b - 2\pw)$.
A point of $\op{Sym}^3_M(\Ctilde )$ is a divisor $\tilde{y} =
\tilde{y}_1 + \tilde{y}_2 + \tilde{y}_3$ on $\Ctilde$, such that the
image divisor $y_1+y_2+y_3$ on $C$ is in $|M|$. This yields a line bundle
$$
U _3(\tilde{y}) := \pi ^{*}(A)(\tilde{y}_1 + \tilde{y}_2 + \tilde{y}_3)
= \Oo _{\Ctilde}(\tilde{a} + \tilde{a}' + \tilde{b} + \tilde{b}' +
\tilde{y}_1 + \tilde{y}_2 + \tilde{y}_3 - 2\pw ' - 2 \pw '')
$$
on $\Ctilde$, whose norm to $C$
gives the line bundle $A^{\otimes 2}\otimes \Oo _C(y_1 + y_2 + y_3
-2\pw ) \cong \Oo _C(\pw )$, so $U _3(\tilde{t}) $ is a point of
$\Prym _3$. This describes the map $\op{Sym}^3_M(\Ctilde ) \rightarrow
\Prym _3$ that lifts in a unique way to a map to the blow-up $Y_1$.

The line bundle $\Lprym_1$ extends to a line bundle $\Lprym _{1,
  \op{Jac}}$ on the Jacobian $\op{Jac}^3(\Ctilde )$, as may be seen by
its definition. The above description gives a map
$$
\op{Sym}^3(\Ctilde ) \rightarrow \op{Jac}^3(\Ctilde ).
$$
Furthermore, compose with the translation $\op{Jac}^3(C) \rightarrow
\op{Jac}^0(C)$ that relates $\Lprym _1$ with $\Lprym$, to get a map
$$
U _{ \op{Jac}}: \op{Sym}^3(\Ctilde ) \rightarrow \op{Jac}^0(\Ctilde )
$$
defined by 
$$
U _{ \op{Jac}}(\tilde{y}) = \tilde{a} + \tilde{a}' + \tilde{b} + \tilde{b}' +
\tilde{y}_1 + \tilde{y}_2 + \tilde{y}_3 - 4\pw ' - 3 \pw '' .
$$
The pullback of $\Lprym$ by this map is the tensor product of the
values of $\Lprym$ on each of the points translated back to
$\op{Jac}^0(\Ctilde )$. All the points except those of $\tilde{y}$ are
constant, so tensoring with those values leads to just tensoring with
constant lines.

Let $s:\Ctilde ^3 \rightarrow \op{Sym}^3(\Ctilde )$ be the projection,
and let $\specN _0 := j^{*}\Lprym$ denote the pullback of $\Lprym$ to
$\Ctilde$ along the map $j:\Ctilde \rightarrow Jac^0(\Ctilde )$ that
sends $\tilde{t}$ to $\Oo _{\Ctilde }(\tilde{t}-\pw ')$.

We conclude that 
$$
s^{*} U _{ \op{Jac}}^{*}(\Lprym ) \cong \op{pr}_1^{*}(\specN _0)
\otimes \op{pr}_2^{*}(\specN _0) \otimes \op{pr}_3^{*}(\specN _0).
$$
This in turn implies that 
$$
U _{ \op{Jac}}^{*}(\Lprym )  \cong (\specN _0)^{(3)}.
$$
Therefore, the restriction of $\Lprym _1$ to $\op{Sym}^3_M(\Ctilde )$ is
the same as the restriction of $(\specN _0)^{(3)}$.

The restriction of the spectral line bundle $\LY _1$ to
$\op{Sym}^3_M(\Ctilde )$ is thus $(\specN _0)^{(3)} \otimes
\Oo_{\ell}(1)$.  We have $\specN = \specN_0 \otimes \Oo _{\Ctilde}
(\pw ' + \pw '')$.

The line bundle $\Oo _{\Ctilde}(\pw ' + \pw '') ^{(3)}$ on
$\op{Sym}^3(\Ctilde )$ is given by a divisor whose pullback to
$\Ctilde ^3$ is the set of points $(\tilde{t}_1, \tilde{t}_2,
\tilde{t}_3)$ such that one of the coordinates is either $\pw '$ or
$\pw ''$. The divisor in $\op{Sym}^3(\Ctilde )$ is the set of sums
$\tilde{t}_1 + \tilde{t}_2 + \tilde{t}_3$ containing either $\pw '$ or
$\pw ''$. This condition is equivalent to the condition that the sum
$t_1+t_2+t_3$ contains $\pw$. Restrict now to $\op{Sym}^3_M(\Ctilde
)$.  For $A$ hence $M$ general, there is a unique sum that may be
written as $\pw + t_2 + t_3$, so our divisor consists of all the lifts
of these points to $\Ctilde$. That, in turn, is the fiber of
$\op{Sym}^3_M(\Ctilde )\rightarrow \ell = |M|$ over the point $\pw +
t_2 + t_3\in |M|$ that is also described as the image of $\pw$ under
the trigonal map. We have now shown that the restriction of $\Oo
_{\Ctilde}(\pw ' + \pw '') ^{(3)}$ to $\op{Sym}^3_M(\Ctilde )$ is
isomorphic to the pullback of $\Oo _{\ell}(1)$.

\

Now, the restriction of $\LY _1$ to $\op{Sym}^3_M(\Ctilde )$ is 
$$
(\specN _0)^{(3)} \otimes \Oo _{\ell}(1) = 
(\specN _0)^{(3)} \otimes \Oo _{\Ctilde}(\pw ' + \pw '') ^{(3)}
= \specN ^{(3)} .
$$
This completes the proof of Lemma \ref{comparison38}, tying up what was needed
for the proof of Theorem \ref{tensor-descrip-line}.
\end{proof}

\

\

We now sketch how to go from here to a global description over
$X_1$. The proofs are left to the reader.  Let $P\rightarrow X_1$ be
the degree $4$ map whose fiber over $x\in X_1$ is the set of four
points corresponding to four lines through $x$. Thus $P\rightarrow
\op{Jac}(C)$ is a $\pp^1$-bundle.

\

\

\begin{lemma}
Given any pair of two lines $\ell _1,\ell _2$ passing through $x$, we
obtain a point of $C$.  If $\ell _3,\ell _4$ is the opposite pair of
lines (so that altogether these are the four lines through $x$) then
the corresponding point of $C$ is the conjugate by the hyperelliptic
involution $\hi_C$.
\end{lemma}

\

Given a point of $P$ lying over $x\in X_1$ it is one of the lines, so
there are three pairs of lines containing that one; this gives a point
of $\op{Sym}^3(C)$. Choosing a different point over $x$ results in changing
two of these three by the hyperelliptic involution, and altogether the
four lines yield three pairs of pairs of points, hence three pairs of
conjugate points in $C$, hence three points of $\pp^1$. We get a
diagram
$$
\xymatrix{
P \ar[r] \ar[d] & \op{Sym}^3(C) \ar[d] \\
X_1 \ar[r] &  \op{Sym}^3(\pp^1)
}
$$
The hyperelliptic involution $\hi_{C}$ acts on $\op{Sym}^3(C)$ by acting on
all three of these points, so we can factor and obtain a diagram:
$$
\xymatrix{
P \ar[r] \ar[d] & \op{Sym}^3(C)/\hi_{C} \ar[d] \\
X_1 \ar[r] & \op{Sym}^3(\pp^1).
}
$$
Here the vertical maps have degree $4$ and the horizontal maps have
degree $32$.

\

\

\begin{proposition} \label{prop:bottommap}
The bottom map is the map given by squaring coordinates in $\pp^5$,
from $X_1$ to $\pp^3 = \op{Sym}^3(\pp^1)$. There are $6$ planes in $\pp^3$
corresponding to the six Weierstrass points of $C$, and the map
$X_1\rightarrow \pp^3$ has ramification of order $2$ along these. In
particular it maps to a covering $X'\rightarrow \pp^3$ of degree $2$
ramified on these $6$ planes.
\end{proposition}

\

\

\begin{remark}
It would be good to look more closely at the various group actions of
$(\zz / 2)^n$ for $n=4,5,6$ as well as the Heisenberg group. We will
not carry out this analysis here.
\end{remark}

\

\

\begin{theorem} \label{theo:discrwobbly}
The space $X'$ is the moduli space of $\pp GL_2$ bundles of odd
degree. Let $\Delta \subset \pp^3$ be the discriminant. Then the
pullback of $\Delta$ to $X'$ and $X_1$ are the wobbly loci of those
spaces respectively.
\end{theorem}

\

We can form the local system $\Lambda ^{(\boxtimes 3)}$ on $Sym^3(C)$.
There are two ways of descending $\Lambda$ to a local system $\Lambda
_{\pp^1}$ on $\pp^1$ with order two monodromy at the $6$ points. For
each of these we get local systems $\Lambda ^{(\boxtimes 3)}
_{\pp^1}$. The rank $8$ local system $\Lambda ^{(\boxtimes 3)}$ has
singularities on the big diagonal of $\op{Sym}^3(C)$.

The rank $8$ local system $\Lambda ^{(\boxtimes 3)} _{\pp^1}$ has
singularities on the discriminant and the six planes in $\pp^3$, and
its pullbacks to $X'$ and $X_1$ have singularities (generically finite
of order $2$) on the wobbly loci.

\

\

\begin{theorem}
The pullback of our Hecke eigenvector local system of $X_1$ to $P$ is
isomorphic to the pullback of $\Lambda ^{(\boxtimes 3)}$ to $P$. The
pullbacks of $\Lambda ^{(\boxtimes 3)} _{\pp^1}$ to $X'$ or $X_1$ are
isomorphic to our Hecke eigenvector local system.
\end{theorem}

\

\

We note that the discriminant $\Delta$ has a cuspidal locus along the
small diagonal; the pullback will give the cuspidal locus of the
wobbly locus.

\newpage

\addcontentsline{toc}{section}{References}


\begin{thebibliography}{FGOPN23}

\bibitem[AG15]{ArinkinGaitsgory-nilp}
D.~Arinkin and D.~Gaitsgory.
\newblock Singular support of coherent sheaves and the geometric {L}anglands
  conjecture.
\newblock {\em Selecta Math. (N.S.)}, 21(1):1--199, 2015.

\bibitem[AGK{\etalchar{+}}22a]{agkrrv3}
D.~Arinkin, D.~Gaitsgory, D.~Kazhdan, S.~Raskin, N.~Rozenblyum, and
  Y.~Varshavsky.
\newblock Automorphic functions as the trace of {F}robenius, 2022, 2102.07906.

\bibitem[AGK{\etalchar{+}}22b]{agkrrv2}
D.~Arinkin, D.~Gaitsgory, D.~Kazhdan, S.~Raskin, N.~Rozenblyum, and
  Y.~Varshavsky.
\newblock Duality for automorphic sheaves with nilpotent singular support,
  2022, 2012.07665.

\bibitem[AGK{\etalchar{+}}22c]{agkrrv1}
D.~Arinkin, D.~Gaitsgory, D.~Kazhdan, S.~Raskin, N.~Rozenblyum, and
  Y.~Varshavsky.
\newblock The stack of local systems with restricted variation and geometric
  {L}anglands theory with nilpotent singular support, 2022, 2010.01906.

\bibitem[Ari01]{Arinkin}
D.~Arinkin.
\newblock Orthogonality of natural sheaves on moduli stacks of
  {$SL(2)$}-bundles with connections on {${\mathbb P}^1$} minus $4$ points.
\newblock {\em Selecta Mathematica}, 7(2):213, 2001.

\bibitem[Ati55]{Ati55}
M.~F. Atiyah.
\newblock Complex fibre bundles and ruled surfaces.
\newblock {\em Proceedings of the London Mathematical Society}, 3(4):407--434,
  1955.

\bibitem[BB73]{BB}
A.~Bialynicki-Birula.
\newblock Some theorems on actions of algebraic groups.
\newblock {\em Ann. of Math. (2)}, 98:480--497, 1973.

\bibitem[BC22]{Beraldo-Chen}
D.~Beraldo and L.~Chen.
\newblock Automorphic gluing, 2022, 2204.09141.

\bibitem[BD97]{BD-Hitchin}
A.~Beilinson and V.~Drinfeld.
\newblock Quantization of {H}itchin’s integrable system and {H}ecke
  eigensheaves,.
\newblock http://www.math.
  uchicago.edu/$\sim$drinfeld/langlands/Quantization.Hitchin.pdf, 1997.

\bibitem[Bea96]{BeauvilleSurfaces}
A.~Beauville.
\newblock {\em Complex algebraic surfaces}.
\newblock Number~34 in London Mathematical Society Student Texts. Cambridge
  University Press, 1996.

\bibitem[Bea06]{Beauville-genus2and3}
A.~Beauville.
\newblock Vector bundles and theta functions on curves of genus 2 and 3.
\newblock {\em Amer. J. Math.}, 128(3):607--618, 2006.

\bibitem[Ber19]{Beraldo-Whittaker}
D.~Beraldo.
\newblock On the extended {W}hittaker category.
\newblock {\em Selecta Math. (N.S.)}, 25(2):Paper No. 28, 55, 2019.

\bibitem[Ber20]{Beraldo-spectral}
D.~Beraldo.
\newblock The spectral gluing theorem revisited.
\newblock {\em \'{E}pijournal G\'{e}om. Alg\'{e}brique}, 4:Art. 9, 34, 2020.

\bibitem[Ber21]{Beraldo}
D.~Beraldo.
\newblock On the geometric {R}amanujan conjecture.
\newblock {\em arXiv preprint arXiv:2103.17211}, 2021.

\bibitem[BH95]{BodenHu}
H.U. Boden and Y.~Hu.
\newblock Variations of moduli of parabolic bundles.
\newblock {\em Mathematische Annalen}, 302:539--560, 1995.

\bibitem[Bis97a]{BiswasChern}
I.~Biswas.
\newblock Chern classes for parabolic bundles.
\newblock {\em Journal of Mathematics of Kyoto University}, 37(4):597--613,
  1997.

\bibitem[Bis97b]{BiswasOrbifold}
I.~Biswas.
\newblock Parabolic bundles as orbifold bundles.
\newblock {\em Duke Math. J.}, 90(1):305--325, 1997.

\bibitem[BK22]{bk-survey}
A.~Braverman and D.~Kazhdan.
\newblock Automorphic functions on moduli spaces of bundles on curves over
  local fields: a survey, 2022, 2112.08139.

\bibitem[BL04]{bl}
C.~Birkenhake and H.~Lange.
\newblock {\em Complex abelian varieties}, volume 302 of {\em Grundlehren der
  mathematischen Wissenschaften [Fundamental Principles of Mathematical
  Sciences]}.
\newblock Springer-Verlag, Berlin, second edition, 2004.

\bibitem[BNR89]{BNR}
A.~Beauville, M.~S. Narasimhan, and S.~Ramanan.
\newblock Spectral curves and the generalised theta divisor.
\newblock {\em J. Reine Angew. Math.}, 398:169--179, 1989.

\bibitem[Bod91]{Boden}
H.U. Boden.
\newblock Representations of orbifold groups and parabolic bundles.
\newblock {\em Commentarii Mathematici Helvetici}, 66:389--447, 1991.

\bibitem[Bog78]{bogomolov-tensors}
F.~Bogomolov.
\newblock Holomorphic tensors and vector bundles on projective manifolds.
\newblock {\em Izv. Akad. Nauk SSSR Ser. Mat.}, 42(6):1227--1287, 1439, 1978.

\bibitem[Bog94]{bogomolov}
F.~Bogomolov.
\newblock Stable vector bundles on projective surfaces.
\newblock {\em Mat. Sb.}, 185(4):3--26, 1994.

\bibitem[Bor07]{Borne}
N.~Borne.
\newblock Fibr{\'e}s paraboliques et champ des racines.
\newblock {\em International Mathematics Research Notices}, 2007(16), 2007.
\newblock Art. ID rnm049.

\bibitem[Boz22]{Bozec-nilp}
T.~Bozec.
\newblock Irreducible components of the global nilpotent cone.
\newblock {\em Int. Math. Res. Not. IMRN}, 23:19054--19077, 2022.

\bibitem[Bri02]{Brieskorn}
E.~Brieskorn.
\newblock Singularities in the work of {F}riedrich {H}irzebruch.
\newblock {\em Surveys in differential geometry}, 7(1):17--60, 2002.

\bibitem[BV12]{BorneVistoli}
N.~Borne and A.~Vistoli.
\newblock Parabolic sheaves on logarithmic schemes.
\newblock {\em Advances in Mathematics}, 231(3-4):1327--1363, 2012.

\bibitem[BZN18]{BenZviNadlerBetti}
D.~Ben-Zvi and D.~Nadler.
\newblock Betti geometric {L}anglands.
\newblock {\em Algebraic geometry: Salt Lake City 2015}, 97:3--41, 2018.

\bibitem[CW19]{Wentworth}
B.~Collier and R.~Wentworth.
\newblock Conformal limits and the {B}ia\l ynicki-{B}irula stratification of
  the space of {$\lambda$}-connections.
\newblock {\em Adv. Math.}, 350:1193--1225, 2019.

\bibitem[DM96]{ron.eyal-spectral}
R.~Donagi and E.~Markman.
\newblock Spectral covers, algebraically completely integrable, {H}amiltonian
  systems, and moduli of bundles.
\newblock In {\em Integrable systems and quantum groups ({M}ontecatini {T}erme,
  1993)}, volume 1620 of {\em Lecture Notes in Math.}, pages 1--119. Springer,
  Berlin, 1996.

\bibitem[Dol20]{Dolgachev200}
I.~Dolgachev.
\newblock Kummer surfaces: 200 years of study.
\newblock {\em Notices of the American Mathematical Society}, 67(10), 2020.

\bibitem[Don80]{donagi-quadrics}
R.~Donagi.
\newblock Group law on the intersection of two quadrics.
\newblock {\em Ann. Scuola Norm. Sup. Pisa Cl. Sci. (4)}, 7(2):217--239, 1980.

\bibitem[Don95]{ron-spectral}
R.~Donagi.
\newblock Spectral covers.
\newblock In {\em Current topics in complex algebraic geometry ({B}erkeley,
  {CA}, 1992/93)}, volume~28 of {\em Math. Sci. Res. Inst. Publ.}, pages
  65--86. Cambridge Univ. Press, Cambridge, 1995.

\bibitem[DP09]{DP1}
R.~Donagi and T.~Pantev.
\newblock Geometric {L}anglands and non-abelian {H}odge theory.
\newblock In {\em Surveys in differential geometry. {V}ol. {XIII}. {G}eometry,
  analysis, and algebraic geometry: forty years of the {J}ournal of
  {D}ifferential {G}eometry}, volume~13 of {\em Surv. Differ. Geom.}, pages
  85--116. Int. Press, Somerville, MA, 2009.

\bibitem[DP12]{DP2012}
R.~Donagi and T.~Pantev.
\newblock Langlands duality for {H}itchin systems.
\newblock {\em Inventiones mathematicae}, 189(3):653--735, 2012.

\bibitem[DP22]{DonagiPantev}
R.~Donagi and T.~Pantev.
\newblock Parabolic {H}ecke eigensheaves.
\newblock {\em Ast\'{e}risque}, 435:viii+192, 2022.

\bibitem[DPS16]{dirim}
R.~{Donagi}, T.~{Pantev}, and C.~{Simpson}.
\newblock {Direct Images in Non Abelian {H}odge Theory}.
\newblock {\em ArXiv e-prints}, December 2016, 1612.06388.

\bibitem[DR76]{DesaleRamanan}
U.~Desale and S.~Ramanan.
\newblock Classification of vector bundles of rank 2 on hyperelliptic curves.
\newblock {\em Inventiones mathematicae}, 38:161--185, 1976.

\bibitem[Dri80]{Drinfeld-ICM}
V.~Drinfeld.
\newblock Langlands' conjecture for {${\rm GL}(2)$} over functional fields.
\newblock In {\em Proceedings of the {I}nternational {C}ongress of
  {M}athematicians ({H}elsinki, 1978),}, pages 565--574,. ,, 1980.

\bibitem[Dri83]{Drinfeld}
V.~Drinfeld.
\newblock Two-dimensional {$l$}-adic representations of the fundamental group
  of a curve over a finite field and automorphic forms on {$GL(2)$}.
\newblock {\em Amer. J. Math.}, 105(1):85--114, 1983.

\bibitem[Dri84]{Drinfeld-LOMI}
V.~Drinfeld.
\newblock Two-dimensional {$l$}-adic representations of the {G}alois group of a
  global field of characteristic {$p$} and automorphic forms on {${\rm
  GL}(2)$}.
\newblock {\em Zap. Nauchn. Sem. Leningrad. Otdel. Mat. Inst. Steklov. (LOMI)},
  134:138--156, 1984.
\newblock Automorphic functions and number theory, II.

\bibitem[EFDK22]{efk2}
P.~Etingof, E.~Frenkel, and D.~D.~Kazhdan.
\newblock Analytic {L}anglands correspondence for {$PGL_2$} on {$\mathbb{P}^1$}
  with parabolic structures over local fields.
\newblock {\em Geom. Funct. Anal.}, 32(4):725--831, 2022.

\bibitem[EFK21]{efk1}
P.~Etingof, E.~Frenkel, and D.~Kazhdan.
\newblock An analytic version of the {L}anglands correspondence for complex
  curves.
\newblock In {\em Integrability, quantization, and geometry {II}. {Q}uantum
  theories and algebraic geometry}, volume 103.2 of {\em Proc. Sympos. Pure
  Math.}, pages 137--202. Amer. Math. Soc., Providence, RI, [2021] \copyright
  2021.

\bibitem[EFK22]{efk4}
P.~Etingof, E.~Frenkel, and D.~Kazhdan.
\newblock Analytic {L}anglands correspondence for {$\mathbb{P}GL(2)$} on
  {$\mathbb{P}^1$} with parabolic structures over local fields, 2022,
  2106.05243.

\bibitem[EFK23]{efk3}
P.~Etingof, E.~Frenkel, and D.~Kazhdan.
\newblock Hecke operators and analytic {L}anglands correspondence for curves
  over local fields.
\newblock {\em Duke Math. J.}, 172(11):2015--2071, 2023.

\bibitem[FGOPN23]{FGOP}
E.~Franco, P.~Gothen, A.~Oliveira, and A.~Pe{\'o}n-Nieto.
\newblock Narasimhan--{R}amanan branes and wobbly {H}iggs bundles.
\newblock {\em arXiv preprint arXiv:2302.02736}, 2023.

\bibitem[FGV02]{FrenkelGaitsgoryVilonen}
E.~Frenkel, D.~Gaitsgory, and K.~Vilonen.
\newblock On the geometric {L}anglands conjecture.
\newblock {\em Journal of the American Mathematical Society}, 15(2):367--417,
  2002.

\bibitem[FR22]{FaegermanRaskin}
J.~Faergeman and S.~Raskin.
\newblock Non-vanishing of geometric {W}hittaker coefficients for reductive
  groups, 2022, 2207.02955.

\bibitem[Fre14]{frenkelC}
E.~Frenkel.
\newblock Gauge theory and {L}anglands duality.
\newblock talk given at MSRI, Berkeley, in September 2014, available at
  https://youtu.be/NQfeBeRKMiw?t=4190., 2014.

\bibitem[Fæ22]{Faegerman-qt}
J.~Færgeman.
\newblock Quasi-tempered automorphic {$D$}-modules, 2022, 2210.09193.

\bibitem[GAB{\etalchar{+}}24]{GLC}
D.~Gaitsgory, D.~Arinkin, D.~Beraldo, L.~Chen, J.~Faergeman, K.~Lin, S.~Raskin,
  and N.~Rozenblyum.
\newblock Proof of the geometric {L}anglands conjecture, 2024.
\newblock https://people.mpim-bonn.mpg.de/gaitsgde/GLC/.

\bibitem[Gai97]{GaitsgoryThesis}
D.~Gaitsgory.
\newblock {\em Automorphic Sheaves and Eisenstein Series}.
\newblock PhD thesis, Tell Aviv University, 1997.

\bibitem[Gai11]{Gaitsgory-indcoh}
D.~Gaitsgory.
\newblock Ind-coherent sheaves.
\newblock {\em arXiv preprint arXiv:1105.4857}, 2011.

\bibitem[Gai15]{Gaitsgory}
D.~Gaitsgory.
\newblock Outline of the proof of the geometric {L}anglands conjecture for
  {$GL_2$}.
\newblock {\em Ast\'{e}risque}, 370:1--112, 2015.

\bibitem[Gai17]{GaitsgoryBourbaki}
D.~Gaitsgory.
\newblock Progr\`es r\'{e}cents dans la th\'{e}orie de {L}anglands
  g\'{e}om\'{e}trique.
\newblock In {\em S\'{e}minaire Bourbaki. Vol. 2015/2016. Expos\'{e}s
  1104--1119}, volume 390, pages Exp. No. 1109, 139--168. Soci\'{e}t\'{e}
  math\'{e}matique de France, 2017.

\bibitem[GH94]{GH}
P.~Griffiths and J.~Harris.
\newblock {\em Principles of algebraic geometry}.
\newblock Wiley Classics Library. John Wiley \& Sons, Inc., New York, 1994.
\newblock Reprint of the 1978 original.

\bibitem[GM88]{GoldmanMillson}
W.~Goldman and J.~Millson.
\newblock The deformation theory of representations of fundamental groups of
  compact {K}{\"a}hler manifolds.
\newblock {\em Publications Math{\'e}matiques de l'IH{\'E}S}, 67:43--96, 1988.

\bibitem[GnR22]{Gothen-nilp}
P.~Gothen and R.A.~Z\'{u}\ {n}iga Rojas.
\newblock Stratifications on the nilpotent cone of the moduli space of
  {H}itchin pairs.
\newblock {\em Rev. Mat. Complut.}, 35(2):311--321, 2022.

\bibitem[Gol97]{Goldman}
W.~Goldman.
\newblock Ergodic theory on moduli spaces.
\newblock {\em Annals of mathematics}, pages 475--507, 1997.

\bibitem[GR14]{GaitsgoryRozenblyum}
D.~Gaitsgory and N.~Rozenblyum.
\newblock {DG} indschemes.
\newblock {\em Perspectives in representation theory}, 610:139--251, 2014.

\bibitem[Hau21]{Hausel}
T.~Hausel.
\newblock Enhanced mirror symmetry for {L}anglands dual {H}itchin systems.
\newblock {\em arXiv preprint arXiv:2112.09455}, 2021.

\bibitem[Hel21]{Helmann-nilp}
I.~Hellmann.
\newblock The nilpotent cone in the {M}ukai system of rank two and genus two.
\newblock {\em Math. Ann.}, 380(3-4):1687--1711, 2021.

\bibitem[Heu09]{Heu}
V.~Heu.
\newblock Stability of rank 2 vector bundles along isomonodromic deformations.
\newblock {\em Math. Ann.}, 344(2):463--490, 2009.

\bibitem[HH22]{HauselHitchin}
T.~Hausel and N.J. Hitchin.
\newblock Very stable {H}iggs bundles, equivariant multiplicity and mirror
  symmetry.
\newblock {\em Inventiones mathematicae}, 228(2):893--989, 2022.

\bibitem[Hit87a]{Hitchin-selfd}
N.~J. Hitchin.
\newblock The self-duality equations on a {R}iemann surface.
\newblock {\em Proc. London Math. Soc. (3)}, 55(1):59--126, 1987.

\bibitem[Hit87b]{Hitchin-spectral}
N.~J. Hitchin.
\newblock Stable bundles and integrable systems.
\newblock {\em Duke Math. J.}, 54(1):91--114, 1987.

\bibitem[Hit22]{Hit22}
N.~Hitchin.
\newblock Multiplicity algebras for rank $2$ bundles on curves of small genus.
\newblock arXiv:2203.03424, 2022.

\bibitem[HL17]{HeuLoray-Higgs}
V.~Heu and F.~Loray.
\newblock Hitchin {H}amiltonians in genus 2.
\newblock In {\em Analytic and algebraic geometry}, pages 153--172. Hindustan
  Book Agency, New Delhi, 2017.

\bibitem[HL19]{HeuLoray-flat}
V.~Heu and F.~Loray.
\newblock Flat rank two vector bundles on genus two curves.
\newblock {\em Mem. Amer. Math. Soc.}, 259(1247):v+103, 2019.

\bibitem[HT03]{HauselThaddeus}
T.~Hausel and M.~Thaddeus.
\newblock Mirror symmetry, {L}anglands duality, and the {H}itchin system.
\newblock {\em Inventiones mathematicae}, 153:197--229, 2003.

\bibitem[Hud05]{Hudson}
R.W.H.T. Hudson.
\newblock {\em Kummer's Quartic Surface}.
\newblock University Press, 1905.

\bibitem[IS07]{IyerSimpson-dr}
J.~Iyer and C.~Simpson.
\newblock A relation between the parabolic {C}hern characters of the de {R}ham
  bundles.
\newblock {\em Math. Ann.}, 338(2):347--383, 2007.

\bibitem[IS08]{IyerSimpson}
J.~Iyer and C.~Simpson.
\newblock The {C}hern character of a parabolic bundle, and a parabolic
  corollary of {R}eznikov's theorem.
\newblock In {\em Geometry and dynamics of groups and spaces}, volume 265 of
  {\em Progr. Math.}, pages 439--485. Birkh\"auser, Basel, 2008.

\bibitem[Kaw88]{Kawamata-cover}
Y.~Kawamata.
\newblock Crepant blowing-up of {$3$}-dimensional canonical singularities and
  its application to degenerations of surfaces.
\newblock {\em Ann. of Math. (2)}, 127(1):93--163, 1988.

\bibitem[Keu97]{Keum}
J.H. Keum.
\newblock Automorphisms of {J}acobian {K}ummer surfaces.
\newblock {\em Compositio Mathematica}, 107(3):269--288, 1997.

\bibitem[Kle70]{Klein}
F.~Klein.
\newblock Zur theorie der linencomplexe des ersten und zweiten grades.
\newblock {\em Math. Ann.}, 2:198--226, 1870.

\bibitem[Kon93]{Konno}
H.~Konno.
\newblock Construction of the moduli space of stable parabolic {H}iggs bundles
  on a {R}iemann surface.
\newblock {\em Journal of the Mathematical Society of Japan}, 45(2):253--276,
  1993.

\bibitem[KW07]{KapustinWitten}
A.~Kapustin and E.~Witten.
\newblock Electric-magnetic duality and the geometric {L}anglands program.
\newblock {\em Communications in Number Theory and Physics}, 1(1):1--236, 2007.

\bibitem[Laf02]{LafforgueChtoucas}
L.~Lafforgue.
\newblock Chtoucas de {D}rinfeld et correspondance de {L}anglands.
\newblock {\em Inventiones mathematicae}, 147:1--241, 2002.

\bibitem[Lan71]{langlands}
R.~Langlands.
\newblock Automorphic forms on {${\rm GL}(2)$}.
\newblock In {\em Actes du {C}ongr\`es {I}nternational des {M}ath\'{e}maticiens
  ({N}ice, 1970), {T}ome 2}, pages 327--329. Gauthier-Villars \'{E}diteur,
  Paris, 1971.

\bibitem[Lau87]{LaumonDuke}
G.~Laumon.
\newblock Correspondance de {L}anglands g\'{e}om\'{e}trique pour les corps de
  fonctions.
\newblock {\em Duke Math. J.}, 54(2):309--359, 1987.

\bibitem[Lau88]{LaumonVeryStable}
G.~Laumon.
\newblock Un analogue global du c\^{o}ne nilpotent.
\newblock {\em Duke Math. J.}, 57(2):647--671, 1988.

\bibitem[Lau95]{Laumon95}
G.~Laumon.
\newblock Faisceaux automorphes pour {GL(n)}: la premiere construction de
  {D}rinfeld.
\newblock arXiv:alg-geom/9511004, 1995.

\bibitem[Moc06]{Mochizuki-kh1}
T.~Mochizuki.
\newblock Kobayashi-{H}itchin correspondence for tame harmonic bundles and an
  application.
\newblock {\em Ast\'erisque}, 309:viii+117, 2006.

\bibitem[Moc07a]{Mochizuki-D1}
T.~Mochizuki.
\newblock Asymptotic behaviour of tame harmonic bundles and an application to
  pure twistor {$D$}-modules. {I}.
\newblock {\em Mem. Amer. Math. Soc.}, 185(869):xii+324, 2007.

\bibitem[Moc07b]{Mochizuki-D2}
T.~Mochizuki.
\newblock Asymptotic behaviour of tame harmonic bundles and an application to
  pure twistor {$D$}-modules. {II}.
\newblock {\em Mem. Amer. Math. Soc.}, 185(870):xii+565, 2007.

\bibitem[Moc09]{Mochizuki-kh2}
T.~Mochizuki.
\newblock Kobayashi-{H}itchin correspondence for tame harmonic bundles. {II}.
\newblock {\em Geom. Topol.}, 13(1):359--455, 2009.

\bibitem[Mum66]{mumford-equations1}
D.~Mumford.
\newblock On the equations defining abelian varieties. {I}.
\newblock {\em Invent. Math.}, 1:287--354, 1966.

\bibitem[Mum74]{Mumford-Prym}
D.~Mumford.
\newblock Prym varieties. {I}.
\newblock In {\em Contributions to analysis (a collection of papers dedicated
  to {L}ipman {B}ers)}, pages 325--350. Academic Press, New York, 1974.

\bibitem[Mum07a]{mumford-theta1}
D.~Mumford.
\newblock {\em Tata lectures on theta. {I}}.
\newblock Modern Birkh\"{a}user Classics. Birkh\"{a}user Boston, Inc., Boston,
  MA, 2007.
\newblock With the collaboration of C. Musili, M. Nori, E. Previato and M.
  Stillman, Reprint of the 1983 edition.

\bibitem[Mum07b]{mumford-theta3}
D.~Mumford.
\newblock {\em Tata lectures on theta. {III}}.
\newblock Modern Birkh\"{a}user Classics. Birkh\"{a}user Boston, Inc., Boston,
  MA, 2007.
\newblock With collaboration of Madhav Nori and Peter Norman, Reprint of the
  1991 original.

\bibitem[Mum08]{Mumford-Abelian}
D.~Mumford.
\newblock {\em Abelian varieties}, volume~5 of {\em Tata Institute of
  Fundamental Research Studies in Mathematics}.
\newblock Published for the Tata Institute of Fundamental Research, Bombay; by
  Hindustan Book Agency, New Delhi, 2008.
\newblock With appendices by C. P. Ramanujam and Yuri Manin, Corrected reprint
  of the second (1974) edition.

\bibitem[MY92]{MaruyamaYokogawa}
M.~Maruyama and K.~Yokogawa.
\newblock Moduli of parabolic stable sheaves.
\newblock {\em Mathematische Annalen}, 293:77--99, 1992.

\bibitem[New68]{Newstead}
P.E. Newstead.
\newblock Stable bundles of rank {$2$} and odd degree over a curve of genus
  {$2$}.
\newblock {\em Topology}, 7:205--215, 1968.

\bibitem[NR69]{NR}
M.S. Narasimhan and S.~Ramanan.
\newblock Moduli of vector bundles on a compact {R}iemann surface.
\newblock {\em Ann. of Math. (2)}, 89:14--51, 1969.

\bibitem[PN20]{PeonNieto}
A.~Pe{\'o}n-Nieto.
\newblock Wobbly and shaky bundles and resolutions of rational maps.
\newblock {\em arXiv preprint arXiv:2007.13447}, 2020.

\bibitem[PP21]{PalPauly}
S.~Pal and C.~Pauly.
\newblock The wobbly divisors of the moduli space of rank $2$ vector bundles.
\newblock {\em Advances in Geometry}, 21(4):473--482, 2021.

\bibitem[Roz21]{RozenblyumHecke}
N.~Rozenblyum.
\newblock Connections on moduli spaces and infinitesimal {H}ecke modifications,
  2021, 2108.07745.

\bibitem[Sab05]{Sabbah}
C.~Sabbah.
\newblock Polarizable twistor {$\mathscr{D}$}-modules.
\newblock {\em Ast\'{e}risque}, 300:vi+208, 2005.

\bibitem[Saw16]{Sawin}
W.~Sawin.
\newblock What do {H}ecke eigensheaves actually look like?
\newblock MathOverflow, (https://mathoverflow.net/users/18060/will-sawin),
  2016, https://mathoverflow.net/q/236706.
\newblock URL:https://mathoverflow.net/q/236706 (version: 2016-04-20).

\bibitem[Sim92]{SimpsonHiggs}
C.~Simpson.
\newblock Higgs bundles and local systems.
\newblock {\em Inst. Hautes \'Etudes Sci. Publ. Math.}, 75:5--95, 1992.

\bibitem[SW01]{SteerWren}
B.~Steer and A.~Wren.
\newblock The {D}onaldson-{H}itchin-{K}obayashi correspondence for parabolic
  bundles over orbifold surfaces.
\newblock {\em Canad. J. Math.}, 53(6):1309--1339, 2001.

\bibitem[Tah10]{Taher1}
C.H. Taher.
\newblock Calculating the parabolic {C}hern character of a locally abelian
  parabolic bundle.
\newblock {\em Manuscripta Math.}, 132(1-2):169--198, 2010.

\bibitem[Tah13]{Taher2}
C.H. Taher.
\newblock The {C}hern invariants for parabolic bundles at multiple points.
\newblock {\em Adv. Geom.}, 13(3):547--570, 2013.

\bibitem[Tyu64]{Tyurin}
A.N. Tyurin.
\newblock On the classification of 2-dimensional vector bundles over an
  algebraic curve of any genus.
\newblock {\em Izv. Akad. Nauk SSSR Set. Mat}, 28:21--52, 1964.

\bibitem[Val34]{duVal}
P.~Du Val.
\newblock On isolated singularities of surfaces which do not affect the
  conditions of adjunction (part i.).
\newblock In {\em Mathematical Proceedings of the Cambridge Philosophical
  Society}, volume 30 (4), pages 453--459. Cambridge University Press, 1934.

\bibitem[vGP96]{PreviatoVanGeemen}
B.~van Geemen and E.~Previato.
\newblock On the {H}itchin system.
\newblock {\em Duke Math. J.}, 85(3):659--683, 1996.

\bibitem[Wil20]{Wilkin}
G.~Wilkin.
\newblock The reverse {Y}ang-{M}ills-{H}iggs flow in a neighbourhood of a
  critical point.
\newblock {\em J. Differential Geom.}, 115(1):111--174, 2020.

\bibitem[Wit82]{Witten-Morse}
E.~Witten.
\newblock Supersymmetry and {M}orse theory.
\newblock {\em J. Differential Geometry}, 17(4):661--692, 1982.

\bibitem[Zel20]{Zelaci}
H.~Zelaci.
\newblock On very stablity of principal {$G$}-bundles.
\newblock {\em Geometriae Dedicata}, 204(1):165--173, 2020.

\bibitem[Zuc79]{Zucker}
S.~Zucker.
\newblock Hodge theory with degenerating coefficients: {$L^2$} cohomology in
  the {P}oincar{\'e} metric.
\newblock {\em Annals of Mathematics}, pages 415--476, 1979.

\end{thebibliography}

\newcommand{\etalchar}[1]{$^{#1}$}

\

\

\bigskip

\noindent
Ron Donagi, {\sc University of Pennsylvania}, donagi@upenn.edu

\smallskip

\noindent
Tony Pantev, {\sc University of Pennsylvania}, tpantev@math.upenn.edu

\smallskip

\noindent
Carlos Simpson, {\sc CNRS, Universit\'e C\^ote d'Azur, LJAD}, 
carlos.simpson@univ-cotedazur.fr

\end{document}